\documentclass[a4paper,12pt]{amsbook}
\usepackage{amsmath,amsthm,amsfonts,graphicx,psfrag,amssymb,dsfont,mathrsfs,enumitem}
\usepackage{minitoc,color}

\usepackage[colorlinks]{hyperref}
\usepackage[alphabetic,initials]{amsrefs}

\renewcommand{\thesection}{\arabic{chapter}.\arabic{section}}
\renewcommand{\thefigure}{\arabic{chapter}.\arabic{figure}}

\theoremstyle{plain}
\newtheorem{thm}{Theorem}[chapter]
\newtheorem{heorem}[thm]{``Theorem''}
\newtheorem{prop}[thm]{Proposition}
\newtheorem{conj}[thm]{Conjecture}
\newtheorem{cor}[thm]{Corollary}
\newtheorem{lemma}[thm]{Lemma}
\newtheorem*{question}{Question}

\newtheorem*{thmu}{Theorem}
\newtheorem*{heoremu}{``Theorem''}
\newtheorem{assumption}[thm]{Assumption}
\newtheorem*{lemmu}{Lemma}

\theoremstyle{definition}
\newtheorem{example}[thm]{Example}
\newtheorem{exercise}[thm]{Exercise}
\newtheorem{defn}[thm]{Definition}

\newtheorem*{notation}{Notation}

\theoremstyle{remark}
\newtheorem{remark}[thm]{Remark}

\newcommand{\AT}{\operatorname{AT}}
\newcommand{\Aut}{\operatorname{Aut}}
\newcommand{\aut}{\mathfrak{aut}}
\newcommand{\br}{^{\operatorname{br}}}
\newcommand{\CH}{{\operatorname{CH}}}
\newcommand{\CCH}{{\operatorname{CCH}}}
\newcommand{\CR}{\mathcal CR}
\newcommand{\SFT}{{\operatorname{SFT}}}
\newcommand{\RSFT}{{\operatorname{RSFT}}}
\newcommand{\codim}{\operatorname{codim}}
\newcommand{\coker}{\operatorname{coker}}
\newcommand{\contr}{{\operatorname{contr}}}
\newcommand{\cov}{\operatorname{cov}}
\newcommand{\Crit}{\operatorname{Crit}}

\newcommand{\Diff}{\operatorname{Diff}}
\newcommand{\defin}[1]{\textbf{#1}}
\newcommand{\dist}{\operatorname{dist}}
\newcommand{\dR}{{\operatorname{dR}}}
\newcommand{\End}{\operatorname{End}}

\newcommand{\ev}{\operatorname{ev}}

\newcommand{\Fred}{\operatorname{Fred}}

\newcommand{\GT}{{\operatorname{GT}}}
\newcommand{\Hom}{\operatorname{Hom}}

\newcommand{\Id}{{\operatorname{Id}}}

\renewcommand{\Im}{\operatorname{Im}}
\newcommand{\im}{\operatorname{im}}

\newcommand{\ind}{\operatorname{ind}}

\newcommand{\Jref}{J^{\operatorname{ref}}}
\newcommand{\Jfix}{J^{\operatorname{fix}}}
\newcommand{\loc}{{\operatorname{loc}}}
\newcommand{\Lie}{{\mathcal{L}}}

\newcommand{\Lin}{\mathscr L}
\newcommand{\Max}{^{\max}}
\newcommand{\muCZ}{\mu_{\text{CZ}}}
\newcommand{\muspec}{\mu^{\operatorname{spec}}}
\newcommand{\node}{^{\operatorname{nd}}}
\newcommand{\ord}{\operatorname{ord}}
\newcommand{\Order}{{\mathcal O}}
\newcommand{\PD}{\operatorname{PD}}
\newcommand{\pr}{\operatorname{pr}}
\newcommand{\rank}{\operatorname{rank}}
\renewcommand{\Re}{\operatorname{Re}}
\newcommand{\reg}{{\operatorname{reg}}}

\newcommand{\Span}{\operatorname{Span}}
\newcommand{\std}{_{\operatorname{std}}}

\newcommand{\Sym}{\operatorname{Sym}}
\newcommand{\sym}{{\operatorname{sym}}}

\newcommand{\Tref}{\mathbf{T}_{{\operatorname{ref}}}}

\newcommand{\transpose}{{\operatorname{T}}}

\newcommand{\virdim}{\operatorname{vir-dim}}
\newcommand{\Vol}{\operatorname{Vol}}
\newcommand{\vol}{d\operatorname{vol}}
\newcommand{\wind}{\operatorname{wind}}

\newcommand{\GL}{\operatorname{GL}}

\newcommand{\PSL}{\operatorname{PSL}}

\newcommand{\U}{\operatorname{U}}

\newcommand{\SL}{\operatorname{SL}}

\newcommand{\Spp}{\operatorname{Sp}}

\renewcommand{\AA}{{\mathbb A}}
\newcommand{\CC}{{\mathbb C}}
\newcommand{\DD}{{\mathbb D}}
\newcommand{\intDD}{{\mathring{\DD}}}
\newcommand{\intZ}{{\mathring{Z}}}
\newcommand{\FF}{{\mathbb F}}
\newcommand{\HH}{{\mathbb H}}
\newcommand{\NN}{{\mathbb N}}
\newcommand{\QQ}{{\mathbb Q}}
\newcommand{\RR}{{\mathbb R}}
\newcommand{\TT}{{\mathbb T}}
\newcommand{\ZZ}{{\mathbb Z}}
\newcommand{\aA}{{\mathcal A}}
\newcommand{\bB}{{\mathcal B}}

\newcommand{\dD}{{\mathcal D}}
\newcommand{\eE}{{\mathcal E}}

\newcommand{\hH}{{\mathcal H}}
\newcommand{\jJ}{{\mathcal J}}
\newcommand{\kK}{{\mathcal K}}
\newcommand{\lL}{{\mathcal L}}
\newcommand{\mM}{{\mathcal M}}
\newcommand{\nN}{{\mathcal N}}
\newcommand{\oO}{{\mathcal O}}
\newcommand{\pP}{{\mathcal P}}
\newcommand{\tT}{{\mathcal T}}
\newcommand{\rR}{{\mathscr R}}
\newcommand{\uU}{{\mathcal U}}
\newcommand{\vV}{{\mathcal V}}

\newcommand{\Poisson}{\mathfrak{P}}
\newcommand{\Weyl}{\mathfrak{W}}

\newcommand{\1}{\mathds{1}}
\newcommand{\p}{\partial}
\renewcommand{\dbar}{\bar{\partial}}
\newcommand{\Cinftyloc}{C^\infty_{\loc}}

\numberwithin{equation}{chapter}

\definecolor{blue}{rgb}{0,0,1}
\definecolor{red}{rgb}{1,0,0}
\definecolor{green}{rgb}{0,.7,0}

\psfrag{Sigmadot}{$\dot{\Sigma}$}
\psfrag{u}{$u$}
\psfrag{What}{$\widehat{W}$}
\psfrag{positive end}{$((-\epsilon,0] \times M_+,d(e^r \alpha_+))$}
\psfrag{negative end}{$([0,\epsilon) \times M_-,d(e^r \alpha_-))$}
\psfrag{W}{$(W,\omega)$}
\psfrag{collar+}{$((-\epsilon,0] \times M_+,d(e^r \alpha_+))$}
\psfrag{collar-}{$([0,\epsilon) \times M_-,d(e^r \alpha_-))$}
\psfrag{end+}{$([0,\infty) \times M_+,d(e^r \alpha_+))$}
\psfrag{end-}{$((-\infty,0] \times M_-,d(e^r \alpha_-))$}
\psfrag{inftyM}{$\{\infty\} \times M$}
\psfrag{-inftyM}{$\{-\infty\} \times M$}
\psfrag{uk}{$u_k$}
\psfrag{M+}{$(M_+,\xi_+)$}
\psfrag{M-}{$(M_-,\xi_-)$}
\psfrag{v1+}{$v_1^+$}
\psfrag{v0}{$v_0$}
\psfrag{v1-}{$v_1^-$}
\psfrag{v2-}{$v_2^-$}
\psfrag{v3-}{$v_3^-$}
\psfrag{RtimesM+}{$\RR \times M_+$}
\psfrag{RtimesM-}{$\RR \times M_-$}
\psfrag{S1 x}{$S^1 \times$}
\psfrag{(W,dlambda)}{$(W,d\lambda)$}
\psfrag{(V2,xi2)}{$(V_2,\xi_2)$}
\psfrag{(V4,xi3)}{$(V_3,\xi_3)$}

\title{Lectures on Symplectic Field Theory}
\author{Chris Wendl}
\address{Institut f\"ur Mathematik \\ Humboldt-Universit\"at zu Berlin \\ Unter den Linden 6 \\ 10099 Berlin \\ Germany}
\email{wendl@math.hu-berlin.de}
\date{2015--2016}

\begin{document}
\frontmatter

\maketitle


\dominitoc
\tableofcontents

\chapter*{Preface}

This book is a slightly expanded version of the lecture notes I produced
for a two-semester course taught at University College London in 2015--16,
for Ph.D.~students with a background in basic symplectic geometry and
interest in symplectic topology and/or geometric analysis.  I say ``slightly expanded,'' although the 
reader will quickly
notice that most individual chapters contain far more material than can
reasonably fit into a two-hour lecture.  In reality, much of that material
was only sketched or mentioned in passing during lectures, and I ended up
using the notes to discuss everything that I would like to have explained
if I'd had unlimited time.  This includes relatively detailed
discussions of several important technical points (e.g.~the definition of
spectral flow, generic transversality
in symplectizations, the punctured Riemann-Roch formula, finite energy and
asymptotics with arbitrary stable Hamiltonian structures) which are either
incompletely covered by the existing literature or, in my opinion, simply
more difficult to learn from other sources than they should be.
For topics that are on the other hand well covered elsewhere, I have usually
not felt obliged to explain every detail, but have tried always to provide
adequate references.

One of the interesting features of SFT is that its foundations are---at the
time of this writing---not yet complete.  When the original ``propaganda
paper'' \cite{SFT} appeared in 2000, it was widely believed that the technical
details would be filled in within a few years, and several papers introducing
important applications of SFT to contact topology were written under this 
assumption.  Since then, a certain realization has set in that the results
in those papers cannot truly be regarded as ``theorems'' in the sense of 
mathematics, and it has become less socially acceptable to preface statements
of results with caveats of the form, ``this theorem is dependent on the foundations of
SFT''.  At the same time, the need for
a robust perturbation scheme to achieve transversality in SFT spawned
the development of a whole new approach to infinite-dimensional differential 
geometry, the \emph{polyfold} project \cite{Hofer:CDM}, which is intended
for much more general applications but is not yet finished.
Opinions vary among symplectic topologists as to how unsatisfied we should
all be with this state of affairs, and what could be done about it---among
other things, one could make an entire course out of the discussion of such 
issues, but I have not chosen to do that.  My
approach is instead to develop the \emph{classical}\footnote{For the
purposes of this discussion, the word ``classical'' may be defined as
``not involving the words \emph{polyfold}, \emph{virtual} or \emph{Kuranishi}''.}
analysis of pseudoholomorphic
curves in symplectizations and symplectic cobordisms, 
to explain how this would lead to a theory of
algebraic contact invariants if transversality for multiple covers
were not an issue, and then
to use the tools and insights gained from this discussion to prove \emph{rigorous
mathematical theorems} about contact manifolds.  Typically, such theorems
can be regarded informally as consequences of computations in a (not yet
well-defined) theory called SFT, but
in a rigorous sense, they are actually consequences of the methods used in those
computations.  Examples covered in these notes 
include distinguishing tight contact structures on the $3$-torus that are
homotopic but not isomorphic (Lecture~\ref{lec:tight3tori}), and the
nonexistence of symplectic fillings or symplectic cobordisms between certain
pairs of contact manifolds (Lecture~\ref{lec:torsion}).
The choice of applications is of course biased somewhat toward my own research
interests.

\subsection*{Prerequisites}

The stated target audience for the lecture course was
``Ph.D.~students in differential geometry or related fields who are not afraid of analysis''.
More precisely, the notes assume some knowledge of the following topics:
\begin{itemize}
\item Differential geometry: manifolds and vector bundles, differential 
forms and Stokes' theorem, connections, basic familiarity with symplectic manifolds
\item Functional analysis: linear operators on Banach spaces, basics of Sobolev spaces, Fredholm operators
\item Differential topology: smooth mapping degree, intersection numbers, Sard's theorem
\item Algebraic topology: fundamental group, homology and cohomology of manifolds, Poincar\'e duality, first Chern class, homological intersection numbers
\end{itemize}
The following topics are not considered formal prerequisites, but some knowledge of 
them is likely in any case to be helpful to the reader, who may
want to have a good reference for them (as suggested below) within arm's reach:
\begin{itemize}
\item Contact manifolds (e.g.~Geiges \cite{Geiges:book})
\item Differential calculus on Banach spaces and Banach manifolds
(e.g.~these two books by Lang: \cite{Lang:analysis} and \cite{Lang:geometry})
\item Closed pseudoholomorphic curves (e.g.~McDuff-Salamon \cite{McDuffSalamon:Jhol}
or my other book in preparation \cite{Wendl:lecturesV33})
\item Floer homology (e.g.~Salamon \cite{Salamon:Floer} or
Audin-Damian \cite{AudinDamian})
\end{itemize}

\subsection*{Acknowledgements}

I would like to thank the students who sat through the course that gave rise
to these notes, and in particular Alexandru Cioba and Agust\'{i}n Moreno
for their assistance in editing the first several lectures.
My understanding of Taubes's approach to the Riemann-Roch formula
(explained in Lecture~\ref{lec:index}) and its generalization to the punctured
case emerged in part from discussions with Chris Gerig, and I am grateful
also to Tim Perutz for helpful hints about Weitzenb\"ock formulas, and
Patrick Massot for patient discussions of singular integral operators and
elliptic regularity.
Thanks also to Michael Hutchings and Janko Latschev for helping me understand 
the combinatorial factors in Lecture~\ref{lec:H}, to Jo Nelson for
helpful comments on coefficients and orbifold singularities, and to Sam Lisi
and Barney Bramham for advice on the Floer $C_\varepsilon$ space.

\adjustmtc

\chapter*{About the current version}

At the time of posting this on the arXiv, Lectures~\ref{lec:automatic},
\ref{lec:intersections} and \ref{lec:torsion} each consist of messy
handwritten notes that have not yet been typed up, but will eventually
appear in the published version of the book.  The main goal for those lectures is
to carry out some explicit computations of the torsion invariant introduced
at the end of Lecture~\ref{lec:SFT}, and to explain the consequences 
for filling and cobordism obstructions, including for
instance the classic result that overtwistedness implies vanishing contact 
homology and thus obstructs fillability.  In keeping with the spirit of the book,
the theorems about torsion in Lecture~\ref{lec:torsion} will need to be
understood with the usual caveat that they depend on the unfinished foundations
of SFT, but part of the point is also to extract complete and
\emph{rigorous} proofs of the
important consequences regarding symplectic fillings.
Lectures~\ref{lec:automatic} and \ref{lec:intersections} are more
technical in nature, in the spirit of Lectures~\ref{lec:local}
through~\ref{lec:compactness} except that they deal with topics that are
only relevant in low-dimensional settings (and thus significantly
increase the power of the theory in those settings).  Aside from dealing with
topics that are valuable in their own right, they specifically precede
Lecture~\ref{lec:torsion} because they introduce techniques that will be
used in the computations in that lecture.

As far as the rest of the manuscript is concerned, I have tried to
produce something that is relatively well polished, but I admit I have
not tried quite as diligently for that as I do with most of my
research papers.  Trying to produce another one of these lectures every
week while teaching the course was a formidable task, and I had more time
to be careful with it in some weeks than in others.  I have since gone
back and reworked some portions, but not all, so I apologize for any
sloppiness that I may have failed so far to expunge.  All comments and
corrections are welcome,\footnote{especially if those corrections are
received before the book goes to press} 
and may be sent to \texttt{wendl@math.hu-berlin.de}.
Updates on the publication of the book will be posted periodically on
my website at
\begin{center}
\url{https://www.mathematik.hu-berlin.de/~wendl/publications.html#notes}
\end{center}

\adjustmtc

\mainmatter
\pagestyle{myheadings}

\chapter{Introduction}
\label{lec:intro}

\minitoc

\vspace{12pt}

Symplectic field theory is a general framework for defining invariants of
contact manifolds and symplectic cobordisms between them via counts of
``asymptotically cylindrical'' pseudoholomorphic curves.  In this first
lecture, we'll summarize some of the historical background of the subject,
and then sketch the basic algebraic formalism of SFT.

\section{In the beginning, Gromov wrote a paper}
\label{sec:Gromov}

Pseudoholomorphic curves first appeared in symplectic geometry in a
1985 paper of Gromov \cite{Gromov}.  The development was revolutionary for
the field of symplectic topology, but it was not unprecedented: a few
years before this, Donaldson had demonstrated the power of using elliptic
PDEs in geometric contexts to define invariants of smooth $4$-manifolds
(see \cite{DonaldsonKronheimer}).  The PDE that Gromov used was a slight
generalization of one that was already familiar from complex geometry.

Recall that if $M$ is a smooth $2n$-dimensional manifold, an \defin{almost
complex structure} on $M$ is a smooth linear bundle map $J : TM \to TM$
such that $J^2 = -\1$.  This makes the tangent spaces of $M$ into complex
vector spaces and thus induces an orientation on~$M$; the pair $(M,J)$ is
called an \defin{almost complex manifold}.  In this context,
a \defin{Riemann surface} is an almost complex manifold of real dimension~$2$
(hence complex dimension~$1$), and a \defin{pseudoholomorphic curve} 
(also called \defin{$J$-holomorphic}) is a smooth map
$$
u : \Sigma \to M
$$
satisfying the \defin{nonlinear Cauchy-Riemann equation}
\begin{equation}
\label{eqn:nonlinearCR}
Tu \circ j = J \circ Tu,
\end{equation}
where $(\Sigma,j)$ is a Riemann surface and $(M,J)$ is an almost complex
manifold (of arbitrary dimension).  The almost complex structure $J$
is called \defin{integrable} if $M$ is admits the structure of a complex 
manifold such that $J$ is multiplication by~$i$ in holomorphic coordinate
charts.  By a basic theorem of the subject, every almost complex structure
in real dimension two is integrable, hence one can always find local
coordinates $(s,t)$ on neighorhoods in $\Sigma$ such that
$$
j \p_s = \p_t, \qquad j\p_t = -\p_s.
$$
In these coordinates, \eqref{eqn:nonlinearCR} takes the form
$$
\p_s u + J(u) \p_t u = 0.
$$

The fundamental insight of \cite{Gromov} was that solutions to the
equation \eqref{eqn:nonlinearCR} capture information about symplectic
structures on $M$ whenever they are related to $J$ in the following
way.

\begin{defn}
\label{defn:tame}
Suppose $(M,\omega)$ is a symplectic manifold.  An almost complex structure
$J$ on $M$ is said to be \defin{tamed} by $\omega$ if
$$
\omega(X,JX) > 0 \quad \text{ for all $X \in TM$ with $X \ne 0$}.
$$
Additionally, $J$ is \defin{compatible} with $\omega$ if the pairing
$$
g(X,Y) := \omega(X,J Y)
$$
defines a Riemannian metric on~$M$.
\end{defn}

We shall denote by $\jJ(M)$ the space of all smooth almost complex
structures on $M$, with the $\Cinftyloc$-topology, and if $\omega$
is a symplectic form on $M$, let
$$
\jJ_\tau(M,\omega) , \jJ(M,\omega) \subset \jJ(M)
$$
denote the subsets consisting of almost complex structures that are
tamed by or compatible with $\omega$ respectively.  Notice that
$\jJ_\tau(M,\omega)$ is an open subset of $\jJ(M)$, but
$\jJ(M,\omega)$ is not.  A proof of the following may be found in
\cite{Wendl:lecturesV33}*{\S 2.2}, among other places.

\begin{prop}
\label{prop:contractible}
On any symplectic manifold $(M,\omega)$, the spaces $\jJ_\tau(M,\omega)$
and $\jJ(M,\omega)$ are each nonempty and contractible.  \qed
\end{prop}

Tameness implies that the \defin{energy} of a $J$-holomorphic curve
$u : \Sigma \to M$,
$$
E(u) := \int_{\Sigma} u^*\omega,
$$
is always nonnegative, and it is strictly positive unless $u$ is constant.
Notice moreover that if the domain $\Sigma$ is closed, then $E(u)$ depends
only on the cohomology class $[\omega] \in H^2_\dR(M)$ and the homology class
$$
[u] := u_*[\Sigma] \in H_2(M),
$$
so in particular, any family of $J$-holomorphic curves in a fixed homology
class satisfies a uniform energy bound.  This basic observation is one of
the key facts behind Gromov's compactness theorem, which states that
moduli spaces of closed curves in a fixed homology class are compact up to
``nodal'' degenerations.

The most famous application of pseudoholomorphic curves presented in 
\cite{Gromov} is Gromov's \emph{nonsqueezing theorem}, which was the first
known example of an obstruction for embedding symplectic domains that is
subtler than the obvious obstruction defined by volume.  The technology
introduced in \cite{Gromov} also led directly to the development of the
\emph{Gromov-Witten invariants} (see \cites{McDuffSalamon:Jhol,RuanTian,
RuanTian:higherGenus}), which follow the same pattern as Donaldson's
earlier smooth $4$-manifold invariants; they use counts of $J$-holomorphic
curves to define invariants of symplectic manifolds up to symplectic
deformation equivalence.

Here is another sample application from \cite{Gromov}.
We denote by
$$
A \cdot B \in \ZZ
$$
the intersection number between two homology classes $A, B \in H_2(M)$ in
a closed oriented $4$-manifold~$M$.

\begin{thm}
\label{thm:S2xS2}
Suppose $(M,\omega)$ is a closed and connected symplectic $4$-manifold with
the following properties:
\begin{enumerate}
\renewcommand{\labelenumi}{(\roman{enumi})}
\item $(M,\omega)$ does not contain any symplectic submanifold
$S \subset M$ that is diffeomorphic to $S^2$ and satisfies $[S] \cdot [S] = -1$.
\item $(M,\omega)$ contains two symplectic submanifolds $S_1, S_2 \subset M$
which are both diffeomorphic to $S^2$, satisfy
$$
[S_1] \cdot [S_1] = [S_2] \cdot [S_2] = 0,
$$
and have exactly one intersection point with each other, which is transverse
and positive.  
\end{enumerate}
Then $(M,\omega)$ is symplectomorphic to
$(S^2 \times S^2,\sigma_1 \oplus \sigma_2)$, where for $i=1,2$,
the $\sigma_i$ are area forms on $S^2$ satisfying
$$
\int_{S^2} \sigma_i = \langle [\omega] , [S_i] \rangle.
$$
\end{thm}
\begin{proof}[Sketch of the proof]
Since $S_1$ and $S_2$ are both symplectic submanifolds, one can choose a
compatible almost complex structure $J$ on $M$ for which both of them
are the images of embedded $J$-holomorphic curves.  One then considers
the moduli spaces $\mM_1(J)$ and $\mM_2(J)$ of equivalence classes of
$J$-holomorphic spheres homologous to $S_1$ and $S_2$ respectively, where
any two such curves are considered equivalent if one is a reparametrization
of the other (in the present setting this just means they have the same
image).  These spaces are both manifestly nonempty, and one can argue via
Gromov's compactness theorem for $J$-holomorphic curves that both are
compact.  Moreover, an infinte-dimensional version of the implicit function
theorem implies that both are smooth $2$-dimensional manifolds, carrying
canonical orientations, hence both are diffeomorphic to closed surfaces.
Finally, one uses \emph{positivity of intersections} to show that every
curve in $\mM_1(J)$ intersects every curve in $\mM_2(J)$ exactly once, and
this intersection is always transverse and positive; moreover, any two
curves in the same space $\mM_1(J)$ or $\mM_2(J)$ are either identical
or disjoint.  It follows that both moduli spaces are diffeomorphic to $S^2$,
and both consist of smooth families of $J$-holomorphic spheres that
foliate~$M$, hence defining a diffeomorphism
$$
\mM_1(J) \times \mM_2(J) \to M
$$
that sends $(u_1,u_2)$ to the unique point in the intersection
$\im u_1 \cap \im u_2$.  This identifies $M$ with $S^2 \times S^2$ such that
each of the submanifolds $S^2 \times \{*\}$ and $\{*\} \times S^2$ are
symplectic.  The latter observation can be used to determine the symplectic 
form up to deformation, so that by the Moser stability theorem, $\omega$
is determined up to isotopy by its cohomology class $[\omega] \in 
H^2_\dR(S^2 \times S^2)$,
which depends only on the evaluation of $\omega$ on $[S^2 \times \{*\}]$
and $[\{*\} \times S^2] \in H_2(S^2 \times S^2)$.
\end{proof}

For a detailed exposition of the above proof of Theorem~\ref{thm:S2xS2},
see \cite{Wendl:rationalRuled}*{Theorem~E}.

\section{Hamiltonian Floer homology}
\label{sec:Floer}

Throughout the following, we write
$$
S^1 := \RR / \ZZ,
$$
so maps on $S^1$ are the same as $1$-periodic maps on~$\RR$.
One popular version of the \emph{Arnold conjecture} on symplectic fixed 
points can be stated as follows.  Suppose $(M,\omega)$ is a closed symplectic
manifold and $H : S^1 \times M \to \RR$ is a smooth function.  Writing
$H_t := H(t,\cdot) : M \to \RR$, $H$ determines a $1$-periodic time-dependent
Hamiltonian vector field $X_t$ via the relation\footnote{Elsewhere in the
literature, you will sometimes see \eqref{eqn:Hamiltonian} without the
minus sign on the right hand side.  If you want to know why I strongly believe
that the minus sign belongs there, see \cite{Wendl:blogSigns}, but to some
extent this is just a personal opinion.}
\begin{equation}
\label{eqn:Hamiltonian}
\omega(X_t,\cdot) = - dH_t.
\end{equation}

\begin{conj}[Arnold conjecture]
If all $1$-periodic orbits of $X_t$ are nondegenerate, then the number of
these orbits is at least the sum of the Betti numbers of~$M$.
\end{conj}

Here a $1$-periodic orbit $\gamma : S^1 \to M$ of $X_t$ is called 
\defin{nondegenerate} if, denoting the flow of $X_t$ by $\varphi^t$, the
linearized time~$1$ flow
$$
d\varphi^1(\gamma(0)) : T_{\gamma(0)} M \to T_{\gamma(0)} M
$$
does not have $1$ as an eigenvalue.  This can be thought of as a Morse
condition for an action functional on the loop space whose critical points are
periodic orbits; like Morse critical points, nondegenerate periodic orbits
occur in isolation.  To simplify our
lives, let's restrict attention to \emph{contractible} orbits and also
assume that $(M,\omega)$ is \defin{symplectically aspherical}, which means
$$
[\omega]|_{\pi_2(M)} = 0.
$$
Then if $C^\infty_\contr(S^1,M)$ denotes the space of all smoothly contractible
smooth loops in $M$, the \defin{symplectic action functional} can be defined by
$$
\aA_H : C^\infty_\contr(S^1,M) \to \RR : \gamma \mapsto - \int_{\DD}
\bar{\gamma}^*\omega + \int_{S^1} H_t(\gamma(t)) \, dt,
$$
where $\bar{\gamma} : \DD \to M$ is any smooth map on the closed unit disk
$\DD \subset \CC$ satisfying
$$
\bar{\gamma}(e^{2\pi i t}) = \gamma(t),
$$
and the symplectic asphericity condition guarantees that $\aA_H(\gamma)$ does
not depend on the choice of~$\bar{\gamma}$.

\begin{exercise}
\label{EX:firstVariation}
Regarding $C^\infty_\contr(S^1,M)$ as a Fr\'echet manifold with tangent
spaces $T_\gamma C^\infty_\contr(S^1,M) = \Gamma(\gamma^*TM)$, show that the
first variation of the action functional $\aA_H$ is
$$
d\aA_H(\gamma) \eta = \int_{S^1} \left[ \omega(\dot{\gamma},\eta) +
d H_t(\eta) \right] \, dt =
\int_{S^1} \omega(\dot{\gamma} - X_t(\gamma), \eta) \, dt
$$
for $\eta \in \Gamma(\gamma^*TM)$.  In particular, the critical points of
$\aA_H$ are precisely the contractible $1$-periodic orbits of~$X_t$.
\end{exercise}

A few years after Gromov's introduction of pseudoholomorphic curves, Floer
proved the most important cases of the Arnold conjecture by developing a
novel version of infinite-dimensional Morse theory for the functional~$\aA_H$.
This approach mimicked the homological approach to Morse theory which has
since been popularized in books such as \cites{AudinDamian,Schwarz:Morse},
but was apparently only known to experts at the time.  
In \emph{Morse homology}, one considers a smooth Riemannian manifold
$(M,g)$ with a Morse function $f : M \to \RR$, and defines a chain complex
whose generators are the critical points of $f$, graded according to their
Morse index.  If we denote the generator corresponding to a given critical
point $x \in \Crit(f)$ by $\langle x \rangle$, the boundary map on this 
complex is defined by
$$
\p \langle x \rangle = \sum_{\ind(y) = \ind(x) - 1}
\# \left(  \mM(x,y) \big/ \RR \right) \langle y \rangle,
$$
where $\mM(x,y)$ denotes the moduli space of negative gradient flow lines
$u : \RR \to M$, satisfying $\p_s u = -\nabla f (u(s))$, 
$\lim_{s \to -\infty} u(s) = x$ and $\lim_{s \to +\infty} u(s) = y$.
This space admits a natural $\RR$-action by shifting the variable in the
domain, and one can show that for generic choices of $f$ and the metric~$g$, 
$\mM(x,y) / \RR$ is a finite set whenever $\ind(x) - \ind(y) = 1$.
The real magic however is contained in the following statement about the
case $\ind(x) - \ind(y) = 2$:

\begin{prop}
\label{prop:Floer}
For generic choices of $f$ and $g$ and any two critical points
$x, y \in \Crit(f)$ with $\ind(x) - \ind(y) = 2$, $\mM(x,y) / \RR$ is 
homeomorphic to a finite collection of circles and open intervals whose
end points are canonically identified with the finite set
$$
\p \overline{\mM}(x,y) := \bigcup_{\ind(z) = \ind(x)-1} 
\mM(x,z) \times \mM(z,y).
$$
\end{prop}

\begin{figure}
\includegraphics{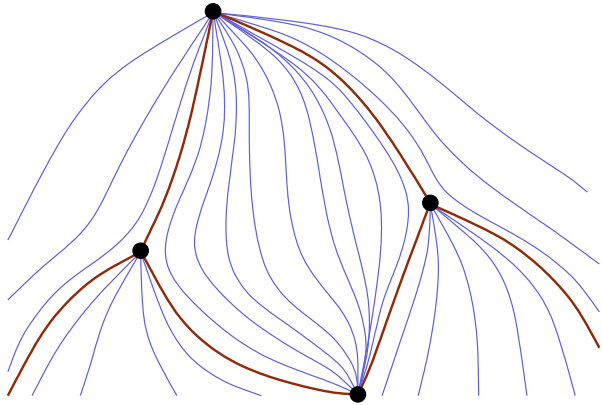}
\caption{\label{fig:brokenMorse} 
One-parameter families of gradient flow lines on a Riemannian manifold
degenerate to broken flow lines.}
\end{figure}

We say that $\mM(x,y)$ has a natural \defin{compatification}
$\overline{\mM}(x,y)$, which has the topology of a compact $1$-manifold
with boundary, and its boundary is the set of all \defin{broken flow lines}
from $x$ to~$y$, cf.~Figure~\ref{fig:brokenMorse}.
This set of broken flow lines is precisely what is counted
if one computes the $\langle y \rangle$ coefficient of $\p^2 \langle x \rangle$,
hence we deduce
$$
\p^2 = 0
$$
as a consequence of the fact that compact $1$-manifolds always have zero
boundary points when counted with appropriate signs.\footnote{Counting with
signs presumes that we have chosen suitable orientations for the moduli
spaces $\mM(x,y)$, and this can always be done.  Alternatively, one can
avoid this issue by counting modulo~$2$ and thus define a homology theory
with $\ZZ_2$ coefficients.}
The homology of the resulting chain complex can be denoted by
$HM_*(M\,;\,g,f)$ and is called the \defin{Morse homology} of~$M$.  The well-known
Morse inequalities can then be deduced from a fundamental theorem stating
that $HM_*(M\,;\,g,f)$ is, for generic $f$ and~$g$, isomorphic to the
singular homology of~$M$.

With the above notion of Morse homology understood, Floer's approach to the
Arnold conjecture can now be summarized as follows:
\begin{itemize}
\item[\textsl{Step~1:}] Under suitable technical assumptions, construct a 
homology theory 
$$
HF_*(M,\omega \,;\,H,\{J_t\}),
$$
depending \emph{a priori} on the choices of a Hamiltonian
$H : S^1 \times M \to \RR$ with all $1$-periodic orbits nondegenerate, and a
generic $S^1$-parametrized family of $\omega$-compatible almost complex 
structures~$\{J_t\}_{t \in S^1}$.  The generators
of the chain complex are the critical points of the symplectic action
functional $\aA_H$, i.e.~$1$-periodic orbits of the Hamiltonian flow,
and the boundary map is defined by counting a suitable notion of gradient
flow lines connecting pairs of orbits (more on this below).
\item[\textsl{Step~2:}] Prove that $HF_*(M,\omega) := HF_*(M,\omega \,;\,H,\{J_t\})$
is a \emph{symplectic invariant}, i.e.~it depends on~$\omega$, but not on the
auxiliary choices $H$ and~$\{J_t\}$.
\item[\textsl{Step~3:}] Show that if $H$ and $\{J_t\}$ are chosen to
be time-independent and $H$ is also $C^2$-small, then the chain complex for
$HF_*(M,\omega \,;\,H,\{J_t\})$ is isomorphic (with a suitable grading shift)
to the chain complex for Morse homology
$HM_*(M\,;\,g,H)$ with $g := \omega(\cdot,J_t\cdot)$.  The isomorphism between
$HM_*(M\,;\,g,H)$ and singular homology thus implies that the Floer complex
must have at least as many generators (i.e.~periodic orbits) as there are
generators of $H_*(M)$, proving the Arnold conjecture.
\end{itemize}

The implementation of Floer's idea required a different type of analysis than
what is needed for Morse homology.  The moduli space $\mM(x,y)$ in Morse
homology is simple to understand as the (generically transverse) intersection 
between the unstable manifold of $x$ and the stable manifold of $y$ with
respect to the negative gradient flow.  Conveniently, both of those are
finite-dimensional manifolds, with their dimensions determined by the
Morse indices of $x$ and~$y$.  We will see in Lecture~\ref{lec:asymptotic}
that no such thing is true
for the symplectic action functional: to the extent that $\aA_H$ can be
thought of as a Morse function on an infinite-dimensional manifold, its
Morse index and its Morse ``co-index'' at every critical point are both
infinite, hence the stable and unstable manifolds are not nearly as nice
as finite-dimensional manifolds, providing no reason to expect that their
intersection should be.  There are additional problems since
$C^\infty_\contr(S^1,M)$ does not have a Banach space topology: in order
to view the negative gradient flow of $\aA_H$ as an ODE and make use of
the usual local existence/uniqueness theorems 
(as in \cite{Lang:geometry}*{Chapter~IV}), one would have to extend to
$\aA_H$ to a smooth function on a suitable Hilbert manifold with a
Riemannian metric.  There is a very limited range of situations in which
one can do this and obtain a reasonable formula for $\nabla \aA_H$,
e.g.~\cite{HoferZehnder}*{\S 6.2} explains the case $M = \TT^{2n}$, in which
$\aA_H$ can be defined on the Sobolev space $H^{1/2}(S^1,\RR^{2n})$ and then
studied using Fourier series.  This approach is very dependent on the fact that
the torus $\TT^{2n}$ is a quotient of~$\RR^{2n}$; for general symplectic 
manifolds $(M,\omega)$, one cannot even 
define $H^{1/2}(S^1,M)$ since functions of class $H^{1/2}$ on $S^1$ need not 
be continuous ($H^{1/2}$ is a ``Sobolev borderline case'' in dimension one).

\begin{figure}
\includegraphics{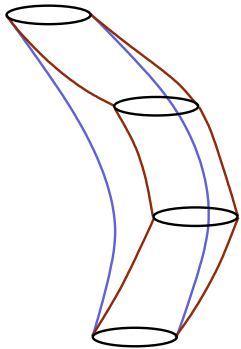}
\caption{\label{fig:brokenFloer} 
A family of smooth Floer trajectories can degenerate into a broken Floer
trajectory.}
\end{figure}

One of the novelties in Floer's
approach was to refrain from viewing the gradient flow as an ODE in a Banach 
space setting, but instead to write down a formal version of the
gradient flow equation and regard it as an elliptic PDE.  To this end,
let us regard $C^\infty_\contr(S^1,M)$ formally as a manifold with
tangent spaces
$$
T_\gamma C^\infty_\contr(S^1,M) := \Gamma(\gamma^*TM),
$$
choose a formal Riemannian metric on this manifold (i.e.~a smoothly varying
family of $L^2$~inner products on the spaces $\Gamma(\gamma^*TM)$) and write 
down the resulting equation for the negative gradient flow.  A suitable
Riemannian metric can be defined by choosing a smooth $S^1$-parametrized
family of compatible almost complex structures
$$
\left\{J_t \in \jJ(M,\omega) \right\}_{t \in S^1},
$$
abbreviated in the following as $\{J_t\}$, and setting
$$
\langle \xi,\eta \rangle_{L^2} := \int_{S^1} \omega(\xi(t),J_t \eta(t))\, dt
$$
for $\xi, \eta \in \Gamma(\gamma^*TM)$.  Exercise~\ref{EX:firstVariation}
then yields the formula
$$
d\aA_H(\gamma) \eta = \langle J_t (\dot{\gamma} - X_t(\gamma)) , \eta \rangle_{L^2},
$$
so that it seems reasonable to define the so-called \emph{unregularized}
gradient of $\aA_H$ by
\begin{equation}
\label{eqn:unregularized}
\nabla \aA_H(\gamma) := J_t (\dot{\gamma} - X_t(\gamma)) \in \Gamma(\gamma^*TM).
\end{equation}
Let us also think of a path $u : \RR \to C^\infty_\contr(S^1,M)$ as a map
$u : \RR \times S^1 \to M$, writing $u(s,t) := u(s)(t)$.
The negative gradient flow equation
$\p_s u + \nabla \aA_H(u(s)) = 0$ then becomes the elliptic PDE
\begin{equation}
\label{eqn:FloerEquation}
\p_s u + J_t(u) \left( \p_t u - X_t(u) \right) = 0.
\end{equation}
This is called the \defin{Floer equation}, and its solutions are often called
\defin{Floer trajectories}.  The relevance of Floer homology to our previous
discussion of pseudoholomorphic curves should now be obvious.  Indeed, the
resemblance of the Floer equation to the nonlinear Cauchy-Riemann equation
is not merely superficial---we will see in Lecture~\ref{lec:cobordisms} that the former
can always be viewed as a special case of the latter.  In any case,
one can use the same set of analytical techniques for both: elliptic regularity
theory implies that Floer trajectories are always smooth, 
Fredholm theory and the implicit function theorem imply that 
(under appropriate assumptions) they form smooth
finite-dimensional moduli spaces.  Most importantly, the same ``bubbling off''
analysis that underlies Gromov's compactness theorem can be used to prove
that spaces of Floer trajectories are compact up to ``breaking'', just as
in Morse homology (see Figure~\ref{fig:brokenFloer})---this is the main
reason for the relation $\p^2 = 0$ in Floer homology.

We should mention one complication that does not arise either in the
study of closed holomorphic curves or in finite-dimensional Morse theory.
Since the gradient flow in Morse homology takes place on a closed manifold,
it is obvious that every gradient flow line asymptotically approaches
critical points at both $-\infty$ and~$+\infty$.  The following example
shows that in the infinite-dimensional setting of Floer theory, this is no
longer true.

\begin{example}
Consider the Floer equation on $M := S^2 = \CC \cup \{\infty\}$ with $H := 0$ and
$J_t$ defined as the standard complex structure~$i$ for every~$t$.
Then the orbits of $X_t$ are all constant, and a map 
$u : \RR \times S^1 \to S^2$ satisfies the Floer equation if and only if it
is holomorphic.  Identifying $\RR \times S^1$ with $\CC^* := \CC \setminus \{0\}$
via the biholomorphic map $(s,t) \mapsto e^{2\pi (s+it)}$, a solution $u$
approaches periodic orbits as $s \to \pm\infty$ if and only if the 
corresponding holomorphic map $\CC^* \to S^2$ extends continuously
(and therefore holomorphically) over $0$ and~$\infty$.  But this is not true
for every holomorphic map $\CC^* \to S^2$, e.g.~take any entire function
$\CC \to \CC$ that has an essential singularity at~$\infty$.
\end{example}
\begin{exercise}
Show that in the above example with an essential singularity at~$\infty$,
the symplectic action $\aA_H(u(s,\cdot))$ is unbounded as $s \to \infty$.
\end{exercise}

\begin{exercise}
\label{EX:FloerEnergy}
Suppose $u : \RR \times S^1 \to M$ is a solution to the Floer equation with
$\lim_{s \to \pm \infty} u(s,\cdot) = \gamma_\pm$ uniformly for a pair
of $1$-periodic orbits $\gamma_\pm \in \Crit(\aA_H)$.  Show that
\begin{equation}
\label{eqn:FloerEnergy}
\aA(\gamma_-) - \aA(\gamma_+) = \int_{\RR \times S^1} 
\omega(\p_s u, \p_t u - X_t(u)) \, ds\, dt =
\int_{\RR \times S^1} \omega(\p_s u, J_t(u) \p_s u) \, ds \, dt.
\end{equation}
\end{exercise}

The right hand side of \eqref{eqn:FloerEnergy} is manifestly nonnegative
since $J_t$ is compatible with~$\omega$, and it is strictly positive unless
$\gamma_- = \gamma_+$.  It is therefore sensible to call this expression
the \defin{energy} $E(u)$ of a Floer trajectory.  The following converse of
Exercise~\ref{EX:FloerEnergy} plays a crucial role in the compactness theory
for Floer trajectories, as it guarantees that all the ``levels'' in a
broken Floer trajectory are asymptotically well behaved.  We will prove
a variant of this result in the SFT context (see Prop.~\ref{prop:HoferEnergy}
below) in Lecture~\ref{lec:compactness}.

\begin{prop}
\label{prop:finiteEnergyFloer}
If $u : \RR \times S^1 \to M$ is a Floer trajectory with $E(u) < \infty$ and
all $1$-periodic orbits of $X_t$ are nonegenerate, then there exist orbits
$\gamma_-,\gamma_+ \in \Crit(\aA_H)$ such that
$\lim_{s \to \pm\infty} u(s,\cdot) = \gamma_\pm$ uniformly.
\qed
\end{prop}

\begin{remark}
It should be emphasized again that we have assumed $[\omega]|_{\pi_2(M)} = 0$
throughout this discussion; Floer homology can also be defined under more
general assumptions, but several details become more complicated.
\end{remark}

For nice comprehensive treatments of Hamiltonian Floer homology---unfortunately
not always with the same sign conventions as used here---see
\cites{Salamon:Floer,AudinDamian}.  Note that this is only one of a few ``Floer homologies''
that were introduced by Floer in the late 80's: the others include
\emph{Lagrangian intersection Floer homology} \cite{Floer:Lagrangian}
(which has since evolved into the \emph{Fukaya category}, see \cite{Seidel:book}), and
\emph{instanton homology} \cite{Floer:instanton}, an extension of 
Donaldson's gauge-theoretic smooth $4$-manifold invariants to dimension three.
The development of new Floer-type theories has since become a major industry.

\section{Contact manifolds and the Weinstein conjecture}
\label{sec:Weinstein}

\begin{figure}
\includegraphics{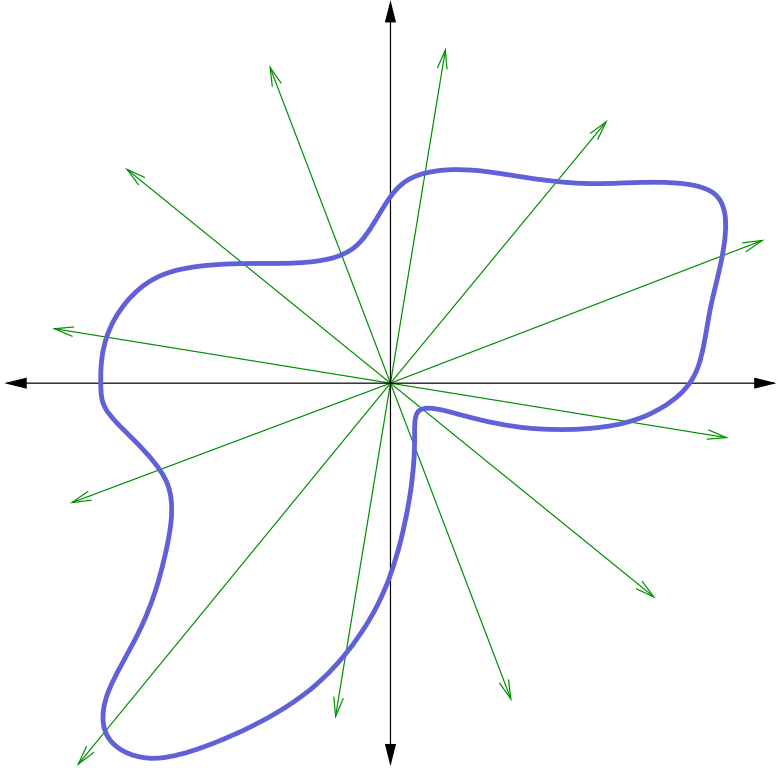}
\caption{\label{fig:starShaped} 
A star-shaped hypersurface in Euclidean space}
\end{figure}

A Hamiltonian system on a symplectic manifold $(W,\omega)$ is called 
\defin{autonomous} if the Hamiltonian $H : W \to \RR$ does
not depend on time.  In this case, the Hamiltonian vector field $X_H$
defined by
$$
\omega(X_H,\cdot) = -dH
$$
is time-independent and its orbits are confined to level sets of~$H$.
The images of these orbits on a given regular level set $H^{-1}(c)$ depend
on the geometry of $H^{-1}(c)$ but not on
$H$ itself, as they are the integral curves (also known as
\defin{characteristics}) of the \defin{characteristic
line field} on $H^{-1}(c)$, defined as the unique direction spanned by a
vector $X$ such that $\omega(X,Y) = 0$ for all $Y$ tangent to $H^{-1}(c)$.
In 1978, Weinstein \cite{Weinstein:convex} and 
Rabinowitz \cite{Rabinowitz:starshaped} proved that certain kinds of regular
level sets in symplectic manifolds are guaranteed to admit closed 
characteristics, hence implying the existence of periodic
Hamiltonian orbits.  In particular, this is true whenever $H^{-1}(c)$ is a
\emph{star-shaped} hypersurface in the standard symplectic~$\RR^{2n}$
(see Figure~\ref{fig:starShaped}).

The following symplectic interpretation of the star-shaped condition provides
both an intuitive reason to believe Rabinowitz's existence result and
motivation for the more general conjecture of Weinstein.  In any symplectic
manifold $(W,\omega)$, a \defin{Liouville vector field} is a smooth vector
field $V$ that satisfies
$$
\Lie_V \omega = \omega.
$$
By Cartan's formula for the Lie derivative, the dual $1$-form $\lambda$
defined by $\lambda := \omega(V,\cdot)$ satisfies $d\lambda = \omega$ if
and only if $V$ is a Liouville vector field; moreover, $\lambda$ then also
satisfies $\Lie_V \lambda = \lambda$, and it is referred to as a
\defin{Liouville form}.  A hypersurface $M \subset (W,\omega)$
is said to be of \defin{contact type} if it is transverse to a Liouville
vector field defined on a neighborhood of~$M$.

\begin{example}
\label{ex:stdLiouville}
Using coordinates $(q_1,p_1,\ldots,q_n,p_n)$ on $\RR^{2n}$, the standard
symplectic form is written as
$$
\omega\std := \sum_{j=1}^n d p_j \wedge d q_j,
$$
and the Liouville form $\lambda\std := \frac{1}{2} \sum_{j=1}^n
(p_j \, d q_j - q_j \, d p_j)$ is dual to the radial Liouville vector field
$$
V\std := \frac{1}{2} \sum_{j=1}^n \left( p_j \frac{\p}{\p p_j} +
q_j \frac{\p}{\p q_j} \right).
$$
Any star-shaped hypersurface is therefore of contact type.
\end{example}

\begin{exercise}
\label{EX:contactType}
Suppose $(W,\omega)$ is a symplectic manifold of dimension~$2n$, 
$M \subset W$ is a smoothly embedded and oriented hypersurface,
$V$ is a Liouville vector field defined near $M$ and $\lambda := 
\omega(V,\cdot)$ is the dual Liouville form.
Define a $1$-form on $M$ by $\alpha := \lambda|_{TM}$.
\begin{enumerate}
\renewcommand{\labelenumi}{(\alph{enumi})}
\item 
Show that $V$ is positively transverse to $M$ if and only if $\alpha$
satisfies
\begin{equation}
\label{eqn:contact}
\alpha \wedge (d\alpha)^{n-1} > 0.
\end{equation}
\item
If $V$ is positively transverse to $M$, choose $\epsilon > 0$
sufficiently small and consider the embedding
$$
\Phi : (-\epsilon,\epsilon) \times M \hookrightarrow W :
(r,x) \mapsto \varphi^r_V(x),
$$
where $\varphi^t_V$ denotes the time~$t$ flow of~$V$.  Show that
$$
\Phi^*\lambda = e^r \alpha,
$$
hence $\Phi^*\omega = d(e^r \alpha)$.
\end{enumerate}
\end{exercise}

The above exercise presents any contact-type hypersurface
$M \subset (W,\omega)$ as one member of a smooth $1$-parameter family
of contact-type hypersurfaces $M_r := \varphi_V^r(M) \subset W$,
each canonically identified with $M$ such that
$\omega|_{TM_r} = e^r \, d\alpha$.  In particular, the characteristic
line fields on $M_r$ are the same for all~$r$, thus the existence of
a closed characteristic on any of these implies that there also exists one
on~$M$.  This observation has sometimes been used to prove such 
existence theorems, e.g.~it is used in \cite{HoferZehnder}*{Chapter~4} to 
reduce Rabinowitz's result to an ``almost existence'' theorem based on 
symplectic capacities.  This discussion hopefully makes the following
conjecture seem believable.

\begin{conj}[Weinstein conjecture, symplectic version]
\label{conj:WeinsteinSymp}
Any closed contact-type hypersurface in a symplectic manifold admits a closed
characteristic.
\end{conj}

Weinstein's conjecture admits a natural rephrasing in the language of contact 
geometry.  A $1$-form $\alpha$ on an oriented $(2n-1)$-dimensional manifold 
$M$ is called a (positive) \defin{contact form} if it satisfies
\eqref{eqn:contact}, and the resulting co-oriented hyperplane field
$$
\xi := \ker \alpha \subset TM
$$
is then called a (positive and co-oriented) \defin{contact structure}.\footnote{The adjective
``positive'' refers to the fact that the orientation of $M$ agrees with
the one determined by the volume form $\alpha \wedge (d\alpha)^{n-1}$; we
call $\alpha$ a \emph{negative} contact form if these two orientations
disagree.  It is also possible in general to define contact structures without
co-orientations, but contact structures of this type will never appear
in these notes; for our purposes, the co-orientation is \emph{always}
considered to be part of the data of a contact structure.}
We call the pair $(M,\xi)$ a \defin{contact manifold}, and refer to a
diffeomorphism $\varphi : M \to M'$ as a \defin{contactomorphism}
from $(M,\xi)$ to $(M',\xi')$ if $\varphi_*$ maps $\xi$ to $\xi'$
and also preserves the respective co-orientations.  Equivalently, if
$\xi$ and $\xi'$ are defined via contact forms $\alpha$ and $\alpha'$
respectively, this means
$$
\varphi^*\alpha' = f \alpha \quad \text{ for some $f \in C^\infty(M,(0,\infty))$}.
$$

Contact topology studies the category of contact manifolds 
$(M,\xi)$ up to contactomorphism.  The following basic result provides 
one good reason to regard $\xi$ 
rather than $\alpha$ as the geometrically meaningful data, as the result
holds for contact \emph{structures}, but not for contact \emph{forms}.

\begin{thm}[Gray's stability theorem]
\label{thm:Gray}
If $M$ is a closed $(2n-1)$-dimensional manifold and $\{\xi_t\}_{t \in [0,1]}$
is a smooth $1$-parameter family of contact structures on $M$, 
then there exists a smooth $1$-parameter family of diffeomorphisms
$\{ \varphi_t \}_{t\in[0,1]}$ such that $\varphi_0 = \Id$ and
$(\varphi_t)_*\xi_0 = \xi_t$.
\end{thm}
\begin{proof}
See \cite{Geiges:book}*{\S 2.2} or \cite{Wendl:lecturesV33}*{Theorem~1.6.12}.
\end{proof}

A corollary is that while the contact form $\alpha$ induced on a contact-type
hypersurface $M \subset (W,\omega)$ via Exercise~\ref{EX:contactType}
is not unique, its induced contact structure is unique up to isotopy.
Indeed, the space of all Liouville vector fields transverse to $M$ is very
large (e.g.~one can add to $V$ any sufficiently small Hamiltonian vector
field), but it is \emph{convex}, hence any two choices of the induced 
contact form $\alpha$ on $M$ are connected by a smooth $1$-parameter family
of contact forms, implying an isotopy of contact structures via
Gray's theorem.

\begin{exercise}
\label{EX:dalpha}
If $\alpha$ is a nowhere zero $1$-form on $M$ and $\xi = \ker \alpha$,
show that $\alpha$ is contact if and only if $d\alpha|_\xi$ defines a
symplectic vector bundle structure on $\xi \to M$.  Moreover, the orientation
of $\xi$ determined by this symplectic bundle structure is compatible with
the co-orientation determined by $\alpha$ and the orientation of $M$ for
which $\alpha \wedge (d\alpha)^{n-1} > 0$.
\end{exercise}

The following definition is based on the fact that since $d\alpha|_\xi$ is
nondegenerate when $\alpha$ is contact, $\ker d\alpha \subset TM$ is
always $1$-dimensional and transverse to~$\xi$.

\begin{defn}
\label{defn:Reeb}
Given a contact form $\alpha$ on $M$, the \defin{Reeb vector field}
is the unique vector field $R_\alpha$ that satisfies
$$
d\alpha(R_\alpha,\cdot) \equiv 0, \quad\text{ and }\quad
\alpha(R_\alpha) \equiv 1.
$$
\end{defn}

\begin{exercise}
\label{EX:ReebFlow}
Show that the flow of any Reeb vector field $R_\alpha$ preserves both
$\xi = \ker\alpha$ and the symplectic vector bundle structure
$d\alpha|_\xi$.
\end{exercise}

\begin{conj}[Weinstein conjecture, contact version]
On any closed contact manifold $(M,\xi)$ with contact form $\alpha$, the
Reeb vector field $R_\alpha$ admits a periodic orbit.
\end{conj}

To see that this is equivalent to the symplectic version of the conjecture,
observe that any contact manifold $(M,\xi = \ker\alpha)$ can be viewed as
the contact-type hypersurface $\{0\} \times M$ in the open symplectic manifold
$$
\left(\RR \times M, d(e^r \alpha) \right),
$$
called the \defin{symplectization} of $(M,\xi)$.

\begin{exercise}
\label{EX:symplectization}
Recall that on any smooth manifold $M$, there is a tautological $1$-form 
$\lambda$ that locally takes the form $\lambda = \sum_{j=1}^n p_j \, d q_j$ 
in any choice of local coordinates $(q_1,\ldots,q_n)$ on a neighbood 
$\uU \subset M$, with $(p_1,\ldots,p_n)$ denoting the induced coordinates on
the cotangent fibers over~$\uU$.  This is a Liouville form, with $d\lambda$
defining the canonical symplectic structure of $T^*M$.  Now if $\xi \subset TM$
is a co-oriented hyperplane field on~$M$, consider the submanifold
$$
S_\xi M := \left\{ p \in T^*M\ \big|\ 
\text{$\ker p = \xi$ and $p(X) > 0$ for any $X \in TM$ pos.~transverse
to~$\xi$} \right\}.
$$
Show that $\xi$ is contact if and only if $S_\xi M$ is a symplectic submanifold
of $(T^*M,d\lambda)$, and the Liouville vector field on $T^*M$ dual to
$\lambda$ is tangent to~$S_\xi M$.  Moreover, if $\xi$ is contact, then any 
choice of contact form for $\xi$ determines a diffeomorphism of
$S_\xi M$ to $\RR \times M$ identifying the Liouville form $\lambda$ along
$S_\xi M$ with $e^r \alpha$.
\end{exercise}
\begin{remark}
Exercise~\ref{EX:symplectization} shows that up to symplectomorphism, 
our definition of the symplectization of $(M,\xi)$ above actually 
depends only on $\xi$ and not on~$\alpha$.
\end{remark}

In 1993, Hofer \cite{Hofer:weinstein} introduced a 
new approach to the Weinstein conjecture that
was based in part on ideas of Gromov and Floer.  Fix a contact manifold
$(M,\xi)$ with contact form~$\alpha$, and let
$$
\jJ(\alpha) \subset \jJ(\RR \times M)
$$
denote the nonempty and contractible space of all almost complex structures 
$J$ on $\RR \times M$ satisfying the following conditions:
\begin{enumerate}
\item The natural translation action on $\RR \times M$ preserves~$J$;
\item $J \p_r = R_\alpha$ and $J R_\alpha = -\p_r$, where $r$ denotes the
canonical coordinate on the $\RR$-factor in $\RR \times M$;
\item $J\xi = \xi$ and $d\alpha(\cdot,J\cdot)|_\xi$ defines a bundle metric
on~$\xi$.
\end{enumerate}
It is easy to check that any $J \in \jJ(\alpha)$ is compatible with the
symplectic structure $d(e^r \alpha)$ on $\RR \times M$.  Moreover, if
$\gamma : \RR \to M$ is any periodic orbit of $R_\alpha$ with period $T > 0$,
then for any $J \in \jJ(\alpha)$, the so-called \defin{trivial cylinder}
$$
u : \RR \times S^1 \to \RR \times M : (s,t) \mapsto (Ts,\gamma(Tt))
$$
is a $J$-holomorphic curve.  Following Floer, one version of Hofer's idea 
would be to look for $J$-holomorphic cylinders that satisfy a finite energy
condition as in Prop.~\ref{prop:finiteEnergyFloer} forcing them to approach
trivial cylinders asymptotically---the existence of such a cylinder would
then imply the existence of a closed Reeb orbit and thus prove the
Weinstein conjecture.  The first hindrance is that the ``obvious'' definition
of energy in this context,
$$
\int_{\RR \times S^1} u^*d(e^r \alpha),
$$
is not the right one: this integral is infinite if $u$ is a trivial cylinder.
To circumvent this, notice that every $J \in \jJ(\alpha)$ is also compatible
with any symplectic structure of the form
$$
\omega_\varphi := d(e^{\varphi(r)} \alpha),
$$
where $\varphi$ is a function chosen freely from the set
\begin{equation}
\label{eqn:tT}
\tT := \left\{ \varphi \in C^\infty(\RR,(-1,1))\ \big|\  \varphi' > 0 \right\}.
\end{equation}
Essentially, choosing $\omega_\varphi$ means identifying $\RR \times M$ with
a subset of the bounded region $(-1,1) \times M$, in which trivial cylinders
have finite symplectic area.  Since there is no preferred choice for the
function $\varphi$, we define the \defin{Hofer energy}\footnote{Strictly
speaking, the energy defined in \eqref{eqn:HoferEnergy} is not identical
to the notion introduced in \cite{Hofer:weinstein} and used in many of
Hofer's papers, but it is equivalent to
it in the sense that uniform bounds on either notion of energy imply
uniform bounds on the other.}
of a $J$-holomorphic curve $u : \Sigma \to \RR \times M$ by
\begin{equation}
\label{eqn:HoferEnergy}
E(u) := \sup_{\varphi \in \tT} \int_{\Sigma} u^*\omega_\varphi.
\end{equation}
This has the desired property of being finite for trivial cylinders, and it
is also nonnegative, with strict positivity whenever $u$ is not constant.

Another useful observation from \cite{Hofer:weinstein} was that if the goal
is to find periodic orbits, then we need not restrict our attention to
$J$-holomorphic \emph{cylinders} in particular.  One can more generally
consider curves defined on an arbitrary \emph{punctured} Riemann surface
$$
\dot{\Sigma} := \Sigma \setminus \Gamma,
$$
where $(\Sigma,j)$ is a closed connected Riemann surface and
$\Gamma \subset \Sigma$ is a finite set of punctures.  For any 
$\zeta \in \Gamma$, one can find coordinates identifying
some punctured neighborhood of $\zeta$ biholomorphically with the
closed punctured disk
$$
\dot{\DD} := \DD \setminus \{0\} \subset \CC,
$$
and then identify this with either the positive or negative half-cylinder
$$
Z_+ := [0,\infty) \times S^1, \qquad Z_- := (-\infty,0] \times S^1
$$
via the biholomorphic maps
$$
Z_+ \to \dot{\DD} : (s,t) \mapsto e^{-2\pi (s+it)}, \qquad
Z_- \to \dot{\DD} : (s,t) \mapsto e^{2\pi (s+it)}.
$$
We will refer to such a choice as a (positive or negative) 
\defin{holomorphic cylindrical coordinate} system near~$\zeta$, and
in this way, we can present $(\dot{\Sigma},j)$ as a \emph{Riemann surface
with cylindrical ends}, i.e.~the union of some compact Riemann surface with
boundary with a finite collection of half-cylinders $Z_\pm$ on which
$j$ takes the standard form $j\p_s = \p_t$.  Note that the standard cylinder
$\RR \times S^1$ is a special case of this, as it can be identified
biholomorphically with $S^2 \setminus \{0,\infty\}$.  Another important
special case is the plane, $\CC = S^2 \setminus \{\infty\}$.

\begin{figure}
\includegraphics{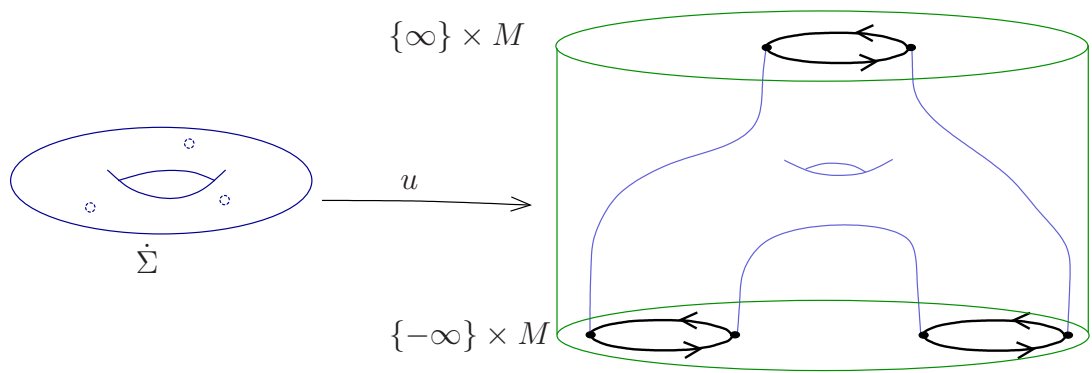}
\caption{\label{fig:sympCurve} 
An asymptotically cylindrical holomorphic curve in a symplectization,
with genus~$1$, one positive puncture and two negative punctures.}
\end{figure}

If $u : (\dot{\Sigma},j) \to (\RR \times M,J)$ is a $J$-holomorphic curve and
$\zeta \in \Gamma$ is one of its punctures, we will say that $u$ is
\defin{positively/negatively asymptotic} to a $T$-periodic Reeb orbit 
$\gamma : \RR \to M$ at $\zeta$ if one can choose 
holomorphic cylindrical coordinates $(s,t) \in Z_\pm$ near~$\zeta$ such that
$$
u(s,t) = \exp_{(Ts,\gamma(Tt))} h(s,t) \quad \text{ for $|s|$ sufficiently large},
$$
where $h(s,t)$ is a vector field along the trivial cylinder satisfying
$h(s,\cdot) \to 0$ uniformly as $|s| \to \infty$, and the exponential map
is defined with respect to any $\RR$-invariant choice of Riemannian metric
on $\RR \times M$.  We say that $u : (\dot{\Sigma},j) \to (\RR \times M,J)$
is \defin{asymptotically cylindrical} if it is (positively or negatively)
asymptotic to some closd Reeb orbit at each of its punctures.  Note that this
partitions the finite set of punctures $\Gamma \subset \Sigma$ into two
subsets,
$$
\Gamma = \Gamma^+ \cup \Gamma^-,
$$
the \emph{positive} and \emph{negative} punctures respectively,
see Figure~\ref{fig:sympCurve}.

\begin{exercise}
\label{EX:HoferEnergy}
Suppose $u : (\dot{\Sigma},j) \to (\RR \times M, J)$ is an asymptotically
cylindrical $J$-holomorphic curve, with the asymptotic orbit at each
puncture $\zeta \in \Gamma^\pm$ denoted by~$\gamma_\zeta$, having period
$T_\zeta > 0$.  Show that
$$
\sum_{\zeta \in \Gamma^+} T_\zeta -
\sum_{\zeta \in \Gamma^-} T_\zeta
= \int_{\dot{\Sigma}} u^*d\alpha \ge 0,
$$
with equality if and only if the image of $u$ is contained in that of a
trivial cylinder.  In particular, $u$ must have at least one positive
puncture unless it is constant.
Show also that $E(u)$ is finite and satisfies an upper
bound determined only by the periods of the positive asymptotic orbits.
\end{exercise}

The following analogue of Prop.~\ref{prop:finiteEnergyFloer} will be proved
in Lecture~\ref{lec:compactness}.
For simplicity, we shall state a weakened version of what Hofer proved
in \cite{Hofer:weinstein}, which did not require any nondegeneracy
assumption.  A $T$-periodic Reeb orbit $\gamma : \RR \to M$ is called
\defin{nondegenerate} if the Reeb flow $\varphi_\alpha^t$ has the
property that its linearization along the contact bundle
(cf.~Exercise~\ref{EX:ReebFlow}),
$$
d\varphi_\alpha^T(\gamma(0))|_{\xi_{\gamma(0)}} : \xi_{\gamma(0)} \to
\xi_{\gamma(0)}
$$
does not have $1$ as an eigenvalue.  Note that since $R_\alpha$ is not
time-dependent, closed Reeb orbits are never completely isolated---they always
exist in $S^1$-parametrized families---but these families are isolated in
the nondegenerate case.

\begin{prop}
\label{prop:HoferEnergy}
Suppose $(M,\xi)$ is a closed contact manifold, with a contact form
$\alpha$ such that all closed Reeb orbits are nondegenerate.  If
$u : (\dot{\Sigma},j) \to (\RR \times M,J)$ is a $J$-holomorphic curve
with $E(u) < \infty$ on a punctured Riemann surface such that none of the
punctures are removable, then $u$ is asymptotically cylindrical.
\qed
\end{prop}

The main results in
\cite{Hofer:weinstein} state that under certain assumptions on a closed
contact $3$-manifold $(M,\xi)$, namely if either $\xi$ is 
\emph{overtwisted} (as defined in \cite{Eliashberg:overtwisted}) or
$\pi_2(M) \ne 0$, one can find for any contact form $\alpha$ on $(M,\xi)$
and any $J \in \jJ(\alpha)$ a finite-energy $J$-holomorphic plane.
By Proposition~\ref{prop:HoferEnergy}, this implies the existence of a
contractible periodic Reeb orbit and thus proves the Weinstein conjecture
in these settings.

\section{Symplectic cobordisms and their completions}
\label{sec:cobordisms}

After the developments described in the previous three sections, it seemed
natural that one might define invariants of contact manifolds via a
Floer-type theory generated by closed Reeb orbits and counting asymptotically
cylindrical holomorphic curves in symplectizations.  This theory is what is
now called SFT, and its basic structure was outlined in a paper 
by Eliashberg, Givental and Hofer \cite{SFT} in 2000, though some of
its analytical foundations remain unfinished in 2016.  The term ``field theory''
is an allusion to ``topological quantum field theories,'' which associate
vector spaces to certain geometric objects and morphisms to cobordisms
between those objects.  Thus in order to place SFT in its proper setting,
we need to introduce symplectic cobordisms between contact manifolds.

Recall that if $M_+$ and $M_-$ are smooth oriented closed manifolds of the
same dimension, an oriented cobordism from $M_-$ to $M_+$ is a
compact smooth oriented manifold $W$ with oriented boundary
$$
\p W = -M_- \sqcup M_+,
$$
where $-M_-$ denotes $M_-$ with its orientation reversed.  Given positive
contact structures $\xi_\pm$ on $M_\pm$, we say that
a symplectic manifold $(W,\omega)$ is a \defin{symplectic cobordism
from $(M_-,\xi_-)$ to $(M_+,\xi_+)$} if $W$ is an oriented 
cobordism\footnote{We assume of course that $W$ is assigned the orientation
determined by its symplectic form.} from
$M_-$ to $M_+$ such that both components of $\p W$ are contact-type
hypersurfaces with induced contact structures isotopic to~$\xi_\pm$.
Note that our chosen orientation conventions imply in this case that
the Liouville vector field chosen near $\p W$ must point \emph{outward}
at $M_+$ and \emph{inward} at~$M_-$; we say in this case that $M_+$ is
a symplectically \defin{convex} boundary component, while $M_-$ is
symplectically \defin{concave}.  As important special cases,
$(W,\omega)$ is a \defin{symplectic filling} of $(M_+,\xi_+)$ if
$M_- = \emptyset$, and it is a \defin{symplectic cap} of $(M_-,\xi_-)$
if $M_+ = \emptyset$.  In the literature, fillings and caps are sometimes
also referred to as \emph{convex fillings} or \emph{concave fillings}
respectively.

\begin{figure}
\includegraphics{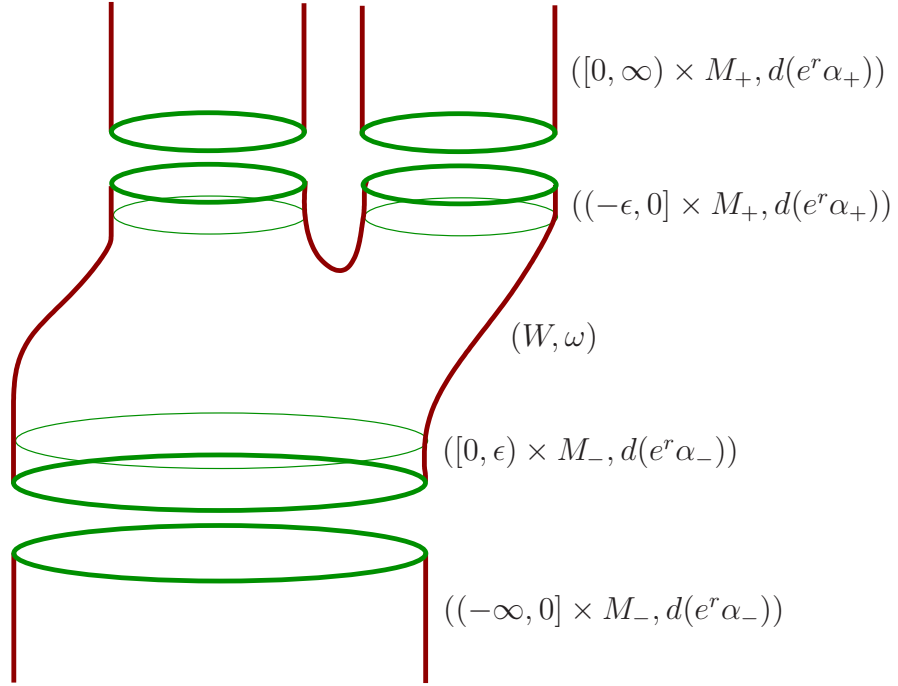}
\caption{\label{fig:completion} The completion of a symplectic cobordism}
\end{figure}

The contact-type condition implies the existence of a Liouville form
$\lambda$ near $\p W$ with $d\lambda = \omega$, such that by
Exercise~\ref{EX:contactType}, neighborhoods of $M_+$ and $M_-$ in $W$ can
be identified with the collars (see Figure~\ref{fig:collars})
$$
(-\epsilon,0] \times M_+ \quad \text{ or } \quad [0,\epsilon) \times M_-
$$
respectively for sufficiently small $\epsilon > 0$, with $\lambda$ taking
the form
$$
\lambda = e^r \alpha_\pm,
$$
where $\alpha_\pm := \lambda|_{TM_\pm}$ are contact forms for~$\xi_\pm$.
The \defin{symplectic completion} of $(W,\omega)$ is the noncompact symplectic
manifold $(\widehat{W},\hat{\omega})$ defined by attaching cylindrical
ends to these collar neighborhoods (Figure~\ref{fig:completion}):
\begin{equation}
\label{eqn:completion}
\begin{split}
(\widehat{W},\hat{\omega}) = \left( (-\infty,0] \times M_-, d(e^r \alpha_-) \right)
&\cup_{M_-} (W,\omega) \\
& \cup_{M_+}
\left( [0,\infty) \times M_+, d(e^r \alpha_+) \right).
\end{split}
\end{equation}
In this context, the symplectization $(\RR \times M,d(e^r \alpha))$ is
symplectomorphic to the completion of the \defin{trivial symplectic cobordism}
$([0,1] \times M, d(e^r \alpha))$ from $(M,\xi=\ker\alpha)$ to itself.
More generally, the object in the following easy exercise can also sensibly 
be called a trivial symplectic cobordism:

\begin{figure}
\includegraphics{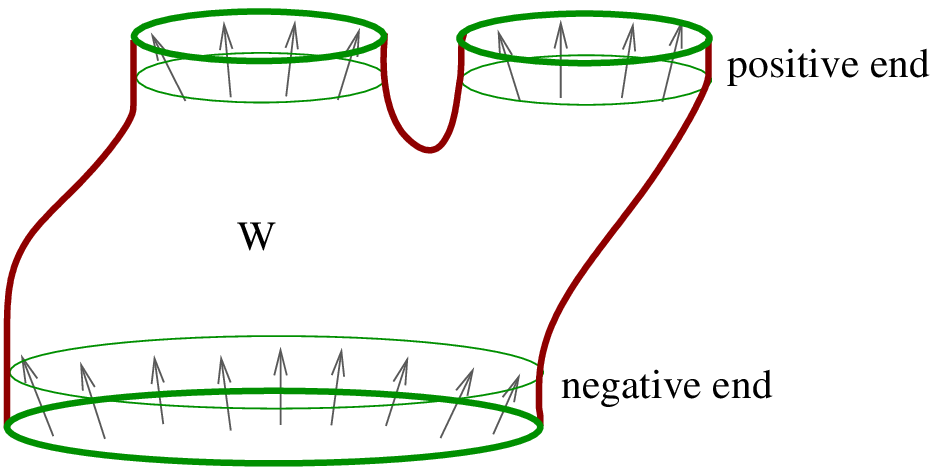}
\caption{\label{fig:collars} A symplectic cobordism with concave boundary
$(M_-,\xi_-)$ and convex boundary $(M_+,\xi_+)$, with symplectic collar
neighborhoods defined by flowing along Liouville vector fields near the boundary.}
\end{figure}

\begin{exercise}
\label{EX:trivialCobordism}
Suppose $(M,\xi)$ is a closed contact manifold with contact form $\alpha$, and
$f_\pm : M \to \RR$ is a pair of functions with $f_- < f_+$ everywhere.
Show that the domain
$$
\left\{ (r,x) \in \RR \times M\ \big|\ f_-(x) \le r \le f_+(x) \right\}
\subset \RR \times M
$$
defines a symplectic cobordism from $(M,\xi)$ to itself, with a global
Liouville form $\lambda = e^r \alpha$ inducing contact forms $e^{f_-} \alpha$ 
and $e^{f_+} \alpha$ on its concave and convex boundaries respectively.
\end{exercise}

We say that $(W,\omega)$ is an \defin{exact symplectic cobordism} or
\defin{Liouville cobordism} if the Liouville form $\lambda$ can be
extended from a neighborhood of $\p W$ to define a global primitive
of $\omega$ on~$W$.  Equivalently, this means that $\omega$ admits a global
Liouville vector field that points inward at $M_-$ and outward at~$M_+$.
An \defin{exact filling} of $(M_+,\xi_+)$ is an exact cobordism whose
concave boundary is empty.  Observe that if $(W,\omega)$ is exact, then
its completion $(\widehat{W},\hat{\omega})$ also inherits a global Liouville
form.

\begin{exercise}
Use Stokes' theorem to show that there is no such thing as an exact
symplectic cap.
\end{exercise}

The above exercise hints at an important difference between cobordisms in
the \emph{symplectic} as opposed to the \emph{oriented smooth} category: 
symplectic cobordisms are not generally reversible.  If $W$ is an oriented 
cobordism from $M_-$
to $M_+$, then reversing the orientation of $W$ produces an oriented
cobordism from $M_+$ to~$M_-$.  But one cannot simply reverse orientations
in the symplectic category, since the orientation is determined by the
symplectic form.  For example, many obstructions to the existence of
symplectic fillings of given contact manifolds are known---some of them
defined in terms of SFT---but we do not know any obstructions at all to 
symplectic caps, in fact it is known that all contact $3$-manifolds admit
them.

\begin{figure}
\includegraphics{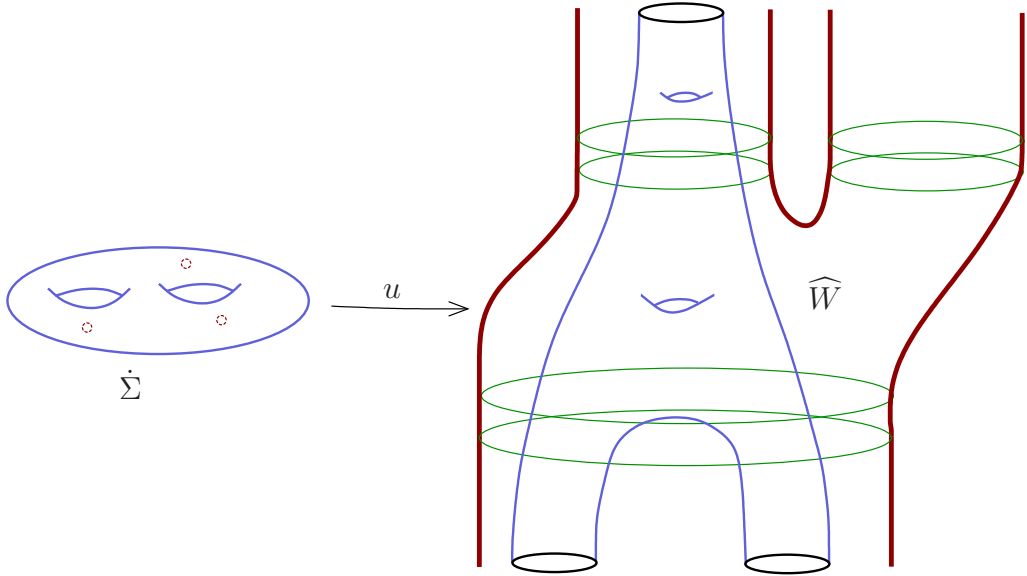}
\caption{\label{fig:asympCyl} 
An asymptotically cylindrical holomorphic curve in a completed symplectic
cobordism, with genus~$2$, one positive puncture and two negative punctures.}
\end{figure}

The definitions for holomorphic curves in symplectizations in the previous
section generalize to completions of symplectic cobordisms in a fairly
straightforward way since these completions look exactly like symplectizations
outside of a compact subset.  Define
$$
\jJ(W,\omega,\alpha_+,\alpha_-) \subset \jJ(\widehat{W})
$$
as the space of all almost complex structures $J$ on $\widehat{W}$ such that
$$
J|_W \in \jJ(W,\omega), \qquad J|_{[0,\infty) \times M_+} \in \jJ(\alpha_+)
\quad \text{ and } \quad
J|_{(-\infty,0] \times M_-} \in \jJ(\alpha_-).
$$
Occasionally it is useful to relax the compatibility condition on $W$ to
tameness,\footnote{It seems natural to wonder whether one could not also
relax the conditions on the cylindrical ends and require $J|_{\xi_\pm}$ to
be tamed by $d\alpha_\pm|_{\xi_\pm}$ instead of compatible with it.
I do not currently know whether this works, but in later lectures we will
see some reasons to worry that it might not.}
i.e.~$J|_W \in \jJ_\tau(W,\omega)$, producing a space that we shall
denote by
$$
\jJ_\tau(W,\omega,\alpha_+,\alpha_-) \subset \jJ(\widehat{W}).
$$
As in Prop.~\ref{prop:contractible}, both of these spaces are 
nonempty and contractible.
We can then consider asymptotically cylindrical $J$-holomorphic curves
$$
u : (\dot{\Sigma} = \Sigma \setminus (\Gamma^+ \cup \Gamma^-),j) \to
(\widehat{W},J),
$$
which are proper maps asymptotic to closed orbits of $R_{\alpha_\pm}$
in $M_\pm$ at punctures in~$\Gamma^\pm$, see Figure~\ref{fig:asympCyl}.

One must again tinker with the
symplectic form on $\widehat{W}$ in order to define a notion of energy
that is finite when we need it to be.  We generalize \eqref{eqn:tT} as
$$
\tT := \left\{ \varphi \in C^\infty(\RR,(-1,1))\ \big| \ 
\text{$\varphi' > 0$ and $\varphi(r) = r$ near $r=0$} \right\},
$$
and associate to each $\varphi \in \tT$ a symplectic form $\hat{\omega}_\varphi$
on $\widehat{W}$ defined by
$$
\hat{\omega}_\varphi :=
\begin{cases}
d(e^{\varphi(r)} \alpha_+) & \text{ on $[0,\infty) \times M_+$},\\
\omega & \text{ on $W$}, \\
d(e^{\varphi(r)} \alpha_-) & \text{ on $(-\infty,0] \times M_-$}.
\end{cases}
$$
One can again check that every $J \in \jJ(W,\omega,\alpha_+,\alpha_-)$ 
or $\jJ_\tau(W,\omega,\alpha_+,\alpha_-)$ is compatible with or,
respectively, tamed by $\hat{\omega}_\varphi$ for every $\varphi \in \tT$.
Thus it makes sense to define the \defin{energy} of
$u : (\dot{\Sigma},j) \to (\widehat{W},J)$ by
$$
E(u) := \sup_{\varphi \in \tT} \int_{\dot{\Sigma}} u^*\hat{\omega}_\varphi.
$$
It will be a straightforward matter to generalize 
Proposition~\ref{prop:HoferEnergy} and show that finite energy implies
asymptotically cylindrical behavior in completed cobordisms.

\begin{exercise}
\label{EX:energyCobordism}
Show that if $(W,\omega)$ is an exact cobordism, then every asymptotically
cylindrical $J$-holomorphic curve in $\widehat{W}$ has at least one
positive puncture.
\end{exercise}

\section{Contact homology and SFT}
\label{sec:CHSFT}

We can now sketch the algebraic structure of SFT.  We shall ignore or
suppress several pesky details that are best dealt with later, some of them
algebraic, others analytical.  Due to analytical problems, some of the 
``theorems'' that we shall (often imprecisely) state in this section are not
yet provable at the current level of technology, though we expect that they
will be soon.  We shall use quotation marks to indicate this caveat wherever
appropriate.

The standard versions of SFT all define homology theories with varying
levels of algebraic structure which are meant to be invariants of a
contact manifold $(M,\xi)$.  The chain complexes always depend on 
certain auxiliary choices, including a nondegenerate contact form
$\alpha$ and a generic $J \in \jJ(\alpha)$.  The generators consist of
formal variables $q_\gamma$, one for each\footnote{Actually I should be
making a distinction here between ``good'' and ``bad'' Reeb orbits, but
let's discuss that later; see Lecture~\ref{lec:orientations}.}
closed Reeb orbit~$\gamma$.  In the most straightforward generalization
of Hamiltonian Floer homology, the chain complex is simply a
graded $\QQ$-vector space generated by the variables $q_\gamma$, and the
boundary map is defined by
$$
\p_{\CCH} q_\gamma = \sum_{\gamma'}
\# \left(  \mM(\gamma,\gamma') \big/ \RR \right) q_{\gamma'},
$$
where $\mM(\gamma,\gamma')$ is the moduli space of $J$-holomorphic cylinders
in $\RR \times M$ with a positive puncture asymptotic to $\gamma$ and a
negative puncture asymptotic to~$\gamma'$, and the sum ranges over all orbits
$\gamma'$ for which this moduli space is $1$-dimensional.  The count
$\# \left( \mM(\gamma,\gamma') / \RR \right)$ is rational, as it includes
rational weighting factors that depend on combinatorial information and are
best not discussed right now.\footnote{Similar combinatorial factors are
hidden behind the symbol ``$\#$'' in our definitions of $\p_\CH$ 
and~$\mathbf{H}$, and will be discussed in earnest in Lecture~\ref{lec:H}.}

\begin{heorem}
\label{thm:CCH}
If $\alpha$ admits no contractible Reeb orbits, then $\p^2_{\CCH} = 0$,
and the resulting homology is independent of the choices of $\alpha$ with
this property and generic $J \in \jJ(\alpha)$.
\end{heorem}

\begin{figure}
\includegraphics[width=5.5in]{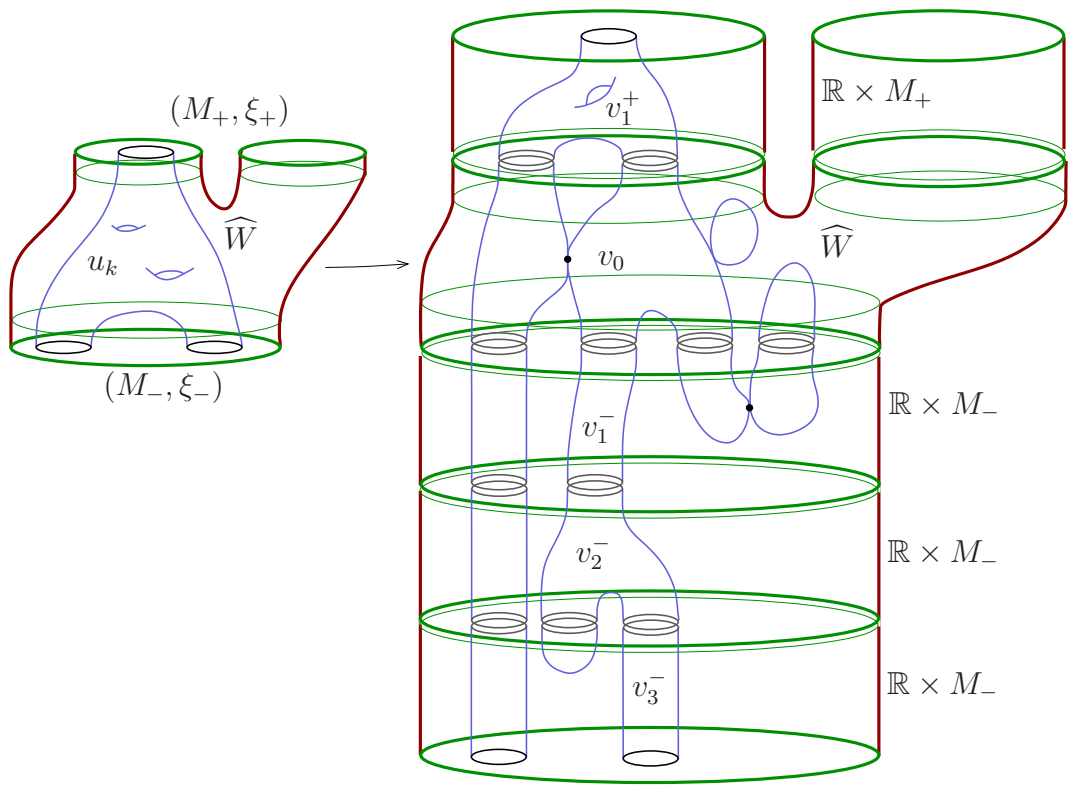}
\caption{\label{fig:SFT} Degeneration of a sequence $u_k$ of finite energy 
punctured holomorphic curves 
with genus~$2$, one positive puncture and two negative punctures in a symplectic cobordism.
The limiting holomorphic building $(v_1^+,v_0,v_1^-,v_2^-,v_3^-)$ in this example
has one upper level living in the symplectization $\RR \times M_+$, 
a main level living in $\widehat{W}$, 
and three lower levels, each of which is a
(possibly disconnected) finite-energy punctured nodal holomorphic curve
in $\RR \times M_-$.
The building has arithmetic genus~$2$ and the same numbers of positive and negative
punctures as~$u_k$.}
\end{figure}

The invariant arising from this result is known as
\defin{cylindrical contact homology}, and it is sometimes quite easy to
work with when it is well defined, though it has the disadvantage of not
always being defined.  Namely,
the relation $\p^2_{\CCH} = 0$ can fail if $\alpha$ admits contractible
Reeb orbits, because unlike in Floer homology,
the compactification of the space of cylinders $\mM(\gamma,\gamma')$
generally includes objects that are not broken cylinders.  In fact,
the objects arising in the ``SFT compactification'' of moduli spaces of
finite-energy curves in completed cobordisms can be quite elaborate,
see Figure~\ref{fig:SFT}.  The combinatorics of the situation are not so
bad however if the cobordism is exact, as is the case for a symplectization:
Exercise~\ref{EX:energyCobordism} then prevents curves without positive
ends from appearing.  The only possible degenerations for cylinders then
consist of broken configurations whose levels each have \emph{exactly one
positive puncture} and arbitrary negative punctures; moreover, all but one
of the negative punctures must eventually be capped off by planes, which
is why ``Theorem''~\ref{thm:CCH} holds in the absence of planes.

If planes do exist, then one can account for them by defining the chain
complex as an \emph{algebra} rather than a vector space, producing the theory
known as \defin{contact homology}.  For this, the chain complex is taken
to be a graded unital algebra over $\QQ$, and we define
$$
\p_{\CH} q_\gamma = \sum_{(\gamma_1,\ldots,\gamma_m)}
\# \left(  \mM(\gamma;\gamma_1,\ldots,\gamma_m) \big/ \RR \right) 
q_{\gamma_1} \ldots q_{\gamma_m},
$$
with $\mM(\gamma;\gamma_1,\ldots,\gamma_m)$ denoting the moduli space of
punctured $J$-holomorphic spheres in $\RR \times M$ with a positive puncture
at $\gamma$ and $m$ negative punctures at the orbits $\gamma_1,\ldots,\gamma_m$,
and the sum ranges over all integers $m \ge 0$ and all $m$-tuples of orbits
for which the moduli space is $1$-dimensional.  
The action of $\p_\CH$ is
then extended to the whole algebra via a graded Leibniz rule
$$
\p_\CH (q_\gamma q_{\gamma'}) := \left(\p_\CH q_\gamma\right) q_{\gamma'} +
(-1)^{|\gamma|} q_\gamma \left( \p_\CH q_{\gamma'} \right).
$$
The general compactness and gluing theory for genus zero curves with one
positive puncture now implies:

\begin{heorem}
$\p^2_\CH = 0$, and the resulting homology is (as a graded unital $\QQ$-algebra)
independent of the choices $\alpha$ and~$J$.
\end{heorem}

Maybe you've noticed the pattern:
in order to accommodate more general classes of
holomorphic curves, we need to add more algebraic structure.  The
\defin{full SFT} algebra counts all rigid holomorphic curves in $\RR \times M$,
including all combinations of positive and negative punctures and all genera.
Here is a brief picture of what it looks like.  Counting all the
$1$-dimensional moduli spaces of $J$-holomorphic curves modulo 
$\RR$-translation in $\RR \times M$ produces a formal power series
$$
\mathbf{H} := \sum \# \left( \mM_g(\gamma_1^+,\ldots,\gamma_{m_+}^+ \,;\,
\gamma_1^-,\ldots,\gamma_{m_-}^-) \Big/ \RR \right) 
q_{\gamma_1^-}\ldots q_{\gamma_{m_-}^-} 
p_{\gamma_1^+}\ldots p_{\gamma_{m_+}^+} \hbar^{g-1},
$$
where the sum ranges over all integers $g, m_+,m_- \ge 0$ and tuples of
orbits, $\hbar$ and $p_\gamma$ (one for each orbit $\gamma$) are additional
formal variables, and 
$$
\mM_g(\gamma_1^+,\ldots,\gamma_{m_+}^+ \,;\,
\gamma_1^-,\ldots,\gamma_{m_-}^-)
$$
denotes the moduli space of $J$-holomorphic
curves in $\RR \times M$ with genus~$g$, $m_+$ positive punctures at the orbits
$\gamma_1^+,\ldots,\gamma_{m_+}^+$, and $m_-$ negative punctures at the orbits
$\gamma_1^-,\ldots,\gamma_{m_-}^+$.  We can regard $\mathbf{H}$ as an
operator on a graded algebra $\Weyl$ of formal power series in the variables 
$\{p_\gamma\}$, $\{q_\gamma\}$ and $\hbar$, equipped with a graded bracket 
operation that satisfies the quantum mechanical commutation relation
$$
[p_\gamma,q_\gamma] = \kappa_\gamma \hbar,
$$
where $\kappa_\gamma$ is a combinatorial factor that is best ignored for now.
Note that due to the signs that accompany the grading, odd elements 
$\mathbf{F} \in \Weyl$ need not satisfy
$[\mathbf{F},\mathbf{F}] = 0$, and $\mathbf{H}$ itself is an odd element,
thus the following statement is nontrivial; in fact, it is the algebraic
manifestation of the general compactness and gluing theory for punctured
holomorphic curves in symplectizations.

\begin{heorem}
$[\mathbf{H},\mathbf{H}] = 0$, hence by the graded Jacobi identity, 
$\mathbf{H}$ determines an operator
$$
D_\SFT : \Weyl \to \Weyl : \mathbf{F} \mapsto [\mathbf{H},\mathbf{F}]
$$
satisfying $D_\SFT^2 = 0$.  The resulting homology depends on $(M,\xi)$ but
not on the auxiliary choices $\alpha$ and~$J$.
\end{heorem}

It takes some time to understand how pictures such as
Figure~\ref{fig:SFT} translate into algebraic relations like
$[\mathbf{H},\mathbf{H}] = 0$, but this is a subject we'll come back to.
There is also an intermediate theory between 
contact homology and full SFT, called \defin{rational SFT}, which counts
only genus zero curves with arbitrary positive and negative punctures.
Algebraically, it is obtained from the full SFT algebra as a
``semiclassical approximation'' by discarding higher-order factors of~$\hbar$
so that the commutation bracket in $\Weyl$ becomes a graded Poisson
bracket.  We will discuss all of this in Lecture~\ref{lec:H}.

\section{Two applications}
\label{sec:applications}

We briefly mention two applications that we will be able to establish rigorously
using the methods developed in this book.  Since SFT itself is not yet
well defined in full generality, this sometimes means using SFT for 
inspiration while proving corollaries via more direct methods.

\subsection{Tight contact structures on $\TT^3$}

The $3$-torus $\TT^3 = S^1 \times S^1 \times S^1$ with coordinates
$(t,\theta,\phi)$ admits a sequence of contact structures
$$
\xi_k := \ker\left( \cos(2\pi k t) \, d\theta + \sin(2\pi k t) \, d\phi \right),
$$
one for each $k \in \NN$.  These cannot be distinguished from each other by
any classical invariants, e.g.~they all have the same Euler class, in fact
they are all homotopic as co-oriented $2$-plane fields.  Nonetheless:

\begin{thm}
For $k \ne \ell$, $(\TT^3,\xi_k)$ and $(\TT^3,\xi_\ell)$ are not
contactomorphic.
\end{thm}

We will be able to prove this in Lecture~\ref{lec:tight3tori}
by rigorously defining and computing
cylindrical contact homology for a suitable choice of contact forms
on $(\TT^3,\xi_k)$.

\subsection{Filling and cobordism obstructions}

\begin{figure}
\includegraphics{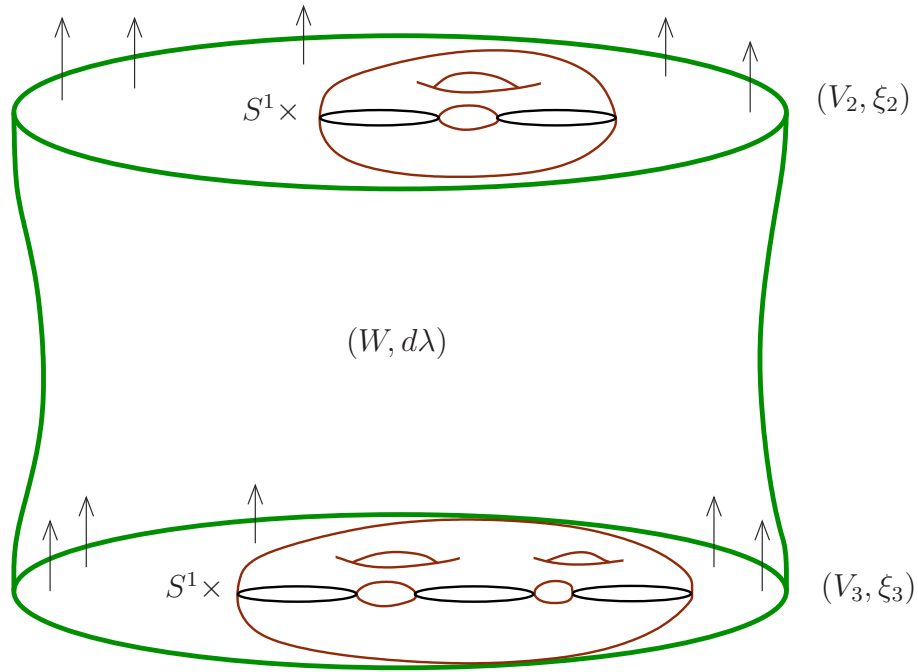}
\caption{\label{fig:torsion} 
This exact symplectic cobordism does not exist.}
\end{figure}

Consider a closed connected and oriented surface 
$\Sigma$ presented as $\Sigma_+ \cup_\Gamma \Sigma_-$, where
$\Sigma_\pm \subset \Sigma$ are each (not necessarily connected) 
compact surfaces with a common boundary~$\Gamma$.  By an old result of
Lutz \cite{Lutz:77}, the $3$-manifold $S^1 \times \Sigma$
admits a unique isotopy class of $S^1$-invariant contact structures $\xi_\Gamma$
such that the loops $S^1 \times \{z\}$ are positively/negatively transverse
to $\xi_\Gamma$ for $z \in \mathring{\Sigma}_\pm$ and tangent to $\xi_\Gamma$
for $z \in \Gamma$.  Now for each $k \in \NN$, define
$$
(V_k,\xi_k) := (S^1 \times \Sigma,\xi_\Gamma)
$$
where $\Sigma = \Sigma_+ \cup_\Gamma \Sigma_-$ is chosen such that
$\Gamma$ has $k$ connected components, $\Sigma_-$ is connected with genus 
zero, and $\Sigma_+$ is connected with positive genus 
(see Figure~\ref{fig:torsion}).

\begin{thm}
\label{thm:torsion}
The contact manifolds $(V_k,\xi_k)$ do not admit any symplectic fillings.
Moreover, if $k > \ell$, then there exists no exact symplectic cobordism from
$(V_k,\xi_k)$ to $(V_\ell,\xi_\ell)$.
\end{thm}

For these examples, one can use explicit constructions 
from \cites{Wendl:cobordisms,Avdek:sums} to show that non-exact cobordisms
from $(V_k,\xi_k)$ to $(V_\ell,\xi_\ell)$ do exist, and so do exact cobordisms
from $(V_\ell,\xi_\ell)$ to $(V_k,\xi_k)$, thus both the directionality of
the cobordism relation and the distinction between exact and non-exact are
crucial.  The proof of the theorem, due to the author with Latschev and 
Hutchings \cite{LatschevWendl}, uses a numerical contact invariant
based on the full SFT algebra---in particular, the curves that cause this
phenomenon have multiple positive ends and are thus not seen by contact
homology.  We will introduce the relevant numerical invariant in Lecture~\ref{lec:SFT}
and compute it for these examples in Lecture~\ref{lec:torsion}.

\chapter{Basics on holomorphic curves}
\label{lec:local}

\minitoc

\vspace{12pt}

In this lecture we begin studying the analysis of $J$-holomorphic curves.
The coverage will necessarily be a bit sparse in some places, but more 
detailed proofs of everything in this lecture 
can be found in \cite{Wendl:lecturesV33}.

\section{Linearized Cauchy-Riemann operators}
\label{sec:linearCR}

In order to motivate the study of linear Cauchy-Riemann type operators,
we begin with a formal discussion of the nonlinear Cauchy-Riemann equation
and its linearization.

Fix a Riemann surface $(\Sigma,j)$ and almost complex manifold $(W,J)$, and
suppose that we wish to understand the structure of some space of the form
\begin{equation}
\label{eqn:spaceOfCurves}
\left\{ u : \Sigma \to W \ |\ Tu \circ j = J \circ Tu 
\text{ plus further conditions} \right\},
\end{equation}
where the ``further conditions'' 
(which we will for now leave unspecified) may impose
constraints on e.g.~the regularity of $u$, as well as its boundary and/or asymptotic
behavior.  The standard approach in global analysis can be summarized as 
follows:
\begin{itemize}
\item[\textsl{Step~1:}]
Construct a smooth Banach manifold $\bB$ of maps $u : \Sigma \to W$ such that
all the solutions we're interested in will be elements of~$\bB$.  The
tangent spaces $T_u \bB$ are then Banach spaces of sections of~$u^*TW$.
\item[\textsl{Step~2:}]
Construct a smooth Banach space bundle $\eE \to \bB$ such that for each
$u \in \bB$, the fiber $\eE_u$ is a Banach space of sections of the vector
bundle
$$
\overline{\Hom}_\CC(T\Sigma,u^*TW) \to \Sigma
$$
of complex-antilinear bundle maps $(T\Sigma,j) \to (u^*TW,J)$.  Since our
purpose is to study a first-order PDE, we need the sections in $\eE_u$ to
be ``one step less regular'' than the maps in $\bB$, e.g.~if $\bB$ consists
of maps of Sobolev class $W^{k,p}$, then the sections in $\eE_u$ should be 
of class $W^{k-1,p}$.
\item[\textsl{Step~3:}]
Show that
$$
\dbar_J : \bB \to \eE : u \mapsto Tu + J(u) \circ Tu \circ j
$$
defines a smooth section of $\eE \to \bB$, whose zero set is precisely the
space of solutions \eqref{eqn:spaceOfCurves}.
\item[\textsl{Step~4:}]
Show that under suitable assumptions (e.g.~on regularity and asymptotic
behavior), one can arrange such that for every
$u \in \dbar_J^{-1}(0)$, the \defin{linearization} of $\dbar_J$,
$$
D\dbar_J(u) : T_u \bB \to \eE_u
$$
is a Fredholm operator and is generically surjective.
(In geometric terms, this would mean that $\dbar_J$ is 
\emph{transverse to the zero section}.)
\item[\textsl{Step~5:}]
Using the implicit function theorem 
in Banach spaces (see \cite{Lang:analysis}), the surjectivity of
$D\dbar_J(u)$ implies
that $\dbar_J^{-1}(0)$ is a smooth finite-dimensional manifold,
with its tangent space at each $u \in \dbar_J^{-1}(0)$ canonically
identified with $\ker D\dbar_J(u)$, hence the 
dimension of $\dbar_J^{-1}(0)$ near~$u$ equals the Fredholm index
of~$D\dbar_J(u)$.
\end{itemize}

Without worrying about the fact that these are actually not Banach spaces, and some Sobolev completion is needed, let us assume, for simplicity, that the bundle $\mathcal{E}\rightarrow \mathcal{B}$ has as base the space $\mathcal{B}=C^{\infty}(\Sigma,W)$ and the fiber over $u \in \mathcal{B}$ is given by 
$\mathcal{E}_u=\overline{\Hom}_\mathbb{C}(T\Sigma,u^*TW)$. The linearization of the section $\overline{\partial}_J$ 
at a point $u \in \dbar_J^{-1}(0)$ should then take the form
$$
\mathbf{D}_u : \Gamma(u^*TW) \to \Omega^{0,1}(\Sigma,u^*TW),
$$
where the right hand side denotes the space of $u^*TW$-valued $(0,1)$-forms
on~$\Sigma$, or equivalently, smooth sections of
$\overline{\Hom}_\CC(T\Sigma,u^*TW)= T^{0,1}\Sigma \otimes_\CC u^*TW$, 
where $T^{0,1}\Sigma$ denotes the $(0,1)$-part of the complexified cotangent bundle.

In order to compute such a linearization, we need to make a choice of ``connection'' on the bundle $\mathcal{E}$. Choose a connection $\nabla$ on $W$, and recall the fact that this naturally induces a connection on the bundles $T^{0,1}\Sigma \otimes_\CC u^*TW$ and $\End(u^*TW)$ by setting $\nabla (\alpha \otimes s)=\alpha \otimes \nabla s$ and $(\nabla J)s=\nabla(J s)-J \nabla s$, for $s \in \Gamma(u^*TW)$, $J \in \End(u^*TW)$ and $\alpha \in \Gamma(T^{0,1}\Sigma)$. We shall make the ansatz that for any smooth $1$-parameter family of maps $u_\rho : \Sigma \to W$ for
$\rho \in (-\epsilon,\epsilon)$ and a section $\eta_\rho \in \mathcal{E}_{u_\rho}$ along the path (i.e a section of the pullback bundle of $\mathcal{E}$ under the map $(-\epsilon,\epsilon)\rightarrow \mathcal{B}$ mapping $\rho$ to $u_\rho$), the connection takes the form $$(\nabla_\rho \eta_\rho)X=\nabla_\rho(\eta_\rho(X)),$$ for $X \in T\Sigma$, where this expression should be interpreted as the pullback connection under the map displayed above. The tensorial property of connections implies that $\nabla_\rho \eta_\rho$ does not depend on the connection at the values $\rho$ for which $\eta_\rho=0$.

Given $u \in \overline{\partial}_J^{-1}(0)$ and $\eta$ in $T_u \mathcal{B}=\Gamma(u^*TW)$, take a one-parameter family $u_\rho \in \mathcal{B}$ with $u_0=u$ and $\left. \partial_\rho u_\rho \right|_{\rho=0}=\eta$. 
We then have that 
$$ 
\mathbf{D}_u \eta = \left. \nabla_\rho \left( \dbar_J(u_\rho) \right) \right|_{\rho=0} = \left. \nabla_\rho(Tu_\rho + J(u_\rho) \circ Tu_\rho \circ j)\right|_{\rho=0}.
$$
Since $\dbar_J u = 0$, this is independent of the connection, and we may therefore choose 
$\nabla$ to be symmetric.

Fix a point $z \in \Sigma$ and choose local holomorphic coordinates $s+it$ around it. The symmetry of the connection implies $\left. \nabla_\rho \partial_s u_\rho \right|_{\rho=0} = \left. \nabla_s \partial_\rho u_\rho \right|_{\rho=0}=\nabla_s \eta$, and similarly for the variable $t$. Observing also that 
$\nabla_\eta J=\nabla_{\left. \partial_\rho u_\rho \right|_{\rho=0}} J 
= \left. \nabla_\rho \left( J(u_\rho) \right) \right|_{\rho=0}$, and using the above ansatz, we obtain 
$$
(\mathbf{D}_u \eta)\partial_s= \left. \nabla_\rho(\partial_s u_\rho + J(u_\rho) \circ \partial_t u_\rho)\right|_{\rho=0}=  \nabla_s \eta + J(u)\nabla_t\eta + (\nabla_\eta J)\partial_t u 
$$
Since $\mathbf{D}_u\eta$ is an antilinear map, and $\partial_t=j\partial_s$, it is therefore determined by its action on $\partial_s$. One can check that the operator on the right hand side below is also antilinear, and thus removing the~$\partial_s$, we obtain 
\begin{equation}
\label{eqn:Du}
\boxed{\mathbf{D}_u \eta = \nabla \eta + J(u) \circ \nabla \eta \circ j +
(\nabla_\eta J) \circ Tu \circ j.}
\end{equation}

\begin{defn}
Fix a complex vector bundle $E$ over a Riemann surface $(\Sigma,j)$.
A (real) linear \defin{Cauchy-Riemann type operator} on $E$ is a real-linear 
first-order differential operator
$$
\mathbf{D} : \Gamma(E) \to \Omega^{0,1}(\Sigma,E)
$$
such that for every $f \in C^\infty(\Sigma,\RR)$ and 
$\eta \in \Gamma(E)$,
\begin{equation}
\label{eqn:Leibniz}
\mathbf{D} (f \eta) = (\dbar f) \eta + f \mathbf{D}\eta,
\end{equation}
where $\dbar f$ denotes the complex-valued $(0,1)$-form
$df + i \, df \circ j$.
\end{defn}

Observe that $\mathbf{D}$ is complex linear if and only if the
Leibniz rule \eqref{eqn:Leibniz} also holds for all
smooth complex-valued functions $f$, not just real-valued.
It is a standard result in complex geometry
that choosing a complex-linear Cauchy-Riemann type 
operator $\mathbf{D}$ on $E$ is equivalent to endowing it with the structure 
of a \emph{holomorphic} vector bundle, where local sections $\eta$ are defined 
to be holomorphic if and only if $\mathbf{D}\eta = 0$.  Indeed, every
holomorphic bundle comes with a canonical Cauchy-Riemann operator that is
expressed as $\dbar$ in holomorphic trivializations, and in the other 
direction, the equivalence follows from a local existence result for solutions
to the equation $\mathbf{D}\eta = 0$, proved in \S\ref{sec:similarity} 
below.\footnote{This statement about the existence of holomorphic vector 
bundle structures is true when the base is a Riemann surface, but not if 
it is a higher-dimensional complex manifold. In higher
dimensions there are obstructions, see e.g.~\cite{Kobayashi}.}

\begin{exercise}
\label{EX:linearCR}
If $\mathbf{D}$ is a linear Cauchy-Riemann type operator on~$E$, 
prove that every other such operator is of the form $\mathbf{D} + A$ where
$A : E \to \overline{\Hom}_\CC(T\Sigma,E)$ is a smooth linear bundle map.
Using this, show that in suitable local trivializations over a subset 
$\uU \subset \Sigma$ identified biholomorphically with an open set in~$\CC$,
every Cauchy-Riemann type operator $\mathbf{D}$ takes the form
$$
\mathbf{D} = \dbar + A : C^\infty(\uU,\CC^m) \to C^\infty(\uU,\CC^m),
$$
where $\dbar = \p_s + i \p_t$ in complex coordinates $z = s + it$ and
$A \in C^\infty(\uU,\End_\RR(\CC^m))$.
\end{exercise}

\begin{exercise}
Verify that the linearized operator $\mathbf{D}_u$ of
\eqref{eqn:Du} is a real-linear Cauchy-Riemann type operator.
\end{exercise}

\section{Some useful Sobolev inequalities}
\label{sec:Sobolev}

In this section, we review a few general properties of Sobolev spaces that
are essential for applications in nonlinear analysis.  The results stated
here are explained in more detail in Appendix~\ref{app:Sobolev}.

Throughout this section we consider functions with values in $\CC$ unless
otherwise specified, and defined on an
open domain~$\uU$ in either $\RR^n$ or a quotient of~$\RR^n$ on which the
Lebesgue measure is well defined.  Certain regularity assumptions must 
generally be placed on the boundary of $\overline{\uU}$ in order for all
the results stated below to hold; we will ignore this detail except to
mention that the necessary assumptions are satisfied for the two classes
of domains that we are most interested in, which are
\begin{equation*}
\begin{split}
\uU &= \intDD \subset \CC, \\
\uU &= (0,L) \times S^1 \subset \CC / \ZZ, \qquad 0 < L \le \infty.
\end{split}
\end{equation*}
Here $\DD$ denotes the closed unit disk and $\intDD$ is its
interior.  Certain results will be specified to hold only for \emph{bounded}
domains, which means in practice that they hold on $\intDD$
and $(0,L) \times S^1$ for any $L > 0$, but not on
$(0,\infty) \times S^1$.

Recall that for $p \in [1,\infty)$ we define the $L^p$ norm of a measurable function 
$f : \uU \rightarrow \RR^m$ 
to be
$$ \|f\|_{L^p} = \left(
\int_\uU |f|^p\right)^{1/p}.$$
For the space $L^\infty$ we define the norm to be the essential supremum of $f$ over $\uU$.\\
Denote by $$C_0^\infty(\uU) \subset C^\infty(\uU)$$ the set of smooth functions with compact support in $\uU$. 
We say a function $f$ has a \textbf{weak $j$-th partial derivative} $g$ if the 
\textit{integration by parts} formula holds for all $\varphi \in C_0^\infty(\uU)$:
$$
\int_\uU g \varphi = - \int_\uU f \, \p_j\varphi.
$$
Equivalently, this means that $g$ is a partial derivative of $f$
\defin{in the sense of distributions} (see e.g.~\cite{LiebLoss}).
Higher order weak partial derivatives are defined similarly: 
recall that for a multiindex $\alpha=(i_1,...i_n)$ we denote 
$$
\p^\alpha f = \frac{\p^{|\alpha|} f}{\p x_1^{i_1}\ldots\p x_n^{i_n}},
$$
where $|\alpha| := \sum_j i_j$. We then write $\p^\alpha f = g$ if
for all $\varphi \in C_0^\infty(\uU)$,
$$\int_\uU g\varphi = (-1)^{|\alpha|} \int_\uU f \, \p^\alpha \varphi.$$
Now we may define $W^{k,p}(\uU)$ to be the set of functions on $\uU$ with weak partial derivatives up to order $k$ lying in $L^p$, and define the norm of such a function by:
$$\|f\|_{W^{k,p}}= \sum_{|\alpha| \leq k} \|\p^\alpha f\|_{L^p}.$$
As $W^{k,p}(\uU)$ can be regarded as a subset of a $k$-fold product of $L^p(\uU)$, it is a Banach space, and it is reflexive and separable for $1<p<\infty$.

While the Sobolev spaces $W^{k,p}(\uU)$ are generally defined on \emph{open}
domains, we often consider the closure $\overline{\uU}$ as the domain for
spaces of differentiable functions 
$C^k(\overline{\uU})$ and $C^\infty(\overline{\uU})$.  For
instance, $C^k(\overline{\uU})$ is the Banach space of $k$-times differentiable
functions on $\uU$ whose derivatives up to order $k$ are bounded and uniformly
continuous on~$\uU$; note that uniform continuity implies the existence of
continuous extensions to the closure~$\overline{\uU}$.
Given suitable regularity assumptions for the boundary
of~$\overline{\uU}$, one can show (with some effort) that $C^k(\overline{\uU})$
is precisely the set 
of functions which admit $k$-times differentiable extensions to some open
set containing~$\overline{\uU}$.

The following two results are special cases of the more general
Theorems~\ref{thm:SobolevEmbedding} and~\ref{thm:Rellich} in
Appendix~\ref{app:Sobolev}, proofs of which may be found
e.g.~in \cite{AdamsFournier}.

\begin{prop}[Sobolev embedding theorem]
\label{prop:SobolevEmb}
Assume $1 \le p < \infty$, $kp > n$ and $d \ge 0$ is an integer.
Then there exists a continuous inclusion
$$
W^{k+d,p}(\uU) \hookrightarrow C^d(\overline{\uU}),
$$
which is compact if $\uU$ is bounded.   \qed
\end{prop}

\begin{prop}[Rellich-Kondrachov compactness theorem]
\label{prop:Rellich}
If $1 \le p < \infty$ and $\uU$ is bounded, then the natural inclusion
$$
W^{k+1,p}(\uU) \hookrightarrow W^{k,p}(\uU)
$$
is compact.    \qed
\end{prop}

\begin{exercise}
Show that Proposition~\ref{prop:Rellich} fails in general for unbounded
domains, e.g.~for~$\RR$.
\end{exercise}

The next three results for the case $kp > n$ are proved in
\S\ref{sec:SobolevProducts} as corollaries of the Sobolev embedding theorem.

\begin{prop}[Banach algebra property]
\label{prop:BanachAlg}
Suppose $1 \le p < \infty$, $kp > n$ and $0 \le m \le k$.  Then 
the product pairing $(f,g) \mapsto fg$ defines a continuous bilinear map
$$
W^{k,p}(\uU) \times W^{m,p}(\uU) \to W^{m,p}(\uU).
$$
In particular, $W^{k,p}(\uU)$ is a Banach algebra.
\qed
\end{prop}

The continuity statements above translate into inequalities between the norms 
in the respective spaces. For example, continuous inclusions
$W^{k+d,p} \hookrightarrow C^d$ and $W^{k+1,p} \hookrightarrow W^{k,p}$
respectively imply that
\begin{equation*}
\begin{split}
\|f\|_{C^d} &\leq c \|f\|_{W^{k+d,p}} \\
\|f\|_{W^{k,p}} &\leq c \|f\|_{W^{k+1,p}} 
\end{split}
\end{equation*}
for some constants $c > 0$ which may depend on $d$, $k$, $p$ or $\uU$, but not $f$. 
Similarly, the Banach algebra property implies
$$
\|fg\|_{W^{m,p}} \leq c \|f\|_{W^{k,p}} \|g\|_{W^{m,p}},
$$
where again, the constant $c$ is independent of $g$ and $f$.

We state the next result only for the case of bounded domains; it does have an
extension to unbounded domains, but the statement becomes more 
complicated (cf.~Theorem~\ref{thm:precomposition}).  
Given an open set $\Omega \subset \RR^n$, we denote
$$
W^{k,p}(\uU,\Omega) := \left\{ u \in W^{k,p}(\uU,\RR^n)\ \Big|\ 
\overline{u(\uU)} \subset \Omega \right\}.
$$
Note that this is an open subset if $kp > n$, due to the
Sobolev embedding theorem.

\begin{prop}[$C^k$-continuity property]
\label{prop:CkContinuity}
Assume $1 \le p < \infty$, $kp > n$, $\uU$ is bounded and $\Omega \subset \RR^n$
is an open set.  Then the map
$$
C^k(\Omega,\RR^N) \times W^{k,p}(\uU,\Omega) \to W^{k,p}(\uU,\RR^N) :
(f,u) \mapsto f \circ u
$$
is well defined and continuous.
\qed
\end{prop}

\begin{remark}
Though we will not yet use it in this lecture, 
Propositions~\ref{prop:SobolevEmb}, \ref{prop:BanachAlg} 
and~\ref{prop:CkContinuity} are the essential conditions needed in order to
define smooth Banach manifold structures on spaces of $W^{k,p}$-smooth
maps from one manifold to another, cf.~\cites{Eliasson,Palais:global}.
This only works under the condition $kp > n$, as the smooth category is
not well equipped to deal with discontinuous maps!
\end{remark}

The following rescaling result will be needed for nonlinear regularity
arguments; see Theorem~\ref{thm:SobolevRescaling} in Appendix~\ref{app:Sobolev}
for a proof.

\begin{prop}
\label{prop:SobolevRescaling}
Assume $p \in [1,\infty)$ and $k \in \NN$ satisfy $kp > n$, let
$\intDD^n$ denote the open unit ball in~$\RR^n$, and
for each $f \in W^{k,p}(\intDD^n)$ and $\epsilon \in (0,1]$, define
$f_\epsilon \in W^{k,p}(\intDD^n)$ by
$$
f_\epsilon(x) := f(\epsilon x).
$$
Then there exist constants $C > 0$ and $r > 0$ such that
for every $f \in W^{k,p}(\intDD^n)$,
$$
\| f_\epsilon - f(0) \|_{W^{k,p}(\intDD^n)} \le C \epsilon^r 
\| f - f(0) \|_{W^{k,p}(\intDD^n)}
\quad\text{ for all $\epsilon \in (0,1]$}.
$$
\qed
\end{prop}

\begin{exercise}
Working on a $2$-dimensional domain with $kp > 2$, 
prove directly that for any multiindex $\alpha$ of positive degree~$k$,
$$
\| \p^\alpha f_\epsilon \|_{L^p(\intDD)} \le
\epsilon^{k - 2/p} \| \p^\alpha f \|_{L^p(\intDD)}
$$
for $f \in W^{k,p}(\intDD)$.  
Find examples (e.g.~in $W^{1,2}(\intDD)$) to show that no estimate of the form
$$
\| \p^\alpha f_\epsilon \|_{L^p(\intDD)} \le
C_\epsilon \| f - f(0) \|_{W^{k,p}(\intDD)}
$$
with $\lim_{\epsilon \to 0^+} C_\epsilon = 0$ is possible when $kp \le 2$.
\end{exercise}

\section{The fundamental elliptic estimate}
\label{sec:estimate}

We will make considerable use of the fact that the linear first-order
differential operator
$$
\dbar := \p_s + i \p_t : C^\infty(\CC,\CC) \to C^\infty(\CC,\CC)
$$
is \defin{elliptic}.  There is no need to discuss here precisely what
ellipticity means in full generality (see 
\cite{Wendl:lecturesV33}*{\S 2.B} if you're curious about this); in
practice, the main consequence is the following pair of analytical
results.

\begin{thm}
\label{thm:rightInverse}
If $1 < p < \infty$, then $\dbar : W^{1,p}(\intDD) \to L^p(\intDD)$ admits
a bounded right inverse $T : L^p(\intDD) \to W^{1,p}(\intDD)$.
\end{thm}

\begin{thm}
\label{thm:estimate}
If $1 < p < \infty$ and $k \in \NN$, then there exists a constant $c > 0$ 
such that for all $f \in W^{k,p}_0(\intDD)$,
$$
\| f \|_{W^{k,p}} \le c \| \dbar f \|_{W^{k-1,p}}.
$$
\end{thm}
Here $W^{k,p}_0(\intDD)$ denotes the $W^{k,p}$-closure of $C_0^\infty(\intDD)$,
the latter being space of smooth functions on $\intDD$ with compact support.

The complete proofs of the two theorems above are rather lengthy, and we
shall refer to \cite{Wendl:lecturesV33}*{\S 2.6 and~2.A} for the details,
but we can at least explain why they hold in the case $p=2$.  First, it is
straightforward to show that the function $K \in L^1_\loc(\CC)$ defined by
$$
K(z) = \frac{1}{2\pi z}
$$
is a \defin{fundamental solution} for the equation $\dbar u = f$, meaning it
satisfies
$$
\dbar K = \delta
$$
in the sense of distributions, where $\delta$ denotes the Dirac
$\delta$-function.  Hence for any $f \in C_0^\infty(\CC)$, one finds a
smooth solution $u : \CC \to \CC$ to the equation $\dbar u = f$ as the
convolution
$$
u(z) = (K*f)(z) := \int_{\CC} K(z - \zeta) f(\zeta) \, d\mu(\zeta),
$$
where $d\mu(\zeta)$ denotes the Lebesgue measure with respect to the
variable $\zeta \in \CC$.  It is not hard to show from this formula that 
whenever $f \in C_0^\infty$, $K*f$ has decaying behavior at infinity
(see \cite{Wendl:lecturesV33}*{Lemma~2.6.13}).  Thus if
$u \in C_0^\infty$ and $\dbar u = f$, it follows that $u - K*f$ is a
holomorphic function on $\CC$ that decays at infinity, hence
$u \equiv K*f$.  Since $C_0^\infty(\intDD)$ is dense in $L^p(\intDD)$ for
all $p < \infty$, Theorem~\ref{thm:rightInverse}
now follows from the claim that for all
$f \in C_0^\infty(\intDD)$, there exist estimates of the form
\begin{equation}
\label{eqn:twoEstimates}
\| K*f \|_{L^p(\intDD)} \le c \| f \|_{L^p(\intDD)}, \qquad
\| \p_j(K*f) \|_{L^p(\intDD)} \le c \| f \|_{L^p(\intDD)},
\end{equation}
with $\p_j = \p_s$ or $\p_t$ for $j=1,2$ respectively,
and the constant $c > 0$ independent of~$f$.

\begin{exercise}
Use Theorem~\ref{thm:rightInverse} and the remarks above to prove 
Theorem~\ref{thm:estimate} for the case $k=1$ with 
$f \in C_0^\infty(\intDD)$, then extend it to $f \in W^{1,p}_0(\intDD)$
by a density argument.  Then extend it to the
general case by differentiating both $f$ and~$\dbar f$.
\end{exercise}

The first estimate in \eqref{eqn:twoEstimates}
is not too hard if you remember your introductory measure theory class: 
it follows from a general ``potential inequality''
for convolution operators (see \cite{Wendl:lecturesV33}*{Lemma~2.6.10}),
similar to Young's inequality, the key points being that $K$ is locally of
class $L^1$ and $\intDD$ has finite measure.  For the second inequality,
observe that $\dbar (K*f) = f$, and the rest of the first derivative
of $K*f$ is determined by $\p (K*f)$, where
$$
\p := \p_s - i \p_t.
$$
Differentiating $K$ in the sense of distributions provides a formula for
$\p (K*f)$ as a principal value integral, namely
$$
\p (K*f)(z) = - \frac{1}{\pi} \lim_{\epsilon \to 0^+}
\int_{|\zeta - z| \ge \epsilon} \frac{f(\zeta)}{(z - \zeta)^2} \, d\mu(\zeta).
$$
This is a so-called \defin{singular integral operator}: it is similar to
our previous convolution operator, but more difficult to handle because the
kernel $\frac{1}{z^2}$ is not of class $L^1_\loc$ on~$\CC$.  The proof of
the estimate $\| \p (K*f) \|_{L^p} \le c \| f \|_{L^p}$ for all
$f \in C_0^\infty(\intDD)$ follows from a rather difficult general estimate 
on singular integral operators, known as the \emph{Calder\'on-Zygmund
inequality}, cf.~\cite{Wendl:lecturesV33}*{\S 2.A} and the references
therein.  The good news however is that the first step in that proof is
not hard: that is the case $p=2$.

As is the case for all elliptic operators with constant coefficients, the
$L^2$-estimate on the fundamental solution of $\dbar$ admits an easy
proof using Fourier transforms:

\begin{prop}
\label{prop:L2estimate}
For all $f \in C_0^\infty(\CC)$, we have
$\| \p(K*f) \|_{L^2} = \| f \|_{L^2}$.
\end{prop}
\begin{proof}
A sufficiently nice function $u : \CC \to \CC$ is related to its Fourier
transform $\hat{u} : \CC \to \CC$ by
$$
u(z) = \int_\CC \hat{u}(\zeta) e^{2\pi i (z \cdot \zeta)} \, d\mu(\zeta)
$$
and thus satisfies the identities
\begin{equation*}
\begin{split}
\widehat{\dbar u}(\zeta) &= 2 \pi i \zeta \hat{u}(\zeta), \\
\widehat{\p u}(\zeta) &= 2\pi i \overline{\zeta} \hat{u}(\zeta).
\end{split}
\end{equation*}
Since $u = K * f$ we have $\hat{u}=\hat{K}\hat{f}$, and since $\dbar K = \delta$, we have $2\pi i \zeta \hat{K}=1$.
Hence we may apply Plancharel's theorem to deduce
\begin{align*}
\|\p (K * f) \|_{L^2} &= \|\p u\|_{L^2}= \|\widehat{\p u}\|_{L^2} 
=\| 2 \pi i \overline{\zeta}\hat{u}\|_{L^2}=\| 2 \pi i \overline{\zeta}\hat{K}\hat{f}\|_{L^2} \\
&=\left\| \frac{\overline{\zeta}}{\zeta} 2 \pi i \zeta \hat{K} \hat{f} \right\|_{L^2}= \left\| \frac{\overline{\zeta}}{\zeta}\hat{f}\right\|_{L^2} 
=\|\hat{f}\|_{L^2}=\|f\|_{L^2}.
\end{align*}
\end{proof}

\section{Regularity}
\label{sec:regularity}

We will now use the estimate $\| u \|_{W^{k,p}} \le c \| \dbar u \|_{W^{k-1,p}}$
from the previous section to prove three types of results about solutions
to Cauchy-Riemann type equations:
\begin{enumerate}
\item All solutions of reasonable Sobolev-type regularity are smooth.
\item Any collection of solutions satisfying uniform bounds in certain
Sobolev norms also locally satisfy uniform $C^\infty$-bounds.
\item All reasonable Sobolev-type topologies on spaces of solutions are
(locally) equivalent to the $C^\infty$-topology.
\end{enumerate}

In the following,
$$
\DD_r \subset \CC
$$
denotes the closed disk of radius $r > 0$, and $\intDD_r$ denotes its
interior.  Note that functions of class $C^\infty(\DD_r)$ are assumed to be
smooth up to the boundary (or equivalently, on some open neighborhood
of $\DD_r$ in~$\CC$), not just on~$\intDD_r$.

\subsection{The linear case}

Recall from Exercise~\ref{EX:linearCR} that every linear Cauchy-Riemann type
operator on a vector bundle of complex rank $n$
locally takes the form $\dbar + A$, where $\dbar = \p_s + i \p_t$,
and $A$ is a smooth function with values in $\End_\RR(\CC^n)$.
Using the Sobolev embedding theorem, the following result implies by induction
that solutions $u \in W^{1,p}$ to the equation $(\dbar + A) u = 0$ 
are always smooth.

\begin{thm}[Linear regularity]
\label{thm:regularityLinear}
Assume $1 < p < \infty$ and $k,m \in \NN$.
\begin{enumerate}
\item \label{item:smoothLinear}
If $u \in W^{k,p}(\intDD)$ satisfies $\dbar u \in W^{m,p}(\intDD)$,
then $u$ is in $W^{m+1,p}$ on every compact subset of~$\intDD$.
\item Suppose $f_\nu \in W^{m,p}(\intDD)$ is a sequence converging in the 
$W^{m,p}$-topology to $f \in W^{m,p}(\intDD)$ as $\nu \to \infty$, and 
$u_\nu \in W^{k,p}(\intDD)$ is a sequence with $\dbar u_\nu = f_\nu$.
\begin{enumerate}
\item \label{item:unifLinear}
If there exist uniform bounds on $\| u_\nu \|_{W^{k,p}}$ and
$\| f_\nu \|_{W^{m,p}}$ over $\intDD$ as $\nu \to \infty$, then
$\| u_\nu \|_{W^{m+1,p}}$ is also uniformly bounded on every compact
subset of~$\intDD$.
\item \label{item:convLinear}
If the sequence $u_\nu$ is $W^{k,p}$-convergent on $\intDD$ to a function
$u \in W^{k,p}(\intDD)$ satisfying $\dbar u = f$, then it is also 
$W^{m+1,p}$-convergent on every compact subset of~$\intDD$.
\end{enumerate}
\end{enumerate}
\end{thm}
\begin{proof}
We begin by proving statement~\eqref{item:unifLinear}, assuming that 
statement~\eqref{item:smoothLinear} is already known, hence 
$u_\nu \in W^{m+1,p}_\loc(\intDD)$ since $f_\nu \in W^{m,p}(\intDD)$.
Assume $m = k$, since there is otherwise nothing to prove. Then by induction, 
it suffices to show that uniform bounds on $\| u_\nu \|_{W^{k,p}(\intDD)}$
and $\| f_\nu \|_{W^{k,p}(\intDD)}$ imply a uniform bound on
$\| u_\nu \|_{W^{k+1,p}(\intDD_r)}$ for any given $r < 1$; equivalently, this 
would mean there is a uniform bound on 
$\| \p_j u_\nu \|_{W^{k,p}(\intDD_r)}$ for $j = 1,2$. 
In order to apply the elliptic estimate, we need to work with functions with 
compact support in~$\intDD$, thus choose a smooth bump function
$$
\beta \in C_0^\infty(\intDD,[0,1])
$$
that satisfies $\beta|_{\DD_r} \equiv 1$.  We then have 
$\beta \, \p_j u_\nu \in C_0^\infty(\intDD)$, so by Theorem~\ref{thm:estimate},
\begin{equation}
\begin{split}
\label{eqn:boundThis}
\| \p_j u_\nu \|_{W^{k,p}(\intDD_r)} &\le
\| \beta \, \p_j u_\nu \|_{W^{k,p}(\intDD)} \le
c \left\| \dbar\left( \beta \, \p_j u_\nu \right) \right\|_{W^{k-1,p}(\intDD)} \\
&\le c \| (\dbar \beta) (\p_j u_\nu) \|_{W^{k-1,p}} +
c \| \beta \, \dbar(\p_j u_\nu) \|_{W^{k-1,p}}.
\end{split}
\end{equation}
The first term on the right hand side is uniformly bounded since
$\dbar\beta$ is smooth and $\| u_\nu \|_{W^{k,p}}$ is uniformly bounded.
To control the second term, we differentiate the equation
$\dbar u_\nu = f_\nu$, giving
$$
\dbar(\p_j u_\nu) = \p_j f_\nu.
$$
This also has a uniformly bounded $W^{k-1,p}$-norm since $\| f_\nu \|_{W^{k,p}}$
is uniformly bounded.  Since $\beta$ is smooth, this bounds the second term 
on the right hand side of \eqref{eqn:boundThis} as $\nu \to \infty$,
and we are done.

Statement~\eqref{item:convLinear} follows by a similar argument bounding
$\| \p_j (u - u_\nu) \|_{W^{k,p}(\intDD_r)}$ in terms of
$\| u - u_\nu \|_{W^{k,p}(\intDD)}$ and $\| f - f_\nu \|_{W^{k,p}(\intDD)}$;
we leave the details as an exercise.

Lastly, we prove statement~\eqref{item:smoothLinear},
where again it suffices to assume $\dbar u = f \in W^{k,p}(\intDD)$ and show
that $u|_{\intDD_r} \in W^{k+1,p}(\intDD_r)$ for some $r < 1$.
The idea is to use the
same argument that was used for statement~\eqref{item:unifLinear}, but with
the partial derivatives $\p_j u$ replaced by the \defin{difference quotients}
$$
D_j^h u(z) := \frac{u(z + h e_j) - u(z)}{h},  \qquad j=1,2,
$$
where $e_1 := \p_s$, $e_2 := \p_t$, and the role of the index $\nu \to \infty$
is now played by the parameter $h \in \RR \setminus \{0\}$ approaching~$0$.
Note that if $u \in W^{k,p}(\intDD)$, then $\beta\, D_j^h u$ is a well-defined
function on $\intDD$ for all $|h| \ne 0$ sufficiently small and belongs to
$W^{k,p}_0(\intDD)$.  The analogue of \eqref{eqn:boundThis} in this context
is then
\begin{equation*}
\begin{split}
\| D_j^h u \|_{W^{k,p}(\intDD_r)} &\le
\| \beta D_j^h u \|_{W^{k,p}(\intDD)} \le
c \left\| \dbar\left( \beta \, D_j^h u \right) \right\|_{W^{k-1,p}(\intDD)} \\
&\le c \| (\dbar \beta) (D_j^h u) \|_{W^{k-1,p}} +
c \| \beta \, \dbar(D_j^h u) \|_{W^{k-1,p}}.
\end{split}
\end{equation*}
The first term is bounded independently of $h$ since 
$\p_j u \in W^{k-1,p}(\intDD)$, 
implying a uniform $W^{k-1,p}$-bound on $D_j^h u$ as $h \to 0$.
To control the second term, we can apply the operator $D_j^h$ to the
equation $\dbar u = f$, giving
$$
\dbar (D_j^h u) = D_j^h (\dbar u) = D_j^h f.
$$
This satisfies a $W^{k-1,p}$-bound that is uniform in~$h$ since
$\p_j f \in W^{k-1,p}(\intDD)$, so we conclude that for all
$|h|$ sufficiently small,
$$
\| D_j^h u \|_{W^{k,p}(\intDD_r)} \le c
$$
for some constant $c > 0$ that does not change as $h \to 0$.  By a
standard application of the Banach-Alaoglu theorem 
(cf.~\cite{Evans}*{\S 5.8.2}), this implies the
existence of a sequence $h_\nu \to 0$ for which $D_j^{h_\nu} u$ is
$W^{k,p}$-convergent on~$\intDD_r$, and its limit is necessarily
$\p_j u$, which therefore belongs to~$W^{k,p}$.  Indeed, if $k=0$,
the uniform $L^p$-bound on $D_j^{h_\nu} u$ over $\intDD_r$ for any
sequence $h_\nu \to 0$ gives rise to a weakly
$L^p$-convergent subsequence via the Banach-Alaoglu theorem.  The limit
of this subsequence belongs to $L^p(\intDD_r)$, and it is straightforward
to show using the definition of weak derivatives that this limit
is~$\p_j u$.  One finds the same result for any $k \in \NN$ by applying
this argument to higher-order derivatives of $\p_j u$.  The conclusion
is that $u$ is in $W^{k+1,p}$ on $\intDD_r$, since $u$ and both of its
first partial derivatives belong to~$W^{k,p}$.  
\end{proof}

\begin{exercise}
Show that all three parts of Theorem~\ref{thm:regularityLinear} continue to 
hold if the operator $\dbar$ is replaced by $\dbar + A$ or $\dbar + A_\nu$,
where $A, A_\nu \in C^\infty(\DD,\End_\RR(\CC^n))$ with $A_\nu \to A$
in $C^\infty$ as $\nu \to \infty$.
\end{exercise}

\begin{exercise}
\label{EX:BRI}
Use Theorem~\ref{thm:regularityLinear}\eqref{item:smoothLinear} to extend
Theorem~\ref{thm:rightInverse} to the existence of a bounded right inverse for
$$
\dbar : W^{k,p}(\intDD) \to W^{k-1,p}(\intDD).
$$
\textsl{Hint: For any $R > 1$, there exists a bounded linear extension
operator $E : W^{k,p}(\intDD) \to W^{k,p}(\intDD_R)$ with the property
$(Ef)|_{\intDD} = f$ for all $f \in W^{k,p}(\intDD)$; see
Theorem~\ref{thm:SobolevExtension} and Corollary~\ref{cor:compactClosureExtension}.}
\end{exercise}

The above exercise can be used to improve the first part
of Theorem~\ref{thm:regularityLinear} to cover weak solutions of
class $L^1_\loc$.  We start with a classical result about ``weakly
holomorphic'' functions:

\begin{lemma}
\label{lemma:Weyl}
If $u \in L^1(\intDD)$ satisfies $\dbar u = 0$ in the sense of distributions,
then $u$ is smooth and holomorphic.
\end{lemma}
\begin{proof}
Taking real and imaginary parts, it suffices to prove that the same
statement holds for the Laplace equation.  By mollification, any weakly
harmonic function can be approximated in $L^1$ with smooth harmonic
functions.  The latter satisfy the mean value property, which behaves well
under $L^1$-convergence, so the result follows from the mean value
characterization of harmonic functions; see 
\cite{Wendl:lecturesV33}*{Lemma~2.6.26} for more details.
\end{proof}

\begin{lemma}
Suppose $1 < p < \infty$, $k \in \NN$, and $u \in L^1(\intDD)$ is a weak
solution to $\dbar u = f$ for some $f \in W^{k,p}(\intDD)$.  Then
$u$ is of class $W^{k+1,p}$ on every compact subset of~$\intDD$.
\end{lemma}
\begin{proof}
Let $T : W^{k,p}(\intDD) \to W^{k+1,p}(\intDD)$ denote a bounded right inverse
of $\dbar : W^{k+1,p}(\intDD) \to W^{k,p}(\intDD)$ as provided by
Exercise~\ref{EX:BRI}.  Then $u - Tf \in L^1(\intDD)$ is a weak solution
to $\dbar (u - Tf) = 0$ and is thus smooth by Lemma~\ref{lemma:Weyl}.
In particular, $u - Tf$ restricts to $\intDD_r$ for every $r < 1$ as a
function of class~$W^{k+1,p}$, implying that $u$ also has a restriction
in $W^{k+1,p}(\intDD_r)$.
\end{proof}

\begin{cor}[Weak linear regularity]
\label{cor:weakReg}
Suppose $1 < p < \infty$.  Then given $A \in C^\infty(\DD,\End_\RR(\CC^n))$, 
every weak solution
$u \in L^p(\intDD,\CC^n)$ of $(\dbar + A) u = 0$ is smooth on~$\intDD$.
\qed
\end{cor}

\subsection{The nonlinear case}

Locally, every $J$-holomorphic curve can be regarded as a map
$u : \intDD \to \CC^n$ satisfying $u(0) = 0$ and
$$
\dbar_J u := \p_s u + J(u) \p_t u = 0,
$$
where $J$ is a smooth almost complex structure on $\CC^n$ satisfying
$J(0) = i$.  Theorem~\ref{thm:regularityLinear} now has the following
analogue.

\begin{thm}[Nonlinear regularity]
\label{thm:regularity}
Assume $1 < p < \infty$ and $k \in \NN$ satisfy $kp > 2$,
and fix a smooth almost complex structure $J$ on $\CC^n$ with
$J(0) = i$.
\begin{enumerate}
\item \label{item:smooth}
Every map $u \in W^{k,p}(\intDD,\CC^n)$ satisfying $u(0) = 0$ and
$\dbar_J u = 0$ is smooth on~$\intDD$.
\item Suppose $J_\nu$ is a sequence of smooth almost complex structures
on $\CC^n$ converging in $C^\infty_\loc$ to $J$ as $\nu \to \infty$, and
$u_\nu \in W^{k,p}(\intDD,\CC^n)$ is a sequence of smooth maps
satisfying $\dbar_{J_\nu} u_\nu = 0$.
\begin{enumerate}
\item \label{item:unif}
If the maps $u_\nu$ are uniformly $W^{k,p}$-bounded on $\intDD$,
then they are also uniformly $C^m$-bounded on compact
subsets of $\intDD$ for every $m \in \NN$.  
\item \label{item:conv}
If the sequence $u_\nu$ is $W^{k,p}$-convergent on~$\intDD$ to a
smooth map $u : \intDD \to \CC^n$, then it
is also $C^\infty$-convergent on every compact subset of~$\intDD$.
\end{enumerate}
\end{enumerate}
\end{thm}

Our proof of this will follow much the same outline as the proof of
Theorem~\ref{thm:regularityLinear}, and indeed, one could use exactly the
same argument if $J$ were identically equal to $i$ (in which case the theorem
can also be deduced from complex analysis).  The reason it works in the
general case is that if we zoom in on a sufficiently small neighborhood of
the origin in $\CC^n$, then $J$ can be viewed as a $C^\infty$-small
perturbation of~$i$.  To make this precise, we shall use the following
rescaling trick.

Associate to any smooth almost complex structure $J$ on $\CC^n$ the function
$$
Q := i - J \in C^\infty(\CC^n,\End_\RR(\CC^n)).
$$
In terms of $Q$, the equation $\p_s u + J(u) \p_t u = 0$ then becomes
\begin{equation}
\label{eqn:rewrittenCR}
\dbar u - (Q \circ u) \p_t u = 0,
\end{equation}
where we are regarding $Q \circ u$ as a function $\intDD \to \End_\RR(\CC^n)$.
Given constants $R \ge 1$ and $\epsilon \in (0,1]$, associate to
$J$ and $u$ the functions
\begin{equation}
\label{eqn:rescaling}
\begin{split}
\widehat{J} : \CC^n \to \End_\RR(\CC^n), &\qquad \widehat{J}(p) := J(p/R), \\
\widehat{Q} : \CC^n \to \End_\RR(\CC^n), &\qquad \widehat{Q}(p) := Q(p/R) = i - \widehat{J}(p), \\
\hat{u} : \intDD \to \CC^n, &\qquad u(z) := R u(\epsilon z).
\end{split}
\end{equation}
Now $u$ satisfies \eqref{eqn:rewrittenCR} if and only if $\hat{u}$ satisfies
\begin{equation}
\label{eqn:rescaledCR}
\dbar \hat{u} - (\widehat{Q} \circ \hat{u}) \p_t \hat{u} = 0.
\end{equation}
The rescaled almost complex structure has the convenient feature that if
$J(0) = i$, then $\widehat{J}$ can be made arbitrarily $C^\infty$-close to~$i$
on the unit disk 
$$
\DD^{2n} \subset \CC^n
$$
by choosing $R$ sufficiently large, which means 
$\| \widehat{Q} \|_{C^m(\DD^{2n})}$ can be made arbitrarily small for
every $m \in \NN$. 
If $u$ is also continuous and satisfies $u(0) = 0$, then after fixing
some large value for~$R$, we can also choose $\epsilon \in (0,1]$ sufficiently
small to ensure $\overline{u(\intDD)} \subset \intDD^{2n}$ and make
$\| \widehat{Q} \circ \hat{u} \|_{C^0(\DD)}$ arbitrarily small.  By 
Propositions~\ref{prop:CkContinuity} and~\ref{prop:SobolevRescaling}, we can
similarly arrange for $\| \widehat{Q} \circ \hat{u} \|_{W^{k,p}}$
to be arbitrarily small if $u$ is of class $W^{k,p}$ with $kp > 2$, and
the same will hold for $\| \widehat{Q}_\nu \circ \hat{u}_\nu \|_{W^{k,p}}$ when
$\nu$ is large if $\| u_\nu \|_{W^{k,p}}$ is uniformly bounded and
$u_\nu(0) \to 0$.  Here of course we abbreviate $Q_\nu := i - J_\nu$ and
$\widehat{Q}_\nu(p) := Q_\nu(p/R)$.
The effect is to make equations such as \eqref{eqn:rescaledCR} 
$W^{k,p}$-close to the linear equation $\dbar\hat{u} = 0$ if $\epsilon > 0$
and $R > 0$ are sufficiently small and large respectively.

The price we pay for this rescaling is that if we are able to prove 
e.g.~a uniform bound on the norms $\| \hat{u}_\nu \|_{W^{k,p}(\intDD)}$ for 
some sequence $u_\nu$, then the resulting $W^{k+1,p}$-bound for $u_\nu$ will 
be valid only on~$\intDD_{\epsilon}$, a very small ball about the origin.  
But this is good enough for obtaining estimates over all compact subsets 
of~$\intDD$: indeed, we can always reparametrize $u : \intDD \to \CC^n$ to
put the origin at some other point and prove suitable estimates near that
point, appealing in the end to the fact that any compact subset of
$\intDD$ is covered by a finite union of small disks about points.

The need to use this rescaling trick is one of a few reasons why the
condition $kp > 2$ is needed in Theorem~\ref{thm:regularity}, while it
was irrelevant in the linear case.

\begin{proof}[Proof of Theorem~\ref{thm:regularity}]
We will prove statement~\eqref{item:unif} and leave the rest as exercises.

By the remarks above, it suffices to prove that if $u_\nu : \intDD \to \CC^n$
are smooth $J_\nu$-holomorphic curves satisfying a uniform bound in
$W^{k,p}(\intDD)$, then for some $r < 1$, the rescaled 
$\widehat{J}_\nu$-holomorphic curves
$\hat{u}_\nu : \intDD \to \CC^n$ defined as in \eqref{eqn:rescaling}
satisfy a uniform $W^{k+1,p}$-bound on~$\intDD_r$.
In fact, it suffices to prove that every subsequence of $u_\nu$ has
a further subsequence for which this is true. Indeed, if the bound for the 
whole sequence did not exist, then we would be able to find a subsequence with 
norms blowing up to infinity, and no further subsequence of this subsequence 
could satisfy a uniform bound. With this understood, we can appeal to the fact 
that $W^{k,p}$-bounded sequences are also $C^0$-bounded for $kp > 2$ and 
thus replace $u_\nu$ with a subsequence (still denoted by $u_\nu$) such that, 
after a suitable change of coordinates on $\CC^n$,
$$
u_\nu(0) \to 0.
$$
Our goal is then to show that for a suitable choice of the rescaling parameters
$\epsilon$ and $R$, this subsequence admits a uniform bound on
$\| \p_j \hat{u}_\nu \|_{W^{k,p}(\intDD_r)}$ for $j = 1,2$	.

The argument begins exactly the same as in the linear case:
choose a smooth bump function
$$
\beta \in C_0^\infty(\intDD,[0,1])
$$
that satisfies $\beta|_{\DD_r} \equiv 1$.  We then have
$\beta \, \p_j \hat{u}_\nu \in C_0^\infty(\intDD)$, so by 
Theorem~\ref{thm:estimate},
\begin{equation}
\label{eqn:boundThis2}
\| \p_j \hat{u}_\nu \|_{W^{k,p}(\intDD_r)} \le
\| \beta \, \p_j \hat{u}_\nu \|_{W^{k,p}(\intDD)} \le
c \left\| \dbar \left( \beta \, \p_j \hat{u}_\nu \right) \right\|_{W^{k-1,p}(\intDD)}.
\end{equation}
Instead of rewriting $\dbar( \beta \, \p_j \hat{u}_\nu)$ as a sum of two
terms, let us derive a PDE satisfied by $\beta \, \p_j \hat{u}_\nu$.
Differentiating the equation
$\dbar \hat{u}_\nu - (\widehat{Q}_\nu \circ \hat{u}_\nu) \p_t \hat{u}_\nu = 0$
gives
$$
\dbar(\p_j \hat{u}_\nu) = \p_j (\dbar \hat{u}_\nu) =
( d \widehat{Q}_\nu \circ \hat{u}_\nu ) \left(\p_j \hat{u}_\nu,
\p_t \hat{u}_\nu \right) + ( \widehat{Q}_\nu \circ \hat{u}_\nu )
\p_j \p_t \hat{u}_\nu,
$$
thus $\beta \, \p_j \hat{u}_\nu$ satisfies
\begin{equation}
\label{eqn:rescaledCRdbeta}
\begin{split}
\dbar (\beta \, \p_j \hat{u}_\nu) - (\widehat{Q}_\nu \circ & \hat{u}_\nu)
\p_t (\beta \, \p_j \hat{u}_\nu) \\
&= \beta (d\widehat{Q}_\nu \circ \hat{u}_\nu)(\p_j\hat{u}_\nu,
\p_t \hat{u}_\nu) + \left( \dbar\beta - (\widehat{Q}_\nu \circ \hat{u}_\nu) 
\p_t \beta \right) \p_j\hat{u}_\nu 
\\ &= (d\widehat{Q}_\nu \circ \hat{u}_\nu)(\beta \, \p_j\hat{u}_\nu,
\p_t \hat{u}_\nu) + \left( \dbar\beta - (\widehat{Q}_\nu \circ \hat{u}_\nu) 
\p_t \beta \right) \p_j\hat{u}_\nu ,
\end{split}
\end{equation}
and combining this with \eqref{eqn:boundThis2} gives
\begin{multline}
\label{eqn:regc1'}
\left\| \beta \, \p_j \hat{u}_\nu \right\|_{W^{k,p}} \le
c \big\| (\widehat{Q}_\nu \circ \hat{u}_\nu) \p_t (\beta \, \p_j \hat{u}_\nu)
\big\|_{W^{k-1,p}} 
+ c \big\| (d\widehat{Q}_\nu \circ \hat{u}_\nu)(\beta \, \p_j\hat{u}_\nu,
\p_t \hat{u}_\nu) \big\|_{W^{k-1,p}} \\
+ c \left\| \left( \dbar\beta - (\widehat{Q}_\nu \circ \hat{u}_\nu) 
\p_t \beta \right) \p_j\hat{u}_\nu \right\|_{W^{k-1,p}}.
\end{multline}
In order to find bounds for the three terms on the right, recall that using
Propositions~\ref{prop:CkContinuity} and~\ref{prop:SobolevRescaling}
and the assumption $u_\nu(0) \to 0$, we can suppose
$$
\left\| \widehat{Q}_\nu \circ \hat{u}_\nu \right\|_{W^{k,p}} \le \delta
$$
for sufficiently large~$\nu$, where $\delta > 0$ is a constant that may be assumed
arbitrarily small via suitable choices of the rescaling parameters
$\epsilon$ and~$R$.
This provides a uniform bound on the third term in \eqref{eqn:regc1'},
as there is also a continuous product pairing $W^{k,p} \times W^{k-1,p} \to
W^{k-1,p}$ by Prop.~\ref{prop:BanachAlg}, giving an estimate
of the form
\begin{equation*}
\begin{split}
\left\| \left( \dbar\beta - (\widehat{Q}_\nu \circ \hat{u}_\nu) 
\p_t \beta \right) \p_j\hat{u}_\nu \right\|_{W^{k-1,p}} &\le
c \left\| \left( \dbar\beta - (\widehat{Q}_\nu \circ \hat{u}_\nu) 
\p_t \beta \right) \right\|_{W^{k,p}} \cdot
\| \p_j\hat{u}_\nu \|_{W^{k-1,p}} \\
&\le c' \| \hat{u}_\nu \|_{W^{k,p}} \le c''.
\end{split}
\end{equation*}
For the first term on the right side of \eqref{eqn:regc1'}, the product
pairing similarly gives
\begin{equation*}
\begin{split}
\big\| (\widehat{Q}_\nu \circ \hat{u}_\nu) \p_t (\beta \, \p_j \hat{u}_\nu)
\big\|_{W^{k-1,p}} &\le c 
\big\| \widehat{Q}_\nu \circ \hat{u}_\nu \big\|_{W^{k,p}}
\cdot \| \p_t (\beta \, \p_j \hat{u}_\nu) \|_{W^{k-1,p}} \\
&\le c \delta \| \beta \, \p_j \hat{u}_\nu \|_{W^{k,p}}.
\end{split}
\end{equation*}
Finally, since $J_\nu \to J$ in $C^{k+1}$ on compact subsets, we are also
free to assume after adjusting the rescaling parameters that
$$
\| d\widehat{Q}_\nu \circ \hat{u}_\nu \|_{W^{k,p}} \le \delta,
$$
so we can apply the product pairing $W^{k,p} \times W^{k-1,p} \to W^{k-1,p}$
twice to estimate
\begin{equation*}
\begin{split}
\big\| (d\widehat{Q}_\nu \circ \hat{u}_\nu)(\beta \, \p_j\hat{u}_\nu,
\p_t \hat{u}_\nu) \big\|_{W^{k-1,p}} &\le 
c \big\| d\widehat{Q}_\nu \circ \hat{u}_\nu \big\|_{W^{k,p}}
\cdot \| \beta \, \p_j \hat{u}_\nu \|_{W^{k,p}} \cdot
\| \p_t \hat{u}_\nu \|_{W^{k-1,p}} \\
&\le c \delta \| \beta \, \p_j \hat{u}_\nu \|_{W^{k,p}} \cdot
\| \hat{u}_\nu \|_{W^{k,p}} \\
&\le c c' \delta 
\| \beta \, \p_j \hat{u}_\nu \|_{W^{k,p}}
=: c'' \delta \| \beta \, \p_j \hat{u}_\nu \|_{W^{k,p}}.
\end{split}
\end{equation*}
Combining the three estimates for the right hand side of
\eqref{eqn:regc1'} now gives
$$
\| \beta \, \p_j \hat{u}_\nu \|_{W^{k,p}} \le c + c \delta 
\| \beta \, \p_j \hat{u}_\nu \|_{W^{k,p}},
$$
so after adjusting the scaling parameters $R$ and $\epsilon$ to ensure
$c \delta < 1$, we obtain the uniform bound
$$
\| \beta \, \p_j \hat{u}_\nu \|_{W^{k,p}} \le \frac{c}{1 - c \delta}.
$$
This provides the desired uniform bound on 
$\| \p_j \hat{u}_\nu \|_{W^{k,p}(\intDD_r)}$.
\end{proof}

\begin{exercise}
Use an analogous argument via difference quotients to prove 
statement~\eqref{item:smooth} in Theorem~\ref{thm:regularity}.
\textsl{Hint: If you're anything like me, you might get stuck trying
to estimate the second term in the difference quotient analogue of
\eqref{eqn:regc1'}.  The difficulty is that this expression was derived
using the chain rule for derivatives, and there is no similarly simple
chain rule for difference quotients.  The trick is to remember that
difference quotients only differ from the corresponding derivatives by a
remainder term.  The remainder will produce an extra term in the difference
quotient version of
\eqref{eqn:regc1'}, but the extra term can be bounded.}
\end{exercise}

\section{Linear local existence and applications}
\label{sec:similarity}

The following lemma can be applied in the case 
$A \in C^\infty(\DD,\End_\CC(\CC^n))$ to prove the aforementioned standard
fact that complex-linear Cauchy-Riemann type operators induce
holomorphic structures on vector bundles.  The version with weakened
regularity will be applied below to prove a useful ``unique continuation''
result about solutions to $(\dbar + A) f = 0$ in the real-linear case.

\begin{lemma}
\label{lemma:linearExistence}
Assume $2 < p < \infty$ and $A \in L^p(\intDD,\End_\RR(\CC^n))$.  
Then for sufficiently small $\epsilon > 0$, the problem
\begin{equation*}
\begin{split}
\dbar u + A u &= 0 \\
u(0) &= u_0
\end{split}
\end{equation*}
has a solution $u \in W^{1,p}(\intDD_\epsilon,\CC^n)$.
\end{lemma}
\begin{remark}
Note that $u : \intDD_\epsilon \to \CC^n$ in the above statement is only a
\emph{weak} solution to $\dbar u + Au = 0$, as it is not necessarily
differentiable, but by the Sobolev embedding theorem, it is at least
continuous.
\end{remark}
\begin{proof}[Proof of Lemma~\ref{lemma:linearExistence}]
The main idea is that if we take $\epsilon > 0$ sufficiently small,
then the restriction of $\dbar + A$ to $\intDD_\epsilon$ can be regarded as a
small perturbation of $\dbar$ in the space of bounded linear operators
$W^{1,p} \to L^p$.  Since the latter has a bounded right inverse
by Theorem~\ref{thm:rightInverse}, the same will be true for the perturbation.

Since $p > 2$, the Sobolev embedding theorem implies that functions
$u \in W^{1,p}$ are also continuous and bounded by $\| u \|_{W^{1,p}}$, 
thus we can define a bounded linear operator
$$
\Phi : W^{1,p}(\intDD) \to L^p(\intDD) \times \CC^n : u \mapsto (\dbar u, u(0)).
$$
Theorem~\ref{thm:rightInverse}
implies that this operator is also surjective and has a bounded right inverse,
namely
$$
L^p(\intDD) \times \CC^n \to W^{1,p}(\intDD) : (f,u_0) \mapsto
Tf - Tf(0) + u_0,
$$
where $T : L^p(\intDD) \to W^{1,p}(\intDD)$ is a right inverse of~$\dbar$.
Thus any operator sufficiently close to $\Phi$ in the norm topology
also has a right inverse.
Now define $\chi_\epsilon : \DD \to \RR$ to be the function that equals $1$ on 
$\DD_\epsilon$ and $0$ outside of it, and let
$$
\Phi_\epsilon : W^{1,p}(\intDD) \to L^p(\intDD) \times \CC^n : u \mapsto
((\dbar + \chi_\epsilon A) u , u(0)).
$$
To see that this is a bounded operator, it suffices to check that
$W^{1,p} \to L^p : u \mapsto A u$ is bounded if $A \in L^p$; indeed,
$$
\| A u \|_{L^p} \le \| A \|_{L^p} \| u \|_{C^0} \le
c \| A \|_{L^p} \| u \|_{W^{1,p}},
$$
again using the Sobolev embedding theorem.
Now by this same trick, we find
$$
\|\Phi_\epsilon u - \Phi u\| = \| \chi_\epsilon A u \|_{L^p(\intDD)} 
\le c \| A \|_{L^p(\intDD_\epsilon)} \| u \|_{W^{1,p}(\intDD)},
$$
thus $\|\Phi_\epsilon - \Phi\|$ is small if $\epsilon$ is small, and it
follows that in this case $\Phi_\epsilon$ is surjective.  Our desired solution
is therefore the restriction of any 
$u \in \Phi_\epsilon^{-1}(0,u_0)$ to~$\intDD_\epsilon$.
\end{proof}

Here is a corollary, which says that every solution to a real-linear
Cauchy-Riemann type equation looks locally like a holomorphic function in
some \emph{continuous} local trivialization.

\begin{thm}[Similarity principle]
\label{thm:similarity}
Suppose $A : \DD \to \End_\RR(\CC^n)$ is smooth and
$u : \intDD \to \CC^n$ satisfies the
equation $\dbar u + Au = 0$ with $u(0) = 0$.  Then for sufficiently 
small $\epsilon > 0$, there exist maps $\Phi \in C^0(\DD_\epsilon,
\End_\CC(\CC^n))$ and $f \in C^\infty(\intDD_\epsilon,\CC^n)$ such that
$$
u(z) = \Phi(z) f(z),
\qquad
\dbar f = 0,
\qquad\text{ and }\qquad
\Phi(0) = \1.
$$
\end{thm}
\begin{proof}
After shrinking the domain if necessary, we may assume without loss of
generality that the smooth solution $u : \intDD \to \CC^n$ is bounded.
Choose a map
$C : \DD \to \End_\CC(\CC^n)$ satisfying $C(z) u(z) = A(z) u(z)$ and
$|C(z)| \le |A(z)|$ for almost every $z \in \DD$.  Then
$C \in L^\infty(\intDD,\End_\CC(\CC^n))$ and $u$ is a weak solution to
$(\dbar + C)u = 0$.  Note that since we do not know anything about the
zero set of $u$, we cannot assume $C$ is continuous, but we have no trouble
assuming $C \in L^p(\intDD)$ for every $p > 2$.

Since $\dbar + C$ is now complex linear, we can use 
Lemma~\ref{lemma:linearExistence} to find a complex
basis of $W^{1,p}$-smooth weak solutions 
to $(\dbar + C) v = 0$ on $\intDD_\epsilon$ that define
the standard basis of $\CC^n$ at~$0$, and these solutions are continuous
by the Sobolev embedding theorem.  This gives rise to a map
$\Phi \in C^0(\intDD_\epsilon,\End_\CC(\CC^n))$ that satisfies
$(\dbar + C)\Phi = 0$ in the sense of distributions 
and $\Phi(0) = \1$.  Since $\Phi$ is continuous,
we can assume without loss of generality that $\Phi(z)$ is
invertible everywhere on~$\intDD_\epsilon$.  Setting $f := \Phi^{-1} u : 
\intDD_\epsilon \to \CC^n$, the Leibniz rule then implies
$$
0 = (\dbar + C) u = (\dbar + C) (\Phi f) = \left[ (\dbar + C) \Phi \right] f
+ \Phi (\dbar f) = \Phi (\dbar f),
$$
thus $\dbar f = 0$, and $f$ is smooth by Lemma~\ref{lemma:Weyl}.
\end{proof}

\begin{cor}[Unique continuation]
Suppose $\mathbf{D}$ is a linear Cauchy-Riemann type operator on a
vector bundle $E$ over a connected Riemann surface, 
and $\eta \in \Gamma(E)$ satisfies 
$\mathbf{D} \eta = 0$.  Then either $\eta$ is identically zero or its
zeroes are isolated.
\end{cor}

The similarity principle also has many nice applications for the nonlinear
Cauchy-Riemann equation.  Here is another ``unique continuation'' type
result for the nonlinear case.

\begin{prop}
\label{prop:uniqueContin}
Suppose $J$ is a smooth almost complex structure on $\CC^n$ and
$u,v : \intDD \to \CC^n$ are smooth $J$-holomorphic curves such that
$u(0) = v(0) = 0$ and $u$ and $v$ have matching partial derivatives 
of all orders at~$0$.  Then $u \equiv v$ on a neighborhood of~$0$.
\end{prop}
\begin{proof}
Let $h = v - u : \intDD \to \CC^n$.  We have
\begin{equation}
\label{eqn:u}
\p_s u + J(u(z)) \p_t u = 0
\end{equation}
and
\begin{equation}
\label{eqn:v}
\begin{split}
\p_s v + J(u(z)) \p_t v &= \p_s v + J(v(z)) \p_t v + \left[
J(u(z)) - J(v(z)) \right] \p_t v \\
&= - \left[ J(u(z) + h(z)) - J(u(z)) \right] \p_t v \\
&= -\left( \int_0^1 \frac{d}{dt} J(u(z) + t h(z)) \,dt \right) \p_t v \\
&= -\left(\int_0^1 dJ(u(z) + t h(z)) \cdot h(z) \,dt \right) 
\p_t v =: -A(z) h(z),
\end{split}
\end{equation}
where the last step defines a smooth family of linear maps $A(z) \in
\End_\RR(\CC^n)$.  Subtracting \eqref{eqn:u} from \eqref{eqn:v} gives
the linear equation
$$
\p_s h(z) + \bar{J}(z) \p_t h(z) + A(z) h(z) = 0,
$$
where $\bar{J}(z) := J(u(z))$.  This is a linear Cauchy-Riemann type
equation on a trivial complex vector bundle over $\intDD$ with complex
structure $\bar{J}(z)$ on the fiber at~$z$.  
The similarity principle thus implies
$h(z) = \Phi(z) f(z)$ near~$0$ for some holomorphic function
$f(z) \in \CC^n$ and some continuous map $\Phi(z) \in \GL(2n,\RR)$ 
representing a change of trivialization.  Now if
$h$ has vanishing derivatives of all orders at~$0$, Taylor's formula
implies 
$$
\lim_{z \to 0} \frac{| \Phi(z) f(z) |}{|z|^k} = 0
$$
for all $k \in \NN$, so $f$ must also have a zero of infinite order
and thus $f \equiv 0$.
\end{proof}

\section{Simple curves and multiple covers}

We now prove a global result about the structure of closed $J$-holomorphic 
curves.  In Lecture~\ref{lec:cobordisms} we will be able to generalize it 
in a straightforward way for
punctured holomorphic curves with asymptotically cylindrical behavior.

\begin{thm}
\label{thm:simple}
Assume $(\Sigma,j)$ is a closed connected Riemann surface,
$(W,J)$ is a smooth almost complex manifold and
$u : (\Sigma,j) \to (W,J)$ is a nonconstant pseudoholomorphic curve.
Then there exists a factorization $u = v \circ \varphi$, where
\begin{itemize}
\item $\varphi : (\Sigma,j) \to (\Sigma',j')$ is a holomorphic map
of positive degree to another closed and connected Riemann surface
$(\Sigma',j')$;
\item $v : (\Sigma',j') \to (W,J)$ is a pseudoholomorphic curve which is
embedded except at a finite set of critical points and self-intersections.
\end{itemize}
\end{thm}

Note that holomorphic maps $(\Sigma,j) \to (\Sigma',j')$ of degree~$1$ are
always diffeomorphisms, so the factorization $u = v \circ \varphi$ in this
case is just a reparametrization, and $u$ is then called a
\defin{simple} curve.  In all other cases, $k := \deg(\varphi) \ge 2$ and
$\varphi$ is in general a branched cover; we then call $u$ a
\defin{$k$-fold branched cover} of the simple curve~$v$.

The main idea in the proof is to construct $\Sigma'$ (minus some punctures) explicitly as 
the image of $u$ after removing finitely many singular points, so that we 
can take $v$ to be the inclusion $\Sigma' \hookrightarrow W$.  The map 
$\varphi : \Sigma \to \Sigma'$ is then uniquely determined.  In order to 
carry out this program, we need some information on what the image of $u$ 
can look like near each of its singularities.  These come in two types, 
each type corresponding to one of the lemmas below, both of which should
seem immediately plausible if your intuition comes from complex analysis.

\begin{lemma}[Intersections]
Suppose $u : (\Sigma,j) \to (W,J)$ and $v : (\Sigma',j') \to (W,J)$ are 
two nonconstant pseudoholomorphic curves with an intersection $u(z) = v(z')$. 
Then there exist neighborhoods $z \in \uU \subset \Sigma$ and 
$z' \in \uU' \subset \Sigma'$ such that
$$
\text{either } \quad
u(\uU) = v(\uU') \qquad \text{ or } \qquad
u(\uU \setminus \{z\}) \cap v(\uU') = u(\uU) \cap v(\uU' \setminus \{z'\})
= \emptyset.
$$
\qed
\end{lemma}

\begin{lemma}[Branching]
Suppose $u : (\Sigma,j) \to (W,J)$ is a nonconstant pseudoholomorphic curve 
and $z_0 \in \Sigma$ is a critical point of~$u$. Then a neighborhood 
$\uU \subset \Sigma$ of $z_0$ can be biholomorphically identified with the 
unit disk $\DD \subset \CC$ such that
$$
u(z) = v(z^k) \quad \text{ for } \quad
z \in \DD = \uU,
$$
where $k \in \NN$, and $v : \DD \to W$ is an injective $J$-holomorphic 
map with no critical points except possibly at the origin.
\qed
\end{lemma}

These two local results follow from a well-known formula of Micallef and White 
\cite{MicallefWhite} describing the local behavior of $J$-holomorphic curves 
near critical points and their intersections.  The proof of that theorem is 
analytically quite involved, but one can also use an easier ``approximate'' 
version, which is proved in \cite{Wendl:lecturesV33}*{\S 2.14}.  Since both
are closely related to the phenomenon of unique continuation, you will not
be surprised to learn that the similarity principle plays a role in the proof:
the main idea is again to exploit the fact that locally $J$ is always a
small perturbation of~$i$, hence the local behavior of $J$-holomorphic curves 
is also similar to the integrable case.

\begin{proof}[Proof of Theorem~\ref{thm:simple}]
Let $\text{Crit}(u) = \{ z \in \Sigma\ |\ du(z) = 0 \}$ denote the set of critical points, and define $\Delta \subset \Sigma$ to be the set of all points $z \in \Sigma$ such that there exists $z' \in \Sigma$ and neighborhoods $z \in {\mathcal U} \subset \Sigma$ and $z' \in {\mathcal U}' \subset \Sigma$ with $u(z) = u(z')$ but $u({\mathcal U} \setminus \{z\}) \cap u({\mathcal U}' \setminus \{z'\}) = \emptyset$.

The lemmas quoted above imply that both of these sets are discrete. Both are therefore finite, and the set $\dot{\Sigma}' = u(\Sigma \setminus (\text{Crit}(u) \cup \Delta)) \subset W$ is then a smooth submanifold of $W$ with $J$-invariant tangent spaces, so it inherits a natural complex structure $j'$ for which the inclusion $(\dot{\Sigma}',j') \hookrightarrow (W,J)$ is pseudoholomorphic. We shall now construct a new Riemann surface $(\Sigma',j')$ from which $(\dot{\Sigma}',j')$ is obtained by removing a finite set of points. Let $\widehat{\Delta} = (\text{Crit}(u) \cup \Delta) / \sim$, where two points in $\text{Crit}(u) \cup \Delta$ are defined to be equivalent whenever they have neighborhoods in $\Sigma$ with identical images under $u$. Then for each $[z] \in \widehat{\Delta}$, the branching lemma provides an injective $J$-holomorphic map $u_{[z]}$ from the unit disk $\mathbb{D}$ onto the image of a neighborhood of $z$ under $u$. We define $(\Sigma',j')$ by $$\Sigma' = \dot{\Sigma}' \cup_\Phi \left( \bigsqcup_{[z] \in \widehat{\Delta}} {\mathbb D} \right),$$ where the gluing map $\Phi$ is the disjoint union of the maps $u_{[z]} : {\mathbb D} \setminus \{0\} \to \dot{\Sigma}'$ for each $[z] \in \widehat{\Delta}$; since this map is holomorphic, the complex structure $j'$ extends from $\dot{\Sigma}'$ to $\Sigma'$. Combining the maps $u_{[z]} : {\mathbb D} \to W$ with the inclusion $\dot{\Sigma}' \hookrightarrow W$ now defines a pseudoholomorphic map $v : (\Sigma',j') \to (W,J)$ which restricts to $\dot{\Sigma}'$ as an embedding and otherwise has at most finitely many critical points and double points. Moreover, the restriction of $u$ to $\Sigma \setminus (\text{Crit}(u) \cup \Delta)$ defines a holomorphic map to $(\dot{\Sigma}',j')$ which extends by removal of singularities to a proper holomorphic map $\varphi : (\Sigma,j) \to (\Sigma',j')$ such that $u = v \circ \varphi$. Its holomorphicity implies that it has positive degree.  
\end{proof}

\chapter{Asymptotic operators}
\label{lec:asymptotic}

\minitoc

\vspace{12pt}

We now begin with the analysis of the particular class of $J$-holomorphic 
curves that are important in SFT.  The next three lectures will focus on the
linearized problem, the goal being to prove that this linearization
is Fredholm and to compute its index.  Using this along with the
implicit function theorem and the Sard-Smale theorem (on genericity of
smooth nonlinear Fredholm maps), we will later be able to show that moduli
spaces of asymptotically cylindrical $J$-holomorphic curves are smooth
finite-dimensional manifolds under suitable genericity assumptions.

\section{The linearization in Morse homology}

Since Morse homology is the prototype for all Floer-type theories, we can
gain useful intuition by recalling how the analysis works for the 
linearization of the gradient flow problem in Morse theory.  The basic
features of the problem were discussed already in \S\ref{sec:Floer}.

Assume $(M,g)$ is a closed $n$-dimensional Riemannian manifold,
$f : M \to \RR$ is a smooth function, and for two critical points
$x_+,x_- \in \Crit(f)$, consider the moduli space of parametrized
gradient flow lines
$$
\mM(x_-,x_+) := \left\{ u \in C^\infty(\RR,M)\ \big|\ 
\dot{u} + \nabla f(u) = 0, \ \lim_{s \to \pm\infty} u(s) = x_\pm \right\}.
$$
The map $\mM(x_-,x_+) \to M : u \mapsto u(0)$ gives a natural
identification of $\mM(x_-,x_+)$ with the intersection between the unstable
manifold of $x_-$ and the stable manifold of $x_+$ for the negative gradient
flow, and we say the pair $(g,f)$ is \defin{Morse-Smale} if $f$ is Morse and
this intersection
is transverse, in which case $\mM(x_-,x_+)$ is a smooth manifold with
$$
\dim \mM(x_-,x_+) = \ind(x_-) - \ind(x_+).
$$
This can all be proved using finite-dimensional differential topology, but
since that approach does not work in the study of Floer trajectories or
holomorphic curves in symplectizations, let us instead see how one proves it
using nonlinear functional analysis.  For more details on the following
discussion, see \cite{Schwarz:Morse}.

Following the strategy laid out in \S\ref{sec:linearCR}, $\mM(x_-,x_+)$ can be
identified with the zero set of a smooth section
$$
\boldsymbol{\sigma} : \bB \to \eE : u \mapsto \dot{u} + 
\nabla f(u),
$$
where $\bB$ is a Banach manifold of maps $u : \RR \to M$ satisfying
$\lim_{s \to \pm\infty} u(s) = x_\pm$, and $\eE \to \bB$ is a smooth
Banach space bundle whose fibers $\eE_u$ contain $\Gamma(u^*TM)$.
The linearization $D\boldsymbol{\sigma}(u) : T_u \bB \to \eE_u$
of this section at a zero $u \in \boldsymbol{\sigma}^{-1}(0)$ defines a
first-order linear differential operator
$$
\mathbf{D}_u : \Gamma(u^*TM) \to \Gamma(u^*TM)
$$
which takes the form
$$
\mathbf{D}_u \eta = \nabla_s \eta + \nabla_\eta \nabla f
$$
for any choice of symmetric connection $\nabla$ on~$M$.  Taking suitable
Sobolev completions of $\Gamma(u^*TM)$, we are therefore led to consider
bounded linear operators\footnote{We are ignoring an analytical subtlety:
since $u^*TM \to \RR$ has no canonical trivialization and $\RR$ is
noncompact, it is not completely obvious what the definition of the
Sobolev space $W^{k,p}(u^*TM)$ should be.  We will return to this issue
in a more general context in the next lecture.}
of the form
\begin{equation}
\label{eqn:MorseLin}
\mathbf{D}_u = \nabla_s + \nabla \nabla f : W^{k,p}(u^*TM) \to 
W^{k-1,p}(u^*TM)
\end{equation}
for $k \in \NN$ and $1 < p < \infty$,
and the first task is to prove that whenever $x_+$ and $x_-$ satisfy the
Morse condition, this is a Fredholm
operator of index $\ind \mathbf{D}_u = \ind(x_-) - \ind(x_+)$.

Choose coordinates near $x_+$ in which $g$ looks like the standard 
Euclidean inner product at~$x_+$.  This induces a trivialization of
$u^*TM$ over $[T,\infty)$ for $T > 0$ sufficiently large, and we are
free to assume that the connection $\nabla$ is the standard one determined
by these coordinates on $[T,\infty)$.
Using the trivialization to identify sections $\beta \in \Gamma(u^*TM)$ 
over $[T,\infty)$ with maps $f : [T,\infty) \to \RR^n$, $\mathbf{D}_u$
now acts on $f$ as
\begin{equation}
\label{eqn:MorseTriv}
(\mathbf{D}_u f)(s) = \p_s f(s) + A(s) f(s),
\end{equation}
where $A(s) \in \RR^{n \times n}$ is the matrix of the linear transformation
$d X(s) : \RR^n \to \RR^n$, with $X(s) \in \RR^n$ being the coordinate
representation of $\nabla f(u(s)) \in T_{u(s)} M$.
As $s \to \infty$, the zeroth-order
term in this expression converges to a symmetric matrix
$$
A_+ := \lim_{s \to \infty} A(s),
$$
which is the coordinate representation of the Hessian $\nabla^2 f(x_+)$.
Any choice of coordinates near $x_-$ produces a similar formula for
$\mathbf{D}_u$ over $(-\infty,-T]$, $A(s)$ converging as $s \to -\infty$ to
another symmetric matrix $A_-$ representing $\nabla^2 f(x_-)$.  Both the
Morse condition and the
dimension $\ind(x_-) - \ind(x_+)$ can now be expressed entirely in terms of
these two matrices: $x_\pm$ is Morse if and only if $A_\pm$ is invertible,
and the Fredholm index of $\mathbf{D}_u$ will then be
$$
\ind(x_-) - \ind(x_+) = \dim E^-(A_-) - \dim E^-(A_+),
$$
where for any symmetric matrix $A$ we denote by $E^-(A)$ the direct sum of
all its eigenspaces with negative eigenvalue.  The main linear functional
analytic result underlying Morse homology can now be stated as
follows (cf.~\cite{Schwarz:Morse}):

\begin{prop}
\label{prop:MorseLinear}
Assume $k \in \NN$ and $1 < p < \infty$.
Suppose $E \to \RR$ is a smooth vector bundle with trivializations fixed
in neighborhoods of $-\infty$ and $+\infty$, and
$\mathbf{D} : W^{k,p}(E) \to W^{k-1,p}(E)$ is a
first-order differential operator which asymptotically takes the form
\eqref{eqn:MorseTriv} near $\pm\infty$ with respect to the chosen
trivializations, where $A(s)$ is a smooth family of $n$-by-$n$ matrices 
with well-defined asymptotic limits $A_\pm := \lim_{s \to \pm\infty} A(s)$
which are symmetric.  If $A_+$ and $A_-$ are also invertible, then
$\mathbf{D}$ is Fredholm and
\begin{equation}
\label{eqn:spectralFlowHint}
\ind(\mathbf{D}) = \dim E^-(A_-) - \dim E^-(A_+).
\end{equation}
\qed
\end{prop}

\begin{remark}
\label{remark:dependsOnApm}
The hypothesis that $A_\pm$ is invertible in Prop.~\ref{prop:MorseLinear}
cannot be lifted: indeed, suppose $\mathbf{D}$ is Fredholm but
e.g.~$A_+$ has $0$ in its spectrum.  Then one can easily perturb $A(s)$
and hence $A_+$ in two distinct ways producing two distinct values of
$\dim E^-(A_+)$, pushing the zero eigenvalue either up or down.
This produces two perturbed Fredholm operators that have different indices
according to \eqref{eqn:spectralFlowHint}, but they also belong to a continuous
family of Fredholm operators, and must therefore have the same index,
giving a contradiction.
\end{remark}

The formula \eqref{eqn:spectralFlowHint} makes sense of course because 
$E^-(A_\pm)$ are both 
finite-dimen\-sion\-al vector spaces, but in Floer-type theories we typically
encounter critical points with infinite Morse index.  With this in mind,
it is useful to note that \eqref{eqn:spectralFlowHint} can be rewritten
without explicitly referencing $E^-(A_+)$ or $E^-(A_-)$.  Indeed,
choose a continuous path of symmetric matrices $\{B_t\}_{t \in [-1,1]}$ 
connecting $B(-1) := A_-$ to $B(1) := A_+$.  The spectrum of $B_t$ varies
continuously with $t$ in the following sense: one can choose a family of
continuous functions 
$$
\{ \lambda_j : [-1,1] \to \RR \}_{j \in I}
$$ 
for the index set ~$I=\{1,\dots,n\}$ such that for every $t \in [-1,1]$, the set of
eigenvalues of $B_t$ counted with multiplicity is
$\{ \lambda_j(t) \}_{j \in I}$.  The \defin{spectral flow} from
$A_-$ to $A_+$ is then defined as a signed count of the number of paths of 
eigenvalues that cross from one side of zero to the other, namely
(cf.~Theorem~\ref{thm:spectralFlow})
$$
\muspec(A_-,A_+) := \# \left\{ j \in I\ \big|\ 
\lambda_j(-1) < 0 < \lambda_j(1) \right\} -
\# \left\{ j \in I\ \big|\ 
\lambda_j(-1) > 0 > \lambda_j(1) \right\}.
$$
The index formula \eqref{eqn:spectralFlowHint} now becomes
$$
\ind(\mathbf{D}) = \muspec(A_-,A_+).
$$
This description of the index has the advantage that it could potentially
make sense and give a well-defined integer even if $A_\pm$ were symmetric 
operators on an infinite-dimensional Hilbert space: they might both have
infinitely many positive and negative eigenvalues, but only finitely many
that change sign along a path from $A_-$ to~$A_+$.  We will make this
discussion precise in the next section.

\section{Spectral flow}
\label{sec:spectral}

We will see in \S\ref{sec:contactHessian} that in Floer-type theories, the role of the
symmetric linear transformation $T_x M \to T_x M$ defined by the
Hessian $\nabla^2 f(x)$ of a Morse function $f : M \to \RR$ at a critical point
is played by a certain class of symmetric operators on the space of
loops $\eta : S^1 \to \RR^{2n}$, namely operators of the form
\begin{equation}
\label{eqn:generalAsympOp}
(\mathbf{A} \eta)(t) := -J_0 \,\p_t \eta(t) - S(t) \eta(t),
\end{equation}
where $J_0$ denotes the standard complex structure on $\RR^{2n} = \CC^n$,
and $S : S^1 \to \End(\RR^{2n})$ is a smooth loop of symmetric matrices.
The goal of this section is to define a
notion of spectral flow for operators of this type.  Regarding
$\mathbf{A}$ as an unbounded linear operator on $L^2(S^1,\RR^{2n})$ with
dense domain $H^1(S^1,\RR^{2n})$, we will see that its
spectrum consists of isolated real eigenvalues with finite multiplicity.
We shall prove:

\begin{thm}
\label{thm:spectralFlow}
Assume $\left\{ S_s : S^1 \to \End(\RR^{2n}) \right\}_{s \in [-1,1]}$ is a 
smooth family of loops of symmetric matrices, and consider the corresponding 
$1$-parameter family of unbounded linear operators
$$
\mathbf{A}_s = -J_0 \p_t - S_s(t) : L^2(S^1,\RR^{2n}) \supset H^1(S^1,\RR^{2n})
\to L^2(S^1,\RR^{2n}).
$$
Then there exists a set of continuous functions
$$
\{ \lambda_j : [-1,1] \to \RR \}_{j \in \ZZ}
$$
such that
for every $s \in [-1,1]$, the spectrum of $\mathbf{A}_s$ consists of
the numbers $\{\lambda_j(s)\}_{j \in \ZZ}$, each of which is an eigenvalue
with finite multiplicity equal to the number of times it is repeated 
as $j$ varies in~$\ZZ$.

Moreover, if additionally $\mathbf{A}_- := \mathbf{A}_{-1}$ and
$\mathbf{A}_+ := \mathbf{A}_1$ both have trivial kernel, then the
number $\muspec(\mathbf{A}_-,\mathbf{A}_+) \in \ZZ$ defined by
$$
\# \left\{ j \in \ZZ\ \big|\ 
\lambda_j(-1) < 0 < \lambda_j(1) \right\} -
\# \left\{ j \in \ZZ\ \big|\ 
\lambda_j(-1) > 0 > \lambda_j(1) \right\}
$$
is well defined and depends only on $\mathbf{A}_-$ and~$\mathbf{A}_+$.
\end{thm}

We will start by giving a more abstract definition of spectral flow as an
intersection number between a path of symmetric index~$0$ Fredholm operators 
and the subvariety of noninvertible operators.  This relies
on the general fact that spaces of operators with kernel and cokernel of
fixed finite dimensions form smooth finite-codimensional submanifolds in the
Banach space of all bounded linear operators.  We explain this fact in
\S\ref{sec:submanifolds}, and then specialize to the case of symmetric 
index~$0$ operators to define the abstract version of spectral flow in 
\S\ref{sec:symmetric}.   In \S\ref{sec:eigenvalues}, we show that the spectra
of such operators vary continuously under small perturbations, 
and in \S\ref{sec:specHomotopy} we specialize further to
operators of the form \eqref{eqn:generalAsympOp} and explain how to 
interpret the abstract definition of spectral flow in terms of eigenvalues 
crossing the origin in~$\RR$, leading to a proof of Theorem~\ref{thm:spectralFlow}.

Spectral flow can be defined more generally for certain classes of self-adjoint
elliptic partial differential operators, see e.g.~\cites{AtiyahPatodiSinger:spectral,RobbinSalamon},
and standard proofs of its existence typically rely on perturbation results as in
\cite{Kato} for the spectra of self-adjoint operators.  In the following
presentation, we have chosen to avoid making explicit use of self-adjointness and instead
focus on the Fredholm property; in this way the discussion is mostly 
self-contained and, in particular, does not require any
results from~\cite{Kato}.

\subsection{Geometry in the space of Fredholm operators}
\label{sec:submanifolds}

Fix a field
$$
\FF := \RR \text{ or } \CC.
$$
Given Banach spaces
$X$ and $Y$ over~$\FF$, denote by $\Lin_\FF(X,Y)$ the Banach space of
bounded $\FF$-linear maps from $X$ to $Y$, with $\Lin_\FF(X) := \Lin_\FF(X,X)$,
and let 
$$
\Fred_\FF(X,Y) \subset \Lin_\FF(X,Y)
$$
denote the open subset consisting of Fredholm operators.
Recall that an operator $\mathbf{T} \in \Lin_\FF(X,Y)$ is \defin{Fredholm}
if its image is closed,\footnote{It is not strictly necessary to require that
$\im \mathbf{T} \subset Y$ be closed, as this follows from the finite-dimensionality
of the kernel and cokernel, cf.~\cite{AbramovichAliprantis}*{Cor.~2.17}.} 
and its kernel and cokernel (i.e.~the quotient $\coker \mathbf{T} := 
Y / \im \mathbf{T}$) are
both finite dimensional.  Its \defin{index} is defined as
$$
\ind_\FF(\mathbf{T}) := \dim_\FF \ker \mathbf{T} - \dim_\FF \coker \mathbf{T}
\in \ZZ.
$$
The index defines a continuous and thus locally constant function
$\Fred_\FF(X,Y) \to \ZZ$, and for each $i \in \ZZ$, we shall denote
$$
\Fred_\FF^i(X,Y) := \left\{ \mathbf{T} \in \Fred_\FF(X,Y)\ \big|\ \ind(\mathbf{T}) = i \right\}.
$$

We will often have occasion to use the following general construction.
Given $\mathbf{T}_0 \in \Fred_\FF(X,Y)$, one can choose splittings into 
closed linear subspaces
$$
X = V \oplus K, \qquad  Y = W \oplus C
$$
such that $K = \ker \mathbf{T}_0$, $W = \im \mathbf{T}_0$,
the quotient projection $\pi_C : Y \to \coker \mathbf{T}_0$ restricts to 
$C \subset Y$ as an isomorphism, and
$\mathbf{T}_0|_{V}$ defines an isomorphism from $V$ to~$W$.  Using these
splittings, any other $\mathbf{T} \in \Fred_\FF(X,Y)$ can be written
in block form as
$$
\mathbf{T} = \begin{pmatrix}
\mathbf{A} & \mathbf{B} \\
\mathbf{C} & \mathbf{D}
\end{pmatrix},
$$
with $\mathbf{T}_0$ itself written in this way as
$\begin{pmatrix} \mathbf{A}_0 & 0 \\ 0 & 0 \end{pmatrix}$ for some
Banach space isomorphism $\mathbf{A}_0 : V \to W$.
Let $\oO \subset \Fred_\FF(X,Y)$ denote the open neighborhood
of $\mathbf{T}_0$ for which the block $\mathbf{A}$ is invertible, and
define a map
\begin{equation}
\label{eqn:Phi}
\Phi : \oO \to \Hom_\FF(\ker \mathbf{T}_0, \coker \mathbf{T}_0) :
\mathbf{T} \mapsto \mathbf{D} - \mathbf{C} \mathbf{A}^{-1} \mathbf{B}.
\end{equation}

\begin{lemma}
\label{lemma:FredCoords}
The map $\Phi$ in \eqref{eqn:Phi} is smooth, and holomorphic in
the case $\FF=\CC$, and its derivative at $\mathbf{T}_0$ defines a
surjective bounded linear operator
$\Lin_\FF(X,Y) \to \Hom_\FF(\ker \mathbf{T}_0,\coker \mathbf{T}_0)$ of the form
$$
d\Phi(\mathbf{T}_0) \mathbf{H} = \pi_C \mathbf{H}|_{\ker \mathbf{T}_0},
$$
where $\pi_C$ denotes the natural projection $Y \to \coker \mathbf{T}_0$.
Moreover, there exists a smooth function $\Psi : \oO \to \Lin_\FF(X)$
such that for every $\mathbf{T} \in \oO$, $\Psi(\mathbf{T}) : X \to X$ maps
$\ker \Phi(\mathbf{T}) \subset \ker \mathbf{T}_0$ isomorphically to
$\ker \mathbf{T}$.
\end{lemma}
\begin{proof}
Smoothness, holomorphicity\footnote{Holomorphicity in this infinite-dimensional
setting means the same thing as usual: $\Lin_\CC(X,Y)$ and 
$\Hom_\CC(\ker \mathbf{T}_0,\coker \mathbf{T}_0)$ both have natural complex
structures if $\mathbf{T}_0 \in \Fred_\CC(X,Y)$, and we require
$d\Phi(\mathbf{T})$ to commute with them for all $\mathbf{T} \in \oO$.}
and the formula for the derivative are easily
verified from the given formula for~$\Phi$; in particular, since
the blocks $\mathbf{B}$ and $\mathbf{C}$ both vanish
for $\mathbf{T} = \mathbf{T}_0$, we have
\begin{equation*}
\begin{split}
d\Phi(\mathbf{T}_0) : \Lin_\FF(X,Y) &\to 
\Hom_\FF(K,C) \\
\begin{pmatrix}
\mathbf{A}' & \mathbf{B}' \\
\mathbf{C}' & \mathbf{D}'
\end{pmatrix}
&\mapsto \mathbf{D}'.
\end{split}
\end{equation*}
The map $\Psi : \oO \to \Lin_\FF(X) = \Lin_\FF(V \oplus K)$ is defined by
$$
\Psi(\mathbf{T}) = \begin{pmatrix}
\1 & - \mathbf{A}^{-1} \mathbf{B} \\
0  & \1
\end{pmatrix}.
$$  
For each $\mathbf{T}$, this is an isomorphism; indeed, its inverse is given by 
$$
\Psi(\mathbf{T})^{-1} = \begin{pmatrix}
\1 & \mathbf{A}^{-1} \mathbf{B} \\
0  & \1
\end{pmatrix}.
$$
Then $\mathbf{T} \Psi(\mathbf{T}) = \begin{pmatrix}
\mathbf{A} & 0 \\
\mathbf{C} & \Phi(\mathbf{T})
\end{pmatrix}$,
and since $\mathbf{A}$ is invertible, 
$\ker\mathbf{T}\Psi(\mathbf{T}) = \{0\} \oplus \ker \Phi(\mathbf{T})$.
\end{proof}

\begin{prop}
\label{prop:submanifolds}
For each $i \in \ZZ$ and each nonnegative integer $k \ge i$, the subset
$$
\Fred_\FF^{i,k}(X,Y) := \left\{ \mathbf{T} \in \Fred_\FF^i(X,Y)\ \big|\ 
\dim_\FF \ker \mathbf{T} = k \text{ and }
\dim_\FF \coker \mathbf{T} = k - i \right\}
$$
admits the structure of a smooth (and complex-analytic if $\FF = \CC$)
finite-codimensional Banach submanifold of $\Lin_\FF(X,Y)$, with
$$
\codim_\FF \Fred_\FF^{i,k}(X,Y) = k (k-i).
$$
\end{prop}
\begin{proof}
Applying the implicit function theorem to the map $\Phi$ from
Lemma~\ref{lemma:FredCoords} endows a neighborhood of $\mathbf{T}_0$ in
$\Phi^{-1}(0) \subset \Fred_\FF(X,Y)$ with the structure of a smooth Banach 
submanifold with
$$
\codim_\FF \Phi^{-1}(0) = \dim_\FF \Hom_\FF(\ker\mathbf{T}_0,\coker\mathbf{T}_0) =
k (k-i).
$$
If $\FF = \CC$, then $\Phi$ is also holomorphic and $\Phi^{-1}(0)$ is thus
a complex-analytic submanifold near~$\mathbf{T}_0$.
Now observe that for every $\mathbf{T} \in \oO$, 
$$
\dim_\FF \ker \mathbf{T} = \dim_\FF \ker \Phi(\mathbf{T}) \le 
\dim_\FF \ker \mathbf{T}_0 = k,
$$
with equality if and only if $\Phi(\mathbf{T}) = 0$, hence, since the index is locally constant, we get
$\Phi^{-1}(0) = \Fred_\FF^{i,k}(X,Y)$ in a neighborhood of~$\mathbf{T}_0$.
\end{proof}

For real-linear operators of index~$0$, one can use 
Prop.~\ref{prop:submanifolds} to define the following ``relative'' invariant.
Given two Banach space isomorphisms $\mathbf{T}_\pm : X \to Y$ that lie
in the same connected component of $\Fred_\RR(X,Y)$, define
$$
\muspec_{\ZZ_2}(\mathbf{T}_-,\mathbf{T}_+) \in \ZZ_2
$$
as the parity of the number of times that a generic smooth path
$[-1,1] \to \Fred_\RR^0(X,Y)$ from $\mathbf{T}_-$ to
$\mathbf{T}_+$ passes through operators with nontrivial kernel.
This is well defined
due to the following consequences of standard transversality theory
(see Exercise~\ref{EX:transversality}):
first, generic paths $\{\mathbf{T}(t) \in \Fred_\RR^0(X,Y)\}_{t \in [-1,1]}$ 
are transverse to $\Fred_\RR^{0,k}(X,Y)$ for every $k \in \NN$, which implies
via the codimension formula in Prop.~\ref{prop:submanifolds} that they
never intersect $\Fred_\RR^{0,k}(X,Y)$ for $k \ge 2$, and their intersections
with $\Fred_\RR^{0,1}(X,Y)$ are transverse and thus isolated.  Second,
transversality also holds for generic homotopies
$$
[0,1] \times [-1,1] \to \Fred_\RR^0(X,Y) : (s,t) \mapsto \mathbf{T}_s(t)
$$
with fixed end points between 
any pair of generic
paths $\mathbf{T}_0(t)$ and $\mathbf{T}_1(t)$, so that the set of 
intersections with $\Fred_\RR^{0,k}(X,Y)$ is again empty for $k \ge 2$
and forms a smooth $1$-dimensional submanifold in $[0,1] \times [-1,1]$ 
for $k=1$.  This submanifold, moreover, is disjoint from
$[0,1] \times \{-1,1\}$ since $\mathbf{T}_s(\pm 1) = \mathbf{T}_\pm$, 
and it is also compact
since the set of $\mathbf{T} \in \Fred_\RR^0(X,Y)$ with nontrivial
kernel is a closed subset.  We therefore obtain a compact $1$-dimensional
cobordism between the intersection sets of $\mathbf{T}_0$ and $\mathbf{T}_1$
respectively with $\Fred_\RR^{0,1}(X,Y)$, implying that the count of
intersections modulo~$2$ does not depend on the choice of generic path.

\begin{exercise}
\label{EX:transversality}
Convince yourself that the standard results (as in e.g.~\cite{Hirsch}*{\S 3.2}
about generic transversality of intersections between smooth maps 
$f : M \to N$ and submanifolds $A \subset N$ continue to hold---with minimal
modifications to the proofs---when $N$ is an infinite-dimensional Banach
manifold and $A \subset N$ has finite codimension.
\end{exercise}

\begin{exercise}
\label{EX:detSign}
For matrices $A_\pm \in \GL(n,\RR)$, show that $\muspec_{\ZZ_2}(A_-,A_+) = 0$
if and only if $\det A_+$ and $\det A_-$ have the same sign.
\end{exercise}

\subsection{Symmetric operators of index zero}
\label{sec:symmetric}

We now add the following assumptions to the setup from the previous subsection:
\begin{itemize}
\item $Y$ is a Hilbert space $\hH$ over~$\FF$, with inner product 
denoted by $\langle\ ,\ \rangle_\hH$;
\item $X$ is an $\FF$-linear subspace $\dD \subset \hH$, carrying a Banach space
structure for which the inclusion $\dD \hookrightarrow \hH$ is a compact linear 
operator.
\end{itemize}
The notation $\dD = X$ is motivated by the fact that if $\mathbf{T} \in
\Lin_\FF(\dD,\hH)$, then we can also regard $\mathbf{T}$ as an
\defin{unbounded operator} on $\hH$ with domain~$\dD$ and thus consider the
spectrum of $\mathbf{T}$, see \S\ref{sec:eigenvalues} below.

Since $\hH$ is a Hilbert space, the space $\Lin_\FF(\hH)$ of bounded linear
operators from $\hH$ to itself contains a distinguished closed linear subspace
$$
\Lin_\FF^\sym(\hH) \subset \Lin_\FF(\hH),
$$
consisting of self-adjoint operators.  For operators that are bounded from
$\dD$ to $\hH$ but not necessarily defined or bounded on~$\hH$, there is
also the space of \defin{symmetric} operators
$$
\Lin_\FF^\sym(\dD,\hH) := \left\{ \mathbf{T} \in \Lin_\FF(\dD,\hH)\ \big|\ 
\langle x , \mathbf{T} y \rangle_\hH = \langle \mathbf{T} x , y \rangle_\hH
\text{ for all $x,y \in \dD$} \right\}.
$$
Important examples of symmetric
operators are those which are self-adjoint (see Remark~\ref{remark:selfAdjoint}
below), though for our purposes, it will suffice to restrict attention
to symmetric operators that are also Fredholm with index~$0$.
It turns out that the space of symmetric operators in
$\Fred_\FF^{0,1}(\dD,\hH)$ is a canonically co-oriented
hypersurface in $\Lin_\FF^\sym(\dD,\hH)$, so that the invariant
$\muspec_{\ZZ_2}(\mathbf{T}_-,\mathbf{T}_+)$ defined above has a natural
integer-valued lift when $\mathbf{T}_\pm$ are symmetric.  
We will need a slightly more specialized version of this statement in order
to give a general definition of spectral flow.

In the following, we let
$$
\Fred_\FF^\sym(\dD,\hH) := \Fred_\FF^0(\dD,\hH) \cap \Lin_\FF^\sym(\dD,\hH)
$$
denote the space of symmetric Fredholm operators with index~$0$, and 
for $k \in \NN$,
$$
\Fred_\FF^{\sym,k}(\dD,\hH) := \Fred_\FF^\sym(\dD,\hH) \cap \Fred_\FF^{0,k}(\dD,\hH).
$$
Given $\Tref \in \Fred_\FF^\sym(\dD,\hH)$, consider the space
$$
\Fred_\FF^{\sym}(\dD,\hH,\Tref) := \left\{ \Tref + \mathbf{K}
: \dD \to \hH \ \big|\ 
\mathbf{K} \in \Lin_\FF^\sym(\hH) \right\}.
$$
Note that the restriction of each $\mathbf{K} \in \Lin_\FF(\hH)$ to $\dD$
is a compact operator $\dD \to \hH$, thus $\Fred_\FF^{\sym}(\dD,\hH,\Tref)$
has a natural continuous inclusion into $\Fred_\RR^\sym(\dD,\hH)$.
It is also an affine space over $\Lin_\FF^\sym(\hH)$ and
can thus be regarded naturally as a smooth Banach manifold locally
modeled on $\Lin_\FF^\sym(\hH)$; in particular, its tangent spaces are
$$
T_{\mathbf{T}} \left(\Fred_\FF^{\sym}(\dD,\hH,\Tref) \right) =
\Lin_\FF^\sym(\hH).
$$
A remark about the case $\FF=\CC$ is in order:
$\Lin_\CC^\sym(\dD,\hH)$ is a \emph{real}-linear and not a
complex subspace of $\Lin_\CC(\dD,\hH)$, thus
$\Fred_\CC^{\sym}(\dD,\hH,\Tref)$ is a real Banach manifold but does not
carry a natural complex structure.

\begin{lemma}
\label{lemma:symmetric}
For any $\mathbf{T} \in \Lin_\FF^\sym(\dD,\hH)$ that is Fredholm with
index~$0$, $\ker \mathbf{T}$ is the
orthogonal complement of $\im \mathbf{T}$ in~$\hH$, hence
there exist splittings into closed linear subspaces
$$
\dD = V \oplus K, \qquad \hH = W \oplus C
$$
where $K = C = \ker \mathbf{T}$, $W = \im \mathbf{T}$ and
$V = W \cap \dD$.
\end{lemma}
\begin{proof}
If $x \in K := \ker \mathbf{T}$, then symmetry implies
$\langle x , \mathbf{T}y \rangle_\hH = \langle \mathbf{T} x , y \rangle_\hH = 0$
for all $y \in \dD$, hence $K \subset W^\perp$, where $W := \im \mathbf{T}$.
But since $\ind \mathbf{T} = 0$, the dimension of $\ker \mathbf{T}$ equals
the codimension of $\im \mathbf{T}$, implying that
$K$ already has the largest possible dimension for a subspace that
intersects $W$ trivially, and therefore $W \oplus K = \hH$.
Since $K$ is also a subspace of~$\dD$ and the latter is a subspace of~$\hH$,
any $x \in \dD$ can be written uniquely as $x = v + k$ where $k \in K$
and $v \in W \cap \dD =: V$.  The continuous inclusion of $\dD$ into $\hH$
and the fact that $W$ is closed in $\hH$ imply that $V$ is a closed
subspace of~$\dD$.
\end{proof}

We now have the following modification of Prop.~\ref{prop:submanifolds}.

\begin{prop}
\label{prop:submanifoldsSym}
For each integer $k \ge 0$, the subset
$$
\Fred_\FF^{\sym,k}(\dD,\hH,\Tref) := \left\{ \mathbf{T} \in \Fred_\FF^{\sym}(\dD,\hH,\Tref)\ \big|\ 
\dim_\FF \ker \mathbf{T} = k \right\}
$$
is a smooth finite-codimensional Banach submanifold
of $\Fred_\FF^{\sym}(\dD,\hH,\Tref)$, with
$$
\codim_\RR \Fred_\FF^{\sym,k}(\dD,\hH,\Tref) = \begin{cases}
\displaystyle k(k+1) / 2 & \text{ if $\FF = \RR$}, \\
k^2 & \text{ if $\FF = \CC$}.
\end{cases}
$$
In particular, $\Fred_\FF^{\sym,1}(\dD,\hH,\Tref)$ is a submanifold of 
$\Fred_\FF^{\sym}(\dD,\hH,\Tref)$ with codimension~$1$, and moreover,
it carries a canonical co-orientation.
\end{prop}
\begin{proof}
Given $\mathbf{T}_0 \in \Fred_\FF^{\sym,k}(\dD,\hH,\Tref)$, fix the splittings
$\dD = V \oplus K$ and $\hH = W \oplus K$ as in Lemma~\ref{lemma:symmetric}.
Using these in the construction of the map $\Phi$ from \eqref{eqn:Phi}
produces a
neighborhood $\oO \subset \Fred_\FF^{\sym}(\dD,\hH,\Tref)$ of 
$\mathbf{T}_0$ such that,
by Lemma~\ref{lemma:FredCoords}, $\{ \mathbf{T} \in \oO\ |\ 
\dim_\FF \ker \mathbf{T} = k \} = \Phi^{-1}(0)$, where
$$
\Phi : \oO \to \End_\FF(K) :
\begin{pmatrix}
\mathbf{A} & \mathbf{B} \\
\mathbf{C} & \mathbf{D}
\end{pmatrix}
\mapsto \mathbf{D} - \mathbf{C} \mathbf{A}^{-1} \mathbf{B}.
$$
Since the splittings are orthogonal, an element $\mathbf{T} = 
\begin{pmatrix} 
\mathbf{A} & \mathbf{B} \\ \mathbf{C} & \mathbf{D} 
\end{pmatrix} \in \oO$ is symmetric if and only if
\begin{equation*}
\begin{split}
\langle x , \mathbf{A} y \rangle_\hH = \langle \mathbf{A} x , y \rangle_\hH
&\quad \text{ for all $x,y \in V$}, \\
\langle x , \mathbf{D} y \rangle_\hH = \langle \mathbf{D} x , y \rangle_\hH
&\quad \text{ for all $x,y \in K$}, \\
\langle x , \mathbf{B} y \rangle_\hH = \langle \mathbf{C} x , y \rangle_\hH
&\quad \text{ for all $x \in V$, $y \in K$}, \\
\langle x , \mathbf{C} y \rangle_\hH = \langle \mathbf{B} x , y \rangle_\hH
&\quad \text{ for all $x \in K$, $y \in V$}, \\
\end{split}
\end{equation*}
and it follows then that $\Phi(\mathbf{T}) \in \End_\FF^\sym(K)$, where
$\End_\FF^\sym(K) \subset \End_\FF(K)$ is the real vector space of symmetric 
(or Hermitian when $\FF = \CC$) linear maps on $(K,\langle\ ,\ \rangle_\hH)$.
We thus have $\oO \cap \Fred_\FF^{\sym,k}(\dD,\hH,\Tref) = \Phi^{-1}(0)$ 
with $\Phi$ regarded as a smooth map
$\oO \cap \Fred_\FF^{\sym}(\dD,\hH,\Tref) \to \End_\FF^\sym(K)$.  The 
derivative at $\mathbf{T}_0$ again takes the form
$$
d\Phi(\mathbf{T}_0) : \Lin_\FF^\sym(\hH) \to \End_\FF^\sym(K) :
\begin{pmatrix}
\mathbf{A}' & \mathbf{B}' \\
\mathbf{C}' & \mathbf{D}'
\end{pmatrix}
\mapsto \mathbf{D}' ,
$$
where now the block matrix represents an element of
$\Lin_\FF^\sym(\hH)$ with respect to the splitting $\hH = W \oplus K$.
This operator is evidently surjective, hence by the implicit function theorem,
$\Phi^{-1}(0)$ is a smooth Banach submanifold with codimension equal to
$\dim_\RR \End_\FF^\sym(K)$.

Finally, we observe that in the case $k=1$, the above identifies
$\Fred_\FF^{\sym,1}(\dD,\hH,\Tref)$ locally with the zero set of a submersion
to $\End_\FF^\sym(K)$, which is a real $1$-dimensional vector space since
$K$ is a $1$-dimensional vector space over~$\FF$.  
The canonical isomorphism
$$
\RR \to \End_\FF^\sym(K) : a \mapsto a \1
$$ 
thus determines a co-orientation on $\Fred_\FF^{\sym,1}(\dD,\hH,\Tref)$.
\end{proof}

The canonical co-orientation of $\Fred_\FF^{\sym,1}(\dD,\hH,\Tref)$ makes it
natural to define signed intersection numbers between $\Fred_\FF^{\sym,1}(\dD,\hH,\Tref)$ 
and smooth paths in the ambient space
$\Fred_\FF^{\sym}(\dD,\hH,\Tref)$.  The codimensions of
$\Fred_\FF^{\sym,k}(\dD,\hH,\Tref)$ for each $k \ge 2$ are still at least~$3$, hence
large enough to ensure that generic paths or homotopies of paths will
never intersect them.  The following notion is therefore independent of
choices.

\begin{defn}
\label{defn:spectralFlow}
Suppose $\mathbf{T}_+: \mathbf{T}_- \in \Fred_\FF^{\sym}(\dD,\hH,\Tref)$ are 
both Banach space isomorphisms $\dD \to \hH$.  The \defin{spectral flow}
$$
\muspec(\mathbf{T}_-,\mathbf{T}_+) \in \ZZ
$$
from $\mathbf{T}_-$ to $\mathbf{T}_+$ is then defined as the signed count
of intersections of $\mathbf{T} : [-1,1] \to \Fred_\FF^{\sym}(\dD,\hH,\Tref)$ with
$\Fred_\FF^{\sym,1}(\dD,\hH,\Tref)$, where the latter is assumed to carry the
co-orientation given by Prop.~\ref{prop:submanifoldsSym}, and
$\mathbf{T} : [-1,1] \to \Fred_\FF^{\sym}(\dD,\hH,\Tref)$
is any smooth path that is transverse to
$\Fred_\FF^{\sym,k}(\dD,\hH,\Tref)$ for every $k \ge 1$ and satisfies
$\mathbf{T}(\pm 1) = \mathbf{T}_\pm$.
\end{defn}

\subsection{Perturbation of eigenvalues}
\label{sec:eigenvalues}

Continuing in the setting of the previous subsection, we shall now regard each
$\mathbf{T} \in \Fred_\FF^{\sym}(\dD,\hH,\Tref)$ as an unbounded 
operator on $\hH$
with domain~$\dD$, see e.g.~\cite{ReedSimon}*{Chapter~VIII}.
Notice that for each scalar $\lambda \in \FF$, the operator
$\mathbf{T} - \lambda$ also belongs to $\Fred_\FF^{\sym}(\dD,\hH,\Tref)$.
The \defin{spectrum}
$$
\sigma(\mathbf{T}) \subset \FF
$$
of $\mathbf{T}$ is defined as the set of all $\lambda \in \FF$ for which
$\mathbf{T} - \lambda : \dD \to \hH$ does not admit a bounded inverse.
In particular, $\lambda \in \sigma(\mathbf{T})$ is an \defin{eigenvalue}
of $\mathbf{T}$ whenever $\mathbf{T} - \lambda : \dD \to \hH$ has nontrivial
kernel, and the dimension of this kernel is called the
\defin{multiplicity} of the eigenvalue.  We call $\lambda$ a
\defin{simple eigenvalue} if it has multiplicity~$1$.
By a standard argument familiar to both mathematicians and
physicists, the eigenvalues of a symmetric complex-linear operator
are always real.

\begin{remark}
\label{remark:selfAdjoint}
If $\dD \subset \hH$ is dense, then the \defin{adjoint} of $\mathbf{T}$ 
is defined as an unbounded operator $\mathbf{T}^*$ with domain $\dD^*$ 
satisfying
$$
\langle x , \mathbf{T} y \rangle_\hH = \langle \mathbf{T}^* x , y \rangle_\hH
\quad \text{ for all $x \in \dD^*$, $y \in \dD$},
$$
where $\dD^*$ is the set of all $x \in \hH$ such that there exists
$z \in \hH$ satisfying $\langle x, \mathbf{T} y \rangle_\hH = 
\langle z,y \rangle_\hH$ for all $y \in \dD$.  One says that $\mathbf{T}$
is \defin{self-adjoint} if $\mathbf{T} = \mathbf{T}^*$, which means both
that $\mathbf{T}$ is symmetric and $\dD = \dD^*$.  In many applications
(e.g.~in Exercise~\ref{EX:selfAdjoint}),
the latter amounts to a condition on ``regularity of weak solutions''.
This condition implies that the inclusion $\ker \mathbf{T} \hookrightarrow
(\im \mathbf{T})^\perp$---valid for all symmetric operators---is also
surjective, so if $\mathbf{T} : \dD \to \hH$ is Fredholm, it is then
automatic that $\ind(\mathbf{T}) = 0$.
\end{remark}

\begin{prop}
\label{prop:eigenvalues}
Assume $\mathbf{T}_0 \in \Fred_\FF^{\sym}(\dD,\hH,\Tref)$.  Then:
\begin{enumerate}
\item \label{item:finite}
Every $\lambda \in \sigma(\mathbf{T}_0)$ is an eigenvalue with finite
multiplicity.
\item \label{item:discrete}
The spectrum $\sigma(\mathbf{T}_0)$ is a discrete subset of~$\RR$.
\item \label{item:perturb}
Suppose $\lambda_0 \in \sigma(\mathbf{T}_0)$ is an eigenvalue with
multiplicity $m \in \NN$ and $\epsilon > 0$ is chosen such that
no other eigenvalues lie in $[\lambda_0 - \epsilon,\lambda_0 + \epsilon]$.
Then $\mathbf{T}_0$ has a neighorhood $\oO \subset \Fred_\FF^{\sym}(\dD,\hH,\Tref)$
such that for all $\mathbf{T} \in \oO$,
$$
\sum_{\lambda \in \sigma(\mathbf{T}) \cap [\lambda_0 - \epsilon,\lambda_0 + \epsilon]}
m(\lambda) = m,
$$
where $m(\lambda) \in \NN$ denotes the multiplicity of 
$\lambda \in \sigma(\mathbf{T})$.
\end{enumerate}
\end{prop}
\begin{proof}
For every $\lambda \in \FF$, $\mathbf{T}_0 - \lambda$ is a Fredholm operator
with index~$0$, so it is a Banach space isomorphism $\dD \to \hH$ and thus
has a bounded inverse if and only if its kernel is trivial.
The Fredholm property also implies that the kernel is finite dimensional
whenever it is nontrivial, so this proves~\eqref{item:finite}.

For \eqref{item:discrete} and \eqref{item:perturb}, let us assume $\FF = \CC$,
as the case $\FF = \RR$ will follow by taking complexifications of real
vector spaces.  We claim therefore that $\sigma(\mathbf{T}_0)$ is a
discrete subset of~$\CC$.  To see this, suppose $\lambda_0 \in \RR$ is an
eigenvalue of $\mathbf{T}_0$ with multiplicity~$m$, so
$$
\mathbf{T}_0 - \lambda_0 \in \Fred_\CC^{\sym,m}(\dD,\hH).
$$
By Lemma~\ref{lemma:symmetric}, there are splittings $\dD = V \oplus K$
and $\hH = W \oplus K$ with $K = \ker (\mathbf{T}_0 - \lambda_0)$,
$W = \im (\mathbf{T}_0 - \lambda_0)$ and $V = W \cap \dD$.  Any scalar
$\lambda \in \CC$ appears in block-diagonal form 
$\begin{pmatrix} \lambda & 0 \\ 0 & \lambda \end{pmatrix}$ with respect
to these splittings, and the block form for $\mathbf{T}_0$ is thus
$$
\mathbf{T}_0 = \begin{pmatrix} \mathbf{A}_0 + \lambda_0 & 0 \\
0 & \lambda_0
\end{pmatrix}
$$
for some Banach space isomorphism $\mathbf{A}_0 : V \to W$.  Writing
nearby operators $\mathbf{T} \in \Fred_\CC(\dD,\hH)$ as
$\begin{pmatrix} \mathbf{A} & \mathbf{B} \\ 
\mathbf{C} & \mathbf{D} \end{pmatrix}$, we can imitate
the construction in \eqref{eqn:Phi} to produce
neighborhoods $\oO(\mathbf{T}_0) \subset \Fred_\CC(\dD,\hH)$ of $\mathbf{T}_0$ 
and $\DD_\epsilon(\lambda_0) \subset \CC$ of $\lambda_0$, admitting a
holomorphic map
$$
\Phi : \oO(\mathbf{T}_0) \times \DD_\epsilon(\lambda_0) \to \End_\CC(K) :
(\mathbf{T},\lambda) \mapsto \left( \mathbf{D} - \lambda \right) -
\mathbf{C} \left( \mathbf{A} - \lambda \right)^{-1} \mathbf{B}
$$
such that $\ker (\mathbf{T} - \lambda) \cong \ker \Phi(\mathbf{T},\lambda)$.
The set of eigenvalues of $\mathbf{T}_0$ near $\lambda_0$ is then the
zero set of the holomorphic function
\begin{equation}
\label{eqn:det}
\DD_\epsilon(\lambda_0) \to \CC : \lambda \mapsto \det \Phi(\mathbf{T}_0,\lambda).
\end{equation}
This function cannot be identically zero since there are no eigenvalues
outside of~$\RR$, thus the zero at $\lambda_0$ is isolated,
proving~\eqref{item:discrete}.

To prove~\eqref{item:perturb}, note finally that if the neighborhood
$\oO(\mathbf{T}_0) \subset \Fred_\CC(\dD,\hH)$ of $\mathbf{T}_0$ is
sufficiently small, then for every $\mathbf{T} \in \oO(\mathbf{T}_0)$,
the holomorphic function
$$
f_\mathbf{T} :
\DD_\epsilon(\lambda_0) \to \CC : \lambda \mapsto \det \Phi(\mathbf{T},\lambda)
$$
has the same algebraic count of zeroes in $\DD_\epsilon(\lambda_0)$,
all of which lie in $[\lambda_0 - \epsilon,\lambda_0 + \epsilon]$ if
$\mathbf{T}$ is symmetric.  Observe moreover that since
$$
\p_\lambda \Phi(\mathbf{T}_0,\lambda_0) = -\1 \in \End_\CC(K),
$$
we are free to assume after possibly shrinking 
$\epsilon$ and $\oO(\mathbf{T}_0)$ that $\p_\lambda \Phi(\mathbf{T},\lambda)$
is always a nonsingular transformation in $\End_\CC(K)$.
Since $\Phi(\mathbf{T},\lambda)$ is in $\End_\CC^\sym(K)$ and thus
diagonalizable whenever $\mathbf{T}$ is symmetric and $\lambda \in \RR$,
it follows via Exercise~\ref{EX:geomMult} below that the order of any
zero $f_\mathbf{T}(\lambda) = 0$ is precisely the multiplicity of $\lambda$
as an eigenvalue of~$\mathbf{T}$.

\end{proof}

\begin{exercise}
\label{EX:geomMult}
Suppose $\uU \subset \CC$ is an open subset, $A : \uU \to \CC^{n \times n}$
is a holomorphic map and $z_0 \in \uU$ is a point at which
$A(z_0)$ is noninvertible but diagonalizable, and $A'(z_0) \in \GL(n,\CC)$.  
Show that $\dim_\CC \ker A(z_0)$ is the order of the zero of the holomorphic
function $\det A : \uU \to \CC$ at~$z_0$.
\end{exercise}

The next result implies that for a generic path of symmetric index~$0$ 
operators as appears in our definition of $\muspec(\mathbf{T}_-,\mathbf{T}_+)$,
the spectral flow is indeed a signed count of eigenvalues crossing~$0$.

\begin{prop}
\label{prop:crossings}
Suppose $\left\{\mathbf{T}_t \in \Fred_\FF^{\sym}(\dD,\hH,\Tref)\right\}_{t \in (-1,1)}$ is a 
smooth path and $\lambda_0 \in \RR$ is a simple eigenvalue of~$\mathbf{T}_0$.  
Then:
\begin{enumerate}
\item \label{item:perturbEigenvalues}
For sufficiently small $\epsilon > 0$, there exists a unique smooth
function $\lambda : (-\epsilon,\epsilon) \to \RR$ such that
$\lambda(0) = \lambda_0$ and $\lambda(t)$ is a simple eigenvalue of
$\mathbf{T}_t$ for each $t \in (-\epsilon,\epsilon)$.
\item \label{item:sign}
The derivative $\lambda'(0)$ is nonzero 
if and only if the intersection of the path
$\left\{\mathbf{T}_t - \lambda_0 \in \Fred_\FF^{\sym}(\dD,\hH,\Tref) \right\}_{t \in (-1,1)}$
with $\Fred_\FF^{\sym,1}(\dD,\hH,\Tref)$ at $t=0$ is transverse, and the sign
of $\lambda'(0)$ is then the sign of the intersection.
\end{enumerate}
\end{prop}
\begin{proof}
Using the same construction as in the proof of 
Proposition~\ref{prop:eigenvalues}, we can find small numbers $\epsilon > 0$
and $\delta > 0$ such that
$$
\left\{ (t,\lambda) \in (-\epsilon,\epsilon) \times (\lambda_0 - \delta, \lambda_0 + \delta)\ \big|\ 
\lambda \in \sigma(\mathbf{T}_t) \right\} = \Phi^{-1}(0),
$$
where
$$
\Phi : (-\epsilon,\epsilon) \times (\lambda_0 - \delta, \lambda_0 + \delta) \to
\End_\FF^\sym(K) : (t,\lambda) \mapsto 
\left( \mathbf{D}_t - \lambda \right) - \mathbf{C}_t \left( \mathbf{A}_t - \lambda \right)^{-1}
\mathbf{B}_t,
$$
and we write $\mathbf{T}_t = \begin{pmatrix} \mathbf{A}_t & \mathbf{B}_t \\
\mathbf{C}_t & \mathbf{D}_t \end{pmatrix}$
with respect to splittings $\dD = V \oplus K$ and $\hH = W \oplus K$ with
$K = \ker (\mathbf{T}_0 - \lambda_0)$, $W = \im (\mathbf{T}_0 - \lambda_0)$
and $V = W \cap \dD$.  In saying this, we've implicitly used the assumption that
$\lambda_0$ is a simple eigenvalue, as it follows that
$\dim_\FF \ker (\mathbf{T} - \lambda)$ cannot be larger than $1$ for any
$\mathbf{T}$ near $\mathbf{T}_0$ and $\lambda$ near~$\lambda_0$, so that
$\Phi^{-1}(0)$ catches all nearby eigenvalues.
Simplicity also means that $\End_\FF^\sym(K)$ is real $1$-dimensional,
and we have
$$
\p_t\Phi(0,\lambda_0) = \p_t \mathbf{D}_t|_{t=0}, \qquad
\p_\lambda \Phi(0,\lambda_0) = -1.
$$
The implicit function theorem thus gives $\Phi^{-1}(0)$ near $(0,\lambda_0)$
the structure of a smooth $1$-manifold with tangent space at $(0,\lambda_0)$
spanned by the vector
$$
\p_t + \left(\p_t \mathbf{D}_t|_{t=0}\right) \p_\lambda,
$$
where we are identifying $\p_t \mathbf{D}_t|_{t=0} \in \End_\FF^\sym(K)$ with
a real number via the natural isomorphism $\End_\FF^\sym(K) = \RR$. 
Therefore $\Phi^{-1}(0)$ can be written as the graph of a uniquely determined smooth function $\lambda$, whose derivative at zero is a multiple of $\p_t \mathbf{D}_t|_{t=0}$.
This proves both statements in the proposition, since by the proof of
Proposition~\ref{prop:submanifoldsSym}, the intersection of
$\{ \mathbf{T}_t \}_{t \in (-1,1)}$ with $\Fred_\FF^{\sym,1}(\dD,\hH,\Tref)$ is
transverse if and only if $\p_t \mathbf{D}_t|_{t=0} \ne 0$, and its sign is
then the sign of $\p_t \mathbf{D}_t|_{t=0}$.
\end{proof}

The purpose of the next lemma is to prevent
eigenvalues from escaping to $\pm\infty$ under
smooth families of operators in $\Fred_\FF^\sym(\dD,\hH,\Tref)$.

\begin{lemma}
\label{lemma:bounded}
Suppose 
$\left\{ \mathbf{K}_t \in \Lin_\FF^\sym(\hH) \right\}_{t \in (a,b)}$
is a smooth path of symmetric bounded linear operators, and
$\lambda : (a,b) \to \RR$ is a smooth function such that 
for every $t \in (a,b)$, $\lambda(t)$ is a simple eigenvalue of 
$\mathbf{T}_t := \Tref + \mathbf{K}_t \in \Fred_\FF^{\sym}(\dD,\hH,\Tref)$.
Then 
$$
|\dot{\lambda}(t)| \le \| \p_t \mathbf{K}_t \|_{\Lin(\hH)} 
\quad \text{ for all $t \in (a,b)$}.
$$ 
\end{lemma}
\begin{proof}
Since $\{\mathbf{T}_t - \lambda(t) \in \Fred_\FF^{\sym}(\dD,\hH,\Tref) \}_{t \in (a,b)}$ is 
a smooth family of operators in $\Fred_\FF(\dD,\hH)$ with $1$-dimensional
kernel, one can use the local families of isomorphisms
$\Psi(\mathbf{T}_t - \lambda(t)) \in \Lin_\FF(\dD)$ from 
Lemma~\ref{lemma:FredCoords} to find a smooth family of eigenvectors
$x(t) \in \ker(\mathbf{T}_t - \lambda(t))$ for $t \in (a,b)$.
Normalize these so that $\| x(t) \|_\hH = 1$ for all~$t$.  Then
$0 = \p_t \langle x(t),x(t) \rangle_\hH = \langle \dot{x}(t),x(t) \rangle_\hH
+ \langle x(t),\dot{x}(t) \rangle_\hH$ and
$\lambda(t) = \langle x(t) , \mathbf{T}_t x(t) \rangle_\hH$, so writing
$\dot{\mathbf{K}}_t := \p_t \mathbf{K}_t = \p_t \mathbf{T}_t$, we have
\begin{equation*}
\begin{split}
\dot{\lambda}(t) &= \p_t \langle x(t) , \mathbf{T}_t x(t) \rangle_\hH
= \langle x(t) , \dot{\mathbf{K}}_t x(t) \rangle_\hH +
\langle \dot{x}(t) , \mathbf{T}_t x(t) \rangle_\hH + 
\langle x(t) , \mathbf{T}_t \dot{x}(t) \rangle_\hH \\
&= \langle x(t) , \dot{\mathbf{K}}_t x(t) \rangle_\hH,
\end{split}
\end{equation*}
as the last two terms in the first line become
$\lambda(t) \left[ \langle \dot{x}(t) , x(t) \rangle_\hH +
\langle x(t) , \dot{x}(t) \rangle_\hH \right] = 0$ since $\mathbf{T}_t$
is symmetric and $\mathbf{T}_t x(t) = \lambda(t) x(t)$.  We obtain
$$
|\dot{\lambda}(t)| \le \| x(t) \|_\hH \| \dot{\mathbf{K}}_t \|_{\Lin(\hH)} \| x(t) \|_\hH
= \| \dot{\mathbf{K}}_t \|_{\Lin(\hH)}.
$$
\end{proof}

\subsection{Homotopies of eigenvalues}
\label{sec:specHomotopy}

Specializing further, we now set $\hH$ and $\dD$ equal to the specific real
Hilbert spaces
$$
\hH := L^2(S^1,\RR^{2n}), \qquad \dD := H^1(S^1,\RR^{2n}),
$$
and set $\Tref := -J_0\, \p_t$, where $J_0$ denotes the standard complex 
structure on $\RR^{2n} = \CC^n$.
Observe that any bounded linear operator on $L^2$ determines
a compact operator $H^1 \to L^2$ via composition with the compact inclusion.
In particular, we shall consider compact perturbations of $- J_0\, \p_t$ 
in the form
\begin{equation}
\label{eqn:A}
\mathbf{A} = -J_0 \, \p_t - S(t)
\end{equation}
with $S : S^1 \to \End_\RR^\sym(\RR^{2n})$ smooth.  It is straightforward to
check that this operator is symmetric with respect to the $L^2$-product
since $S(t)$ is symmetric for every~$t$.  The following then implies that
$\mathbf{A} \in \Fred_\RR^\sym(\dD,\hH,\Tref)$.

\begin{lemma}
\label{lemma:itIsFredholm}
The operator $-J_0\, \p_t : H^1(S^1,\RR^{2n}) \to L^2(S^1,\RR^{2n})$
is Fredholm with index~$0$. 
\end{lemma}
\begin{proof}
Since $J_0$ defines an isomorphism, it suffices actually to show that the ordinary differential operator
$$
\p_t : H^1(S^1,\RR^{2n}) \to L^2(S^1,\RR^{2n})
$$
is Fredholm with index~$0$. The kernel of this operator is the space of constant functions $S^1 \rightarrow \RR^{2n}$, which has dimension $2n$. To compute the dimension of the cokernel, we observe that if $f=\p_t F$ lies in the image of this operator, we have $\int_{S^1}f(t) \, dt=0$ since $F$ is periodic in $t$. Conversely, 
if $\int_{S^1}f(t) \, dt=0$
with $f \in L^2(S^1,\RR^{2n})$, then the function $F(s)=\int_0^s f(t)\, dt$ is periodic in $s$
and defines an element of $H^1(S^1,\RR^{2n})$ satisfying
$\p_t F = f$. Hence the image of $\p_t$ is exactly the set 
$$
\im(\p_t)=\left\{f \in L^2(S^1,\RR^{2n}) \ \bigg|\ \int_{S^1}f(t) \, dt=0\right\},
$$ 
which has codimension $2n$.
\end{proof}

The proof of Theorem~\ref{thm:spectralFlow} requires only one more technical
ingredient, whose proof is given in Appendix~\ref{app:spectral} and should
probably be skipped on first reading unless you have already read
Lecture~\ref{lec:transversality} or seen similar applications of the
Sard-Smale theorem.  You might however find the result 
plausible in accordance with the notion that maps from $2$-dimensional domains,
such as a map of the form
$$
(-1,1) \times \RR \to \Fred_\RR^{\sym}(\dD,\hH,\Tref) : (t,\lambda) \mapsto
\mathbf{T}_t - \lambda
$$
should \emph{generically} not intersect submanifolds that have 
codimension~$3$ or more, such as $\Fred_\RR^{\sym,k}(\dD,\hH,\Tref)$ when
$k \ge 2$.

\begin{lemma}
\label{lemma:generic}
Fix a smooth map $S : [-1,1] \times S^1 \to \End_\RR^\sym(\RR^{2n})$ and
consider the $1$-parameter family of unbounded linear operators
$$
\mathbf{A}_s := -J_0 \, \p_t - S(s,\cdot) : L^2(S^1,\RR^{2n}) \supset
H^1(S^1,\RR^{2n}) \to L^2(S^1,\RR^{2n})
$$
for $s \in [-1,1]$.  One can arrange after a $C^\infty$-small perturbation of 
$S$ fixed at $s= \pm 1$ that the following conditions hold:
\begin{enumerate}
\item For each $s \in (-1,1)$, all eigenvalues
of $\mathbf{A}_s$ are simple.
\item All intersections of the path
$$
(-1,1) \to \Fred_\RR^\sym(\dD,\hH,\Tref) : s \mapsto \mathbf{A}_s
$$
with $\Fred_\RR^{\sym,1}(\dD,\hH,\Tref)$ are transverse.
\end{enumerate}
\qed
\end{lemma}

\begin{proof}[Proof of Theorem~\ref{thm:spectralFlow}]
Given a smooth family $\{ \mathbf{A}_s \}_{s \in [-1,1]}$ as stated
in the theorem, use Lemma~\ref{lemma:generic} to obtain a $C^\infty$-small 
perturbation for which the eigenvalues are simple for $s \in (-1,1)$ and
all intersections with $\Fred_\RR^{\sym,1}(\dD,\hH)$ are transverse.
Proposition~\ref{prop:crossings} then implies that the eigenvalues
depend smoothly on $s$, and Lemma~\ref{lemma:bounded} imposes a uniform
bound on their derivatives with respect to~$s$ so that each one varies only
in a bounded subset of $\RR$ for $s \in (-1,1)$.  The smooth families of
eigenvalues for $s \in (-1,1)$ therefore extend to continuous families
for $s \in [-1,1]$ since the space of noninvertible Fredholm operators
with index~$0$ is closed.  Proposition~\ref{prop:eigenvalues} ensures
moreover that these continuous families hit every eigenvalue with the 
correct multiplicity at $s = \pm 1$, and by
Proposition~\ref{prop:crossings}, the formula for
$\muspec(\mathbf{A}_-,\mathbf{A}_+)$ stated in the theorem is correct for
the perturbed family with simple eigenvalues and transverse crossings.
To obtain the same result for the original family, suppose we have a
sequence of perturbations $\{ \mathbf{A}_s^\nu \}_{s \in [-1,1]}$
converging in $C^\infty$ as 
$\nu \to \infty$ to $\{ \mathbf{A}_s \}_{s \in [-1,1]}$.
Lemma~\ref{lemma:bounded} then provides a uniform $C^1$-bound for each
sequence of smooth families of eigenvalues, so they have $C^0$-convergent
subsequences as $\nu \to \infty$, giving rise to the continuous families
in the statement of the theorem.
\end{proof}

\begin{remark}
\label{remark:noLift}
It is important to understand that the definition of spectral flow depends
on the particular co-orientation of $\Fred_\FF^{\sym,1}(\dD,\hH,\Tref)$ 
that arose in the proof of Prop.~\ref{prop:submanifoldsSym}; we saw in
Prop.~\ref{prop:crossings} that this is indeed the \emph{right} co-orientation
to use if we want to interpret signed intersections with $\Fred_\FF^{\sym,1}(\dD,\hH,\Tref)$
as signed crossing numbers of eigenvalues.  In the non-symmetric setting of
\S\ref{sec:submanifolds}, one can show that $\Fred_\RR^{0,1}(X,Y)$ is also
co-orientable; this is obvious in the finite-dimensional case since
$\Fred_\RR^{0,1}(\RR^n,\RR^n)$ is then a regular level set of the determinant
function.  Moreover, $\Fred_\RR^{0,1}(\RR^n,\RR^n)$ is connected
(see Exercise~\ref{EX:FredConnected} below), so the
co-orientation is unique up to a sign.  One can therefore lift the
$\ZZ_2$-valued spectral flow of \S\ref{sec:submanifolds} to~$\ZZ$, but as
in Exercise~\ref{EX:detSign}, the result will be a different and much
less interesting invariant than $\muspec(A_-,A_+)$,
as its value will always be either $0$ (if $\det A_-$ and $\det A_+$
have the same sign) or~$\pm 1$ (if they don't).  The reason for the discrepancy
is that the canonical co-orientation of $\Fred_\RR^{\sym,1}(\dD,\hH,\Tref)$
must generally differ on some connected components from any possible
co-orientation of the larger hypersurface
$\Fred_\RR^{0,1}(\dD,\hH) \subset \Fred_\RR^0(\dD,\hH)$.
\end{remark}
\begin{exercise}
\label{EX:FredConnected}
Show that the space $\Fred_\RR^{0,1}(\RR^2,\RR^2)$ of rank~$1$ matrices in $\RR^{2 \times 2}$ is connected,
but the space $\Fred_\RR^{\sym,1}(\RR^2,\RR^2)$ of \emph{symmetric} rank~$1$ matrices is not, and that the
canonical co-orientation of $\Fred_\RR^{\sym,1}(\RR^2,\RR^2)$ coming from
Prop.~\ref{prop:submanifoldsSym} differs on some components from any possible
co-orientation of $\Fred_\RR^{0,1}(\RR^2,\RR^2) \subset \RR^{2 \times 2}$.
\textsl{Hint: A non-symmetric $2$-by-$2$ matrix may have rank~$1$ even if
both of its eigenvalues are~$0$.  For symmetric matrices this cannot happen.}
\end{exercise}
\begin{exercise}
\label{EX:noLift}
Find a smooth path $A : [-1,1] \to \RR^{2 \times 2}$ of symmetric matrices such that
$A_\pm := A(\pm 1)$ are both invertible and $\muspec(A_-,A_+) = 2$, but
$A_+$ and $A_-$ can also be connected by a smooth path of (not necessarily symmetric)
invertible matrices in~$\RR^{2 \times 2}$.
\end{exercise}

\section{The Hessian of the contact action functional}
\label{sec:contactHessian}

Before returning to contact geometry, let's quickly revisit the Floer
homology for a time-dependent Hamiltonian $\{ H_t : M \to \RR \}_{t \in S^1}$
on a symplectic manifold $(M,\omega)$.  In Lecture~\ref{lec:intro}, we introduced the
symplectic action functional $\aA_H : C^\infty_\contr(S^1,M) \to \RR$
and wrote down the formula
$$
\nabla \aA_H(\gamma) = J_t(\gamma) \left( \dot{\gamma} - X_t(\gamma) \right)
\in \Gamma(\gamma^*TM) =: T_\gamma C^\infty_\contr(S^1,M)
$$
for the ``unregularized'' gradient of $\aA_H$ at a contractible loop
$\gamma \in C^\infty_\contr(S^1,M)$.  Here $X_t$ denotes the Hamiltonian
vector field and $J_t$ is a time-dependent family of compatible almost
complex structures, which determines the $L^2$-product
$$
\langle \eta_1 , \eta_2 \rangle_{L^2} = \int_{S^1} \omega(\eta_1(t),J_t \eta_2(t))\, dt.
$$
The critical points of $\aA_H$ are the loops $\gamma$ such that
$\nabla \aA_H(\gamma) = 0$.  Formally, the Hessian of $\aA_H$ at 
$\gamma \in \Crit(\aA_H)$ is the ``linearization of $\aA_H$ at $\gamma$,''
which gives a linear operator
$$
\mathbf{A}_\gamma := \nabla^2 \aA_H(\gamma) : \Gamma(\gamma^*TM) \to \Gamma(\gamma^*TM).
$$
To write it down, one can choose any connection $\nabla$ on $M$, and choose
for $\eta \in \Gamma(\gamma^*TM)$ a smooth family 
$\{ \gamma_\rho : S^1 \to M \}_{\rho \in (-\epsilon,\epsilon)}$ with
$\gamma_0 = \gamma$ and $\p_\rho \gamma_\rho|_{\rho=0} = \eta$, and then compute
$$
\mathbf{A}_\gamma \eta := \left. \nabla_\rho \left[ \nabla \aA_H(\gamma_\rho)
\right] \right|_{\rho=0}.
$$
The result is independent of the choice of connection since
$\nabla \aA_H(\gamma) = 0$.

\begin{exercise}
\label{EX:FloerA}
Show that if the connection $\nabla$ on $M$ is chosen to be symmetric,
then $\mathbf{A}_\gamma \eta = J_t (\nabla_t \eta - \nabla_\eta X_t)$.
\end{exercise}

We now introduce the class of symmetric operators that appear in
asymptotic formulas in SFT.  Fix a $(2n-1)$-dimensional contact manifold
$(M,\xi)$ with contact form $\alpha$, induced Reeb vector field $R_\alpha$,
and a complex structure $J : \xi \to \xi$ compatible with the symplectic
structure $d\alpha|_{\xi}$.  Let
$$
\pi_\xi : TM \to \xi
$$
denote the projection along~$R_\alpha$.
The \defin{contact action functional} is defined by
$$
\aA_\alpha : C^\infty(S^1,M) \to \RR : \gamma \mapsto \int_{S^1} \gamma^*\alpha.
$$
The first variation of this functional for $\gamma \in C^\infty(S^1,M)$ and
$\eta \in \Gamma(\gamma^*TM)$ is
$$
d \aA_\alpha(\gamma) \eta = \int_{S^1} d\alpha(\eta,\dot{\gamma})\, dt =
- \int_{S^1} d\alpha(\pi_\xi \dot{\gamma} , \eta) \, dt.
$$
The functional has a built-in degeneracy since it is parametrization-invariant;
in particular, $d\aA_\alpha(\gamma) \eta = 0$ whenever $\eta$ points in the
direction of the Reeb vector field, a symptom of the fact that closed
Reeb orbits always come in families related to each other by reparametrization.
A loop $\gamma : S^1 \to M$ is critical for $\aA_\alpha$ if and only if 
$\dot{\gamma}$ is everywhere tangent to $R_\alpha$, allowing for an
infinite-dimensional family of distinct perturbations---however, there exist
preferred parametrizations, namely those for which $\dot{\gamma}$ is a
\emph{constant} multiple of $R_\alpha$, meaning
\begin{equation}
\label{eqn:preferred}
\dot{\gamma} = T \cdot R_\alpha(\gamma), \qquad T := \aA_\alpha(\gamma).
\end{equation}
Such a loop corresponds to a $T$-periodic solution $x : \RR \to M$ to
$\dot{x} = R_\alpha(x)$, where $\gamma(t) = x(Tt)$.

The discussion above indicates that we cannot derive a ``Hessian'' of
$\aA_\alpha$ in the same straightforward way as in Floer homology, as the
resulting operator will always have nontrivial kernel due to the degeneracy
in the $R_\alpha$ direction.  To avoid this, we shall consider only
preferred parametrizations $\gamma : S^1 \to M$ of the form 
\eqref{eqn:preferred}, and
perturbations in directions tangent to $\xi$, which is transverse to every
Reeb orbit.  For $\eta \in \Gamma(\gamma^*\xi)$, we then have
$$
d \aA_\alpha(\gamma) \eta = \int_{S^1} d\alpha(-J \pi_\xi \dot{\gamma}, J\eta)\, dt
= \langle -J \pi_\xi \dot{\gamma} , \eta \rangle_{L^2},
$$
where we define an $L^2$-product for sections of $\gamma^*\xi$ by
\begin{equation}
\label{eqn:L2xi}
\langle \eta,\eta' \rangle_{L^2} := \int_{S^1} d\alpha(\eta, J \eta')\, dt.
\end{equation}
It therefore seems sensible to write
$$
\nabla \aA_\alpha(\gamma) := - J \pi_\xi \dot{\gamma} \in \Gamma(\gamma^*\xi),
$$
and we shall define the Hessian at a critical point $\gamma$ as the 
linearization in $\xi$ directions, i.e.
$$
\nabla^2 \aA_\alpha(\gamma) : \Gamma(\gamma^*\xi) \to \Gamma(\gamma^*\xi).
$$
Given $\eta \in \Gamma(\gamma^*\xi)$, choose a smooth family
$\{ \gamma_\rho : S^1 \to M \}_{\rho \in (-\epsilon,\epsilon)}$ with
$\gamma_0 = \gamma$ and $\p_\rho \gamma_\rho|_{\rho=0} = \eta$, and fix
a symmetric connection $\nabla$ on~$M$.  Since $\pi_\xi \dot{\gamma} = 0$,
the covariant derivative of
$\nabla \aA_\alpha(\gamma_\rho)$ at $\rho=0$ is then
\begin{equation*}
\begin{split}
\left. \nabla_\rho\left( -J \pi_\xi \dot{\gamma}_\rho \right) \right|_{\rho=0}
&= -J \left.\nabla_\rho\left(\pi_\xi \dot{\gamma}_\rho \right)\right|_{\rho=0}
= -J \left.\nabla_\rho\left[\dot{\gamma}_\rho - \alpha(\dot{\gamma}_\rho) R_\alpha(\gamma_\rho)\right]\right|_{\rho=0} \\
&= -J \left(\nabla_t \eta - T \nabla_\eta R_\alpha - \p_\rho\left[ 
\alpha(\dot{\gamma}_\rho) \right]|_{\rho=0} \cdot R_\alpha(\gamma) \right).
\end{split}
\end{equation*}
In the last term, we can write $\p_\rho\left[ \alpha(\dot{\gamma}_\rho) \right]|_{\rho=0}
= d\alpha(\eta,\dot{\gamma}) + \p_t\left[ \alpha(\eta) \right] = 0$
since $\dot{\gamma} = T R_\alpha(\gamma)$ and $\alpha(\eta) = 0$ for
$\eta \in \Gamma(\gamma^*\xi)$.  One can now check that the remaining terms
define a section of $\gamma^*\xi$, thus we are led to the following definition.

\begin{defn}
\label{defn:asymptotic}
Given a loop $\gamma : S^1 \to M$ parametrizing a closed Reeb orbit
in $(M,\xi = \ker \alpha)$ with period $T \equiv \alpha(\dot{\gamma})$,
the \defin{asymptotic operator} associated to $\gamma$ is the
first-order differential operator on $\gamma^*\xi$ defined by
$$
\mathbf{A}_\gamma : \Gamma(\gamma^*\xi) \to \Gamma(\gamma^*\xi) :
\eta \mapsto -J (\nabla_t \eta - T \nabla_\eta R_\alpha)
$$
\end{defn}

\begin{exercise}
\label{EX:symmetric}
Show that $\mathbf{A}_\gamma$ is symmetric with respect to
the $L^2$ inner product \eqref{eqn:L2xi} on $\Gamma(\gamma^*\xi)$.
Moreover, $\gamma$ is nondegenerate (see \S\ref{sec:Weinstein}) if and only if
$\ker \mathbf{A}_\gamma$ is trivial.
\textsl{Hint for nondegeneracy: Consider the pullback of $\gamma^*\xi$ via
the cover $\RR \to S^1 = \RR / \ZZ$, and show that solutions to
$\nabla_t \eta - T \nabla_\eta R_\alpha=0$ on the pullback are given by
operating on $\xi_{\gamma(0)}$ with the linearized Reeb flow.  To see this,
try differentiating families of solutions to the equation
$\dot{x} = T R_\alpha(x)$.}
\end{exercise}

\begin{remark}
\label{remark:connection}
Another way of phrasing the hint in the the above exercise is as follows:
$\mathbf{A}_\gamma$ can also be written as
$-J \widehat{\nabla}_t$, where $\widehat{\nabla}_t$ is the unique
\emph{symplectic connection} on $(\gamma^*\xi,d\alpha)$ for which parallel
transport is given by the linearized Reeb flow.
\end{remark}

You might be slightly concerned about the sign difference between
the two formulas we've derived for asymptotic operators in contact geometry
and in Floer homology.  I also find this troubling, but the discrepancy  seems
to originate from the fact that our account of Floer homology has referred
always to the \emph{negative} gradient flow of $\aA_H$, while SFT is
actually defined via the \emph{positive} gradient flow of~$\aA_\alpha$.
The words ``gradient flow'' in SFT must in any case be interpreted very
loosely.  If 
$$
u : [0,\infty) \times S^1 \to \RR \times M
$$
is the cylindrical end of a finite-energy $J$-holomorphic curve for some
$J \in \jJ(\alpha)$ as we described in Lecture~\ref{lec:intro}, then $u(s,t)$ does not
satisfy anything so straightforward as $\p_s - \nabla \aA_\alpha(u(s,\cdot)) = 0$,
but it does satisfy
$$
\pi_\xi \p_s u + J \pi_\xi \p_t u = 0,
$$
which can be interpreted as the projection of a positive gradient flow equation
to the contact bundle.  This observation is a local symptom of a more important
global fact that follows from Stokes' theorem: 
any asymptotically cylindrical $J$-holomorphic curve 
$u : \dot{\Sigma} \to \RR \times M$ with positive and negative punctures
$\Gamma^\pm$ asymptotic to orbits $\{ \gamma_z \}_{z \in \Gamma^\pm}$ satisfies
$$
\sum_{z \in \Gamma^+} \aA_\alpha(\gamma) - \sum_{z \in \Gamma^-} \aA_\alpha(\gamma)
= \int_{\dot{\Sigma}} u^*d\alpha \ge 0.
$$
This generalizes the basic fact in Floer homology that flow lines decrease
action and, conversely, have their energy controlled by the action.

We would now like to develop some of the general properties of
asymptotic operators.  Recall that on any symplectic vector bundle 
$(E,\omega)$, a compatible complex structure $J$ determines a Hermitian
inner product
$$
\langle v , w \rangle = \omega( v , J w) + i \omega(v,w),
$$
and conversely, any Hermitian inner product on a complex vector bundle
determines a symplectic structure via the same relation.  For this reason,
we shall refer to any vector bundle $E$ with a compatible pair $(J,\omega)$
as a \defin{Hermitian vector bundle}.  A \defin{unitary trivialization}
of such a bundle is a trivialization that identifies fibers with
$\RR^{2n}$ such that $J$ and $\omega$ become the standard complex structure
$J_0$ and symplectic structure $\omega_0$ respectively.

\begin{defn}
\label{defn:asymptotic2}
Fix a Hermitian vector bundle $(E,J,\omega)$ over $S^1$.  An
\defin{asymptotic operator} on $(E,J,\omega)$ is any real-linear differential
operator $\mathbf{A} : \Gamma(E) \to \Gamma(E)$ that takes the form
\begin{equation}
\label{eqn:asympTriv}
\mathbf{A} : C^\infty(S^1,\RR^{2n}) \to C^\infty(S^1,\RR^{2n}) :
\eta \mapsto - J_0 \p_t \eta - S(t) \eta
\end{equation}
in unitary trivializations,
where $S : S^1 \to \End(\RR^{2n})$ is a smooth loop of symmetric matrices.

Equivalently, an asymptotic operator on $(E,J,\omega)$ is any operator of
the form $-J \nabla$ where $\nabla$ is a symplectic connection on~$E$.
\end{defn}

\begin{exercise}
Show that any asymptotic operator on a Hermitian vector bundle $(E,J,\omega)$
over $S^1$ is symmetric with respect to the real $L^2$ bundle metric
$$
\langle \eta_1,\eta_2 \rangle_{L^2} := \int_{S^1} \omega(\eta_1(t),J \eta_2(t))\, dt.
$$
\end{exercise}

\begin{exercise}
Show that the asymptotic operator $\mathbf{A}_\gamma$ for a closed Reeb orbit
$\gamma$ is also an asymptotic operator on $(\gamma^*\xi,J,d\alpha)$ in
the sense of Definition~\ref{defn:asymptotic2}.
\end{exercise}

For functional analytic purposes, we shall regard asymptotic operators on
Hermitian bundles $(E,J,\omega)$ as bounded real-linear operators
$$
\mathbf{A} : H^1(E) \to L^2(E).
$$
By Lemma~\ref{lemma:itIsFredholm}, all asymptotic operators are then
Fredholm with index~$0$, and any two such operators on the same bundle
are compact perturbations of each other.  Regarding them alternatively as
unbounded symmetric operators on~$L^2(E)$, the spectral flow 
$$
\muspec(\mathbf{A}_-,\mathbf{A}_+) \in \ZZ
$$
between two such operators $\mathbf{A}_\pm$ with trivial kernel
is defined by choosing any unitary trivialization to write both in the
form $-J_0 \, \p_t - S(t)$, and it is independent of this choice.
The following is what we mean when we say that 
critical points of the action functional have ``infinite Morse index''
and ``infinite Morse co-index'':

\begin{prop}
Every asymptotic operator has infinitely many eigenvalues of both signs.
\end{prop}
\begin{proof}
It is easy to verify that this is true for 
$\mathbf{A}_0 := -J_0 \p_t : H^1(S^1,\RR^{2n}) \to L^2(S^1,\RR^{2n})$;
see the proof of theorem \ref{thm:winding} below.  It is therefore also
true for $\mathbf{A}_0 + \epsilon$ for any $\epsilon \in \RR$, and this
operator has trivial kernel whenever $\epsilon \not\in 2\pi\ZZ$.
For any other trivialized asymptotic operator $\mathbf{A}$ with 
$0 \not\in \sigma(\mathbf{A})$, the result then 
follows  from Theorem~\ref{thm:spectralFlow} since
$\muspec(\mathbf{A}_0 + \epsilon,\mathbf{A})$ is finite, 
and this is precisely the signed count of eigenvalues which change sign.
The condition $0 \not\in \sigma(\mathbf{A})$ can then be lifted by
replacing $\mathbf{A}$ with $\mathbf{A} + \epsilon$.
\end{proof}

\begin{exercise}
\label{EX:selfAdjoint}
Show that asymptotic operators are self-adjoint (as unbounded 
operators on $L^2$ with domain~$H^1$) in the sense of
Remark~\ref{remark:selfAdjoint}.
\end{exercise}

\section{The Conley-Zehnder index}
\label{sec:CZ}

We are now in a position to define a suitable replacement for the
Morse index in the context of SFT.  We shall say that an asymptotic
operator $\mathbf{A}$ is \defin{nondegenerate} whenever
$0 \not\in \sigma(\mathbf{A})$.
We will begin by defining the Conley-Zehnder index as an integer-valued
invariant of homotopy classes of nondegenerate asymptotic operators
on the trivial Hermitian bundle $S^1 \times \RR^{2n}$; the definition
on arbitrary Hermitian bundles will then depend on a choice of trivialization.

It is customary elsewhere in the literature 
(see e.g.~\cite{SalamonZehnder:Morse}) to adopt a somewhat different
perspective on the Conley-Zehnder index, in which it defines an integer-valued
invariant of connected components of the space of ``nondegenerate
symplectic arcs''
$$
\left\{ \Psi \in C^0([0,1],\Spp(2n))\ \big|\ 
\text{$\Psi(0) = \1$ and $1 \not\in \sigma(\Psi(1))$} \right\}.
$$
These are two different perspectives on the same notion.  A dictionary from
ours to the other perspective is provided by
associating to any trivialized nondegenerate 
asymptotic operator $\mathbf{A} = -J_0 \p_t - S(t)$ the symplectic arc $\Psi$
defined by the initial value problem
$$
(-J_0 \p_t - S(t)) \Psi(t) = 0, \qquad
\Psi(0) = \1.
$$
Conversely, any smooth symplectic arc determines via this same formula
a smooth path of symmetric matrices $S : [0,1] \to \End(\RR^{2n})$,
producing a mild generalization of our notion of an asymptotic 
operator.\footnote{If $S(t)$ is not continuous on $S^1$ but is continuous
on $[0,1]$, then $-J_0 \p_t - S(t)$ cannot be regarded as a linear operator
on $C^\infty(S^1,\RR^{2n})$ but is still a very well-behaved symmetric
Fredholm operator from $H^1(S^1)$ to $L^2(S^1)$.  All of the important
functional analytic results in this lecture can thus be generalized to
allow this.}

\begin{defn}
\label{defn:CZtriv}
The \defin{Conley-Zehnder index} associates to every
trivialized nondegenerate asymptotic operator 
$\mathbf{A} : H^1(S^1,\RR^{2n})
\to L^2(S^1,\RR^{2n})$ as in \eqref{eqn:asympTriv} an integer
$$
\muCZ(\mathbf{A}) \in \ZZ
$$
determined uniquely by the following properties:
\begin{enumerate}
\item Set $\muCZ(\mathbf{A}) := 0$ for the operator
$\mathbf{A} = -J_0 \p_t - \begin{pmatrix} \1 & 0 \\ 0 & -\1 \end{pmatrix}$.
\item For any two nondegenerate operators $\mathbf{A}_\pm$, set
$$
\muCZ(\mathbf{A}_-) - \muCZ(\mathbf{A}_+) := \muspec(\mathbf{A}_-,\mathbf{A}_+).
$$
\end{enumerate}
\end{defn}

\begin{defn}
\label{defn:CZ}
Given a nondegenerate asymptotic operator $\mathbf{A}$ on a Hermitian
bundle $(E,J,\omega)$ over $S^1$ and a choice of complex trivialization
$\tau$ for $(E,J)$, the \defin{Conley-Zehnder index} of $\mathbf{A}$
with respect to $\tau$ is the integer
$$
\muCZ^\tau(\mathbf{A}) \in \ZZ
$$
defined by choosing any unitary trivialization homotopic to $\tau$ to
write $\mathbf{A}$ as an operator $H^1(S^1,\RR^{2n}) \to L^2(S^1,\RR^{2n})$
and then plugging in Definition~\ref{defn:CZtriv}.

If $\gamma$ is a nondegenerate Reeb orbit $\gamma$ in a $(2n-1)$-dimensional 
contact manifold $(M,\xi=\ker\alpha)$, then for any 
complex trivialization
$\tau$ of $\gamma^*\xi \to S^1$, the \defin{Conley-Zehnder index} of $\gamma$
relative to $\tau$ is defined as
$$
\muCZ^\tau(\gamma) := \muCZ^\tau(\mathbf{A}_\gamma).
$$
\end{defn}

\begin{remark}
From the perspective of \cite{SalamonZehnder:Morse}, $\muCZ^\tau(\gamma)$ is
the Conley-Zehnder index of the linearized Reeb flow along $\gamma$ restricted
to~$\xi$, expressed via a choice of unitary trivialization 
as a nondegenerate arc in $\Spp(2n-2)$.
\end{remark}

\begin{exercise}
\label{EX:CZsum}
Show that if $\mathbf{A}_1$ and $\mathbf{A}_2$ are nondegenerate
asymptotic operators on
Hermitian bundles $E_1$ and $E_2$ respectively, then
$\mathbf{A}_1 \oplus \mathbf{A}_2$ defines a nondegenerate asymptotic
operator on $E_1 \oplus E_2$, and given trivializations
$\tau_j$ for $j=1,2$,
$$
\muCZ^{\tau_1 \oplus \tau_2}(\mathbf{A}_1 \oplus \mathbf{A}_2) =
\muCZ^{\tau_1}(\mathbf{A}_1) + \muCZ^{\tau_2}(\mathbf{A}_2).
$$
\end{exercise}

The following is a functional-analytic version of the well-known fact that
the Conley-Zehnder index classifies homotopy classes of nondegenerate
symplectic arcs.

\begin{thm}
\label{thm:CZclassification}
On any Hermitian bundle $(E,J,\omega) \to S^1$ with complex trivialization
$\tau$, two nondegenerate asymptotic operators $\mathbf{A_\pm}$ lie in the
same connected component of the space of nondegenerate asymptotic operators
if and only if $\muCZ^\tau(\mathbf{A}_+) = \muCZ^\tau(\mathbf{A}_-)$.
\end{thm}
\begin{proof}
Trivializing the bundle, we need to show that if $\mathbf{A}_\pm =
-J_0 \p_t - S_\pm(t)$ satisfy $\muspec(\mathbf{A}_-,\mathbf{A}_+) = 0$,
then there exists a path of asymptotic operators between them for which
no eigenvalues cross~$0$.  To see this, we can first choose any path
$\{\mathbf{A}_t\}_{t \in [-1,1]}$ of asymptotic operators with
$\mathbf{A}_{\pm 1} = \mathbf{A}_\pm$, and then use Lemma~\ref{lemma:generic}
to add generic compact perturbations producing a family 
$$
\left\{\mathbf{A}_t' \in \Fred_\RR^{\sym}(H^1,L^2,\mathbf{A}_+) \right\}_{t \in [-1,1]}
$$
whose intersections with $\Fred_\RR^{\sym,k}(H^1,L^2,\mathbf{A}_+)$ are transverse
for every $k \ge 1$, hence only simple eigenvalues cross~$0$ and they
cross transversely.  Any neighboring pair of crossings with opposite
signs can then be eliminated by changing $\{\mathbf{A}_t'\}_{t \in [-1,1]}$ to
$\{\mathbf{A}_t' + c(t)\}_{t \in [-1,1]}$ for a suitable choice of smooth
function $c : [-1,1] \to \RR$.  Since the spectral flow is zero, one can
repeat this modification until one obtains a path of perturbed operators 
with no crossings,
and it is a small perturbation of the path of asymptotic operators
$\{\mathbf{A}_t + c(t)\}_{t \in [-1,1]}$.  Since $\mathbf{A}_\pm$ are both
nondegenerate, one can assume moreover that
all eigenvalues of $\mathbf{A}_t + c(t)$ stay a fixed distance 
$\delta > 0$ away from~$0$, where $\delta$ is
independent of the perturbation.  One can therefore
``turn off the perturbation'' as in the proof of Theorem~\ref{thm:spectralFlow},
i.e.~there exists a sequence of perturbed paths $\{\mathbf{A}_t^\nu\}_{t \in [-1,1]}$
converging to $\{\mathbf{A}_t + c(t)\}$ whose eigenvalues stay a fixed 
distance away from~$0$, and the same is therefore true for the continuous
families of eigenvalues of $\mathbf{A}_t + c(t)$ obtained as $\nu \to \infty$.
\end{proof}

To compute Conley-Zehnder indices, Exercise~\ref{EX:CZsum} shows that it
suffices if we know how to compute them for operators on Hermitian
line bundles.  The next two theorems provide a tool for handling the latter.

\begin{thm}
\label{thm:winding}
Let $\mathbf{A} = -J_0 \p_t - S(t) : H^1(S^1,\RR^2) \to L^2(S^1,\RR^2)$,
where $S(t)$ is a smooth loop of symmetric $2$-by-$2$ matrices.
For each $\lambda \in \sigma(\mathbf{A})$, denote the corresponding
eigenspace by $E_\lambda \subset H^1(S^1,\RR^2)$.
\begin{enumerate}
\item \label{item:nowherezero}
Every nontrivial eigenfunction $e_\lambda \in E_\lambda$ is nowhere zero
and thus has a well-defined winding number $\wind(e_\lambda) \in \ZZ$.
\item \label{item:invariance}
Any two nontrivial eigenfunctions in the same eigenspace $E_\lambda$ have
the same winding number.
\item \label{item:monotonicity}
If $\lambda, \mu \in \sigma(\mathbf{A})$ satisfy $\lambda < \mu$, then any
two nontrivial eigenfunctions $e_\lambda \in E_\lambda$ and
$e_\mu \in E_\mu$ satisfy $\wind(e_\lambda) \le \wind(e_\mu)$.
\item \label{item:multiplicity}
For every $k \in \ZZ$, $\mathbf{A}$ has exactly two eigenvalues (counting
multiplicity) for which the corresponding eigenfunctions have
winding number equal to~$k$.
\end{enumerate}
\end{thm}
\begin{proof}
We follow the proof given in \cite{HWZ:props2}.
 
Observe first that (\ref{item:nowherezero}) follows from the fact that nontrivial eigenfunctions are solutions to an ODE, for which classical existence and uniqueness results are available. Since the trivial map is a solution, every eigenfunction which vanishes at a point must be itself trivial, by uniqueness.

To prove (\ref{item:invariance}), let $\nu_0$ and $\nu_1$ be 
nontrivial eigenfunctions for the same eigenvalue $\lambda$.  If their
winding numbers are different, then there exists $t_0 \in S^1$ at which
$\nu_1(t_0)$ is a nonzero real multiple of~$\nu_0(t_0)$, so after rescaling,
we can assume $\nu_0(t_0) = \nu_1(t_0)$.  But $\nu_0$ and $\nu_1$ are both
solutions to the same linear ODE, so this implies $\nu_0(t) = \nu_1(t)$ for
all~$t$ and thus contradicts the assumption on the winding numbers.
 
We first prove the rest for the case $S=0$ and the operator $\mathbf{A}_0 = -J_0 \p_t$. Given $\nu \in H^1(S^1,\mathbb{R}^2)$, written as $\nu(t)=(x(t),y(t))$, we have that $\nu$ is an element of $E_\lambda$ for the operator $\mathbf{A}_0$ if and only $(\dot{y},-\dot{x})=\lambda(x,y)$. This has solutions of the form 
\[ 
\left \{
  \begin{array}{c}
  x(t)=A\cos(\lambda t)-B\sin(\lambda t)\\
  y(t)=B\cos(\lambda t)+A\sin(\lambda t)\\
  \end{array}   
\right., 
\] for some constants $A,B \in \RR$, which are defined on $S^1$ as long as $\lambda \in 2\pi\mathbb{Z}$.
In other words, the spectrum of this operator is $\sigma(\mathbf{A}_0)=2\pi \mathbb{Z}$. 
Hence $\nu(t)=\nu(0)e^{i\lambda t}$, which has winding number $$\wind(\nu)=\frac{\lambda}{2\pi}$$
Statements (\ref{item:invariance}) and (\ref{item:monotonicity}) are now
obvious, and (\ref{item:multiplicity}) follows from the observation that $E_\lambda$ is two-dimensional, so in this case each eigenvalue is to be counted with multiplicity two.

For the general case, consider the path of asymptotic operators given by
$$
\left\{ \mathbf{A}_\tau=-J_0\p_t - \tau S(t) \right\}_{\tau \in [0,1]}.
$$ 
Theorem~\ref{thm:spectralFlow} gives continuous families $\{\lambda_j:[0,1]\rightarrow \RR\}_{j \in \mathbb{Z}}$ and $\{\nu_j:[0,1] \rightarrow H^1(S^1,\mathbb{R}^2)\}_{j \in \mathbb{Z}}$ such that for every $\tau \in [0,1]$, $\nu_j(\tau)$ is an eigenfunction for the operator $\mathbf{A}_\tau$ with eigenvalue $\lambda_j(\tau)$, whose multiplicity is given by the number of $i \in \mathbb{Z}$ for which $\lambda_i(\tau)=\lambda_j(\tau)$, and such that $\lambda_{2n+k}(0)=2\pi n$, for $k=0,1$ (this eigenvalue has multiplicity 2). Now, since the winding number is a homotopy invariant (hence invariant under deformations), we have 
$$
\wind(\nu_{2n+k}(\tau))= \wind(\nu_{2n+k}(0))=n,
$$ 
for $k=0,1$. Moreover, since the winding 
only depends on the eigenvalue, the only paths that can possibly meet are $\lambda_{2n}$ and $\lambda_{2n+1}$, which implies that the multiplicity of every eigenvalue $\lambda_i(\tau)$ is at most two, with equality where these two ``branches''
meet. Hence (\ref{item:monotonicity}) and (\ref{item:multiplicity}) follow, where equality in (\ref{item:monotonicity}) holds if and only if the two branches of paths of eigenvalues with same winding number end up at different points. 
\end{proof}

The theorem implies the existence of a well-defined and nondecreasing function
$$
\sigma(\mathbf{A}) \to \ZZ : \lambda \mapsto \wind(\lambda),
$$
where $\wind(\lambda)$ is defined as $\wind(e_\lambda)$ for any
nontrivial $e_\lambda \in E_\lambda$, and this function
attains every value exactly twice (counting multiplicity of eigenvalues).
Since eigenvalues of $\mathbf{A}$ are isolated, we can therefore associate 
to any nondegenerate asymptotic operator $\mathbf{A}$
on the trivial Hermitian line bundle its \defin{extremal winding numbers} 
and its \defin{parity},
\begin{equation}
\label{eqn:alphas}
\begin{split}
\alpha_+(\mathbf{A}) &= \min_{\lambda \in \sigma(\mathbf{A}) \cap (0,\infty)} 
\wind(\lambda) \in \ZZ, \\
\alpha_-(\mathbf{A}) &= \max_{\lambda \in \sigma(\mathbf{A}) \cap (-\infty,0)}
\wind(\lambda) \in \ZZ, \\
p(\mathbf{A}) &= \alpha_+(\mathbf{A}) - \alpha_-(\mathbf{A}) \in \{0,1\}.
\end{split}
\end{equation}

\begin{thm}
\label{thm:CZwinding}
If $\mathbf{A}$ is a nondegenerate asymptotic operator on the trivial
Hermitian line bundle $S^1 \times \RR^2 \to S^1$, then
$$
\muCZ(\mathbf{A}) = 2\alpha_-(\mathbf{A}) + p(\mathbf{A}) =
2\alpha_+(\mathbf{A}) - p(\mathbf{A}).
$$
\end{thm}
\begin{proof}
The operator 
$\mathbf{A}_0 = -J_0 \p_t - \begin{pmatrix} 1 & 0 \\ 0 & -1 \end{pmatrix}$
satisfies $\muCZ(\mathbf{A}_0) = 0$ by definition, and it has two constant
eigenfunctions with eigenvalues of opposite signs, hence
$$
\alpha_-(\mathbf{A}_0) = \alpha_+(\mathbf{A}_0) = 0,
$$
consistent with the stated formula.  The general case then follows by
computing the spectral flow from $\mathbf{A}_0$ to any other nondegenerate
operator $\mathbf{A}$, and observing that the winding number associated
to any continuous family of eigenvalues (as in
Theorem~\ref{thm:spectralFlow}) for a path 
$\{\mathbf{A}_t\}_{t \in [-1,1]}$ of asymptotic operators cannot change.
\end{proof}

For any Hermitian line bundle $(E,J,\omega)$ over $S^1$ with a nondegenerate
asymptotic operator $\mathbf{A}$, we can similarly
choose a complex trivialization $\tau$ to define the winding numbers
$\alpha_\pm^\tau(\mathbf{A}) \in \ZZ$ and parity
$p(\mathbf{A}) = \alpha^\tau_+(\mathbf{A}) - \alpha^\tau_-(\mathbf{A})
\in \{0,1\}$;
note that the dependence on $\tau$ cancels out in the last formula,
so that $p(\mathbf{A})$ is independent of choices.  We then can associate
to any nondegenerate Reeb orbit $\gamma$ in a contact $3$-manifold 
$(M,\xi=\ker\alpha)$ with a trivialization $\tau$
of $\gamma^*\xi$ the integers $\alpha^\tau_\pm(\gamma)$ and $p(\gamma)$,
such that
$$
\muCZ^\tau(\gamma) = 2\alpha_-^\tau(\gamma) + p(\gamma) =
2\alpha_+^\tau(\gamma) - p(\gamma)
$$
holds.

\begin{exercise}
\label{EX:changeTriv}
Given a Hermitian vector bundle $(E,J,\omega) \to S^1$ with two complex
trivializations $\tau_j : E \to S^1 \times \RR^{2n}$ for $j=1,2$, denote by
$$
\deg(\tau_1 \circ \tau_2^{-1}) \in \ZZ
$$
the winding number of $\det g : S^1 \to \CC \setminus \{0\}$, where
$g : S^1 \to \GL(n,\CC)$ is the transition map appearing in the formula
$\tau_1 \circ \tau_2^{-1}(t,v) = (t,g(t)v)$.  Show that for any
asymptotic operator $\mathbf{A}$ on $(E,J,\omega)$,
$$
\muCZ^{\tau_2}(\mathbf{A}) = \muCZ^{\tau_1}(\mathbf{A}) + 2 \deg(\tau_2
\circ \tau_1^{-1}).
$$
\end{exercise}

Exercise~\ref{EX:changeTriv} provides the useful formula
$$
\muCZ^{\tau_2}(\gamma) = \muCZ^{\tau_1}(\gamma) + 
2 \deg(\tau_2 \circ \tau_1^{-1})
$$
for any two trivializations $\tau_1,\tau_2$ of $\xi$ along a 
nondegenerate Reeb orbit~$\gamma$.  
In particular, this shows that the \defin{parity}
$$
\muCZ^{\ZZ_2}(\gamma) := [\muCZ^\tau(\gamma)] \in \ZZ_2
$$
of the orbit does not depend on a choice of trivialization.  We sometimes
refer to \defin{even orbits} and \defin{odd orbits} accordingly.

\begin{exercise}
Show that if a Reeb orbit $\gamma : S^1 \to M$ in a contact $3$-manifold
$(M,\xi=\ker\alpha)$ is nondegenerate and has even parity,
then the same is true for all of its multiple covers
$$
\gamma^k : S^1 \to M : t \mapsto \gamma(kt), \qquad k \in \NN.
$$
\end{exercise}

\chapter{Fredholm theory with cylindrical ends}
\label{lec:Fredholm}

\minitoc

\vspace{12pt}

In this lecture we will study the class of linear Cauchy-Riemann type
operators that arise by linearizing the nonlinear equation for moduli
spaces in SFT.  We saw in the previous lecture that linearizing PDEs
over domains with cylindrical ends naturally leads one to consider
certain symmetric \emph{asymptotic operators} (e.g.~the Hessian of a
Morse function at its critical points), which have trivial kernel
if and only if a nondegeneracy (i.e.~Morse) condition is satisfied.
Our goal in this lecture is to write down the SFT version of this story
and show that the linear Cauchy-Riemann type operators are Fredholm
if their asymptotic operators are nondegenerate.

\section{Cauchy-Riemann operators with punctures}
\label{sec:setup}

The setup throughout this lecture will be as follows.

Assume $(\Sigma,j)$ is a closed connected Riemann surface of genus $g \ge 0$,
$\Gamma \subset \Sigma$ is a finite set partitioned into two subsets
$$
\Gamma = \Gamma^+ \cup \Gamma^-,
$$
and $\dot{\Sigma} := \Sigma \setminus \Gamma$ denotes the resulting punctured
Riemann surface.  We shall fix a choice of \defin{holomorphic cylindrical
coordinate} near each puncture $z \in \Gamma^\pm$, meaning the following.
Given $R \ge 0$, let $(Z_\pm^R,i)$ denote the half-cylinders
$$
Z_+^R := [R,\infty) \times S^1, \qquad Z_-^R := (-\infty,-R] \times S^1,
\qquad Z_\pm := Z_\pm^0,
$$
with complex structure $i\p_s = \p_t$, $i\p_t = -\p_s$ in coordinates
$(s,t) \in \RR \times S^1$.  The standard half-cylinders $Z_\pm$ are each
biholomorphically equivalent to
the punctured disk $\dot{\DD} := \DD \setminus \{0\}$ via the maps
$$
\psi_\pm : Z_\pm \to \dot{\DD} : (s,t) \mapsto e^{\mp 2\pi(s+it)}.
$$
For $z \in \Gamma^\pm$, we choose a closed neighborhood 
$\uU_z \subset \Sigma$ of $z$ with a biholomorphic map
$$
\varphi_z : (\dot{\uU}_z,j) \to (Z_\pm,i),
$$
where $\dot{\uU}_z := \uU_z \setminus \{z\}$,
such that $\psi_\pm \circ \varphi_z : \dot{\uU}_z \to \dot{\DD}$
extends holomorphically to $\uU_z \to \DD$ with $z \mapsto 0$.  One can
always find such coordinates by choosing holomorphic coordinates 
near~$z$.  We can thus view the punctured neighborhoods
$\dot{\uU}_z \subset \dot{\Sigma}$ as \defin{cylindrical ends}~$Z_\pm$.

Suppose $(E,J)$ is a smooth complex vector bundle of rank $m$
over $(\dot{\Sigma},j)$.  An \defin{asymptotically Hermitian structure}
on $(E,J)$ is a choice of Hermitian vector bundles $(E_z,J_z,\omega_z)$ of
rank~$m$ associated to each puncture $z \in \Gamma^\pm$, together with 
choices of complex bundle isomorphisms
$$
E|_{\dot{\uU}_z} \to \pr_2^*E_z
$$
covering $\varphi_z : \dot{\uU}_z \to Z_\pm$, where
$\pr_2 : Z_\pm \to S^1$ denotes the natural projection to the $S^1$ factor.
This isomorphism induces from any unitary trivialization $\tau$ of
$(E_z,J_z,\omega_z)$ a complex trivialization 
\begin{equation}
\label{eqn:asympTriv2}
\tau : E|_{\dot{\uU}_z} \to Z_\pm \times \RR^{2m}
\end{equation}
over the cylindrical end, which we will call an \defin{asymptotic
trivialization} near~$z$.  The bundle $(E_z,J_z,\omega_z)$ will be referred
to as the \defin{asymptotic bundle} associated to $(E,J)$ near~$z$.

Fixing asymptotic trivializations near every puncture,
we can now define Sobolev spaces of sections of $E$ by
$$
W^{k,p}(E) := \left\{ \eta \in W^{k,p}_\loc(E)\ \Big|\ 
\text{$\eta_z \in W^{k,p}(\intZ_\pm,\RR^{2m})$ for every $z \in \Gamma^\pm$}
\right\},
$$
where $\eta_z : Z_\pm \to \RR^{2m}$ denotes the expression of 
$\eta|_{\dot{\uU}_z}$ in terms of the asymptotic trivialization, and
we use the standard area form $ds \wedge dt$ on $Z_\pm$ to define the
norm.  Since $S^1$ is compact, different choices of asymptotic trivialization
give rise to equivalent norms, however:

\begin{exercise}
Convince yourself that different choices of asymptotically Hermitian
structure on $E \to \dot{\Sigma}$ can give rise to \emph{inequivalent}
$W^{k,p}$-norms.
\end{exercise}

Any linear Cauchy-Riemann type operator on $E$ has as its target the complex
vector bundle
$$
F := \overline{\Hom}_\CC(T\dot{\Sigma},E),
$$
so sections of $F$ are the same thing as $E$-valued $(0,1)$-forms.
An asymptotic trivialization $\tau$ as in \eqref{eqn:asympTriv2} then
also induces a complex trivialization 
$$
F|_{\dot{\uU}_z} \to Z_\pm \times \RR^{2m}
: \lambda \mapsto \tau(\lambda(\p_s)),
$$
where $\p_s$ is the vector field on $\dot{\uU}_z$ arising from its
identification with~$Z_\pm$.  This trivialization yields a corresponding
definition for the Sobolev spaces 
$W^{k,p}(F)$, which depend on the
asymptotically Hermitian structure of $E$ but not on the choices of
asymptotic trivializations.
Having made these choices, a
Cauchy-Riemann type operator $\mathbf{D} : \Gamma(E) \to \Gamma(F)$ always
appears over $\dot{\uU}_z$
as a linear map on $C^\infty(Z_\pm,\RR^{2m})$ of the form
\begin{equation}
\label{eqn:CRasymp}
\mathbf{D} \eta(s,t) = \dbar \eta(s,t) + S(s,t) \eta(s,t),
\end{equation}
where $\dbar := \p_s + J_0 \p_t$ and $S \in C^\infty(Z_\pm , \End(\RR^{2m}))$.

\begin{defn}
\label{defn:asympCR}
Suppose $\mathbf{A}_z$ is an asymptotic operator on $(E_z,J_z,\omega_z)$
and $\mathbf{D}$ is a linear Cauchy-Riemann
type operator on $(E,J)$.  We say that $\mathbf{D}$ is \defin{asymptotic to}
$\mathbf{A}_z$ at $z$ if $\mathbf{D}$ appears in the form
\eqref{eqn:CRasymp} with respect to an asymptotic trivialization near~$z$, with
$$
\| S - S_\infty \|_{C^k(Z_\pm^R)} \to 0 \quad \text{ as } \quad R \to \infty
$$
for all $k \in \NN$, where $S_\infty(s,t) := S_\infty(t)$ is a smooth
loop of symmetric matrices such that $\mathbf{A}_z$ appears in the
corresponding unitary trivialization of $(E_z,J_z,\omega_z)$ as
$-J_0 \p_t - S_\infty$.
\end{defn}

Recall that an asymptotic operator is called \defin{nondegenerate} if 
$0$ is not in its spectrum, which means it defines an isomorphism $H^1 \to L^2$.
The objective of this lecture will be to prove the following:

\begin{thm}
\label{thm:Fredholm}
Suppose $(E,J)$ is an asymptotically Hermitian vector bundle
over $(\dot{\Sigma},j)$, $\mathbf{A}_z$ is a nondegenerate
asymptotic operator on the associated asymptotic bundle $(E_z,J_z,\omega_z)$
for each $z \in \Gamma$, and $\mathbf{D}$ is a linear Cauchy-Riemann type
operator asymptotic to $\mathbf{A}_z$ at each puncture~$z$.
Then for every $k \in \NN$ and $1 < p < \infty$,
$$
\mathbf{D} : W^{k,p}(E) \to W^{k-1,p}(F)
$$
is Fredholm.  Moreover, $\ind \mathbf{D}$ and $\ker \mathbf{D}$ are each
independent of $k$ and $p$, the latter being
a space of smooth sections whose derivatives of all orders
decay to $0$ at infinity.
\end{thm}

\begin{remark}
The asymptotic decay conditions on $S(s,t)$ in Definition~\ref{defn:asympCR}
can be relaxed at the cost of limiting the range of $k \in \NN$ for
which Theorem~\ref{thm:Fredholm} is valid.  To prove that
$\mathbf{D} : W^{1,p} \to L^p$ is Fredholm, it suffices to assume
$S(s,\cdot) \to S_\infty$ uniformly as $|s| \to \infty$.
\end{remark}

The index of $\mathbf{D}$ is determined by a generalization of the Riemann-Roch
formula involving the Conley-Zehnder indices $\muCZ^\tau(\mathbf{A}_z)$
that were introduced in the previous lecture.
We will postpone serious discussion of the index formula until the next lecture, 
but here is the statement:

\begin{thm}
\label{thm:RiemannRochPreview}
In the setting of Theorem~\ref{thm:Fredholm},
$$
\ind \mathbf{D} = m \chi(\dot{\Sigma}) + 2 c_1^\tau(E) +
\sum_{z \in \Gamma^+} \muCZ^\tau(\mathbf{A}_z) -
\sum_{z \in \Gamma^-} \muCZ^\tau(\mathbf{A}_z),
$$
where $\tau$ is an arbitrary choice of asymptotic trivializations, 
$c_1^\tau(E) \in \ZZ$ is the \emph{relative first Chern number} of $E$ with
respect to $\tau$, and the sum is independent of this choice.
\end{thm}

For the rest of this lecture, we maintain as standing assumptions
that $k \in \NN$, $1 < p < \infty$, and
$\mathbf{D}$ is a linear Cauchy-Riemann type operator on $E$
asymptotic at the punctures to a fixed set of asymptotic operators 
$\{\mathbf{A}_z\}_{z \in \Gamma}$.  We will not always need to assume
that the $\mathbf{A}_z$ are nondegenerate, so this condition
will be specified whenever it is relevant.
For subdomains $\Sigma_0 \subset \dot{\Sigma}$, we will sometimes denote
the $W^{k,p}$-norm on sections of $E$ restricted to $\Sigma_0$ by
$$
\| \eta \|_{W^{k,p}(\Sigma_0)} := \| \eta \|_{W^{k,p}(E|_{\Sigma_0})},
$$
and we will use the same notation for sections of other bundles 
such as $F = \overline{\Hom}_\CC(T\dot{\Sigma},E)$ over this
domain when there is no danger of confusion.  The space
$$
W^{k,p}_0(\Sigma_0) \subset W^{k,p}(E)
$$
is defined in this case as the $W^{k,p}$-closure of the space of smooth
sections of $E$ with compact support in~$\Sigma_0 \setminus \p\Sigma_0$.  
For some background discussion on Sobolev spaces of sections of vector bundles, 
see Appendix~\ref{app:Sobolev}.

\section{A global weak regularity result}
\label{sec:regularityGlobal}

In Lecture~\ref{lec:local} we proved that for $1 < p < \infty$, weak solutions of class
$L^p_\loc$ to linear Cauchy-Riemann type equations are always smooth.
Here is a global version of that result.

\begin{prop}
\label{prop:globalReg}
Suppose $1 < p < \infty$ and $k \in \NN$.
If $\eta \in L^p(E)$ weakly satisfies
$$
\mathbf{D} \eta \in W^{k-1,p}(F),
$$
then $\eta \in W^{k,p}(E)$.
\end{prop}
\begin{proof}
By induction,
it suffices to show that if $\eta \in W^{k-1,p}$ and
$\mathbf{D}\eta \in W^{k-1,p}$ then $\eta \in W^{k,p}$.
We already know that this is true locally, so the task is to
bound the $W^{k,p}$-norm of $\eta$ on the cylindrical ends.  
Pick an asymptotic trivialization and write $\mathbf{D}$ on one of the
ends $Z_\pm \cong \dot{\uU}_z$ as $\dbar + S(s,t)$.
Let us assume for concreteness that the puncture is a positive one, and now
consider the $W^{k,p}$-norm of $\eta$ on 
$(N,N+1) \times S^1 \subset \dot{\uU}_z$ for $N \in \NN$.
Choosing a smooth bump function $\beta : \RR \times S^1 \to [0,1]$
supported in $(N-1,N+2) \times S^1$ with $\beta = 1$
on $[N,N+1] \times S^1$, we can use the usual elliptic estimate to write
\begin{equation*}
\begin{split}
\| \eta \|_{W^{k,p}((N,N+1) \times S^1)} &\le
\| \beta \eta \|_{W^{k,p}((N-1,N+2) \times S^1)} \le
c \| \dbar(\beta \eta) \|_{W^{k-1,p}((N-1,N+2) \times S^1)} \\
&\le c \| \eta \|_{W^{k-1,p}((N-1,N+2) \times S^1)} +
c \| \dbar \eta \|_{W^{k-1,p}((N-1,N+2) \times S^1)} \\
&= c \| \eta \|_{W^{k-1,p}((N-1,N+2) \times S^1)} +
c \| \mathbf{D}\eta - S\eta \|_{W^{k-1,p}((N-1,N+2) \times S^1)} \\
&\le c' \| \eta \|_{W^{k-1,p}((N-1,N+2) \times S^1)} +
c' \| \mathbf{D}\eta \|_{W^{k-1,p}((N-1,N+2) \times S^1)}.
\end{split}
\end{equation*}
An important detail here is that the constants in these estimates can be
assumed independent of~$N$: indeed, one can use shifts of the same cutoff
function for any $N$, and the $C^{k-1}$-norm of $S$ on $[N-1,N+2] \times S^1$
is also bounded uniformly in $N$ since $S(s,t)$ converges asymptotically
to some~$S_\infty(t)$.  We can therefore take the sum of this estimate for
all $N \in \NN$, producing
$$
\| \eta \|_{W^{k,p}(\intZ_+^1)} \le c \| \eta \|_{W^{k-1,p}(\intZ_+)} +
c \| \mathbf{D}\eta \|_{W^{k-1,p}(\intZ_+)}.
$$
\end{proof}

\begin{cor}
For $1 < p < \infty$, any weak solution $\eta \in L^p(E)$ of 
$\mathbf{D}\eta = 0$ is smooth,
with derivatives of all orders decaying to $0$ at infinity.
\end{cor}
\begin{proof}
Proposition~\ref{prop:globalReg} implies $\eta \in W^{k,p}(E)$ for every
$k \in \NN$, so smoothness follows from the Sobolev embedding theorem.
Moreover, suppose $k$ and $p$ are large enough to have a continuous
inclusion $W^{k,p} \hookrightarrow C^m$ for some $m \in \NN$.
Then the finiteness of the $W^{k,p}$-norm also implies that for each
end $\dot{\uU}_z = Z_\pm$,
$$
\| \eta \|_{C^m(Z_\pm^R)} \le c \| \eta \|_{W^{k,p}(\intZ_\pm^R)} \to 0
\quad \text{ as } \quad R \to \infty.
$$
\end{proof}

\section{Elliptic estimates on cylindrical ends}
\label{sec:estimatesEnds}

The local elliptic estimates for $\dbar = \p_s + J_0 \p_t$ in Lecture~\ref{lec:local} 
applied to functions on $\intDD \subset \CC$ with compact support.
Using a finite open covering with a subordinate partition of unity,
it is a straightforward matter to turn these local estimates into
the following global result (cf.~\cite{Wendl:lecturesV33}*{Lemma~3.3.2}):

\begin{prop}
\label{prop:globalCpct}
If $\Sigma_0 \subset \dot{\Sigma}$ is a compact $2$-dimensional submanifold 
with boundary, then there exists a constant $c > 0$ such that
$$
\| \eta \|_{W^{k,p}(\Sigma_0)} \le c \| \mathbf{D}\eta \|_{W^{k-1,p}(\Sigma_0)}
+ c \| \eta \|_{W^{k-1,p}(\Sigma_0)}
$$
for all $\eta \in W^{k,p}_0(\Sigma_0)$.
\qed
\end{prop}

This unfortunately is unsufficient for the global problem under consideration,
since one has to chop off the cylindrical ends of $\dot{\Sigma}$ in order to
obtain a compact domain.  We therefore supplement the previous local estimates
with an asymptotic estimate.

\begin{prop}
\label{prop:ellipticEnd}
Suppose $z \in \Gamma^\pm$ is a puncture such that the asymptotic operator
$\mathbf{A}_z$ is nondegenerate.  Then on $Z_\pm^R \subset \dot{\uU}_z$
for sufficiently large $R \ge 0$, there exists a constant $c > 0$ such that
$$
\| \eta \|_{W^{k,p}(\intZ_\pm^R)} \le c \| \mathbf{D} \eta \|_{W^{k-1,p}(\intZ_\pm^R)}
\qquad \text{ for all } \qquad 
\eta \in W^{k,p}_0(\intZ_\pm^R).
$$
\end{prop}
\begin{remark}
Recall that $W^{k,p}_0(\intZ_\pm^R)$ denotes the $W^{k,p}$-closure of
$C_0^\infty(\intZ_\pm^R)$, so such functions remain in $W^{k,p}$ if they
are extended as zero to larger domains containing $\intZ_\pm^R$.  Note that
functions of class $W^{k,p}_0$ on $\intZ_\pm^R$ need not actually have
compact support; in fact $C_0^\infty$ is dense in $W^{k,p}(\RR \times S^1)$,
see~\S\ref{sec:SobolevEnds}.
\end{remark}

The proof of this requires a basic result about 
translation-invariant Cauchy-Riemann type operators on the cylinder.
Other than the elliptic estimates we discussed in Lecture~\ref{lec:local}, this is the
main analytical ingredient that makes all Floer-type theories in symplectic geometry
work.

\begin{thm}
\label{thm:invertible}
Suppose $k \in \NN$, $1 < p < \infty$, and $\mathbf{A} = -J_0 \p_t - S(t)$
is a nondegenerate asymptotic operator on the trivial Hermitian vector
bundle $S^1 \times \RR^{2n} \to S^1$.  Then the operator
$$
\p_s - \mathbf{A} = \p_s + J_0 \p_t + S(t) : W^{k,p}(\RR\times S^1,\RR^{2n})
\to W^{k-1,p}(\RR \times S^1,\RR^{2n})
$$
is an isomorphism.
\qed
\end{thm}

A detailed proof of this result for $k=1$ can be found in 
\cite{Salamon:Floer}*{Lemma~2.4}, and the general result follows easily
from this using regularity (Proposition~\ref{prop:globalReg}).
I will not attempt to reproduce the proof in Salamon's notes here since it is 
somewhat involved, but let us informally sketch the first step, which is the 
interesting part.  The goal is to prove that $\mathbf{D}_0 := \p_s
- \mathbf{A}$ is an invertible operator from $H^1(\RR \times S^1)$ 
to $L^2(\RR \times S^1)$.
To gain some intuition on this, consider the special case where the asymptotic 
operator is of the form $\mathbf{A} = - i \p_t - C$ for some constant 
$C \in \RR$.
One can then write down an inverse of $\mathbf{D}_0$
explicitly by combining a Fourier transform in the $s$ variable with a
Fourier series in the $t$ variable.  That is, sufficiently nice functions
$u$ on $\RR \times S^1$ can be expressed as
$$
u(s,t) = \sum_{k \in \ZZ} \int_\RR \hat{u}_k(\sigma) e^{2\pi i \sigma s}
e^{2\pi i kt}\, d\sigma,
$$
where the hybrid Fourier transform/series $\hat{u}$ depends on a continuous
variable $\sigma \in \RR$ and a discrete variable~$k \in \ZZ$.  One can then
obtain $\hat{u}$ from $u$ by
$$
\hat{u}_k(\sigma) = \int_{\RR \times S^1} u(s,t) e^{-2\pi i \sigma s}
e^{-2\pi i kt}\, ds \, dt,
$$
and we have the usual derivative formulas $\widehat{\p_s u}_k(\sigma) = 
2\pi i \sigma \hat{u}_k(\sigma)$ and
$\widehat{\p_t u}_k(\sigma) = 2\pi i k \hat{u}_k(\sigma)$.
The relation $(\p_s + i \p_t + C)u = f$ therefore produces an inversion 
formula of the form
$$
\hat{u}_k(\sigma) = \frac{\hat{f}_k(\sigma)}{2\pi i \sigma - 2\pi k + C}.
$$
This is a nice formula and produces from any $f \in L^2$ an element
$u \in H^1$ unless $C \in 2\pi \ZZ$, in which case the denominator has a
singularity.  This condition means $C$ must not be an eigenvalue of
$-i \p_t$, or in other words, $\mathbf{A} = -i \p_t - C$ is nondegenerate.
One can perhaps imagine carrying out a similar argument in the general
case using an orthonormal set of eigenfunctions\footnote{Recall from
Lecture~\ref{lec:asymptotic} that the spectrum $\sigma(\mathbf{A})$ of an arbitrary
asymptotic operator $\mathbf{A}$ always consists only of isolated real
eigenvalues, thus one can find $\lambda \in \RR$ for which 
$\lambda - \mathbf{A} : H^1(S^1) \to L^2(S^1)$ is invertible.  Its inverse, also
known as the \defin{resolvent}, then defines a \emph{compact} self-adjoint
operator $(\lambda - \mathbf{A})^{-1} : L^2(S^1) \to L^2(S^1)$ due to the compact
inclusion $H^1(S^1) \hookrightarrow L^2(S^1)$.  The spectral theorem for compact
self-adjoint operators now provides an orthonormal basis of $L^2(S^1)$
consisting of eigenfunctions of $(\lambda - \mathbf{A})^{-1}$, which are
also eigenfunctions of~$\mathbf{A}$.}
for $\mathbf{A}$ in place
of the functions $e^{2\pi i k t}$; this is presumably part of the idea
behind the actual proof in \cite{Salamon:Floer}, which uses strongly continuous
semigroups generated by the self-adjoint operator~$\mathbf{A}$.

\begin{proof}[Proof of Proposition~\ref{prop:ellipticEnd}]
Write $\mathbf{D} = \p_s + J_0 \p_t + S(s,t)$ 
and $\mathbf{D}_0 = \p_s + J_0 \p_t + S_\infty(t)$
in an asymptotic trivialization 
on $\dot{\uU}_z = Z_\pm$, where the nondegenerate asymptotic operator
is $\mathbf{A} = -J_0 \p_t - S_\infty(t)$ and we assume
$$
\| S - S_\infty \|_{C^{k-1}(Z^R_\pm)} \to 0 \quad \text{ as } \quad
R \to \infty.
$$
For $\eta \in W^{k,p}_0(\intZ_\pm^R)$, there is a canonical extension
$\eta \in W^{k,p}(\RR \times S^1)$ that equals zero outside
$Z_\pm^R$, so by Theorem~\ref{thm:invertible} we have
$$
\| \eta \|_{W^{k,p}(\intZ^R_\pm)} = \| \eta \|_{W^{k,p}(\RR \times S^1)}
\le c \| \mathbf{D}_0 \eta \|_{W^{k-1,p}(\RR \times S^1)} =
c \| \mathbf{D}_0 \eta \|_{W^{k-1,p}(\RR \times S^1)}.
$$
Rewriting this in terms of $\mathbf{D}$ gives
$$
\| \eta \|_{W^{k,p}(\intZ^R_\pm)} \le c \| \mathbf{D}\eta \|_{W^{k-1,p}(Z^R_\pm)}
+ c \| (S_\infty - S) \eta \|_{W^{k-1,p}(Z^R_\pm)},
$$
where the constants $c > 0$ do not depend on~$R$.  For this reason, we are
free to make $R \ge 0$ large enough to make the $C^{k-1}$-norm of
$S_\infty - S$ on $Z^R_\pm$ less than an arbitrarily small number
$\delta > 0$, in which case the above gives
$$
\| \eta \|_{W^{k,p}(\intZ^R_\pm)} \le c \| \mathbf{D}\eta \|_{W^{k-1,p}(Z^R_\pm)}
+ c \delta \| \eta \|_{W^{k-1,p}(Z^R_\pm)},
$$
and thus by the inclusion $W^{k-1,p} \hookrightarrow W^{k,p}$,
$$
\| \eta \|_{W^{k,p}(\intZ^R_\pm)} \le \frac{c}{1 - c\delta}
\| \mathbf{D}\eta \|_{W^{k-1,p}(Z^R_\pm)}.
$$
\end{proof}

\section{The semi-Fredholm property}

The standard approach for proving that elliptic operators are Fredholm
begins by proving that they are \defin{semi-Fredholm}, meaning
$\dim \ker \mathbf{D} < \infty$ and $\im \mathbf{D}$ is closed.  
In most settings, it is not hard to show that local elliptic estimates
give rise to global estimates of the form $\| \eta \|_{W^{k,p}} \le
c \| \mathbf{D}\eta \|_{W^{k-1,p}} + \| \eta \|_{W^{k-1,p}}$.
The step from these estimates to the semi-Fredholm property is then provided
by the following lemma.

\begin{lemma}
\label{lemma:semiFred}
Suppose $X$, $Y$ and $Z$ are Banach spaces, $\mathbf{T} \in \Lin(X,Y)$,
$\mathbf{K} \in \Lin(X,Z)$ is compact, and there is a constant $c > 0$ such that for all
$x \in X$,
\begin{equation}
\label{eqn:generalEstimate}
\| x \|_X \le c \| \mathbf{T}x \|_Y + c \| \mathbf{K}x \|_Z.
\end{equation}
Then $\ker \mathbf{T}$ is finite dimensional and $\im \mathbf{T}$ is closed.
\end{lemma}
\begin{proof}
A vector space is finite dimensional if and only if the unit ball in that
space is a compact set, so we begin by proving the latter holds for $\ker \mathbf{T}$.
Suppose $x_k \in \ker \mathbf{T}$ is a bounded sequence.  Then since $\mathbf{K}$ is a
compact operator, $\mathbf{K} x_k$ has a convergent subsequence in $Z$, which is
therefore Cauchy.  But \eqref{eqn:generalEstimate} 
then implies that the corresponding
subsequence of $x_k$ in $X$ is also Cauchy, and thus converges.

Since we now know $\ker \mathbf{T}$ is finite dimensional, we also know there is a
closed complement $V \subset X$ with $\ker \mathbf{T} \oplus V = X$.  
Then the restriction
$\mathbf{T}|_V$ has the same image as $\mathbf{T}$, thus if $y \in \overline{\im \mathbf{T}}$, there
is a sequence $x_k \in V$ such that $\mathbf{T} x_k \to y$.  We claim that $x_k$
is bounded.  If not, then $\mathbf{T}(x_k / \|x_k\|_X) \to 0$ and $\mathbf{K}(x_k / \|x_k\|_X)$
has a convergent subsequence, so \eqref{eqn:generalEstimate} implies that
a subsequence of $x_k / \|x_k\|_X$ also converges to some $x_\infty \in V$
with $\|x_\infty\| = 1$ and $\mathbf{T} x_\infty = 0$, a contradiction.
But now since $x_k$ is bounded, $\mathbf{K} x_k$ also has a convergent subsequence 
and $\mathbf{T} x_k$ converges by assumption, thus \eqref{eqn:generalEstimate}
yields also a convergent subsequence of $x_k$, whose limit $x$ satisfies
$\mathbf{T}x = y$.  This completes the proof that $\im \mathbf{T}$ is closed.
\end{proof}

In the analysis of closed $J$-holomorphic curves, one makes use of the above
lemma by placing the inclusion $W^{k-1,p} \hookrightarrow W^{k,p}$ in the
role of the compact operator~$\mathbf{K}$.  Unfortunately, 
$W^{k-1,p} \hookrightarrow W^{k,p}$ is not compact when the domain 
$\dot{\Sigma}$ has cylindrical ends; in contrast to the case of a compact
domain, there is no way to write the norm on the ends as a finite sum 
of norms for functions on domains of finite measure.  To circumvent this
problem, let
$$
\Sigma^R \subset \dot{\Sigma}
$$
denote the compact complement of the ends $\intZ^R_\pm \subset \dot{\uU}_z$ for
all $z \in \Gamma$.

\begin{lemma}
\label{lemma:semiFredEstimate}
Fix $k \in \NN$ and $1 < p < \infty$, and assume all the $\mathbf{A}_z$
are nondegenerate.  Then for sufficiently large $R > 0$, 
there exists a constant $c > 0$ such that
$$
\| \eta \|_{W^{k,p}(\dot{\Sigma})} \le 
c \| \mathbf{D} \eta \|_{W^{k-1,p}(\dot{\Sigma})} +
c \| \eta \|_{W^{k-1,p}(\Sigma^R)}
$$
for all $\eta \in W^{k,p}(E)$.
\end{lemma}
\begin{proof}
Fix a smooth cutoff function $\beta \in C_0^\infty(\Sigma^R)$
such that $\beta|_{\Sigma^{R-1}} \equiv 1$, and write
$$
\dot{\uU}_\Gamma^R \subset \dot{\Sigma}
$$
for the union of all the ends $\intZ_\pm^R \subset \dot{\uU}_z$ for
$z \in \Gamma^+ \cup \Gamma^-$.  Then we can write any
$\eta \in W^{k,p}(E)$ as $\eta = \beta \eta + (1 - \beta) \eta$ so that
$\beta \eta \in W^{k,p}_0(\Sigma^R)$ and $(1 - \beta) \eta
\in W^{k,p}_0(\dot{\uU}_\Gamma^{R-1})$.  Choosing $R$ large enough to make
Proposition~\ref{prop:ellipticEnd} valid, we can apply this together with
Proposition~\ref{prop:globalCpct} to show
\begin{equation*}
\begin{split}
\| \eta \|_{W^{k,p}(\dot{\Sigma})} &\le \| \beta \eta \|_{W^{k,p}(\Sigma^R)}
+ \| (1-\beta) \eta \|_{W^{k,p}(\dot{\uU}_\Gamma^{R-1})} \\
&\le c \| \mathbf{D}(\beta \eta) \|_{W^{k-1,p}(\Sigma^R)} +
c \| \beta \eta \|_{W^{k-1,p}(\Sigma^R)}  +
\| \mathbf{D} \left[ (1 - \beta) \eta \right] \|_{W^{k-1,p}(\dot{\uU}_\Gamma^{R-1})}.
\end{split}
\end{equation*}
After applying the Leipbniz rule and absorbing the norms of $\beta$ and
$\dbar\beta$ into the constants, this produces the
stated inequality since the term involving the $W^{k-1,p}$-norm of $\eta$
on the cylindrical ends includes $\dbar(1 - \beta)$, which vanishes
outside of~$\Sigma^R$.
\end{proof}

Lemma~\ref{lemma:semiFred} is now applicable since the operator
$$
W^{k,p}(\dot{\Sigma}) \to W^{k-1,p}(\Sigma^R) : \eta \mapsto 
\eta|_{\Sigma^R}
$$
involves the compact inclusion $W^{k,p}(\Sigma^R) \hookrightarrow
W^{k-1,p}(\Sigma^R)$ and is thus compact.

\begin{cor}
If all the $\mathbf{A}_z$ are nondegenerate, then
$$
\mathbf{D} : W^{k,p}(E) \to W^{k-1,p}(F)
$$
is semi-Fredholm.
\qed
\end{cor}

\section{Formal adjoints and proof of the Fredholm property}

In order to show that $\coker \mathbf{D}$ is also finite dimensional,
we will apply the above arguments to the formal adjoint of
$\mathbf{D}$, an operator whose kernel is naturally isomorphic to
the cokernel of~$\mathbf{D}$.  Let us choose Hermitian bundle
metrics $\langle\ , \rangle_E$ on $E$ and $\langle\ ,\ \rangle_F$ on~$F$,
and fix an area form $\vol$ on $\dot{\Sigma}$ that takes the form
$\vol = ds \wedge dt$ on the cylindrical ends.  The \defin{formal adjoint}
of $\mathbf{D}$ is then defined as the unique first-order linear
differential operator
$$
\mathbf{D}^* : \Gamma(F) \to \Gamma(E)
$$
that satisfies the relation
$$
\langle \lambda , \mathbf{D} \eta \rangle_{L^2(F)} =
\langle \mathbf{D}^* \lambda , \eta \rangle_{L^2(E)} \quad \text{ for all }\quad
\eta \in C_0^\infty(E),\ \lambda \in C_0^\infty(F),
$$
where we use the real-valued $L^2$-pairings
\begin{equation*}
\begin{split}
\langle \eta,\xi \rangle_{L^2(E)} &:= \Re \int_{\dot{\Sigma}} 
\langle \eta,\xi \rangle_E \, \vol, \quad \text{ for } \quad
\eta,\xi \in \Gamma(E), \\
\langle \alpha,\lambda \rangle_{L^2(F)} &:= \Re \int_{\dot{\Sigma}} 
\langle \alpha,\lambda \rangle_F \, \vol, \quad \text{ for } \quad
\alpha,\lambda \in \Gamma(F).
\end{split}
\end{equation*}
The word ``formal'' refers to the fact that we are not viewing $\mathbf{D}^*$
as the adjoint of an unbounded operator on a Hilbert 
space (cf.~\cite{ReedSimon}); that would be a stronger condition.

\begin{exercise}
\label{EX:antiCR}
Show that $\mathbf{D}^*$ is well defined and, for suitable choices of 
complex local trivializations of $E$ and $F$ and holomorphic coordinates on
open subsets $\uU \subset \dot{\Sigma}$, can be written locally as
$$
\mathbf{D}^* = -\p + A : C^\infty(\uU,\RR^{2n}) \to C^\infty(\uU,\RR^{2n})
$$
for some $A \in C^\infty(\uU,\End(\RR^{2n}))$, where $\p := \p_s - J_0 \p_t$.
\end{exercise}

The formula in the above exercise reveals that $\mathbf{D}^*$ is also an
elliptic operator\footnote{Technically, this property of the formal adjoint
is part of the definition of ellipticity: we call a differential operator
elliptic whenever (1)~it has the properties necessary for proving fundamental 
estimates using Fourier transforms as we did with $\dbar$ in 
\S\ref{sec:estimate}, and (2)~its formal adjoint also has this property.  The former
requires the principal symbol of the operator to be everywhere injective, and 
the latter requires it to be surjective.} 
and thus has the same local properties as~$\mathbf{D}$; 
indeed, $-\p + A$ can be transformed into $\dbar + B$ for some zeroth-order
term $B$ if we
conjugate it by a suitable complex-antilinear change of trivialization.
In particular, our local estimates for $\mathbf{D}$ and their consequences,
notably Proposition~\ref{prop:globalCpct}, are all equally 
valid for~$\mathbf{D}^*$.

To obtain suitable asymptotic estimates for $\mathbf{D}^*$, let us fix
asymptotic trivializations $\tau$ of $E$, use the corresponding
trivializations of $F$ over the ends as described in \S\ref{sec:setup},
and choose the bundle metrics such that both appear standard in these
trivializations over the ends.  We will say that the bundle metrics
are \defin{compatible with the asymptotically Hermitian structure} of $E$
whenever they are chosen in this way outside of a
compact subset of~$\dot{\Sigma}$.  We can then express $\mathbf{D}$ as
$\dbar + S(s,t)$ on $\dot{\uU}_z = Z_\pm$, and integrate by parts to obtain
$$
\mathbf{D}^* = -\p + S(s,t)^\transpose.
$$
To identify this expression with a Cauchy-Riemann type operator,
let $C := \begin{pmatrix} \1 & 0 \\ 0 & -\1 \end{pmatrix}$ denote the
$\RR$-linear transformation on $\RR^{2n} = \CC^n$ representing complex 
conjugation.  Then since $C$ anticommutes with $J_0$, we have
\begin{equation*}
\begin{split}
(C^{-1}\mathbf{D}^*C) \eta &= - C \p_s (C\eta) + C J_0 \p_t (C\eta) + 
C S(s,t)^\transpose C \eta \\
&= - \p_s \eta - J_0 \p_t \eta + C S(s,t)^\transpose C \eta 
= - (\dbar \eta - C S(s,t)^\transpose C \eta) \\
&=: - (\dbar + \bar{S}(s,t)) \eta,
\end{split}
\end{equation*}
where we've defined $\bar{S}(s,t) := - C S(s,t)^\transpose C$.
Now if the asymptotic operator $\mathbf{A}_z$ at $z \in \Gamma^\pm$ is
written in the chosen trivialization as $\mathbf{A} := -J_0 \p_s - S_\infty(t)$, 
the asymptotic convergence of $S(s,t)$ implies that similarly
$$
\| \bar{S} - \bar{S}_\infty \|_{C^k(Z_\pm^R)} \to 0 \quad \text{ as } 
\quad R \to \infty
$$
for all $k \in \NN$, where
$$
\bar{S}_\infty(t) := - C S_\infty(t) C.
$$
This defines a trivialized asymptotic operator 
$\overline{\mathbf{A}} = -J_0 \p_t - \bar{S}_\infty(t)$ to which
$- \mathbf{D}^*$ is (after a suitable change of trivialization) asymptotic
at the puncture~$z$; in particular, our proof of the global regularity result,
Proposition~\ref{prop:globalReg}, now also works for~$\mathbf{D}^*$.
Finally, notice that $\mathbf{A}$ and $-\overline{\mathbf{A}}$ are conjugate: 
indeed,
$$
(C^{-1} \overline{\mathbf{A}} C ) \eta = - C J_0 \p_t (C\eta) + C C S_\infty(t) C (C\eta)
= J_0 \p_t \eta + S_\infty(t) \eta = - \mathbf{A} \eta.
$$
This implies that $\mathbf{A}$ is nondegenerate if and only if 
$\overline{\mathbf{A}}$ is; applying this assumption for all of the 
$\mathbf{A}_z$, the proofs of Proposition~\ref{prop:ellipticEnd} and
Lemma~\ref{lemma:semiFredEstimate} now also go through for~$\mathbf{D}^*$.

We've proved:

\begin{prop}
\label{prop:adjointSemi}
Suppose $\mathbf{D}^*$ is defined with respect Hermitian bundle metrics on
$E$ and $F = \overline{\Hom}_\CC(T\dot{\Sigma},E)$ that are compatible with
the asymptotically Hermitian structure of~$E$.
If additionally all the asymptotic operators $\mathbf{A}_z$ are nondegenerate, 
then
$$
\mathbf{D}^* : W^{k,p}(F) \to W^{k-1,p}(E)
$$
is semi-Fredholm, and its kernel is a space of smooth sections contained
in $W^{m,q}(F)$ for all $m \in \NN$ and $q \in (1,\infty)$.
\qed
\end{prop}

Since $\ker \mathbf{D}^*$ is now known to be finite dimensional, the next
result completes the proof of the Fredholm property for $\mathbf{D}$ by
showing that its image has finite codimension:

\begin{lemma}
\label{lemma:codim}
Under the same assumptions as in Proposition~\ref{prop:adjointSemi},
$$
W^{k-1,p}(F) = \im \mathbf{D} + \ker \mathbf{D}^*.
$$
\end{lemma}
\begin{proof}
Consider first the case $k=1$.
Since $\mathbf{D} : W^{1,p}(E) \to L^p(F)$ is semi-Fredholm, its image is 
closed, hence $\im \mathbf{D} + \ker \mathbf{D}^*$ is a closed subspace 
of~$L^p(F)$.  Then if $\im \mathbf{D} + \ker \mathbf{D}^*
\ne L^p(F)$, the Hahn-Banach theorem\footnote{In the case $p=2$, one can
forego the Hahn-Banach theorem and simply take an $L^2$-orthogonal complement.}
provides a nontrivial element $\alpha \in \left(L^p(F)\right)^* \cong L^q(F)$
for $\frac{1}{p} + \frac{1}{q} = 1$ such that
\begin{equation}
\label{eqn:conditions}
\langle \mathbf{D}\eta + \lambda , \alpha \rangle_{L^2(F)} = 0 \quad 
\text{ for all } \quad
\eta \in W^{1,p}(E),\ \lambda \in \ker \mathbf{D}^*.
\end{equation}
Choosing $\lambda=0$, this implies in particular
$$
\langle \mathbf{D}\eta , \alpha \rangle_{L^2(F)} = 0 \quad
\text{ for all } \quad
\eta \in C_0^\infty(E),
$$
which means that $\alpha$ is a weak solution of class $L^q$ to the formal
adjoint equation $\mathbf{D}^*\alpha = 0$.  By Proposiiton~\ref{prop:globalReg},
$\alpha$ is therefore smooth and belongs to $\ker\mathbf{D}^*$.
But this contradicts \eqref{eqn:conditions} if we plug in $\eta=0$ and
$\lambda = \alpha$, so this completes the proof for $k=1$.

For $k \ge 2$, suppose $\alpha \in W^{k-1,p}(F) \subset L^p(F)$ is given: 
then the case $k=1$ provides elements $\eta \in W^{1,p}(E)$ and $\lambda \in
\ker \mathbf{D}^*$ such that $\mathbf{D}\eta + \lambda = \alpha$.  Since
Proposition~\ref{prop:globalReg} implies $\lambda \in W^{m,q}(F)$ for all
$m \in \NN$ and $q \in (1,\infty)$, we have $\mathbf{D}\eta = \alpha - \lambda
\in W^{k-1,p}(F)$ and thus, by Prop.~\ref{prop:globalReg} again,
$\eta \in W^{k,p}(E)$, completing the proof for all $k \in \NN$.
\end{proof}

The proof of Theorem~\ref{thm:Fredholm} is now complete, but as long as
we're talking about the formal adjoint, let us take note of a few more
properties that will be useful in the future.  Assume from now on that
all the assumptions of Proposition~\ref{prop:adjointSemi} are satisfied.
We can now strengthen Lemma~\ref{lemma:codim} as follows.

\begin{prop}
\label{prop:duality}
$W^{k-1,p}(F) = \im \mathbf{D} \oplus \ker \mathbf{D}^*$ and
$W^{k-1,p}(E) = \im \mathbf{D}^* \oplus \ker \mathbf{D}$.  In particular,
the projections defined by these splittings give isomorphisms
$$
\coker \mathbf{D} \cong \ker \mathbf{D}^* \quad\text{ and }\quad
\coker \mathbf{D}^* \cong \ker \mathbf{D},
$$
thus $\mathbf{D}^* : W^{k,p}(F) \to W^{k-1,p}(E)$ is a Fredholm operator with
$$
\ind \mathbf{D}^* = - \ind \mathbf{D}.
$$
\end{prop}
\begin{proof}
By Lemma~\ref{lemma:codim}, the first splitting follows if we can show
that $\im \mathbf{D} \cap \ker \mathbf{D}^* = \{0\}$.
Recall first (see \S\ref{sec:SobolevEnds})
that $C_0^\infty(\dot{\Sigma})$ is dense in
$W^{k,p}(\dot{\Sigma})$ for every $k \ge 0$ and $p \in [1,\infty)$, so the
definition of the formal adjoint implies via density and
H\"older's inequality that if $1 < p,q < \infty$ and 
$\frac{1}{p} + \frac{1}{q} = 1$,
\begin{equation}
\label{eqn:adjointHoelder}
\langle \lambda,\mathbf{D} \eta \rangle_{L^2(F)} = 
\langle \mathbf{D}^*\lambda , \eta \rangle_{L^2(E)} \quad \text{ for all }
\quad \eta \in W^{1,p}(E),\ \lambda \in W^{1,q}(F).
\end{equation}
Now suppose $\lambda \in \im \mathbf{D} \cap \ker \mathbf{D}^*$ and write
$\lambda = \mathbf{D}\eta$, assuming $\eta \in W^{k,p}(E)$.  Regularity
implies that since $\mathbf{D}^*\lambda = 0$, $\lambda \in W^{1,q}(F)$,
where $q$ can be chosen to satisfy $\frac{1}{p} + \frac{1}{q} = 1$.
We can therefore apply \eqref{eqn:adjointHoelder} and obtain
$$
\langle \lambda , \lambda \rangle_{L^2(F)} = 
\langle \lambda , \mathbf{D}\eta \rangle_{L^2(F)} =
\langle \mathbf{D}^*\lambda , \eta \rangle_{L^2(E)} = 0,
$$
hence $\lambda = 0$.

The proof that $W^{k-1,p}(E) = \im \mathbf{D}^* \oplus \ker \mathbf{D}$
is analogous.
\end{proof}

This result hints at the fact that $\mathbf{D}^*$ is in fact---under some
natural extra assumptions---globally equivalent to another Cauchy-Riemann 
type operator.  To see this, let us impose a further constraint on the
relation between the Hermitian bundle metrics $\langle\ ,\ \rangle_E$ and
$\langle\ ,\ \rangle_F$.  Note that since the area form $\vol$ is necessarily
$j$-invariant, it induces a Hermitian structure on $T\dot{\Sigma}$, namely
$$
\langle X,Y \rangle_\Sigma := \vol(X,jY) + i \, \vol(X,Y),
$$
which matches the standard bundle metric in the trivializations over the
ends defined via the cylindrical coordinates.
This induces real-linear isomorphisms from $T\dot{\Sigma}$ to the
complex-linear and -antilinear parts of the complexified cotangent bundle,
\begin{equation*}
\begin{split}
T\dot{\Sigma} \to \Lambda^{1,0}T^*\dot{\Sigma} : X \mapsto X^{1,0} := \langle X,\cdot\rangle_\Sigma, \\
T\dot{\Sigma} \to \Lambda^{0,1}T^*\dot{\Sigma} : X \mapsto X^{0,1} := \langle \cdot,X \rangle_\Sigma,
\end{split}
\end{equation*}
where the first isomorphism is complex antilinear and the second is complex linear.
We use these to define Hermitian bundle metrics on $\Lambda^{1,0}T^*\dot{\Sigma}$
and $\Lambda^{0,1}T^*\dot{\Sigma}$ in terms of the metric on $T\dot{\Sigma}$;
note that this is a straightforward definition for $\Lambda^{0,1}T^*\dot{\Sigma}$,
but since the isomorphism to $\Lambda^{1,0}T^*\dot{\Sigma}$ is
complex \emph{antilinear}, we really mean
$$
\langle X^{1,0} , Y^{1,0} \rangle_\Sigma := \langle Y,X \rangle_\Sigma \quad
\text{ for } \quad X,Y \in T\dot{\Sigma}.
$$
Now observe that as a vector bundle with complex structure
$\lambda \mapsto J \circ \lambda$,
$F = \overline{\Hom}_\CC(T\dot{\Sigma},E)$ is naturally isomorphic to the
complex tensor product
$$
F = \Lambda^{0,1}T^*\Sigma \otimes E.
$$
We can therefore make a natural choice for $\langle\ ,\ \rangle_F$ as the
tensor product metric determined by $\langle\ ,\ \rangle_\Sigma$
and $\langle\ ,\ \rangle_E$.  It is easy to check that this choice is 
compatible with the asymptotically Hermitian structure of~$E$.

Next, we notice that the area form $\vol$ also induces a natural complex
bundle isomorphism
$$
E \to \Hom_\CC(T\dot{\Sigma},F).
$$
Indeed, the right hand side is canonically isomorphic to the complex tensor
product
$$
\Hom_\CC(T\dot{\Sigma},F) = \Lambda^{1,0}T^*\dot{\Sigma} \otimes F =
\Lambda^{1,0}T^*\dot{\Sigma} \otimes \Lambda^{0,1}T^*\dot{\Sigma} \otimes E,
$$
and $\Lambda^{1,0}T^*\dot{\Sigma} \otimes \Lambda^{0,1}T^*\dot{\Sigma}$
is isomorphic to the trivial complex
line bundle $\epsilon^1 := \dot{\Sigma} \times \CC \to \dot{\Sigma}$ via
$$
\Lambda^{1,0}T^*\dot{\Sigma} \otimes \Lambda^{0,1}T^*\dot{\Sigma} \to 
\epsilon^1 : X^{1,0} \otimes Y^{0,1} \mapsto X^{1,0}(Y) = 
\langle X,Y \rangle_\Sigma.
$$

\begin{exercise}
Assuming $\langle\ ,\ \rangle_F$ is chosen as the tensor product metric
described above, show that under the natural identification of $E$ with
$\Hom_\CC(T\dot{\Sigma},F)$,
$$
-\mathbf{D}^* : \Gamma(F) \to \Omega^{1,0}(\dot{\Sigma},F)
$$
satisfies the Leibniz rule
$$
-\mathbf{D}^*(f \lambda) = (\p f) \lambda + f (-\mathbf{D}^*\lambda)
$$
for all $f \in C^\infty(\dot{\Sigma},\RR)$, where
$\p f \in \Omega^{1,0}(\dot{\Sigma})$ denotes the complex-valued
$(1,0)$-form $df - i\, df \circ j$.
\end{exercise}

We might summarize this exercise by saying that $-\mathbf{D}^*$ is 
an ``anti-Cauchy-Riemann type'' operator on~$F$.
But such an object is easily transformed into an honest Cauchy-Riemann type
operator: let $\bar{F}$ denote the \defin{conjugate bundle} to $F$,
which we define as the same real vector bundle $F$ but with the sign of its
complex structure reversed, so
$\lambda \mapsto -J \circ \lambda$.  Now there is a canonical isomorphism
$$
\Hom_\CC(T\dot{\Sigma},F) = \overline{\Hom}_\CC(T\dot{\Sigma},\bar{F}),
$$
and the same operator defines a real-linear map
$$
-\mathbf{D}^* : \Gamma(\bar{F}) \to \Omega^{0,1}(\dot{\Sigma},\bar{F})
$$
which satisfies our usual Leibniz rule for Cauchy-Riemann type operators.

Its asymptotic behavior also fits into the scheme we've been describing:
we have already seen this by computing $\mathbf{D}^*$ on the ends with respect
to asymptotic trivializations.  To express this in trivialization-invariant
language, observe that each of the Hermitian bundles $(E_z,J_z,\omega_z)$
over $S^1$ for $z \in \Gamma$ has a conjugate bundle
$\bar{E}_z$ with complex structure $-J_z$ and symplectic structure
$-\omega_z$; its natural Hermitian inner product is then the complex conjugate
of the one on~$E_z$.  The asymptotic operator $\mathbf{A}_z$ on $E_z$
can be expressed as $-J_z \widehat{\nabla}_t$, where $\widehat{\nabla}_t$
is a symplectic connection on $(E_z,\omega_z)$.  Then $\widehat{\nabla}_t$
is also a symplectic connection on $(\bar{E}_z,-\omega_z)$, so we
naturally obtain an asymptotic operator on $\bar{E}_z$ in the form
\begin{equation}
\label{eqn:conjAsymp}
\overline{\mathbf{A}}_z := - \mathbf{A}_z : \Gamma(\bar{E}_z) \to
\Gamma(\bar{E}_z),
\end{equation}
where the sign reversal arises from the reversal of the complex structure.
One can check that if we choose a unitary trivialization of $E_z$ and
the conjugate trivialization of $\bar{E}_z$, this relationship between
$\mathbf{A}_z$ and $\overline{\mathbf{A}}_z$ produces precisely the
relationship between $\mathbf{A} = -J_0 \p_t - S_\infty(t)$ and
$\overline{\mathbf{A}} = -J_0 \p_t - \bar{S}_\infty(t)$ that we saw
previously, with $\bar{S}_\infty(t) = -C S_\infty(t) C$.  Let us summarize
all this with a theorem.

\begin{thm}
\label{thm:adjoint}
Assume $\langle\ ,\ \rangle_F$ is chosen to be the tensor product metric
on $F = \Lambda^{0,1}T^*\Sigma \otimes E$ induced by $\langle\ ,\ \rangle_E$
and the area form~$\vol$.  Then under the isomorphism induced by $\vol$
from $E$ to $\Hom_\CC(T\dot{\Sigma},F)$ and the natural identification of
the latter with its conjugate $\overline{\Hom}_\CC(T\dot{\Sigma},\bar{F})$,
the operator $-\mathbf{D}^* : \Gamma(F) \to \Gamma(E)$ defines a linear
Cauchy-Riemann type operator on the conjugate bundle~$\bar{F}$,
$$
-\mathbf{D}^* : \Gamma(\bar{F}) \to \Omega^{0,1}(\dot{\Sigma},\bar{F}),
$$
and it is asymptotic at each puncture $z \in \Gamma$ to the conjugate
asymptotic operator \eqref{eqn:conjAsymp}.
\qed
\end{thm}

\chapter{The index formula}
\label{lec:index}

\minitoc

\vspace{12pt}

\section{Riemann-Roch with punctures}
\label{sec:RiemannRoch}

As in the previous lecture, let $\mathbf{D}$ denote a linear Cauchy-Riemann
type operator on an asymptotically Hermitian vector bundle $E$ of complex
rank~$m$ over a 
punctured Riemann surface $(\dot{\Sigma} = \Sigma \setminus (\Gamma^+
\cup \Gamma^-),j)$, and assume that $\mathbf{D}$ is asymptotic at each
puncture $z \in \Gamma$ to a nondegenerate asymptotic operator $\mathbf{A}_z$
on the asymptotic bundle $(E_z,J_z,\omega_z)$ over~$S^1$.  
Writing
$$
F := \overline{\Hom}_\CC(T\dot{\Sigma},E)
$$
for the bundle of complex-antilinear homomorphisms $T\dot{\Sigma} \to E$,
the main result of the previous lecture was that
$$
\mathbf{D} : W^{k,p}(E) \to W^{k-1,p}(F)
$$
is Fredholm for any $k \in \NN$ and $p \in (1,\infty)$, and its kernel
and index do not depend on $k$ or~$p$.  The main goal of this lecture is
to compute $\ind(\mathbf{D}) \in \ZZ$.

The index will depend on the Conley-Zehnder indices
$\muCZ^\tau(\mathbf{A}_z) \in \ZZ$ introduced in Lecture~\ref{lec:asymptotic}, but
since these depend on arbitrary choices of unitary trivializations~$\tau$,
we need a way of selecting preferred trivializations.  The most natural
condition is to require that every $(E_z,J_z,\omega_z)$ be endowed with
a unitary trivialization such that the corresponding asymptotic 
trivializations of $(E,J)$ extend to a global 
trivialization\footnote{Note that $(E,J)$ is always globally 
trivializable unless $\Gamma = \emptyset$, as
a punctured surface can be retracted to its 1-skeleton.}; if there is
only one puncture~$z$, for instance, then this condition determines
$\muCZ^\tau(\mathbf{A}_z)$ uniquely.  This convention has been used to state
the formula for $\ind(\mathbf{D})$ in several of the standard references,
e.g.~in \cite{HWZ:props3}.  We would prefer
however to state a formula which is also valid when $\Gamma = \emptyset$
and $E \to \Sigma$ is nontrivial.  One way to do this is by allowing
completely arbitrary asymptotic trivializations, but introducing a topological
invariant to measure their failure to extend globally over~$E$.

\begin{defn}
Fix a compact oriented surface $S$ with boundary.  The \defin{relative first Chern
number} associates to every complex vector bundle $(E,J)$ over $S$ and
trivialization $\tau$ of $E|_{\p S}$ an integer
$$
c_1^\tau(E) \in \ZZ
$$
satisfying the following properties:
\begin{enumerate}
\item
If $(E,J) \to S$ is a line bundle, then $c_1^\tau(E)$ is the signed count
of zeroes for a generic smooth section $\eta \in \Gamma(E)$ that appears
as a nonzero constant at $\p S$ with respect to~$\tau$.
\item
For any two bundles $(E_1,J_1)$ and $(E_2,J_2)$ with trivializations
$\tau_1$ and $\tau_2$ respectively over $\p S$,
$$
c_1^{\tau_1 \oplus \tau_2}(E_1 \oplus E_2) = c_1^{\tau_1}(E_1) + c_1^{\tau_2}(E_2).
$$
\end{enumerate}
\end{defn}

These two conditions uniquely determine $c_1^\tau(E)$ for all complex vector
bundles since bundles of higher rank can always be split into direct sums
of line bundles.  The definition clearly matches the usual first Chern
number $c_1(E)$ when $\p S = \emptyset$, and it extends in an obvious way
to the category of asymptotically Hermitian vector bundles with asymptotic
trivializations.

\begin{exercise}
\label{EX:c1Change}
Given two distinct choices of asymptotic trivializations $\tau_1$ and $\tau_2$
for an asymptotically Hermitian bundle $E$ of rank~$m$, show that
$$
c_1^{\tau_2}(E) = c_1^{\tau_1}(E) - \deg(\tau_2 \circ \tau_1^{-1}),
$$
where $\deg(\tau_2 \circ \tau_1^{-1}) \in \ZZ$ denotes the sum over all punctures
of the winding numbers of the determinants of the transition maps
$S^1 \to \U(m)$.\footnote{Caution: to compute this winding number at a
negative puncture using cylindrical coordinates $(s,t) \in (-\infty,0] \times S^1$,
one must traverse $\{-s\} \times S^1$ for $s \gg 1$ in the \emph{wrong direction}, as
this is consistent with the orientation induced on $\{-s\} \times S^1$ as a 
boundary component of a large compact subdomain of~$\dot{\Sigma}$.}
\end{exercise}

\begin{exercise}
\label{EX:sumIndep}
Combining Exercise~\ref{EX:c1Change} above with Exercise~\ref{EX:changeTriv},
show that for our asymptotically Hermitian vector bundle
$E$ with Cauchy-Riemann type operator $\mathbf{D}$ and asymptotic
operators $\mathbf{A}_z$, the number
$$
2 c_1^\tau(E) + \sum_{z \in \Gamma^+} \muCZ^\tau(\mathbf{A}_z) -
\sum_{z \in \Gamma^-} \muCZ^\tau(\mathbf{A}_z)
$$
is independent of the choice of asymptotic trivializations~$\tau$.
\end{exercise}

The above exercise shows that the right hand side of the following index
formula is independent of all choices.

\begin{thm}
\label{thm:RiemannRoch}
The Fredholm index of $\mathbf{D}$ is given by
$$
\ind \mathbf{D} = m \chi(\dot{\Sigma}) + 2 c_1^\tau(E) +
\sum_{z \in \Gamma^+} \muCZ^\tau(\mathbf{A}_z) -
\sum_{z \in \Gamma^-} \muCZ^\tau(\mathbf{A}_z),
$$
where $m = \rank_\CC E$ and
$\tau$ is an arbitrary choice of asymptotic trivializations.
\end{thm}

\begin{notation}
Throughout this lecture, we shall denote the integer on the right hand
side in Theorem~\ref{thm:RiemannRoch} by
$$
I(\mathbf{D}) := m \chi(\dot{\Sigma}) + 2 c_1^\tau(E) +
\sum_{z \in \Gamma^+} \muCZ^\tau(\mathbf{A}_z) -
\sum_{z \in \Gamma^-} \muCZ^\tau(\mathbf{A}_z) \in \ZZ.
$$
Our goal is thus to prove that $\ind(\mathbf{D}) = I(\mathbf{D})$.
\end{notation}

When $\Gamma = \emptyset$, Theorem~\ref{thm:RiemannRoch} is equivalent to
the classical Riemann-Roch formula, which is more often stated for
\emph{holomorphic} vector bundles over a closed Riemann surface $(\Sigma,j)$
with genus $g$~as
\begin{equation}
\label{eqn:RiemannRochClassical}
\ind_\CC(\mathbf{D}_0) = m(1 - g) + c_1(E).
\end{equation}
This formula assumes that the Cauchy-Riemann type operator $\mathbf{D}_0$ is 
complex linear, but an arbitrary
real-linear Cauchy-Riemann operator is then of the form
$\mathbf{D} = \mathbf{D}_0 + B$, where the zeroth-order term 
$B \in \Gamma(\Hom_\RR(E,F))$
defines a compact perturbation since the inclusion $W^{k,p}(\Sigma) 
\hookrightarrow W^{k-1,p}(\Sigma)$ is compact.  It follows that $\mathbf{D}$ 
has the same \emph{real} Fredholm index as $\mathbf{D}_0$, namely twice the
complex index shown on the 
right hand side of \eqref{eqn:RiemannRochClassical}, which matches what we
see in Theorem~\ref{thm:RiemannRoch}.

\begin{remark}
Now seems a good moment to clarify explicitly that all dimensions (and therefore
also Fredholm indices) in this lecture are \emph{real} dimensions, not
complex dimensions, unless otherwise stated.
\end{remark}

Reduction to the complex-linear case does not work in general if there
are punctures: it remains true that arbitrary Cauchy-Riemann type operators
can be written as $\mathbf{D} = \mathbf{D}_0 + B$ where $\mathbf{D}_0$ is
complex linear, but the perturbation introduced by the zeroth-order term $B$
is not compact since $W^{k,p}(\dot{\Sigma}) \hookrightarrow 
W^{k-1,p}(\dot{\Sigma})$ is not compact when $\Gamma \ne \emptyset$.
Another indication that this idea cannot work is the fact that while
the formula in Theorem~\ref{thm:RiemannRoch} always gives an \emph{even}
integer when $\Gamma = \emptyset$, it can be odd when there are punctures,
in which case $\mathbf{D}$ clearly cannot have the same index is any
complex-linear operator.  Our proof will therefore have to deal with more
than just the complex category.

The punctured version of Theorem~\ref{thm:RiemannRoch} was first proved by
Schwarz in his thesis \cite{Schwarz}, its main purpose at the time being to
help define algebraic
operations (notably the \emph{pair-of-pants product}) in Hamiltonian
Floer homology.  Schwarz's proof used a ``linear gluing'' construction
that gives a relation between indices of operators on bundles over surfaces 
obtained
by gluing together constituent surfaces along matching cylindrical ends.
Since any surface with ends can be ``capped off'' to form a closed surface,
one obtains the general index formula if one already knows how to compute
it for closed surfaces and for planes (i.e.~caps).  For the latter, it is
simple enough to write down model Cauchy-Riemann operators on planes and
compute their kernels and cokernels explicitly, so in this way the general
case is reduced to the classical Riemann-Roch formula.  An analogous
linear gluing argument for compact surfaces with boundary is used in
\cite{McDuffSalamon:Jhol}*{Appendix~C} to reduce the general Riemann-Roch 
formula to an explicit computation for Cauchy-Riemann operators on the disk 
with a totally real boundary condition.

In this lecture, we will follow a different path and use an argument that was
first sketched by Taubes for the closed case in \cite{Taubes:counting}*{\S 7}, 
with an additional argument
for the punctured case suggested by Chris Gerig \cite{Gerig:thesis}.
The argument is (in my opinion)
analytically somewhat easier than the more standard approaches, and in
addition to proving the formula we need for punctured surfaces, it
produces a new proof in the closed case without assuming the classical
Riemann-Roch formula.  It also provides a gentle preview of two analytical
phenomena that will later assume prominent roles in our discussion of SFT:
\emph{bubbling} and \emph{gluing}.

To see the idea behind Taubes's argument, we can start by noticing an apparent
numerical coincidence in the closed case.  Assume $(E,J)$ is a complex line
bundle over a closed Riemann surface $(\Sigma,j)$, and 
$\mathbf{D} : \Gamma(E) \to \Gamma(F) = \Omega^{0,1}(\Sigma,E)$ 
is a Cauchy-Riemann type operator.  We know that $\ind(\mathbf{D}) =
\ind(\mathbf{D} + B)$ for any zeroth-order term 
$B \in \Gamma(\Hom_\RR(E,F))$.  But $E$ and
$F$ are both complex vector bundles, so $B$ can always be split uniquely into
its complex-linear and complex-antilinear parts, i.e.~there is a natural
splitting of $\Hom_\RR(E,F)$ into a direct sum of complex line 
bundles\footnote{Here the complex structure
on $\Hom_\RR(E,F)$ and its subbundles is defined in terms of the complex 
structure of $F$, i.e.~it sends 
$B \in \Hom_\RR(E,F)$ to $J \circ B \in \Hom_\RR(E,F)$.}
$$
\Hom_\RR(E,F) = \Hom_\CC(E,F) \oplus \overline{\Hom}_\CC(E,F).
$$
Out of curiosity, let's compute the first Chern number of the second factor;
this will be the signed count of zeroes of a generic complex-\emph{antilinear}
zeroth-order perturbation.  To start with, note that
$$
\overline{\Hom}_\CC(E,F) = \overline{\Hom}_\CC(E,\CC) \otimes F,
$$
and then observe that $\overline{\Hom}_\CC(E,\CC)$ and $E$ are isomorphic:
indeed, any Hermitian bundle metric $\langle\ ,\ \rangle_E$ on $E$ gives rise
to a bundle isomorphism\footnote{We are assuming as usual that Hermitian
inner products are complex antilinear in the first argument and linear in
the second.}
$$
E \to \overline{\Hom}_\CC(E,\CC) : \eta \mapsto \langle \cdot,\eta \rangle_E.
$$
We thus have $\overline{\Hom}_\CC(E,F) \cong E \otimes F$, so
$c_1(\overline{\Hom}_\CC(E,F)) = c_1(E) + c_1(F)$.  We can compute $c_1(F)$
by the same trick since
$$
F = \overline{\Hom}_\CC(T\Sigma,E) = \overline{\Hom}_\CC(T\Sigma,\CC) \otimes E
\cong T\Sigma \otimes E,
$$
so $c_1(F) = c_1(T\Sigma) + c_1(E) = \chi(\Sigma) + c_1(E)$, and thus
$$
c_1(\overline{\Hom}_\CC(E,F)) = \chi(\Sigma) + 2 c_1(E).
$$
Since we're looking at a line bundle over a surface without punctures,
this number is the same as $I(\mathbf{D})$.  This coincidence is too 
improbable to ignore, and indeed, it turns out not to be coincidental.
Here is an informal statement of a result that we will later prove
a more precise version of in order to deduce Theorem~\ref{thm:RiemannRoch}.

\begin{heoremu}
Given a Cauchy-Riemann type operator $\mathbf{D} : H^1(E) \to L^2(F)$ 
on a line bundle $(E,J)$
over a closed Riemann surface $(\Sigma,j)$, choose a complex-antilinear
zeroth-order perturbation $B \in \Gamma(\overline{\Hom}_\CC(E,F))$
whose zeroes are all nondegenerate.  Then for sufficiently large $\sigma > 0$,
$\ker (\mathbf{D} + \sigma B)$ is approximately spanned by $1$-dimensional
spaces of sections with support localized near the positive zeroes of~$B$.
In particular, $\dim \ker(\mathbf{D} + \sigma B)$ equals the number of
positive zeroes of~$B$.
\end{heoremu}

To deduce $\ind(\mathbf{D}) = I(\mathbf{D})$ from this, we need
to apply the same trick to the formal adjoint~$\mathbf{D}^*$.
As we will review in \S\ref{sec:adjoint}, $-\mathbf{D}^*$ can be regarded 
under certain natural assumptions as a 
Cauchy-Riemann type operator on the bundle $\bar{F}$ conjugate to~$F$,
and the formal adjoint of $\mathbf{D} + \sigma B$ then gives rise to a
Cauchy-Riemann type operator of the form
$$
-\mathbf{D}^* + \sigma B' : \Gamma(\bar{F}) \to \Gamma(\bar{E}) = 
\Omega^{0,1}(\Sigma,\bar{F}),
$$
where $B' : \bar{F} \to \bar{E}$ is also complex antilinear and has the
same zeroes as~$B$, but with opposite signs.  Applying the above
``theorem'' to $-\mathbf{D}^*$ thus identifies 
$\ker (\mathbf{D} + \sigma B)^*$ for sufficiently large $\sigma > 0$ with a
space whose dimension equals the number of \emph{negative} zeroes of~$B$.
This gives
\begin{equation*}
\begin{split}
\ind(\mathbf{D}) &= \ind(\mathbf{D} + \sigma B) = \dim \ker (\mathbf{D} + \sigma B) - 
\dim \ker (\mathbf{D} + \sigma B)^* \\ &= c_1(\overline{\Hom}_\CC(E,F))
= I(\mathbf{D}).
\end{split}
\end{equation*}

It's worth mentioning that the ``large perturbation'' argument we've just 
sketched is only one simple example of an idea with a long
and illustrious history: another simple example is the observation by
Witten \cite{Witten:Morse} that after choosing a Morse function on a
Riemannian manifold, certain large deformations of the
de Rham complex lead to an approximation of the
Morse complex, with generators of the de Rham complex having support
concentrated near the critical points of the Morse function---this yields
a somewhat novel proof of de Rham's theorem.
A much deeper example is Taubes's isomorphism 
\cite{Taubes:SWtoGr} between the
Seiberg-Witten invariants of symplectic $4$-manifolds and certain
holomorphic curve invariants: here also, the idea is to consider a
large compact perturbation of the Seiberg-Witten equations and show that,
in the limit where the perturbation becomes infinitely large, solutions
of the Seiberg-Witten equations localize near $J$-holomorphic curves.
For a more recent exploration of this idea in the context of Dirac operators,
see~\cite{Maridakis}.

Before proceeding with the details, let us fix two simplifying assumptions that
can be imposed without loss of generality:

\begin{assumption}
$(E,J)$ has complex rank~$1$.
\end{assumption}

Indeed, an asymptotically Hermitian bundle $E$ of complex rank $m \in \NN$ 
always admits a decomposition into
asymptotically Hermitian line bundles $E = E_1 \oplus \ldots \oplus E_m$,
producing a corresponding splitting of the target bundle
$F = F_1 \oplus \ldots \oplus F_m$.  The operator $\mathbf{D}$ need not
respect these splittings, but it is always \emph{homotopic through Fredholm
operators} to one that does: we saw in Theorem~\ref{thm:CZclassification} that
the asymptotic operators $\mathbf{A}_z$
are homotopic through nondegenerate asymptotic operators to any other
operators $\mathbf{A}_z'$ that have the same Conley-Zehnder
indices, so one can choose $\mathbf{A}_z'$ to respect the splitting.  Any
homotopy of Cauchy-Riemann operators following such a homotopy of 
nondegenerate asymptotic operators then produces a continuous family of
Fredholm operators by the main result of Lecture~\ref{lec:Fredholm}, implying that their
indices do not change.  The general index formula then follows from the
line bundle case since any two Cauchy-Riemann type Fredholm operators
$\mathbf{D}_1$ and $\mathbf{D}_2$ over the same Riemann surface satisfy
$$
\ind(\mathbf{D}_1 \oplus \mathbf{D}_2) = \ind(\mathbf{D}_1) + \ind(\mathbf{D}_2)
\quad\text{ and }\quad
I(\mathbf{D}_1 \oplus \mathbf{D}_2) = I(\mathbf{D}_1) + I(\mathbf{D}_2).
$$

\begin{assumption}
$k=1$ and $p=2$.
\end{assumption}
This means we will concretely be considering the operator
$$
\mathbf{D} : H^1(E) \to L^2(F),
$$
where $H^1$ as usual is an abbreviation for~$W^{1,2}$.  This assumption is
clearly harmless since we know that $\ind \mathbf{D}$ does not depend on
the choice of $k$ and~$p$.

\section{Some remarks on the formal adjoint}
\label{sec:adjoint}

For the beginning of this section we can drop the assumption that $(E,J)$ is
a line bundle and assume $\rank_\CC E = m \in \NN$, though later we will
again set $m=1$.

Recall from the end of Lecture~\ref{lec:Fredholm} that if we fix global Hermitian structures
$\langle\ ,\ \rangle_E$ and $\langle\ ,\ \rangle_F$ on $(E,J)$ and $(F,J)$ 
respectively and an area form $\vol$ on $\dot{\Sigma}$ that
matches $ds \wedge dt$ on the cylindrical ends, then $\mathbf{D}$ has a
\emph{formal adjoint}
$$
\mathbf{D}^* : \Gamma(F) \to \Gamma(E)
$$
satisfying
$$
\langle \lambda,\mathbf{D}\eta \rangle_{L^2(F)} = \langle \mathbf{D}^*\lambda,
\eta \rangle_{L^2(E)} \quad\text{ for all }\quad
\eta \in H^1(E),\ \lambda \in H^1(F).
$$
Here the real-valued $L^2$~pairings are defined by
$$
\langle \eta,\xi \rangle_{L^2(E)} := \Re \int_{\dot{\Sigma}}
\langle \eta,\xi \rangle_E \, \vol \quad \text{ for }\quad
\eta,\xi \in \Gamma(E),
$$
and similarly for sections of~$F$.  The essential features of the formal
adjoint are that $\ker \mathbf{D}^* \cong \coker \mathbf{D}$ and 
$\coker \mathbf{D}^* \cong \ker \mathbf{D}$, hence
$\ind(\mathbf{D}^*) = - \ind(\mathbf{D})$.  
Recall moreover that $\vol$ induces a natural Hermitian bundle metric
on $\dot{\Sigma}$ by
$$
\langle\cdot ,\cdot \rangle_\Sigma = \vol(\cdot,j\cdot) + i \, \vol(\cdot,\cdot),
$$
which determines a bundle isomorphism
$$
T\dot{\Sigma} \to \Lambda^{0,1}T^*\dot{\Sigma} : X \mapsto X^{0,1} :=
\langle \cdot,X \rangle_\Sigma,
$$
as well as a complex-\emph{antilinear} isomorphism
$$
T\dot{\Sigma} \to \Lambda^{1,0}T^*\dot{\Sigma} : X \mapsto X^{1,0} :=
\langle X,\cdot \rangle_\Sigma.
$$
If $\langle\ ,\ \rangle_F$ is then chosen to be the tensor product 
metric determined via the natural isomorphism
$$
F = \overline{\Hom}_\CC(T\dot{\Sigma},E) = \Lambda^{0,1}T^*\dot{\Sigma} \otimes E =
T\dot{\Sigma} \otimes E,
$$
then $E$ admits a natural isomorphism to $\Lambda^{1,0}T^*\dot{\Sigma} \otimes F$
such that
$$
-\mathbf{D}^* : \Gamma(F) \to \Gamma(E) = \Omega^{1,0}(\dot{\Sigma},F)
$$
becomes an \emph{anti-Cauchy-Riemann} type operator, i.e.~it satisfies the
Leibniz rule
$$
-\mathbf{D}^* (f\lambda) = (\p f) \lambda + f (-\mathbf{D}^* \lambda)
$$
for all $f \in C^\infty(\dot{\Sigma},\RR)$, with $\p f := df - i\, df \circ j
\in \Omega^{1,0}(\dot{\Sigma})$.  Equivalently, $-\mathbf{D}^*$ defines a
Cauchy-Riemann type operator on the \defin{conjugate} bundle
$\bar{F} \to \dot{\Sigma}$, defined as the real bundle $F \to \dot{\Sigma}$
but with the sign of its complex structure reversed; we shall distinguish
this Cauchy-Riemann operator from $-\mathbf{D}^*$ by writing it as
$$
-\overline{\mathbf{D}}^* : \Gamma(\bar{F}) \to \Omega^{0,1}(\dot{\Sigma},\bar{F}),
$$
though it is technically the same operator.
Recall that the identity map defines a natural complex-antilinear isomorphism
between any complex vector bundle and its conjugate bundle; we shall denote
this isomorphism generally by
$$
E \to \bar{E} : v \mapsto \bar{v},
$$
so in particular it satisfies $\overline{c v} = \bar{c} \bar{v}$ for all
scalars $c \in \CC$, and similarly
$$
\overline{\mathbf{D}}^* \bar{\lambda} = \overline{\mathbf{D}^* \lambda}
$$
for $\lambda \in \Gamma(F)$.  The asymptotic operators for
$-\overline{\mathbf{D}}^*$ are
$$
\overline{\mathbf{A}}_z = -\mathbf{A}_z : \Gamma(\bar{E}_z) \to \Gamma(\bar{E}_z).
$$

\begin{lemma}
\label{lemma:conjTriv}
If $\tau$ is a choice of asymptotic trivialization on $E$ and
$\bar{\tau}$ denotes the \emph{conjugate} asymptotic trivialization\footnote{If
$\tau : E|_{\uU} \to \uU \times \CC^m$ is a local trivialization of $E$ with
$\tau(v) = (z,w)$, the conjugate trivialization $\bar{\tau} : \bar{E}|_{\uU}
\to \uU \times \CC^m$ is defined by $\bar{\tau}(\bar{v}) = (z,\bar{w})$.}, then
$$
c_1^{\bar{\tau}}(\bar{E}) = - c_1^\tau(E), \quad \text{ and } \quad
\muCZ^{\bar{\tau}}(\overline{\mathbf{A}}_z) = - \muCZ^\tau(\mathbf{A}_z)
\text{ for all $z \in \Gamma$}.
$$
\end{lemma}
\begin{proof}
Assuming $E$ is a line bundle, suppose $\eta$ is a generic section of
$E$ that matches a nonzero constant with respect to $\tau$ on the cylindrical
ends, so $c_1^\tau(E)$ is the signed count of zeroes of~$\eta$.  Then
$\bar{\eta} \in \Gamma(\bar{E})$ is similarly a nonzero constant on the ends
with respect to $\bar{\tau}$, but the signs of its zeroes are opposite
those of $\eta$ because they are defined as winding numbers with respect to
\emph{conjugate} local trivializations.  This proves $c_1^{\bar{\tau}}(\bar{E}) =
- c_1^\tau(E)$.  

The Conley-Zehnder indices can be computed from the formula
$$
\muCZ^\tau(\mathbf{A}_z) = \alpha_+^\tau(\mathbf{A}_z) + 
\alpha_-^\tau(\mathbf{A}_z),
$$
see Theorem~\ref{thm:CZwinding}.  Here $\alpha_-^\tau(\mathbf{A}_z)$ is the
largest possible winding number relative to $\tau$ of an eigenfunction
for $\mathbf{A}_z$ with negative eigenvalue, and $\alpha_+^\tau(\mathbf{A}_z)$
is the smallest possible winding number with positive eigenvalue.
The eigenfunctions of $\overline{\mathbf{A}}_z = - \mathbf{A}_z$ are the same,
but the signs of their eigenvalues are reversed, and the signs of their
winding numbers are also reversed because they must be measured relative
to the conjugate trivialization, thus
$$
\alpha_\pm^{\bar{\tau}}(\overline{\mathbf{A}}_z) = - 
\alpha_\mp^\tau(\mathbf{A}_z),
$$
implying
$$
\muCZ^{\bar{\tau}}(\overline{\mathbf{A}}_z) = 
\alpha_+^{\bar{\tau}}(\overline{\mathbf{A}}_z) +
\alpha_-^{\bar{\tau}}(\overline{\mathbf{A}}_z) =
- \alpha_-^\tau(\mathbf{A}_z) - \alpha_+^\tau(\mathbf{A}_z) =
-\muCZ^\tau(\mathbf{A}_z).
$$

The above calculations are all valid for line bundles, but the general case
follows by taking direct sums.
\end{proof}

We are now able to show that Theorem~\ref{thm:RiemannRoch} is consistent
with what we already know about the formal adjoint.

\begin{prop}
\label{prop:adjointIndex}
$I(-\overline{\mathbf{D}}^*) = - I(\mathbf{D})$.
\end{prop}
\begin{proof}
Under the isomorphism $F = \Lambda^{0,1}T^*\dot{\Sigma} \otimes E =
T\dot{\Sigma} \otimes E$, an asymptotic trivialization $\tau$ on $E$
induces an asymptotic trivialization $\p_s \otimes \tau$ on $F$,
where $\p_s$ denotes the asymptotic trivialization of $T\dot{\Sigma}$
defined via an outward pointing vector field on the cylindrical ends.
Counting zeroes of vector fields then proves
$c_1^{\p_s}(T\dot{\Sigma}) = \chi(\dot{\Sigma})$, so
$$
c_1^{\p_s \otimes \tau}(F) = c_1^{\p_s \otimes \tau}(T\dot{\Sigma} \otimes E) =
m c_1^{\p_s}(T\dot{\Sigma}) + c_1^\tau(E)
= m \chi(\dot{\Sigma}) + c_1^\tau(E).
$$
Applying Lemma~\ref{lemma:conjTriv} to the conjugate bundle then gives
$$
c_1^{\overline{\p_s \otimes \tau}}(\bar{F}) = - m \chi(\dot{\Sigma}) -
c_1^\tau(E).
$$
The unitary trivializations of the asymptotic bundles $\bar{E}_z$
corresponding to $\overline{\p_s \otimes \tau}$ are simply $\bar{\tau}$,
thus using Lemma~\ref{lemma:conjTriv} again for the Conley-Zehnder terms,
\begin{equation*}
\begin{split}
I(-\overline{\mathbf{D}}^*) &= m\chi(\dot{\Sigma}) + 2 c_1^{\overline{\p_s \otimes \tau}}(\bar{F})
+ \sum_{z \in \Gamma^+} \muCZ^{\bar{\tau}}(\overline{\mathbf{A}}_z) -
\sum_{z \in \Gamma^-} \muCZ^{\bar{\tau}}(\overline{\mathbf{A}}_z) \\
&= -m \chi(\dot{\Sigma}) - 2 c_1^\tau(E) 
- \sum_{z \in \Gamma^+} \muCZ^{\tau}(\mathbf{A}_z) +
\sum_{z \in \Gamma^-} \muCZ^{\tau}(\mathbf{A}_z) \\
&= - I(\mathbf{D}).
\end{split}
\end{equation*}
\end{proof}

We next consider the effect of an antilinear zeroth-order perturbation on
the formal adjoint.  By ``antilinear zeroth-order perturbation,'' we generally
mean a smooth section
$$
B \in \Gamma(\overline{\Hom}_\CC(E,F)).
$$
It is perhaps easier to understand $B$ in terms of the conjugate bundle
$\bar{E}$: indeed, there exists a unique
$$
\beta \in \Gamma(\Hom_\CC(\bar{E},F))
$$
such that
$$
B\eta = \beta \bar{\eta},
$$
and this correspondence defines a bundle isomorphism
$\overline{\Hom}_\CC(E,F) = \Hom_\CC(\bar{E},F)$.  

\begin{exercise}
\label{EX:conjugates}
Assume $X$ and $Y$ are complex vector bundles over the same base.
\begin{enumerate}
\renewcommand{\labelenumi}{(\alph{enumi})}
\item
Show that $\bar{X} \otimes \bar{Y}$ is canonically isomorphic to the
conjugate bundle of $X \otimes Y$.
\item
Show that $\Hom_\CC(\bar{X},\bar{Y})$ is canonically isomorphic to the
conjugate bundle of $\Hom_\CC(X,Y)$, and
$\overline{\Hom}_\CC(\bar{X},\bar{Y})$ is canonically isomorphic to the
conjugate bundle of $\overline{\Hom}_\CC(X,Y)$.
\item
Show that $\Lambda^{0,1}X := \overline{\Hom}_\CC(X,\CC)$ is canonically 
isomorphic to the conjugate bundle of $\Lambda^{1,0}X := \Hom_\CC(X,\CC)$.
\end{enumerate}
\end{exercise}

Define the Cauchy-Riemann type operator
$$
\mathbf{D}_B := \mathbf{D} +  B : \Gamma(E) \to \Gamma(F)
= \Omega^{0,1}(\dot{\Sigma},E),
$$
so $\mathbf{D}_B \eta = \mathbf{D}\eta + \beta \bar{\eta}$.  
To write down $\mathbf{D}_B^*$, observe that
since $\beta : \bar{E} \to F$ is a complex-linear bundle map between
Hermitian bundles, it has a complex-linear adjoint
$$
\beta^\dagger : F \to \bar{E} \quad \text{ such that } \quad
\langle \beta^\dagger \lambda, \bar{\eta} \rangle_{\bar{E}} = 
\langle \lambda , \beta \bar{\eta} \rangle_F \text{ for 
$\lambda \in F$, $\bar{\eta} \in \bar{E}$}.
$$
Here the bundle metric on $\bar{E}$ is defined by
$\langle \bar{\eta},\bar{\xi} \rangle_{\bar{E}} := 
\langle \xi , \eta \rangle_E$.  We then have
\begin{equation*}
\begin{split}
\Re \langle \lambda , B \eta \rangle_F &=
\Re \langle \lambda, \beta \bar{\eta} \rangle_F = 
\Re \langle \beta^\dagger \lambda , \bar{\eta} \rangle_{\bar{E}} =
\Re \langle \eta , \overline{\beta^\dagger \lambda} \rangle_E =
\Re \langle \overline{\beta^\dagger \lambda} , \eta \rangle_E \\
&= \Re \langle \overline{\beta^\dagger} \bar{\lambda} , \eta \rangle_E,
\end{split}
\end{equation*}
where $\overline{\beta^\dagger} \in \Gamma(\Hom_\CC(\bar{F},E))$ denotes the 
image of $\beta^\dagger \in \Gamma(\Hom_\CC(F,\bar{E}))$ under the 
complex-antilinear identity map from $\Hom_\CC(F,\bar{E})$ to its conjugate
bundle (see Exercise~\ref{EX:conjugates}).
The formal adjoint of $\mathbf{D}_B$ is thus
$$
\mathbf{D}_B^* = \mathbf{D}^* +  B^* : \Gamma(F) \to \Gamma(E),
$$
where $B^* : F \to E$ is defined by
$$
B^* \lambda := \overline{\beta^\dagger} \bar{\lambda}.
$$
To write down the resulting Cauchy-Riemann type operator on $\bar{F}$,
we replace $B^* : F \to E$ with $\overline{B}^* : \bar{F} \to \bar{E}$,
defined by
$$
\overline{B}^* \bar{\lambda} := \overline{B^*\lambda} = \beta^\dagger \lambda,
$$
giving a Cauchy-Riemann operator
$$
-\overline{\mathbf{D}}_B^* = -\overline{\mathbf{D}}^* + 
(-\overline{B}^*) : \Gamma(\bar{F}) \to \Gamma(\bar{E}) = 
\Omega^{0,1}(\dot{\Sigma},\bar{F}).
$$
The point of writing down this formula is to make the following observations:

\begin{lemma}
\label{lemma:adjointFun}
The zeroth-order perturbation $-\overline{B}^* : \bar{F} \to \bar{E}$ appearing
in $-\overline{\mathbf{D}}_B^*$ has the following properties:
\begin{enumerate}
\item $-\overline{B}^* : \bar{F} \to \bar{E}$ is complex antilinear;
\item There is a natural complex bundle isomorphism
$\overline{\Hom}_\CC(\bar{F},\bar{E}) = \Hom_\CC(F,\bar{E})$ that identifies
$-\overline{B}^*$ with $-\beta^\dagger$;
\item If $m=1$ and $B \in \Gamma(\overline{\Hom}_\CC(E,F))$ has only nondegenerate 
zeroes, then $-\overline{B}^* \in \Gamma(\overline{\Hom}_\CC(\bar{F},\bar{E}))$
has the same zeroes but with opposite signs.
\end{enumerate}
\end{lemma}
\begin{proof}
The first two statements follow immediately from the fact that 
$-\overline{B}^*$ is
the composition of the canonical conjugation map $\bar{F} \to F$ with the 
complex-linear bundle map $-\beta^\dagger : F \to \bar{E}$.  For the third,
it suffices to compare what $\beta \in \Gamma(\Hom_\CC(\bar{E},F))$ 
and $-\beta^\dagger : \Gamma(\Hom_\CC(F,\bar{E}))$ look like in local
trivializations near a zero: one is minus the complex conjugate of the other,
hence their zeroes count with opposite signs.
\end{proof}

\section{The index zero case on a torus}
\label{sec:torus}

As a warmup for the general case, we now fill in the details of Taubes's proof of
Theorem~\ref{thm:RiemannRoch} in the case 
$$
\dot{\Sigma} = \TT^2 := \CC \setminus (\ZZ \oplus i\ZZ)
$$
and $E = \TT^2 \times \CC$, i.e.~a trivial line bundle.  In this case
$I(\mathbf{D}) = \chi(\TT^2) + 2 c_1(E) = 0$, so our aim is to prove
$\ind(\mathbf{D}) = 0$.  What we will show in fact is that $\mathbf{D}$ is
homotopic through a continuous family of Fredholm operators to one that
is an isomorphism.  Since $E$ and $F$ are now both trivial, it will suffice
to consider the operator
$$
\mathbf{D} := \dbar = \p_s + i \p_t : H^1(\TT^2,\CC) \to L^2(\TT^2,\CC),
$$
whose formal adjoint is $\mathbf{D}^* := -\p = -\p_s + i\p_t$.
An antilinear zeroth-order perturbation is then equivalent to a choice of
function $\beta : \TT^2 \to \CC$, giving rise to a family of operators
$$
\mathbf{D}_\sigma \eta := \dbar \eta + \sigma \beta \bar{\eta}
$$
for $\sigma \in \RR$, where $\bar{\eta} : \TT^2 \to \CC$ now denotes the
straightforward complex conjugate of~$\eta$.  Let us assume that
$\beta : \TT^2 \to \CC$ is nowhere zero; note that this would not be
possible in more general situations, but is possible here because
$\Hom_\CC(\bar{E},F)$ is a trivial bundle.

\begin{lemma}
\label{lemma:T2}
$\mathbf{D}_\sigma$ is injective for all $\sigma > 0$ sufficiently large.
\end{lemma}
\begin{proof}
Elliptic regularity implies any $\eta \in \ker \mathbf{D}_\sigma$ is smooth,
so we shall restrict our attention to smooth functions $\eta : \TT^2 \to \CC$.
We start by comparing the two second-order differential operators
$$
\mathbf{D}^*\mathbf{D} \text{ and } \mathbf{D}_\sigma^* \mathbf{D}_\sigma 
: C^\infty(\TT^2,\CC) \to C^\infty(\TT^2,\CC).
$$
Both are nonnegative $L^2$-symmetric operators, and in fact the first
is simply the Laplacian
$$
\mathbf{D}^*\mathbf{D} = - \p \dbar = (-\p_s + i \p_t) (\p_s + i \p_t) =
-\p_s^2 - \p_t^2 = -\Delta.
$$
The formal adjoint of $\mathbf{D}_\sigma$ takes the form
$$
\mathbf{D}_\sigma^* \eta = \mathbf{D}^*\eta + \sigma B^*\eta =
\mathbf{D}^*\eta + \sigma \beta \bar{\eta},
$$
thus for any $\eta \in C^\infty(\TT^2,\CC)$,
\begin{equation}
\label{eqn:babyWeitzenbock}
\begin{split}
\mathbf{D}_\sigma^*\mathbf{D}_\sigma \eta &= (\mathbf{D}^* + \sigma B^*)
(\mathbf{D} + \sigma B)\eta \\
&= \mathbf{D}^*\mathbf{D} \eta + \sigma \left( \beta \overline{\dbar \eta} 
- \p (\beta \bar{\eta}) \right) + \sigma^2 B^*B \eta \\
&= \mathbf{D}^*\mathbf{D} \eta + \sigma \left( \beta \p \bar{\eta}
- (\p\beta) \bar{\eta} - \beta \p \bar{\eta} \right) + \sigma^2 B^*B \eta\\
&= \mathbf{D}^*\mathbf{D} \eta + \sigma^2 B^*B \eta - \sigma (\p\beta) \bar{\eta}.
\end{split}
\end{equation}
This is a \emph{Weitzenb\"ock formula}: its main message is that the 
Laplacian $\mathbf{D}^*\mathbf{D}$ and the related operator
$\mathbf{D}_\sigma^*\mathbf{D}_\sigma$ differ from each other only by
a zeroth-order term that will be positive definite if $\sigma$ is sufficiently
large.  Indeed, since $\beta$ is nowhere zero, we have
$|B\eta| \ge c |\eta|$ for some constant $c > 0$, thus
\begin{equation*}
\begin{split}
\| \mathbf{D}_\sigma \eta \|_{L^2}^2 &= \langle \eta , \mathbf{D}_\sigma^*
\mathbf{D}_\sigma \eta \rangle_{L^2} =
\langle \eta , \mathbf{D}^*\mathbf{D} \eta \rangle_{L^2}
+ \sigma^2 \langle \eta , B^*B\eta \rangle_{L^2} 
- \sigma \langle \eta , (\p\beta) \bar{\eta} \rangle_{L^2} \\
&= \| \mathbf{D}\eta \|_{L^2}^2 + \sigma^2 \| B\eta \|_{L^2}^2 -
\sigma \langle \eta , (\p\beta) \bar{\eta} \rangle_{L^2} \\
&\ge \left( \sigma^2 c^2 - \sigma \| \p\beta \|_{C^0} \right) \| \eta \|_{L^2}^2.
\end{split}
\end{equation*}
We conclude that as soon as $\sigma > 0$ is large enough to make the quantity
in parentheses positive, $\mathbf{D}_\sigma \eta$ cannot vanish unless
$\| \eta \|_{L^2} = 0$.
\end{proof}

\begin{proof}[Proof of Theorem~\ref{thm:RiemannRoch} for $E = \TT^2 \times \CC$]
The lemma above shows that one can add a large antilinear perturbation to
$\mathbf{D} = \dbar$ making the deformed operator $\mathbf{D}_\sigma$
injective.  By Lemma~\ref{lemma:adjointFun}, the same argument applies to
the formal adjoint $\mathbf{D}^*$, implying that for sufficiently large
$\sigma > 0$, $\mathbf{D}_\sigma^*$ is injective and thus $\mathbf{D}_\sigma$
is also surjective, and therefore an isomorphism.  This proves
$\ind(\mathbf{D}) = \ind(\mathbf{D}_\sigma) = 0$.
\end{proof}

Let's consider which particular details of the setup made the proof above 
possible.

First, the zeroth-order perturbation is complex antilinear.  We used this,
if only implicitly, in deriving the Weitzenb\"ock formula
\eqref{eqn:babyWeitzenbock}: the key step is in the third line, where
the two terms involving $\p\bar{\eta}$ cancel each other out and leave
nothing but zeroth-order terms remaining.  This would
not have happened if e.g.~$B : E \to F$ had been complex linear---we would
then have seen terms depending on the first derivative of $\eta$ in
$\mathbf{D}_\sigma^*\mathbf{D}_\sigma \eta - \mathbf{D}^*\mathbf{D} \eta$,
and this would have killed the whole argument.  The fact that this
cancellation happens when the perturbation is antilinear probably looks
like magic at this point, but there is a principle behind it; we will 
discuss it further in \S\ref{sec:Weitzenbock} below,
see Remark~\ref{remark:why}.

The second crucial fact we used was that $\beta : \TT^2 \to \CC$ is nowhere
zero, in order to obtain the lower bound on $\| B\eta \|_{L^2}$
in terms of $\| \eta \|_{L^2}$.  This cannot always be achieved---it is
possible in this special case only because $E$ and $F$ are both trivial
bundles and thus so is $\Hom_\CC(\bar{E},F)$.  On more general
bundles, the best we could hope for would be to pick 
$\beta \in \Gamma(\Hom_\CC(\bar{E},F))$ with finitely many zeroes, all
nondegenerate.  In this case the above argument fails, but it still tells
us something.  Suppose $\Sigma_\epsilon \subset \TT^2$ is a region disjoint from
the isolated zeroes of~$\beta$.  Then there exists a constant $c_\epsilon > 0$,
dependent on the region $\Sigma_\epsilon$, such that
$$
\| \beta\bar{\eta} \|_{L^2(\TT^2)}^2 \ge \| \beta\bar{\eta} 
\|_{L^2(\Sigma_\epsilon)}^2 \ge c_\epsilon \| \eta \|_{L^2(\Sigma_\epsilon)}^2,
$$
so instead of the estimate at the end of the proof above implying
$\mathbf{D}_\sigma$ is injective, we obtain one of the form
$$
\| \mathbf{D}_\sigma \eta \|_{L^2(\TT^2)}^2 \ge c_\epsilon \sigma^2 \| \eta \|_{L^2(\Sigma_\epsilon)}
- c \sigma \| \eta \|_{L^2(\TT^2)}^2.
$$
To see what this means, imagine we have sequences $\sigma_\nu \to \infty$ and
$\eta_\nu \in \ker \mathbf{D}_{\sigma_\nu}$, normalized so that
$\| \eta_\nu \|_{L^2} = 1$ for all~$\nu$.  The estimate above then implies
$$
\| \eta_\nu \|_{L^2(\Sigma_\epsilon)}^2 \le \frac{c}{c_\epsilon \sigma_\nu} \to 0
\quad \text{ as } \quad \nu \to \infty,
$$
so while all sections $\eta_\nu$ have the same amount of ``energy'' (as
measured via their $L^2$-norms), the energy is escaping from
$\Sigma_\epsilon$ as $\sigma_\nu$ increases.  This is true for \emph{any}
domain $\Sigma_\epsilon$ disjoint from the zeroes, so we
conclude that in the limit as $\sigma \to \infty$, sections in 
$\ker \mathbf{D}_\sigma$ have their energy concentrated in infinitesimally
small neighborhoods of the zeroes of~$\beta$.  We will see in the following
how to extract useful information from this concentration of energy.

\section{A Weitzenb\"ock formula for Cauchy-Riemann operators}
\label{sec:Weitzenbock}

The Weitzenb\"ock formula \eqref{eqn:babyWeitzenbock} can be generalized
to a useful relation between any two Cauchy-Riemann type operators that
differ by an \emph{antilinear} zeroth-order term.  To see this, we start with
a short digression on holomorphic and antiholomorphic vector bundles.

A smooth function
$f : \CC \supset \uU \to \CC$ is called \defin{antiholomorphic} if it 
satisfies
$(\p_s - i\p_t) f = 0$, which means its differential anticommutes with the
complex structure on~$\CC$.  The class of antiholomorphic functions is not
closed under composition, but it is closed under products, hence one can
define an \defin{antiholomorphic structure} on a complex vector bundle to be
a system of local trivializations for which all transition maps are
antiholomorphic.  Given the standard correspondence between holomorphic 
structures and Cauchy-Riemann type operators, it is easy to establish a 
similar correspondence between aniholomorphic structures and (complex-linear)
\defin{anti-Cauchy-Riemann type} operators, i.e.~those which satisfy
$$
\mathbf{D}(f\eta) = (\p f) \eta + f \mathbf{D}\eta
$$
for all $f \in C^\infty(\dot{\Sigma},\CC)$, where $\p f := df - i\, df \circ j
\in \Omega^{1,0}(\dot{\Sigma})$.  We've seen one important example of such
an operator already: if $\mathbf{D} : \Gamma(E) \to \Gamma(F)$ is complex linear, 
then $-\mathbf{D}^*$ is a complex-linear anti-Cauchy-Riemann operator on $F$
and thus endows $F$ with an antiholomorphic structure.
Another natural example occurs naturally on conjugate bundles: if
$E$ has a holomorphic structure, then $\bar{E}$ inherits from
this an antiholomorphic structure.  This is immediate from the fact that
$f : \CC \supset \uU \to \CC$ is holomorphic if and only if
$\bar{f} : \uU \to \CC$ is antiholomorphic.  If 
$\mathbf{D} : \Gamma(E) \to \Gamma(F) = \Omega^{0,1}(\dot{\Sigma},E)$ is the
corresponding complex-linear Cauchy-Riemann type operator on $E$, we shall
denote the resulting anti-Cauchy-Riemann operator by
$$
\overline{\mathbf{D}} : \Gamma(\bar{E}) \to \Gamma(\bar{F}) = 
\Omega^{1,0}(\dot{\Sigma},\bar{E}),
$$
where by definition $\overline{\mathbf{D}} \bar{\eta} = 
\overline{\mathbf{D}\eta}$.

\begin{exercise}
\label{EX:anti}
Show that if $X$ and $Y$ are antiholomorphic vector bundles over the same base,
then $X \otimes Y$ and $\Hom_\CC(X,Y)$ both naturally inherit antiholomorphic
bundle structures such that the obvious Leibniz rules are satisfied.
\textsl{Remark: the proof of this is exactly the same as for holomorphic
bundles, one only needs to change some signs.}
\end{exercise}
\begin{exercise}
\label{EX:real}
Suppose $X$ and $Y$ are complex vector bundles over the same base,
carrying real-linear anti-Cauchy-Riemann operators $\p_X$ and $\p_Y$
respectively.  Show that $H := \Hom_\RR(X,Y)$ then admits a real-linear
anti-Cauchy-Riemann operator $\p_H$ such that for all $\Phi \in \Gamma(H)$
and $\eta \in \Gamma(X)$,
$$
\p_Y (\Phi \eta) = (\p_H \Phi) \eta + \Phi (\p_X \eta).
$$
\textsl{Hint: write $\p_X$ and $\p_Y$ as complex-linear operators with
real-linear zeroth-order perturbations, and apply Exercise~\ref{EX:anti}.
Show moreover that any $C^k$-bounds satisfied by the zeroth-order terms
in $\p_X$ and $\p_Y$ are inherited by the zeroth-order term in~$\p_H$.}
\end{exercise}

The setup for the next result is as follows.  We assume again $m=1$, so
$E$ and $F$ are line bundles.
Fix $\beta \in \Gamma(\Hom_\CC(\bar{E},F))$, define $B \in 
\Gamma(\overline{\Hom}_\CC(E,F))$ by $B\eta := \beta\bar{\eta}$,
and use this to define the perturbed Cauchy-Riemann type operator
$$
\mathbf{D}_B := \mathbf{D} + B : \Gamma(E) \to \Gamma(F),
$$
whose formal adjoint is $\mathbf{D}_B^* = \mathbf{D}^* + B^*$ with
$B^* \lambda := \overline{\beta^\dagger} \bar{\lambda}$.

\begin{prop}
\label{prop:Weitzenbock}
The second-order differential
operators $\mathbf{D}^*\mathbf{D}$ and $\mathbf{D}_B^*\mathbf{D}_B$ on $E$
are related by
$$
\mathbf{D}_B^*\mathbf{D}_B \eta = \mathbf{D}^*\mathbf{D}^*\eta +
B^*B \eta - (\p_H \beta) \bar{\eta},
$$
where $\p_H$ is a real-linear anti-Cauchy-Riemann type operator on
$\Hom_\CC(\bar{E},F)$.  Moreover, if $\beta$ is $C^1$-bounded on~$\dot{\Sigma}$,
then $\p_H \beta$ is $C^0$-bounded.
\end{prop}
\begin{proof}
We have real-linear anti-Cauchy-Riemann operators $\overline{\mathbf{D}}$
and $-\mathbf{D}^*$ on $\bar{E}$ and $F$ respectively, so 
Exercise~\ref{EX:real} produces an operator $\p_H$ on $\Hom_\CC(\bar{E},F)$
for which the Leibniz rule is satisfied.  We can then write
\begin{equation*}
\begin{split}
\mathbf{D}_B^*\mathbf{D}_B \eta &= (\mathbf{D}^* + B^*)(\mathbf{D} + B) \eta \\
&= \mathbf{D}^*\mathbf{D} \eta + \overline{\beta^\dagger}\overline{\mathbf{D}\eta}
- (- \mathbf{D}^*)(\beta \bar{\eta}) + B^*B \eta \\
&= \mathbf{D}^*\mathbf{D} \eta + \overline{\beta^\dagger}\overline{\mathbf{D}}
\bar{\eta} - (\p_H\beta) \bar{\eta} - \beta \overline{\mathbf{D}}\bar{\eta}
+ B^*B \eta \\
&= \mathbf{D}^*\mathbf{D} \eta + B^*B\eta - (\p_H \beta) \bar{\eta} +
\left( \overline{\beta^\dagger} - \beta \right) \overline{\mathbf{D}}\bar{\eta}.
\end{split}
\end{equation*}
Here $\beta$ and $\overline{\beta^\dagger}$ are both viewed as complex-linear
bundle maps $\bar{F} \to E$, the latter in the obvious way, and the former
acting as $\1 \otimes \beta$ on $\bar{F} = \Lambda^{1,0}T^*\dot{\Sigma} \otimes
\bar{E}$ with target $\Lambda^{1,0}T^*\dot{\Sigma} \otimes F =
\Lambda^{1,0}T^*\dot{\Sigma} \otimes \Lambda^{0,1}T^*\dot{\Sigma} \otimes E = E$.
Choosing unitary local trivializations, $\beta$ and $\overline{\beta^\dagger}$
are represented by the same complex-valued function: indeed, the latter is
the transpose of the former as $m$-by-$m$ complex matrices, but since $m=1$,
this means they are identical.

Finally, we observe that the asymptotic convergence conditions satisfied by
$\mathbf{D}$ on the cylindrical ends imply similar conditions for
all other Cauchy-Riemann and anti-Cauchy-Riemann operators in this picture,
yielding an estimate of the form $\|\p_H\beta\|_{C^0} \le c 
\|\beta\|_{C^1}$ globally on~$\dot{\Sigma}$.
\end{proof}

\begin{remark}
The above proof used the assumption $m=1$ in order to
conclude $\overline{\beta^\dagger} - \beta \equiv 0$.  For higher rank
bundles, this imposes a nontrivial condition that must be satisfied in
order for the Weitzenb\"ock formula to hold, cf.~\cite{GerigWendl}.
\end{remark}

\begin{remark}
\label{remark:why}
We can now pick out a geometric reason for the miraculous cancellation in
the Weitzenb\"ock formula: the perturbation $B$ is described by a complex
bundle map $\bar{E} \to F$, where $\bar{E}$ and $F$ both have natural
antiholomorphic bundle structures defined via the complex-linear parts of
$\overline{\mathbf{D}}$ and $-\mathbf{D}^*$ respectively.  A complex-linear
perturbation $B : E \to F$ would not work because $E$ is holomorphic rather
than antiholomorphic: while $\overline{\mathbf{D}}$ can be fit into the
same Leibniz rule with $-\mathbf{D}^*$, the same is not true
of~$\mathbf{D}$.
\end{remark}

\section{Large antilinear perturbations and energy concentration}
\label{sec:concentration}

We continue in the setting of Proposition~\ref{prop:Weitzenbock} and set
$$
\mathbf{D}_\sigma := \mathbf{D} + \sigma B : \Gamma(E) \to \Gamma(F)
$$
for $\sigma > 0$.  After a compact perturbation of $\mathbf{D}$, we can
without loss of generality also impose the following
assumptions on $\mathbf{D}$, $\beta \in \Gamma(\Hom_\CC(\bar{E},F))$ and the 
area form~$\vol$:
\begin{enumerate}
\renewcommand{\labelenumi}{(\roman{enumi})}
\item All zeroes of $\beta$ are nondegenerate.
\item Both $|\beta|$ and $1 / |\beta|$ are bounded outside of a compact subset
of~$\dot{\Sigma}$.
\item Near each point $\zeta \in \dot{\Sigma}$ with $\beta(\zeta) = 0$, there 
exists a neighborhood $\dD(\zeta) \subset \dot{\Sigma}$ of~$\zeta$, a
holomorphic coordinate chart identifying $(\dD(\zeta),j,\zeta)$ with the
unit disk $(\DD,i,0)$, and a local trivialization of $E$ over $\dD(\zeta)$
that identifies $\mathbf{D}$ with $\dbar = \p_s + i\p_t : C^\infty(\DD,\CC)
\to C^\infty(\DD,\CC)$ and $\beta$ with one of the functions
$$
\beta(z) = z \quad\text{ or }\quad \beta(z) = \bar{z},
$$
the former if $\zeta$ is a positive zero and the latter if it is negative.
\item In the holomorphic coordinate on $\dD(\zeta)$ described above,
$\vol$ is the standard Lebesgue measure.
\end{enumerate}

As in the torus case discussed in \S\ref{sec:torus}, we will see that the
Weitzenb\"ock formula implies a concentration of energy near the
zeroes of $\beta$ for sections $\eta \in \ker \mathbf{D}_\sigma$ as
$\sigma \to \infty$.  To understand what really happens in this limit,
we will use a rescaling trick.  Denote the zero set of $\beta$ by
$$
Z(\beta) = Z^+(\beta) \cup Z^-(\beta) \subset \dot{\Sigma},
$$
partitioned into the positive and negative zeroes.  For any
$\eta \in \Gamma(E)$, $\zeta \in Z^\pm(\beta)$ and $\sigma > 0$, we then
define a rescaled function
$$
\eta^{(\zeta,\sigma)} : \DD_{\sqrt{\sigma}} \to \CC : z \mapsto
\frac{1}{\sqrt{\sigma}} \eta(z / \sqrt{\sigma}),
$$
where the right hand side denotes the local representation of $\eta$ 
on $\dD(\zeta)$ in the chosen coordinate and trivialization.
Notice that the equation $\mathbf{D}_\sigma \eta = 0$ appears in this
local representation as either $\dbar \eta + \sigma z \bar{\eta} = 0$ or
$\dbar \eta + \sigma \bar{z} \bar{\eta} = 0$ depending on the sign of~$\zeta$,
and the function $f := \eta^{(\zeta,\sigma)}$ then satisfies
$$
\dbar f+ z \bar{f} = 0 \quad\text{ or } \quad
\dbar f + \bar{z} \bar{f} = 0
\qquad \text{ on $\DD_{\sqrt{\sigma}}$}.
$$
We will take a closer look at these two PDEs in \S\ref{sec:plane} below.
But first, observe that by change of variables,
$$
\left\| \eta^{(\zeta,\sigma)} \right\|_{L^2(\DD_{\sqrt{\sigma}})} =
\| \eta \|_{L^2(\dD(\zeta))}.
$$

\begin{lemma}
\label{lemma:concentration}
Assume $\sigma_\nu \to \infty$, and $\eta_\nu \in \ker \mathbf{D}_{\sigma_\nu}$ is
a sequence satisfying a uniform $L^2$-bound.  Then after passing to
a subsequence, the rescaled functions $\eta_\nu^\zeta := 
\eta_\nu^{(\zeta,\sigma_\nu)} : \DD_{\sqrt{\sigma_\nu}} \to \CC$ for each
$\zeta \in Z^\pm(\beta)$ converge in $C^\infty_{\loc}(\CC)$ to smooth functions
$\eta_\infty^\zeta \in L^2(\CC)$ satisfying
\begin{equation*}
\begin{split}
\dbar \eta_\infty^\zeta + z \overline{\eta_\infty^\zeta} =0\quad & 
\text{ if $\zeta \in Z^+(\beta)$},\\
\dbar \eta_\infty^\zeta + \bar{z} \overline{\eta_\infty^\zeta} =0\quad &
\text{ if $\zeta \in Z^-(\beta)$}.
\end{split}
\end{equation*}
Moreover, if $\xi_\nu \in \ker \mathbf{D}_{\sigma_\nu}$ is another sequence
with these same properties and convergence $\xi_\nu^\zeta \to \xi_\infty^\zeta$,
then
$$
\lim_{\nu \to \infty} \langle \eta_\nu , \xi_\nu \rangle_{L^2(E)} =
\sum_{\zeta \in Z(\beta)} \langle \eta_\infty^\zeta , \xi_\infty^\zeta \rangle_{L^2(\CC)}.
$$
\end{lemma}
\begin{proof}
The uniform $L^2$-bound implies uniform bounds on 
$\|\eta_\nu^\zeta\|_{L^2(\DD_R)}$ for every $R > 0$, where $\nu$ here is
assumed sufficiently large so that $R < \sqrt{\sigma_\nu}$.  Since
$\eta_\nu^\zeta$ satisfies a Cauchy-Riemann type equation on $\DD_R$, the
usual elliptic estimates (see Lecture~\ref{lec:local}) then imply uniform
$H^k$-bounds for every $k \in \NN$ on every compact subset in the interior
of~$\DD_R$, hence $\eta_\nu^\zeta$ has a $C^\infty_{\loc}$-convergent 
subsequence on~$\CC$, and the limit $\eta_\infty^\zeta$ clearly satisfies
the stated PDE.  The uniform $L^2$-bound also implies a uniform bound on
$\| \eta_\nu^\zeta \|_{L^2(\DD_{\sqrt{\sigma_\nu}})}$ and thus an
$R$-independent uniform bound on $\| \eta_\nu^\zeta \|_{L^2(\DD_R)}$
as $\nu \to \infty$, implying that $\eta_\infty^\zeta$ is in~$L^2(\CC)$.

The limit of $\langle \eta_\nu,\xi_\nu \rangle_{L^2(E)}$ is now
proved using the Weitzenb\"ock formula.  Let
$$
\dot{\Sigma}_\epsilon := \dot{\Sigma} \setminus \bigcup_{\zeta \in Z(\beta)}
\dD(\zeta),
$$
so there exists a constant $c > 0$ such that $\beta$ satisfies 
$|\beta(z) \bar{v}| \ge c |v|$ for all $v \in E_z$, 
$z \in \dot{\Sigma}_\epsilon$.  (Note that this depends on the assumption
of $1 / |\beta|$ being bounded outside of a compact subset.)  Now
by Proposition~\ref{prop:Weitzenbock},
\begin{equation*}
\begin{split}
0 &= \| \mathbf{D}_{\sigma_\nu} \eta_\nu \|_{L^2(\dot{\Sigma})}^2 = 
\langle \eta_\nu , \mathbf{D}_{\sigma_\nu}^*\mathbf{D}_{\sigma_\nu} \eta_\nu \rangle_{L^2(\dot{\Sigma})} \\
&= \langle \eta_\nu , \mathbf{D}^*\mathbf{D} \eta_\nu \rangle_{L^2(\dot{\Sigma})}
+ \sigma_\nu^2 \langle \eta_\nu , B^*B \eta_\nu \rangle_{L^2(\dot{\Sigma})} -
\sigma_\nu \langle \eta_\nu , (\p_H \beta) \bar{\eta}_\nu \rangle_{L^2(\dot{\Sigma})} \\
&\ge \| \mathbf{D} \eta_\nu \|_{L^2(\dot{\Sigma})}^2 + 
\sigma_\nu^2 c^2 \| \eta_\nu \|_{L^2(\dot{\Sigma}_\epsilon)}^2 
- \sigma_\nu c' \| \eta_\nu \|_{L^2(\dot{\Sigma})}^2 \\
&\ge \sigma_\nu^2 c^2 \| \eta_\nu \|_{L^2(\dot{\Sigma}_\epsilon)}^2 
- \sigma_\nu c' \| \eta_\nu \|_{L^2(\dot{\Sigma})}^2
\end{split}
\end{equation*}
for some constant $c' > 0$ independent of~$\nu$.  This implies
$$
\| \eta_\nu \|_{L^2(\dot{\Sigma}_\epsilon)}^2 \le \frac{c'}{c^2 \sigma_\nu} 
\| \eta_\nu \|_{L^2(\dot{\Sigma})}^2 \to 0 \quad \text{ as } \quad
\nu \to \infty
$$
since $\| \eta_\nu \|_{L^2(\dot{\Sigma})}$ is uniformly bounded.
The same estimate applies to $\xi_\nu$, so that
$\langle \eta_\nu , \xi_\nu \rangle_{L^2(\dot{\Sigma}_\epsilon)} \to 0$
and thus by change of variables,
\begin{equation*}
\begin{split}
\lim_{\nu \to \infty} \langle \eta_\nu,\xi_\nu \rangle_{L^2(\dot{\Sigma})} &=
\lim_{\nu \to \infty} \sum_{\zeta \in Z(\beta)}
\langle \eta_\nu,\xi_\nu \rangle_{L^2(\dD(\zeta))} 
= \lim_{\nu \to \infty} \sum_{\zeta \in Z(\beta)}
\langle \eta_\nu^\zeta,\xi_\nu^\zeta \rangle_{L^2(\DD_{\sqrt{\sigma_\nu}})} \\
&= \sum_{\zeta \in Z(\beta)} \langle \eta_\infty^\zeta,\xi_\infty^\zeta \rangle_{L^2(\CC)}.
\end{split}
\end{equation*}
\end{proof}

\section{Two Cauchy-Riemann type problems on the plane}
\label{sec:plane}

The rescaling trick in the previous section produced smooth solutions
$f : \CC \to \CC$ of class $L^2(\CC)$ to the two equations
$$
\dbar f + z \bar{f} = 0, \qquad \dbar f + \bar{z} \bar{f} = 0.
$$
It turns out that we can say precisely what all such solutions are.
Write $\mathbf{D}_+ f := \dbar f + z \bar{f}$ and
$\mathbf{D}_- f := \dbar f + \bar{z} \bar{f}$.  Both operators
differ from $\dbar$ by antilinear perturbations, so they satisfy
Weitzenb\"ock formulas relating $\mathbf{D}^*_\pm \mathbf{D}_\pm$ to the
Laplacian $-\Delta = \dbar^*\dbar = -\p_s^2 - \p_t^2$.  Indeed, repeating
Proposition~\ref{prop:Weitzenbock} in these special cases gives
$$
\mathbf{D}^*_+ \mathbf{D}_+ f = -\Delta f + |z|^2 f - 2 \bar{f}
\qquad\text{ and }\qquad
\mathbf{D}^*_- \mathbf{D}_- f = -\Delta f + |z|^2 f.
$$
To make use of this, recall that a smooth function $u : \uU \to \RR$ on an 
open subset $\uU \subset \CC$ is called \defin{subharmonic} if it satisfies
$$
- \Delta u \le 0.
$$
Subharmonic functions satisfy a \defin{mean value property}:
$$
-\Delta u \le 0 \text{ on $\uU$} \qquad \Rightarrow \qquad
u(z_0) \le \frac{1}{\pi r^2} \int_{\DD_r(z_0)} u(z) \, d\mu(z)
\quad \text{ for all } \quad \DD_r(z_0) \subset \uU,
$$
where $\DD_r(z_0) \subset \CC$ denotes the disk of radius $r > 0$ about a point
$z_0 \in \uU$, and $d\mu(z)$ is the Lebesgue measure on~$\CC$;
see e.g.~\cite{Evans}*{p.~85}.

\begin{exercise}
\label{EX:Kato}
Show that for any smooth complex-valued function $f$ on an open subset of~$\CC$,
$$
\Delta |f|^2 = 2 \Re \langle f , \Delta f \rangle + 2 |\nabla f|^2,
$$
where $\langle\ ,\ \rangle$ denotes the standard Hermitian inner product
on $\CC$ and $|\nabla f|^2 := |\p_s f|^2 + |\p_t f|^2$.
\end{exercise}

\begin{prop}
\label{prop:D-}
The equation $\dbar f + \bar{z} \bar{f}=0$ does not admit any nontrivial 
smooth solutions $f \in L^2(\CC,\CC)$.
\end{prop}
\begin{proof}
If $f : \CC \to \CC$ is smooth with $\mathbf{D}_- f = 0$, then the 
Weitzenb\"ock formula for $\mathbf{D}_-$ implies $\Delta f = |z|^2 f$.
Then by Exercise~\ref{EX:Kato},
$$
\Delta |f|^2 = 2 \Re \langle f , |z|^2 f \rangle + 2 |\nabla f|^2 =
2 |z|^2 |f|^2 + 2 |\nabla f|^2,
$$
implying that $|f|^2 : \CC \to \RR$ is subharmonic.  Now if $f(z_0) \ne 0$
for some $z_0 \in \CC$, the mean value property implies
$$
\int_{\DD_r(z_0)} |f(z)|^2 \, d\mu(z) \ge \pi r^2 |f(z_0)|^2 \to \infty
\quad \text{ as }\quad r \to \infty,
$$
so $f \not\in L^2(\CC)$.
\end{proof}

\begin{prop}
\label{prop:D+}
Every smooth solution $f \in L^2(\CC,\CC)$ to the equation
$\dbar f + z \bar{f} = 0$ is a constant real multiple of
$f_0(z) := e^{-\frac{1}{2} |z|^2}$.
\end{prop}
\begin{proof}
We claim first that every smooth solution in $L^2(\CC,\CC)$ of
$\mathbf{D}_+ f = 0$ is purely real valued.  The Weitzenb\"ock formula for
this case gives $\Delta f = |z|^2 f - 2\bar{f}$, and taking the difference
between this equation and its complex conjugate then implies that
$u := \Im f : \CC \to \RR$ satisfies
$$
\Delta u = (|z|^2 + 2) u.
$$
Now by Exercise~\ref{EX:Kato},
$$
\Delta(u^2) = 2 |\nabla u|^2 + 2 (|z|^2 + 2) u^2 \ge 0,
$$
so $u^2 : \CC \to \RR$ is subharmonic, and the mean value property implies
as in the proof of Prop.~\ref{prop:D-} that $u \not\in L^2(\CC)$ and hence
$f \not\in L^2(\CC)$ unless $u \equiv 0$.  This proves the claim.

It is easy to check however that $f_0$ is a solution and is in~$L^2(\CC)$.
Since it is also nowhere zero, every other solution $f$ must then take
the form $f(z) = v(z) f_0(z)$ for some \emph{real-valued} function
$v : \CC \to \RR$.  Since $\mathbf{D}_+$ is a Cauchy-Riemann type operator,
the Leibniz rule then implies $\dbar v \equiv 0$.  But the only globally
holomorphic functions with trivial imaginary parts are constant.
\end{proof}

\section{A linear gluing argument}

Now we're getting somewhere.

\begin{lemma}
\label{lemma:largePerturbation}
Suppose the assumptions of \S\ref{sec:concentration} hold and
$\beta \in \Gamma(\Hom_\CC(\bar{E},F))$ has $I_+ \ge 0$ positive and 
$I_- \ge 0$ negative zeroes.  Then for all $\sigma > 0$ sufficiently large,
$$
\dim \ker \mathbf{D}_\sigma \le I_+ \quad\text{ and }\quad
\dim \coker \mathbf{D}_\sigma \le I_-.
$$
In particular, for sufficiently large~$\sigma$, $\mathbf{D}_\sigma$ is
injective if all zeroes of $\beta$ are negative and surjective if
all zeroes are positive.
\end{lemma}
\begin{proof}
Arguing by contradiction, suppose there exists a sequence
$\sigma_\nu \to \infty$ such that $\dim \ker \mathbf{D}_{\sigma_\nu} > I_+$,
and pick $(I_+ + 1)$ sequences of sections 
$\eta_\nu^1,\ldots,\eta_\nu^{I_+ + 1} \in \ker \mathbf{D}_{\sigma_\nu}$
which form $L^2$-orthonormal sets for each~$\nu$.  By 
Lemma~\ref{lemma:concentration}, we can then extract a subsequence such that
rescaling near the zeroes of $\beta$ produces $C^\infty_\loc$-convergent
sequences whose limits form an $(I_+ + 1)$-dimensional orthonormal set in
$$
\bigoplus_{\zeta \in Z(\beta)} L^2(\CC,\CC),
$$
where the component functions $f \in L^2(\CC,\CC)$ for $\zeta \in Z^+(\zeta)$ satisfy
$\dbar f + z \bar{f} = 0$, while those for $\zeta \in Z^-(\zeta)$ satisfy
$\dbar f + \bar{z} \bar{f} = 0$.  Proposition~\ref{prop:D-} now implies
that the component functions for $\zeta \in Z^-(\zeta)$ are all trivial,
and by Proposition~\ref{prop:D+}, the components for
$\zeta \in Z^+(\zeta)$ belong to $1$-dimensional subspaces
$\ker \mathbf{D}_+ \subset L^2(\CC)$ generated by the function
$e^{-\frac{1}{2}|z|^2}$.
We conclude that the limiting orthonormal set lives in a precisely
$I_+$-dimensional subspace
$$
\bigoplus_{\zeta \in Z^+(\beta)} \ker \mathbf{D}_+ \subset
\bigoplus_{\zeta \in Z(\beta)} L^2(\CC,\CC),
$$
and this is a contradiction since there are $I_+ + 1$ elements in the set.

Applying the same argument to the formal adjoint implies similarly
$\dim \ker \mathbf{D}_\sigma^* \le I_-$ for $\sigma$ sufficiently large.
\end{proof}

We would next like to turn the two inequalities in the above lemma into
equalities, which means showing that the $I_+$-dimensional subspace of
$\bigoplus_{\zeta \in Z^+(\beta)} L^2(\CC,\CC)$ generated by solutions
of $\dbar f + z\bar{f} = 0$ is isomorphic to $\ker \mathbf{D}_\sigma$ for
$\sigma$ sufficiently large.  This requires a simple example of a \emph{linear
gluing} argument, the point of which is to reverse the ``convergence after
rescaling'' process that we saw in Lemma~\ref{lemma:concentration}.
The first step is a \defin{pregluing} construction which turns elements
of $\bigoplus_{\zeta \in Z^+(\beta)} \ker \mathbf{D}_+$ into 
\emph{approximate} solutions to $\mathbf{D}_\sigma \eta = 0$
for large~$\sigma$.  To this end, fix a smooth bump function
$$
\rho \in C_0^\infty(\intDD,[0,1]), \qquad \rho|_{\DD_{1/2}} \equiv 1
$$
and define for each $\zeta \in Z^+(\beta)$ and $\sigma > 0$ a linear map
$$
\Phi_\sigma^\zeta : \ker \mathbf{D}_+ \to \Gamma(E)
$$
such that $\Phi_\sigma^\zeta(f)$ is a section with support in $\dD(\zeta)$
whose expression in our fixed coordinate and trivialization on that 
neighborhood is the function
$$
f_\sigma^\zeta(z) = \rho(z) \sqrt{\sigma} f(\sqrt{\sigma} z).
$$
Adding up the $\Phi_\sigma^\zeta$ for all $\zeta \in Z^+(\beta)$ then produces
a linear map
$$
\Phi_\sigma : \bigoplus_{\zeta \in Z^+(\beta)} \ker \mathbf{D}_+ \to
\Gamma(E)
$$
whose image consists of sections supported near $Z^+(\beta)$, each a linear
combination of cut-off Gaussians with energy concentrated in smaller
neighborhoods of $Z^+(\beta)$ for larger~$\sigma$.  These sections are
manifestly not in $\ker \mathbf{D}_\sigma$ since they vanish on open subsets
and thus violate unique continuation, but they are close, in a quantitative
sense:

\begin{lemma}
\label{lemma:PhiBounds}
For each $\sigma > 0$, there exists a constant $c_\sigma > 0$ such that
$$
\| \mathbf{D}_\sigma \Phi_\sigma(f) \|_{L^2} \le c_\sigma \| f \|_{L^2}
\quad \text{ for all } \quad
f \in \bigoplus_{\zeta \in Z^+(\beta)} \ker \mathbf{D}_+, 
$$
and $c_\sigma \to 0$ as $\sigma \to \infty$.
Moreover, for every pair $f,g \in \bigoplus_{\zeta \in Z^+(\beta)}
\ker \mathbf{D}_+$,
$$
\langle \Phi_\sigma(f) , \Phi_\sigma(g) \rangle_{L^2} \to \langle f,g \rangle_{L^2}
$$
as $\sigma \to \infty$.
\end{lemma}
\begin{proof}
First, observe that any $f \in \bigoplus_{\zeta \in Z^+(\beta)} \ker \mathbf{D}_+$ 
is described by a collection of functions 
$\{ f_\zeta \in L^2(\CC) \}_{\zeta \in \beta^+(Z)}$ which take the form
$$
f_\zeta(z) = K_\zeta e^{-\frac{1}{2}|z|^2},
$$
for some constants $K_\zeta \in \RR$.  Since each $f_\zeta$ is in
$\ker \mathbf{D}_+$, we find
\begin{equation}
\begin{split}
\mathbf{D}_\sigma\left(\Phi_\sigma(f)\vert_{\dD(\zeta)}\right)(z)&=\overline{\partial}\rho(z)\sqrt{\sigma}f_\zeta(\sqrt{\sigma}z)+\rho(z)\sigma\overline{\partial}f_\zeta(\sqrt{\sigma}z)\\
& \quad +\sigma z\rho(z)\sqrt{\sigma}f_\zeta(\sqrt{\sigma}z)\\
&=\overline{\partial}\rho(z)\sqrt{\sigma}f_\zeta(\sqrt{\sigma}z)+\rho(z)\sigma(\mathbf{D}_+f_\zeta)(\sqrt{\sigma}z)\\
&=\overline{\partial}\rho(z)\sqrt{\sigma} K_\zeta e^{-\frac{1}{2} \sigma |z|^2}.
\end{split}
\end{equation}
Now since $\overline{\partial}\rho=0$ in $\DD_{1/2}$, we obtain
\begin{equation*}
\begin{split}
\|\mathbf{D}_\sigma\Phi_\sigma(f)\|^2_{L^2} &= 
\sum_{\zeta \in Z^+(\beta)} \int_{\dD(\zeta)}|\mathbf{D}_\sigma\Phi_\sigma(f)(z)|^2 \, d\mu(z)\\
&=\sum_{\zeta \in Z^+(\beta)} \int_{\DD \setminus \DD_{1/2}}|\overline{\partial}\rho(z)|^2\sigma K_\zeta^2 e^{-\sigma |z|^2}\, d\mu(z) \\
& \leq I \sigma e^{-\sigma/4}\sum_{\zeta \in Z^+(\beta)} K_\zeta^2, 
\end{split}
\end{equation*}
where we abbreviate $I := \int_{\DD \setminus \DD_{1/2}} \left| \dbar\rho(z) \right|^2 \, d\mu(z)$.
The norm of $f$ is given by 
$$
\|f\|^2_{L^2}= \sum_{\zeta \in Z^+(\beta)} \int_{\CC} K_\zeta^2 e^{-|z|^2}\, d\mu(z) 
= \left( \int_{\CC} e^{-|z|^2}\, d\mu(z) \right) \sum_{\zeta \in Z^+(\beta)} K_\zeta^2.
$$
We conclude that there is a bound of the form
$$
\|\mathbf{D}_\sigma\Phi_\sigma(f)\|_{L^2}\leq C \sqrt{\sigma} e^{-\sigma/2} \| f \|_{L^2},
$$
which proves the first statement since $\sqrt{\sigma} e^{-\sigma/2} \to 0$
as $\sigma \to \infty$.

The second statement follows by a change of variable, since
\begin{equation*}
\begin{split}
 \langle \Phi_\sigma(f) , \Phi_\sigma(g) \rangle_{L^2} &= \sum_{\zeta\in Z^+(\beta)} \langle \Phi_\sigma(f)\vert_{\dD(\zeta)} , \Phi_\sigma(g)\vert_{\dD(\zeta)} \rangle_{L^2(\dD(\zeta))}\\&=
\sum_{\zeta\in Z^+(\beta)}\int_{\DD}\rho^2(z)\sigma f_\zeta(\sqrt{\sigma}z) g_\zeta(\sqrt{\sigma}z) \, d\mu(z)\\&= \sum_{\zeta\in Z^+(\beta)}\int_{\DD_{\sqrt{\sigma}}}\rho^2\left(\frac{z}{\sqrt{\sigma}}\right)f_\zeta(z)g_\zeta(z) \, d\mu(z)
\end{split}
\end{equation*}
The functions $f_\zeta$ and $g_\zeta$ are both real multiples of
$e^{-\frac{1}{2}|z|^2}$, so
this last integral for each $\zeta \in Z^+(\beta)$ is bounded between
$\int_{\DD_{\sqrt{\sigma}/2}} f_\zeta(z) g_\zeta(z) \, d\mu(z)$ and 
$\int_{\DD_{\sqrt{\sigma}}} f_\zeta(z) g_\zeta(z) \, d\mu(z)$, both of which
converge to $\int_\CC f_\zeta(z) g_\zeta(z) \, d\mu(z)$ as $\sigma \to \infty$,
thus
$$
\lim_{\sigma \rightarrow \infty} \langle \Phi_\sigma(f) , \Phi_\sigma(g) \rangle_{L^2}=\langle f,g\rangle_{L^2}.
$$
\end{proof}

To turn approximate solutions into actual solutions, let
$$
\Pi_\sigma : L^2(E) \to \ker \mathbf{D}_\sigma
$$
denote the orthogonal projection.  We will prove:

\begin{prop}
\label{prop:inj}
If all zeroes of $\beta$ are positive, then the linear map 
$$
\Pi_\sigma \circ \Phi_\sigma : \bigoplus_{\zeta \in Z^+(\beta)} \ker \mathbf{D}_+ \to
\ker \mathbf{D}_\sigma
$$
is injective for all $\sigma > 0$ sufficiently large.
\end{prop}

This statement says in effect that whenever $\sigma > 0$ is large enough and
$\eta := \Phi_\sigma(f) \in \Gamma(E)$ is in the image of the pregluing map,
with $f$ normalized by $\| f \|_{L^2} = 1$, we can find a ``correction''
$\xi \in (\ker \mathbf{D}_\sigma)^\perp$ such that
$$
\eta + \xi \ne 0 \quad\text{ but }\quad \mathbf{D}_\sigma(\eta + \xi) = 0.
$$
An element $\xi \in (\ker \mathbf{D}_\sigma)^\perp$ with the second
property certainly exists, and in fact it's unique: indeed, the assumption 
$Z^-(\beta) = \emptyset$ implies via Lemma~\ref{lemma:largePerturbation} that
$\mathbf{D}_\sigma$ is surjective and thus restricts to an isomorphism
from $(\ker \mathbf{D})^\perp \cap H^1(E)$ to~$L^2(F)$, with a bounded
right inverse
$$
\mathbf{Q}_\sigma : L^2(F) \to H^1(E) \cap (\ker \mathbf{D})^\perp,
$$
hence $\xi := -\mathbf{Q}_\sigma(\mathbf{D}_\sigma \eta)$.
We know moreover from Lemma~\ref{lemma:PhiBounds} that $\| \eta \|_{L^2}$ is 
close to $\| f \|_{L^2} = 1$, so to prove $\eta + \xi \ne 0$, it would suffice
to show $\| \xi \|_{L^2}$ is small, which sounds likely since we also know
$\|\mathbf{D}_\sigma \eta\|_{L^2}$ is small and $\mathbf{Q}_\sigma$ is a
bounded operator.  To make this reasoning precise, we just need to have some 
control over $\| \mathbf{Q}_\sigma \|$ as $\sigma \to \infty$, or equivalently,
a quantitative measure of the injectivity of 
$\mathbf{D}_\sigma|_{(\ker \mathbf{D}_\sigma)^\perp \cap H^1(E)}$.
This requires one last appeal to the Weitzenb\"ock formula.

\begin{lemma}
\label{lemma:injBound}
Assume all zeroes of $\beta$ are positive.  Then
there exist constants $c > 0$ and $\sigma_0$ such that for all
$\sigma > \sigma_0$,
$$
\| \eta \|_{L^2} \le c \| \mathbf{D}_\sigma \eta \|_{L^2} 
\quad \text{ for all }\quad
\eta \in H^1(E) \cap (\ker \mathbf{D}_\sigma)^\perp.
$$
\end{lemma}
\begin{proof}
Let us instead prove that if zeroes of $\beta$ are all \emph{negative},
then the same bound holds for all $\eta \in H^1(E)$.  The stated
result follows from this by considering the formal adjoint and using
Exercise~\ref{EX:adjointBound} below.  Note that by density, it suffices
to prove the estimate holds for all $\eta \in C_0^\infty(E)$.

Assume therefore that $Z^+(\beta) = \emptyset$ and, arguing by contradiction,
suppose there exist sequences $\sigma_\nu \to \infty$ and 
$\eta_\nu \in C_0^\infty(E)$ with $\| \eta_\nu \|_{L^2} = 1$ and
$$
\| \mathbf{D}_{\sigma_\nu} \eta_\nu \|_{L^2} \to 0.
$$
The usual rescaling trick and application of the Weitzenb\"ock formula
then produces for each $\zeta \in Z^-(\beta)$ a sequence of functions
$\eta_\nu^\zeta := \eta_\nu^{(\zeta,\sigma_\nu)} : \DD_{\sqrt{\sigma_\nu}} \to \CC$
which satisfy
$$
\sum_{\zeta \in Z^-(\beta)} \| \eta_\nu^\zeta \|_{L^2(\DD_{\sqrt{\sigma_\nu}})}^2 \to 1 \quad \text{ and }\quad
\| \mathbf{D}_- \eta_\nu^\zeta \|_{L^2(\DD_{\sqrt{\sigma_\nu}})} \to 0
$$
as $\nu \to \infty$.
Indeed, defining $\dot{\Sigma}_\epsilon$ as in the proof of 
Lemma~\ref{lemma:concentration}, a similar application of the Weitzenb\"ock 
formula yields
$$
\|\mathbf{D}_{\sigma_\nu}\eta_\nu\|^2_{L^2(\dot{\Sigma})}\geq 
\sigma_\nu^2 c^2\|\eta_\nu\|^2_{L^2(\dot{\Sigma}_\epsilon)} -
\sigma_\nu c^\prime\|\eta_\nu\|^2_{L^2(\dot{\Sigma})}=
\sigma_\nu^2 c^2\|\eta_\nu\|^2_{L^2(\dot{\Sigma}_\epsilon)}-\sigma_\nu c^\prime,
$$ 
for some $c^\prime>0$. Thus we obtain 
$$
\|\eta_\nu\|^2_{L^2(\dot{\Sigma}_\epsilon)}\leq 
\frac{\|\mathbf{D}_{\sigma_\nu}\eta_\nu\|^2_{L^2(\dot{\Sigma})}}{c^2 \sigma_\nu^2}+\frac {c^\prime}{\sigma_\nu c^2}
\to 0 \quad\text{ as }\quad \nu \to \infty,
$$
so there is again concentration of energy near the 
zeroes of the antilinear perturbation: in particular,
\begin{equation*}
\begin{split} 1&=\lim_{\nu \rightarrow \infty}\|\eta_\nu\|^2_{L^2(\dot{\Sigma})}\\&=\lim_{\nu \rightarrow \infty}\|\eta_\nu\|^2_{L^2(\dot{\Sigma}_\epsilon)}+\lim_{\nu \rightarrow \infty}\sum_{\zeta \in Z^-(\beta)}\|\eta_\nu\|^2_{L^2(\dD(\zeta))}\\&=\lim_{\nu \rightarrow \infty}\sum_{\zeta \in Z^-(\beta)}\|\eta^\zeta_\nu\|^2_{L^2(\DD_{\sqrt{\sigma_\nu}})}.
\end{split}
\end{equation*}
Moreover, we have
$$
\mathbf{D}_-\eta_\nu^\zeta(z)= \frac{1}{\sigma_\nu}\overline{\partial}\eta_\nu\left(\frac{z}{\sqrt{\sigma_\nu}}\right)+\frac{\bar{z}}{\sqrt{\sigma_\nu}}\bar{\eta}_\nu\left(\frac{z}{\sqrt{\sigma_\nu}}\right)=\frac{1}{\sigma_\nu}\mathbf{D}_{\sigma_\nu} \eta_\nu \left(\frac{z}{\sqrt{\sigma_\nu}}\right).
$$ 
Taking the square of the norms on each side, we may integrate and use
change of variables to obtain
$$
\|\mathbf{D}_-\eta_\nu^\zeta\|_{L^2(\DD_{\sqrt{\sigma_\nu}})}=\frac{1}{\sqrt{\sigma_\nu}}\|\mathbf{D}_{\sigma_\nu}\eta_\nu\|_{L^2(\dD(\zeta))}
\to 0 \quad \text{ as } \quad \nu \to \infty.
$$
 
The elliptic estimates from Lecture~\ref{lec:local} now provide
uniform $H^k$-bounds for each $\eta_\nu^\zeta$ on compact subsets of $\CC$
for every $k \in \NN$, so that a subsequence converges in $C^\infty_{\loc}(\CC)$
to a smooth map $\eta_\infty^\zeta \in L^2(\CC,\CC)$ satisfying
$\mathbf{D}_- \eta_\infty^\zeta = 0$.  But $\sum_{\zeta \in Z^-(\beta)}
\| \eta_\infty^\zeta \|_{L^2(\CC)}^2 = 1$, so at least one of these
solutions is nontrivial and thus contradicts Proposition~\ref{prop:D-}.
\end{proof}

\begin{exercise}
\label{EX:adjointBound}
Show that for any Fredholm Cauchy-Riemann type operator $\mathbf{D}$ on~$E$,
the following two estimates are equivalent, with the same constant
$c > 0$ in both:
\begin{enumerate}
\renewcommand{\labelenumi}{(\roman{enumi})}
\item $\| \eta \|_{L^2(E)} \le c \| \mathbf{D} \eta \|_{L^2(F)}$ for all
$\eta \in H^1(E) \cap (\ker \mathbf{D})^\perp$;
\item $\| \lambda \|_{L^2(F)} \le c \| \mathbf{D}^*\lambda \|_{L^2(E)}$ for all
$\lambda \in H^1(F) \cap (\ker \mathbf{D}^*)^\perp$.
\end{enumerate}
\textsl{Hint: elliptic regularity implies that for $\mathbf{D}$ and 
$\mathbf{D}^*$ as bounded linear operators $H^1 \to L^2$, $(\ker\mathbf{D})^\perp
= \im \mathbf{D}^*$ and $(\ker \mathbf{D}^*)^\perp = \im \mathbf{D}$.}
\end{exercise}

\begin{proof}[Proof of Proposition~\ref{prop:inj}]
If the statement is not true, then there exist sequences $\sigma_\nu \to \infty$
and 
$$
f_\nu \in \bigoplus_{\zeta \in Z^+(\beta)} \ker \mathbf{D}_+
$$
such that $\| f_\nu \|_{L^2} = 1$ and $\eta_\nu := \Phi_{\sigma_\nu}(f_\nu) \in
(\ker \mathbf{D}_{\sigma_\nu})^\perp$ for all~$\nu$.
Lemmas~\ref{lemma:PhiBounds} and~\ref{lemma:injBound} then provide estimates
of the form
\begin{itemize}
\item $\| \eta_\nu \|_{L^2} \to 1$,
\item
$\| \mathbf{D}_{\sigma_\nu} \eta_\nu \|_{L^2} \to 0$, and
\item $\| \eta_\nu \|_{L^2} \le c \| \mathbf{D}_{\sigma_\nu} \eta_\nu \|_{L^2}$
\end{itemize}
as $\nu \to \infty$, with $c > 0$ independent of~$\nu$.  These imply:
$$
1 = \lim_{\nu \to \infty} \| \eta_\nu \|_{L^2} \le 
\lim_{\nu \to \infty} c \| \mathbf{D}_{\sigma_\nu} \eta_\nu \|_{L^2} = 0.
$$
\end{proof}

We've proved:

\begin{prop}
\label{prop:largePerturbation}
Suppose the assumptions of \S\ref{sec:concentration} hold and that the
section $\beta \in \Gamma(\Hom_\CC(\bar{E},F))$ has $I_+ \ge 0$ positive and 
$I_- \ge 0$ negative zeroes.  If $I_- = 0$,
then $\mathbf{D}_\sigma$ is surjective with
$\dim \ker \mathbf{D}_\sigma = I_+$ for all $\sigma > 0$ sufficiently large.
If $I_+ = 0$, then $\mathbf{D}_\sigma$ is injective with
$\dim \coker \mathbf{D}_\sigma = I_-$ for all $\sigma > 0$ sufficiently large.
In either case,
$$
\ind(\mathbf{D}_\sigma) = I_+ - I_-
$$
for all $\sigma > 0$ sufficiently large.
\qed
\end{prop}

\section{Antilinear deformations of asymptotic operators}
\label{sec:asymptotic}

Proposition~\ref{prop:largePerturbation} suffices to prove the index formula
in the closed case, but there is an additional snag if $\Gamma \ne \emptyset$:
since $H^1(\dot{\Sigma}) \hookrightarrow L^2(\dot{\Sigma})$ is not a
compact inclusion, we have no guarantee that $\mathbf{D}$ and
$\mathbf{D}_\sigma := \mathbf{D} + \sigma B$ will have the same index, and
generally they will not.  A solution to this problem has been pointed out by
Chris Gerig \cite{Gerig:thesis}, using a special class of asymptotic operators
that also originate in the work of Taubes (see \cite{Taubes:ECH=SWF1}*{Lemma~2.3}).

In general, the only obvious way to guarantee
$\ind(\mathbf{D}) = \ind(\mathbf{D}_\sigma)$ for large $\sigma > 0$ is if we
can arrange for \emph{every} operator in the family
$\{ \mathbf{D}_\sigma \}_{\sigma \ge 0}$ to be Fredholm, which is not
automatic since the zeroth-order perturbation $B : E \to F$ is required
to be bounded away from zero near $\infty$ and must therefore change
the asymptotic operators at the punctures.  We are therefore led to ask:

\begin{question}
For what nondegenerate asymptotic operators 
$\mathbf{A} : H^1(E) \to L^2(E)$ on a Hermitian line
bundle $(E,J,\omega) \to S^1$ can one find complex-antilinear bundle maps
$B : E \to E$ such that
$$
\mathbf{A}_\sigma := \mathbf{A} - \sigma B : H^1(E) \to L^2(E)
$$
is an isomorphism for every $\sigma \ge 0$?
\end{question}

It turns out that it will suffice to find, for each unitary trivialization
$\sigma$ and every $k \in \ZZ$, a particular pair $(\mathbf{A}_k,B_k)$
such that $\mathbf{A}_k - \sigma B_k$ is nondegenerate for all $\sigma \ge 0$
and $\muCZ^\tau(\mathbf{A}_k) = k$.  To see why, let us proceed under the
assumption that such pairs can be found, and use them to compute the index:

\begin{lemma}
\label{lemma:precomputation}
Given $\mathbf{D}$ as in Theorem~\ref{thm:RiemannRoch}, fix asymptotic 
trivializations $\tau$ and suppose that for each puncture $z \in \Gamma$
there exists an asymptotic operator $\mathbf{A}_z'$ on $(E_z,J_z,\omega_z)$
with $\muCZ^\tau(\mathbf{A}_z') = \muCZ^\tau(\mathbf{A}_z)$, such that if
$\mathbf{A}_z'$ is written with respect to $\tau$ as
$-J_0 \p_t - S_z(t)$, then the deformed asymptotic operator
\begin{equation}
\label{eqn:deformedAsymp}
C^\infty(S^1,\RR^2) \to C^\infty(S^1,\RR^2) : \eta \mapsto -J_0 \p_t \eta -
S_z(t)\eta - \sigma \beta_z(t) \bar{\eta}
\end{equation}
is nondegenerate for some loop $\beta_z : S^1 \to \CC \setminus \{0\}$
and every $\sigma \ge 0$.  Then
$$
\ind(\mathbf{D}) = \chi(\dot{\Sigma}) + 2 c_1^\tau(E) + \sum_{z \in \Gamma^+}
\wind(\beta_z) - \sum_{z \in \Gamma^-} \wind(\beta_z).
$$
\end{lemma}
\begin{proof}
Since $\muCZ^\tau(\mathbf{A}_z) = \muCZ^\tau(\mathbf{A}_z')$, we can deform
$\mathbf{A}_z$ to $\mathbf{A}_z'$ continuously through a family of
nondegenerate asymptotic operators.  It follows that we can deform
$\mathbf{D}$ through a continuous family of Fredholm Cauchy-Riemann type
operators to a new operator $\mathbf{D}'$ whose asymptotic operators
are $\mathbf{A}_z'$ for $z \in \Gamma$, and $\ind(\mathbf{D}') =
\ind(\mathbf{D})$.  We are free to assume in fact that $\mathbf{D}'$
is written with respect to the trivialization $\tau$ on the cylindrical
end near $z \in \Gamma^\pm$ as
$$
\p_s + J_0 \p_t + S_z(t).
$$
Now choose $\beta \in \Gamma(\Hom_\CC(\bar{E},F))$ with nondegenerate zeroes
such that the deformed
operators $\mathbf{D}_\sigma \eta := \mathbf{D}'\eta + \sigma \beta \bar{\eta}$
appear in trivialized form on the cylindrical end near $z \in \Gamma^\pm$ as
$$
\mathbf{D}_\sigma \eta = \p_s \eta + J_0 \p_t \eta + S_z(t) \eta + \sigma \beta_z(t) \bar{\eta}.
$$
This means $\mathbf{D}_\sigma$ is asymptotic at $z$ to \eqref{eqn:deformedAsymp},
which is nondegenerate for every $\sigma \ge 0$, implying $\mathbf{D}_\sigma$
is Fredholm for every $\sigma \ge 0$ and thus
$$
\ind(\mathbf{D}) = \ind(\mathbf{D}_\sigma).
$$

The trivializations $\tau$ induce trivializations over the cylindrical ends
for $\bar{E}$ and $F = \Lambda^{0,1}T^*\dot{\Sigma} \otimes E$, and the
expression for $\beta$ in the resulting asymptotic trivialization of
$\Hom_\CC(\bar{E},F)$ near $z \in \Gamma$ is~$\beta_z(t)$.  It follows that
the signed count of zeroes of $\beta$ is
\begin{equation*}
\begin{split}
i(\mathbf{D}) &:= c_1^\tau(\Hom_\CC(\bar{E},F)) + 
\sum_{z \in \Gamma^+} \wind(\beta_z) - \sum_{z \in \Gamma^-} \wind(\beta_z) \\
&= \chi(\dot{\Sigma}) + 2 c_1^\tau(E) +
\sum_{z \in \Gamma^+} \wind(\beta_z) - \sum_{z \in \Gamma^-} \wind(\beta_z),
\end{split}
\end{equation*}
where the computation $c_1^\tau(\Hom_\CC(\bar{E},F)) = \chi(\dot{\Sigma}) +
2 c_1^\tau(E)$ follows from the natural isomorphism
\begin{equation*}
\begin{split}
\Hom_\CC(\bar{E},F) &= \bar{E}^* \otimes F = E \otimes F =
E \otimes \Lambda^{0,1}T^*\dot{\Sigma} \otimes E =
\Lambda^{0,1}T^*\dot{\Sigma} \otimes E \otimes E \\
&= T\dot{\Sigma} \otimes E \otimes E.
\end{split}
\end{equation*}
We are free to assume that all zeroes of $\beta$ are either positive or
negative, depending on the sign of $i(\mathbf{D})$.
Proposition~\ref{prop:largePerturbation} then implies
$\ind(\mathbf{D}_\sigma) = i(\mathbf{D})$ for large~$\sigma$.
\end{proof}

Notice that instead of nondegenerate families 
$\mathbf{A} - \sigma B$
parametrized by $\sigma \in [0,\infty)$, it is just as well to find such families
which are nondegenerate and have the right Conley-Zehnder index for all
$\sigma > 0$, as the $\sigma \ge 1$ portion of this family can be rewritten as
$(\mathbf{A} - B) - \sigma B$ for $\sigma \ge 0$.
The following lemma thus completes the proof of Theorem~\ref{thm:RiemannRoch}.

\begin{lemma}
For every $k \in \ZZ$, the trivial Hermitian line bundle over $S^1$ admits
an asymptotic operator $\mathbf{A}_k$ and a loop
$\beta_k : S^1 \to \CC \setminus \{0\}$ such that the deformed asymptotic
operators 
$$
\mathbf{A}_{k,\sigma} \eta := \mathbf{A}_k \eta - \sigma \beta_k \bar{\eta}
$$
are nondegenerate for every $\sigma > 0$ and satisfy
$$
\muCZ(\mathbf{A}_{k,\sigma}) = \wind(\beta_k) = k.
$$
\end{lemma}
\begin{proof}
We claim that the choices
$$
\mathbf{A}_k \eta := -J_0 \p_t \eta - \pi k \eta \quad \text{ and }\quad
\beta_k(t) := e^{2\pi i kt}
$$
do the trick.  We prove this in three steps.

\textsl{Step~1: $k=0$.}
The above formula gives $\mathbf{A}_{0,\sigma} = -J_0 \p_t \eta - \sigma \bar{\eta}$,
in which the $\sigma=1$ case is precisely the operator that we used in
Lecture~\ref{lec:asymptotic} to normalize the Conley-Zehnder index, hence
$\muCZ(\mathbf{A}_{0,1}) = 0$ by definition.  More generally, all of these
operators can be expressed in the form $\mathbf{A} := -J_0 \p_t - S$
where $S \in \End_\RR(\RR^2)$ is a constant nonsingular $2$-by-$2$ symmetric
matrix that anticommutes with~$J_0$.  We claim that \emph{all} asymptotic 
operators of this form are nondegenerate.  Indeed, the conditions
$S^\transpose = S$ and $S J_0 = -J_0 S$ for
$J_0 = \begin{pmatrix} 0 & -1 \\ 1 & 0 \end{pmatrix}$ imply that
$S$ takes the form 
$\begin{pmatrix} a & b \\ b & -a \end{pmatrix}$ with $\det S = -a^2 - b^2 \ne 0$, 
and moreover $S$ is of this form if and only if $J_0 S$ also is.  In particular,
$J_0 S$ is traceless, symmetric, and nonsingular.  Solutions of 
$\mathbf{A} \eta = 0$ then satisfy 
$\dot{\eta} = J_0 S \eta$, which has no periodic solutions since
$J_0 S$ has one positive and one negative eigenvalue, hence
$\ker \mathbf{A} = \{0\}$.

\textsl{Step~2: even~$k$.}
There is a cheap trick to deduce the case $k=2m$ for any $m \in \NN$ from
the $k=0$ case.  Recall that by Exercise~\ref{EX:changeTriv} in Lecture~\ref{lec:asymptotic}, 
conjugating $\mathbf{A}_{0,\sigma}$ by a change of trivialization changes its 
Conley-Zehnder index by twice the degree of that change.  In particular,
the operator
$$
\tilde{\mathbf{A}}_{0,\sigma} \eta := 
e^{2\pi i m t} \mathbf{A}_{0,\sigma} (e^{-2\pi i m t} \eta)
$$
is also a nondegenerate asymptotic operator, but with
$\muCZ(\tilde{\mathbf{A}}_{0,\sigma}) = 
\muCZ(\mathbf{A}_{0,\sigma}) + 2m = k$.  Explicitly, we compute
$$
\tilde{\mathbf{A}}_{0,\sigma} \eta = -J_0 \p_t \eta - \pi k \eta 
- \sigma k e^{2\pi i kt} \bar{\eta},
$$
so $\mathbf{A}_{k,\sigma} = \tilde{\mathbf{A}}_{0,\sigma/k}$ is also
nondegenerate for every $\sigma > 0$.

\textsl{Step~3: odd~$k$.}
Another cheap trick relates each $\mathbf{A}_{k,\sigma}$ to
$\mathbf{A}_{2k,\sigma}$ after an adjustment in~$\sigma$.  Given an arbitrary
asymptotic operator $\mathbf{A} = -J_0 \p_t - S(t)$ and $m \in \NN$, define
$$
\mathbf{A}^m := -J_0 \p_t - m S(mt).
$$
Geometrically, if $\mathbf{A}$ is a trivialized representation for the
asymptotic operator of a Reeb orbit $\gamma : S^1 \to M$, then
$\mathbf{A}^m$ is the operator for the $m$-fold covered orbit
$\gamma^m : S^1 \to M : t \mapsto \gamma(mt)$.  It is easy to check in
particular that if we define $\eta^m(t) := \eta(mt)$ for any given
loop $\eta : S^1 \to \RR^2$, then
$$
\mathbf{A}^m \eta^m = m (\mathbf{A} \eta)^m,
$$
so this gives an embedding of $\ker \mathbf{A}$ into $\ker \mathbf{A}^m$,
implying that whenever $\mathbf{A}^m$ is nondegenerate for some $m \in \NN$,
so is $\mathbf{A}$.  To make use of this, observe that
$$
\mathbf{A}_{k,\sigma}^2 \eta = -J_0 \p_t \eta - \pi 2k \eta - 2\sigma
e^{4\pi i kt} \bar{\eta} = \mathbf{A}_{2k,2\sigma}\eta,
$$
so $\mathbf{A}_{k,\sigma}^2$ is nondegenerate for all $\sigma > 0$ by
Step~2, and therefore so is~$\mathbf{A}_{k,\sigma}$.
\end{proof}

The proof of Theorem~\ref{thm:RiemannRoch} is now complete.

\begin{exercise}
Derive a Weitzenb\"ock formula for asymptotic operators and use it to show
that for any asymptotic operator $\mathbf{A}$ on the trivial Hermitian
line bundle and any smooth $\beta : S^1 \to \CC \setminus \{0\}$, 
the deformed operators
$\mathbf{A}_\sigma \eta := \mathbf{A}\eta - \sigma \beta \bar{\eta}$ are all
nondegenerate for $\sigma > 0$ sufficiently large.
Deduce from this that $\muCZ(\mathbf{A}_\sigma) = \wind(\beta)$ for large 
$\sigma > 0$.
\end{exercise}

\psfrag{What}{$\widehat{W}$}
\psfrag{positive end}{$((-\epsilon,0] \times M_+,d(r \lambda_+) + \omega_+)$}
\psfrag{negative end}{$([0,\epsilon) \times M_-, d(r \lambda_-) + \omega_-)$}
\psfrag{W}{$(W,\omega)$}
\psfrag{collar+}{$((-\epsilon,0] \times M_+, d(r \lambda_+) + \omega_+)$}
\psfrag{collar-}{$([0,\epsilon) \times M_-, d(r \lambda_-) + \omega_-)$}
\psfrag{end+}{$([0,\infty) \times M_+, d(\varphi(r)\lambda_+) + \omega_+)$}
\psfrag{end-}{$((-\infty,0] \times M_-, d(\varphi(r)\lambda_-) + \omega_-)$}
\psfrag{Sigmadot}{$\dot{\Sigma}$}
\psfrag{u}{$u$}

\chapter{Symplectic cobordisms and moduli spaces}
\label{lec:cobordisms}

\minitoc

\vspace{12pt}

In this lecture we introduce the moduli spaces of holomorphic curves
that are used to define SFT.

\section{Stable Hamiltonian structures and their symplectizations}

In Lecture~\ref{lec:intro}, we motivated the notion of a contact manifold by considering
hypersurfaces $M$ in a symplectic manifold $(W,\omega)$ that satisfy
a \emph{convexity} (also known as ``contact type'') condition.  The point
of that condition was that it presents $M$ as one member of a smooth
$1$-parameter family of hypersurfaces that all have the same Hamiltonian
dynamics; that $1$-parameter family furnishes the basic model of
what we call the \emph{symplectization} of $M$ with its induced contact
structure.  A useful generalization of this notion was introduced in
\cite{HoferZehnder} and was later recognized to be the most natural 
geometric setting for punctured holomorphic curves.  It has the advantage
of allowing us to view seemingly distinct theories such as Hamiltonian Floer
homology as special cases of SFT---and even if we are only interested in
contact manifolds, the generalization sometimes makes computations easier
than they might be in a purely contact setting.

Recall that every smooth hypersurface $M$ in a $2n$-dimensional symplectic
manifold $(W,\omega)$ has a \defin{characteristic line field}
$$
\ker \left( \omega|_{TM} \right) \subset TM,
$$
whose integral curves are the orbits on $M$ of any Hamiltonian vector field
generated by a function $H : W \to \RR$ that has $M$ as a regular level set.
We say that $M \subset (W,\omega)$ is \defin{stable} if a neighborhood of 
$M$ admits a \defin{stabilizing
vector field}~$V$: this means that $V$ is transverse to $M$ and the
$1$-parameter family of hypersurfaces
$$
M_t := \varphi_V^t(M), \qquad -\epsilon < t < \epsilon
$$
generated by the flow $\varphi_V^t$ of $V$ has the property that each of
the diffeomorphisms $M \to M_t$ defined by flowing along $V$ preserves
characteristic line fields.

\begin{exercise}
\label{EX:SHS}
Show that if $V$ is a stabilizing vector field for $M \subset (W,\omega)$,
then the $2$-form and $1$-form pair $(\Omega,\Lambda)$ defined on $M$ by
$$
\Omega := \omega|_{TM}, \qquad \Lambda := \iota_V \omega|_{TM}
$$
has the following properties:
\begin{enumerate}[label=(\roman*)]
\item $\Omega|_{\ker \Lambda}$ is nondegenerate;
\label{item:nondeg}
\item $\ker \Omega \subset \ker d\Lambda$.
\label{item:kernel}
\end{enumerate}
Show moreover that if $M$ is assigned the orientation for which $V$ is
\emph{positively} transverse to $M$ and $\xi := \ker\Lambda \subset TM$ is
assigned the natural co-orientation determined by $\Lambda$, then the
induced orientation of $\xi$ matches the orientation determined by
the symplectic vector bundle structure~$\Omega|_{\xi}$, hence 
condition~\ref{item:nondeg} can equivalently be written as
\begin{enumerate}[label=(\roman*)]
\setcounter{enumi}{2}
\item $\Lambda \wedge \Omega^{n-1} > 0$
\label{item:positive}
\end{enumerate}
where $\dim W = 2n$.
\end{exercise}

A \defin{stable Hamiltonian structure} (or ``SHS'' for short)
on an arbitrary oriented 
$(2n-1)$-dimensional manifold $M$ is a pair
$(\Omega,\Lambda)$ consisting of a closed $2$-form $\Omega$ and $1$-form
$\Lambda$ such that properties~\ref{item:kernel} and~\ref{item:positive}
in Exercise~\ref{EX:SHS} are satisfied.

\begin{exercise}
Show that if $(\Omega,\Lambda)$ is a stable Hamiltonian structure, then
$$
\omega := d(r\Lambda) + \Omega
$$
is a symplectic form on $(-\epsilon,\epsilon) \times M$
for $\epsilon > 0$ sufficiently small, where $r$ denotes the coordinate
on $(-\epsilon,\epsilon)$; moreover, $\{0\} \times M$ is a stable
hypersurface in $((-\epsilon,\epsilon) \times M,\omega)$.
\end{exercise}

\begin{example}
If $M \subset (W,\omega)$ is a contact type hypersurface, then a
Liouville vector field $V$ transverse to $M$ is a stabilizing vector field,
and the induced stable Hamiltonian structure is $(d\alpha,\alpha)$,
where $\alpha := \lambda|_{TM}$ with $\lambda := \omega(V,\cdot)$.
We will refer to this example henceforward as the \defin{contact case}.
\end{example}

\begin{prop}
\label{prop:collar}
Suppose $M \subset (W,\omega)$ is a closed stable hypersurface with stabilizing
vector field $V$ and induced stable Hamiltonian structure $(\Omega,\Lambda)$
where $\Omega = \omega|_{TM}$ and $\Lambda = \iota_V \omega|_{TM}$.
Then a neighborhood of $M$ in $(W,\omega)$ admits a symplectomorphism
to $((-\epsilon,\epsilon) \times M , d(r\Lambda) + \Omega)$ for some
$\epsilon > 0$, identifying $M \subset W$ with $\{0\} \times M \subset
(-\epsilon,\epsilon) \times M$.
\end{prop}
\begin{proof}
By the smooth tubular neighbourhood theorem and the preceeding exercise, we can view $\omega_0=d(r\Lambda) + \Omega$ as a symplectic form in some neighbourhood $\uU_0 \cong((-\epsilon,\epsilon) \times M)$ of $M$. In this neighbourhood,
$$
(\omega_0-\omega)|_M=0
$$
by definition of $\omega_0$ and thus
$$
\omega_0-\omega=d\mu 
$$
for some 1-form $\mu$ such that $\mu|_M=0$. Now define 
$$\omega_t = \omega +t \, d\mu
$$ 
and observe that it is a closed 2-form which can be assumed to be non-degenerate for a small enough choice of $\uU_0$.
Solving the Moser equation $$
\iota_{v_t} \omega_t = -\mu$$
yields a well-defined, time-dependent vector field $v_t$ with the property that $v_t|_M=0$. Working back we produce an isotopy as follows:
$$d\iota_{v_t} \omega_t = -d\mu \Rightarrow$$
$$\lL_{v_t}\omega_t=d\iota_{v_t} \omega_t +\iota_{v_t} d\omega_t=d\iota_{v_t} \omega_t=-d\mu = -\frac{d\omega_t}{dt} \Rightarrow$$
$$\frac{d}{dt}(\rho^*_t\omega_t)=\lL_{v_t}\omega_t+\frac{d\omega_t}{dt}=0$$
where $\rho^*_t$ is the flow of $v_t$. Then $$\rho^*_t \omega_t= \rho^*_0 \omega = \omega$$ since $\rho_0$ is the identity. The required symplectomorphism is then 
$$\rho_1:\rho_1^{-1}\uU_0 \to \uU_0$$
and the fact that $M$ is fixed under the isotopy follows from $v_t|_M=0$. 
\end{proof}

\begin{example}
In the contact case $(\Omega,\Lambda) = (d\alpha,\alpha)$, the symplectic form
on the collar neighborhood in Proposition~\ref{prop:collar} can be rewritten as 
$d(e^t \alpha)$ by defining the coordinate $t := \ln (r+1)$.  The proposition
is easier to prove in this case: one can construct the collar neighborhood simply
by flowing along~$V$, with no need for the Moser isotopy trick.
\end{example}

A stable Hamiltonian structure $\hH = (\Omega,\Lambda)$ 
gives rise to two important additional objects: a co-oriented
hyperplane distribution 
$$
\xi := \ker\Lambda,
$$
and a positively transverse vector field~$R$ determined by the conditions
$$
\Omega(R,\cdot) \equiv 0 \quad\text{ and }\quad
\Lambda(R) \equiv 1.
$$
By analogy with the contact case, we will refer to $R$ as the
\defin{Reeb vector field} of~$\hH$.  The condition 
$\ker \Omega \subset \ker d\Lambda$ implies that it reduces to the usual
contact notion of the Reeb vector field for $\Lambda$ whenever the latter
happens also to be a contact form.

The \defin{symplectization of $(M,\hH)$} for any stable Hamiltonian structure
$\hH = (\Omega,\Lambda)$ can be defined by choosing suitable
diffeomorphisms of $(-\epsilon,\epsilon) \times M$ with $\RR \times M$:
equivalently, this means we consider $\RR \times M$ with the family of
symplectic forms $\omega_\varphi$ defined by
\begin{equation}
\label{eqn:SHSsymplectization}
\omega_\varphi := 
d\left( \varphi(r) \Lambda \right) + \Omega
\end{equation}
where $\varphi$ is chosen arbitrarily from the set
\begin{equation}
\label{eqn:tT2}
\tT := \left\{ \varphi \in C^\infty(\RR,(-\epsilon,\epsilon))\ \big|\ 
\varphi' > 0 \right\}.
\end{equation}

\begin{example}
\label{ex:FloerSHS}
The following stable Hamiltonian structure places Hamiltonian Floer homology
into the setting of SFT.  Suppose $(W,\omega)$ is a closed symplectic manifold
and $H : S^1 \times W \to \RR$ is a smooth function, and denote
$H_t := H(t,\cdot) : W \to \RR$.  The time-dependent Hamiltonian vector
field $X_t$ defined by $d H_t = -\omega(X_t,\cdot)$ can then be viewed as
defining a \emph{symplectic connection} on the trivial symplectic fiber bundle
$$
M := S^1 \times W \stackrel{t}{\longrightarrow} S^1,
$$
i.e.~the flow of $R(t,x) := \p_t + X_t(x)$ defines symplectic parallel
transport maps between fibers.  The horizontal subbundle for this connection is the 
``symplectic complement'' of the vertical subbundle with respect to the
closed $2$-form
$$
\Omega = \omega + dt \wedge dH.
$$
In other words, $\Omega$ restricts to the fibers of $M \to S^1$ as $\omega$
and the subbundle 
$\{ X \in TM\ |\ \omega(X,\cdot)|_{T(\{\text{const}\} \times W)} \}$
is generated by~$R$, so $\Omega$ is the \defin{connection $2$-form} defining
the connection, cf.~\cite{McDuffSalamon:ST}.  Setting $\Lambda := dt$ then makes
$\hH := (\Omega,\Lambda)$ a stable Hamiltonian structure with Reeb vector
field~$R$, and its closed orbits in homotopy classes that project to $S^1$ 
with degree one are in $1$-to-$1$ correspondence
with the $1$-periodic Hamiltonian orbits on~$W$.  Notice that this is very
different from the contact case: $\xi = \ker dt$ is as far as possible from
being a contact structure, it is instead an integrable distribution whose
integral submanifolds are the fibers of $M \to S^1$.
\end{example}

\begin{exercise}
Show that for any stable Hamiltonian structure $\hH = (\Omega,\Lambda)$,
the flow of $R$ preserves $\xi = \ker\Lambda$ along with its symplectic
bundle structure $\Omega|_\xi$.
\end{exercise}

\begin{defn}
A $T$-periodic orbit $x : \RR \to M$ of $R$ is called \defin{nondegenerate}
if $1$ is not an eigenvalue of $d\varphi^T|_{\xi_{x(0)}} :
\xi_{x(0)} \to \xi_{x(0)}$, where $\varphi^t$ denotes the flow of~$R$.
\end{defn}

\begin{exercise}
Show that in Example~\ref{ex:FloerSHS}, the notions of nondegeneracy for
closed Reeb orbits on $M$ and for $1$-periodic Hamiltonian orbits
on $W$ (see Lecture~\ref{lec:intro}) coincide.
\end{exercise}

If $\gamma : S^1 \to M$ parametrizes a $T$-periodic orbit of $R$ with
$\dot{\gamma} = T \cdot R(\gamma)$, then the formula of Lecture~\ref{lec:asymptotic} for the
\defin{asymptotic operator}
$$
\mathbf{A}_\gamma \eta = -J (\nabla_t \eta - T \nabla_\eta R)
$$
still makes sense in this more general context, and it
defines an $L^2$-symmetric operator on the Hermitian
vector bundle $(\gamma^*\xi,J,\Omega)$ over~$S^1$.  It can also be interpreted
as a Hessian at a critical point, though for an action functional that is only
locally defined: indeed, while $\Omega$ need not be globally exact, it is
necessarily exact on a neighborhood of $\gamma_0(S^1)$ for any given 
loop $\gamma_0 : S^1 \to M$, so one can pick any
primitive $\lambda$ of $\Omega$ on this neighborhood and, for a sufficiently small 
neighborhood $\uU(\gamma_0) \subset C^\infty(S^1,M)$ of~$\gamma_0$, consider
the action functional
\begin{equation}
\label{eqn:SHSaction}
\aA_\hH : \uU(\gamma_0) \to \RR : \gamma \mapsto \int_{S^1} \gamma^*\lambda.
\end{equation}
Its first variation at $\gamma \in \uU(\gamma_0)$ in the direction 
$\eta \in \Gamma(\gamma^*\xi)$ is then
$$
d\aA_\hH(\gamma) \eta = - \int_{S^1} \Omega(\dot{\gamma},\eta) \, dt =
\langle - J \pi_\xi \dot{\gamma} , \eta \rangle_{L^2},
$$
where $\pi_\xi : TM \to \xi$ denotes the projection along~$R$ and the
$L^2$-pairing on $\gamma^*\xi$ is defined via the bundle metric
$\Omega(\cdot,J\cdot)|_\xi$.  This leads us to interpret 
$-J \pi_\xi \dot{\gamma}$ as a ``gradient'' $\nabla\aA_\hH(\gamma)$, and
if $\dot{\gamma} = T \cdot R(\gamma)$, then differentiating this gradient
in the direction of $\eta \in \Gamma(\gamma^*\xi)$ gives
$\mathbf{A}_\gamma \eta$.  As one would expect, nondegeneracy of $\gamma$
is then equivalent to the condition $\ker \mathbf{A}_\gamma = \{0\}$,
and one can in this case define the Conley-Zehnder index
$\muCZ^\tau(\gamma) \in \ZZ$ as in Lecture~\ref{lec:asymptotic}, relative to a choice of
unitary trivialization $\tau$ for $(\xi,J,\Omega)$.

\begin{exercise}
In the setting of Example~\ref{ex:FloerSHS}, work out the relationship
between $\aA_\hH$ and the symplectic action functional for Hamiltonian
systems that we discussed in Lecture~\ref{lec:intro}.  (Try not to worry too much about
signs.)
\end{exercise}

\begin{defn}
Given a stable Hamiltonian structure $\hH = (\Omega,\Lambda)$, denote by
$$
\jJ(\hH) \subset \jJ(\RR \times M)
$$
the space of smooth almost complex structures $J$ on $\RR \times M$ with the
following properties:
\begin{itemize}
\item $J$ is invariant under the $\RR$-action on $\RR \times M$ by translation
of the first factor;
\item $J \p_r = R$ and $J R = -\p_r$, where $r$ denotes the natural
coordinate on the first factor;
\item $J(\xi) = \xi$ and $J|_{\xi}$ is compatible with the symplectic vector
bundle structure $\Omega|_\xi$.
\end{itemize}
\end{defn}

Notice that if $\hH = (d\alpha,\alpha)$ for a contact form $\alpha$, then
$\jJ(\hH)$ matches the space $\jJ(\alpha)$ defined in Lecture~\ref{lec:intro}.

\begin{exercise}
\label{EX:positiveEnergy}
Show that every $J \in \jJ(\hH)$ is tamed by all of the symplectic
structures $\omega_\varphi$ as defined in \eqref{eqn:SHSsymplectization}
for $\varphi \in \tT$.
\end{exercise}

Given $J \in \jJ(\hH)$, we define the \defin{energy} of a $J$-holomorphic
curve $u : (\Sigma,j) \to (\RR \times M,J)$ by
$$
E(u) := \sup_{\varphi \in \tT} \int_\Sigma u^*\omega_\varphi.
$$
Exercise~\ref{EX:positiveEnergy} above implies that $E(u) \ge 0$, with
equality if and only if $u$ is constant.
In the contact case, this notion of energy is not identical to the
``Hofer energy'' that we defined in Lecture~\ref{lec:intro}, nor to Hofer's original
definition from \cite{Hofer:weinstein}, but all three are equivalent for
our purposes since uniform bounds on any of them imply uniform bounds
on the others.

Just as in the contact case, the simplest example of a finite-energy
$J$-holomorphic curve is a \defin{trivial cylinder}
$$
u_\gamma : \RR \times S^1 \to \RR \times M : (s,t) \mapsto (Ts,\gamma(t)),
$$
where $\gamma : S^1 \to M$ is a ``constant velocity'' parametrization of
a $T$-periodic orbit of $R$, i.e.~$\dot{\gamma} = T \cdot R(\gamma)$.
More generally, given a punctured Riemann surface 
$(\dot{\Sigma} = \Sigma \setminus \Gamma,j)$ with 
$\Gamma = \Gamma^+ \cup \Gamma^-$, we consider \defin{asymptotically 
cylindrical} $J$-holomorphic curves $u : (\dot{\Sigma},j) \to (\RR \times M,J)$,
which are assumed to have the property that for each $z \in \Gamma^\pm$,
there exist holomorphic cylindrical coordinates identifying a punctured
neighborhood $\dot{\uU}_z \subset \dot{\Sigma}$
of $z$ with $Z_+ = [0,\infty) \times S^1$ or
$Z_- = (-\infty,0] \times S^1$ respectively, and a trivial cylinder
$u_{\gamma_z} : \RR \times S^1 \to \RR \times M$ such that
$$
u(s,t) = \exp_{u_{\gamma_z}(s,t)} h_z(s,t) \quad \text{ for $|s|$ sufficiently
large},
$$
where $h_z(s,t)$ is a vector field along $u_{\gamma_z}$ satisfying
$|h_z(s,\cdot)| \to 0$ uniformly as $s \to \pm\infty$.  As usual, both the
norm $|h_z(s,t)|$ and the exponential
map here are assumed to be defined with respect to a translation-invariant
choice of Riemannian metric on $\RR \times M$.  The vector fields $h_z$
along $u_{\gamma_z}$ for each $z \in \Gamma$ are sometimes called 
\defin{asymptotic representatives} of $u$ near~$z$.  

Asymptotic representatives satisfy a regularity estimate that will be
important to know about, though its proof (given originally in
\cite{HWZ:props1}) would be too lengthy to present here.  The methods
behind the following statement involve a combination of nonlinear
regularity arguments as in Lecture~\ref{lec:local} with the asymptotic elliptic 
estimates from Lecture~\ref{lec:Fredholm}.  To prepare for the statement, note that $\hH$
induces a splitting of complex vector bundles
\begin{equation}
\label{eqn:splitting}
T(\RR \times M) = \epsilon \oplus \xi,
\end{equation}
where $\epsilon$ denotes the trivial complex line bundle generated by the
vector field $\p_r$, or equivalently, the Reeb vector field.  It follows
that if $\gamma : S^1 \to M$ is a Reeb orbit
and $u_\gamma : \RR \times S^1 \to \RR \times M$ is the corresponding
trivial cylinder, then any unitary trivialization $\tau$ of the Hermitian
bundle $(\gamma^*\xi,J,\Omega)$ naturally induces a trivialization of
$u_\gamma^*T(\RR \times M)$.

\begin{prop}[\cite{HWZ:props1}]
\label{prop:expDecay}
Assume $J \in \jJ(\hH)$, $u : (\dot{\Sigma},j) \to (\RR \times M,J)$ is
$J$-holomorphic and asymptotically cylindrical, and its asymptotic orbit
$\gamma_z$ at $z \in \Gamma^\pm$ is nondegenerate.  Let
$h(s,t) \in \CC^n$ denote the asymptotic representative of $u$ near~$z$
expressed via the trivialization induced by a choice of unitary trivialization
for $(\gamma_z^*\xi,J,\Omega)$.  If $\delta > 0$ is small enough 
so that the asymptotic operator $\mathbf{A}_{\gamma_z}$ has no eigenvalues
in the closed interval between $0$ and $\mp \delta$, then
$$
h(s,t) = e^{\mp \delta s} g(s,t)
$$
for some bounded function $g(s,t) \in \CC^n$ whose derivatives of all orders 
are bounded as $s \to \pm\infty$.
\end{prop}

\begin{remark}
The range of $\delta > 0$ for which Prop.~\ref{prop:expDecay} holds is
open, thus by adjusting $\delta$ slightly, one can equivalently say that
$h(s,t) = e^{\mp \delta s} g(s,t)$ where the derivatives of all orders 
of $g(s,t)$ \emph{decay to zero} as $s \to \pm\infty$.
\end{remark}

\begin{exercise}
Convince yourself that the analogue of Proposition~\ref{prop:expDecay}
in Morse theory is true.  Namely, suppose $(M,g)$ is a Riemannian manifold,
$f : M \to \RR$ is smooth and $u : \RR \to M$ is a solution to
$\dot{u} + \nabla f(u) = 0$ with $\lim_{s \to \pm\infty} u(s) = x_\pm
\in \Crit(f)$, where $x_\pm$ are nondegenerate critical points.  We can write
$u(s)$ asymptotically as
$$
u(s) = \exp_{x_\pm} h_\pm(s)
$$
for some functions $h_\pm(s) \in T_{x_\pm}M$ that are defined for $s$ close
to $\pm\infty$ and satisfy $|h_\pm(s)| \to 0$ as $s \to \pm\infty$.
Show that if $\delta > 0$ is small enough so that $\nabla^2 f(x_\pm)$ has
no eigenvalue in the closed interval between $0$ and $\pm \delta$, then
$$
h_\pm(s) = e^{\mp \delta s} g_\pm(s)
$$
for some functions $g_\pm(s)$ with bounded derivatives of all orders
as $s \to \pm\infty$.\footnote{The apparent discrepancy in signs between
this and Proposition~\ref{prop:expDecay} is due to the fact that
$u(s)$ satisfies a \emph{negative} gradient flow equation, whereas the
nonlinear Cauchy-Riemann equation in symplectizations is interpreted loosely
as a \emph{positive} gradient flow equation.}
\textsl{Hint: fix local coordinates identifying $x_\pm$ with $0 \in \RR^n$
and first consider the case where $\nabla f(x)$ in these coordinates
depends linearly on~$x$.  Then try to compare $u(s)$ with solutions of this
idealized equation.}
\end{exercise}

\begin{example}
In the setting of Example~\ref{ex:FloerSHS}, a choice of $J \in \jJ(\hH)$
is equivalent to a choice of smooth $S^1$-parametrized family of
compatible almost complex structures $\{J_t\}_{t \in S^1}$ on $(W,\omega)$,
and $J$-holomorphic curves $u : (\dot{\Sigma},j) \to (\RR \times M,J)$ can
be written as
$$
u = (f,v) : \dot{\Sigma} \to \left( \RR \times S^1 \right) \times W,
$$
where $f : (\dot{\Sigma},j) \to (\RR \times S^1,i)$ is holomorphic.  In
particular, if $(\dot{\Sigma},j) = (\RR \times S^1,i)$ and $f$ is taken to
have an extension to $S^2 \to S^2$ of degree one, then $u$ can be 
reparametrized so that $f$ is the identity map, hence
$u = (\Id,v) : \RR \times S^1 \to (\RR \times S^1) \times W$ is a section
of the trivial fiber bundle $(\RR \times S^1) \times W \to \RR \times S^1$,
and one can check that the equation satisfied by $v : \RR \times S^1 \to W$
is precisely the Floer equation
$$
\p_s v + J_t(v) (\p_t v - X_t(v)) = 0.
$$
\end{example}

\section{Symplectic cobordisms with stable boundary}
\label{sec:stableBoundary}

We discussed symplectic cobordisms between contact manifolds in
Lecture~\ref{lec:intro}.  Let us now generalize this notion in the context of stable
Hamiltonian structures.

\begin{figure}
\includegraphics{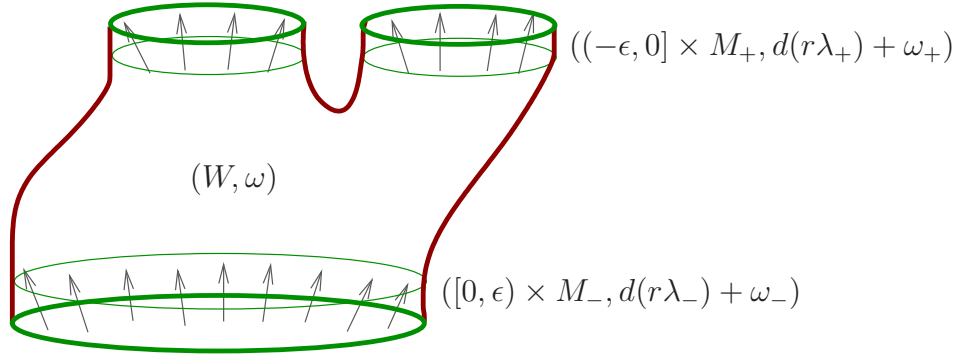}
\caption{\label{fig:collars2} A symplectic cobordism with stable boundary
components $\p W = - M_- \sqcup M_+$ and symplectic collar neighborhoods 
induced by the stable Hamiltonian structures 
$\hH_\pm = (\omega_\pm,\lambda_\pm)$ on~$M_\pm$.}
\end{figure}

A \defin{symplectic cobordism with stable boundary} is a compact symplectic
manifold $(W,\omega)$ with boundary
$\p W = -M_- \sqcup M_+$, equipped with a stabilizing
vector field $V$ that points transversely inward at $M_-$ and outward at~$M_+$.
This induces stable Hamiltonian structures $\hH_\pm = (\omega_\pm,\lambda_\pm)$
on $M_\pm$, where
$$
\omega_\pm := \omega|_{T M_\pm}, \qquad \lambda_\pm := (\iota_V \omega)|_{T M_\pm},
$$
and observe that the orientation conventions for $M_+$ and $M_-$ (with the
latter carrying the opposite of the natural boundary orientation) have been 
chosen such that if $\dim W = 2n$,
$$
\lambda_\pm \wedge \omega_\pm^{n-1} > 0
\quad \text{ on $M_\pm$}.
$$
We can now identify neighborhoods of $M_\pm$ in $(W,\omega)$ symplectically with
collars of the form
\begin{equation*}
\begin{split}
& \left( [0,\epsilon) \times M_+ , d\left( r \lambda_+ \right) + \omega_+ \right), \\
& \left( (-\epsilon,0] \times M_- , d\left( r \lambda_- \right) + \omega_- \right),
\end{split}
\end{equation*}
see Figure~\ref{fig:collars2}.

\begin{figure}
\includegraphics{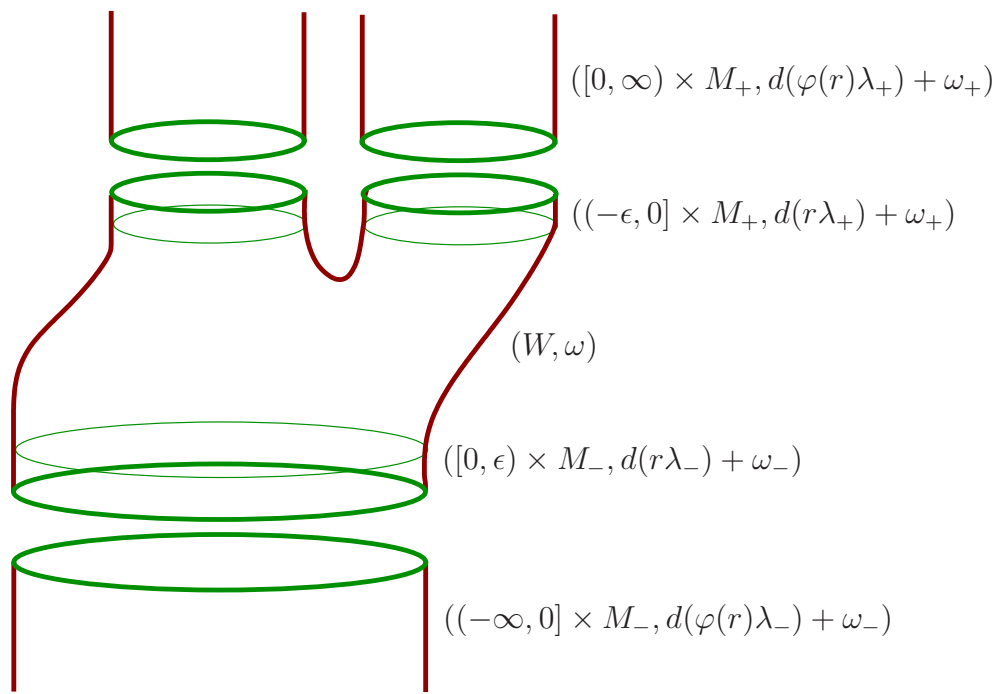}
\caption{\label{fig:completion2} The completion $(\widehat{W},\omega_\varphi)$
of a symplectic cobordism with stable boundary.}
\end{figure}

Modifying \eqref{eqn:tT2} by
\begin{equation}
\label{eqn:tTcob}
\tT_0 := \left\{ \varphi \in C^\infty(\RR,(-\epsilon,\epsilon))\ \big|\ 
\varphi' > 0 \text{ and $\varphi(r) = r$ for $r$ near~$0$} \right\},
\end{equation}
we can use any $\varphi \in \tT_0$ to define a \defin{symplectic completion}
$(\widehat{W},\omega_\varphi)$ of $(W,\omega)$ by
$$
\widehat{W} := \big( (-\infty,0] \times M_- \big) \cup_{M_-} W
\cup_{M_+} \big( [0,\infty) \times M_+ \big),
$$
where the above collar neighborhoods are used to glue the pieces together
smoothly and the symplectic form is defined by
$$
\omega_\varphi := \begin{cases}
d\left( \varphi(r) \lambda_- \right) + \omega_- & 
\text{ on $(-\infty,0] \times M_-$},\\
\omega & \text{ on $W$},\\
d\left( \varphi(r) \lambda_+ \right) + \omega_+ & 
\text{ on $[0,\infty) \times M_+$},
\end{cases}
$$
see Figure~\ref{fig:completion2}.  For each $r_0 \ge 0$, we define the compact
submanifold
$$
W^{r_0} := \left( [-r_0,0] \times M_- \right) \cup_{M_-} W
\cup_{M_+} \left( [0,r_0] \times M_+ \right),
$$
and observe that $(W^{r_0},\omega_\varphi)$ is also a symplectic cobordism with
stable boundary for every $\varphi \in \tT_0$.

Since $\widehat{W}$ is noncompact, almost complex structures $J$ on 
$\widehat{W}$ will need to satisfy conditions near infinity in order for 
moduli spaces of $J$-holomorphic curves to be well behaved, but we would 
like to preserve the freedom of choosing arbitrary compatible or tame
almost complex structures in compact subsets.

\begin{defn}
\label{defn:stableCobJ}
Given $\psi \in \tT_0$ and $r_0 \ge 0$, let
$$
\jJ_\tau(\omega_\psi,r_0,\hH_+,\hH_-) \subset \jJ(\widehat{W})
$$
denote the space of smooth almost complex structures $J$ on 
$\widehat{W}$ such that:
\begin{itemize}
\item $J$ on $[r_0,\infty) \times M_+$ matches an element of $\jJ(\hH_+)$;
\item $J$ on $(-\infty,-r_0] \times M_-$ matches an element of $\jJ(\hH_-)$;
\item $J$ on $W^{r_0}$ is tamed by $\omega_\psi$.
\end{itemize}
Let
$$
\jJ(\omega_\psi,r_0,\hH_+,\hH_-) \subset \jJ_\tau(\omega_\psi,r_0,\hH_+,\hH_-)
$$
denote the subset for which $J$ is additionally compatible with
$\omega_\psi$ on~$W^{r_0}$.
\end{defn}

\begin{figure}
\includegraphics{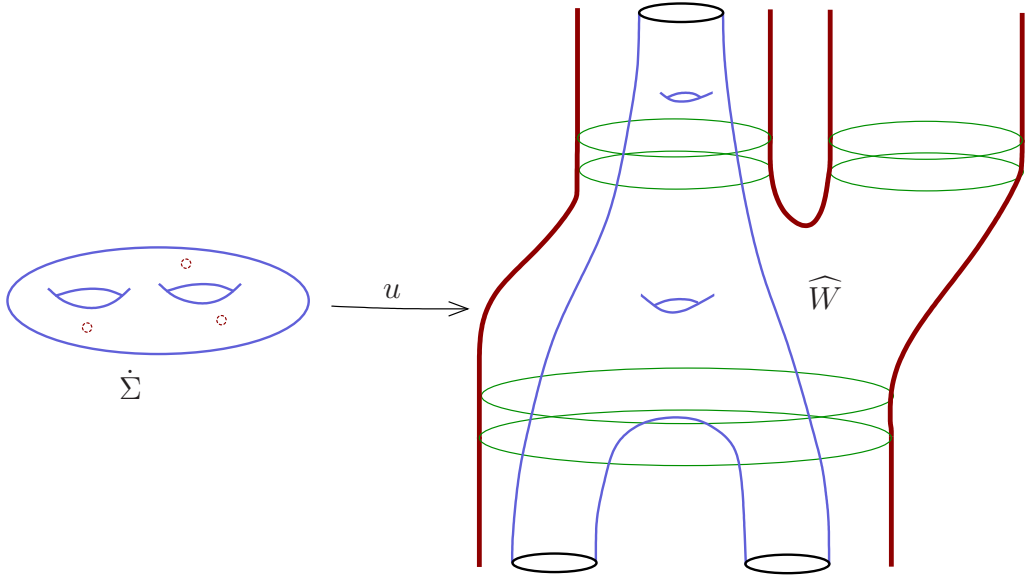}
\caption{\label{fig:asympCyl2} 
An asymptotically cylindrical holomorphic curve in $(\widehat{W},J)$
with genus~$2$, one positive puncture and two negative punctures.}
\end{figure}

Setting
\begin{equation}
\label{eqn:tTcobBig}
\tT(\psi,r_0) := \left\{ \varphi \in \tT_0 \ \big|\ 
\text{$\varphi \equiv \psi$ on $[-r_0,r_0]$} \right\},
\end{equation}
Exercise~\ref{EX:positiveEnergy} implies that every 
$J \in \jJ(\omega_\psi,r_0,\hH_+,\hH_-)$ is tamed by $\omega_\varphi$ 
for every $\varphi \in \tT(\psi,r_0)$.  It is therefore sensible to define the 
energy of a $J$-holomorphic curve $u : (\Sigma,j) \to (\widehat{W},J)$ by
$$
E(u) := \sup_{\varphi \in \tT(\psi,r_0)} \int_\Sigma u^*\omega_\varphi.
$$
The notion of asymptotically cylindrical $J$-holomorphic curves extends in a
straightforward way to the setting of $(\widehat{W},J)$: such curves
are proper maps whose positive/negative punctures are asymptotic to
closed orbits of the Reeb vector field $R_\pm$ induced by $\hH_\pm$ on
$\{\pm\infty\} \times M_\pm$, see Figure~\ref{fig:asympCyl2}.
The exponential decay estimate in Proposition~\ref{prop:expDecay} is also
immediately applicable in this more general setting since asymptotically
cylindrical curves in $\widehat{W}$ are indistinguishable near their
punctures from curves in the symplectizations $\RR \times M_\pm$.

It is easy to check that asymptotically cylindrical $J$-holomorphic curves
always have finite energy.  We will prove in Lecture~\ref{lec:Dragnev} that the converse
is also true whenever the Reeb orbits are nondegenerate.

\begin{remark}
Strictly speaking, the ``trivial stable cobordism'' 
$$
([0,1] \times M,d(\varphi(r)\Lambda,\Omega))
$$
induces different stable
Hamiltonian structures at $M_- := \{0\} \times M$ and $M_+ := \{1\} \times M$,
thus one cannot technically regard $\jJ(\hH)$ as contained in any space
of the form $\jJ(\omega_\psi,r_0,\hH_+,\hH_-)$ without inventing questionable
new notions such as the ``infinitesimal trivial cobordism''
$[0,0] \times M$ (whose completion would be the symplectization of
$(M,\hH)$).  It is nonetheless true for fairly trivial reasons that most
results about $\jJ(\omega,r_0,\hH_+,\hH_-)$ apply equally well to $\jJ(\hH)$,
and we shall use this fact in the following without always mentioning it.
\end{remark}

Every asymptotically cylindrical curve $u : \dot{\Sigma} \to \widehat{W}$ has
a well-defined \defin{relative homology class}, meaning the following.
Denote the asymptotic orbits of $u$ at its punctures $z \in \Gamma^\pm$
by $\gamma_z$, and let $\bar{\boldsymbol{\gamma}}^\pm \subset M_\pm$ denote
the closed $1$-dimensional submanifold defined as the union over 
$z \in \Gamma^\pm$ of the images of the orbits~$\gamma_z$.
Let $\overline{\Sigma}$ denote the compact oriented topological surface with
boundary obtained from $\dot{\Sigma}$ by appending $\{\pm\infty\} \times S^1$
to each of its cylindrical ends, and let $\pi : \widehat{W} \to W$ denote
the retraction defined as the identity on $W$ and $\pi(r,x) = x \in M_\pm 
\subset \p W$ for $(r,x)$ in $[0,\infty) \times M_+$ or 
$(-\infty,0] \times M_-$.  Then $\pi \circ u : \dot{\Sigma} \to W$ has a
natural continuous extension
$$
\bar{u} : (\overline{\Sigma},\p\overline{\Sigma}) \to (W,
\bar{\boldsymbol{\gamma}}^+ \cup \bar{\boldsymbol{\gamma}}^-)
$$
and thus represents a relative homology class
$$
[u] \in H_2(W,\bar{\boldsymbol{\gamma}}^+ \cup \bar{\boldsymbol{\gamma}}^-).
$$

\section{Moduli spaces of unparametrized holomorphic curves}
\label{sec:moduliSpaces}

We continue in the setting of a completed symplectic cobordism $\widehat{W}$
with fixed choices of $\psi \in \tT_0$, $r_0 \ge 0$ and
$J \in \jJ(\omega_\psi,r_0,\hH_+,\hH_-)$.  We shall denote by $\xi_\pm$ and
$R_\pm$ the hyperplane distribution and Reeb vector field respectively
determined by the stable Hamiltonian structure 
$\hH_\pm = (\omega_\pm,\lambda_\pm)$.

Fix integers $g, m, k_+,k_- \ge 0$ along with ordered sets of Reeb orbits
$$
\boldsymbol{\gamma}^\pm = (\gamma_1^\pm,\ldots,\gamma_{k_\pm}^\pm),
$$
where each $\gamma_i^\pm$ is a closed orbit of $R_\pm$ in~$M_\pm$.
Denote the union of the images of the $\gamma_i^\pm$ by
$\bar{\boldsymbol{\gamma}}^\pm \subset M_\pm$, and choose a relative homology
class
$$
A \in H_2(W, \bar{\boldsymbol{\gamma}}^+ \cup \bar{\boldsymbol{\gamma}}^-)
$$
whose image under the boundary map 
$H_2(W, \bar{\boldsymbol{\gamma}}^+ \cup \bar{\boldsymbol{\gamma}}^-) 
\stackrel{\p}{\longrightarrow}
H_1(\bar{\boldsymbol{\gamma}}^+ \cup \bar{\boldsymbol{\gamma}}^-)$ defined via
the long exact sequence of the pair $(W,\bar{\boldsymbol{\gamma}}^+ \cup \bar{\boldsymbol{\gamma}}^-)$
is
$$
\p A = \sum_{i=1}^{k_+} [\gamma_i^+] - \sum_{i=1}^{k_-} [\gamma_i^-] \in
H_1(\bar{\boldsymbol{\gamma}}^+ \cup \bar{\boldsymbol{\gamma}}^-).
$$
The \defin{moduli space} of \defin{unparametrized} $J$-holomorphic curves 
of \defin{genus~$g$} with \defin{$m$ marked points}, 
\defin{homologous to $A$} and 
\defin{asymptotic to $(\boldsymbol{\gamma}^+,\boldsymbol{\gamma}^-)$} is
then defined as a set of equivalence classes of tuples
$$
\mM_{g,m}(J,A,\boldsymbol{\gamma}^+,\boldsymbol{\gamma}^-) = 
\left\{ (\Sigma,j,\Gamma^+,\Gamma^-,\Theta,u) \right\} \big/ \sim,
$$
where:
\begin{enumerate}
\item
$(\Sigma,j)$ is a closed connected Riemann surface of genus~$g$;
\item
$\Gamma^+ = (z_1^+,\ldots,z_{k_+}^+)$, 
$\Gamma^- = (z_1^-,\ldots,z_{k_-}^-)$ 
and $\Theta = (\zeta_1,\ldots,\zeta_m)$ are disjoint ordered sets of 
distinct points in~$\Sigma$;
\item 
$u : (\dot{\Sigma} := \Sigma \setminus (\Gamma^+\cup \Gamma^-),j) \to
(\widehat{W},J)$ is an asymptotically cylindrical $J$-holomorphic map
with $[u] = A$, asymptotic at $z_i^\pm \in \Gamma^\pm$ to $\gamma_i^\pm$
for $i=1,\ldots,k_\pm$;
\item
Equivalence
$$
(\Sigma_0,j_0,\Gamma^+_0,\Gamma^-_0,\Theta_0,u_0) \sim
(\Sigma_1,j_1,\Gamma^+_1,\Gamma^-_1,\Theta_1,u_1)
$$
means the existence of a biholomorphic map $\psi : (\Sigma_0,j_0) \to
(\Sigma_1,j_1)$, taking $\Gamma_0^\pm$ to $\Gamma_1^\pm$ and
$\Theta_0$ to $\Theta_1$ with the ordering preserved, such that
$$
u_1 \circ \psi = u_0.
$$
\end{enumerate}

We shall usually abuse notation by abbreviating elements
$[(\Sigma,j,\Gamma^+,\Gamma^-,\Theta,u)]$ in this moduli space by
$$
u \in \mM_{g,m}(J,A,\boldsymbol{\gamma}^+,\boldsymbol{\gamma}^-).
$$
The \defin{automorphism group}
$$
\Aut(u) = \Aut(\Sigma,j,\Gamma^+,\Gamma^-,\Theta,u)
$$
of $u$ is defined as the group of biholomorphic maps
$\psi : (\Sigma,j) \to (\Sigma,j)$ which act as the identity on
$\Gamma^+ \cup \Gamma^- \cup \Theta$ and satisfy $u = u \circ \psi$.
Clearly the isomorphism class of this group depends only on the
equivalence class 
$[(\Sigma,j,\Gamma^+,\Gamma^-,\Theta,u)] \in 
\mM_{g,m}(J,A,\boldsymbol{\gamma}^+,\boldsymbol{\gamma}^-)$, and we will
see in \S\ref{sec:covers} below that it is always finite unless
$u : \dot{\Sigma} \to \widehat{W}$ is constant.  The significance of
the marked points is that they determine an \defin{evaluation map}
$$
\ev : \mM_{g,m}(J,A,\boldsymbol{\gamma}^+,\boldsymbol{\gamma}^-) \to
\widehat{W}^m : [(\Sigma,j,\Gamma^+,\Gamma^-,\Theta,u)] \mapsto
(u(\zeta_1),\ldots,u(\zeta_m))
$$
where $\Theta = (\zeta_1,\ldots,\zeta_m)$.  For most of our applications we
will be free to assume $m=0$, as marked points are not needed for defining
the most basic versions of SFT; the evaluation map does play a prominent
role however in more algebraically elaborate versions of the theory, and
especially in the Gromov-Witten invariants (the ``closed case'' of SFT).

We will assign a topology to
$\mM_{g,m}(J,A,\boldsymbol{\gamma}^+,\boldsymbol{\gamma}^-)$ in the next
lecture by locally identifying it with subsets of certain manifolds of
maps $\dot{\Sigma} \to \widehat{W}$ with Sobolev-type regularity
and exponential decay conditions at the ends.  In reality, this topology
admits a simpler description: one can define convergence of a sequence
$$
[(\Sigma_\nu,j_\nu,\Gamma^+_\nu,\Gamma^-_\nu,\Theta_\nu,u_\nu)] \to
[(\Sigma,j,\Gamma^+,\Gamma^-,\Theta,u)]
$$
to mean that for sufficiently large $\nu$, the equivalence classes in the
sequence admit representatives of the form
$(\Sigma,j_\nu',\Gamma^+,\Gamma^-,\Theta,u_\nu')$ such that
\begin{enumerate}
\item $j_\nu' \to j$ in $C^\infty$;
\item $u_\nu' \to u$ in $C^\infty_{\loc}(\dot{\Sigma},\widehat{W})$;
\item $\bar{u}_\nu' \to \bar{u}$ in $C^0(\overline{\Sigma},W)$.
\end{enumerate}
The proof that this topology matches what we will define in the next lecture
in terms of weighted Sobolev spaces requires asymptotic elliptic regularity
arguments along the lines of Proposition~\ref{prop:expDecay}.

\section{Simple curves and multiple covers}
\label{sec:covers}

In Lecture~\ref{lec:local}, we proved that closed $J$-holomorphic curves are all either
embedded in the complement of a finite set or are multiple covers of
curves with this property.  The same thing holds in the punctured case:

\begin{thm}
\label{thm:simple2}
Assume $u : (\dot{\Sigma},j) \to (\widehat{W},J)$ is a nonconstant
asymptotically cylindrical $J$-holomorphic curve whose asymptotic orbits
are all nondegenerate, where 
$\dot{\Sigma} = \Sigma \setminus \Gamma$ for some closed Riemann surface
$(\Sigma,j)$ and finite subset $\Gamma \subset \Sigma$.  Then there
exists a factorization $u = v \circ \varphi$, where
\begin{itemize}
\item $\varphi : (\Sigma,j) \to (\Sigma',j')$ is a holomorphic map
of positive degree to another closed and connected Riemann surface
$(\Sigma',j')$;
\item $v : (\dot{\Sigma}',j') \to (\widehat{W},J)$ is an asymptotically 
cylindrical $J$-holomorphic curve which is
embedded except at a finite set of critical points and self-intersections,
where $\dot{\Sigma}' := \Sigma' \setminus \Gamma'$ with
$\Gamma' := \varphi(\Gamma)$ and $\Gamma = \varphi^{-1}(\Gamma')$.
\end{itemize}
\end{thm}

As in the closed case, we call $u$ a \defin{simple} curve if the
holomorphic map $\varphi : (\Sigma,j) \to (\Sigma',j')$ is a diffeomorphism,
and $u$ is otherwise a \defin{$k$-fold multiple cover} of $v$ with
$k := \deg(\varphi) \ge 2$.

The proof of this theorem is an almost verbatim repeat of the proof of
Theorem~\ref{thm:simple} in Lecture~\ref{lec:local}, but with one new ingredient added.  Recall that
in the closed case, our proof required two lemmas which described the local
picture of a $J$-holomorphic curve $u : \dot{\Sigma} \to \widehat{W}$ near 
either a double point $u(z_0) = u(z_1)$ for $z_0 \ne z_1$
or a critical point $du(z_0) = 0$.  Both statements were completely local
and thus equally valid for non-closed curves, but we now need similar
statements to describe what kinds of singularities can appear in the
neighborhood of a puncture.  The following lemma is due to Siefring
\cite{Siefring:asymptotics} and follows from a ``relative asymptotic
formula'' analogous to Proposition~\ref{prop:expDecay}.

\begin{lemma}[Asymptotics]
\label{lemma:asymptotics}
Assume $u : (\dot{\Sigma} = \Sigma \setminus \Gamma,j) \to (\widehat{W},J)$ 
is asymptotically cylindrical and is asymptotic at $z_0 \in \Gamma$ to
a nondegenerate Reeb orbit.  Then a punctured neighborhood
$\dot{\uU}_{z_0} \subset \dot{\Sigma}$ of $z_0$ can be identified
biholomorphically with the punctured disk $\dot{\DD} = \DD \setminus \{0\}$
such that
$$
u(z) = v(z^k) \quad \text{ for } \quad
z \in \dot{\DD} = \dot{\uU}_{z_0},
$$
where $k \in \NN$ and $v : (\dot{\DD},i) \to (\widehat{W},J)$ is an embedded
and asymptotically cylindrical $J$-holomorphic curve.  Moreover, 
if $u' : (\dot{\Sigma}' = \Sigma' \setminus \Gamma',j') \to (\widehat{W},J)$
is another asymptotically cylindrical curve with a puncture
$z_0' \in \Gamma'$, then the images of $u$ near $z_0$ and $u'$ near $z_0'$
are either identical or disjoint.
\qed
\end{lemma}

\begin{exercise}
With Lemma~\ref{lemma:asymptotics} in hand, adapt the proof of
Theorem~\ref{thm:simple} in Lecture~\ref{lec:local} to prove Theorem~\ref{thm:simple2}.
If you get stuck, see \cite{Nelson:Abendblatt}*{\S 3.2}.
\end{exercise}

\begin{prop}
\label{prop:finite}
If $[(\Sigma,j,\Gamma^+,\Gamma^-,\Theta,u)] \in 
\mM_{g,m}(J,A,\boldsymbol{\gamma}^+,\boldsymbol{\gamma}^-)$ is represented
by a simple curve, then $\Aut(u)$ is trivial.  If it is represented by a
$k$-fold cover of a simple curve, then $|\Aut(u)| \le k$.
In particular, $\Aut(u)$ is always finite unless $u$ is constant.
\end{prop}
\begin{proof}
If $u$ is simple, then it is a diffeomorphism onto its image in a small 
neighbourhood of some point, and any map $\varphi$ satisfying 
$u=u \circ \varphi$ would be the identity on such a neighbourhood.
By unique continuation, we conclude that $\Aut(u)$ is trivial. 
In general if $u =v \circ \varphi$ for some simple 
$$v: \Sigma' \to W$$ and
$$\varphi: \Sigma \to \Sigma'$$ 
a $k$-fold branched cover, we have $$\Aut(u)=\{f:\Sigma \to \Sigma\ |\ v \circ \varphi \circ f = v \circ \varphi \}.$$
By a similar argument as in the previous case, knowing that $v$ is simple implies we only need to look at solutions to $$\varphi \circ f = \varphi.$$
Remove the set of branch points $B$ from $\Sigma'$ together with the set $\varphi^{-1}(B)$ from $\Sigma$,
so that $\varphi$ becomes an honest covering map. Any $\varphi \in \Aut(u)$ then defines a biholomorphic deck
transformation of the cover, so it remains to argue that there are at most $k$ of them. 
In fact, there is at most one transformation that takes $w_1$ to $w_2$ for any two given points
$w_1,w_2 \in \varphi^{-1}(x)$. If there were two such transformations $f$ and $g$, 
then $f \circ g^{-1}$ would be the identity on an open neighbourhood and would thus be globally
the identity by unique continuation.
\end{proof}

\section{A local structure result}

The following statement, which we will prove in the next lecture, is the
main goal of most of the analysis we have discussed recently.
It is essentially an application of the implicit function theorem for a
smooth nonlinear Fredholm section of a Banach space bundle.  The implicit
function theorem (see \cite{Lang:analysis}) implies in particular that if 
$F$ is a smooth map between Banach spaces such that $F(x_0) = 0$ and
$dF(x_0)$ is a surjective Fredholm operator, then $F^{-1}(0)$ is a smooth
manifold near $x_0$ with its dimension equal to the Fredholm index
of~$dF(x_0)$.  Surjectivity is an extra hypothesis, referred to in the
statement below as ``Fredholm regularity,'' a notion that we will define
precisely in the next lecture.  The dimension formula should look familiar,
but is only an \emph{indirect} consequence of the index formula for
Cauchy-Riemann type operators that we proved
in Lecture~\ref{lec:index}; one also needs to account for the fact that in defining our
moduli space 
$\mM_{g,m}(J,A,\boldsymbol{\gamma}^+,\boldsymbol{\gamma}^-)$, we did not
fix the complex structures on our domain curves, hence they are free to
move about in the moduli space of Riemann surfaces, whose dimension 
therefore plays a role in determining the dimension of
$\mM_{g,m}(J,A,\boldsymbol{\gamma}^+,\boldsymbol{\gamma}^-)$.

\begin{thm}
\label{thm:moduliDimension}
The set of \emph{Fredholm regular} curves forms an open subset
$$
\mM^\reg_{g,m}(J,A,\boldsymbol{\gamma}^+,\boldsymbol{\gamma}^-)
\subset \mM_{g,m}(J,A,\boldsymbol{\gamma}^+,\boldsymbol{\gamma}^-)
$$
which naturally
admits the structure of a smooth finite-dimensional orbifold of dimension
\begin{equation*}
\begin{split}
\dim \mM^\reg_{g,m}(J,A,\boldsymbol{\gamma}^+,\boldsymbol{\gamma}^-) &=
(n-3) (2 - 2g - k_+ - k_-) + 2 c_1^\tau(A) \\
& \quad + \sum_{i=1}^{k_+} \muCZ^\tau(\gamma_i^+) -
\sum_{i=1}^{k_-} \muCZ^\tau(\gamma_i^-) + 2m,
\end{split}
\end{equation*}
where $\dim W = 2n$, $\tau$ is a choice of unitary trivialization for
$(\xi_\pm,J,\omega_\pm)$ along each of the asymptotic orbits $\gamma_i^\pm$,
and $c_1^\tau(A)$ denotes the normal first Chern number of the 
complex vector bundle $(u^*T\widehat{W},J) \to \dot{\Sigma}$
with respect to the asymptotic trivialization determined by $\tau$ and
the splitting $T(\RR \times M_\pm) = \epsilon \oplus \xi_\pm$
(cf.~\eqref{eqn:splitting}).  The local isotropy group of
$\mM^\reg_{g,m}(J,A,\boldsymbol{\gamma}^+,\boldsymbol{\gamma}^-)$
at $u$ is $\Aut(u)$, hence the moduli space is a manifold near any
regular element with trivial automorphism group.
\end{thm}

\begin{exercise}
Verify that the number in the above index formula is independent of the
choice of trivializations~$\tau$, and that $c_1^\tau(u^*T\widehat{W})$
depends only on the relative homology class~$A$.
\end{exercise}

\psfrag{zero}{$0$}
\psfrag{one}{$1$}
\psfrag{param}{$s$}
\psfrag{F0}{$\mM(J_0)$}
\psfrag{F1}{$\mM(J_1)$}

\chapter{Smoothness of the moduli space}
\label{lec:transversality}

\minitoc

\vspace{12pt}

In this lecture, we continue the study of the moduli space
$$
\mM(J) := \mM_{g,m}(J,A,\boldsymbol{\gamma}^+,\boldsymbol{\gamma}^-).
$$
We assume as before that $(W,\omega)$ is a $2n$-dimensional symplectic
cobordism with stable boundary $\p W = -M_- \sqcup M_+$ inheriting stable
Hamiltonian structures $\hH_\pm = (\omega_\pm,\lambda_\pm)$ with induced 
Reeb vector fields $R_\pm$ and hyperplane distributions
$\xi_\pm = \ker \lambda_\pm$, $g , m , k_+,k_- \ge 0$ are integers,
$\boldsymbol{\gamma}^\pm = (\gamma_1^\pm,\ldots,\gamma_{k_\pm}^\pm)$ are ordered
sets of periodic $R_\pm$-orbits in $M_\pm$, and $A \in H_2(W,\bar{\boldsymbol{\gamma}}^+
\cup \bar{\boldsymbol{\gamma}}^-)$ is a relative homology class with
$\p A = \sum_i [\gamma_i^+] - \sum_i [\gamma_i^-] \in 
H_1(W,\bar{\boldsymbol{\gamma}}^+ \cup \bar{\boldsymbol{\gamma}}^-)$.  The noncompact
completion of $(W,\omega)$ is denoted by $(\widehat{W},\omega_\psi)$ for some
fixed function $\psi : \RR \to (-\epsilon,\epsilon)$ that scales the symplectic
form on the cylindrical ends, and $r_0 \ge 0$ is a fixed constant which determines
the size of the ends $[r_0,\infty) \times M_+$ and $(-\infty,-r_0] \times M_-$ on
which we require our almost complex structures $J \in \jJ(\omega_\psi,r_0,\hH_+,\hH_-)$
to be $\RR$-invariant.  The complement of these ends has closure
$$
W^{r_0} := \left([-r_0,0] \times M_-\right) \cup_{M_-} W \cup_{M_+}
\left([0,r_0] \times M_+ \right).
$$
We will often make use of the fact that since $J$ matches translation-invariant
almost complex structures in $\jJ(\hH_\pm)$ outside of $W^{r_0}$, there are
natural complex vector bundle splittings
$$
T(\RR \times M_\pm) = \epsilon \oplus \xi_\pm,
$$
where $\epsilon$ denotes the canonically trivial line bundle spanned by $\p_r$ and
the Reeb vector field.

\section{Transversality theorems in cobordisms}

We concluded the previous lecture with the statement of the following theorem.

\begin{thm}
\label{thm:moduliDimension2}
If the orbits $\gamma_i^\pm$ are all nondegenerate and
$J \in \jJ(\omega_\psi,r_0,\hH_+,\hH_-)$, then
the moduli space $\mM(J)$ contains an open subset
$$
\mM^\reg(J) \subset \mM(J)
$$
consisting of so-called \emph{Fredholm regular} curves, which naturally
admits the structure of a smooth finite-dimensional orbifold of dimension
\begin{equation*}
\begin{split}
\dim \mM^\reg(J) &= (n-3) (2 - 2g - k_+ - k_-) + 2 c_1^\tau(A) \\ &\quad
+ \sum_{i=1}^{k_+} \muCZ^\tau(\gamma_i^+) - \sum_{i=1}^{k_-} \muCZ^\tau(\gamma_i^-) + 2m,
\end{split}
\end{equation*}
where $\dim W = 2n$, $\tau$ is a choice of unitary trivialization for
$(\xi_\pm,J,\omega_\pm)$ along each of the asymptotic orbits $\gamma_i^\pm$,
and $c_1^\tau(A)$ denotes the normal first Chern number of the 
complex vector bundle $(u^*T\widehat{W},J) \to \dot{\Sigma}$
with respect to the asymptotic trivialization determined by $\tau$ and
the splitting $T(\RR \times M_\pm) = \epsilon \oplus \xi_\pm$.
The local isotropy group of
$\mM^\reg(J)$ at $u$ is $\Aut(u)$, hence the moduli space is a manifold near any
regular element with trivial automorphism group.
\end{thm}

The integer in the above dimension formula is often called the 
\defin{virtual dimension} of  $\mM(J)$ and denoted by
\begin{equation*}
\begin{split}
\virdim \mM(J) &:= (n-3) (2 - 2g - k_+ - k_-) + 2 c_1^\tau(A) \\
& \quad + \sum_{i=1}^{k_+} \muCZ^\tau(\gamma_i^+) 
- \sum_{i=1}^{k_-} \muCZ^\tau(\gamma_i^-) + 2m.
\end{split}
\end{equation*}
Ignoring the marked points, the virtual dimension of a space
$\mM_{g,0}(J,A,\boldsymbol{\gamma}^+,\boldsymbol{\gamma}^-)$ containing
a curve $u : (\dot{\Sigma},j) \to (\widehat{W},J)$ with punctures 
$z \in \Gamma^\pm$ and nondegenerate asymptotic orbits
$\{ \gamma_z \}_{z \in \Gamma^\pm}$ is sometimes also called 
the \defin{index} of~$u$,
$$
\ind(u) := (n-3) \chi(\dot{\Sigma}) + 2 c_1^\tau(u^*T\widehat{W}) +
\sum_{z \in \Gamma^+} \muCZ^\tau(\gamma_z) - 
\sum_{z \in \Gamma^-} \muCZ^\tau(\gamma_z) \in \ZZ,
$$
and we will see that it is in fact the Fredholm index of an operator
closely related to the linearized Cauchy-Riemann operator
$\mathbf{D}_u$ at~$u$.
The word ``virtual'' refers to the fact that in general, the regularity
condition may fail and thus 
$\mM(J)$ might not
be smooth, or if it is, it might actually be of a different dimension
(see Example~\ref{ex:noTransversality} below), but in an ideal world where
transversality is always satisfied, its dimension would be
$\virdim \mM(J)$.
This notion makes sense in finite-dimensional contexts as well: if
$f : \RR^n \to \RR^m$ is a smooth map, then we would say that $f^{-1}(0)$
has virtual dimension $n-m$, even though $f^{-1}(0)$ might in general be
all sorts of strange things other than a smooth $(n-m)$-dimensional manifold.
In particular, $n-m$ could be negative, in which case $f^{-1}(0)$ would be
empty if transversality were satisfied, but in general this need not
be the case.  It is true however that $f$ can always be \emph{perturbed}
to a map whose zero set is an $(n-m)$-dimensional manifold (or empty if
$n - m < 0$).  The same is true in principle of the nonlinear Cauchy-Riemann equation,
but in general it is a formidably difficult problem to find perturbations
that respect all symmetries inherent in the setup as well as the extra
structure provided by the \emph{compatification} of $\mM(J)$, which is
usually crucial for meangingful applications.
Such issues require more sophisticated methods than we will discuss here,
but a good place to read about them is \cite{FabertFishGolovkoWehrheim}.

The first goal of this lecture is to define the notion ``Fredholm regular''
and prove Theorem~\ref{thm:moduliDimension2}.  In practice, however,
Fredholm regularity is a technical condition that can rarely be directly
checked.  To remedy this, we will also prove a genericity result
for \emph{somewhere injective} $J$-holomorphic curves.  A smooth map
$u : \dot{\Sigma} \to \widehat{W}$ is said to have an \defin{injective point}
$z \in \dot{\Sigma}$ if
$$
du(z) : T_z \dot{\Sigma} \to T_{u(z)}\widehat{W} \text{ is injective}
\quad \text{ and } \quad
u^{-1}(u(z)) = \{z\}.
$$
If $u$ is a proper map, then it is easy to see that the set of injective
points is open in $\dot{\Sigma}$, though in general it could also be empty;
this is the case e.g.~for multiply covered $J$-holomorphic curves.
We say $u$ is \defin{somewhere injective} if its set of injective points is
nonempty; for asymptotically cylindrical $J$-holomorphic curves with nondegenerate
asymptotic orbits, Theorem~\ref{thm:simple2} implies that somewhere injectivity
is equivalent to being \emph{simple}, i.e.~not multiply covered.

Recall that if $X$ is a topological space, a subset $Y \subset X$ is called
\defin{comeager} if it contains a countable intersection of open and dense sets.\footnote{Elsewhere
in the symplectic literature, comeager subsets are sometimes referred to as
``sets of second category,'' which is unfortunately slightly at odds with the
standard meaning of ``second category,'' though it is accurate to say that
the \emph{complement} of a comeager subset (also known as a ``meager'' subset)
is a set of first category.  The term \emph{Baire subset} is also sometimes
used as a synonym for ``comeager subset''.}
If $X$ is complete, then the Baire category theorem implies that comeager
subsets are always dense; moreover, any countable intersection of comeager
subsets is also comeager and therefore dense.  Comeager subsets often play
the role in infinite dimensions that the term ``almost everywhere'' plays
in finite dimensions.  Informally, we often say that a given statement 
dependent on a choice of auxiliary data (living in a complete metric space) is
true \defin{generically}, or ``for \defin{generic} choices,'' if it is true
whenever the data are chosen from some comeager subset of the space of 
all possible data.

\begin{thm}
\label{thm:genericCobordism}
Fix the same data as in Theorem~\ref{thm:moduliDimension2}, an almost complex
structure $\Jfix \in \jJ(\omega_\psi,r_0,\hH_+,\hH_-)$ and an open subset
$$
\uU \subset W^{r_0}.
$$
Then there exists a comeager subset
$$
\jJ_\uU^\reg \subset \left\{ J \in \jJ(\omega_\psi,r_0,\hH_+,\hH_-) \ \big|\ 
\text{$J = \Jfix$ on~$\widehat{W} \setminus \uU$} \right\},
$$
such that for every $J \in \jJ_\uU^\reg$, every curve $u \in \mM(J)$
that has an injective point mapped into $\uU$ is Fredholm regular.
In particular, the curves with this property define an open subset of
$\mM(J)$ that is a smooth manifold with dimension equal to its virtual
dimension.
\end{thm}

\begin{remark}
\label{remark:whatever}
Since $\uU \subset \widehat{W}$ has compact closure, the set 
$$
\left\{ J \in \jJ(\omega_\psi,r_0,\hH_+,\hH_-) \ \big|\ 
\text{$J = \Jfix$ on~$\widehat{W} \setminus \uU$} \right\}
$$
has a natural $C^\infty$-topology that
makes it a Fr\'echet manifold and thus a complete metric space, hence
comeager subsets of it are dense.
\end{remark}

\begin{figure}
\includegraphics{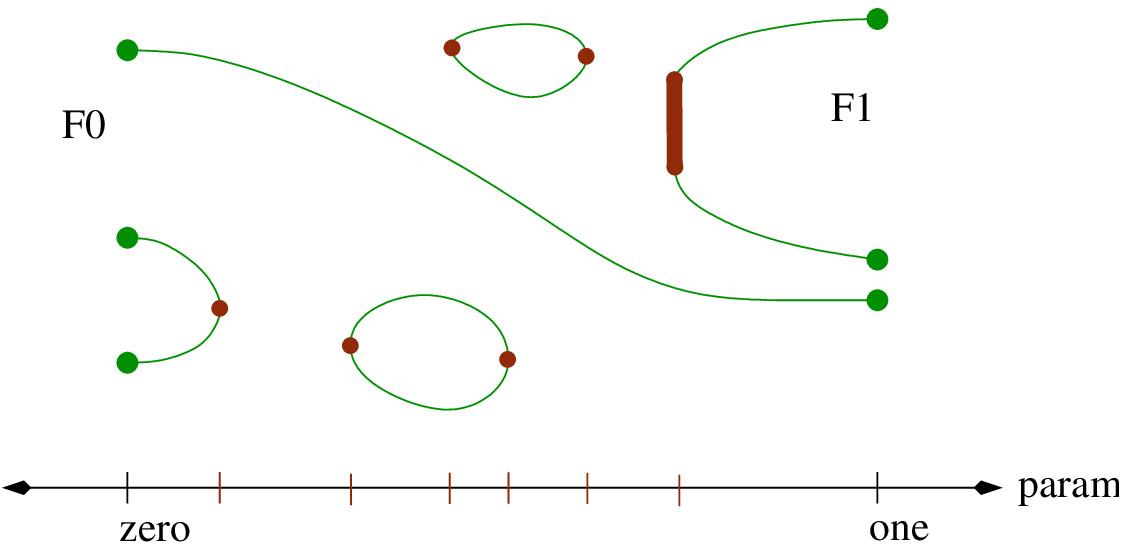}
\caption{\label{fig:parametric} The picture shows a smooth parametric
moduli space $\mM(\{J_s\}_{s \in [0,1]})$ and its projection
$\mM(\{J_s\}) \to [0,1] : (s,u) \mapsto s$ in a case where
$\virdim \mM(J_s) = 0$.  The parametric moduli space is
$1$-dimensional and the spaces $\mM(J_s)$ are regular and
$0$-dimensional for almost every $s \in [0,1]$, but
this need not hold when $s$ is a critical value of the projection; 
in the picture, one such space $\mM(J_s)$ contains a
$1$-dimensional component consisting of non-regular curves, so its dimension 
differs from its virtual dimension.}
\end{figure}

\begin{remark}
\label{remark:parametric}
Both of the above theorems admit easy extensions to the study of moduli
spaces dependent on finitely many parameters.  Concretely, suppose
$P$ is a smooth finite-dimensional manifold and $\{ J_s \}_{s \in P}$ is
a smooth family of almost complex structures satisfying the usual conditions.
One can then define a \emph{parametric moduli space}
$$
\mM(\{J_s\}_{s \in P}) = \left\{ (s,u)\ \big| \ s \in P, \ 
u \in \mM(J_s) \right\}
$$
and a notion of \emph{parametric regularity} for pairs
$(s,u) \in \mM(\{J_s\})$, which is again an open condition, such that the
space $\mM^\reg(\{J_s\})$ of parametrically regular elements will be an orbifold of dimension
$$
\dim \mM^\reg(\{J_s\}) = \virdim \mM(J) + \dim P.
$$
Similarly, one can show that if the family
$\{J_s\}_{s \in P}$ is allowed to vary on an open subset
$\uU \subset W^{r_0}$ for $s$ lying in some precompact open subset
$\vV \subset P$, then all elements $(s,u)$ for which $s \in \vV$ and 
$u$ has an injective point mapping to $\uU$ will be parametrically regular.
See \cite{Wendl:lecturesV33}*{\S 4.5} for details in the closed case, which
is not fundamentally different from the punctured case.
The standard and most important example is $P = [0,1]$ with $\vV = (0,1)$, 
so we consider \emph{generic homotopies} of almost complex structures.  
Here it is important
to observe that while regularity in the sense of Theorem~\ref{thm:moduliDimension2}
always implies parametric regularity, the converse is false: there can exist
parametrically regular pairs $(s,u) \in \mM(\{J_s\})$ for which $u$ is \emph{not}
a Fredholm regular element of~$\mM(J_s)$, hence $\mM(\{J_s\})$ may be smooth
even if $\mM(J_s)$ is not smooth for some $s \in P$.  This can happen
in particular whenever $s$ is a critical value of the projection map
$$
\mM(\{J_s\}) \to P : (s,u) \mapsto s,
$$
see Figure~\ref{fig:parametric}.
In general these cannot be excluded by making generic choices of the homotopy, 
though it is possible in certain cases
using ``automatic'' transversality results, which guarantee regularity
for all $J_s$ with no need for genericity (cf.~\cite{Wendl:automatic}).
\end{remark}

\begin{example}
\label{ex:noTransversality}
It is not hard to imagine situations in which transversality \emph{must}
fail generically for multiply covered curves.  Suppose for instance that
$(W,\omega)$ is an $8$-dimensional symplectic manifold with compatible
almost complex structure~$J_0$, and $u_0 : S^2 \to W$ is a simple
$J_0$-holomorphic sphere with no punctures and $[u_0] = A \in H_2(W)$,
where $c_1(A) = -1$.  This means $u_0$ represents an element of a
moduli space $\mM_{0,0}(J_0,A)$ with
$$
\virdim \mM_{0,0}(J_0,A) = 2 - 2g + 2 c_1(A) = 0.
$$
In particular if $u_0$ is regular and $\{J_s \in \jJ(\omega) \}_{s \in \RR^k}$
is a smooth $k$-parameter family of compatible almost complex structures 
including~$J_0$, then Remark~\ref{remark:parametric} implies that a
neighborhood of $(0,u_0)$ in the parametric moduli space
$\mM(\{J_s\}) = \{ (s,u)\ |\ s \in P,\ u \in \mM_{0,0}(J_s,A) \}$ is a 
smooth $k$-dimensional manifold, and this will be true no matter how the
family $\{J_s\}$ is chosen.  But for each of the elements
$(s,u) \in \mM(\{J_s\})$ parametrized by a $J$-holomorphic map
$u : (S^2 = \CC \cup \{\infty\},i) \to (W,J_s)$, there is also a double cover
$$
u' : S^2 \to W : z \mapsto u(z^2),
$$
with $[u'] = 2A$, so $u' \in \mM_{0,0}(J_s,2A)$ and
$$
\virdim \mM_{0,0}(J_s,2A) = 2 - 2g + 2 c_1(2A) = -2.
$$
Negative virtual dimension means that $\mM_{0,0}(J_0,2A)$ should be empty
whenever Fredholm regularity is achieved, but this is clearly impossible,
even generically, since elements of $\mM_{0,0}(J_s,A)$ always have double
covers belonging to $\mM_{0,0}(J_s,2A)$.
\end{example}

\begin{remark}
The most common way to apply Theorem~\ref{thm:genericCobordism} is by setting
$\uU$ equal to the interior of~$W^{r_0}$,
so generic perturbations of $J$ are allowed everywhere except on the
regions where it is required to be $\RR$-invariant.  The theorem then
achieves transversality for all simple curves that are not confined to the
$\RR$-invariant regions.  We will show in the next lecture that
transversality for all curves of the latter type can also be achieved by
generic perturbations within the spaces $\jJ(\hH_\pm)$ of compatible
$\RR$-invariant almost complex structures on the symplectizations 
$\RR \times M_\pm$, hence generic choices in $\jJ(\omega_\psi,r_0,\hH_+,\hH_-)$
do achieve transversality for all simple curves.
\end{remark}

Our proofs of Theorems~\ref{thm:moduliDimension2} and~\ref{thm:genericCobordism}
will mostly follow the same line of argument that is carried out for the
closed case in \cite{Wendl:lecturesV33}*{Chapter~4}, thus we will not
discuss every detail but will instead emphasize aspects which are unique to
the punctured case.

\section{Functional analytic setup}
\label{sec:setup2}

Fix $k \in \NN$ and $p \in (1,\infty)$ with $kp > 2$, a small
number $\delta \ge 0$, and a Riemannian metric on $\widehat{W}$ that
is translation-invariant in the cylindrical ends.  Fix also
a closed connected surface $\Sigma$ of genus~$g$, and disjoint finite 
ordered sets of distinct points
$$
\Gamma^\pm = (z_1^\pm,\ldots,z_{k_\pm}^\pm), \qquad
\Theta = (\zeta_1,\ldots,\zeta_m)
$$
in $\Sigma$, together with disjoint neighborhoods 
$$
\uU_j^\pm \subset \Sigma
$$
of each $z_j^\pm \in \Gamma^\pm$ with complex structures
$j_\Gamma$ and biholomorphic identifications of $(\uU_j^\pm,j_\Gamma,z_j)$
with $(\DD,i,0)$ for each $j=1,\ldots,k_\pm$.  This determines holomorphic
cylindrical coordinates identifying each of the punctured neighborhoods
$$
\dot{\uU}_j^\pm \subset \dot{\Sigma} := \Sigma \setminus (\Gamma^+ \cup \Gamma^-)
$$
biholomorphically with the half-cylinder~$Z_\pm$.

For reasons that will become clear when we study the linearized Cauchy-Riemann
operator in the punctured setting, we will need to consider exponentially
weighted Sobolev spaces.  Suppose $E \to \dot{\Sigma}$ is an asymptotically
Hermitian vector bundle: then the Banach space
$$
W^{k,p,\delta}(E) \subset W^{k,p}_{\loc}(E)
$$
is defined to consist of sections $\eta \in W^{k,p}_{\loc}(E)$ whose 
representatives $f : Z_\pm \to \CC^m$ in cylindrical coordinates $(s,t) \in Z_\pm$
and asymptotic trivializations at the ends satisfy
\begin{equation}
\label{eqn:weights}
\| e^{\pm \delta s} f \|_{W^{k,p}(Z_\pm)} < \infty.
\end{equation}
The norm of a section $\eta \in W^{k,p,\delta}(E)$ is defined by adding 
the $W^{k,p}$-norm of $\eta$ over a large compact subdomain in $\dot{\Sigma}$
to the weighted norms \eqref{eqn:weights} for each cylindrical end.
If $\delta=0$, this just produces the usual $W^{k,p}(E)$, but for $\delta > 0$,
sections in $W^{k,p,\delta}(E)$ are guaranteed to have exponential decay
at infinity.

\begin{remark}
\label{remark:negativeExp}
It is occasionally useful to observe that the definition of
$W^{k,p,\delta}(E)$ also makes sense when $\delta < 0$.  In this case,
sections in $W^{k,p,\delta}(E)$ are of class $W^{k,p}_{\loc}$ but need not be
globally in~$W^{k,p}(E)$, as they are also allowed to have exponential
\emph{growth} at infinity.
\end{remark}

We now want to define a Banach manifold of maps $u : \dot{\Sigma} \to \widehat{W}$
that will contain all the asymptotically cylindrical $J$-holomorphic
curves with our particular choice of asymptotic orbits.  Recall that the
asymptotically cylindrical condition means
\begin{equation}
\label{eqn:asympCyl}
u(s,t) = \exp_{(T_j^\pm s,\gamma_j^\pm(t))} h(s,t) \quad \text{ for
sufficiently large~$|s|$}
\end{equation}
in suitable cylindrical coordinates $(s,t) \in Z_\pm$ near each
puncture $z_j^\pm \in \Gamma^\pm$, where $T_j^\pm > 0$ is the period of
the orbit $\gamma_j^\pm : S^1 \to M_\pm$ and $h(s,t)$ is a vector field
along the trivial cylinder that decays as $s \to \pm\infty$.  The catch is
that this definition was not formulated with respect to a \emph{fixed}
choice of the holomorphic cylindrical coordinates $(s,t)$; in general
the coordinates in which \eqref{eqn:asympCyl} is valid may
depend on~$u$, and different choices of coordinates might be required for
different maps.  One can show however that any two distinct choices of
holomorphic cylindrical coordinates are related to each other
by a transformation
that converges asymptotically to a constant shift, which implies that for our \emph{fixed}
choice of coordinates $(s,t)$, every asymptotically cylindrical map can
be assumed to satisfy
$$
u(s,t) = \exp_{(T_j^\pm s + a,\gamma_j^\pm(t + b))} h(s,t), \qquad
\lim_{s \to \pm\infty} h(s,t) = 0
$$
for some constants $a \in \RR$ and $b \in S^1$.
We therefore define the space
$$
\bB^{k,p,\delta} := W^{k,p,\delta}(\dot{\Sigma},\widehat{W} \,;\,
\boldsymbol{\gamma}^+,\boldsymbol{\gamma}^-) \subset C^0(\dot{\Sigma},\widehat{W})
$$
to consist of all continuous maps $u : \dot{\Sigma} \to \widehat{W}$ of the form
$$
u = \exp_f h,
$$
where:
\begin{itemize}
\item $f : \dot{\Sigma} \to \widehat{W}$ is smooth and, in our fixed cylindrical
coordinates $(s,t) \in Z_\pm$ on neighborhoods of the punctures $z_j^\pm \in
\Gamma^\pm$, takes the form
$$
f(s,t) = (T_j^\pm s + a , \gamma_j^\pm(t + b)) \quad \text{ for 
$|s|$ sufficiently large,}
$$
where $a \in \RR$ and $b \in S^1$ are arbitrary constants and
$T_j^\pm > 0$ is the period of the Reeb orbit $\gamma_j^\pm : S^1 \to M_\pm$;
\item $h \in W^{k,p,\delta}(f^*T\widehat{W})$.
\end{itemize}

Though it is not immediate since $\dot{\Sigma}$ is noncompact, 
one can generalize the ideas in \cite{Eliasson} to give $\bB^{k,p,\delta}$
the structure of a smooth, separable and metrizable Banach manifold.  The
key point is the condition $kp > 2$, which guarantees the continuous
inclusion $W^{k,p,\delta}(f^*T\widehat{W}) \hookrightarrow
C^0(f^*T\widehat{W})$ as well as Banach algebra and $C^k$-continuity
properties, cf.~Propositions~2.4, 2.7 and~2.8 in Lecture~\ref{lec:local}.  These properties
are needed in order to show that the transition maps between pairs of
charts of the form $\exp_f h \mapsto h$ are smooth.

The tangent space to $\bB^{k,p,\delta}$ at $u \in \bB^{k,p,\delta}$ can be 
written as
$$
T_u \bB^{k,p,\delta} = W^{k,p,\delta}(u^*T\widehat{W}) \oplus V_\Gamma,
$$
where $V_\Gamma \subset \Gamma(u^*T\widehat{W})$ is a non-canonical choice
of a $2(k_+ + k_-)$-dimensional vector space of smooth sections asymptotic
at the punctures to constant linear combinations of the vector fields spanning 
the canonical trivialization of the first factor in
$T(\RR \times M_\pm) = \epsilon \oplus \xi_\pm$, i.e.~they point in the
$\RR$- and $R_\pm$-directions.  The space $V_\Gamma$ appears due to the
fact that two distinct elements of $\bB^{k,p,\delta}$ are generally 
asymptotic to collections of trivial cylinders that differ from each other 
by $k_+ + k_-$ pairs of constant shifts $(a,b) \in \RR \times S^1$.

Fix $J \in \jJ(\omega_\psi,r_0,\hH_+,\hH_-)$ and a smooth complex
structure $j$ on $\Sigma$ that matches $j_\Gamma$ in the 
neighborhoods~$\uU_j^\pm$ of the punctures.
The nonlinear Cauchy-Riemann operator is then defined as a smooth section
$$
\dbar_{j,J} : \bB^{k,p,\delta} \to \eE^{k-1,p,\delta} : u \mapsto 
Tu + J \circ Tu \circ j
$$
of a Banach space bundle
$$
\eE^{k-1,p,\delta} \to \bB^{k,p,\delta}
$$
with fibers
$$
\eE^{k-1,p,\delta}_u = W^{k-1,p,\delta}(\overline{\Hom}_\CC(T\dot{\Sigma},u^*T\widehat{W})).
$$
The zero set of $\dbar_{j,J}$ is the set of all maps $u \in \bB^{k,p,\delta}$
that are pseudoholomorphic from $(\dot{\Sigma},j)$ to $(\widehat{W},J)$.
Note that the smoothness of $\dbar_{j,J}$ depends mainly on the fact that
$J$ is smooth.  Indeed, in local coordinates $\dbar_{j,J}$ looks like
$u \mapsto \p_s u + (J \circ u) \p_t u$, in which the most obviously nonlinear
ingredient is $u \mapsto J \circ u$.  If $J$ were only of class $C^k$, then the
$C^k$-continuity property would imply that the map
$u \mapsto J \circ u$
sends maps of class $W^{k,p}$ continuously to maps of class $W^{k,p}$,
and one can use an inductive argument to show that this map then becomes
$r$-times differentiable if $J$ is of class $C^{k+r}$,
see \cite{Wendl:lecturesV33}*{Lemma~2.12.5}.  Moreover, the fact that
$\dbar_{j,J} u$ satisfies the same exponential weighting condition as $u$
at the cylindrical ends depends on the fact that $J$ is $\RR$-invariant
near infinity.

For $u \in \dbar_{j,J}^{-1}(0)$, the linearization 
$D\dbar_{j,J}(u) : T_u \bB^{k,p,\delta} \to \eE_u^{k-1,p,\delta}$
defines a bounded linear operator
$$
\mathbf{D}_u : W^{k,p,\delta}(u^*T\widehat{W}) \oplus V_\Gamma \to
W^{k-1,p,\delta}(\overline{\Hom}_\CC(T\dot{\Sigma},u^*T\widehat{W})).
$$
We derived a formula for this operator in Lecture~\ref{lec:local} and showed that it is
of Cauchy-Riemann type.  Since $V_\Gamma$ is finite dimensional,
$\mathbf{D}_u$ will be Fredholm if and only if its restriction to the
first factor is Fredholm; denote this restriction by
$$
\mathbf{D}_\delta : W^{k,p,\delta}(u^*T\widehat{W}) \to
W^{k-1,p,\delta}(\overline{\Hom}_\CC(T\dot{\Sigma},u^*T\widehat{W})),
$$
where we've chosen the notation to emphasize the dependence of this
operator on the choice of exponential weight $\delta \ge 0$ in the
definition of our Banach space.  We will see presently why it's important
to pay attention to this detail.

To see whether $\mathbf{D}_\delta$ is Fredholm, consider first the special
case where $u$ is a trivial cylinder
$$
u_\gamma : \RR \times S^1 \to \RR \times M : (s,t) \mapsto (Ts,\gamma(t))
$$
over some Reeb orbit $\gamma : S^1 \to M$ with period $T > 0$ in
$M$ with stable Hamiltonian structure~$\hH = (\omega,\lambda)$ on~$M$.  
In this case, there is a
more convenient way to write down $\mathbf{D}_{u_\gamma}$ than the formula
from Lecture~\ref{lec:local}.  To start with, we use the splitting $T(\RR \times M)
= \epsilon \oplus \xi$ to decompose $u_\gamma^*T(\RR \times M) = u_\gamma^*\epsilon \oplus
u_\gamma^*\xi$ and thus write $\mathbf{D}_{u_\gamma}$ in block form
$$
\mathbf{D}_{u_\gamma} = \begin{pmatrix} \mathbf{D}_{u_\gamma}^\epsilon & \mathbf{D}_{u_\gamma}^{\epsilon\xi} \\
\mathbf{D}_{u_\gamma}^{\xi\epsilon} & \mathbf{D}_{u_\gamma}^\xi \end{pmatrix}.
$$

\begin{exercise}
\label{EX:split}
Suppose $\mathbf{D} : \Gamma(E) \to \Omega^{0,1}(\dot{\Sigma},E)$ is a
linear Cauchy-Riemann type operator on a vector bundle $E$ with a
complex-linear splitting $E = E_1 \oplus E_2$, and
$$
\mathbf{D} = \begin{pmatrix}
\mathbf{D}_{11} & \mathbf{D}_{12} \\
\mathbf{D}_{21} & \mathbf{D}_{22}
\end{pmatrix}
$$
is the resulting block decomposition of~$\mathbf{D}$.  Use the Leibniz
rule satisfied by $\mathbf{D}$ to show that
$\mathbf{D}_{11}$ and $\mathbf{D}_{22}$ are also Cauchy-Riemann type operators
on $E_1$ and $E_2$ respectively, while the off-diagonal terms are tensorial,
i.e.~they commute with multiplication by smooth real-valued functions and
thus define bundle maps $\mathbf{D}_{12} : E_2 \to \Lambda^{0,1}T^*\dot{\Sigma} \otimes E_1$
and $\mathbf{D}_{21} : E_1 \to \Lambda^{0,1}T^*\dot{\Sigma} \otimes E_2$.
\end{exercise}

Now observe that if $u = (u_\RR,u_M) : \RR \times S^1 \to \RR \times M$ is 
another cylinder near~$u_\gamma$, the nonlinear operator
$(\dbar_{j,J}u) \p_s = \p_s u + J\, \p_t u \in \Gamma(u^*T(\RR \times M)) = \Gamma(u^*\epsilon \oplus u^*\xi)$
takes the form
$$
(\dbar_{j,J} u) \p_s = \begin{pmatrix}
\p_s u_\RR - \lambda(\p_t u_M) + i \left( \p_t u_\RR + \lambda(\p_t u_M) \right) \\
\pi_\xi\, \p_s u_M + J \pi_\xi \, \p_t u_M
\end{pmatrix},
$$
where we are using the canonical trivialization of $u^*\epsilon$ via $\p_r$ and $R$
to express the top block as a complex-valued function.  As we observed in
Lecture~\ref{lec:asymptotic}, the bottom block of this expression can be interpreted in terms
of the gradient flow of an action functional $\aA_\hH : C^\infty(S^1) \to \RR$,
with $\nabla \aA_\hH(\gamma) = - J \pi_\xi \, \p_t\gamma$.
Linearizing in the direction of a section $\eta^\xi \in \Gamma(u_\gamma^*\xi)$
and taking the $\xi$ component thus yields an expression involving the
Hessian of $\aA_\hH$ at the critical point~$\gamma$, namely
$$
(\mathbf{D}_{u_\gamma}^\xi \eta^\xi) \p_s = (\p_s - \mathbf{A}_\gamma) \eta^\xi.
$$
To compute the blocks $\mathbf{D}_{u_\gamma}^\epsilon$ and 
$\mathbf{D}_{u_\gamma}^{\xi\epsilon}$, notice that $\mathbf{D}_{u_\gamma} \eta^\epsilon = 0$
whenever $\eta^\epsilon$ is a constant linear combination of $\p_r$ and $R$,
as $\eta^\epsilon$ is then the derivative of a smooth family of
$J$-holomorphic reparametrizations of~$u_\gamma$.  This is enough to prove
$\mathbf{D}_{u_\gamma}^{\xi\epsilon} = 0$ since the latter is tensorial
by Exercise~\ref{EX:split},
and expressing arbitrary sections of $u_\gamma^*\epsilon$ as $f\p_r + gR$,
we can apply the Leibniz rule for $\mathbf{D}_{u_\gamma}^\epsilon$ and conclude
$$
(\mathbf{D}_{u_\gamma}^{\epsilon} \eta^\epsilon) \p_s = (\p_s + i \, \p_t) \eta^\epsilon
$$
in the canonical trivialization.  To compute the remaining off-diagonal term,
one needs to compute $dr(\mathbf{D}_{u_\gamma} \eta^\xi)$ and
$\lambda(\mathbf{D}_{u_\gamma} \eta^\xi)$ for an arbitrary section
$\eta^\xi \in \Gamma(u_\gamma^*\xi)$, e.g.~by picking a smooth family 
$u_\rho : \RR \times S^1 \to \RR \times M$ with $\p_\rho u_\rho|_{\rho=0}
= \eta^\xi$ and a connection $\nabla$ and computing
$$
dr\big( \left.\nabla_\rho (\dbar_{j,J} u_\rho) \right|_{\rho=0} \big)
\quad\text{ and }\quad
\lambda\big( \left.\nabla_\rho( \dbar_{j,J} u_\rho )\right|_{\rho=0} \big).
$$
This calculation is straightforward but unenlightening, so I will leave it as
an exercise for now---in the next lecture we'll derive a general formula
(see Lemma~\ref{lemma:thisLemma}), which implies that since $\pi_\xi \, \p_s u_\gamma \equiv
\pi_\xi \, \p_t u_\gamma \equiv 0$ in the present setting,
$\mathbf{D}_{u_\gamma}^{\epsilon\xi} = 0$.  All this leads to the formula
$$
(\mathbf{D}_{u_\gamma} \eta) \p_s = \left( \p_s - \begin{pmatrix} -i \p_t & 0 \\ 0 & \mathbf{A}_\gamma \end{pmatrix} \right) \eta.
$$
Here the upper left block is the ``trivial'' asymptotic operator acting on the
trivial line bundle over~$S^1$.  Since every asymptotically 
cylindrical curve approximates a trivial cylinder near infinity, one can
deduce from this calculuation the following:

\begin{prop}
\label{prop:asympOperators}
The Cauchy-Riemann type operator $\mathbf{D}_u$ on $u^*T\widehat{W}$
is asymptotic at its punctures $z_j^\pm \in \Gamma^\pm$ for $j=1,\ldots,k_\pm$
to the asymptotic operators $(-i\p_t) \oplus \mathbf{A}_{\gamma_j^\pm}$ on
$(\gamma_j^\pm)^*(\epsilon \oplus \xi_\pm)$.
\end{prop}

Perhaps you can now see a problem: even if the orbits $\gamma_j^\pm$
are all nondegenerate, the asymptotic operators 
$(-i\p_t) \oplus \mathbf{A}_\gamma$ are degenerate, as they have nontrivial
kernel consisting of constant sections in the first (trivial) factor of
$(\gamma_j^\pm)^*(\epsilon \oplus \xi_\pm)$.  This implies in particular that
$$
\mathbf{D}_0 : W^{k,p}(u^*T\widehat{W}) \to W^{k-1,p}(\overline{\Hom}_\CC(T\dot{\Sigma},u^*T\widehat{W}))
$$
is \emph{not} Fredholm, except of course in the special case where there
are no punctures.  

The situation is saved by the exponential weight:

\begin{lemma}
\label{lemma:expRescue}
For every $\delta > 0$ sufficiently small, the operator
$\mathbf{D}_\delta$ is Fredholm and has index
$$
\ind(\mathbf{D}_\delta) = n \chi(\Sigma) - (n+1)\#\Gamma +
2 c_1^\tau(u^*T\widehat{W}) + \sum_{j=1}^{k_+} \muCZ^\tau(\gamma_j^+) -
\sum_{j=1}^{k_-} \muCZ^\tau(\gamma_j^-).
$$
Moreover, every element of
$\mM(J)$ can be represented by a map $u \in \bB^{k,p,\delta}$.
\end{lemma}
\begin{proof}
The second claim follows from the exponential decay estimate of
Hofer-Wysocki-Zehnder \cite{HWZ:props1} mentioned in the previous lecture,
see Proposition~\ref{prop:expDecay}.

To see that $\mathbf{D}_\delta : W^{k,p,\delta} \to W^{k-1,p,\delta}$ is
Fredholm and to compute its index, we can identify it with a Cauchy-Riemann
type operator from $W^{k,p}$ to $W^{k-1,p}$.  Indeed, pick any smooth
function $f : \dot{\Sigma} \to \RR$ with $f(s,t) = \mp \delta s$ on the
cylindrical ends near~$\Gamma^\pm$, define Banach space isomorphisms
\begin{equation*}
\begin{split}
\Phi_\delta : W^{k,p} &\to W^{k,p,\delta} : \eta \mapsto e^f \eta, \\
\Psi_\delta : W^{k-1,p} &\to W^{k-1,p,\delta} : \theta \mapsto e^f \theta,
\end{split}
\end{equation*}
and consider the bounded linear map
$$
\mathbf{D}_\delta' := \Psi_\delta^{-1} \mathbf{D}_\delta \Phi_\delta :
W^{k,p}(u^*T\widehat{W}) \to W^{k-1,p}(\overline{\Hom}_\CC(T\dot{\Sigma},
u^*T\widehat{W})).
$$
Using the Leibniz rule for $\mathbf{D}_\delta$, it is straightforward to
show that $\mathbf{D}_\delta'$ is also a linear Cauchy-Riemann type operator.
Moreover, suppose $\mathbf{D}_\delta$ takes the form $\dbar + S(s,t)$ in
coordinates and trivialization on the cylindrical end near $z_j^\pm$, 
where $S(s,t) \to S_\infty(t)$ as $s \to \pm\infty$ and
$\mathbf{A}_{\gamma_j^\pm} = -i \p_t - S_\infty(t)$.  Then $\mathbf{D}_\delta'$
on this same end takes the form
$$
\mathbf{D}_\delta' \eta = e^{\pm \delta s} (\dbar + S(s,t))(e^{\mp \delta s} \eta)
= \dbar \eta + (S(s,t) \mp \delta) \eta
$$
and is therefore asymptotic to the perturbed asymptotic operator
$$
\tilde{\mathbf{A}}_j^\pm := \left( (-i\p_t) \oplus \mathbf{A}_{\gamma_j^\pm} \right)
\pm \delta.
$$
The latter is the direct sum of two asymptotic operators
$-i\p_t \pm \delta$ on the trivial line bundle and
$\mathbf{A}_{\gamma_j^\pm} \pm \delta$ on $(\gamma_j^\pm)^*\xi_\pm$
respectively.  Since $\gamma_j^\pm$ is nondegenerate by assumption and the
spectrum of $\mathbf{A}_{\gamma_j^\pm}$ is discrete, we can
assume $\ker (\mathbf{A}_{\gamma_j^\pm} \pm \delta)$ remains trivial if
$\delta > 0$ is sufficiently small, and the Conley-Zehnder index of this
perturbed operator will be the same as without the perturbation.
On the other hand, the spectrum
of $-i\p_t$ consists of the integer multiples of $2\pi$, thus
$-i \p_t \pm \delta$ also becomes nondegenerate for any $\delta > 0$ small.
Its Conley-Zehnder index can be deduced from the winding numbers of its
eigenfunctions using Theorem~\ref{thm:CZwinding} in Lecture~\ref{lec:asymptotic}: 
$-i \p_t$ has a $2$-dimensional nullspace consisting of sections with
winding number~$0$, and this becomes an eigenspace for the smallest
positive eigenvalue if the puncture is positive or the largest negative
eigenvalue if the puncture is negative.  Theorem~\ref{thm:CZwinding} thus gives
$$
\muCZ(-i \p_t \pm \delta) = \mp 1,
$$
and therefore,
$$
\muCZ^\tau(\tilde{\mathbf{A}}_j^\pm) = \mp 1 + \muCZ^\tau(\gamma_j^\pm).
$$
Plugging this into the general index formula from Lecture~\ref{lec:index} then gives
the stated result.
\end{proof}

Putting back the missing $2(\#\Gamma)$ dimensions in the domain of
$\mathbf{D}_u$, we have:

\begin{cor}
\label{cor:index}
For all $\delta > 0$ sufficiently small, the linearized Cauchy-Riemann
operator
$\mathbf{D}_u : T_u \bB^{k,p,\delta} \to \eE_u^{k-1,p,\delta}$ is Fredholm
with index
$$
\ind(\mathbf{D}_u) = 
n \chi(\Sigma) - (n-1)\#\Gamma +
2 c_1^\tau(u^*T\widehat{W}) + \sum_{j=1}^{k_+} \muCZ^\tau(\gamma_j^+) -
\sum_{j=1}^{k_-} \muCZ^\tau(\gamma_j^-).
$$
\end{cor}

\section{Teichm\"uller slices}

Since the moduli space $\mM(J)$ is not defined with reference to any
fixed complex structure on the domains $\dot{\Sigma}$, we must build
this freedom into the setup.  For a more detailed version of the following
discussion, see \cite{Wendl:lecturesV33}*{\S 4.2.1}.

For any integers $g,\ell \ge 0$, the
\defin{moduli space of Riemann surfaces} of genus~$g$ with $\ell$ marked points
is a space of equivalence classes
$$
\mM_{g,\ell} = \left\{ (\Sigma,j,\Theta) \right\} \big/ \sim
$$
where $(\Sigma,j)$ is a compact connected surface with genus~$g$,
$\Theta \subset \Sigma$ is an ordered set of $\ell$ points and equivalence
is defined via biholomorphic maps that preserve the marked points with their
ordering.  This space has been studied extensively in algebraic geometry,
though it can also be understood using the same global analytic methods
that we have been applying for $\mM(J)$.  It is known in particular that
$\mM_{g,\ell}$ is always a smooth orbifold, and for any
$[(\Sigma,j,\Theta)] \in \mM_{g,\ell}$, it satisfies
\begin{equation}
\label{eqn:indexSurface}
\dim \Aut(\Sigma,j,\Theta) - \dim \mM_{g,\ell} = 3\chi(\Sigma) - 2\ell,
\end{equation}
where $\Aut(\Sigma,j,\Theta)$ is the group of biholomorphic transformations
of $(\Sigma,j)$ that fix the points in~$\Theta$.  This group is finite
whenever $(\Sigma,j,\Theta)$ is \defin{stable}, meaning
$\chi(\Sigma \setminus \Theta) < 0$, and in that case
\eqref{eqn:indexSurface} turns into the well-known dimension formula
$$
\dim \mM_{g,\ell} = -3\chi(\Sigma) + 2 \ell = 6g - 6 + 2 \ell.
$$
This is also the dimension of the \defin{Teichm\"uller space}
$$
\tT(\Sigma,\Theta) := \jJ(\Sigma) / \Diff_0(\Sigma,\Theta),
$$
where $\jJ(\Sigma)$ denotes the space of all smooth complex structures
on $\Sigma$ compatible with its orientation, and $\Diff_0(\Sigma,\Theta)$
is the identity component of the group of diffeomorphisms that fix~$\Theta$.
It is a classical result that $\tT(\Sigma,\Theta)$ is a smooth manifold
of the same dimension as $\mM_{g,\ell}$, and indeed, the latter can be
presented as the quotient of the former by the discrete action of the
mapping class group of $(\Sigma,\Theta)$.

Equation~\eqref{eqn:indexSurface} is actually a formula for a Fredholm index.
To see how this works, consider first the case $\ell = 0$.  The right hand
side is then $\chi(\Sigma) + 2 c_1(T\Sigma)$, which is, according to
Riemann-Roch, the index of the natural Cauchy-Riemann operator on $T\Sigma$
that defines its holomorphic structure.  This operator can also be
interpreted as the linearization at the identity map of the \emph{nonlinear}
Cauchy-Riemann operator for holomorphic maps $(\Sigma,j) \to (\Sigma,j)$,
so its kernel is naturally isomorphic to $T_{\Id} \Aut(\Sigma,j)$.
Similarly, one can show that the cokernel of this operator is naturally 
isomorphic to $T_{[j]} \tT(\Sigma)$.  This discussion remains valid if
marked points are included: the main difference is then that the Cauchy-Riemann
operator on $T\Sigma$ should be restricted to a space of vector fields that
vanish at~$\Theta$, defining a $2\ell$-codimensional subspace as the domain
and thus reducing the index by~$2\ell$.

For a proof of the following, see \cite{Wendl:lecturesV33}*{Chapter~4}
and \cite{Wendl:automatic}*{\S 3.1}.

\begin{prop}
\label{prop:Teichmueller}
Given a closed Riemann surface $(\Sigma,j)$ with a finite ordered set
$\Theta \subset \Sigma$, there exists a smooth finite-dimensional submanifold
$\tT \subset \jJ(\Sigma)$ with the following properties:
\begin{enumerate}
\item The map $\tT \to \tT(\Sigma,\Theta) : j' \mapsto [j']$ is bijective
onto a neighborhood of $[j]$ in $\tT(\Sigma,\Theta)$;
\item The subspace $T_j \tT \subset \Gamma(\overline{\End}_\CC(T\Sigma))$
is complementary in $W^{k-1,p}(\overline{\End}_\CC(T\Sigma))$ to the
image of the standard Cauchy-Riemann operator of $T\Sigma$ acting on the
domain $\{ X \in W^{k,p}(T\Sigma)\ |\ X|_{\Theta} = 0 \}$;
\item Every $j' \in \tT$ equals $j$ near~$\Theta$ and is
invariant under the action of $\Aut(\Sigma,j,\Theta)$ by diffeomorphisms
on~$\Sigma$.
\end{enumerate}
\qed
\end{prop}

We will refer to the family $\tT \subset \jJ(\Sigma)$ in this proposition
as a \defin{Teichm\"uller slice through~$j$}.

\section{Fredholm regularity and the implicit function theorem}
\label{sec:IFT}

We are now in a position to define the necessary regularity condition and
prove that a neighborhood of any given regular element
$[(\Sigma,j_0,\Gamma^+,\Gamma^-,\Theta,u_0)]$ in $\mM(J)$ is an orbifold
of the stated dimension.  After reparametrizing, we can assume without
loss of generality that $\Sigma$, $\Gamma^\pm$ and $\Theta$ are precisely
the data that were fixed in \S\ref{sec:setup2}, and $j_0 \in \jJ(\Sigma)$
matches $j_\Gamma$ on our fixed coordinate neighborhoods of~$\Gamma^\pm$.
We can then choose a Teichm\"uller slice
$$
\tT \subset \jJ(\Sigma)
$$
through $j_0$ as provided by Prop.~\ref{prop:Teichmueller}, but with 
$j$ in that statement replaced by $j_0$ and $\Theta$ 
replaced by $\Gamma^+ \cup \Gamma^- \cup \Theta$.  In particular, 
$\tT$ is invariant under the action of the group
$$
G_0 := \Aut(\Sigma,j_0,\Gamma^+ \cup \Gamma^- \cup \Theta),
$$
and \eqref{eqn:indexSurface} now becomes
\begin{equation}
\label{eqn:indexSurface2}
\dim G_0 - \dim \tT = 3\chi(\Sigma) - 2(k_+ + k_- + m).
\end{equation}

There is a natural extension of the nonlinear operator $\dbar_{j,J}$
in \S\ref{sec:setup2} to a smooth section
$$
\dbar_J : \tT \times \bB^{k,p,\delta} \to \eE^{k-1,p,\delta} : (j,u) \mapsto
Tu + J \circ Tu \circ j
$$
of a Banach space bundle $\eE^{k-1,p,\delta} \to \tT \times \bB^{k,p,\delta}$
with fibers
$$
\eE^{k-1,p,\delta}_{(j,u)} = W^{k-1,p,\delta}\big(\overline{\Hom}_\CC(
(T\dot{\Sigma},j),(u^*T\widehat{W},J))\big).
$$
The zero set $\dbar_J^{-1}(0) \subset \tT \times \bB^{k,p,\delta}$ consists of
pairs $(j,u)$ for which $u : (\dot{\Sigma},j) \to (\widehat{W},J)$ is
pseudoholomorphic, and it contains $(j_0,u_0)$ by construction.  It also
admits a natural action of the automorphism group $G_0$,
$$
G_0 \times \dbar_J^{-1}(0) \to \dbar_J^{-1}(0) : (\varphi,(j,u)) \mapsto
(\varphi^*j,u \circ \varphi),
$$
whose stabilizer at $(j_0,u_0)$ is $\Aut(u_0)$, a finite group whenever
$u_0$ is not constant.  Observe that any two elements in the same
$G_0$-orbit of $\dbar_J^{-1}(0)$ define equivalent elements of the moduli
space~$\mM(J)$, as they are related to each other by a biholomorphic
reparametrization that fixes the punctures and marked points.

\begin{lemma}
\label{lemma:bijection}
The map
$$
\dbar_J^{-1}(0) \big/ G_0 \to \mM(J) : [(j,u)] \mapsto [(\Sigma,j,\Gamma^+,\Gamma^-,\Theta,u)]
$$
is a homeomorphism between open neighborhoods of 
$[(j_0,u_0)]$ and $[(\Sigma,j_0,\Gamma^+,\Gamma^-,\Theta,u_0)]$.
\end{lemma}
\begin{proof}
This depends fundamentally on the same fact underlying the smoothness of
Teichm\"uller space: the action of $\Diff_0(\Sigma,\Gamma^+ \cup \Gamma^- \cup \Theta)$
on $\jJ(\Sigma)$ is free and proper.\footnote{This is true at least in the
stable case, i.e.~when $\chi(\dot{\Sigma} \setminus \Theta) < 0$.  There are
finitely many cases not satisfying this hypothesis, for which the lemma can
be proved by more direct arguments since explicit descriptions of both
Teichm\"uller space and the automorphism groups of Riemann surfaces are
available; see \cite{Wendl:automatic}*{\S 3.1 and \S 3.2} for more details.}
See the proof of \cite{Wendl:lecturesV33}*{Theorem~4.3.6}.
\end{proof}

\begin{defn}
We say that $[(\Sigma,j_0,\Gamma^+,\Gamma^-,\Theta,u_0)]$ is
\defin{Fredholm regular} if there exists a choice of Teichm\"uller slice
$\tT$ through $j_0$ such that the linearization
$$
D\dbar_J(j_0,u_0) : T_{j_0}\tT \oplus T_{u_0} \bB^{k,p,\delta} \to
\eE_{(j_0,u_0)}^{k-1,p,\delta}
$$
is surjective.
\end{defn}

One can show that the surjectivity condition in this definition does not
actually depend on the choice of Teichm\"uller slice.  This follows from the
identification of $T_{j_0}\tT$ with the cokernel of the
natural Cauchy-Riemann operator on $T\dot{\Sigma}$; see
\cite{Wendl:lecturesV33}*{Lemma~4.3.2}.

\begin{proof}[Proof of Theorem~\ref{thm:moduliDimension2}]
The fact that $\mM(J)$ is an orbifold in a neighborhood of 
$[(\Sigma,j_0,\Gamma^+,\Gamma^-,\Theta,u_0)]$ with isotropy group
$\Aut(u_0)$ follows from Lemma~\ref{lemma:bijection} and the implicit 
function theorem, which gives $\dbar_J^{-1}(0)$ the structure of a
finite-dimensional manifold near $(j_0,u_0)$ if Fredholm regularity is
satisfied.  There is a bit of work to be done in showing that
transition maps relating any two overlapping charts that arise in this way 
from the implicit function theorem are smooth; for this, we refer again
to the proof of Theorem~4.3.6 in \cite{Wendl:lecturesV33} and merely comment
that the key ingredient is elliptic regularity.

The dimension of $\mM(J)$ is
$$
\dim \mM(J) = \dim \dbar_J^{-1}(0) - \dim G_0 = \ind D\dbar_J(j_0,u_0) - \dim G_0.
$$
The restriction of $D\dbar_J(j_0,u_0)$ to $T_{u_0}\bB^{k,p,\delta}$ is
the operator $\mathbf{D}_{u_0}$ that we studied in \S\ref{sec:setup2}, hence
$$
\ind D\dbar_J(j_0,u_0) = \dim \tT + \ind \mathbf{D}_{u_0}.
$$
Using \eqref{eqn:indexSurface2} to replace $\dim \tT - \dim G_0$ and combining
this with Corollary~\ref{cor:index} now gives
the stated formula for $\dim \mM(J)$.
\end{proof}

\section{A universal moduli space}
\label{sec:universal}

The remainder of this lecture is devoted to the proof of 
Theorem~\ref{thm:genericCobordism}.  The main tool for this purpose is the
Sard-Smale theorem \cite{Smale:Sard}, an infinite-dimensional version of
Sard's theorem stating that the regular values of a smooth nonlinear
Fredholm map between separable Banach spaces (i.e.~a smooth map whose derivative 
at every point is a Fredholm operator) form a comeager subset of the 
target space.  In order to incorporate perturbations of the almost complex
structure into our functional analytic setup, we need to choose a suitable
Banach manifold of almost complex structures.  All known ways of doing
this are in some sense non-ideal, e.g.~one could take almost complex structures
of class $C^k$ or $W^{k,p}$, but this necessarily introduces non-smooth
almost complex structures into the picture, with the consequence that the
nonlinear Cauchy-Riemann operator has only finitely many derivatives.
That is not the end of the world, and indeed, this is the approach taken
in \cite{McDuffSalamon:Jhol}, but I will instead present an approach that
was introduced by Floer in \cite{Floer:action}, in terms of what is now
called the ``Floer $C_\varepsilon$ space''.  The idea is to work with a Banach
manifold that continuously embeds into the space of smooth almost complex
structures, so that the nonlinear Cauchy-Riemann operator will always be
smooth.  It's a nice trick, but the catch is that we obtain a space
that is strictly smaller than the actual space of smooth almost complex 
structures we're interested in, and has a much stronger topology.
The $C_\varepsilon$ space should be viewed as a useful tool but not a
deeply meaningful object---you might notice that while some of the 
intermediate results stated below depend on its (somewhat ad hoc)
definition, Theorem~\ref{thm:genericCobordism} does not.  This is due to a general
trick described in \S\ref{sec:Taubes} below for turning results about
$C_\varepsilon$ into results about~$C^\infty$.

As in the statement of Theorem~\ref{thm:genericCobordism}, assume
$\uU \subset W^{r_0}$ is open and
$\Jfix \in \jJ(\omega_\psi,r_0,\hH_+,\hH_-)$.  Let
$$
\jJ_\uU := \left\{ J \in \jJ(\omega_\psi,r_0,\hH_+,\hH_-) \ \Big|\ 
\text{$J = \Jfix$ on $\widehat{W} \setminus \uU$} \right\},
$$
and choose any almost complex structure
$$
\Jref \in \jJ_\uU.
$$
We can regard $\jJ_\uU$ as a smooth Fr\'echet manifold with tangent spaces
$$
T_{\Jref} \jJ_\uU = \left\{ Y \in \Gamma\big(\overline{\End}_\CC(T\widehat{W},\Jref)\big)\ \Big|\ 
Y|_{\widehat{W}\setminus \uU} \equiv 0 \text{ and }
\omega_\psi(\cdot,Y\cdot) + \omega_\psi(Y\cdot,\cdot) \equiv 0 \right\},
$$
where the antilinearity of $Y \in T_{\Jref} \jJ_\uU$ means that $Y$ is tangent
to the space almost complex structures, and the condition relating it to
$\omega_\psi$ means that these structures are compatible with~$\omega_\psi$.
One can check that the map
$$
Y \mapsto J_Y := \left( \1 + \frac{1}{2} \Jref Y \right) \Jref
\left( \1 + \frac{1}{2} \Jref Y \right)^{-1}
$$
maps a neighborhood of $0 \in T_{\Jref} \jJ_\uU$ bijectively to
a neighborhood of $\Jref$ in~$\jJ_\uU$.  We thus fix a sufficiently small
constant $c > 0$ and define the space of ``$C_\varepsilon$-small perturbations
of~$\Jref$'' by
$$
\jJ_\uU^\varepsilon := \left\{ J_Y \in \jJ_\uU \ \bigg|\ 
Y \in T_{\Jref} \jJ_\uU \text{ with } \sum_{\ell=0}^\infty \varepsilon_\ell \| Y \|_{C^\ell(\overline{\uU})} < c \right\},
$$
where $\varepsilon := (\varepsilon_\ell)_{\ell=0}^\infty$ is a fixed sequence
of positive numbers with $\varepsilon_\ell \to 0$ as $\ell \to \infty$.
The sum
$$
\| Y \|_{C_\varepsilon} := \sum_{\ell=0}^\infty \varepsilon_\ell \| Y \|_{C^\ell(\overline{\uU})}
$$
defines a norm, and the space of smooth sections $Y \in T_{\Jref} \jJ_\uU$
for which this norm is finite is then a separable Banach space; see
Appendix~\ref{app:Cepsilon} for a proof of this statement.
This makes $\jJ_\uU^\varepsilon$ a separable and metrizable Banach manifold, 
as the map $J_Y \mapsto Y$ can be
viewed as a chart identifying it with an open subset of the aforementioned
Banach space.  Not every $J \in \jJ_\uU$ near $\Jref$ belongs to
$\jJ_\uU^\varepsilon$, but there is a continuous inclusion
$$
\jJ_\uU^\varepsilon \hookrightarrow \jJ_\uU,
$$
where the latter carries its usual $C^\infty$-topology and $\jJ_\uU^\varepsilon$
carries the topology induced by the $C_\varepsilon$-norm.
By a lemma due to Floer, choosing a sequence $\varepsilon_\ell$ that decays
sufficiently fast makes $\jJ_\uU^\varepsilon$ large enough to contain
perturbations in arbitrary directions with arbitrarily small support
near arbitrary points in~$\uU$;
see Theorem~\ref{thm:bump} in Appendix~\ref{app:Cepsilon} for a precise
version of this statement and its proof.  We will assume from now on
that a suitably fast decaying sequence has been fixed.

We now define a \defin{universal moduli space}
\begin{equation*}
\begin{split}
\mM^*(\jJ_\uU^\varepsilon) := \big\{ (u,J) \ \big|\ 
&\text{$J \in \jJ_\uU^\varepsilon$, $u \in \mM(J)$ and}\\
&\text{$u$ has an injective point mapped into~$\uU$} \big\}.
\end{split}
\end{equation*}
The terminology is somewhat unfortunate, as $\mM^*(\jJ_\uU^\varepsilon)$
depends on many auxiliary choices such as $\Jref$ and $(\varepsilon_\ell)_{\ell=0}^\infty$
and thus should not really be thought of as a ``universal'' object.
Nonetheless:

\begin{lemma}
\label{lemma:universal}
The universal moduli space $\mM^*(\jJ_\uU^\varepsilon)$ is a smooth separable Banach
manifold, and the projection $\mM^*(\jJ_\uU^\varepsilon) \to \jJ_\uU^\varepsilon :
(u,J) \mapsto J$ is smooth.
\end{lemma}
\begin{proof}
As in the proof of Theorem~\ref{thm:moduliDimension2}, one can identify
$\mM^*(\jJ_\uU^\varepsilon)$ locally with the zero set of a smooth section
of a Banach space bundle.  Suppose $J_0 \in \jJ_\uU^\varepsilon$ and
$[(\Sigma,j_0,\Gamma^+,\Gamma^-,\Theta,u_0)] \in \mM(J_0)$ where
$u_0 : \dot{\Sigma} \to \widehat{W}$ has an injective point $z_0$ with
$u_0(z_0) \in \uU$.  Choose a Teichm\"uller slice $\tT$ through $j_0$
as in Proposition~\ref{prop:Teichmueller} and consider the smooth section
$$
\dbar : \tT \times \bB^{k,p,\delta} \times \jJ_\uU^\varepsilon \to
\eE^{k-1,p,\delta} : (j,u,J) \mapsto Tu + J \circ Tu \circ j,
$$
where $\eE^{k-1,p,\delta}$ is the obvious extension of our previous Banach
space bundle to a bundle over $\tT \times \bB^{k,p,\delta} \times \jJ_\uU^\varepsilon$.
We're assuming as before that $k \in \NN$, $1 < p < \infty$, $kp > 2$,
and $\delta > 0$ is small.
A neighborhood of $(u_0,J_0)$ in $\mM^*(\jJ_\uU^\varepsilon)$ can then be
identified with a neighborhood of $[(j_0,u_0,J_0)]$ in
$$
\dbar^{-1}(0) \big/ G_0,
$$
where $G_0 := \Aut(\Sigma,j_0,\Gamma^+ \cup \Gamma^- \cup \Theta)$ acts
on $\dbar^{-1}(0)$ by $\varphi \cdot (j,u,J) := (\varphi^*j,u\circ \varphi,J)$.
Since $u_0$ has an injective point, $\Aut(u_0)$ is trivial and the
$G_0$-action at $(j_0,u_0,J_0)$ is therefore free;
hence it suffices to show that $\dbar^{-1}(0)$ is a smooth Banach manifold
near $(j_0,u_0,J_0)$.  This follows from the implicit function theorem if
we can show that
$$
D\dbar(j_0,u_0,J_0) : T_{j_0}\tT \oplus T_{u_0}\bB^{k,p,\delta} \oplus T_{J_0}\jJ_\uU^\varepsilon
\to \eE^{k-1,p,\delta}_{(j_0,u_0,J_0)}
$$
is surjective; indeed, the infinite-dimensional implicit function theorem
(see \cite{Lang:analysis}) requires the additional hypothesis that 
$D\dbar(j_0,u_0,J_0)$ has a bounded right inverse, but this is immediate
since the restriction of this operator to the factor $T_{u_0}\bB^{k,p,\delta}$
is Fredholm (see Exercise~\ref{EX:bri} below).  We claim in fact that
\begin{equation*}
\begin{split}
T_{u_0}\bB^{k,p,\delta} \oplus T_{J_0}\jJ_\uU^\varepsilon &\to \eE^{k-1,p,\delta}_{(j_0,u_0,J_0)} \\
(\eta,Y) &\mapsto D\dbar(j_0,u_0,J_0)(0,\eta,Y) = \mathbf{D}_{u_0} \eta +
Y \circ Tu_0 \circ j_0
\end{split}
\end{equation*}
is surjective.  Consider first the case $k=1$,\footnote{Since the present
discussion is purely linear, it does not require the assumption $kp > 2$.} 
so we are looking at a bounded linear map
$$
W^{1,p,\delta}(u_0^*T\widehat{W}) \oplus V_\Gamma \oplus T_{J_0}\jJ_\uU^\varepsilon
\to L^{p,\delta}(\overline{\Hom}_\CC(T\dot{\Sigma},u_0^*T\widehat{W})).
$$
Note that the dual of any space of sections of class $L^{p,\delta}$ can be 
identified with sections of class $L^{q,-\delta}$ for 
$\frac{1}{p} + \frac{1}{q} = 1$ (recall Remark~\ref{remark:negativeExp}).
Indeed, choosing a suitable $L^2$-pairing defines a bounded bilinear map
\begin{equation}
\label{eqn:LpdeltaPairing}
\langle\ ,\ \rangle_{L^2} : L^{p,\delta} \times L^{q,-\delta} \to \RR ,
\end{equation}
and one can use isomorphisms of the form $L^p \to L^{p,\delta} : \eta \mapsto
e^f \eta$ as in the proof of Lemma~\ref{lemma:expRescue} to prove 
$(L^{p,\delta})^* \cong L^{q,-\delta}$ as a corollary of the standard fact
that $(L^p)^* \cong L^q$.  With this understood, observe that
since $\mathbf{D}_{u_0} : W^{1,p,\delta} \oplus V_\Gamma \to L^{p,\delta}$ 
is Fredholm, we know by Exercise~\ref{EX:closedRange} below that the map 
under consideration has closed range.
Thus if it is not surjective, the Hahn-Banach theorem provides a nontrivial
element $\theta \in L^{q,-\delta}(\overline{\Hom}_\CC(T\dot{\Sigma},u_0^*T\widehat{W}))$
that annihilates its image under the pairing \eqref{eqn:LpdeltaPairing},
which amounts to the two conditions
\begin{equation}
\label{eqn:twoConditions}
\begin{split}
\langle \mathbf{D}_{u_0} \eta , \theta \rangle_{L^2} = 0 & \text{ for all
$\eta \in W^{1,p,\delta}(u_0^*T\widehat{W}) \oplus V_\Gamma$}, \\
\langle Y \circ Tu_0 \circ j_0 , \theta \rangle_{L^2} = 0 & \text{ for all
$Y \in T_{J_0}\jJ_\uU^\varepsilon$}.
\end{split}
\end{equation}
The first relation is valid in particular for all smooth sections $\eta$
with compact support and thus means that $\theta$ is a weak solution to
the formal adjoint equation $\mathbf{D}_{u_0}^*\theta = 0$; applying
elliptic regularity and the similarity principle, $\theta$ is therefore
smooth and has only isolated zeroes.  We will see however that this
contradicts the second relation as long as there exists an injective point
$z_0 \in \dot{\Sigma}$ with $u_0(z_0) \in \uU$.  Indeed, since the set of
injective points with this property is open and zeroes of $\theta$ are
isolated, let us assume without loss of generality that $\theta(z_0) \ne 0$.
Then by a standard lemma in symplectic linear algebra (see \cite{Wendl:lecturesV33}*{Lemma~4.4.12}),
one can find a smooth section $Y \in T_{J_0}\jJ_\uU$ whose value at
$u_0(z_0)$ is chosen such that $Y \circ Tu_0 \circ j_0 = \theta$ at
$z_0$, so their pointwise inner product is positive in some neighborhood
of~$z_0$.  But by Theorem~\ref{thm:bump}, one can multiply a small
perturbation of $Y$ by a bump function to produce a section (still denoted by~$Y$) 
of class $C_\varepsilon$ so that the pointwise inner product of
$Y \circ Tu_0 \circ j_0$ with $\theta$ is positive near~$z_0$ but
vanishes everywhere else; note that this requires the assumption
$u_0^{-1}(u_0(z_0)) = \{z_0\}$, so that the value of $Y$ near $u_0(z_0)$
affects the value of $Y \circ Tu_0 \circ j_0$ near $z_0$ but nowhere
else.  This violates the second condition in \eqref{eqn:twoConditions}
and thus completes the proof for $k=1$.
In the general case, suppose $\theta \in W^{k-1,p,\delta}(\overline{\Hom}_\CC(T\dot{\Sigma},u_0^*T\widehat{W}))$.
Then $\theta$ is also of class $L^{p,\delta}$, so surjectivity in the $k=1$
case implies the existence of $\eta \in W^{1,p,\delta}$ and 
$Y \in T_{J_0}\jJ_\uU^\varepsilon$ with $\mathbf{D}_{u_0} \eta +
Y \circ Tu_0 \circ j_0 = \theta$.  Since $Y \circ Tu_0 \circ j_0$ is smooth
with compact support, one can then use elliptic regularity
to show $\eta \in W^{k,p,\delta}$, and this proves surjectivity for arbitrary
$k \in \NN$ and $p \in (1,\infty)$.

The implicit function theorem now implies that whenever $kp > 2$ so that
$\bB^{k,p,\delta}$ is a well-defined Banach manifold,
$\dbar^{-1}(0)$ is a smooth
Banach submanifold of $\tT \times \bB^{k,p,\delta} \times \jJ_\uU^\varepsilon$
in a neighborhood of $(j_0,u_0,J_0)$.  The projection map
$$
\dbar^{-1}(0) \to \jJ_\uU^\varepsilon : (j,u,J) \mapsto J
$$
is also smooth since it is the restriction to a smooth submanifold of the
obviously smooth projection map $\tT \times \bB^{k,p,\delta} \times
\jJ_\uU^\varepsilon \to \jJ_\uU^\varepsilon$.  Since $G_0$ acts freely
and properly on $\dbar^{-1}(0)$, the quotient $\dbar^{-1} / G_0$ then
inherits a smooth Banach manifold structure for which the projection is still
smooth, and this quotient is identified locally with $\mM^*(\jJ_\uU^\varepsilon)$.
Smoothness of transition maps is shown via the same regularity arguments as
in the proof of Theorem~\ref{thm:moduliDimension2}.
\end{proof}

\begin{exercise}
\label{EX:closedRange}
Show that if $X$, $Y$ and $Z$ are Banach spaces, $\mathbf{T} : X \to Y$ is a 
Fredholm operator and $\mathbf{A} : Z \to Y$ is a bounded linear operator,
then the linear map 
$$
\mathbf{L} : X \oplus Z \to Y : (x,z) \mapsto \mathbf{T}x + \mathbf{A}z
$$
has closed range.
\textsl{Hint: it might help to write $X = V \oplus \ker \mathbf{T}$ and
$Y = W \oplus \coker C$ so that $C \cong \coker\mathbf{T}$ and
$V \stackrel{\mathbf{T}}{\longrightarrow} W$ is an isomorphism.}
\end{exercise}

\begin{exercise}
\label{EX:bri}
Under the same assumptions as in Exercise~\ref{EX:closedRange}, show that if
$\mathbf{T}$ is surjective, then $\mathbf{L}$ has a bounded right inverse.
\end{exercise}

\section{Applying the Sard-Smale theorem}

We claim now that the smooth map
\begin{equation}
\label{eqn:projection}
\mM^*(\jJ_\uU^\varepsilon) \to \jJ_\uU^\varepsilon :
(u,J) \mapsto J
\end{equation}
is a nonlinear Fredholm map, i.e.~its derivative at every point is a
Fredholm operator.  Using the local identification of
$\mM^*(\jJ_\uU^\varepsilon)$ with $\dbar^{-1}(0) / G_0$ as in the proof of
Lemma~\ref{lemma:universal} and lifting the projection to $\dbar^{-1}(0)$,
the derivative of $\dbar^{-1}(0) \to \jJ_\uU^\varepsilon$ at $(j_0,u_0,J_0)$
takes the form
$$
\ker D\dbar(j_0,u_0,J_0) \to T_{J_0}\jJ_\uU^\varepsilon : (y,\eta,Y) \mapsto Y.
$$
The Fredholm property for this projection is a consequence of the Fredholm
property for $\mathbf{D}_{u_0}$ via the following general lemma,
whose proof is a routine matter of linear algebra 
(cf.~\cite{Wendl:lecturesV33}*{Lemma~4.4.13}):

\begin{lemma}
\label{lemma:linearAlg}
Under the assumptions of Exercise~\ref{EX:closedRange}, suppose
$\mathbf{L}$ is surjective.  Then the projection
$$
\boldsymbol{\Pi} : \ker \mathbf{L} \to Z : (x,z) \mapsto z
$$
has kernel and cokernel isomorphic to the kernel and cokernel respectively of
$\mathbf{T} : X \to Y$.
\qed
\end{lemma}

By the Sard-Smale theorem, the set of regular values of the projection
\eqref{eqn:projection} is a comeager subset
$$
\jJ_\uU^{\varepsilon,\reg} \subset \jJ_\uU^\varepsilon,
$$
and by Lemma~\ref{lemma:linearAlg}, every $(u_0,J_0) \in \mM^*(\jJ_\uU^\varepsilon)$
with $J \in \jJ_\uU^{\varepsilon,\reg}$ then has the property that
$$
D\dbar_{J_0}(j_0,u_0) : T_{j_0}\tT \oplus T_{u_0}\bB^{k,p,\delta} \to
\eE_{(j_0,u_0)}^{k-1,p,\delta}
$$
is surjective, which means $u_0$ represents a Fredholm regular element
of $\mM(J_0)$.

\section{From $C_\varepsilon$ to~$C^\infty$}
\label{sec:Taubes}

The arguments above would constitute a proof of Theorem~\ref{thm:genericCobordism}
if we were allowed to replace the space of smooth almost complex structures
$\jJ_\uU$ with the space $\jJ_\uU^\varepsilon$ of $C_\varepsilon$-small
perturbations of~$\Jref$.  Let us \emph{define}
$$
\jJ_\uU^\reg \subset \jJ_\uU
$$
to be the space of all $J \in \jJ_\uU$ with the property that all
curves in $\mM(J)$ that have injective points mapping to~$\uU$ are Fredholm
regular.  The theorem claims that this set is comeager in~$\jJ_\uU$.
We can already see at this point that it is dense: indeed, the Baire
category theorem implies that $\jJ_\uU^{\varepsilon,\reg}$ is dense in
$\jJ_\uU^\varepsilon$, so in particular there exists a sequence
$J_\nu \in J_\uU^{\varepsilon,\reg}$ that converges in to $\Jref$ in
the $C_\varepsilon$-topology and therefore also in the $C^\infty$-topology.
The choice of $\Jref \in \jJ_\uU$ in this discussion was arbitrary, so this
proves density.

To prove that $\jJ_\uU^\reg$ is not only dense but also contains a countable
intersection of \emph{open} and dense sets in $\jJ_\uU$, we can adapt an
argument originally due to Taubes.  The idea is to present the sets of
somewhere injective curves in $\mM(J)$ as countable unions of compact 
subsets $\mM_N^*(J)$
for $N \in \NN$, and thus present $\jJ_\uU^\reg$ as a corresponding
countable intersection of spaces $\jJ_\uU^{\reg,N}$ that achieve regularity
only for the elements in $\mM_N^*(J)$.  The compactness of $\mM_N^*(J)$ will
then permit us to prove that $\jJ_\uU^{\reg,N}$ is not only dense but also
open.

The definition of $\mM_N^*(J)$ is motivated in part by the knowledge that
spaces of $J$-holomorphic curves have natural compactifications.  We have
not yet discussed the compactification $\overline{\mM}(J)$ of $\mM(J)$,
but we have covered enough of the analytical techniques behind this 
construction to suffice for the present discussion.  Recall first that
the moduli space of Riemann surfaces $\mM_{g,\ell}$ of genus $g$ with
$\ell$ marked points also has a natural compactification whenever
$2g + \ell \ge 3$, known as the \defin{Deligne-Mumford compactification}
$$
\overline{\mM}_{g,\ell} \supset \mM_{g,\ell}.
$$
The space $\overline{\mM}_{g,\ell}$ consists of ``nodal'' Riemann surfaces,
which can be understood as objects that arise from smooth Riemann surfaces
with pair-of-pants decompositions in the limit where some of the lengths
of the circles separating two pairs of pants from each other may degenerate
to~$0$ (see e.g.~\cite{SeppalaSorvali}).  We will discuss this in a bit 
more detail in Lecture~\ref{lec:compactness}; for now, all you really need to know is that
$\overline{\mM}_{g,\ell}$ is a compact and metrizable topological space
that contains $\mM_{g,\ell}$ as an open subset.  Let us fix a metric on
$\mM_{g,\ell}$ and denote the distance function by $\dist(\ ,\ )$.

Similarly, fix Riemannian metrics on $\widehat{W}$ and $\dot{\Sigma}$ 
with translation-invariance on the cylindrical ends and use
$\dist(\ ,\ )$ to denote the distance functions.  For $N \in \NN$ and
$J \in \jJ_\uU$, we define
$$
\mM_N^*(J) \subset \mM(J)
$$
to be the set of equivalence classes admitting representatives
$(\Sigma,j,\Gamma^+,\Gamma^-,\Theta,u)$ with the following properties:
\begin{itemize}
\item The equivalence class in $\mM_{g,k_+ + k_- + m}$ represented by
$(\Sigma,j,\Gamma^+ \cup \Gamma^- \cup \Theta)$ lies at a distance of at
most $1/N$ from $\overline{\mM}_{g,k_+ + k_- + m} \setminus \mM_{g,k_+ + k_- + m}$;\footnote{If
the stability condition $2g + k_+ + k_- + m \ge 3$ is not satisfied, one 
should amend this by asking for the distance
condition to hold for some tuple $(\Sigma,j,\Gamma^+ \cup \Gamma^-,\Theta')$,
where $\Theta'$ is the union of $\Theta$ with enough extra marked points to
achieve stability.}
\item $\sup_{z \in \dot{\Sigma}} |du(z)| \le N$;
\item There exists $z_0 \in \dot{\Sigma}$ such that
$$
\dist(u(z_0) , \widehat{W} \setminus \uU) \ge \frac{1}{N}, \qquad
|du(z_0)| \ge \frac{1}{N},
$$
and
$$
\inf_{z \in \dot{\Sigma} \setminus \{z_0\}} 
\frac{\dist(u(z_0),u(z))}{\dist(z_0,z)} \ge \frac{1}{N}.
$$
\end{itemize}
We observe that every element of $\mM_N^*(J)$ has an injective point mapped
into~$\uU$, and conversely, every asymptotically cylindrical $J$-holomorphic
curve with that property belongs to
$\mM_N^*(J)$ for $N \in \NN$ sufficiently large.  It is crucial to observe
that all three conditions in this definition are \emph{closed} conditions: 
morally, we are defining $\mM_N^*(J)$ to be a closed subset in the
compactification of $\mM(J)$, and it will therefore be compact.

Define
$$
\jJ_\uU^{\reg,N} \subset \jJ_\uU
$$
as the set of all $J \in \jJ_\uU$ for which every element of
$\mM_N^*(J)$ is Fredholm regular.

\begin{lemma}
For every $N \in \NN$, $\jJ_\uU^{\reg,N}$ is open and dense.
\end{lemma}
\begin{proof}
Density is immediate, since we've seen already that every $J \in \jJ_\uU$
admits a $C^\infty$-small perturbation that achieves regularity for all
curves in $\bigcup_{N \in \NN} \mM_N^*(J)$.  For openness, suppose the
contrary: then there exists $J_\infty \in \jJ_\uU^{\reg,N}$ and a sequence
$J_\nu \in \jJ_\uU \setminus \jJ_\uU^{\reg,N}$ with $J_\nu \to J_\infty$
in the $C^\infty$-topology.  There must also exist a sequence of curves
$u_\nu \in \mM_N^*(J_\nu)$ that are not Fredholm regular.  By the
definition of $\mM_N^*(J_\nu)$, they have domains that are uniformly
bounded away from the singular part of the Deligne-Mumford space of Riemann
surfaces, so we can extract a subsequence for which these domains converge.
Similarly, the first derivatives of $u_\nu$ are uniformly bounded,
implying in particular a uniform $W^{1,p}$-bound locally for some $p > 2$,
and elliptic regularity (Theorem~\ref{thm:regularity} in Lecture~\ref{lec:local}) turns this into 
uniform $C^\infty$-bound and thus a $C^\infty$-convergent subsequence
$u_\nu \to u_\infty \in \mM_N^*(J_\infty)$.  But $u_\infty$ must then be
Fredholm regular, which is an open condition, implying that $u_\nu$ is
also regular for $\nu$ sufficiently large, and this is a contradiction.
\end{proof}

\begin{proof}[Proof of Theorem~\ref{thm:genericCobordism}]
Since the space of all curves in $\mM(J)$ with injective points mapped into
$\uU$ is the union of the spaces $\mM_N^*(J)$ for $N \in \NN$, we have
$$
\jJ_\uU^\reg = \bigcap_{N \in \NN} \jJ_\uU^{\reg,N},
$$
which is a countable intersection of open and dense sets.
\end{proof}

\chapter{Transversality in symplectizations}
\label{lec:Dragnev}

\minitoc

\vspace{12pt}

This lecture is an addendum to the transversality discussion in Lecture~\ref{lec:transversality}:
we need to prove that Fredholm regularity can also be achieved for generic
\emph{translation-invariant} almost complex structures on symplectizations.

\section{Statement of the theorem and discussion}

Theorem~\ref{thm:genericCobordism} in the previous lecture stated that generic perturbations
of $J$ in a precompact open subset $\uU$ of a completed symplectic cobordism
suffice to achieve regularity for all simple holomorphic curves that pass
through that subset.  In the more specialized setting of a symplectization
$\RR \times M$ with an $\RR$-invariant almost complex structure
$J \in \jJ(\hH)$, we need a more specialized transversality result, as the 
generic perturbation from Theorem~\ref{thm:genericCobordism} cannot be expected to stay in the
space $\jJ(\hH)$, in particular it will usually not be $\RR$-invariant.  
The following statement refers to a stable Hamiltonian
structure $\hH = (\omega,\lambda)$ with induced hyperplane distribution
$\xi = \ker \lambda$ and Reeb vector field $R$, and we denote by
$$
\pi_\xi : T(\RR \times M) \to \xi
$$
the projection along the trivial subbundle generated by $\p_r$ and~$R$.  
We assume as usual that $\mM(J)$ denotes a
moduli space of asymptotically cylindrical $J$-holomorphic curves with a 
fixed genus and number of marked points, representing a fixed relative
homology class and asymptotic to fixed sets of nondegenerate Reeb orbits
at its positive and negative punctures.

\begin{thm}
\label{thm:genericRinvt}
Suppose $M$ is a closed $(2n-1)$-dimensional manifold carrying a stable 
Hamiltonian structure $\hH = (\omega,\lambda)$, $\Jfix \in \jJ(\hH)$, and
$$
\uU \subset M
$$
is an open subset.  Then there exists a comeager subset
$$
\jJ_\uU^\reg \subset \left\{ J \in \jJ(\hH) \ \big|\ 
\text{$J = \Jfix$ on $\RR \times (M \setminus \uU)$} \right\}
$$
such that for every $J \in \jJ_\uU^\reg$, every curve $u \in \mM(J)$
with a representative $u : \dot{\Sigma} \to \RR \times M$ that has
an injective point $z \in \dot{\Sigma}$ satisfying
\begin{enumerate}[label=(\roman{enumi})]
\item $u(z) \in \RR \times \uU$,
\item $\pi_\xi \circ du(z) \ne 0$, and
\item $\im \left( \pi_\xi \circ du(z) \right) \cap \ker \left( d\lambda|_\xi \right) = \{0\}$
\end{enumerate}
is Fredholm regular.
\end{thm}

This result is applied most frequently with $\uU = M$, in which case the 
condition $u(z) \in \RR \times \uU$ is vacuous.  The second and third
conditions on the injective point $z$ can be rephrased by asking for the 
linear map
$$
d\lambda(\pi_\xi\, Tu(X),\cdot)|_{\xi_{u(z)}} : \xi_{u(z)} \to \RR
$$
to be nontrivial for every nonzero $X \in T_z \dot{\Sigma}$.  
If $\lambda$ is contact,
then this is immediate whenever $\pi_\xi \, Tu(X) \ne 0$ since
$d\lambda|_{\xi}$ is nondegenerate, and the condition $\pi_\xi \, Tu(X) \ne 0$
is also easy to achieve:

\begin{prop}
\label{prop:piTu}
If $J \in \jJ(\hH)$, then for any connected $J$-holomorphic curve
$u : (\dot{\Sigma},j) \to (\RR \times M,J)$, the section
$$
\pi_\xi \circ du \in \Gamma(\Hom_\CC(T\dot{\Sigma},u^*\xi))
$$
either is identically zero or has only isolated zeroes.
\end{prop}
As you might guess, this result is a consequence of the similarity
principle; see \S\ref{sec:injective} for a proof.
Notice that if $\pi_\xi \circ du \equiv 0$, then $u$ is everywhere tangent
to the vector fields $\p_r$ and $R$, so if it is asymptotically cylindrical,
then it can only be a trivial cylinder or a cover thereof.

\begin{prop}
\label{prop:trivial}
All trivial cylinders over nondegenerate Reeb orbits have 
index~$0$ and are Fredholm regular.
\end{prop}
\begin{proof}
Let $u_\gamma : \RR \times S^1 \to \RR \times M$ denote the trivial
cylinder over an orbit $\gamma : S^1 \to M$.  The virtual dimension
formula proved in Lecture~\ref{lec:transversality} gives
\begin{equation*}
\begin{split}
\ind(u_\gamma) &= (n-3) \chi(\RR \times S^1) + 2 c_1^\tau(u_\gamma^*T(\RR \times M))
+ \muCZ^\tau(\gamma) - \muCZ^\tau(\gamma) \\
&= 2 c_1^\tau(u_\gamma^*T(\RR \times M)) = 0
\end{split}
\end{equation*}
since the asymptotic trivialization $\tau$ has an obvious extension to a
global trivialization of $u_\gamma^*\xi$, and $u_\gamma^*T(\RR \times M)$
is globally the direct sum of the latter with the trivial line bundle spanned by
$\p_r$ and~$R$.  Using this splitting, the linearized Cauchy-Riemann 
operator $\mathbf{D}_{u_\gamma}$ can be identified with
$\dbar \oplus (\p_s - \mathbf{A}_\gamma)$, where 
$$
\dbar = \p_s + i\p_t : W^{k,p,\delta}(\RR \times S^1,\CC) \oplus V_\Gamma \to
W^{k-1,p,\delta}(\RR \times S^1,\CC)
$$
and
$$
\p_s - \mathbf{A}_\gamma : W^{k,p,\delta}(u_\gamma^*\xi) \to
W^{k-1,p,\delta}(u_\gamma^*\xi).
$$
Here we are assuming without loss of generality that $V_\Gamma$ is a
complex $2$-dimensional space of smooth sections of the trivial line bundle
spanned by $\p_r$ and~$R$ that are constant near infinity, and we are
identifying this with a space of smooth complex-valued functions
on $\RR \times S^1$.
Nondegeneracy implies that $\p_s - \mathbf{A} : W^{k,p} \to W^{k-1,p}$
is an isomorphism, recall Theorem~\ref{thm:invertible} in Lecture~\ref{lec:Fredholm}.  Using weight
functions as in the proof of Lemma~\ref{lemma:expRescue} to define isomorphisms between
$W^{k,p,\delta}$ and $W^{k,p}$, one can identify
$\p_s - \mathbf{A}_\gamma : W^{k,p,\delta} \to W^{k-1,p,\delta}$ with a
small perturbation of the same operator $W^{k,p} \to W^{k-1,p}$,
hence it is also an isomorphism for $\delta > 0$ sufficiently small.
To see that $\dbar : W^{k,p,\delta} \oplus V_\Gamma \to W^{k-1,p,\delta}$
is also surjective, observe first that its index is~$2$; this follows from
our calculation of $\ind(u_\gamma)$ and corresponds to the fact that
$\dim \Aut(\RR \times S^1,i) = 2$.  The kernel of this operator consists of
bounded holomorphic $\CC$-valued functions on $\RR \times S^1$, so it is
precisely the real $2$-dimensional space of constant functions, implying
$$
\dim_\RR \coker(\dbar) = \dim_\RR \ker(\dbar) - \ind_\RR(\dbar) = 2 - 2 = 0,
$$
so $\mathbf{D}_{u_\gamma}$ is surjective.
\end{proof}

\begin{cor}
\label{cor:contactRegular}
For any contact form $\alpha$ on a closed manifold $M$, there exists a 
comeager subset $\jJ^\reg(\alpha) \subset \jJ(\alpha)$ such that
for every $J \in \jJ^\reg(\alpha)$, all somewhere injective asymptotically
cylindrical $J$-holomorphic
curves in $\RR \times M$ are Fredholm regular.
\qed
\end{cor}

Note that in the setting of Corollary~\ref{cor:contactRegular},
a curve that is not a cover of a trivial cylinder always belongs to a
smooth $1$-parameter family of curves related to each other by
$\RR$-translation, so that the kernel of the linearized Cauchy-Riemann
operator automatically has kernel of dimension at least~$1$.  This
precludes Fredholm regularity for curves of index~$0$, thus:

\begin{cor}
If $\alpha$ is a contact form and $J \in \jJ^\reg(\alpha)$, then all
simple asymptotically cylindrical $J$-holomorphic curves 
$u : (\dot{\Sigma},j) \to (\RR \times M,J)$ other than trivial cylinders
satisfy
$$
\ind(u) \ge 1.
$$
\qed
\end{cor}

The following example shows that the third condition on the
injective point in Theorem~\ref{thm:genericRinvt} cannot be fully
removed in general.

\begin{example}[cf.~Examples~6.6 and~6.16 in Lecture~\ref{lec:cobordisms}]
Assume $(W,\omega)$ is a closed symplectic manifold of dimension $2n-2$
with a periodic time-dependent Hamiltonian $H : S^1 \times W \to \RR$,
and $M := S^1 \times W$ is assigned the stable Hamiltonian structure
$(\Omega,\Lambda) := (\omega + dt \wedge dH,dt)$.  A choice of
$J \in \jJ(\hH)$ is then equivalent to a choice of $t$-dependent family
of $\omega$-compatible almost complex structures $\{J_t \in \jJ(W,\omega)\}_{t \in S^1}$,
and for any $t \in S^1$ and $s \in \RR$, $J_t$-holomorphic curves
$u : (\Sigma,j) \to (W,J_t)$ give rise to $J$-holomorphic curves
$$
\bar{u} : (\Sigma,j) \to (\RR \times M,J) : z \mapsto (s,t,u(z)).
$$
In particular, when $n=2$ one can consider the example where $W = \Sigma$
is a closed surface, so curves of this form exist for any choice of
$J \in \jJ(\hH)$, no matter how generic (remember that the domain complex
structure $j$ is arbitrary, it is not fixed in advance).  If 
$\Sigma$ has genus~$g$ and the map $u : \Sigma \to \Sigma$ has degree~$1$, 
then since $\bar{u}$ has no punctures and satisfies
$c_1([\bar{u}]) = c_1(\bar{u}^*T(\RR\times S^1 \times \Sigma)) =
c_1(T\Sigma) = \chi(\Sigma)$, the index of $\bar{u}$ is
$$
\ind(\bar{u}) = (n-3)\chi(\Sigma) + 2 \chi(\Sigma) = \chi(\Sigma) = 2 - 2g.
$$
This shows that $\bar{u}$ cannot be Fredholm regular unless $g=0$.
\end{example}

Theorem~\ref{thm:genericRinvt} appeared for the first time in the contact
case in \cite{Dragnev}, and alternative proofs have since appeared in
the appendix of \cite{Bourgeois:homotopy} (for cylinders in the contact
case) and in \cite{Wendl:blogTransversality} (under slightly different
assumptions in the stable Hamiltonian setting).  What I will describe
below is a generalization of Bourgeois's proof.

\section{Injective points of the projected curve}
\label{sec:injective}

One point of difficulty in proving transversality in $\RR \times M$
is that in contrast to the setting of Theorem~\ref{thm:genericCobordism}, generic perturbations 
within $\jJ(\hH)$
can never be truly local, i.e.~if you perturb $J$ near a point
$(r,x) \in \RR \times M$, then you are also perturbing it in a neighborhood
of the entire line $\RR \times \{x\}$.  We therefore need to know that we
can find a point $z \in \dot{\Sigma}$ that is the \emph{only} point
where $u : \dot{\Sigma} \to \RR \times M$ passes through such a line;
put another way, we need to know that not only $u = (u_\RR,u_M) :
\dot{\Sigma} \to \RR \times M$ but also the projected map $u_M : \dot{\Sigma} \to M$
is somewhere injective.  The first step in showing this is
Proposition~\ref{prop:piTu} above, as the zeroes of the section
$$
\pi_\xi \circ du \in \Gamma(\Hom_\CC(T\dot{\Sigma},u^*\xi))
$$
are precisely the critical points of $u_M : \dot{\Sigma} \to M$;
everywhere else, $u_M$ is an immersion transverse to the Reeb vector field.
To prove Proposition~\ref{prop:piTu}, we shall use the fact that the vector
fields $\p_r$ and $R$ generate an integrable $J$-invariant distribution
on $\RR \times M$.  Indeed, the zeroes of $\pi_\xi \circ du$ are the points
of tangency with this distribution, hence the result is an immediate 
consequence of the following statement:

\begin{lemma}
Suppose $(W,J)$ is an almost complex manifold,
$\Xi \subset TW$ is a smooth integrable $J$-invariant distribution and
$u : (\Sigma,j) \to (W,J)$ is a connected pseudoholomorphic curve whose image
is not contained in a leaf of the foliation generated by~$\Xi$.
Then all points $z \in \Sigma$ with $\im du(z) \subset \Xi$ are isolated
in~$\Sigma$.
\end{lemma}
\begin{proof}
Statement is local, so assume $(\Sigma,j) = (\DD,i)$ with coordinates
$s + it$, $W = \CC^n$, and $u(0) = 0$.
Let $2m$ denote the real dimension of~$\Xi$, and observe that since $\Xi$
is integrable, we can change coordinates near~$0$ and assume without loss of 
generality that at every point $p \in \CC^n$ near~$0$, $\Xi_p = 
\CC^m \oplus \{0\} \subset \CC^n = T_p \CC^n$.  The $J$-invariance of
$\Xi$ then implies that in coordinates $(w,\zeta) \in \CC^m \times \CC^{n-m}$,
$J$ takes the form
$$
J(w,\zeta) = \begin{pmatrix} J_1(w,\zeta) & Y(w,\zeta) \\ 0 & J_2(w,\zeta) \end{pmatrix},
$$
where $J_1^2$ and $J_2^2$ are both $-\1$, and $J_1 Y + Y J_2 = 0$.
Writing $u(z) = (f(z),v(z)) \in \CC^m \times \CC^{n-m}$, the Cauchy-Riemann
equation $\p_s u + J(u) \p_t u = 0$ is then equivalent to the two equations
\begin{equation}
\label{eqn:twoEquations}
\begin{split}
\p_s f + J_1(f,v) \, \p_t f + Y(f,v) \, \p_t v &= 0, \\
\p_s v + J_2(f,v) \, \p_t v &= 0.
\end{split}
\end{equation}
We have $\im du(z) \subset \Xi$ wherever $\p_s v = \p_t v = 0$; notice that
it suffices to consider the condition $\p_s v = 0$ since
$\p_t v = J_2(f,v)\, \p_s v$.  Differentiating the second equation in
\eqref{eqn:twoEquations} with respect to $s$ gives
$$
\p_s (\p_s v) + J_2(f,v) \, \p_t (\p_s v) + \p_s\left[ J_2(f,v) \right] J_2(f,v)\, \p_s v
= 0,
$$
where in the last term we've substituted $J_2(f,v) \, \p_s v$ for $\p_t v$.
Setting $\bar{J}(z) := J_2(f(z),v(z))$ and $A(z) := \p_s\left[ J_2(f(z),f(z) \right] J_2(f(z),v(z))$,
this becomes a linear Cauchy-Riemann type equation
$\p_s (\p_s v) + \bar{J} \, \p_t (\p_s v) + A (\p_s v) = 0$, so the
similarity principle implies that zeroes of $\p_s v$ are isolated unless
it is identically zero.  The latter would mean $v$ is constant, so $u$
is contained in a leaf of~$\Xi$.
\end{proof}

\begin{lemma}
\label{lemma:noIntOrbits}
Suppose $J \in \jJ(\hH)$, $\gamma : S^1 \to M$ is a closed Reeb orbit,
and $u = (u_\RR,u_M) : (\dot{\Sigma},j) \to (\RR \times M,J)$
is an asymptotically cylindrical $J$-holomorphic curve that is not a cover
of a trivial cylinder.  Then all intersections of the map
$u_M : \dot{\Sigma} \to M$ with the image of the orbit $\gamma$ are isolated.
\end{lemma}
\begin{proof}
The trivial cylinder over $\gamma$ is a $J$-holomorphic curve, so the
statement follows from the fact that two asymptotically cylindrical
$J$-holomorphic curves can only
have isolated intersections unless both are covers of the same simple
curve.
\end{proof}

We can now prove the statement we need about somewhere injectivity
for $u_M : \dot{\Sigma} \to M$.  This result first appeared in
\cite{HWZ:props3}*{Theorem~1.13}.

\begin{prop}
Suppose $J \in \jJ(\hH)$ and 
$$
u = (u_\RR,u_M) : (\dot{\Sigma},j) \to (\RR \times M,J)
$$
is a simple asymptotically cylindrical $J$-holomorphic curve which is not
a trivial cylinder and has only nondegenerate asymptotic orbits.
Then the set of injective points $z \in \dot{\Sigma}$ of the 
map $u_M : \dot{\Sigma} \to M$ for which $u_M(z)$ is not contained in any
of the asymptotic orbits of $u$ is open and dense.
\end{prop}
\begin{proof}
Openness is clear, so our main task is to prove density.
The idea is first to show via elementary topological arguments that if the 
set of injective
points is not dense, then $\dot{\Sigma}$ contains two disjoint open sets on
which $u_M$ is an embedding with identical images.  We will then conclude
from this that if $u$ is simple, it must be equivalent to one of its
nontrivial $\RR$-translations, and the latter is impossible for an
asymptotically cylindrical curve.

\textbf{Step~1}: We begin by harmlessly removing some discrete sets of
points in $\dot{\Sigma}$ that would make the subsequent arguments more
complicated.  Let 
$$
P \subset M
$$
denote the union of the images of the asymptotic orbits
of~$u$, a finite disjoint union of circles.  Lemma~\ref{lemma:noIntOrbits}
implies that $u_M^{-1}(P)$ is a discrete subset of~$\dot{\Sigma}$.
By Proposition~\ref{prop:piTu}, there is also a discrete set
$Z \subset \dot{\Sigma} \setminus u_M^{-1}(P)$ containing all points 
$z \not\in u_M^{-1}(P)$ where $\pi_\xi \circ du(z) = 0$, and we claim that
$$
Z' := u_M^{-1}(u_M(Z))
$$
is a discrete subset of $\dot{\Sigma} \setminus u_M^{-1}(P)$.  Indeed,
$u_M(Z)$ is a discrete subset of $M \setminus P$ since the points in
$Z$ can only accumulate at infinity,\footnote{Actually the asymptotic
formula of \cite{HWZ:props1} implies that both $Z$ and $u_M^{-1}(P)$ are
always finite for curves that are not covers of trivial cylinders, but we
do not need to use that here.} hence accumulation points of $u_M(Z) \subset M$
can occur only in~$P$.  For each individual point $p \in u_M(Z)$, the fact
that $p \not\in P$ implies $u_M^{-1}(p)$ is compact, and it consists of a
discrete (and therefore finite) set of points with
$\pi_\xi \circ du(z) = 0$, plus possibly some other points where 
$\pi_\xi \circ du(z) \ne 0$, but $u_M$ is an embedding near each point of
the latter type, so that these points of $u_M^{-1}(p)$ must always be
isolated and are therefore also finite in number.  This proves the claim,
and we conclude that
$$
\ddot{\Sigma} := \dot{\Sigma} \setminus \left(u_M^{-1}(P) \cup Z' \right)
$$
an open and dense subset of $\dot{\Sigma}$, as it
is obtained by removing a discrete subset from the open and dense subset
$\dot{\Sigma} \setminus u_M^{-1}(P)$.  To prove the proposition, it will now
suffice to prove that
the set of points $z \in \ddot{\Sigma}$ which are injective points of
$u_M : \dot{\Sigma} \to M$ is dense in~$\ddot{\Sigma}$.  We shall
argue by contradiction and assume from now on that density fails.

\textbf{Step~2}: We will find two open subsets
$\uU , \vV \subset \dot{\Sigma}$ such that $u_M$ restricts to an embedding
on both, but
$$
\uU \cap \vV = \emptyset \quad \text{ and }\quad u_M(\uU) = u_M(\vV).
$$
Indeed, assume the set of injective points of $u_M$ lying in $\ddot{\Sigma}$
is not dense in~$\ddot{\Sigma}$.  Then
there exists a point $z_0 \in \ddot{\Sigma}$
with a closed neighborhood $\dD(z_0) \subset \ddot{\Sigma}$ such that
no $z \in \dD(z_0)$ is an injective point.  Since $z \in \ddot{\Sigma}$
implies $\pi_\xi \circ du(z) \ne 0$, this means that for every
$z \in \dD(z_0)$, there exists $\zeta \in \dot{\Sigma} \setminus \{z\}$
with $u_M(z) = u_M(\zeta)$, and the definition of $\ddot{\Sigma}$ implies
$\zeta$ is also in~$\ddot{\Sigma}$, hence $\pi_\xi \circ du(\zeta) \ne 0$ and
$u_M$ is a local embedding near~$\zeta$.  Since $u(z) \not\in P$ and
$u_M$ maps $\dot{\Sigma} \setminus u_M^{-1}(P)$ properly to $M \setminus P$,
we also conclude that $u_M^{-1}(u_M(z))$ is finite.  Now suppose
$u_M^{-1}(u_M(z_0)) = \{z_0,\zeta_1,\ldots,\zeta_m\}$, and let
$\dD(\zeta_j) \subset \ddot{\Sigma}$ for $j=1,\ldots,m$ denote closed
neighborhoods on which $u_M$ is an embedding.  We claim that after
possibly shrinking $\dD(z_0)$, we can assume
$$
u_M(\dD(z_0)) \subset \bigcup_{j=1}^m u_M(\dD(\zeta_j).
$$
Let us first shrink $\dD(z_0)$ so that $u_M$ is an embedding on $\dD(z_0)$,
which is possible since $\pi_\xi \circ du(z_0) \ne 0$.
Then if the claim is false, there exists a sequence $z_\nu \in \dD(z_0)$ of
noninjective points with $z_\nu \to z_0$, hence there is also a sequence
$z_\nu' \in \ddot{\Sigma} \setminus \dD(z_0)$ with $u_M(z_\nu) = u_M(z_\nu')$ 
but $z_\nu'$ not converging to any of
$\zeta_1,\ldots,\zeta_m$.  But since $u_M(z_\nu') \to u_M(z_0) \not\in P$,
the points $z_\nu'$ are confined to a compact subset of $\dot{\Sigma}$
and therefore have a subsequence $z_\nu' \to z_\infty' \in \dot{\Sigma}$
with $u_M(z_\infty') = u_M(z_0)$.  The limit cannot be $z_0$ itself since
$z_\nu' \not\in \dD(z_0)$,
thus $z_\infty'$ must be one of the $\zeta_1,\ldots,\zeta_m$, and we have
a contradiction.  We claim next that at least one of the sets 
$u_M(\dD(z_0)) \cap u_M(\dD(\zeta_j))$ has nonempty interior.  This is a 
simple exercise in metric space topology: it can be reduced to the fact 
that if $X$ is a metric space with closed subsets $V , W \subset X$ that
both have empty interior (meaning no open subset of $X$ is contained in
$V$ or $W$), then $V \cup W$ also has empty interior.  Since the subsets
$u_M(\dD(z_0)) \cap u_M(\dD(\zeta_j)) \subset u_M(\dD(z_0))$ for 
$j=1,\ldots,m$ are all closed but their union is $u_M(\dD(z_0))$, they
cannot all have empty interior.  This achieves the goal of Step~2.

\textbf{Step~3}: We show that $u$ is biholomorphically equivalent to one of
its $\RR$-translations
$$
\tau \cdot u := (u_\RR + \tau,u_M) : \dot{\Sigma} \to \RR \times M
$$
for $\tau \in \RR \setminus \{0\}$.  To see this, note that for 
$J \in \jJ(\hH)$, the nonlinear Cauchy-Riemann equation $Tu \circ j =
J(u) \circ Tu$ is equivalent to the two equations
\begin{equation}
\label{eqn:twoMoreEquations}
\begin{split}
d u_\RR &= u_M^*\lambda \circ j,\\
\pi_\xi \circ T u_M \circ j &= J(u_M) \circ \pi_\xi \circ T u_M.
\end{split}
\end{equation}
Since $\pi_\xi \circ T u_M : \dot{\Sigma} \to u_M^*\xi$ is injective
everywhere on the neighborhoods $\uU$ and $\vV$, the second equation
determines $j$ in terms of $J$ on each of these regions; in particular,
the identification of $u_M(\uU)$ with $u_M(\vV)$ provides a biholomorphic
map of $\vV$ to $\uU$ so that $u|_{\uU}$ and $u|_{\vV}$ may be regarded as
two $J$-holomorphic maps from the same Riemann surface which differ only
in the $\RR$-factor.  But with $j$ and $u_M$ both fixed, the first
equation in \eqref{eqn:twoMoreEquations} determines $du_\RR$ and thus
determines $u_\RR$ up to the addition of a constant $\tau \in \RR$.
If $\tau = 0$, this means $u$ has two disjoint regions on which its
images are identical, contradicting the assumption that $u$ is simple.
Thus $\tau \ne 0$, and since two distinct simple curves can only intersect
each other at isolated points, we conclude $u = \tau \cdot u$ up to
parametrization.

\textbf{Step~4}: We now derive a contradiction.  The relation
$u = \tau \cdot u$ implies that in fact $u = k\tau \cdot u$ for every
$k \in \ZZ$, so we obtain a diverging sequence of $\RR$-translations
$\tau_k \to \infty$ such that $u$ and $\tau_k \cdot u$ always have identical
images in $\RR \times M$.  It follows that for some point $z \in \dot{\Sigma}$ 
with $u(z) = (r,x)$ where $x$ is not contained in any of the asymptotic
orbits of~$u$, the points
$(r - \tau_k,x)$ are all in the image of $u$ as $\tau_k \to \infty$.
But this contradicts the asymptotically cylindrical behavior of~$u$.
\end{proof}

\section{Smoothness of the universal moduli space}

The overall outline of the proof of Theorem~\ref{thm:genericRinvt}
is the same as for Theorem~\ref{thm:genericCobordism}: one needs to define a suitable space
$\jJ_\uU^\varepsilon$ of perturbed almost complex structures, giving rise to 
a universal moduli space $\mM^*(\jJ_\uU^\varepsilon)$ that is a smooth 
Banach manifold, and then apply the Sard-Smale theorem to conclude that
generic elements of $\jJ_\uU^\varepsilon$ are regular values of
the projection $\mM^*(\jJ_\uU^\varepsilon) \to \jJ_\uU^\varepsilon :
(u,J) \mapsto J$.  If $\jJ_\uU^\varepsilon$ is a space of
$C_\varepsilon$-perturbed almost complex structures, then in the final step
one can use the Taubes trick as in \S\ref{sec:Taubes} to transform the 
genericity result in $\jJ_\uU^\varepsilon$ into a genericity result
within the space $\jJ(\hH)$ of smooth almost complex structures.
The only step that differs meaningfully from what we've already discussed
is the smoothness of the universal moduli space, so let us focus on this
detail.

Assume $\Jref \in \jJ(\hH)$ with $\Jref = \Jfix$ outside $\RR \times \uU$,
and $\jJ_\uU^\varepsilon$ is a Banach manifold of $C_\varepsilon$-small
perturbations of $\Jref$ in $\jJ(\hH)$ that are also fixed outside of 
$\RR \times \uU$.
The relevant universal moduli space is then defined by
\begin{equation*}
\begin{split}
\mM^*(\jJ_\uU^\varepsilon) := \big\{ (u,J) \ \big|\ 
&\text{$J \in \jJ_\uU^\varepsilon$, $u \in \mM(J)$ and}\\
&\text{$u : \dot{\Sigma} \to \RR \times M$ has an injective point 
$z \in \dot{\Sigma}$ with} \\
&\text{$u(z) \in \RR \times \uU$ and $\im\left(\pi_\xi \circ du(z)\right)
\cap \ker \left( d\lambda|_\xi \right) = \{0\}$} \big\}.
\end{split}
\end{equation*}
Notice that both of the constraints satisfied by $u$ at the injective point
are open.  The local structure of $\mM^*(\jJ_\uU^\varepsilon)$ near an
element $(u_0,J_0)$ with representative $u_0 : (\dot{\Sigma},j_0) \to
(\RR \times M,J_0)$ can again
be described via the zero set of a smooth section
$$
\dbar : \tT \times \bB^{k,p,\delta} \times \jJ_\uU^\varepsilon \to
\eE^{k-1,p,\delta} : (j,u,J) \mapsto Tu \circ J \circ Tu \circ j,
$$
where $\tT$ is a Teichm\"uller slice through $j_0$, and it suffices to
show that the linearization
$$
\mathbf{L} : T_{u_0}\bB^{k,p,\delta} \oplus T_{J_0} \jJ_\uU^\varepsilon \to
\eE_{(j_0,u_0,J_0)}^{k-1,p,\delta} : (\eta,Y) \mapsto
\mathbf{D}_{u_0} \eta + Y \circ Tu_0 \circ j_0
$$
is always surjective.  As usual, here we're assuming $k \in \NN$,
$1 < p < \infty$, and the exponential weight $\delta > 0$ is small but positive
so that $\mathbf{D}_{u_0}$ is Fredholm.  The image of $\mathbf{L}$ is then
closed, and focusing on the $k=1$ case, if $\mathbf{L}$ is not surjective
then there exists a nontrivial
element $\theta \in L^{q,-\delta}(\overline{\Hom}_\CC(T\dot{\Sigma},u_0^*T(\RR \times M)))$
such that
\begin{equation}
\label{eqn:twoConditions2}
\begin{split}
\langle \mathbf{D}_{u_0} \eta , \theta \rangle_{L^2} = 0 & \text{ for all
$\eta \in W^{1,p,\delta}(u_0^*T(\RR \times M)) \oplus V_\Gamma$}, \\
\langle Y \circ Tu_0 \circ j_0 , \theta \rangle_{L^2} = 0 & \text{ for all
$Y \in T_{J_0}\jJ_\uU^\varepsilon$}.
\end{split}
\end{equation}
The first condition implies via elliptic regularity and the similarity
principle that $\theta$ is smooth and has only isolated zeroes.
So far this is all the same as in the proof of Theorem~\ref{thm:genericCobordism}, but the next
step is trickier: since perturbing $J_0$ within $\jJ(\hH)$ only changes
the action of the almost complex structure on $\xi$ but not on the
trivial subbundle generated by $\p_r$ and $R$, it is not clear whether the
range of values allowed for $Y$ is large enough to force
$\langle Y \circ Tu_0 \circ j_0 , \theta \rangle_{L^2} > 0$.

To overcome this, let us decompose everything in this picture with respect
to the natural splitting
$$
T(\RR \times M) = \epsilon \oplus \xi,
$$
where $\epsilon$ denotes the trivial line bundle spanned by $\p_r$ and~$R$.
In particular, the domain and target bundles of the Cauchy-Riemann type
operator $\mathbf{D}_{u_0}$ now split as
\begin{equation*}
\begin{split}
u_0^*T(\RR \times M) &= u_0^*\epsilon \oplus u_0^*\xi, \\
\overline{\Hom}_\CC(T\dot{\Sigma},u_0^*T(\RR \times M)) &=
\overline{\Hom}_\CC(T\dot{\Sigma},u_0^*\epsilon) \oplus
\overline{\Hom}_\CC(T\dot{\Sigma},u_0^*\xi),
\end{split}
\end{equation*}
and we shall write $\eta = (\eta^\epsilon,\eta^\xi)$ and
$\theta = (\theta^\epsilon,\theta^\xi)$ accordingly.  This gives a block
decomposition of $\mathbf{D}_{u_0}$ as
$$
\mathbf{D}_{u_0} \eta = \begin{pmatrix} (\mathbf{D}_{u_0}\eta)^\epsilon \\
(\mathbf{D}_{u_0}\eta)^\xi \end{pmatrix}
= \begin{pmatrix}
\mathbf{D}_{u_0}^\epsilon & \mathbf{D}_{u_0}^{\epsilon\xi} \\
\mathbf{D}_{u_0}^{\xi\epsilon} & \mathbf{D}_{u_0}^\xi
\end{pmatrix}
\begin{pmatrix} \eta^\epsilon \\ \eta^\xi \end{pmatrix}.
$$
It is easy to verify that $\mathbf{D}_{u_0}^\epsilon$ and
$\mathbf{D}_{u_0}^\xi$ each satisfy suitable Leibniz rules and are thus
Cauchy-Riemann type operators on $u_0^*\epsilon$ and $u_0^*\xi$ respectively,
while the off-diagonal terms are both tensorial, i.e.~zeroth-order
operators.  Since perturbations of $J_0$ in $\jJ(\hH)$ only change its
action on $\xi$, $Y \in T_{J_0}\jJ_\uU^\varepsilon$ now takes the block form
$$
Y = \begin{pmatrix} 0 & 0 \\
0 & Y^\xi \end{pmatrix},
$$
where $Y^\xi$ is a $C_\varepsilon$-small section of the bundle
$\overline{\End}_\CC(\xi,J_0)$ over~$M$.  
Assuming the $L^2$-pairings are defined so as to respect these splittings,
the second condition in \eqref{eqn:twoConditions2} now becomes
$$
\langle Y^\xi \circ \pi_\xi \circ Tu_0 \circ j_0 , \theta^\xi \rangle_{L^2} = 0,
$$
and given any injective point $z_0 \in \dot{\Sigma}$ of
$(u_0)_M : \dot{\Sigma} \to M$ satisfying $u_0(z_0) \in \RR \times \uU$, we have
enough freedom to choose $Y^\xi$ near $\RR \times \{u_0(z_0)\}$ such that this
pairing becomes positive unless
$$
\theta^\xi = 0 \quad \text{ near~$z_0$}.
$$
It remains to show that $\theta^\epsilon$ also vanishes near~$z_0$, which will 
contradict the fact that $\theta$ only has isolated zeroes.  To this end,
notice that the first condition in \eqref{eqn:twoConditions2} implies
via separate choices of the components $\eta^\epsilon$ and $\eta^\xi$ 
with support near $z_0$ that
\begin{equation}
\label{eqn:anotherTwoConditions}
\begin{split}
\langle \mathbf{D}_{u_0}^\epsilon \eta^\epsilon , \theta^\epsilon \rangle_{L^2} = 0 & \text{ for all
$\eta^\epsilon$ supported near~$z_0$}, \\
\langle \mathbf{D}_{u_0}^{\epsilon\xi} \eta^\xi , \theta^\epsilon \rangle_{L^2} = 0 &
\text{ for all $\eta^\xi$ supported near~$z_0$}.
\end{split}
\end{equation}
The first of these two conditions gives no new
information, since we already know that $\theta = (\theta^\epsilon,0)$
solves an anti-Cauchy-Riemann equation.  To get some information out of the
second condition, we will need an explicit formula for
$\mathbf{D}_{u_0}^{\epsilon\xi}$.

\begin{lemma}
\label{lemma:thisLemma}
The tensorial operator
$\mathbf{D}_{u_0}^{\epsilon\xi} : u_0^*\xi \to 
\overline{\Hom}_\CC(T\dot{\Sigma},u_0^*\epsilon)$ takes the form
$$
\mathbf{D}_{u_0}^{\epsilon\xi} \eta^\xi = \left[ -d\lambda\big(\eta^\xi , 
J_0^\xi \circ \pi_\xi \circ Tu(\cdot)\big) \right] \p_r +
\left[ d\lambda \big( \eta^\xi , \pi_\xi \circ Tu(\cdot)\big) \right] R.
$$
\end{lemma}
\begin{proof}
As a preliminary step, notice that $-dr \circ J = \lambda$ for any
$J \in \jJ(\hH)$; indeed, the conditions $J(\xi) = \xi \subset \ker dr$
and $J\p_r = R$ imply that these two $1$-forms have matching values on
$\p_r$, $R$ and~$\xi$.  As a consequence, $\lambda \circ J_0 = dr$, so
in particular $\lambda \circ J_0$ is closed.

Choosing local holomorphic coordinates $(s,t)$ in an arbitrary neighborhood
in $\dot{\Sigma}$, we have
$$
(\mathbf{D}_{u_0}^{\epsilon\xi} \eta^\xi) \p_s =
dr\big( (\mathbf{D}_{u_0} \eta^\xi)\p_s \big) \, \p_r +
\lambda\big( (\mathbf{D}_{u_0} \eta^\xi) \p_s \big) \, R.
$$
Extend $u_0 : \dot{\Sigma} \to \RR \times M$ to a 
smooth $1$-parameter family of maps
$\{u_\rho : \dot{\Sigma} \to \RR \times M\}_{\rho \in \RR}$ with $\p_\rho u_\rho|_{\rho=0} =
\eta^\xi \in \Gamma(u_0^*\xi)$.  Then by the definition of the linearized
Cauchy-Riemann operator,
$$
(\mathbf{D}_{u_0} \eta^\xi) \p_s = \left. \nabla_\rho \left(
\p_s u_\rho + J_0(u_\rho) \p_t u_\rho \right)\right|_{\rho=0},
$$
for any choice of connection $\nabla$ on $\RR \times M$.  Since
$\p_s u_0 + J_0(u_0) \p_t u_0 = 0$, we find
\begin{equation*}
\begin{split}
\lambda\big( (\mathbf{D}_{u_0} \eta^\xi) \p_s \big) &=
\lambda\big( \nabla_\rho \left.\left( \p_s u_\rho + J_0(u_\rho) \p_t u_\rho\right)\right|_{\rho=0} \big)
= \left.\p_\rho \left[ \lambda(\p_s u_\rho + J_0(u_\rho) \p_t u_\rho) \right]\right|_{\rho=0} \\
&= \left.\p_\rho \left[ \lambda(\p_s u_\rho) \right]\right|_{\rho=0}
+ \left.\p_\rho \left[ (\lambda \circ J_0)(\p_t u_\rho) \right]\right|_{\rho=0} \\
&= d\lambda(\eta^\xi, \p_s u) + d(\lambda \circ J_0)(\eta^\xi, \p_t u) \\
&= d\lambda(\eta^\xi, \pi_\xi \p_s u),
\end{split}
\end{equation*}
where we've used the formula
$$
d\lambda(X,Y) = \Lie_X\left[ \lambda(Y) \right] - \Lie_Y\left[ \lambda(X) \right]
- \lambda([X,Y])
$$
and eliminated several terms using the fact that $\lambda(\eta^\xi) = 
\lambda(J_0 \eta^\xi) = 0$ since $\eta^\xi$ is valued in~$\xi$,
plus $d(\lambda \circ J_0) = 0$.  A similar computation gives
$$
dr\big( (\mathbf{D}_{u_0} \eta^\xi) \p_s \big) = 
-d\lambda(\eta^\xi, \pi_\xi \p_t u) = -d\lambda(\eta^\xi , J_0 \circ 
\pi_\xi \p_s u),
$$
so removing the local coordinates from the picture produces the stated
formula.
\end{proof}

The following exercise in symplectic linear algebra shows that this bundle
map $u_0^*\xi \to \overline{\Hom}_\CC(T\dot{\Sigma},u_0^*\epsilon)$ 
is surjective on any fiber over a point
$z$ with $\pi_\xi \circ du_0(z) \ne 0$.  (If you have no patience for the
exercise, just convince yourself that it's true whenever
$d\lambda|_\xi$ is nondegenerate and tames $J|_\xi$, i.e.~the contact case.)

\begin{exercise}
\label{Ex:fibSurj}
Assume $V$ is a finite-dimensional vector space, $X, Y \subset V$ are 
linearly independent vectors, and $\Omega$ is an alternating bilinear form
on~$V$.  Show that the real-linear map
$$
A : V \to \CC : v \mapsto \Omega(v,X) + i \Omega(v,Y)
$$
is surjective if and only if $\Span(X,Y) \cap \ker \Omega = \{0\}$. \\
\textsl{Hint: Under the latter condition, one loses no generality by 
replacing $V$ with a subspace that is complementary to $\ker \Omega$ and 
contains $\Span(X,Y)$, in which case $(V,\Omega)$ becomes a symplectic 
vector space.  Now consider the restriction of $A$ to a $2$-dimensional 
subspace transverse to the symplectic complement of $\Span(X,Y)$.}
\end{exercise}

The conclusion of this discussion is that unless $\theta^\epsilon$ vanishes
near~$z_0$, $\eta^\xi$ can be chosen
with support near~$z_0$ so that 
$\langle \mathbf{D}_{u_0} \eta^\xi,\theta^\epsilon \rangle_{L^2} > 0$,
violating the second condition in \eqref{eqn:anotherTwoConditions}.
This proves that $\theta$ vanishes altogether near~$z_0$ and thus,
by unique continuation, $\theta \equiv 0$, a contradiction.

We've proved that the universal moduli space is smooth as claimed.
Since the rest of the proof of Theorem~\ref{thm:genericRinvt} is the same
as in the non-$\RR$-invariant case, we leave those details to the reader.

\begin{remark}
You may have noticed that in both this and the previous lecture, our proof
that the universal moduli space is smooth relied on a surjectivity result
that was actually stronger than needed: in both cases, we needed to prove
that an operator of the form
$$
T_{j_0}\tT \oplus T_{u_0}\bB^{k,p,\delta} \oplus T_{J_0}\jJ_\uU^\varepsilon
\stackrel{\mathbf{L}}{\longrightarrow} \eE_{(j_0,u_0,J_0)}^{k-1,p,\delta}
$$
was surjective, but we ended up proving that its restriction to the
smaller domain $T_{u_0}\bB^{k,p,\delta} \oplus T_{J_0}\jJ_\uU^\varepsilon$
is already surjective.  This technical detail hints at a stronger result
that can be proved using these methods: one can show that not only is
$\mM^*(\jJ_\uU^\varepsilon)$ smooth but also the \defin{forgetful map}
\begin{equation*}
\begin{split}
\mM^*(\jJ_\uU^\varepsilon) &\to \mM_{g,k_+ + k_- + m} \\
([(\Sigma,j,\Gamma^+,\Gamma^-,\Theta,u)],J) &\mapsto
[(\Sigma,j,\Gamma^+ \cup \Gamma^- \cup \Theta)]
\end{split}
\end{equation*}
sending a $J$-holomorphic curve to its underlying domain in the moduli space
of Riemann surfaces is a submersion, cf.~the blog post
\cite{Wendl:blogForgetful} and its sequel.  One can use this to prove generic
transversality results for spaces of $J$-holomorphic curves whose domains
are constrained within the moduli space of Riemann surfaces, which can be
used to define more elaborate algebraic structures on SFT, e.g.~this idea
plays a very prominent role in the study of Gromov-Witten invariants.
\end{remark}

\psfrag{Br(p)}{$B_r(p)$}
\psfrag{p}{$p$}
\psfrag{u(S)}{$u(\Sigma)$}
\psfrag{u(S')}{$u(\Sigma')$}

\chapter{Asymptotics and compactness}
\label{lec:compactness}

\minitoc

\vspace{12pt}

Moduli spaces of pseudoholomorphic curves are generally not compact, but they
have natural \emph{compactifications}, obtained by allowing certain types
of curves with singular behavior.  For closed holomorphic curves, this 
fact is known as \emph{Gromov's compactness theorem}, and our main goal in
this lecture is to state its generalization to punctured curves, which
is usually called the \emph{SFT compactness theorem}.  The theorem was
first proved in \cite{SFTcompactness} (see also \cite{CieliebakMohnke:compactness}
for an alternative approach), and we do not have space here to present a
complete proof, but we can still describe the main geometric and
analytical ideas behind it.

The overarching theme of this lecture is the notion of \emph{bubbling},
of which we will see several examples.  Bubbling arises in a natural way
from elliptic regularity: recall that in Lecture~\ref{lec:local}, we proved that
whenever $kp > 2$, any uniformly $W^{k,p}$-bounded sequence $u_\nu$ 
of holomorphic curves is also uniformly $C^m_\loc$-bounded for every
$m \ge \NN$ (cf.~Theorem~\ref{thm:regularity}).  The Arzel\`a-Ascoli theorem implies that
such sequences have $C^\infty_\loc$-convergent subsequences, and this is true
in particular whenever $u_\nu$ is uniformly $C^1$-bounded, as a $C^1$-bound
implies a $W^{1,p}$-bound with $p > 2$.  Let us take note of this fact
for future use:

\begin{prop}
\label{prop:C1bound}
If $(W,J_\nu)$ is a sequence of almost complex manifolds with $J_\nu \to J$
in $C^\infty$, then any uniformly $C^1$-bounded sequence of 
$J_\nu$-holomorphic maps $u_\nu : \DD \to W$ has a subsequence convergent
in $C^\infty_\loc$ on~$\intDD$.
\end{prop}

If one wants to prove
compactness for a moduli space of $J$-holomorphic curves, it therefore
suffices in general to establish a $C^1$-bound.  The catch is, of course, that the
first derivatives of $u_\nu$ might \emph{not} be uniformly bounded, and
this is when interesting things are seen to happen: while the sequence
$u_\nu$ is not compact, it turns out that it becomes compact after removing
finitely many points from its domain, and near those points one can take
a sequence of reparametrizations to find additional nontrivial
holomorphic curves in the limit, the so-called ``bubbles''.  This is one of
the ways that the ``nodal'' curves in Gromov's compactness
theorem can arise, and we will see the same phenomenon at work in
several other contexts as well.

\section{Removal of singularities}
\label{sec:singularities}

As an important tool for use in the rest of this lecture, we begin with the
following result from \cite{Gromov}:

\begin{thm}[Gromov's removable singularity theorem]
\label{thm:removable}
Assume $(W,\omega)$ is a symplectic manifold with a tame almost complex
structure $J$, and $u : \DD \setminus \{0\} \to W$ is a $J$-holomorphic
curve that has its image contained in a compact subset of $W$ and satisfies
$$
\int_{\DD \setminus\{0\}} u^*\omega < \infty.
$$
Then $u$ admits a smooth extension to~$\DD$.
\end{thm}

\begin{figure}
\includegraphics{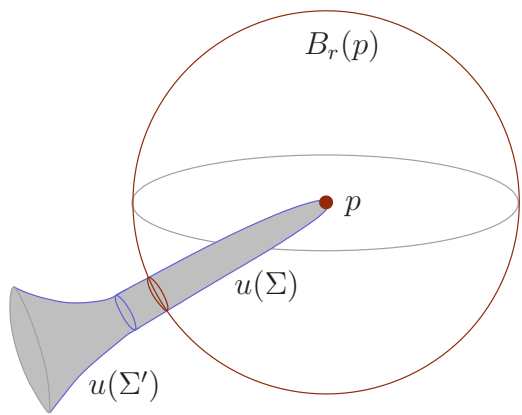}
\caption{\label{fig:monotonicity} The intersection of a $J$-holomorphic curve
$u$ with an open ball $B_r(p)$ defines a proper map $\Sigma \to B_r(p)$.
The monotonicity lemma prevents this map from having
arbitrarily small area if it passes through~$p$.}
\end{figure}

We will prove the slightly weaker statement that $u$ has a \emph{continuous}
extension.  If $\dim_\RR W = 2$, then the smooth extension follows from
this by classical complex analysis; in higher dimensions, one can instead
apply results on local elliptic regularity, see e.g.~\cite{McDuffSalamon:Jhol}.
We will use as a black box the following additional result from \cite{Gromov},
which is closely related to a standard result about minimal surfaces:

\begin{thmu}[Gromov's monotonicity lemma \cite{Gromov}]
Suppose $(W,\omega)$ is a compact symplectic manifold (possibly with boundary),
$J$ is an $\omega$-tame almost complex structure, and $B_r(p) \subset W$ 
denotes the open ball of radius $r > 0$ about $p \in W$ with respect to the
Riemannian metric $g(X,Y) := \frac{1}{2} \omega(X,JY) + \frac{1}{2} \omega(Y,JX)$.
Then there exist constants $c , R > 0$ such that for all $r \in (0,R)$ and
$p \in W$ with $B_r(p) \subset W$, every proper non-constant
$J$-holomorphic curve $u : (\Sigma,j) \to (B_r(p),J)$ passing through $p$ satisfies
$$
\int_\Sigma u^*\omega \ge c r^2.
$$
\end{thmu}
In the above statement, $(\Sigma,j)$ is assumed to be an arbitrary 
(generally noncompact) Riemann surface \emph{without boundary}.  In 
applications, one typically has a larger (e.g.~closed or punctured) domain 
$\Sigma'$ in the picture, and $\Sigma$ is defined to be the connected
component of $u^{-1}(B_r(p)) \subset \Sigma'$ containing some point 
$z \in u^{-1}(p)$.  The main message of the theorem is that $u$ must use up
at least a certain amount of energy for every ball whose center it passes
through, so e.g.~the portion of the curve passing through $B_r(p)$ cannot
become arbitrarily ``thin'' as in Figure~\ref{fig:monotonicity}.

Returning to the removable singularity theorem, we shall use
the biholomorphic map
$$
Z_+ := [0,\infty) \times S^1 \to \DD \setminus \{0\} : (s,t) \mapsto
e^{-2\pi(s+it)}
$$
to transform $J$-holomorphic maps $\DD \setminus \{0\} \to W$ into maps
$Z_+ \to W$, and the goal will be to show that whenever
such a map $u$ has precompact
image and satisfies $\int_{Z_+} u^*\omega < \infty$, there exists
a point $p \in W$ such that
\begin{equation}
\label{eqn:unifConv}
u(s,\cdot) \to p \qquad \text{ in } \quad \text{$C^\infty(S^1,W)$ as $s \to \infty$}.
\end{equation}
Fix the obvious
flat metric on $Z_+$ and any Riemannian metric on $W$ in order to define
norms such as $|du(s,t)|$ for $(s,t) \in Z_+$.

\begin{lemma}
\label{lemma:boundedDeriv}
There exists a constant $C > 0$ such that $|du(s,t)| \le C$ for all
$(s,t) \in Z_+$.
\end{lemma}
\begin{proof}[Proof, part~1]
\renewcommand{\qedsymbol}{\textsc{proof to be continued\ldots}}
Arguing by contradiction, suppose there exists a sequence $z_k = (s_k,t_k)
\in Z_+$ with $|du(z_k)| =: R_k \to \infty$.  Choose a sequence of positive
numbers $\epsilon_k > 0$ that converge to zero but not too fast, so that
$\epsilon_k R_k \to \infty$.  We then consider the sequence of reparametrized
maps
$$
v_k : \DD_{\epsilon_k R_k} \to W : z \mapsto u(z_k + z/R_k).
$$
These are also $J$-holomorphic since $z \mapsto z_k + z / R_k$ is holomorphic,
and the values of $v_k$ depend only on the values of $u$ over the
$\epsilon_k$-disk about~$z_k$.  Notice that since $s_k \to \infty$ and
$\epsilon_k \to 0$, we are free to assume that all of these $\epsilon_k$-disks
are disjoint; moreover, tameness of $J$ implies $u^*\omega \ge 0$ and
$v_k^*\omega \ge 0$, thus
$$
\sum_k \int_{\DD_{\epsilon_k R_k}} v_k^*\omega = \sum_k
\int_{\DD_{\epsilon_k}(z_k)} u^*\omega \le \int_{Z_+} u^*\omega < \infty,
$$
implying
\begin{equation}
\label{eqn:energyDisappearing}
\int_{\DD_{\epsilon_k R_k}} v_k^*\omega \to 0 \quad \text{ as }\quad
k \to \infty.
\end{equation}
We would now like to say something about a limit of the maps $v_k$ as
$k \to \infty$, but this will require a brief pause in the proof, as we
don't yet have quite enough information to do so.  We know that the $v_k$
are uniformly $C^0$-bounded since $u(Z_+)$ is contained in a compact subset.
It would be ideal if we also had a uniform $C^1$-bound, as then elliptic
regularity (Prop.~\ref{prop:C1bound}) would give a $C^\infty_\loc$ convergent
subsequence on the union of all the domains $\DD_{\epsilon_k R_k}$,
i.e.~on the entire plane.  We have
$$
d v_k(z) = \frac{1}{R_k} du(z_k + z/R_k),
$$
hence $|d v_k(0)| = 1$, but we will need to know more about
$|du|$ on the rest of $\DD_{\epsilon_k}(z_k)$ in order to deduce a
$C^1$-bound for $v_k$ on all of~$\DD_{\epsilon_k R_k}$.  We'll come back
to this in a moment.
\end{proof}

Here is the auxiliary lemma that is needed to complete the above proof:

\begin{lemma}[Hofer]
\label{lemma:Hofer}
Suppose $(X,d)$ is a complete metric space, $g : X \to [0,\infty)$ is
continuous, $x_0 \in X$ and $\epsilon_0 > 0$.  Then there exist
$x \in X$ and $\epsilon > 0$ such that,
\begin{enumerate}
\renewcommand{\labelenumi}{(\alph{enumi})}
\item $\epsilon \le \epsilon_0$,
\item $g(x) \epsilon \ge g(x_0) \epsilon_0$,
\item $d(x,x_0) \le 2\epsilon_0$, and
\item $g(y) \le 2 g(x)$ for all $y \in \overline{B_\epsilon(x)}$.
\end{enumerate}
\end{lemma}
\begin{proof}
If there is no $x_1 \in \overline{B_{\epsilon_0}(x_0)}$ such that
$g(x_1) > 2 g(x_0)$, then we can set $x = x_0$ and $\epsilon = \epsilon_0$
and are done.  If such a point~$x_1$ does exist, then we
set $\epsilon_1 := \epsilon_0 / 2$ and repeat the above process for the
pair $(x_1,\epsilon_1)$: that is, if there is no $x_2 \in
\overline{B_{\epsilon_1}(x_1)}$ with $g(x_2) > 2 g(x_1)$, we set
$(x,\epsilon) = (x_1,\epsilon_1)$ and are finished, and 
otherwise define $\epsilon_2 = \epsilon_1 / 2$ and repeat for
$(x_2,\epsilon_2)$.  This process must eventually terminate, as otherwise
we obtain a Cauchy sequence $x_n$ with $g(x_n) \to \infty$, which is
impossible if~$X$ is complete.
\end{proof}

\begin{proof}[Proof of Lemma~\ref{lemma:boundedDeriv}, part~2]
Applying Lemma~\ref{lemma:Hofer} to $X = Z_+$ with $g(z) = |du(z)|$,
we can replace the original sequences $\epsilon_k$ and $z_k$ with
new sequences for which all the previously stated properties still hold,
but additionally,
$$
|du(z)| \le 2 |du(z_k)| \quad \text{ for all } \quad
z \in \DD_{\epsilon_k}(z_k).
$$
Our sequence of reparametrizations $v_k$ then satisfies
$$
| d v_k(z) | \le 2 \quad \text{ for all } \quad
z \in \DD_{\epsilon_k R_k},
$$
so by elliptic regularity, $v_k$ has a subsequence convergent in
$C^\infty_\loc(\CC)$ to a $J$-holomorphic map
$$
v_\infty : \CC \to W
$$
which is not constant since $|d v_\infty(0)| = \lim_{k \to \infty} |d v_k(0)| = 1$.
Informally, we say that the blow-up of the derivatives at $z_k$ has caused
a plane to ``bubble off''.
However, \eqref{eqn:energyDisappearing} implies that for every $R > 0$, one
can write $\epsilon_k R_k \ge R$ for $k$ sufficiently large and thus
$$
\int_{\DD_R} v_\infty^*\omega = \lim_{k \to \infty} \int_{\DD_R} v_k^*\omega
\le \lim_{k \to \infty} \int_{\DD_{\epsilon_k R_k}} v_k^*\omega = 0,
$$
implying $\int_\CC v_\infty^*\omega = 0$.  It follows that $v_\infty$ must
be constant, so we have a contradiction.
\end{proof}

To obtain the uniform limit of $u(s,\cdot)$ as $s \to \infty$, we now 
pick any sequence of nonnegative numbers $s_k \to \infty$ and consider 
the sequence of $J$-holomorphic half-cylinders
$$
u_k : [-s_k,\infty) \times S^1 \to W : (s,t) \mapsto u(s + s_k,t).
$$
By Lemma~\ref{lemma:boundedDeriv}, these maps are uniformly
$C^1$-bounded, 
so elliptic regularity gives a subsequence converging
in $C^\infty_\loc$ on $\RR \times S^1$ to a $J$-holomorphic cylinder
$$
u_\infty : \RR \times S^1 \to W.
$$
Observe that for any $c > 0$, we can write
$-s_k/2 \le -c$ for sufficiently large $k$ and thus compute
\begin{equation*}
\begin{split}
\int_{[-c,c] \times S^1} u_\infty^*\omega &=
\lim_{k \to \infty} \int_{[-c,c,] \times S^1} u_k^*\omega \le
\lim_{k \to \infty} \int_{[-s_k/2,\infty) \times S^1} u_k^*\omega \\
&= \lim_{k \to \infty} \int_{[s_k/2,\infty) \times S^1} u^*\omega = 0
\end{split}
\end{equation*}
since $\int_{Z_+} u^*\omega < \infty$.  This implies $\int_{\RR \times S^1}
u_\infty^*\omega = 0$, so $u_\infty$ is a constant map to some
point $p \in W$, hence after replacing $s_k$ with a subsequence,
$$
u(s_k,\cdot) = u_k(0,\cdot) \to p \quad \text{ in } \quad
\text{$C^\infty(S^1,W)$ as $k \to \infty$}.
$$

To finish the proof of \eqref{eqn:unifConv}, we need to show that one 
cannot find two sequences
$s_k \to \infty$ and $s_k' \to \infty$ such that $u(s_k,\cdot) \to p$ and
$u(s_k',\cdot) \to p'$ for distinct points $p \ne p' \in W$.  This is an
easy consequence of the monotonicity lemma: indeed, if two such sequences
exist, then we can find a sequence $s_k'' \to \infty$ for which the loops
$u(s_k'',\cdot)$ alternate between arbitrarily small neighborhoods of
$p$ and $p'$.  Since $u$ is continuous, it must then pass through 
$\p B_{2r}(p)$ infinitely many times for $r > 0$ sufficiently small,
and in fact there exists an infinite sequence of pairwise disjoint 
neighborhoods $\uU_k \subset Z_+$ such that each
$$
u|_{\uU_k} : \uU_k \to B_r(q_k)
$$
is a proper map passing through some point $q_k \in \p B_{2r}(p)$.  
The monotonicity lemma then implies
$$
\int_{Z_+} u^*\omega \ge \sum_k \int_{\uU_k} u^*\omega \ge \sum_k c r^2 = \infty,
$$
a contradiction.

\begin{exercise}
\label{EX:essential}
Given an area form $\omega$ on $S^2 = \CC \cup \{\infty\}$ and a finite subset
$\Gamma \subset S^2$, show that a
holomorphic function $f : S^2 \setminus \Gamma \to \CC$ has an essential 
singularity at one of its punctures if and only if
$\int_\CC f^*\omega = \infty$.
\end{exercise}

\section{Finite energy and asymptotics}
\label{sec:asymptotics}

As further preparation for the compactness discussion, we now prove the
long-awaited converse of the fact that asymptotically cylindrical curves
have finite energy.  We work in the setting described in \S\ref{sec:stableBoundary}:
$(W,\omega)$ is a symplectic cobordism with stable boundary 
$\p W = -M_- \sqcup M_+$ carrying stable Hamiltonian structures
$\hH_\pm = (\omega_\pm,\lambda_\pm)$ with induced hyperplane distributions
$\xi_\pm = \ker \lambda_\pm$ and Reeb vector fields~$R_\pm$.  The completion
$(\widehat{W},\omega_h)$ carries the symplectic structure
$$
\omega_h := \begin{cases}
d\left( h(r) \lambda_+ \right) + \omega_+ & 
\text{ on $[0,\infty) \times M_+$} \\
\omega & \text{ on $W$},\\
d\left( h(r) \lambda_- \right) + \omega_- & 
\text{ on $(-\infty,0] \times M_-$},
\end{cases}
$$
for some $C^0$-small smooth function $h(r)$ with $h' > 0$ that is the 
identity near $r=0$, and for a fixed constant $r_0$,
we define a compact subset
$$
W^{r_0} := \left( [-r_0,0] \times M_- \right) \cup_{M_-} W 
\cup_{M_+} \left([0,r_0] \times M_+ \right) \subset \widehat{W},
$$
outside of which our $\omega_h$-tame almost complex structures
$J \in \jJ_\tau(\omega_h,r_0,\hH_+,\hH_-)$ are required to be 
translation-invariant and compatible with~$\hH_\pm$.  The
\defin{energy} of a $J$-holomorphic curve $u : (\dot{\Sigma},j) \to
(\widehat{W},J)$ is defined by
$$
E(u) := \sup_{f \in \tT(h,r_0)} \int_{\dot{\Sigma}} u^*\omega_f,
$$
where
$$
\tT(h,r_0) := \left\{ f \in C^\infty(\RR,(-\epsilon,\epsilon))\ \big|\ 
\text{$f' > 0$ and $f \equiv h$ near $[-r_0,r_0]$} \right\}.
$$
The constant $\epsilon > 0$ should always be assumed sufficiently small so that
if $J_\pm \in \jJ(\hH_\pm)$ and $X \in \xi_\pm$,
\begin{equation}
\label{eqn:epsilonSmall}
(\omega_\pm + \kappa \, d\lambda_\pm)(X,J_\pm X) > 0 \quad \text{ whenever } \quad
X \ne 0 \text{ and } \kappa \in (-2\epsilon,2\epsilon).
\end{equation}
This condition implies that every
$J \in \jJ_\tau(\omega_h,r_0,\hH_+,\hH_-)$ is tamed by every
$\omega_f$ for every $f \in \tT(h,r_0)$,
thus all $J$-holomorphic curves
satisfy $E(u) \ge 0$, with equality if and only if $u$ is constant.

\begin{thm}
\label{thm:asymptotics}
Assume all closed Reeb orbits in $(M_+,\hH_+)$ and $(M_-,\hH_-)$ are nondegenerate,
$J \in \jJ_\tau(\omega_h,r_0,\hH_+,\hH_-)$,
$(\Sigma,j)$ is a closed Riemann surface with $\dot{\Sigma} =
\Sigma \setminus \Gamma$ for some finite subset $\Gamma \subset \Sigma$,
and $u : (\dot{\Sigma},j) \to (\widehat{W},J)$ is a $J$-holomorphic curve such
that none of the singularities in $\Gamma$ are removable and
$E(u) < \infty$.  Then $u$ is asymptotically cylindrical.
\end{thm}

\begin{remark}
\label{remark:energySymp}
The theorem also holds in the setting
of a symplectization $(\RR \times M,J)$ with $J \in \jJ(\hH)$ for a
stable Hamiltonian structure $\hH = (\omega,\lambda)$ on~$M$.  The only real 
difference in this case is the slightly simpler definition of energy,
$$
E(u) = \sup_{f \in \tT} \int_{\dot{\Sigma}} u^*\omega_f,
$$
where $\omega_f := d\big( f(r) \lambda\big) + \omega$ and
$$
\tT = \left\{ f \in C^\infty(\RR,(-\epsilon,\epsilon))\ \big|\ 
f' > 0 \right\}.
$$
This change necessitates a few trivial modifications to the proof of
Theorem~\ref{thm:asymptotics} given below.
\end{remark}

Like removal of singularities, Theorem~\ref{thm:asymptotics} is really a 
local result, so let us formulate a more precise and more general
statement in these terms.  Let
$$
\dot{\DD} := \DD \setminus \{0\} \subset \CC
$$
and define the two biholomorphic maps
\begin{equation}
\label{eqn:halfCylinders}
\begin{split}
& \varphi_+ : Z_+ := [0,\infty) \times S^1 \to \dot{\DD} : (s,t) \mapsto e^{-2\pi(s+it)} \\
& \varphi_- : Z_- := (-\infty,0] \times S^1 \to \dot{\DD} : (s,t) \mapsto e^{2\pi(s+it)}.
\end{split}
\end{equation}

\begin{thm}
\label{thm:asympLocal}
Suppose $J \in \jJ_\tau(\omega_h,r_0,\hH_+,\hH_-)$ and $u : \dot{\DD} \to
\widehat{W}$ is a $J$-holomorphic map with $E(u) < \infty$.  Then either
the singularity at $0 \in \DD$ is removable or $u$ is a proper map.  In the
latter case the puncture is either positive or negative, meaning that $u$
maps neighborhoods of $0$ to neighborhoods of $\{\pm\infty\} \times M_\pm$,
and the puncture has a well-defined \defin{charge}, defined as
$$
Q = \lim_{\epsilon \to 0^+} \int_{\p \DD_\epsilon} u^*\lambda_\pm,
$$
which satisfies $\pm Q > 0$.  Moreover, the map
$$
(u_\RR(s,t),u_M(s,t)) := u \circ \varphi_\pm(s,t) \in \RR \times M_\pm
\quad \text{ for $(s,t) \in Z_\pm$ near infinity}
$$
satisfies 
$$
u_\RR(s,\cdot) - Ts \to c \quad \text{ in } \quad
\text{$C^\infty(S^1)$ as $s \to \pm\infty$}
$$
for $T := |Q|$ and a constant $c \in \RR$, while for every
sequence $s_k \to \pm\infty$, one can restrict to a subsequence such that
$$
u_M(s_k,\cdot) \to \gamma(T\cdot)
\quad \text{ in } \quad
\text{$C^\infty(S^1,M_\pm)$ as $k \to \infty$}
$$
for some $T$-periodic Reeb orbit $\gamma : \RR / T\ZZ \to M_\pm$.
If $\gamma$ is nondegenerate or Morse-Bott, then in fact
$$
u_M(s,\cdot) \to \gamma(T\cdot) \quad \text{ in } \quad
\text{$C^\infty(S^1,M_\pm)$ as $s \to \pm\infty$}
$$
\end{thm}

We will not prove this result in its full strength, as in particular the
last step (when $\gamma$ is nondegenerate or Morse-Bott) requires some
asymptotic elliptic regularity results that we do not have space to
explain here.  Note however that most of the above statement does not
require any nondegeneracy assumption at all.  The price for this level of
generality is that if
$s_k, s_k' \to \pm\infty$ are two distinct sequences, then
we have no guarantee in general that the two Reeb orbits obtained as
limits of subsequences of $u_M(s_k,\cdot)$ and $u_M(s_k',\cdot)$ will be the
same; at present, neither an example of this rather unpleasant possibility
nor any general argument to rule it out is known.  If one of these orbits
is assumed to be isolated, however---which is always true when
the Reeb vector field is nondegenerate---then we will be able
to show that both are the same up to parametrization, hence
\emph{geometrically}, $u_M(s,t)$ lies in arbitrarily
small neighborhoods of the orbit $\gamma$ as $s \to \pm\infty$.
This turns out to be also true in the more general Morse-Bott setting,
though it is then much harder to prove
since $\gamma$ need not be isolated.  
Once $u_M(s,\cdot)$ is localized near~$\gamma$, one can use the nondegeneracy
condition as we did in the Fredholm theory of Lecture~\ref{lec:Fredholm}
to develop asymptotic regularity results that give
much finer control over the behavior of $u_M$ as $s \to \pm\infty$,
implying in particular that $u_M(s,\cdot) \to \gamma(T\cdot)$ in
$C^\infty(S^1,M_\pm)$.  For details on this step,
we refer to the original sources: \cites{HWZ:props1,HWZ:FIMpreprint}
for the nondegenerate case, and \cites{HWZ:props4,Bourgeois:thesis} when
the Reeb vector field is Morse-Bott.  Those papers deal exclusively with
the contact case, but the setting of general stable Hamiltonian structures
is also dealt with in \cite{Siefring:asymptotics}.

Ignoring the final step for now,
the proof of Theorem~\ref{thm:asympLocal} will reuse most of the techniques 
that we already saw in our proof of
removal of singularities in \S\ref{sec:singularities}.  The main idea is to
use a combination of the monotonicity lemma and bubbling analysis to
show that unless $u$ has a removable singularity, it is a proper map,
and for any sequence $s_k \to \pm\infty$,
the holomorphic half-cylinders defined by
$$
u_k(s,t) = u \circ \varphi_\pm(s + s_k,t)
$$
on a sequence of increasingly large half-cylinders
must have a subsequence converging in $C^\infty_\loc(\RR \times S^1)$
to either a constant map or a trivial cylinder.  The first case will
turn out to mean (as in Theorem~\ref{thm:removable}) that the puncture
is removable, and the second implies asymptotic convergence to a closed
Reeb orbit.  

One major difference between the proof of Theorem~\ref{thm:asympLocal}
and removal of singularities is that since $\widehat{W}$ is noncompact,
sequences of curves in $\widehat{W}$ with uniformly bounded first derivatives 
need not be locally $C^0$-bounded.  This issue will arise both in the
bubbling argument to prove $|d u_k(s,t)| \le C$ and in the analysis of the
sequence $u_k$ itself.  In such cases, one can
use the $\RR$-translation action
\begin{equation}
\label{eqn:Raction}
\tau_c : \RR \times M_\pm \to \RR \times M_\pm : (r,x) \mapsto (r + c,x)
\quad \text{ for }\quad c \in \RR
\end{equation}
on suitable subsets of the cylindrical ends to replace unbounded sequences
with uniformly $C^1$-bounded sequences of curves mapping into $\RR \times M_+$ or 
$\RR \times M_-$.  These $\RR$-translations are the reason why our definition
of energy needs to be something slightly
more complicated than just the symplectic area $\int_{\dot{\Sigma}} u^*\Omega$
for a single choice of symplectic form.  To understand bubbling in the
presence of arbitrarily large $\RR$-translations, we will need the following
lemma.

\begin{lemma}
\label{lemma:noContactArea}
Suppose $J \in \jJ(\hH)$ for some stable Hamiltonian structure
$\hH = (\omega,\lambda)$ on an odd-dimensional manifold~$M$, 
and $u : (\dot{\Sigma},j) \to (\RR \times M,J)$
is a $J$-holomorphic curve satisfying
$$
E(u) < \infty \quad \text{ and }\quad \int_{\dot{\Sigma}} u^*\omega = 0.
$$
If $\dot{\Sigma} = \CC$, then $u$ is constant.  If $\dot{\Sigma} = \RR \times S^1$,
then $u$ either is constant or is biholomorphically equivalent to a
trivial cylinder over a closed Reeb orbit.
\end{lemma}
\begin{proof}
Denote $\xi = \ker \lambda$ and let
$$
\pi_\xi : T(\RR \times M) \to \xi
$$
denote the projection along the subbundle spanned by $\p_r$ (the unit vector field in the
$\RR$-direction) and the Reeb vector field~$R$.  Then since $\omega$ annihilates
both $\p_r$ and $R$, for any local 
holomorphic coordinates $(s,t)$ on a subset 
of~$\dot{\Sigma}$, the compatibility of $J|_\xi$ with $\omega|_\xi$ implies
$$
u^*\omega(\p_s,\p_t) = \omega(\p_s u,\p_t u) = \omega(\p_s u,J\p_s u) =
\omega( \pi_\xi \p_s u , J \pi_\xi \p_s u) \ge 0,
$$
hence $\int_{\dot{\Sigma}} u^*\omega \ge 0$ for every $J$-holomorphic curve, and
equality means that $u$ is everywhere tangent to the subbundle spanned by
$\p_r$ and~$R$.  This implies that $\im u$ is contained in the image of
some $J$-holomorphic plane of the form
$$
u_\gamma : \CC \to \RR \times M : s + it \mapsto (s,\gamma(t)),
$$
where $\gamma : \RR \to M$ is a (not necessarily periodic) orbit of~$R$.
If $\gamma$ is not periodic, then $u_\gamma$ is embedded, hence there exists a
unique (and necessarily holomorphic) map $\Phi : (\dot{\Sigma},j) \to (\CC,i)$
such that $u = u_\gamma \circ \Phi$.  If on the other hand $\gamma$ is
periodic with minimal period $T > 0$, then $u_\gamma$ descends to an
embedding of the cylinder
$$
\hat{u}_\gamma : \CC / iT \ZZ \to \RR \times M,
$$
and we can view $u_\gamma$ as a covering map to this embedded cylinder.
Now there exists a unique holomorphic map $\Phi : \dot{\Sigma} \to
\CC / iT \ZZ$ such that $u = \hat{u}_\gamma \circ \Phi$.  If
$\dot{\Sigma} = \CC$, then since $\pi_1(\CC) = 0$ implies that $\Phi$ can be lifted to
a (necessarily holomorphic) map $\widetilde{\Phi} : \CC \to \CC$ with
$u_\gamma \circ \widetilde{\Phi} = u$.  Relabeling symbols, we conclude that
in general if $\dot{\Sigma} = \CC$, then 
$u = u_\gamma \circ \Phi$ for a holomorphic map $\Phi : \CC \to \CC$.

Let us consider all cases in which the factorzation $u = u_\gamma \circ \Phi$
exists, where $\Phi : (\dot{\Sigma},j) \to (\CC,i)$ is holomorphic and
$\dot{\Sigma} = \Sigma \setminus \Gamma$ for a closed Riemann surface $(\Sigma,j)$.
We will now use the removable singularity theorem for $\Phi : \dot{\Sigma}
\to S^2 \setminus \{0\}$ to show that unless $\Phi$ 
is constant, $\int_{\dot{\Sigma}} u^*\omega_f = \infty$ for suitable choices
of $f \in \tT$.  This integral can be rewritten as
\begin{equation}
\label{eqn:integrals}
\int_{\dot{\Sigma}} u^*\omega_f = \int_{\dot{\Sigma}} \Phi^* u_\gamma^*\omega_f =
\int_{\dot{\Sigma}} \Phi^* d\left(f(s) \, dt \right) =
\int_{\dot{\Sigma}} \Phi^*\left( f'(s) \, ds \wedge dt \right)
\end{equation}
since $\omega_f = d\big(f(r) \, \lambda\big) + \omega$ and
$u_\gamma(s,t) = (s,\gamma(t))$.  Since $f' > 0$, $f'(s)\, ds \wedge dt$
is an area form on $\CC$ with infinite area.  We claim now that for suitable
choices of $f \in \tT$, one can find an area form $\Omega$ on $S^2
= \CC \cup \{\infty\}$ such that $\Omega \le f'(s)\, ds \wedge dt$.
To see this, let us change coordinates so that $\infty$ becomes~$0$: setting
$\Psi : \CC^* \to \CC^* : z \mapsto 1/z$, a slightly tedious but
straightforward computation gives
\begin{equation}
\label{eqn:areaForm}
\begin{split}
\Psi^*\left(f'(s)\, ds \wedge dt\right) &= f'(s / |z|^2)
\frac{1}{|z|^4} \left( 1 + \frac{(2 st)^2}{|z|^4} \right) \, ds \wedge dt \\
&\ge f'(s / |z|^2) \frac{1}{|z|^4} \, ds \wedge dt
\quad \text{ for } \quad z = s + it \in \CC \setminus \{0\}.
\end{split}
\end{equation}
We need to show that this $2$-form can be bounded away from~$0$ as
$z \to 0$.  Let us choose $f \in \tT$ such that
\begin{equation}
\label{eqn:fCondition}
f(r) = \pm \left(\epsilon - \frac{\epsilon}{2r}\right) \quad \text{ for } \quad
\pm r \ge 1
\end{equation}
and extend $f$ arbitrarily to $[-1,1]$ such that $f' > 0$.  We can then
find a constant $c > 0$ such that $f'$ satisfies
$$
f'(r) > \min\left\{ c , \frac{\epsilon}{2 r^2} \right\} \quad \text{ for all }
\quad r \in \RR.
$$
Plugging this into \eqref{eqn:areaForm} gives
$$
\Psi^*\left(f'(s)\, ds \wedge dt\right) \ge 
\min\left\{ \frac{c}{|z|^4} , \frac{\epsilon}{2 s^2} \right\} \, 
ds \wedge dt,
$$
which clearly blows up as $|z| \to 0$.  With this established, we
observe that for any number $C > 0$, the fact that $f'(s)\, ds \wedge dt$
has infinite area implies we can choose an area form $\Omega$ on $S^2$ with
$$
\Omega \le f'(s) \, ds \wedge dt \text{ on $S^2 \setminus \{\infty\}$} 
\qquad \text{ and } \quad
\int_{S^2} \Omega > C.
$$
We now have two possibilities:
\begin{enumerate}
\item If $\int_{\dot{\Sigma}} \Phi^*\Omega < \infty$, then 
Theorem~\ref{thm:removable} 
implies that the singularities of $\Phi : \dot{\Sigma} \to \CC$ at $\Gamma$
are all removable, i.e.~$\Phi$ extends to a
holomorphic map $(\Sigma,j) \to (S^2,i)$, which has
a well-defined mapping degree $k \ge 0$.  Then
$$
\int_{\dot{\Sigma}} u^*\omega_f = \int_{\dot{\Sigma}} 
\Phi^*\left( f'(s) \, ds \wedge dt \right)
\ge \int_{\dot{\Sigma}} \Phi^*\Omega = \int_{\Sigma} \Phi^*\Omega = 
k \int_{S^2} \Omega > k C.
$$
Since $C > 0$ can be chosen arbitrarily large, this implies
$\int_{\dot{\Sigma}} u^*\omega_f = \infty$ unless $k=0$, meaning 
$\Phi$ is constant.
\item If $\int_{\dot{\Sigma}} \Phi^*\Omega = \infty$ (meaning there is an
essential singularity, cf.~Exercise~\ref{EX:essential}), then since
$\Phi^*\left( f'(s) \, ds \wedge dt \right) \ge \Phi^*\Omega$, 
\eqref{eqn:integrals} implies $\int_\CC u^*\omega_f = \infty$.
\end{enumerate}
Since $u$ is constant whenever $\Phi$ is,
this completes the proof for $\dot{\Sigma} = \CC$.

If $\dot{\Sigma} = \RR \times S^1$, then it remains to deal with the case
where the factorization $u = u_\gamma \circ \Phi$ does not exist because
$\gamma$ is periodic.  If the minimal period is $T > 0$, then let us in this
case redefine $u_\gamma$ as an embedded $J$-holomorphic trivial cylinder
$$
u_\gamma : \RR \times S^1 : (s,t) \mapsto (Ts , \gamma(Tt)).
$$
Since the new $u_\gamma$ is embedded, we can now write 
$u = u_\gamma \circ \Phi$ for a unique holomorphic map
$\Phi : \RR \times S^1 \to \RR \times S^1$.  Identifying $\RR \times S^1$
biholomorphically with $S^2 \setminus \{0,\infty\}$, we claim that $\Phi$
extends to a holomorphic map $S^2 \to S^2$.  Indeed, by the removable 
singularity theorem, this is true if and only if $\int_{\RR \times S^1} \Phi^*\Omega 
< \infty$ for some area form $\Omega$ on $S^2$.  Notice that
$u_\gamma^*\omega_f = T^2 \cdot f'(Ts)\, ds \wedge dt$, defines an
area form on $\RR \times S^1$ with finite area for any $f \in \tT$
since $\int_{-\infty}^\infty
f'(s) \, ds < \infty$; this is equivalent to the observation that
trivial cylinders always have finite energy.  Using the biholomorphic map
$(s,t) \mapsto e^{2\pi (s+it)}$ to identify $\RR \times S^1$ with
$\CC^* = S^2 \setminus \{0,\infty\}$ and using coordinates $z = x+iy$ on
the latter, another tedious but straightforward computation gives
$$
u_\gamma^*\omega_f = \frac{T^2}{4\pi^2} 
\frac{f'\left(\frac{T}{2\pi} \log |z| \right)}{|z|^2} \, dx \wedge dy \quad
\text{ for } \quad z = x + i y \in \CC^*.
$$
Now suppose $f \in \tT$ is chosen as in
\eqref{eqn:fCondition}.  Then one can check
that the positive function in front of $dx \wedge dy$ in the above formula
goes to $+\infty$ as $|z| \to 0$; this means that one can find an area form
$\Omega$ on $\CC$ with $\Omega \le u_\gamma^*\omega_f$ on~$\CC^*$.
The singularity at $+\infty \in S^2$ can be handled in a similar way, thus
we can find an area form $\Omega$ on $S^2$ such that
$\Omega \le u_\gamma^*\omega_f$ on $\RR \times S^1$.
Now since $E(u) < \infty$, we have
$$
\int_{\RR \times S^1} \Phi^*\Omega \le \int_{\RR \times S^1} 
\Phi^*u_\gamma^*\omega_f =
\int_{\RR \times S^1} u^*\omega_f < \infty,
$$
so by Theorem~\ref{thm:removable}, $\Phi$ has a holomorphic extension
$S^2 \to S^2$, which is then a map of degree $k \ge 0$ with
$\Phi^{-1}(\{0,\infty\}) \subset \{0,\infty\}$.  If $k=0$ then
$\Phi$ is constant, and so is~$u$.  Otherwise, $\Phi$ is surjective and
thus hits both $0$ and $\infty$, but it can only do this at either $0$ or
$\infty$, thus it either fixes both or interchanges them.
After composing with a biholomorphic map of $S^2$
preserving $\RR \times S^1$, we may assume without loss of generality that
$\Phi(0) = 0$ and $\Phi(\infty) = \infty$.  This makes $\Phi$ a polynomial
with only one zero, hence as a map on $\CC \cup \{\infty\}$,
$\Phi(z) = c z^k$ for some $c \in \CC^*$.  Up to biholomorphic equivalence,
$\Phi(z)$ is then $z^k$, which appears in cylindrical coordinates as the
map $(s,t) \mapsto (ks,kt)$, so $u$ is now the trivial cylinder
$$
u(s,t) = u_\gamma(ks,kt) = (kT s , \gamma(kT t))
$$
over the $k$-fold cover of~$\gamma$.
\end{proof}
\begin{remark}
It may be useful for some applications to observe that 
Lemma~\ref{lemma:noContactArea} does not require $M$ to be compact.
In contrast, the compactness arguments in this lecture almost always
depend on the assumption that $W$ and $M_\pm$ are compact---without this,
one would need add some explicit assumption to guarantee local
$C^0$-bounds on sequences of holomorphic curves, e.g.~the assumption in
Theorem~\ref{thm:removable} that $u(\DD \setminus \{0\})$ is contained
in a compact subset.
\end{remark}

Before continuing, it is worth noting that neither of the two definitions of 
energy stated above
(one for curves in $\widehat{W}$ and the other for symplectizations)
is unique, i.e.~each can be tweaked in various ways such that the results of
this section still hold.  Indeed, the original definitions appearing in
\cites{Hofer:weinstein,SFTcompactness} are slightly different, but
equivalent to these.  The next lemma illustrates one further example of this
freedom, which will be useful in some of the arguments below.

\begin{lemma}
\label{lemma:Renergy}
Given a stable Hamiltonian structure $\hH = (\omega,\lambda)$ on $M$, a
sufficiently small constant $\epsilon > 0$ as in \eqref{eqn:epsilonSmall},
and $J \in \jJ(\hH)$, consider the alternative notion of energy for $J$-holomorphic
curves $u : (\dot{\Sigma},j) \to (\RR \times M,J)$ defined by
$$
E_0(u) = \sup_{f \in \tT_0} \int_{\dot{\Sigma}} u^*\omega_f
$$
where $\omega_f = d\left( f(r)\, \lambda \right) + \omega$ and
$$
\tT_0 = \left\{ f \in C^\infty(\RR,(a,b))\ \big|\ f' > 0 \right\}
$$
for some constants $-\epsilon \le a < b \le \epsilon$.  Then if $E(u)$
denotes the energy as written in Remark~\ref{remark:energySymp}, there exists
a constant $c > 0$, depending on the data $a$, $b$, $\epsilon$ and $\hH$
but not on~$u$, such that
$$
c E(u) \le E_0(u) \le E(u).
$$
\end{lemma}
\begin{proof}
The second of the two inequalities is immediate since $\tT_0 \subset \tT$.
For the first inequality, note that since $\epsilon > 0$ is small, we can
assume there exists a constant $c > 1$ such that for every  
$X \in T(\RR \times M)$ and every $\kappa \in [-\epsilon,\epsilon]$,
\begin{equation}
\label{eqn:omegaEstimate}
\frac{1}{c} (\omega + \kappa \, d\lambda)(X,JX) \le \omega(X,JX) \le
c (\omega + \kappa \, d\lambda)(X,JX).
\end{equation}
This uses \eqref{eqn:epsilonSmall} and the fact that $d\lambda$ annihilates
$\ker \omega$.
Now suppose $f \in \tT$, choose a constant $\delta \in (0,b-a]$
and define $\tilde{f} \in \tT_0$ by
$$
\tilde{f}(r) = \frac{\delta}{2 \epsilon} f(r) + \frac{a+b}{2}.
$$
Then $\tilde{f}'(r) = \frac{\delta}{2\epsilon} f'(r)$, and
given a $J$-holomorphic curve $u : \dot{\Sigma} \to \RR \times M$, we
can write $\omega_f = \omega + f(r)\, d\lambda +
f'(r) \, dr \wedge \lambda$ and use \eqref{eqn:omegaEstimate} to
estimate
\begin{equation*}
\begin{split}
\int_{\dot{\Sigma}} u^*\omega_f &=
\int_{\dot{\Sigma}} u^*\left( \omega + f(r)\, d\lambda \right) +
\int_{\dot{\Sigma}} u^* \left( f'(r) \, dr \wedge \lambda \right) \\
&\le c \int_{\dot{\Sigma}} u^*\omega + 
\frac{2\epsilon}{\delta} \int_{\dot{\Sigma}} u^*\left( 
\tilde{f}'(r) \, dr \wedge \lambda \right) \\
&\le c^2 \int_{\dot{\Sigma}} u^*\left( \omega + \tilde{f}(r) \, d\lambda \right)
+ \frac{2\epsilon}{\delta} \int_{\dot{\Sigma}} u^*\left(
\tilde{f}'(r) \, dr \wedge \lambda \right) .
\end{split}
\end{equation*}
If $c^2 \ge \frac{2\epsilon}{b-a}$, then we can choose 
$\delta := 2\epsilon / c^2 \le b-a$ and rewrite the last expression as
\begin{equation*}
\begin{split}
c^2 \int_{\dot{\Sigma}} & u^*\left( \omega + \tilde{f}(r) \, d\lambda \right)
+ \frac{2\epsilon}{\delta} \int_{\dot{\Sigma}} u^*\left(
\tilde{f}'(r) \, dr \wedge \lambda \right) \\
&= c^2 \int_{\dot{\Sigma}}
u^*\left( \omega + \tilde{f}(r)\, d\lambda + \tilde{f}'(r)\, dr \wedge \lambda \right)
= c^2 \int_{\dot{\Sigma}} u^*\omega_{\tilde{f}} \le c^2 E_0(u).
\end{split}
\end{equation*}
On the other hand if $c^2 < \frac{2\epsilon}{b-a}$, we can set
$\delta := b-a$ and write
\begin{equation*}
\begin{split}
c^2 \int_{\dot{\Sigma}} u^*\left( \omega + \tilde{f}(r) \, d\lambda \right)
&+ \frac{2\epsilon}{\delta} \int_{\dot{\Sigma}} u^*\left(
\tilde{f}'(r) \, dr \wedge \lambda \right) \\
&\le \frac{2\epsilon}{b-a} \int_{\dot{\Sigma}}
u^*\left( \omega + \tilde{f}(r)\, d\lambda + \tilde{f}'(r)\, dr \wedge \lambda \right) \\
&= \frac{2\epsilon}{b-a} \int_{\dot{\Sigma}} u^*\omega_{\tilde{f}}
\le \frac{2\epsilon}{b-a} E_0(u).
\end{split}
\end{equation*}
\end{proof}

With this preparation out of the way, we now begin in earnest with the proof
of Theorem~\ref{thm:asympLocal}.  Assume 
$u : \dot{\DD} \to \widehat{W}$ is a $J$-holomorphic punctured disk satisfying
$E(u) < \infty$.  Using the maps $\varphi_\pm : Z_\pm \to \dot{\DD}$
defined in \eqref{eqn:halfCylinders}, we shall write
$$
u_\pm := u \circ \varphi_\pm : Z_\pm \to \widehat{W}
$$
and observe that these reparametrizations have no impact on the energy, i.e.
$$
E(u_\pm) = \sup_{f \in \tT(h,r_0)} \int_{Z_\pm} (u \circ \varphi_\pm)^*\omega_f
= \sup_{f \in \tT(h,r_0)} \int_{\dot{\DD}} u^*\omega_f = E(u).
$$
Fix a Riemannian metric on $\widehat{W}$ that is translation-invariant on the
cylindrical ends, and fix the standard metric on the half-cylinders~$Z_\pm$.  
We will use these metrics implicitly whenever referring to quantities 
such as~$|du_\pm(z)|$.

\begin{lemma}
\label{lemma:boundedDeriv2}
There exists a constant $C > 0$ such that $|du_+(s,t)| \le C$ 
for all $(s,t) \in Z_+$.
\end{lemma}
\begin{proof}
We use a bubbling argument as in the proof of Lemma~\ref{lemma:boundedDeriv}.
Suppose the contrary, so there exists a sequence $z_k = (s_k,t_k) \in Z_+$
with $R_k := |du_+(z_k)| \to \infty$.  Choose a sequence 
$\epsilon_k > 0$ with $\epsilon_k \to 0$ but $\epsilon_k R_k \to \infty$, and 
using Lemma~\ref{lemma:Hofer}, assume without loss of generality that
$$
|du_+(z)| \le 2 R_k \quad \text{ for all } \quad 
z \in \DD_{\epsilon_k}(z_k).
$$
Define a rescaled sequence of $J$-holomorphic disks by
$$
v_k : \DD_{\epsilon_k R_k} \to \widehat{W} : z \mapsto 
u \circ \varphi_+(z_k + z / R_k).
$$
These satisfy $|d v_k| \le 2$ on their domains, but they are not
necessarily $C^1$-bounded since their images may escape to infinity.
We distinguish three possibilities, at least one of which must hold:

\textsl{Case~1: $v_k(0)$ has a bounded subsequence.} \\
Then the corresponding subsequence of $v_k : \DD_{\epsilon_k R_k}
\to \widehat{W}$ is uniformly $C^1$-bounded on every compact subset and thus
(by elliptic regularity) has a further subsequence convergent in
$C^\infty_\loc(\CC)$ to a $J$-holomorphic plane
$$
v_\infty : \CC \to \widehat{W}
$$
with $|d v_\infty(0)| = \lim_{k \to \infty} |d v_k(0)| = 1$.  But by the
same argument we used in the proof of Lemma~\ref{lemma:boundedDeriv},
the fact that $\int_{Z_+} u_+^*\omega_f < \infty$ for any choice of
$f \in \tT(h,r_0)$ implies
$$
\int_{\CC} v_\infty^*\omega_f = 0,
$$
hence $v_\infty$ is constant, and this is a contradiction.

\textsl{Case~2: $v_k(0)$ has a subsequence diverging to $\{+\infty\} \times M_+$}. \\
Restricting to this subsequence, suppose
$$
v_k(0) \in \{r_k\} \times M_+,
$$
so $r_k \to \infty$, and assume without loss of generality that
$r_k > r_0$ for all~$k$.  Let $\tilde{R}_k \in (0,\epsilon_k R_k]$ for each $k$
denote the largest radius such that $v_k(\DD_{\tilde{R}_k}) \subset
(r_0,\infty) \times M_+$.  Then $\tilde{R}_k \to \infty$ since $|dv_k|$
is bounded.  Now using the $\RR$-translation maps $\tau_r$ defined in
\eqref{eqn:Raction}, define
$$
\tilde{v}_k := \tau_{-r_k} \circ v_k|_{\DD_{\tilde{R}_k}} : \DD_{\tilde{R}_k}
\to \RR \times M_+.
$$
Since we're using a translation-invariant metric on $[r_0,\infty) \times M_+$,
$\tilde{v}_k$ is now a uniformly $C^1_\loc$-bounded sequence of maps
into $\RR \times M_+$.  Elliptic regularity thus provides a subsequence 
convergent in $C^\infty_\loc(\CC)$ to a plane
$$
v_\infty : \CC \to \RR \times M_+,
$$
which is $J_+$-holomorphic, where $J_+ \in \jJ(\hH_+)$ denotes the restriction
of $J$ to $[r_0,\infty) \times M_+$, extended over $\RR \times M_+$ by
$\RR$-invariance.  We claim,
\begin{equation}
\label{eqn:energyClaim}
E(v_\infty) < \infty \quad \text{ and } \quad \int_\CC v_\infty^*\omega_+ = 0,
\end{equation}
where $E(v_\infty)$ is now defined as in Remark~\ref{remark:energySymp}.  
By Lemma~\ref{lemma:Renergy}, the first part of the claim will follow if we
can fix a constant $a \in (-\epsilon,\epsilon)$ and 
establish a uniform bound
$$
\int_\CC v_\infty^*\Omega^+_f \le C,
$$
with $\Omega^+_f := \omega_+ + d\big(f(r)\, \lambda_+\big)$,
for all smooth and strictly increasing functions $f : \RR \to
(a,\epsilon)$.  For convenience in the following, we shall assume
$a > h(r_0)$.  Now if $f$ is such a function, then for any
$R > 0$,
$$
\int_{\DD_R} v_\infty^*\Omega^+_f = \lim_{k\to \infty}
\int_{\DD_R} v_k^*\tau_{-r_k}^*\Omega^+_f =
\lim_{k \to \infty} \int_{\DD_R} v_k^*\Omega^+_{f_k},
$$
where $f_k(r) := f(r - r_k)$.  Notice that the dependence of
the last integral on $f_k$ is limited to the interval $(r_0,\infty)$
since $v_k(\DD_R) \subset (r_0,\infty) \times M_+$.  Then
since $f > a > h(r_0)$
by assumption, there exists for each $k$ a function $h_k \in \tT(h,r_0)$
that matches $f_k$ outside some neighborhood of $(-\infty,r_0]$ and
thus satisfies
$$
\int_{\DD_R} v_k^*\Omega^+_{f_k} = \int_{\DD_R} v_k^*\omega_{h_k}
\le \int_{\DD_{\epsilon_k R_k}} v_k^*\omega_{h_k} =
\int_{\DD_{\epsilon_k}(z_k)} u_+^*\omega_{h_k} \le
\int_{Z_+} u_+^*\omega_{h_k} \le E(u).
$$
This is true for every $R > 0$ and thus proves the first part of
\eqref{eqn:energyClaim}.  To establish the second part, fix $R > 0$ again
and pick any $f \in \tT(h,r_0)$.  Observe that since we can assume 
(after perhaps passing to a subsequence) the disks $\DD_{\epsilon_k}(z_k)$ 
are all disjoint,
\begin{equation*}
\begin{split}
0 &= \lim_{k \to \infty} \int_{\DD_{\epsilon_k}(z_k)} u_+^*\omega_f =
\lim_{k \to \infty} \int_{\DD_{\epsilon_k R_k}} v_k^*\omega_f =
\lim_{k \to \infty} \int_{\DD_{\epsilon_k R_k}} \tilde{v}_k^* \tau_{r_k}^* \omega_f \\
&\ge \lim_{k \to \infty} \int_{\DD_R} \tilde{v}_k^*\tau_{r_k}^*\omega_f 
= \lim_{k \to \infty} \int_{\DD_R} \tilde{v}_k^*\Omega^+_{f_k},
\end{split}
\end{equation*}
where now $f_k(r) := f(r + r_k)$.
Writing $\Omega^+_{f_k} = \omega_+ + d\big( f_k(r) \, \lambda_+\big)
= \omega_+ + f_k(r)\, d\lambda_+ + f_k'(r)\, dr \wedge \lambda_+$,
we can choose $f$ such that $f'(r) = f'(r + r_k) \to 0$
as $k \to \infty$, so the third term contributes nothing to the integral.
For the second term, let $f_+ := \lim_{k \to \infty} f_k(r) =
\lim_{r \to \infty} f(r)$, so the calculation above becomes
$$
0 \ge \int_{\DD_R} v_\infty^*\left( \omega_+ + f_+ \, d\lambda_+ \right).
$$
Now observe that since $f_+ \in [-\epsilon,\epsilon]$, 
condition~\eqref{eqn:epsilonSmall}
implies that the $2$-form $\omega_+ + f_+\, d\lambda_+$ is nondegenerate on
$\xi_+$, and it also annihilates $\p_r$ and~$R_+$, so the vanishing of this 
integral implies that $v_\infty$ is everywhere tangent to $\p_r$ and $R_+$
over~$\DD_R$.  But $R > 0$ was arbitrary, so this is true on the whole
plane, which is equivalent to $\int_\CC v_\infty^*\omega_+ = 0$.
With the claim established, we apply Lemma~\ref{lemma:noContactArea} and
conclude that $v_\infty$ is constant, contradicting the fact that
$|dv_\infty(0)| = 1$.

\textsl{Case~3: $v_k(0)$ has a subsequence diverging to $\{-\infty\} \times M_-$}. \\
This is simply the mirror image of case~2: writing the restriction of
$J$ to $(-\infty,-r_0] \times M_-$ as $J_-$, one can follow the same
bubbling argument but translate up and instead of down, giving rise to a
limiting nonconstant $J_-$-holomorphic plane 
$v_\infty : \CC \to \RR \times M_-$ that has finite energy but
$\int_\CC v_\infty^*\omega_- = 0$, in contradiction to
Lemma~\ref{lemma:noContactArea}.
\end{proof}

Consider now a sequence $s_k \to \infty$ and construct the $J$-holomorphic
half-cylinders
$$
u_k : [-s_k,\infty) \times S^1 \to \widehat{W} : (s,t) \mapsto 
u_+(s + s_k,t).
$$
The derivatives $|du_k|$ are uniformly bounded due to 
Lemma~\ref{lemma:boundedDeriv2}, though again, $u_k$ might fail to be
uniformly bounded in~$C^0$.  We distinguish three cases.

\textsl{Case~1: $u_k(0,0)$ has a bounded subsequence.}\\
Then the corresponding subsequence of $u_k$ is uniformly $C^1$-bounded
on compact subsets and
thus has a further subsequence converging in $C^\infty_\loc(\RR \times S^1)$
to a $J$-holomorphic cylinder
$$
u_\infty : \RR \times S^1 \to \widehat{W}.
$$
For any $f \in \tT(h,r_0)$ and any $c > 0$, we have
\begin{equation}
\label{eqn:cylinderZero}
\begin{split}
\int_{[-c,c] \times S^1} u_\infty^*\omega_f &=
\lim_{k \to \infty} \int_{[-c,c] \times S^1} u_k^*\omega_f \le
\lim_{k \to \infty} \int_{[-s_k/2,\infty) \times S^1} u_k^*\omega_f \\
&= \lim_{k \to \infty} \int_{[s_k/2,\infty) \times S^1} u_+^*\omega_f = 0
\end{split}
\end{equation}
since $\int_{Z_+} u_+^*\omega_f < \infty$.  It follows that
$\int_{\RR \times S^1} u_\infty^*\omega_f = 0$, so $u_\infty$ is a
constant map to some point $p \in \widehat{W}$, implying that after passing
to a subsequence of~$s_k$,
$$
u_+(s_k,\cdot) \to p \quad \text{ in $C^\infty(S^1,\widehat{W})$ } \quad 
\text{ as $k \to \infty$}.
$$

\textsl{Case~2: $u_k(0,0)$ has a subsequence diverging to $\{+\infty\} \times M_+$.}\\
Passing to the corresponding subsequence of $u_k$, suppose
$$
u_k(0,0) \in \{r_k\} \times M_+,
$$
so $r_k \to \infty$.  Since the derivatives 
$|d u_k|$ are uniformly bounded, we can then find a sequence of intervals 
$[-R_k^-,R_k^+] \subset [-s_k,\infty)$ such that
$$
u_k([-R_k^-,R_k^+] \times S^1) \subset [r_0,\infty) \times M_+ \quad \text{ and }\quad
R_k^\pm \to \infty.
$$
Now the translated sequence 
$$
\tau_{-r_k} \circ u_k|_{[-R_k^-,R_k^+] \times S^1} :
[-R_k^-,R_k^+] \times S^1 \to \RR \times M_+
$$
is uniformly $C^1$-bounded on compact subsets and thus has a subsequence
coverging in $C^\infty_\loc$ to a $J_+$-holomorphic cylinder
$$
u_\infty : \RR \times S^1 \to \RR \times M_+,
$$
where $J_+$ again denotes the restriction of $J$ to $[r_0,\infty) \times M_+$,
extended over $\RR \times M_+$ by $\RR$-translation.  We claim that this
cylinder satisfies
$$
E(u_\infty) < \infty \quad \text{ and } \quad 
\int_{\RR \times S^1} u_\infty^*\omega_+ = 0.
$$
The proof of this should be an easy exercise if you understood the proofs of
\eqref{eqn:energyClaim} and \eqref{eqn:cylinderZero} above, so I will leave
it as such.  Lemma~\ref{lemma:noContactArea} now implies that $u_\infty$
is either constant or is a reparametrization of a trivial cylinder
$$
u_\gamma : \RR \times S^1 \to \RR \times M_+ : (s,t) \mapsto (Ts,\gamma(Tt))
$$
for some Reeb orbit $\gamma : \RR / T\ZZ \to M_+$ with period $T > 0$.
More precisely, all the biholomorphic reparametrizations of $\RR \times S^1$
are of the form $(s,t) \mapsto (\pm s + a,\pm t + b)$, thus after
shifting the parametrization of $\gamma$, we can write $u_\infty$
without loss of generality in the form
\begin{equation}
\label{eqn:paramTrivial}
u_\infty(s,t) = \left( \pm T s + a , \gamma(\pm T t) \right)
\end{equation}
for some constant $a \in \RR$ and a choice of signs to be determined below
(see Lemma~\ref{lemma:chargePosNeg}).

\textsl{Case~3: $u_k(0,0)$ has a subsequence diverging to $\{-\infty\} \times M_-$.}\\
Writing $J_- := J|_{(-\infty,-r_0] \times M_-} \in \jJ(\hH_-)$ and
imitating the argument for case~2, we suppose $u_k(0,0) \in \{-r_k\} \times M_-$
with $r_k \to \infty$ and obtain a subsequence for which
$\tau_{r_k} \circ u_k$ converges in $C^\infty_\loc(\RR \times S^1)$ to a
$J_-$-holomorphic cylinder $u_\infty : \RR \times S^1 \to \RR \times M_-$,
where $u_\infty$ is either a constant or takes the form
\eqref{eqn:paramTrivial} for some orbit Reeb $\gamma : \RR / T\ZZ \to M_-$
of period $T > 0$.

Here is one easy consequence of the discussion so far.  Use the Riemannian
metric on $\widehat{W}$ to define a metric $\dist_{C^0}(\cdot,\cdot)$ on
the space of continuous loops $S^1 \to \widehat{W}$.

\begin{lemma}
\label{lemma:dist}
Given $\delta > 0$, there exists $s_0 \ge 0$ such that for every $s \ge s_0$,
the loop $u_+(s,\cdot) : S^1 \to \widehat{W}$ satisfies
$$
\dist_{C^0}(u_+(s,\cdot) , \ell) < \delta,
$$
where $\ell : S^1 \to \widehat{W}$ either is constant or is a loop of the
form $\ell(t) = (r,\gamma(\pm Tt))$ in $[r_0,\infty) \times M_+$ or
$(-\infty,r_0] \times M_-$ for some constant $r \in \RR$ and
Reeb orbit $\gamma : \RR / T\ZZ \to M_\pm$ of period $T > 0$.
\end{lemma}
\begin{proof}
If not, then there exists a sequence $s_k \to \infty$ such that each of the
loops $u_+(s_k,\cdot)$ lies at $C^0$-distance at least $\delta$ 
away from any loop of the above form.  However, the preceding discussion then 
gives a subsequence for which $u(s_k,\cdot)$ becomes arbitrarily 
$C^\infty$-close to such a loop, so this is a contradiction.
\end{proof}

\begin{lemma}
\label{lemma:proper}
If $u : \dot{\DD} \to \widehat{W}$ is not bounded, then it is proper.
\end{lemma}
\begin{proof}
We use the monotonicity lemma.  Suppose there exists a sequence
$(s_k,t_k) \in Z_+$ such that $u_+(s_k,t_k)$ diverges to
$\{+\infty\} \times M_+$.  This implies $s_k \to \infty$, and we claim then
that for every $R \ge r_0$, there exists $s_0 \ge 0$ such that
$$
u_+((s_0,\infty) \times S^1) \subset (R,\infty) \times M_+.
$$
If not, then we find $R \ge r_0$ and a sequence $(s_k',t_k') \in Z_+$ with 
$s_k' \to \infty$ such that 
$u_+(s_k',t_k') \not\in (R,\infty) \times M_+$
for every~$k$.  By continuity, we are free to suppose 
$u_+(s_k',t_k') \in \{R\} \times M_+$ for all~$k$ since 
Lemma~\ref{lemma:dist}
implies $u_+(\{s_k\} \times S^1) \subset (2R,\infty) \times M_+$
for $k$ sufficiently large.
Using Lemma~\ref{lemma:dist} again, we also have
$$
u_+(\{s_k'\} \times S^1) \subset (R-1,R+1) \times M_+
$$
for all $k$ large.  Assuming $2R > R + 2$ without loss of generality,
we can therefore find infinitely many pairwise disjoint
annuli of the form $[s_k',s_j] \times S^1 \subset Z_+$ containing open
sets that $u$ maps properly to small balls centered at points
in $\{R + 2\} \times M_+$.  Choosing any $f \in \tT(h,r_0)$, the monotonicity
lemma implies that each of these contributes at least some fixed amount
to $\int_{Z_+} u_+^*\omega_f$, contradicting the
assumption that $E(u) < \infty$.\footnote{The fact that $\widehat{W}$ is
noncompact is not a problem for this application of the monotonicity lemma,
as we are only using it in the compact subset $W^{2R} \subset \widehat{W}$.}

A similar argument works if $u_+(s_k,t_k)$ diverges to
$\{-\infty\} \times M_-$, proving that for every $R \ge r_0$, there exists
$s_0 \ge 0$ with
$$
u_+((s_0,\infty) \times S^1) \subset (-\infty,-R) \times M_-.
$$
\end{proof}

If $u$ is bounded, then the singularity at $0$ is removable by
Theorem~\ref{thm:removable}.  If not, then Lemma~\ref{lemma:proper} implies
that it maps neighborhoods of the puncture to neighborhoods of either
$\{+\infty\} \times M_+$ or $\{-\infty\} \times M_-$, and we shall refer to
the puncture as \emph{positive} or \emph{negative} accordingly.  

\begin{lemma}
\label{lemma:charge}
If the puncture is positive/negative, then the limit
$$
Q := \lim_{s \to \infty} \int_{S^1} u_+(s,\cdot)^*\lambda_\pm \in \RR
$$
exists.
\end{lemma}
\begin{proof}
If the puncture is positive, fix $s_0 \ge 0$ such that
$u_+([s_0,\infty) \times S^1) \subset [r_0,\infty) \times M_+$.
Then by Stokes' theorem, it suffices to show that the integral
$\int_{[s_0,\infty) \times S^1} u_+^*d\lambda_+$ exists,
which is true if
\begin{equation}
\label{eqn:L1estimate}
\int_{[s_0,\infty) \times S^1} \left| u_+^*d\lambda_+ \right| < \infty.
\end{equation}
We claim first that $\int_{[s_0,\infty) \times S^1} u_+^*\omega_+ < \infty$.
Indeed, for any $s > s_0$ and $f \in \tT(h,r_0)$, we have
$$
E(u) \ge \int_{[s_0,s] \times S^1} u_+^*\omega_f =
\int_{[s_0,s] \times S^1} u_+^*\omega_+ +
\int_{[s_0,s] \times S^1} u_+^*d\left( f(r) \, \lambda_+ \right).
$$
Applying Stokes' theorem, the second term becomes the sum of some number
not dependent on~$s$ and the integral
$$
\int_{S^1} u_+(s,\cdot)^*\left( f(r)\, \lambda_+\right) =
\int_{S^1} [f \circ u_+(s,\cdot)] \, u_+(s,\cdot)^*\lambda_+,
$$
which is bounded as $s \to \infty$ since $f$ and $|du_+|$ are both bounded.
This proves that $\int_{[s_0,s] \times S^1} u_+^*\omega_+$
is also bounded as $s \to \infty$, and since $u_+^*\omega_+ \ge 0$,
the claim follows.  Now observe that since $d\lambda_+$ annihilates the
kernel of $\omega_+$ and the latter tames $J$ on~$\xi_+$, there exists a
constant $c > 0$ such that $| u_+^*d\lambda_+ | \le
c | u_+^*\omega_+ |$, implying \eqref{eqn:L1estimate}.

An analogous argument works if the puncture is negative.
\end{proof}

The number $Q \in \RR$ defined in the above lemma matches what we referred
to in the statement of Theorem~\ref{thm:asympLocal} as the \defin{charge}
of the puncture.

\begin{lemma}
\label{lemma:chargePosNeg}
If the puncture is nonremovable and $Q \ne 0$, then the puncture is 
positive/negative if and only if $Q > 0$ or $Q < 0$ respectively.  In either
case, given any 
sequence $s_k \to \infty$ with
$u_+(s_k,0) \in \{\pm r_k\} \times M_\pm$, 
one can find a sequence $R_k \in [0,s_k]$ with $R_k \to \infty$ such that
$u_+$ maps $[s_k - R_k,\infty) \times S^1$ into the positive/negative
cylindrical end for every $k$, and the sequence of half-cylinders
$$
u_k : [-R_k,\infty) \times S^1 \to \RR \times M_+ \quad \text{ or }\quad
u_k : (-\infty,R_k] \times S^1 \to \RR \times M_-
$$
defined by $u_k(s,t) = \tau_{\mp r_k} \circ u_\pm(s \pm s_k,t)$ has a
subsequence convergent in $C^\infty_\loc(\RR \times S^1)$ to a 
$J_\pm$-holomorphic cylinder of the form
$$
u_\infty : \RR \times S^1 \to \RR \times M_\pm : (s,t) \mapsto
(Ts + a, \gamma(Tt))
$$
for some constant $a \in \RR$ and Reeb orbit $\gamma : \RR / T\ZZ \to M_\pm$
with period $T := \pm Q$.
\end{lemma}
\begin{proof}
Assume the puncture is either positive or negative and $Q \ne 0$.
In the discussion preceding Lemma~\ref{lemma:dist},
we showed that the sequence 
$u'(s,t) := \tau_{\mp r_k} \circ u_+(s + s_k,t)$ defined on 
$[-R_k,\infty) \times S^1$ has a subsequence convergent in $C^\infty_\loc$
to a $J_\pm$-holomorphic cylinder $u_\infty' : \RR \times S^1 \to
\RR \times M_\pm$ which is either constant or of the form
\begin{equation}
\label{eqn:aTrivCyl}
u_\infty'(s,t) = (\sigma T s + a, \gamma(\sigma Tt))
\end{equation}
for some $a \in \RR$, $\sigma = \pm 1$ and a Reeb orbit $\gamma : \RR / T\ZZ
\to M_\pm$ of period $T > 0$.  We then have
$$
0 \ne Q = \lim_{s \to \infty} \int_{S^1} u_+(s,\cdot)^*\lambda_\pm = 
\lim_{k \to \infty} \int_{S^1} u_k'(0,\cdot)^*\lambda_\pm =
\int_{S^1} u_\infty'(0,\cdot)^*\lambda_\pm,
$$
so $u_\infty'$ cannot be constant, and from \eqref{eqn:aTrivCyl} we deduce
$Q = \sigma T$, hence $u_\infty'(s,t) = (Qs + a, \gamma(Qt))$.
Writing $u_+(s,t) = (u_\RR(s,t),u_M(s,t)) \in \RR \times M_\pm$ for
$s$ sufficiently large, it follows that every sequence $s_k \to \infty$
admits a subsequence for which
$$
\p_s u_\RR(s_k,\cdot) \to Q \quad \text{ in } \quad
C^\infty(S^1,\RR),
$$
and consequently $\p_s u_\RR(s,\cdot) \to Q$ in $C^\infty(S^1,\RR)$ as
$s \to \infty$.  This proves that the sign of $Q$ matches the sign of
the puncture whenever $Q \ne 0$.  The stated formula for $u_\infty$
now follows by adjusting all the appropriate signs in the case $Q < 0$.
\end{proof}

\begin{lemma}
\label{lemma:Q0}
If the puncture is nonremovable, then $Q \ne 0$.
\end{lemma}
\begin{proof}
Assume on the contrary that $u$ is a proper map, say with a positive puncture,
but $Q = 0$.  In this case, the argument of the previous lemma shows that
the limiting map $u_\infty : \RR \times S^1 \to \RR \times M_+$ will always be
\emph{constant}, thus for every sequence $s_k \to \infty$, there exists
a point $p \in M_+$ such that $u_+(s_k,0) \in \{r_k\} \times M_+$ with
$r_k \to \infty$ and
$$
\tau_{-r_k} \circ u_+(s_k,\cdot) \to (0,p) \in \RR \times M_+ \quad 
\text{ in } \quad \text{$C^\infty(S^1,\RR \times M_+)$ as $k \to \infty$}.
$$
In particular, this implies that all derivatives of $u_+$ decay to $0$ as
$s \to \infty$.  Intuitively, this should suggest to you that portions
of $u_+$ near infinity will have improbably small symplectic area, perhaps
violating the monotonicity lemma---this will turn out to be true, but we have 
to be a bit clever with our argument since $u_+$ is unbounded.  We will make
this argument precise by translating pieces of $u_+$ downward so that we
only compute its symplectic area in $[0,2] \times M_+$.  Fix
a function $f : \RR \to (-\epsilon,\epsilon)$ with $f' > 0$ and set
$\Omega^+_f = \omega_+ + d\left( f(r)\, \lambda_+ \right)$.

Given a small number $\delta > 0$, we can find $s_0 \ge 0$ such that
$|du_+(s,t)| < \delta$ for all $s \ge s_0$ and each of the loops
$u_+(s,\cdot)$ for $s \ge s_0$ is $\delta$-close to a constant in~$C^1(S^1)$.  
Assume $u_+(s_0,0) \in \{R\} \times M_+$ and choose $s_1 > s_0$ such that
$u_+(s_1,0) \in \{R + 2\} \times M_+$, which is possible since 
$u_+(s,t) \to \{+\infty\} \times M_+$ as $s \to \infty$.  Now consider the
$J_+$-holomorphic annulus
$$
v_\delta := \tau_{-R} \circ u_+|_{[s_0,s_1] \times S^1} : [s_0,s_1] \times S^1 \to
\RR \times M_+.
$$
We claim that $\int_{[s_0,s_1] \times S^1} v_\delta^*\Omega^+_f$ can be made
arbitrarily small by choosing $\delta$ suitably small.  Indeed, we can
use Stokes' theorem to write this integral as
\begin{equation*}
\begin{split}
\int_{[s_0,s_1] \times S^1} v_\delta^*\Omega^+_f &=
\int_{[s_0,s_1] \times S^1} v_\delta^*\omega_+ + 
\int_{[s_0,s_1] \times S^1} v_\delta^*d\left( f(r)\, \lambda_+ \right) \\
&= \int_{[s_0,s_1] \times S^1} v_\delta^*\omega_+ +
\int_{S^1} \left[ v_\delta(s_1,\cdot)^*\left( f(r)\, \lambda_+ \right) -
v_\delta(s_0,\cdot)^*\left( f(r) \, \lambda_+ \right) \right].
\end{split}
\end{equation*}
The second term is small because $f(r)$ is bounded and $|v_\delta(s,\cdot)^*\lambda_+|$
is small in proportion to $|dv_\delta(s,t)| = |d u_+(s,t)|$ for $s \ge s_0$.
For the first term, observe that since both of the loops $v_\delta(s_i,\cdot)$ for
$i=0,1$ are nearly constant, they are contractible and can be filled in with
disks $\bar{v}_i : \DD \to \RR \times M_+$ for which $\left| \int_\DD \bar{v}_i^*\omega_+ \right|$
may be assumed arbitrarily small.  Moreover, since all of the loops
$v_\delta(s,\cdot)$ are similarly contractible, the union of these two disks
with the annulus $v_\delta$ defines a closed cycle in $M_+$ that is trivial
in $H_2(M_+)$, hence the integral of the closed $2$-form $\omega_+$ over
this cycle vanishes, implying
$$
\int_{[s_0,s_1] \times S^1} v_\delta^*\omega_+ = 
\int_\DD \bar{v}_1^*\omega_+ - \int_\DD \bar{v}_0^*\omega_+,
$$
which is therefore arbitrarily small, and this proves the claim.

To finish, notice that since $v_\delta$ maps its boundary components to small
neighborhoods of $\{0\} \times M_+$ and $\{2\} \times M_+$, one can fix
a suitable choice of radius $r_1 > 0$ such that
$v_\delta$ must pass through a point in $p \in \{1\} \times M_+$ for which
the boundary of $v_\delta$ is outside the ball $B_{r_1}(p)$.  The monotonicity
lemma then bounds the symplectic area of $v_\delta$ from below by a constant
times $r_1^2$, but since we can also make this area arbitrarily small
by choosing $\delta$ smaller, this is a contradiction.

As usual, the case of a negative puncture can be handled similarly.
\end{proof}

We've now proved every statement in Theorem~\ref{thm:asympLocal} up to the
final detail about the case where the asymptotic orbit is nondegenerate
or Morse-Bott.  The complete proof of this part requires delicate analytical
results from \cites{HWZ:props1,HWZ:FIMpreprint,HWZ:props4,Bourgeois:thesis},
but we can explain the first step for the nondegenerate case.  In the following,
we say that a closed Reeb orbit $\gamma : \RR / T\ZZ \to M_\pm$ is 
\defin{isolated} if, after rescaling the domain to write it as an element
of $C^\infty(S^1,M_\pm)$, there exists a neighborhood
$\gamma \in \uU \subset C^\infty(S^1,M_\pm)$ such that all closed Reeb
orbits in $\uU$ are reparametrizations of~$\gamma$.

\begin{lemma}
\label{lemma:isolated}
Suppose the puncture is nonremovable, write 
$$
u_+(s,t) = (u_\RR(s,t),u_M(s,t)) \in \RR \times M_\pm
$$
for $s \ge 0$ sufficiently
large, and suppose $s_k \to \infty$ is a sequence and 
$\gamma : \RR / T\ZZ \to M_\pm$ is a Reeb orbit such that
$$
u_M(s_k,\cdot) \to \gamma(T\cdot) \quad \text{ in }\quad
C^\infty(S^1,M_\pm).
$$
If $\gamma$ is isolated, then for every neighborhood $\uU \subset C^\infty(S^1,M_\pm)$
of the set of parametrizations $\{ \gamma(\cdot + \theta)\ |\ \theta \in S^1 \}$,
we have $u_M(s,\cdot) \in \uU$ for all sufficiently large~$s$.
\end{lemma}
\begin{proof}
Note first that if $\gamma$ is isolated, then its image admits a neighborhood
$\im \gamma \subset \vV \subset M_\pm$ such that no point in
$\vV \setminus \im \gamma$ is contained in another Reeb orbit of period~$T$.
Indeed, we could otherwise find a sequence of $T$-periodic Reeb orbits
passing through a sequence of points in $\vV \setminus \im \gamma$ that
converge to a point in $\im \gamma$.  Since their derivatives are
determined by the Reeb vector field and are therefore bounded,
the Arzel\`a-Ascoli theorem then
gives a subsequence of these orbits converging to a reparametrization
of~$\gamma$, contradicting the assumption that $\gamma$ is isolated.

Arguing by contradiction,
suppose now that there exists a sequence $s_k' \to \infty$ with $u_M(s_k,\cdot) \not\in
\uU$ for all~$k$.  We can nonetheless restrict to a 
subsequence for which $u_M(s_k',\cdot)$
converges to some Reeb orbit $\tilde{\gamma} : \RR / T\ZZ \to M_\pm$.
Then $\tilde{\gamma}$ is disjoint from $\gamma$, and by continuity, one
can find a sequence $s_k'' \to \infty$ for which each $u_M(s_k'',0)$ lies
in the region $\vV$ some fixed distance away from $\im \gamma$.
There must then be a subsequence for which $u_M(s_k'',\cdot)$ converges to
another $T$-periodic orbit, but this is impossible since no such orbits
exist in $\vV \setminus \im \gamma$.
\end{proof}

\psfrag{What}{$\widehat{W}$}
\psfrag{RtimesM+}{$\RR \times M_+$}
\psfrag{RtimesM-}{$\RR \times M_-$}
\psfrag{z1}{$\zeta_1$}
\psfrag{z2}{$\zeta_2$}
\psfrag{z3}{$\zeta_3$}
\psfrag{uk}{$u_k$}
\psfrag{uinfty}{$u_\infty$}
\psfrag{vinfty1}{$v_\infty^1$}
\psfrag{vinfty3}{$v_\infty^3$}
\psfrag{u1+}{$v^1_+$}
\psfrag{u2+}{$v^2_+$}
\psfrag{u3+}{$v^3_+$}
\psfrag{u-}{$v_-$}
\psfrag{v1+}{$v^1_+$}
\psfrag{v2+}{$v^2_+$}
\psfrag{v0}{$v_0$}
\psfrag{v-}{$v_-$}
\psfrag{=}{$=$}
\psfrag{W}{$\widehat{W}$}
\psfrag{M+}{$M_+$}
\psfrag{M-}{$M_-$}

\section{Degenerations of holomorphic curves}
\label{sec:degenerations}

To motivate the SFT compactness theorem, we shall now discuss three examples
of phenomena that can prevent a sequence of holomorphic curves from having
a compact subsequence.  The theorem will then tell us that these three things
are, in essence, the only things that can go wrong.

Throughout this section and the next, 
assume $J_k \to J \in \jJ_\tau(\omega_h,r_0,\hH_+,\hH_-)$
is a $C^\infty$-convergent sequence of tame almost complex structures
on the completed cobordism~$\widehat{W}$.  More generally, one can also
allow the data $\omega$, $h$ and $\hH_\pm$ to vary in $C^\infty$-convergent
sequences, but let's not clutter the notation too much.  We shall denote
the restrictions of $J$ to the cylindrical ends by
$$
J_+ := J|_{[r_0,\infty) \times M_+} \in \jJ(\hH_+), \qquad
J_- := J|_{(-\infty,-r_0] \times M_-} \in \jJ(\hH_-).
$$
Suppose
$$
u_k := [(\Sigma_k,j_k,\Gamma_k^+,\Gamma_k^-,\Theta_k,u_k)] \in 
\mM_{g,m}(J_k,A_k,\boldsymbol{\gamma}^+,\boldsymbol{\gamma}^-)
$$
is a sequence of $J_k$-holomorphic curves in $\widehat{W}$ with fixed genus
$g \ge 0$ and $m \ge 0$ marked points, 
varying relative homology classes $A_k \in H_2(W,\bar{\boldsymbol{\gamma}}^+
\cup \bar{\boldsymbol{\gamma}}^-)$ and fixed collections of asymptotic orbits
$\boldsymbol{\gamma}^\pm = (\gamma_1^\pm,\ldots,\gamma^\pm_{m_\pm})$.
Observe that the energies $E(u_k)$ depend only on the orbits
$\boldsymbol{\gamma}^\pm$ and relative homology classes $A_k$, so in 
particular, $E(u_k)$ is uniformly bounded whenever the relative homology
class is also fixed.  The fundamental question of this section is:

\begin{question}
If $E(u_k)$ is uniformly bounded and no subsequence of $u_k$ 
converges to an element of $\mM_{g,m}(J,A,\boldsymbol{\gamma}^+,\boldsymbol{\gamma}^-)$
for any $A \in H_2(W,\bar{\boldsymbol{\gamma}}^+ \cup \bar{\boldsymbol{\gamma}}^-)$,
what can happen?
\end{question}

\subsection{Bubbling}
\label{sec:bubbling}

Suppose $(\Sigma_k,j_k,\Gamma_k^+,\Gamma_k^-,\Theta_k) =
(\Sigma,j,\Gamma^+,\Gamma^-,\theta)$ is a fixed sequence of domains,
and choose Riemannian metrics on $\dot{\Sigma} = \Sigma \setminus \Gamma$ 
and $\widehat{W}$ that are
translation-invariant on the cylindrical ends of both.  Suppose there
exists a point $\zeta_0 \in \dot{\Sigma}$ such that $u_k(\zeta_0)$ is
contained in a compact subset for all~$k$.  Suppose also that the maps
$u_k : \dot{\Sigma} \to \widehat{W}$ are locally $C^1$-bounded outside some
finite subset 
$$
\Gamma' = \{\zeta_1,\ldots,\zeta_N\} \subset \dot{\Sigma},
$$
i.e.~for every compact set
$K \subset \dot{\Sigma} \setminus \Gamma'$, there exists a constant
$C_K > 0$ independent of~$k$ such that
$$
|du_k| \le C_K \quad \text{ on~$K$}.
$$
Then elliptic regularity gives a subsequence that converges in 
$C^\infty_\loc(\dot{\Sigma} \setminus \Gamma')$ to a $J$-holomorphic curve
$$
u_\infty : \dot{\Sigma} \setminus \Gamma' \to \widehat{W}
$$
with $E(u_\infty) \le \limsup E(u_k) < \infty$, thus all the punctures
$\Gamma^+ \cup \Gamma^- \cup \Gamma'$ of $u_\infty$ are either
removable or positively or negatively asymptotic to Reeb orbits.  We cannot
be sure that the asymptotic behavior of $u_\infty$ at $\Gamma^\pm$ is the
same as for $u_k$, but let's assume this for now (\S\ref{sec:breaking}
below discusses some things that can happen if this does not hold).  Then to complete the
picture, we need to understand not only what $u_\infty$ is doing at
the additional punctures $\Gamma'$, but also what is happening to $u_k$
near these points as its first derivative blows up.  For this we can apply
the familiar rescaling trick: choose for each $\zeta_i$ a sequence 
$z_k^i \to \zeta_i$ such that $| d u_k(z_k^i) | =: R_k \to \infty$, 
along with a sequence $\epsilon_k \to 0$ with $\epsilon_k R_k \to \infty$,
and using Lemma~\ref{lemma:Hofer}, assume without loss of generality that
$|d u_k(z)| \le 2 R_k$ for all $z$ in the $\epsilon_k$-ball about~$z_k^i$.
For convenience, we can choose a holomorphic coordinate system identifying
a neighborhood of $\zeta_i$ with $\DD \subset \CC$ and placing $\zeta_i$
at the origin, so $z_k^i \to 0$ in these coordinates,
and assume without loss of generality that they identify our chosen
metric near~$\zeta_i$ with the Euclidean metric.
Now setting
$$
v_k^i(z) = u(z_k^i + z/R_k) \quad \text{ for } \quad z \in \DD_{\epsilon_k R_k}
$$
gives a sequence of $J_k$-holomorphic maps $v_k^i : \DD_{\epsilon_k R_k} \to
\widehat{W}$ whose energies and first derivatives are both uniformly
bounded.  As in the arguments of \S 2, we now have three possibilities:
\begin{itemize}
\item If $u_k^i(z_k^i)$ has a bounded subsequence, then the corresponding 
subsequence of $v_k^i$ converges in $C^\infty_\loc(\CC)$
to a $J$-holomorphic plane $v_\infty^i : \CC \to \widehat{W}$ with finite 
energy.
\item If $u_k(z_k^i)$ has a subsequence diverging to 
$\{\pm\infty\} \times M_\pm$, then translating $v_k^i$ by the
$\RR$-action produces a limiting finite-energy plane
$v_\infty^i$ in the positive/negative symplectization $\RR \times M_\pm$.
\end{itemize}
Viewing $\CC$ as the punctured sphere $S^2 \setminus \{\infty\}$,
the singularity of $v_\infty^i$ at $\infty$ may be removable,
in which case $v_\infty^i$ extends to a $J$-holomorphic sphere and we say
that $u_k$ has ``bubbled off a sphere'' at~$\zeta_i$.  Alternatively,
$v_\infty^i$ may be positively or negatively asymptotic to a Reeb orbit
at~$\infty$.  

\begin{figure}
\includegraphics{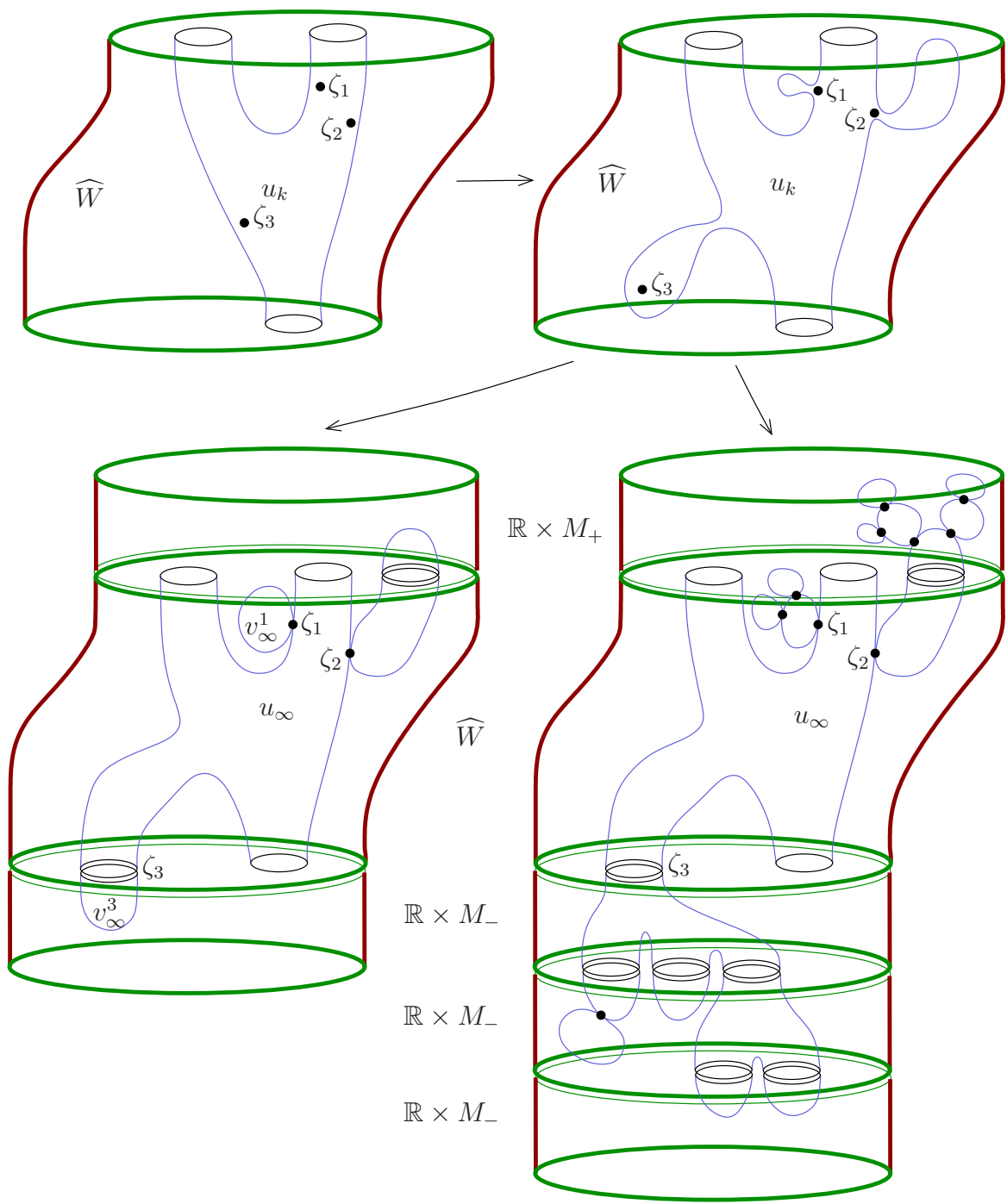}
\caption{\label{fig:bubbling} Two possible pictures of spheres and/or
planes that can bubble off when the first derivative blows up near three 
points.}
\end{figure}

Figure~\ref{fig:bubbling} shows two scenarios
that could occur for a sequence in which $| d u_k |$ blows up at
three points $\Gamma' = \{\zeta_1,\zeta_2,\zeta_3\}$.
Both scenarios show $u_\infty$ with $\zeta_1$ and $\zeta_2$ as removable
singularities and $\zeta_3$ as a negative puncture, but the behavior
of the various $v_\infty^i$ reveals a wide spectrum of possibilities.
In the lower-left picture, the points $u_k(z_k^1)$ 
are bounded and bubble off a sphere $v_\infty^1 : S^2 \to \widehat{W}$.
The picture shows that $v_\infty^1$ passes 
through $u_\infty(\zeta_1)$ at some point; this does
not follow from our argument so far, but in this situation one can use
a more careful analysis of $u_k$ near~$\zeta_1$ to show that it must
be true, i.e.~``bubbles connect''.  At $\zeta_3$, we have $u_k(z_k^3) \to
\{-\infty\} \times M_-$ and $v_\infty^3$ is a plane in $\RR \times M_-$ 
with a positive puncture
asymptotic to the same orbit as~$\zeta_3$; the coincidence of these orbits
is another detail that does not follow from the analysis above but
turns out to be true in the general picture.  The situation at $\zeta_2$
allows two different interpretations:
$v_\infty^2$ could be the plane with negative end
in $\RR \times M_+$, meaning $u_k(z_k^2) \to \{+\infty\} \times M_+$, and
the picture then shows an additional plane in $\widehat{W}$ with a positive
end approaching the same asymptotic orbit as $v_\infty^2$ as well as a point
passing through $u_\infty(\zeta_2)$.  One would need to choose a different
rescaled sequence near~$\zeta_2$ to find this extra plane, but as we will
see, the SFT compactness theorem dictates that some such object must be there.
Alternatively, $u_k(z_k^2)$ could also be bounded at $\zeta_2$, in which case
$v_\infty^2$ must be the plane in $\widehat{W}$ with positive end, and the
extra plane above this is something that one could find via a different
choice of rescaled sequence.  
In general, the range of actual possibilities
can involve arbitrarily many additional curves that could be discovered
via different choices of rescaled sequences:
e.g.~there could be entire
``bubble trees'' as shown in the lower-right picture, where each
$v_\infty^i$ is only one of several curves
that arise as limits of different parametrizations of~$u_k$ near~$\zeta_i$.
One good place to read about the analysis of bubble trees is
\cite{HWZ:foliations}*{\S 4}.

\subsection{Breaking}
\label{sec:breaking}

\begin{figure}
\includegraphics{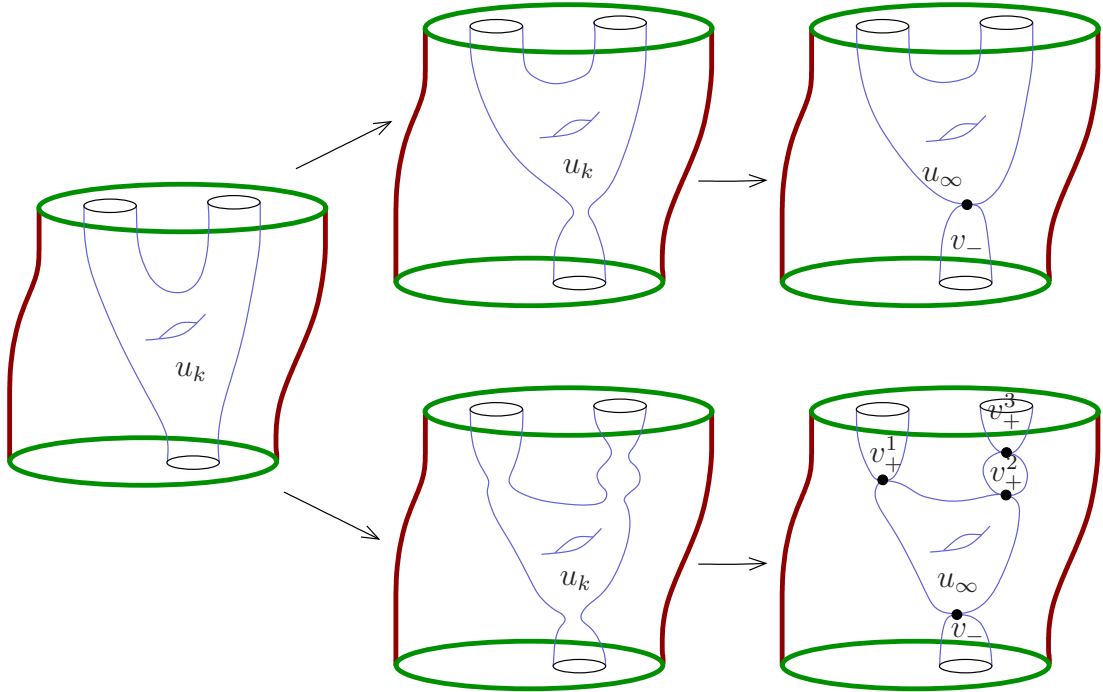}
\caption{\label{fig:breaking1} Even with fixed conformal structures on the 
domains and without bubbling, a sequence of
punctured holomorphic curves in $\widehat{W}$ can break to produce 
multiple curves in $\widehat{W}$ with extra removable punctures.
The picture shows two such scenarios.}
\end{figure}

\begin{figure}
\includegraphics{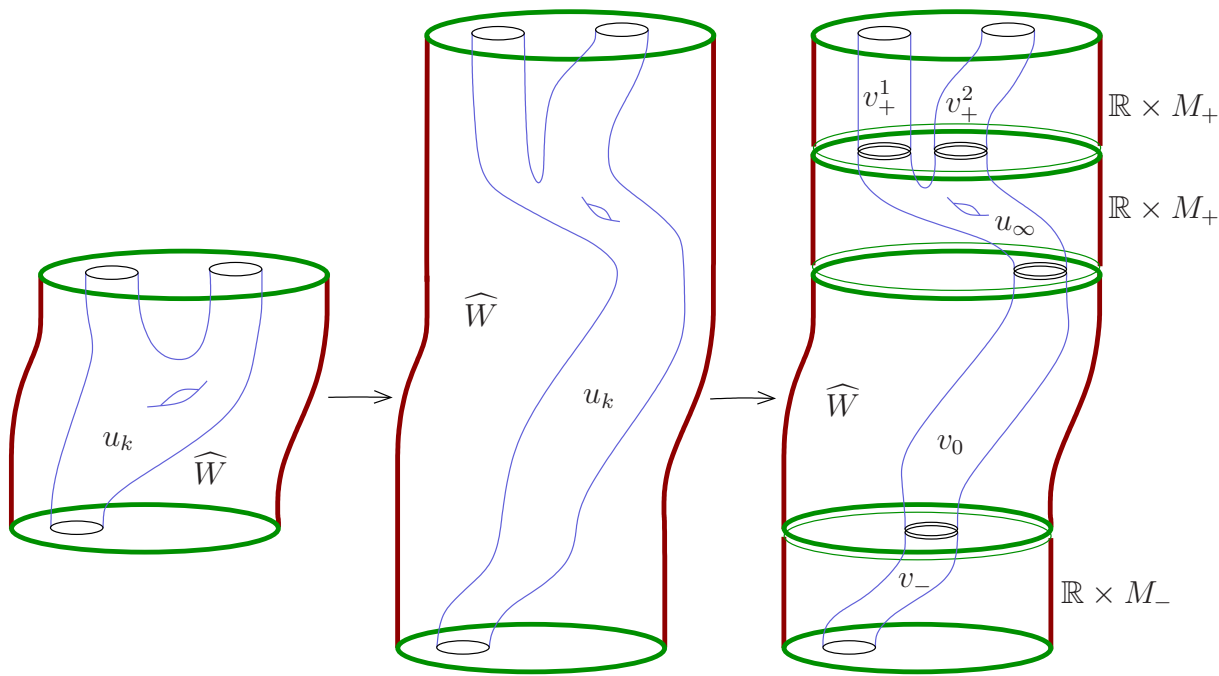}
\caption{\label{fig:breaking2} Different portions of a breaking sequence of
curves may also become infinitely far apart in the limit, so that some
live in~$\widehat{W}$ while others live in the symplectization
of $M_+$ or~$M_-$.}
\end{figure}

Figure~\ref{fig:bubbling} already shows some phenomena that could be interpreted
as ``breaking'' in the Floer-theoretic sense, but breaking can also happen
when no derivatives are blowing up, simply due to the fact that our domains
are noncompact.  Figures~\ref{fig:breaking1} and~\ref{fig:breaking2} show
three such scenarios, where we assume again that
$(\Sigma_k,j_k,\Gamma_k^+,\Gamma_k^-,\Theta_k) =
(\Sigma,j,\Gamma^+,\Gamma^-,\Theta)$ is a fixed sequence of domains, and
$\dot{\Sigma} = \Sigma \setminus \Gamma$ and $\widehat{W}$ carry Riemannian
metrics that are translation-invariant on the cylindrical ends such that
$$
|du_k| \le C \quad \text{ everywhere on $\dot{\Sigma}$}
$$
for some constant $C > 0$ independent of~$k$.  This is a stronger
condition than we had in \S\ref{sec:bubbling}, and if there exists a point
$\zeta_0 \in \dot{\Sigma}$ such that $u_k(\zeta_0)$ is bounded, it implies
that $u_\infty$ converges in $C^\infty_\loc(\dot{\Sigma})$ to a
$J$-holomorphic map
$$
u_\infty : \dot{\Sigma} \to \widehat{W}
$$
with $E(u_\infty) \le \limsup E(u_k) < \infty$.  Convergence in
$C^\infty_\loc$ is, however, not very strong: there may in general be no
relation between the asymptotic behavior of $u_\infty$ and $u_k$ at
corresponding punctures, e.g.~the top scenario in Figure~\ref{fig:breaking1} 
shows a case in which a negative puncture of $u_k$ becomes a removable
singularity of~$u_\infty$.  Whenever this happens, there must be more to
the story: in this example, one can choose holomorphic cylindrical coordinates
$(s,t) \in (-\infty,0] \times S^1 \subset \dot{\Sigma}$ 
near the negative puncture of $u_k$
and find a sequence $s_k \to \infty$ such that the sequence of
half-cylinders
$$
(-\infty,s_k] \times S^1 \to \widehat{W} : (s,t) \mapsto
u_k(s - s_k,t)
$$
is uniformly $C^1$-bounded and thus converges in $C^\infty_\loc(\RR \times S^1)$
to a finite-energy $J$-holomorphic cylinder
$v_- : \RR \times S^1 \to \widehat{W}$.  In the picture, $v_-$ turns out to
have a removable singularity at $+\infty$ mapping to the same point as
the removable singularity of $u_\infty$, and its negative puncture approaches
the same orbit as the negative puncture of~$u_k$.  

More complicated things
can happen in general: the bottom scenario in this same figure shows a
case where all three singularities of $u_\infty$ are removable, thus it
extends to a closed curve, while at one of the positive cylindrical ends
$[0,\infty) \times S^1 \subset \dot{\Sigma}$ of $u_k$, we can find a sequence $s_k \to \infty$
such that the half-cylinders
$$
[-s_k,\infty) \times S^1 \to \widehat{W} : (s,t) \mapsto
u_k(s + s_k,t)
$$
are uniformly $C^1$-bounded and converge in $C^\infty_\loc(\RR \times S^1)$
to a $J$-holomorphic cylinder $v_+^1 : \RR \times S^1 \to \widehat{W}$
with one removable singularity and one positive puncture.
At the other positive end, we can perform the same trick in two distinct
ways for two sequences $s_k \to \infty$, one diverging faster than the other:
the result is a pair of $J$-holomorphic cylinders
$v_+^2 , v_+^3 : \RR \times S^1 \to \widehat{W}$, the former with both
singularities removable (thus forming a holomorphic sphere in the picture),
and the latter with one removable singularity and one positive puncture.

It can get weirder.  Remember that $\widehat{W}$ is also noncompact!

In each of the above scenarios, we tacitly assumed that all of the various
sequences obtained by reparametrizing portions of $u_k$ were locally
$C^0$-bounded, thus all of the limits were curves in~$\widehat{W}$.  But
it may also happen that some of these sequences are $C^0_\loc$-bounded
while others locally diverge toward $\{\pm\infty\} \times M_\pm$; in fact,
two such sequences that both diverge toward, say, $\{+\infty\} \times M_+$,
might even locally diverge infinitely far from \emph{each other},
meaning one of them approaches $\{+\infty\} \times M_+$ quantitatively
faster than the other.  This phenomenon leads to the notion of limiting 
curves with multiple \emph{levels}.

In Figure~\ref{fig:breaking2}, we see a scenario in which $u_k$ satisfies
the same conditions as above, except that instead of $u_k(\zeta_0)$ being
bounded, it diverges to $\{+\infty\} \times M_+$.  It follows that after
applying suitable $\RR$-translations, a subsequence converges in
$C^\infty_\loc(\dot{\Sigma})$ to a $J_+$-holomorphic curve
$$
u_\infty : \dot{\Sigma} \to \RR \times M_+
$$
with finite energy.  In the example, all three of its punctures are 
nonremovable, but two of them approach orbits that have nothing to do with the
asymptotic orbits of~$u_k$.  Now observe that since $u_k$ has a negative
cylindrical end $(-\infty,0] \times S^1 \subset \dot{\Sigma}$, one can 
necessarily find a sequence $s_k \to \infty$ such that $u_k(-s_k,0)$
is bounded, and the sequence of half-cylinders
$$
(-\infty,s_k] \times S^1 \to \widehat{W} : (s,t) \mapsto
u_k(s - s_k,t)
$$
is then uniformly $C^1$-bounded and thus has a subsequence convergent 
in $C^\infty_\loc(\RR \times S^1)$ to a finite-energy $J$-holomorphic cylinder
$v_0 : \RR \times S^1 \to \widehat{W}$.  In the picture, $v_0$ has both
a positive and a negative puncture, but its negative end again approaches a
different Reeb orbit from the negative ends of~$u_k$, so one can deduce that
there must be still more happening near $-\infty$: there exists another
sequence $s_k' \to \infty$ with $s_k' - s_k \to \infty$ such that
suitable $\RR$-translations of the half-cylinders
$$
(-\infty,s_k] \times S^1 \to (-\infty,-r_0] \times M_- :
(s,t) \mapsto u_k(s - s_k',t)
$$
define uniformly $C^1$-bounded maps into $\RR \times M_-$, giving a
subsequence that converges in $C^\infty_\loc(\RR \times S^1)$ to a
finite-energy $J_-$-holomorphic cylinder
$$
v_- : \RR \times S^1 \to \RR \times M_-.
$$
Finally, the fact that $u_\infty$ has a positive asymptotic orbit different
from those of $u_k$ indicates that something more must also be happening 
near~$+\infty$: in the example, one of the positive ends 
$[0,\infty) \times S^1 \subset
\dot{\Sigma}$ admits a sequence $s_k \to \infty$ such that
$u_k(s_k,0) \in \{r_k\} \times M_+$ for some $r_k \to \infty$, and
suitable $\RR$-translations of
$$
[-s_k,\infty) \times S^1 \to [r_0,\infty) \times M_+ :
(s,t) \mapsto u_k(s + s_k,t)
$$
become a uniformly $C^1$-bounded sequence of half-cylinders in $\RR \times M_+$,
with a subsequence converging in $C^\infty_\loc(\RR \times S^1)$ to a
finite-energy $J_+$-holomorphic cylinder
$$
v_+^2 : \RR \times S^1 \to \RR \times M_+
$$
that connects the errant asymptotic orbit of $u_\infty$ to the corresponding
orbit of~$u_k$.  One can now perform the same trick at the other positive
end of $\dot{\Sigma}$, as there necessarily also exists a sequence
$s_k' \to \infty$ in this end such that $u_k(s_k',0) \in \{r_k\} \times M_+$
for the same sequence $r_k \to \infty$ as in the above discussion.
The resulting limit curve $v_+^1 : \RR \times S^1 \to \RR \times M_+$
however is not guaranteed to be interesting: in the picture, it turns out
to be a trivial cylinder.

The type of degeneration shown in Figure~\ref{fig:breaking2} happens
whenever the sequence $u_k$ does interesting things in multiple regions of
its domain that are sent increasingly far away from each other in the
image.  The usual picture of $\widehat{W}$ that collapses the cylindrical
ends to a finite size therefore becomes increasingly inadequate for
visualizing $u_k$ as $k \to \infty$: the middle picture in Figure~\ref{fig:breaking2}
deals with this by expanding the scale of the cylindrical ends so that the
convergence to upper and lower levels becomes visible.

\subsection{The Deligne-Mumford space of Riemann surfaces}
\label{sec:DM}

We next need to relax the assumption that the Riemann surfaces
$(\Sigma_k,j_k,\Gamma_k^+ \sqcup \Gamma_k^- \sqcup \Theta_k)$ are fixed.  Recall that for integers
$g \ge 0$ and $\ell \ge 0$, the moduli space of pointed Riemann surfaces is
the space of equivalence classes
$$
\mM_{g,\ell} = \left\{ (\Sigma,j,\Theta) \right\} \big/ \sim,
$$
where $(\Sigma,j)$ is a closed connected Riemann surface of genus~$g$,
$\Theta \subset \Sigma$ is an ordered set of $\ell$ distinct points,
and $(\Sigma,j,\Theta) \sim (\Sigma',j',\Theta')$ whenever there exists a
biholomorphic map $\varphi : (\Sigma,j) \to (\Sigma',j')$ taking
$\Theta$ to $\Theta'$ with the ordering preserved.  This space is fairly easy 
to understand in the finitely many cases with $2g + \ell < 3$, 
e.g.~$\mM_{0,\ell}$ is a one-point space for each $\ell \le 3$.  
We say that $(\Sigma,j,\Theta)$ is \defin{stable} whenever
$\chi(\Sigma \setminus \Theta) < 0$, which means $2g + \ell \ge 3$.
In the stable case, one can show that every pointed Riemann surface has
a finite automorphism group, and $\mM_{g,\ell}$ is a smooth orbifold of
dimension $6g - 6 + 2 \ell$.  It is generally not compact, but it admits
a natural compactification
$$
\overline{\mM}_{g,\ell} \supset \mM_{g,\ell},
$$
known as the \defin{Deligne-Mumford compactification}.  We shall now give
a sketch of this construction from the perspective of hyperbolic geometry;
for more details, see \cites{Hummel,SeppalaSorvali}.

We recall first the following standard result.

\begin{thmu}[Uniformization theorem]
Every simply connected Riemann surface is biholomorphically equivalent
to either the Riemann sphere $S^2 = \CC \cup \{\infty\}$, the complex
plane $\CC$ or the upper half plane $\HH = \{ \Im z > 0 \} \subset \CC$.
\end{thmu}

The uniformization theorem implies that every Riemann surface can be 
presented as a quotient of either $(S^2,i)$, $(\CC,i)$ or $(\HH,i)$ by some
freely acting discrete group of biholomorphic transformations.  The only
punctured surface $\dot{\Sigma} = \Sigma \setminus \Theta$ that has
$S^2$ as its universal cover is $S^2$ itself.  It is almost as easy to see 
which surfaces are covered by~$\CC$, as the only biholomorphic transformations 
on $(\CC,i)$ with no fixed points are the translations, so every freely 
acting discrete subgroup
of $\Aut(\CC,i)$ is either trivial, a cyclic group of translations or a
lattice.  The resulting quotients are, respectively, $(\CC,i)$,
$(\RR\times S^1,i) \cong (\CC\setminus\{0\},i)$ and the unpunctured
tori $(T^2,j)$.  All \emph{stable} pointed Riemann surfaces are thus
quotients of~$(\HH,i)$.

\begin{prop}
\label{prop:Poincare}
There exists on $(\HH,i)$ a complete Riemannian metric $g_P$ of constant
curvature~$-1$ that defines the same conformal structure as~$i$ and has
the property that all conformal transformations on $(\HH,i)$ are also
isometries of~$(\HH,g_P)$.
\end{prop}
\begin{proof}
We define $g_P$ at $z = x + iy \in \HH$ by
$$
g_P = \frac{1}{y^2} g_E,
$$
where $g_E$ is the Euclidean metric.  The conformal transformations
on $(\HH,i)$ are given by fractional linear transformations
\begin{equation*}
\begin{split}
\Aut(\HH,i) &= 
\left\{ \varphi(z) = \frac{a z + b}{c z + d} \ \Big|\ a,b,c,d \in \RR,
\quad ad - bc = 1 \right\} \bigg/ \{\pm 1\} \\
&= \SL(2,\RR) / \{\pm 1\} =: \PSL(2,\RR), \\
\end{split}
\end{equation*}
and one can check that each of these defines an isometry with respect
to~$g_P$.  One can also compute that $g_P$ has curvature~$-1$, and
the geodesics of~$g_P$ are precisely
the lines and semicircles that meet $\RR$ orthogonally, parametrized so that
they exist for all forward and backward time, thus $g_P$ is complete.
For more details on all of this, the book by Hummel \cite{Hummel} is
highly recommended.
\end{proof}
By lifting to universal covers, this implies the following.
\begin{cor}
\label{cor:Poincare}
For every pointed Riemann surface $(\Sigma,j,\Theta)$ such that
$\chi(\Sigma\setminus\Theta) < 0$, the punctured Riemann surface
$(\Sigma\setminus\Theta,j)$ admits a complete Riemannian metric $g_j$
of constant curvature~$-1$ that defines the same conformal structure as~$j$,
and has the property that all biholomorphic transformations on
$(\Sigma\setminus\Theta,j)$ are also isometries of $(\Sigma\setminus\Theta,g_j)$.
\end{cor}
The metric $g_j$ in this corollary is often called the 
\defin{Poincar\'{e} metric}.  It is uniquely determined by~$j$.

\begin{figure}
\includegraphics{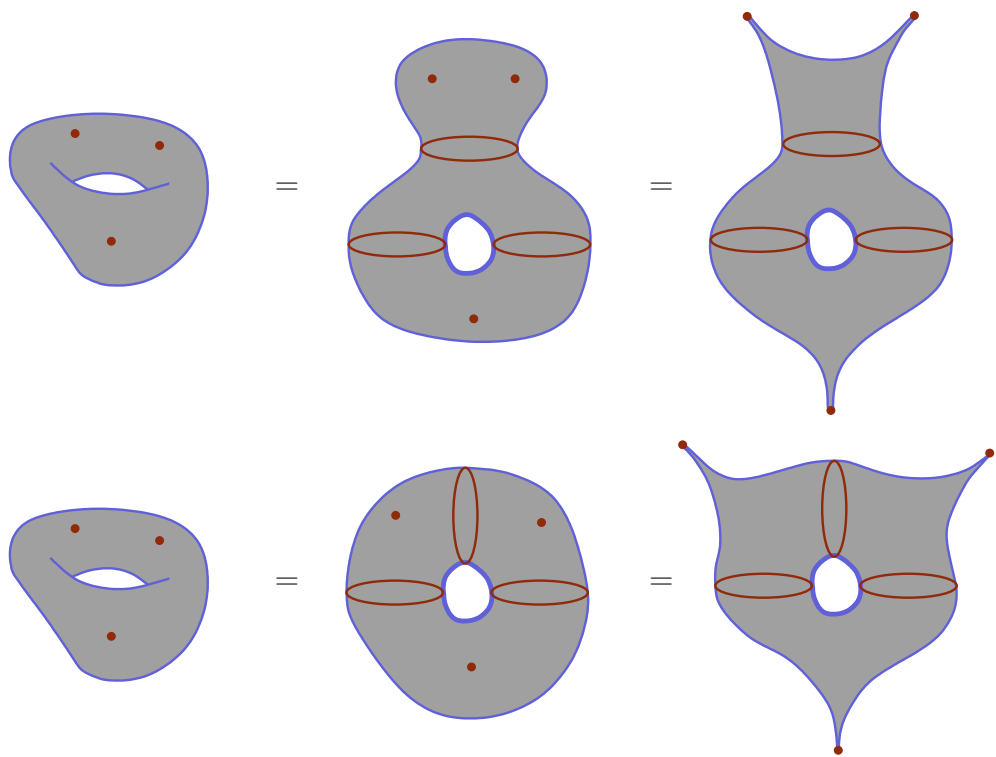}
\caption{\label{fig:POP} Two distinct pair-of-pants decompositions for the
same genus~$1$ Riemann surface with three marked points.  The decompositions
are shown from two perspectives: the pictures at the right are meant to give
a more accurate impression of the Poincar\'e metric, which becomes singular
and forms a cusp at each marked point.}
\end{figure}

Every class in $\pi_1(\dot{\Sigma})$ contains a unique geodesic for~$g_j$. 
Now suppose $C \subset \dot{\Sigma}$ is a union of disjoint 
embedded geodesics such that each connected component of 
$\dot{\Sigma} \setminus C$ has the
homotopy type of a disk with two holes.  The
components are then called \defin{singular pairs of pants}, and the result
is called a \defin{pair-of-pants decomposition} of $(\dot{\Sigma},j)$.
Two examples for the case $g=1$ and $\ell=3$ are shown in
Figure~\ref{fig:POP}.

A pair-of-pants decomposition for $(\Sigma,j,\Theta)$ gives rise to a local 
parametrization of $\mM_{g,\ell}$ near $[(\Sigma,j,\Theta)]$, known as
the \emph{Fenchel-Nielsen coordinates}.  These consist of two parameters
that can be associated to each of the geodesics $\gamma \subset \Sigma$
in the decomposition, namely the length $\ell(\gamma) > 0$ of the geodesic
and a \emph{twist} parameter $\theta(\gamma) \in S^1$, which describes how
the two neighboring pairs of pants are glued together along~$\gamma$.
Note that by computing Euler characteristics, there are always exactly
$-\chi(\Sigma \setminus \Theta) = 2g - 2 + \ell$ pairs of pants in a
decomposition, so that the total number of geodesics involved is
$\left[ 3 (2g - 2 + \ell) - \ell \right] / 2 = 3g - 3 + \ell$, thus one
can read off the formula $\dim \mM_{g,\ell} = 6g - 6 + 2\ell$ from this
geometric picture.

One can also see the noncompactness of $\mM_{g,\ell}$ in this picture
quite concretely: the twist parameters belong to a compact space, but
each length parameter can potentially shrink to~$0$ or blow up to $\infty$
as $j$ (and hence~$g_j$) is deformed.  It turns out that the latter possibility
is an illusion, but one may need to switch to a different pair-of-pants
decomposition to see why:

\begin{thmu}
For every pair of integers $g \ge 0$ and $\ell \ge 0$ with $2g + \ell \ge 3$,
there exists a constant $C = C(g,\ell) > 0$ such that every
$[(\Sigma,j,\Theta)] \in \mM_{g,\ell}$ admits a pair-of-pants decomposition
in which all geodesics bounding the pairs of pants have length at most~$C$.
\end{thmu}

\begin{figure}
\includegraphics{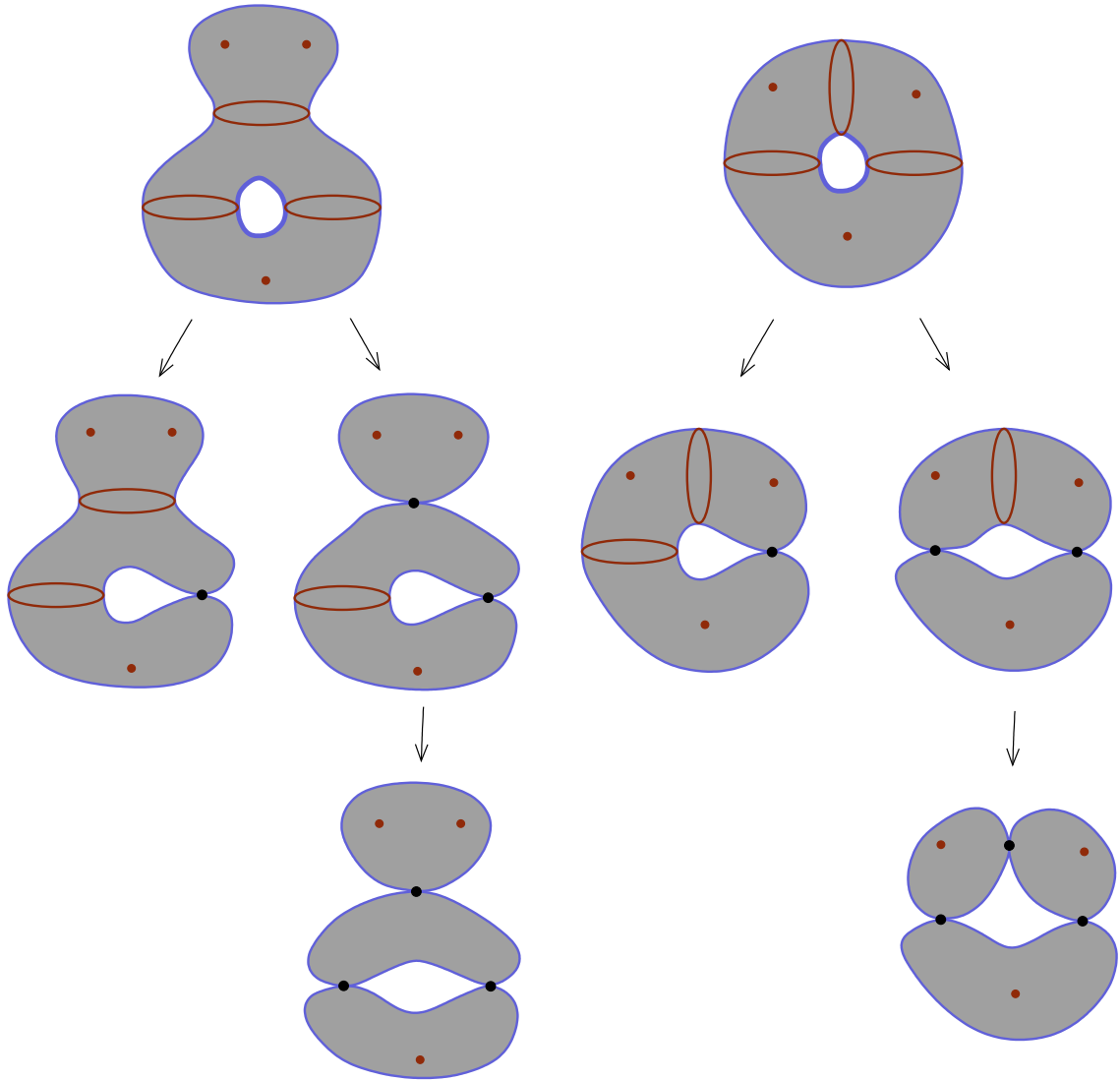}
\caption{\label{fig:POPsing} Starting from each of the pair-of-pants
decompositions for the $g=1$ and $\ell=3$ case from Figure~\ref{fig:POP},
shrinking geodesic lengths to zero produces various examples of stable nodal
Riemann surfaces belonging to~$\overline{\mM}_{1,3}$.}
\end{figure}

This theorem implies that from a hyperbolic perspective, the only 
meaningful way for stable pointed Riemann surfaces to degenerate is when
some of the bounding geodesics in a pair-of-pants decomposition shrink to
length zero.  Figure~\ref{fig:POPsing} shows several examples
of degenerate Riemann surfaces that can arise in this way for $g=1$
and $\ell=3$, giving elements of the space that we will now define
as~$\overline{\mM}_{1,3}$.

\begin{defn}
\label{defn:nodalRS}
A \defin{nodal Riemann surface} with $\ell \ge 0$ marked points 
and $N \ge 0$ \defin{nodes} is a tuple 
$(S,j,\Theta,\Delta)$ consisting of:
\begin{itemize}
\item A closed but not necessarily connected Riemann surface $(S,j)$;
\item An ordered set of $\ell$ points $\Theta \subset S$;
\item An unordered set of $2N$ points $\Delta \subset S \setminus \Theta$
equipped with an involution $\sigma : \Delta \to \Delta$.  Each pair
$\{z,\sigma(z)\}$ for $z \in \Delta$ is referred to as a \defin{node}.
\end{itemize}
Let $\widehat{S}$ denote the closed surface obtained by performing
connected sums on $S$ at each node $\{ z^+,z^- \} \subset \Delta$.
We then say that $(S,j,\Theta,\Delta)$ is \defin{connected}
if and only if $\widehat{S}$ is connected, and the genus of
$\widehat{S}$ is called the \defin{arithmetic genus} of
$(S,j,\Theta,\Delta)$.  We say that $(S,j,\Theta,\Delta)$ is \defin{stable}
if every connected component of $S \setminus (\Theta \cup \Delta)$ has
negative Euler characteristic.  Finally, two nodal Riemann surfaces
$(S,j,\Theta,\Delta)$ and $(S',j',\Theta',\Delta')$ are considered
\defin{equivalent} if there exists a biholomorphic map $\varphi : (S,j) \to
(S',j')$ taking $\Theta$ to $\Theta'$ with the ordering preserved and
taking $\Delta$ to $\Delta'$ such that nodes are mapped to nodes.
\end{defn}

The nodes $\{z^+,z^-\} \subset \Delta$ are typically represented in
pictures as self-intersections of~$S$, cf.~Figure~\ref{fig:POPsing}.
We can think of the \emph{stable} nodal surfaces as precisely those which
admit (possibly singular) pair-of-pants decompositions.
All nodal Riemann surfaces we consider will be assumed connected in the
sense defined above unless otherwise noted; note that $S$ itself
can nonetheless be disconnected, as is the case in four out of the six
nodal surfaces shown in Figure~\ref{fig:POPsing}.

We now introduce some further terminology and notation that will be useful
in the next section as well.  Whenever $\dot{\Sigma} = \Sigma \setminus \Gamma$
is obtained by puncturing a Riemann surface $(\Sigma,j)$ at finitely many
points $\Gamma \subset \Sigma$, we shall define the \defin{circle compactification}
$$
\overline{\Sigma} := \dot{\Sigma} \cup \bigcup_{z \in \Gamma} \delta_z,
$$
where for each $z \in \Gamma$, the circle $\delta_z$ is defined as a
``half-projectivization'' of the tangent space at~$z$:
$$
\delta_z := \left( T_z \Sigma \setminus \{0\} \right) \Big/ \RR_+^*,
$$
with the positive real numbers $\RR_+^*$ acting by scalar multiplication.
To understand the topology of $\overline{\Sigma}$, one can equivalently
define it by choosing holomorphic cylindrical coordinates
$[0,\infty) \times S^1 \subset \dot{\Sigma}$ near each~$z$, and replacing
the open half-cylinder with $[0,\infty] \times S^1$, where $\delta_z$ is
now the \defin{circle at infinity} $\{\infty\} \times S^1$.  There is no
natural choice of global smooth structure on $\overline{\Sigma}$, but it is
homeomorphic to an oriented surface with boundary and carries both smooth and 
conformal structures on its interior, due to the obvious identification
$$
\dot{\Sigma} = \overline{\Sigma} \setminus \bigcup_{z \in \Gamma} \delta_z
\subset \overline{\Sigma}.
$$
The conformal structure of $\Sigma$ at each $z \in \Gamma$ does induce on
each of the circles $\delta_z$ an \defin{orthogonal structure}, meaning a 
preferred class of homeomorphisms to $S^1$ that are all related to each other
by rotations.  One can therefore speak of \defin{orthogonal maps}
$\delta_z \to \delta_{z'}$ for $z , z' \in \Gamma$, which are always
homeomorphisms and can either preserve or reverse orientation.

Now if $(S,j,\Theta,\Delta)$ is a nodal Riemann surface, we let
$\dot{S} = S \setminus \Delta$ and form the circle compactification
$\overline{S}$, which has the topology of a compact oriented surface
with boundary.  Given a node $\{z^+,z^-\} \subset \Delta$, a
\defin{decoration} for $\{z^+,z^-\}$ is a choice of orientation reversing
orthogonal map 
$$
\Phi : \delta_{z^+} \to \delta_{z^-}.
$$
We say that $(S,j,\Theta,\Delta)$ is a \defin{decorated nodal surface}
if it is equipped with a choice of decoration $\Phi$ for every node, or
\defin{partially decorated} if $\Phi$ is defined for some subset of the
nodes.  A partial decoration $\Phi$ gives rise to another compact
oriented surface
$$
\widehat{S}_\Phi := \overline{S} \big/ \sim,
$$
where the equivalence relation identifies $\delta_{z^+}$ with
$\delta_{z^-}$ via~$\Phi$ for each decorated node $\{z^+,z^-\} \subset \Delta$.
Note that if every node is decorated, then $\widehat{S}_\Phi$ has the topology
of a closed
connected and oriented surface whose genus defines the arithmetic genus
of $(S,j,\Theta,\Delta)$ according to Definition~\ref{defn:nodalRS}.
We shall denote the collection of special circles in $\widehat{S}_\Phi$
where boundray components $\delta_{z^+},\delta_{z^-} \subset \p \overline{S}$
have been identified by
$$
C_\Phi \subset \widehat{S}_\Phi.
$$
Since $\widehat{S}_\Phi \setminus (\p \widehat{S}_\Phi \cup C_\Phi)$
has a natural identification with $\dot{S}$, it inherits smooth and
conformal structures which degenerate along $C_\Phi$ and~$\p\widehat{S}_\Phi$.
We will say that two partially decorated nodal Riemann surfaces
$(S,j,\Theta,\Delta,\Phi)$ and $(S',j',\Theta',\Delta',\Phi')$ are 
\defin{equivalent} if $(S,j,\Theta,\Delta)$ and $(S',j',\Theta',\Delta')$
are equivalent via a biholomorphic map $\varphi : (S,j) \to (S',j')$
that extends continuously from $\dot{S} \to \dot{S}'$
to a homeomorphism $\widehat{S}_\Phi \to \widehat{S}'_{\Phi'}$.

Now if $2g + \ell \ge 3$, define $\overline{\mM}_{g,\ell}$ as the set
of equivalence classes of stable nodal Riemann surfaces with $\ell$ marked 
points and arithmetic genus~$g$.  There is a natural inclusion 
$$
\mM_{g,\ell} \subset \overline{\mM}_{g,\ell}
$$
by regarding each pointed Riemann surface $(\Sigma,j,\Theta)$ as a nodal
Riemann surface $(\Sigma,j,\Theta,\Delta)$ with $\Delta = \emptyset$.
The most important property of $\overline{\mM}_{g,\ell}$ is that it admits
the structure of a compact metrizable topological space for which
the inclusion $\mM_{g,\ell} \hookrightarrow \overline{\mM}_{g,\ell}$ is 
continuous onto an open subset.  Rather than formulating all of this in 
precise terms, let us state the main corollary that is important to know
in practice.

\begin{thm}
\label{thm:DM}
Fix $g \ge 0$ and $\ell \ge 0$ with $2g + \ell \ge 3$.  Then for any
sequence $[(\Sigma_k,j_k,\Theta_k)] \in \mM_{g,\ell}$, there exists a
stable nodal Riemann surface
$[(S,j,\Theta,\Delta)] \in \overline{\mM}_{g,\ell}$ such that
after restricting to a subsequence,
$$
[(\Sigma_k,j_k,\Theta_k)] \to [(S,j,\Theta,\Delta)]
$$
in the following sense: $(S,j,\Theta,\Delta)$ admits a decoration~$\Phi$
such that for sufficiently large~$k$, there are homeomorphisms
$$
\varphi : \widehat{S}_\Phi \to \Sigma_k,
$$
smooth outside of $C_\Phi$, which map $\Theta$ to $\Theta_k$ preserving
the ordering and satisfy
$$
\varphi^*j_k \to j \quad\text{ in } \quad
C^\infty_\loc(\widehat{S}_\Phi \setminus C_\Delta).
$$
\end{thm}

As one might gather from the above statement, one could just as well define
a compact metrizable topology on the space of equivalence classes of
\emph{decorated} nodal Riemann surfaces and then characterize the topology 
of $\overline{\mM}_{g,\ell}$ via the natural projection that forgets
the decorations.

\begin{exercise}
The space $\mM_{0,4}$ has a natural identification with 
$S^2 \setminus \{0,1,\infty\}$, defined by choosing the unique
identification of any $4$-pointed Riemann sphere
$(S^2,j,(z_1,z_2,z_3,z_4))$ with $\CC \cup \{\infty\}$ such that
$z_1,z_2,z_3$ are identified with $0,1,\infty$ respectively, while 
$z_4$ is sent to some point in $S^2 \setminus \{0,1,\infty\}$.  Show that
this extends continuously to an identification of $\overline{\mM}_{0,4}$
with~$S^2$.  What do the three nodal curves in $\overline{\mM}_{0,4}
\setminus \mM_{0,4}$ look like in terms of pair-of-pants decompositions?
\end{exercise}

\section{The SFT compactness theorem}

We now introduce the natural compactification of 
$\mM_{g,m}(J,A,\boldsymbol{\gamma}^+,\boldsymbol{\gamma}^-)$.

\subsection{Nodal curves}
\label{sec:nodal}

A punctured $J$-holomorphic \defin{nodal curve} in $(\widehat{W},J)$ 
with $m \ge 0$ marked points consists of 
the data $(S,j,\Gamma^+,\Gamma^-,\Theta,\Delta,u)$, where
\begin{itemize}
\item $(S,j,\Gamma^+ \sqcup \Gamma^- \sqcup \Theta,\Delta)$ is a 
nodal Riemann surface, with $|\Theta| = m$;
\item $u : (\dot{S},j) \to (\widehat{W},J)$ for 
$\dot{S} := S \setminus (\Gamma^+ \cup \Gamma^-)$
is an asymptotically cylindrical $J$-holomorphic map with positive
punctures $\Gamma^+$ and negative punctures $\Gamma^-$ such that for each
node $\{z^+,z^-\} \subset \Delta$, $u(z^+) = u(z^-)$.
\end{itemize}
Equivalence of two nodal curves
$$
(S_0,j_0,\Gamma^+_0,\Gamma^-_0,\Theta_0,\Delta_0,u_0) \sim
(S_1,j_1,\Gamma^+_1,\Gamma^-_1,\Theta_1,\Delta_1,u_1)
$$
is defined as the existence of an equivalence of nodal Riemann surfaces
$\varphi : (S_0,j_0,\Gamma^+_0 \sqcup \Gamma^-_0 \sqcup \Theta_0,\Delta_0) \to
(S_1,j_1,\Gamma^+_1 \sqcup \Gamma^-_1 \sqcup \Theta_1,\Delta_1)$ such that
$u_0 = u_1 \circ \varphi$.  We say that 
$(S,j,\Gamma^+,\Gamma^-,\Theta,\Delta,u)$ is \defin{connected} if and only
if the nodal Riemann surface $(S,j,\Gamma^+ \sqcup \Gamma^- \sqcup \Theta,\Delta)$
is connected, and its \defin{arithmetic genus} is then defined to be
the arithmetic genus of the latter.  We say that
$(S,j,\Gamma^+,\Gamma^-,\Theta,\Delta,u)$ is
\defin{stable} if every
connected component of $S \setminus (\Gamma^+ \cup \Gamma^- \cup \Theta \cup \Delta)$
on which $u$ is constant has negative Euler characteristic.  Note that the
underlying nodal Riemann surface $(S,j,\Gamma^+ \sqcup \Gamma^- \sqcup \Theta,\Delta)$
need not be stable in general.

Nodal curves are sometimes also referred to as \emph{holomorphic buildings
of height~1}.  These are the objects that form the \emph{Gromov 
compactification} of $\mM_{g,m}(J,A)$ when $W$ is a closed symplectic manifold.
One can now roughly imagine how the compactness theorem in that setting is
proved: given a converging sequence of almost complex structures $J_k \to J$
and a sequence $[(\Sigma_k,j_k,\Theta_k,u_k)] \in
\mM_{g,m}(J_k,A_k)$ with uniformly bounded energy, we can first add some
auxiliary marked points if necessary to assume that $2g + m \ge 3$.
Now a subsequence of the domains $[(\Sigma_k,j_k,\Theta_k)] \in \mM_{g,m}$
converges to an element of the Deligne-Mumford space
$[(S,j,\Theta,\Delta)] \in \overline{\mM}_{g,m}$.  Concretely, this means
that for large~$k$,
our sequence in $\mM_{g,m}(J_k,A_k)$ admits representatives
$(\Sigma,j_k',\Theta,u_k')$, with $\Sigma$ a fixed surface with fixed marked
points $\Theta \subset \Sigma$, and $(S,j,\Theta,\Delta)$ admits decorations
$\Phi$ so that one can identify $\widehat{S}_\Phi$ with $\Sigma$ and find
$$
j_k' \to j \quad \text{ in } \quad C^\infty_\loc(\Sigma \setminus C)
$$
for some collection of disjoint circles $C \subset \Sigma$.
The connected components of $(\Sigma \setminus C,j)$ are then biholomorphically
equivalent to the connected components of $(S \setminus \Delta,j)$, and
if the newly reparametrized maps $u_k' : \Sigma \to W$ are
uniformly $C^1_\loc$-bounded on $\Sigma \setminus C$, then a subsequence 
converges in $C^\infty_\loc(\Sigma \setminus C)$ to a limiting finite-energy
$J$-holomorphic map
$u_\infty : (S \setminus \Delta,j) \to (W,J)$, whose singularities 
at $\Delta$ are removable.  In particularly nice cases, this may be the
end of the story, and our subsequence of $[(\Sigma_k,j_k,\Theta_k,u_k)]
\in \mM_{g,m}(J_k,A_k)$ converges to the nodal curve
$[(S,j,\Theta,\Delta,u_\infty)]$; in particular the domain 
$[(S,j,\Theta,\Delta)]$ in this case is stable and is thus an
element of $\overline{\mM}_{g,m}$.  But more complicated things can also
happen, e.g.~$u_k'$ might not be $C^1$-bounded, in which case there is
bubbling.  The bubbles that arise will be either planes or spheres, 
so they produce extra domain
components with nonnegative Euler characteristic, but since
they are never constant, the limiting nodal curve is still considered
stable.  Similarly, since $\Sigma \setminus C$ is not compact, there can
also be breaking as in Figure~\ref{fig:breaking1}, producing more 
non-stable domain components which can be cylinders in addition to planes
and spheres---but again, the limiting map on these components will never
be constant.

\subsection{Holomorphic buildings}
\label{sec:buildings}

Only a small subset of the phenomena observed in \S\ref{sec:degenerations}
can be described via nodal curves: we've seen that in general, we also have
to allow ``broken'' curves with multiple ``levels''.  This notion can be
formalized as follows.

Given integers
$g,m,N_+,N_- \ge 0$, a \defin{holomorphic building of height $N_-|1|N_+$}
with arithmetic genus $g$ and $m$ marked points is a tuple
$$
\mathbf{u} = (S,j,\Gamma^+,\Gamma^-,\Theta,\Delta\node,\Delta\br,L,\Phi,u),
$$
with the various data defined as follows:
\begin{itemize}
\item The \defin{domain} 
$(S,j,\Gamma^+ \sqcup \Gamma^- \sqcup \Theta,\Delta\node \sqcup \Delta\br)$ is
a connected but not necessarily stable nodal Riemann surface of
arithmetic genus~$g$, where $|\Theta| = m$, and the involution on
$\Delta\node \sqcup \Delta\br$ is assumed to preserve the subsets $\Delta\node$
and~$\Delta\br$.  Matched pairs in these subsets are called the \defin{nodes}
and \defin{breaking pairs} respectively of~$\mathbf{u}$.  The \defin{marked
points} of $\mathbf{u}$ are the points in $\Theta$, while
$\Gamma^+$ and $\Gamma^-$ are its positive and negative \defin{punctures}
respectively.
\item The \defin{level structure} is a locally constant function
$$
L : S \to \{-N_-,\ldots,-1,0,1,\ldots,N_+\}
$$
that attains every value in $\{-N_-,\ldots,N_+\}$ except possibly~$0$, and
satisfies:
\begin{enumerate}
\item
$L(z^+) = L(z^-)$ for each node $\{ z^+,z^- \} \subset \Delta\node$;
\item 
Each breaking pair $\{ z^+,z^- \} \subset \Delta\br$ can be labelled
such that $L(z^+) - L(z^-) = 1$;
\item
$L(\Gamma^+) = \{N_+\}$ and $L(\Gamma^-) = \{-N_-\}$.
\end{enumerate}
\item
The \defin{decoration} is a choice of orientation-reversing orthogonal
map 
$$
\delta_{z^+} \stackrel{\Phi}{\longrightarrow} \delta_{z^-}
$$ 
for each breaking pair $\{ z^+,z^- \} \subset \Delta\br$.
\item 
The \defin{map} is an asymptotically cylindrical pseudoholomorphic curve
$$
u : (\dot{S} := S \setminus (\Gamma^+ \cup \Gamma^- \cup \Delta\br),j) \to
\bigsqcup_{N \in \{-N_-,\ldots,N_+\}} (\widehat{W}_N,J_N),
$$
where
$$
(\widehat{W}_N,J_N) := 
\begin{cases}
(\RR \times M_+,J_+) & \text{ for $N \in \{1,\ldots,N_+\}$},\\
(\widehat{W},J)      & \text{ for $N = 0$},\\
(\RR \times M_-,J_-) & \text{ for $N \in \{-N_-,\ldots,-1\}$},
\end{cases}
$$
and $u$ sends $\dot{S} \cap L^{-1}(N)$ into $\widehat{W}_N$ for each~$N$,
with positive punctures at $\Gamma^+$ and negative punctures at~$\Gamma^-$.
Moreover,
$$
u(z^+) = u(z^-) \quad \text{ for every node $\{z^+,z^-\} \subset \Delta\node$},
$$
and for each breaking pair $\{z^+,z^-\} \subset \Delta\br$ labelled with
$L(z^+) - L(z^-) = 1$, $u$ has a positive puncture at $z^-$ and a negative
puncture at $z^+$ asymptotic to the same orbit, such that
if $u_+ : \delta_{z^+} \to M_\pm$ and $u_- : \delta_{z^-} \to M_\pm$
denote the induced asymptotic parametrizations of the orbit, then
$$
u_+ = u_- \circ \Phi : \delta_{z^+} \to M_\pm.
$$
\end{itemize}

The following additional notation and terminology for the building
$\mathbf{u}$ will be useful to keep in mind.  For each
$N \in \{-N_-,\ldots,0,\ldots,N_+\}$, denote
$$
\dot{S}_N := \left( S \setminus (\Gamma^+ \cup \Gamma^- \cup \Delta\br) \right)
\cap L^{-1}(N),
$$
and denote the restriction of $u$ to this subset by
$$
u^N : \dot{S}_N \to \begin{cases}
\RR \times M_+ & \text{ if $N > 0$},\\
\widehat{W} & \text{ if $N=0$},\\
\RR \times M_- & \text{ if $N < 0$}.
\end{cases}
$$
Including $\Theta \cap L^{-1}(N)$ and $\Delta\node \cap L^{-1}(N)$ in the data 
defines $u^N$ as a (generally disconnected) nodal curve with 
marked points, whose positive punctures are in bijective correspondence with
the negative punctures of $u^{N+1}$ if $N < N_+$.  We call $u_N$ the
\defin{$N$th level} of~$\mathbf{u}$, and all it an \defin{upper} or \defin{lower}
level if $N > 0$ or $N < 0$ respectively, and the \defin{main level} if
$N=0$.  By convention, every holomorphic building
in $\widehat{W}$ has exactly one main level (which lives in $\widehat{W}$ itself)
and arbitrary nonnegative numbers
of upper and lower levels (which live in the symplectizations $\RR \times M_\pm$).
One slightly subtle detail is that it is possible for the main
level to be \emph{empty}, meaning $0$ is not in the image of the level
function~$L$.  The requirement
that $L$ should attain every other value from $-L_-$ to~$L_+$ is a convention
to ensure that upper and lower levels are not empty, so e.g.~if a building 
has an empty main level and no lower levels, 
then the lowest nonempty upper level is always labelled~$1$ instead of 
something arbitrary.

The positive
punctures of the topmost level of $\mathbf{u}$ are $\Gamma^+$, and the 
negative punctures of the bottommost level are $\Gamma^-$, so these give rise 
to lists of positive/negative asymptotic orbits 
$\boldsymbol{\gamma}^\pm = (\gamma_1^\pm,\ldots,\gamma_{k_\pm}^\pm)$ 
in~$M_\pm$.  There is also a relative homology class
$$
[\mathbf{u}] \in H_2(W,\bar{\boldsymbol{\gamma}}^+ \cup \bar{\boldsymbol{\gamma}}^-).
$$
To define this, recall from \S\ref{sec:stableBoundary} how it was defined for smooth
curves $u : \dot{\Sigma} \to \widehat{W}$:
we considered the retraction $\pi : \widehat{W} \to W$ that collapses each
cylindrical end to $M_\pm \subset \p W$, and noted that since $u$ is
asymptotically cylindrical, the map $\pi \circ u : \dot{\Sigma} \to W$
extends to a continuous map on the circle compactification,
$$
\bar{u} : \overline{\Sigma} \to W,
$$
whose relative homology class gives the definition of~$[u]$.
The conditions on nodes and breaking orbits allow us to perform a
similar trick for the building~$\mathbf{u}$, using the map
$$
\pi : \bigsqcup_{N \in \{-N_-,\ldots,N_+\}} \widehat{W}_N \to W
$$ 
which acts as the identity on~$W$ but collapses cylindrical ends of
$\widehat{W}$ to $\p W$ and similarly collapses each copy of
$\RR \times M_\pm$ to $M_\pm \subset \p W$.  Extending the decorations $\Phi$
arbitrarily to decorations of the nodes $\Delta\node$, one can then take
the circle compactification of $\dot{S} := S \setminus (\Gamma^+ \cup \Gamma^- \cup
\Delta\node \cup \Delta\br)$ and glue matching boundary components together
along~$\Phi$ to form a compact surface with boundary $\overline{S}_\Phi$
such that $\pi \circ u : \dot{S} \to W$ extends to a continuous map
$$
\bar{u} : \overline{S}_\Phi \to W.
$$
Its relative homology class defines
$[\mathbf{u}] \in H_2(W,\bar{\boldsymbol{\gamma}}^+ \cup \bar{\boldsymbol{\gamma}}^-)$.

We say that the building $\mathbf{u}$ is \defin{stable} if two properties hold:
\begin{enumerate}
\item Every connected component of $S \setminus (\Gamma^+ \cup \Gamma^- \cup \Theta \cup \Delta\node \cup \Delta\br)$ on
which the map $u$ is constant has negative Euler characteristic;
\item There is no $N \in \{-N_-,\ldots,N_+\}$ for which the $N$th level
consists entirely of a disjoint union of
trivial cylinders without any marked points or nodes.
\end{enumerate}
An \defin{equivalence} between two holomorphic buildings
$$
\mathbf{u}_i = (S_i,j_i,\Gamma^+_i,\Gamma^-_i,\Theta_i,\Delta\node_i,\Delta\br_i,L_i,\Phi_i,u_i), \qquad
i=0,1
$$
is defined as an equivalence of partially decorated nodal Riemann surfaces
$$
(S_0,j_0,\Gamma^+_0 \sqcup \Gamma^+_0 \sqcup \Theta_0,\Delta\node_0 \sqcup \Delta\br_0,\Phi_0) 
\stackrel{\varphi}{\longrightarrow}
(S_1,j_1,\Gamma^+_1 \sqcup \Gamma^+_1 \sqcup \Theta_1,\Delta\node_1 \sqcup \Delta\br_1,\Phi_1)
$$
such that $\varphi(\Gamma^\pm_0) = \Gamma^\pm_1$, $\varphi(\Theta_0) = \Theta_1$,
$\varphi(\Delta\node_0) = \Delta\node_1$, $\varphi(\Delta\br_0) = \Delta\br_1$,
$L_1 \circ \varphi = L_0$, and 
$$
u_1^0 \circ \varphi = u_0^0,
$$
while
$$
u_1^N \circ \varphi = u_0^N \text{ up to $\RR$-translation} \quad 
\text{for each $N \ne 0$}.
$$

Given lists of orbits $\boldsymbol{\gamma}^\pm$ and a relative homology
class~$A$, the set of equivalence classes of stable holomorphic buildings 
in $(\widehat{W},J)$ with 
arithmetic genus~$g$ and $m$ marked points, positively/negatively asymptotic
to $\boldsymbol{\gamma}^\pm$ and homologous to $A$ will be denoted by
$$
\overline{\mM}_{g,m}(J,A,\boldsymbol{\gamma}^+,\boldsymbol{\gamma}^-).
$$
Observe that for any $A \ne 0$, there is a natural inclusion
$\mM_{g,m}(J,A,\boldsymbol{\gamma}^+,\boldsymbol{\gamma}^-) \subset
\overline{\mM}_{g,m}(J,A,\boldsymbol{\gamma}^+,\boldsymbol{\gamma}^-)$
defined by regarding $J$-holomorphic curves in
$\mM_{g,m}(J,A,\boldsymbol{\gamma}^+,\boldsymbol{\gamma}^-)$ as buildings with
no upper or lower levels and no nodes.  Such buildings are always stable
if $A \ne 0$ because they are not constant.

\subsection{Convergence}
\label{sec:convergence}

\begin{figure}
\includegraphics[width=6in]{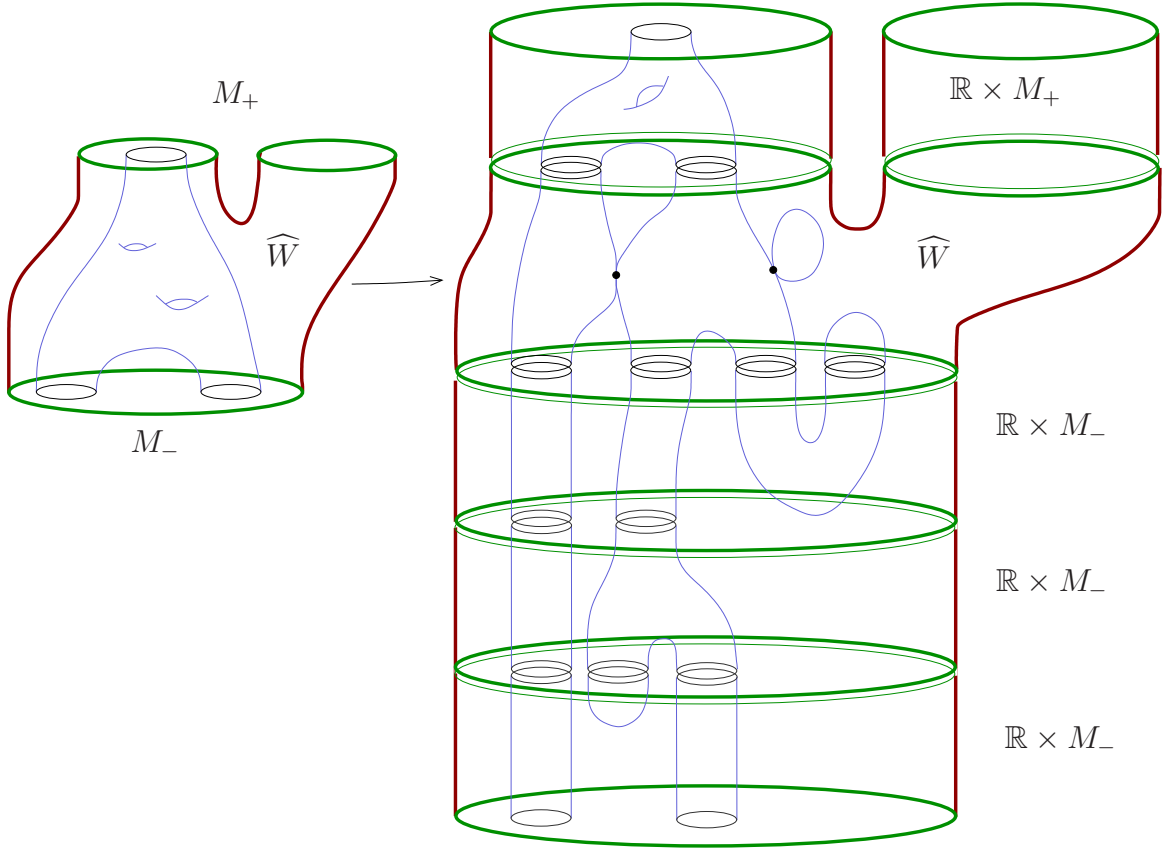}
\caption{\label{fig:SFTcomp} Convergence to a building with arithmetic
genus~$2$, one upper level and three lower levels.}
\end{figure}

For a general definition of the topology of
$\overline{\mM}_{g,m}(J,A,\boldsymbol{\gamma}^+,\boldsymbol{\gamma}^-)$
and the proof that it is both compact and metrizable, we refer to
\cite{SFTcompactness} or the more comprehensive treatment
in \cite{Abbas:book}.  The following statement contains all the details
about the topology that one usually needs to know in practice
(see Figure~\ref{fig:SFTcomp}).

\begin{thm}
\label{thm:SFTcompactness}
Fix integers $g \ge 0$ and $m \ge 0$, and assume all Reeb orbits in
$(M,\hH_+)$ and $(M,\hH_-)$ are nondegenerate.  Then for any
sequence
$$
[(\Sigma_k,j_k,\Gamma_k^+,\Gamma_k^-,\Theta_k,u_k)] \in 
\mM_{g,m}(J_k,A_k,\boldsymbol{\gamma}^+,\boldsymbol{\gamma}^-)
$$
of nonconstant $J_k$-holomorphic curves in $\widehat{W}$ with
uniformly bounded energy~$E(u_k)$, there exists a stable holomorphic
building
$$
[\mathbf{u}_\infty] = [(S,j,\Gamma^+,\Gamma^-,\Theta,\Delta\node,\Delta\br,L,\Phi,u_\infty)]
\in \overline{\mM}_{g,m}(J,A,\boldsymbol{\gamma}^+,\boldsymbol{\gamma}^-)
$$
such that after restricting to a subsequence, 
$[(\Sigma_k,j_k,\Gamma_k^+,\Gamma_k^-,\Theta_k,u_k)] \to [\mathbf{u}_\infty]$
in the following sense.  
The decorations $\Phi$ at $\Delta\br$ can be
extended to decorations at $\Delta\node$ so that if $\widehat{S}_\Phi$ denotes
the closed oriented topological $2$-manifold obtained from $S \setminus (\Delta\node \cup \Delta\br)$
by gluing circle compactifications along~$\Phi$, then for $k$ sufficiently
large, there exist homeomorphisms
$$
\varphi_k : \widehat{S}_\Phi \to \Sigma_k
$$
that are smooth outside of $C_\Phi$, map
$\Gamma^+ \sqcup \Gamma^- \sqcup \Theta$ to $\Gamma^+_k \sqcup
\Gamma^-_k \sqcup \Theta_k$ with the ordering preserved, and satisfy
$$
\varphi_k^*j_k \to j \quad \text{ in } \quad
C^\infty_\loc(\widehat{S}_\Phi \setminus C_\Phi).
$$
Moreover for $N=\{-N_-,\ldots,0,\ldots,N\}$, let
$$
v_k^N := u_k \circ \varphi_k|_{\ddot{S}_N} : \ddot{S}_N \to \widehat{W},
$$
with $\ddot{S}_N := \left(S \setminus (\Gamma^+ \cup \Gamma^- \cup \Delta\node \cup \Delta\br)\right) \cap L^{-1}(N)$
regarded as a subset of $\widehat{S}_\Phi \setminus C_\Phi$.
Then:
\begin{enumerate}
\item $v_k^0 \to u_\infty^N$ in $C^\infty_\loc(\ddot{S}_N,\widehat{W})$;
\item For each $\pm N > 0$, $v_k^N$ has image in the positive/negative 
cylindrical end
for all $k$ sufficiently large, and there exists a sequence 
$r_k^N \to \pm\infty$ such that
the resulting $\RR$-translations converge:
$$
\tau_{-r_k^N} \circ v_k^N \to u_\infty^N \quad \text{ in }\quad
C^\infty_\loc(\ddot{S}_N,\RR \times M_\pm).
$$
\end{enumerate}
The rates of divergence of the sequences $r_k^N \to \pm\infty$ are related by
$$
r_k^{N+1} - r_k^N \to +\infty \quad \text{ for all $N < N_+$}.
$$
Finally, let $\overline{S}_\Phi$ denote the compact topological surface with
boundary defined as the circle compactification of $\widehat{S}_\Phi
\setminus (\Gamma^+ \cup \Gamma^-)$, and let $\overline{\Sigma}_k$ denote the
circle compactification of $\dot{\Sigma}_k := \Sigma_k \setminus (\Gamma_k^+ \cup \Gamma_k^-)$.
Then for all $k$ large, $\varphi_k$ extends to a continuous map
$$
\bar{\varphi}_k : \overline{S}_\Phi \to \overline{\Sigma}_k
$$
such that
$$
\bar{u}_k \circ \bar{\varphi}_k \to \bar{u}_\infty \quad \text{ in } \quad
C^0(\overline{S}_\Phi,W).
$$
\end{thm}

\begin{remark}
The theorem is also true under the more general hypothesis that
the Reeb vector fields are Morse-Bott.  In this case, one can also allow
the asymptotic Reeb orbits of the sequence to vary, as long as the sum of
their periods is uniformly bounded---such a bound plays the role of an
energy bound and guarantees a convergent subsequence of orbits via the
Arzel\`a-Ascoli theorem.
\end{remark}

\begin{remark}
Stability of the limit in Theorem~\ref{thm:SFTcompactness} is guaranteed
for the same reasons as in our discussion of Gromov compactness in
\S\ref{sec:nodal}: stable domains degenerate to stable nodal domains as
geodesics in pair-of-pants decompositions shrink to zero length, while
bubbling and breaking produce additional domain components that are not
stable but on which the maps are never trivial.  Moreover, stability
guarantees the \emph{uniqueness} of the limiting building for any
convergent sequence, i.e.~it is the reason why
$\overline{\mM}_{g,m}(J,A,\boldsymbol{\gamma}^+,\boldsymbol{\gamma}^-)$
is a Hausdorff space.  Indeed, if $u_k$ converges to a stable building
$\mathbf{u}_\infty$, then under the notion of convergence described in the
theorem, it will also converge to a building $\mathbf{u}'_\infty$ constructed
out of $\mathbf{u}_\infty$ by adding to $S$ an extra spherical component,
attaching it to the rest by a single node and extending the map $u_\infty$ to 
be constant on the extra component.  One can also insert extra levels into 
$\mathbf{u}_\infty$ that consist only of trivial cylinders, and $u_k$ will still
converge to the resulting building.  But these modifications produce buildings
that are not stable and thus are not elements of
$\overline{\mM}_{g,m}(J,A,\boldsymbol{\gamma}^+,\boldsymbol{\gamma}^-)$.
\end{remark}

\subsection{Symplectizations, stretching and so forth}

A few minor modifications to the above discussion are necessary to compactify
the moduli space of curves in a symplectization $(\RR \times M,J)$ for
$J \in \jJ(\hH)$.  It is possible to view this as a special case of a
completed symplectic cobordism, but this perspective produces a certain
amount of extraneous data that is not meaningful.  The key observation is
that in the presence of an $\RR$-action, one should really compactify
$\mM_{g,m}(J,A,\boldsymbol{\gamma}^+,\boldsymbol{\gamma}^-) \big/ \RR$
instead of $\mM_{g,m}(J,A,\boldsymbol{\gamma}^+,\boldsymbol{\gamma}^-)$.
The compactification 
$\overline{\mM}_{g,m}(J,A,\boldsymbol{\gamma}^+,\boldsymbol{\gamma}^-)$ then
consists of holomorphic buildings as defined in \S\ref{sec:buildings},
but since all levels live in the same symplectization $\RR \times M$,
there is no longer a distinguished \emph{main level} or any meaningful
notion of \emph{upper} vs.~\emph{lower} levels; the level structure is simply a
function $L : S \to \{1,\ldots,N\}$ for some $N \in \NN$, and equivalence
of buildings must permit $\RR$-translations within each level.
For these reasons, the SFT compactness theorem in symplectizations has a few
qualitative differences, but is still very much analogous to
Theorem~\ref{thm:SFTcompactness}.

To complete the picture, we should mention one more type of compactness
theorem that appears in \cite{SFTcompactness}, which is colloquially described
as \emph{stretching the neck}.  The geometric idea is as follows: suppose
$(W,\omega)$ is a closed symplectic manifold and $M \subset W$ is a
stable hypersurface that separates $W$ into two pieces $W = W_- \cup_M W_+$,
with an induced stable Hamiltonian structure $\hH = (\omega,\lambda)$ that
orients $M$ as the boundary of~$W_-$.\footnote{The assumption that $M \subset W$ 
separates $W$ is inessential, but makes certain details in this discussion
more convenient.}
A neighborhood of $M$ in $(W,\omega)$ can then 
be identified symplectically with
$$
(\nN_\epsilon,\omega_\epsilon) := \left( (-\epsilon,\epsilon) \times M,
d(r \lambda) + \omega \right)
$$
for sufficiently small $\epsilon > 0$.  The idea now is to replace
$\nN_\epsilon$ with larger collars of the form
$$
\left( (-T,T) \times M , d\left( f(r) \lambda\right) + \omega \right),
$$
with $C^0$-small functions $f$ chosen with $f' > 0$ so that the collar
can be glued in smoothly to replace $(\nN_\epsilon,\omega_\epsilon)$.
This collar looks like a piece of the symplectization of $(M,\hH)$,
thus we are free to choose tame almost complex structures whose restrictions
to the inserted collar belong to~$\jJ(\hH)$.  Symplectic manifolds
constructed in this way are all symplectomorphic, but their almost complex
structures degenerate as one takes $T \to \infty$.  Given a sequence
$T_k \to \infty$ and a corresponding degenerating sequence~$J_k$, 
a sequence $u_k$ of $J_k$-holomorphic
curves with bounded energy converges to yet another form of holomorphic
building, this time involving a bottom level in $\widehat{W}_- :=
W_- \cup_M \left([0,\infty) \times M \right)$ with positive 
punctures approaching orbits in~$M$, some finite number of middle levels that 
live in the symplectization of~$M$, and a top level that lives in
$\widehat{W}_+ := \left((-\infty,0] \times M\right) \cup_M W_+$ with negative
punctures approaching~$M$.  

A very popular example for applications arises from Lagrangian submanifolds
$L \subset W$.  By the Weinstein neighborhood theorem, $L$ always has a
neighborhood $W_-$ symplectomorphic to a neighborhood of the zero-section
in~$T^*L$, so $M := \p W_-$ is a contact-type hypersurface contactomorphic
to the unit cotangent bundle of~$L$.  Stretching the neck then yields
$T^*L$ as the completion of~$W_-$, and $W \setminus L$ as the completion
of $W_+ := W \setminus \mathring{W}_-$.  This construction has often been
used in order to study Lagrangian submanifolds via SFT-type methods,
see e.g.~\cite{SFT}*{Theorem~1.7.5} and 
\cites{Evans:delPezzo,CieliebakMohnke:Audin}.

\psfrag{rho}{$\rho$}
\psfrag{phi}{$\phi$}
\psfrag{theta}{$\theta$}

\chapter{Cylindrical contact homology and the tight $3$-tori}
\label{lec:tight3tori}

\minitoc

\vspace{12pt}

We've now developed enough of the technical machinery of holomorphic curves
to be able to give a rigorous construction of the most basic version of SFT 
and apply it to a problem in contact topology.

\section{Contact structures on $\TT^3$ and Giroux torsion}

\begin{figure}
\includegraphics{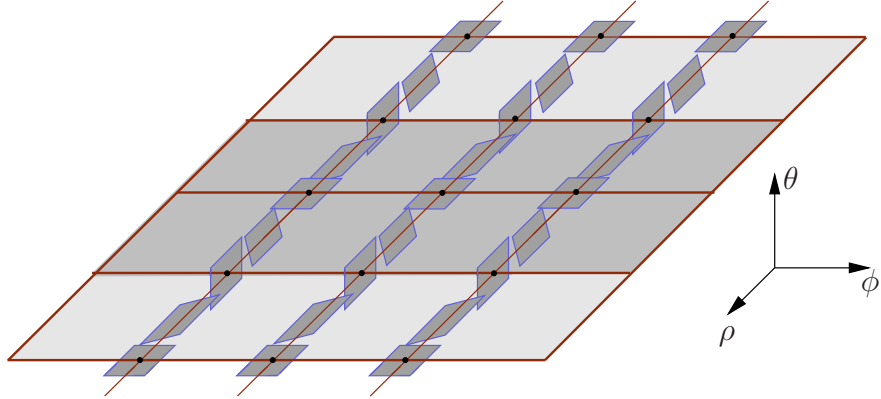}
\caption{\label{fig:T3} The contact structures $\xi_k$ on $\TT^3$ can be
constructed by gluing $k$ copies of the same model $[0,1] \times \TT^2$
to each other cyclically.}
\end{figure}

As a motivating goal in this lecture, we will prove a result about the 
classification of contact structures on $\TT^3 = S^1 \times S^1 \times S^1$.  
Denote the three global coordinates on $\TT^3$ valued in $S^1 = \RR / \ZZ$
by $(\rho,\phi,\theta)$, and
for any $k \in \NN$, consider the contact structure
$$
\xi_k := \ker \alpha_k, \quad \text{ where } \quad \alpha_k := \cos(2\pi k \rho)\, d\theta + \sin(2\pi k \rho)\, d\phi.
$$
It is an easy exercise to verify that these all satisfy the contact condition
$\alpha_k \wedge d\alpha_k > 0$; see Figure~\ref{fig:T3} for a visual
representation.
The following result is originally due to Giroux \cite{Giroux:plusOuMoins}
and Kanda \cite{Kanda:torus}.

\begin{thm}
\label{thm:T3}
For each pair of positive integers $k \ne \ell$, the contact manifolds
$(\TT^3,\xi_k)$ and $(\TT^3,\xi_\ell)$ are not contactomorphic.
\end{thm}

One of the reasons this result is interesting is that it cannot be proved
using any so-called ``classical'' invariants, i.e.~invariants coming from
algebraic topology.  An example of a classical invariant would be the
Euler class of the oriented vector bundle $\xi_k \to \TT^3$, or
anything else that depends only on the isomorphism class of this bundle.
The following observation shows that such invariants will never distinguish
$\xi_k$ from~$\xi_\ell$.

\begin{prop}
\label{prop:homotopy}
For every $k, \ell \in \NN$, $\xi_k$ and $\xi_\ell$ are homotopic through a
smooth family of oriented $2$-plane fields on~$\TT^3$.
\end{prop}
\begin{proof}
In fact, all the $\xi_k$ can be deformed smoothly to $\ker d\rho$, 
via the homotopy
$$
\ker\left[ (1-s) \, \alpha_k + s\, d\rho \right], \qquad s \in [0,1].
$$
\end{proof}
\begin{remark}
One can check in fact that the $1$-form in the homotopy given above is
contact for every $s \in [0,1)$, so Gray's stability theorem implies that
every $\xi_k$ is isotopic to an arbitrarily small perturbation of the
foliation $\ker d\rho$.  In \cite{Giroux:plusOuMoins}, Giroux used this
observation to show that all of them are what we now call \emph{weakly
symplectically fillable}.  If $\ker d\rho$ were also contact, then Gray's
theorem would imply that $\xi_k$ and $\xi_\ell$ are always isotopic.
Thus Theorem~\ref{thm:T3} indicates the
impossibility of modifying a homotopy from $\xi_k$ to $\xi_\ell$ into one
that passes only through contact structures.
\end{remark}

Let us place this discussion in a larger context.
Using the coordinates $(\rho,\phi,\theta)$ on $\RR \times \TT^2$, a pair of
smooth functions $f, g : \RR \to \RR$ gives rise to a contact form
$$
\alpha = f(\rho)\, d\theta + g(\rho)\, d\phi
$$
whenever the function $D(\rho) := f(\rho) g'(\rho) - f'(\rho) g(\rho)$ 
is everywhere positive.  Indeed, we have $\alpha \wedge d\alpha =
D(\rho)\, d\rho \wedge d\phi \wedge d\theta$, and one easily derives a
similar formula for the Reeb vector field,
$$
R_\alpha = \frac{1}{D(\rho)} \left[ g'(\rho) \, \p_\theta - f'(\rho)\, \p_\phi \right].
$$
The condition $D > 0$ means geometrically that the path $(f,g) : \RR \to \RR^2$
winds counterclockwise around the origin with its angular coordinate
strictly increasing.  The simplest special case is the contact form
$$
\alpha_\GT := \cos(2\pi \rho) \, d\theta + \sin(2\pi \rho)\, d\phi,
$$
which matches the formula for $\alpha_1$ on $\TT^3$ given above.
Let $\xi_\GT := \ker \alpha_\GT$ on $\RR \times \TT^2$.

\begin{defn}
\label{defn:GT}
The \defin{Giroux torsion} $\GT(M,\xi) \in \NN \cup \{0,\infty\}$
of a contact $3$-manifold $(M,\xi)$ is the
supremum of the set of positive integers $k$ such that there exists a
contact embedding
$$
\left( [0,k] \times \TT^2, \xi_\GT \right) \hookrightarrow (M,\xi).
$$
We write $\GT(M,\xi) = 0$ if no such embedding exists for any~$k$,
and $\GT(M,\xi) = \infty$ if it exists for all~$k$.
\end{defn}

\begin{example}
The tori $(\TT^3,\xi_k)$ for $k \ge \ZZ$ are contactomorphic to
$(\RR \times \TT^2, \xi_\GT) / k\ZZ$, with $k\ZZ$ acting by translation of
the $\rho$-coordinate.  Thus $\GT(\TT^3,\xi_k) \ge k - 1$.
\end{example}

A $2$-torus $T \subset (M,\xi)$ embedded in a contact $3$-manifold is called
\defin{pre-Lagrangian} if a neighborhood of $T$ in $(M,\xi)$ admits a
contactomorphism to a neighborhood of $\{0\} \times \TT^2$ in
$(\RR \times \TT^2,\xi_\GT)$, identifying $T$ with $\{0\} \times \TT^2$.
The neighborhood in $\RR \times \TT^2$ can be arbitrarily small, thus
the existence of a pre-Lagrangian torus does not imply $\GT(M,\xi) > 0$;
in fact, pre-Lagrangian tori always exist in abundance, e.g.~as boundaries
of neighborhoods of transverse knots (using the contact model provided by
the transverse neighborhood theorem).  But given any pre-Lagrangian torus
$T \subset (M,\xi)$, one can make a local modification of $\xi$ near $T$
to produce a new contact structure (up to isotopy) 
with positive Giroux torsion.  Define $(M',\xi')$ from $(M,\xi)$ by replacing
the small neighborhood $((-\epsilon,\epsilon) \times \TT^2,\xi_\GT)$ with
$((-\epsilon,1 + \epsilon) \times \TT^2,\xi_\GT)$, then identify $M'$
with $M$ by a choice of compactly supported diffeomorphism $(-\epsilon,1+\epsilon)
\to (-\epsilon,\epsilon)$.  There is now an obvious contact embedding of
$([0,1] \times \TT^2,\xi_\GT)$ into $(M,\xi')$, hence $\GT(M,\xi') \ge 1$.
Moreover, one can adapt the proof of Prop.~\ref{prop:homotopy} above to show
that $\xi'$ is homotopic to~$\xi$ through a smooth family of oriented
$2$-plane fields.  The operation changing $\xi$ to $\xi'$ is known as a
\defin{Lutz twist} along~$T$.  In this language, we see that for each
$k \in \NN$, $(\TT^3,\xi_{k+1})$ is obtained from $(\TT^3,\xi_k)$ by
performing a Lutz twist along $\{0\} \times \TT^2$.

The invariant $\GT(M,\xi)$ is easy to define, but hard to compute in general.
The natural guess, 
$$
\GT(\TT^3,\xi_k) = k-1,
$$
turns out to be correct, as was
shown in \cite{Giroux:bifurcations}, so this is one way to prove
Theorem~\ref{thm:T3}, but not the approach we will take.  The following
example shows that one must in any case be careful with such guesses.

\begin{example}
\label{ex:S1S2}
For each $k \in \NN$, define a model of $S^1 \times S^2$ by
$$
S^1 \times S^2 \cong \left( [0,k+ 1/2] \times \TT^2 \right) \big/ \sim
$$
where the equivalence relation identifies
$(\rho,\phi,\theta) \sim (\rho,\phi',\theta)$ for $\rho \in \{0,k+1/2\}$
and every $\theta,\phi,\phi' \in S^1$.  Near $\rho=0$ and $\rho=k+1/2$, this 
means thinking of $(\rho,\phi)$ as polar coordinates, so the two subsets
$\{\rho=0\}$ and $\{\rho=k+1/2\}$ become circles of the form 
$S^1 \times \{\text{const}\}$ embedded in $S^1 \times S^2$.
Since the $\phi$-coordinate is singular at these two circles, the contact
form $\alpha_\GT$ needs to be modified slightly in this region before it
will descend to a smooth contact form on $S^1 \times S^2$: this can be
done by a $C^0$-small modification of the 
form $f(\rho)\, d\theta + g(\rho)\, d\phi$, and the resulting contact
structure is then uniquely determined up to isotopy.  We shall call this
contact manifold
$$
(S^1 \times S^2, \xi_k).
$$
Now observe that for each $k \in \NN$, $(S^1 \times S^2,\xi_{k+1})$ is
obtained from $(S^1 \times S^2,\xi_k)$ by a Lutz twist.  However, both
contact manifolds are also \defin{overtwisted}: recall that a contact
$3$-manifold $(M,\xi)$ is overtwisted whenever it contains an embedded
closed $2$-disk $\dD \subset M$ such that $T(\p\dD) \subset \xi$ but
$T\dD|_{\p\dD} \pitchfork \xi$.  (Exercise: find a disk with this property
in $(S^1 \times S^2,\xi_k)$!)  Eliashberg's flexibility theorem for
overtwisted contact structures \cite{Eliashberg:overtwisted} implies
that whenever $\xi$ and $\xi'$ are two contact structures
on a closed $3$-manifold that are both overtwisted and are homotopic as
oriented $2$-plane fields, they are actually isotopic.  As a consequence,
the contact structures $\xi_k$ on $S^1 \times S^2$ defined above for every
$k \in \NN$ are all isotopic to each other.  As tends to be the case with
most interesting h-principles, the isotopy is very hard to see concretely,
but it must exist.
\end{example}

\begin{exercise}
\label{EX:GTinfty}
Show that if $(M,\xi)$ is a closed overtwisted contact $3$-manifold, then
$\GT(M,\xi) = \infty$.
\end{exercise}

In contrast to the $S^1 \times S^2$ example above, the contact manifolds
$(\TT^3,\xi_k)$ are not overtwisted, they are \defin{tight}---in fact,
the classification of contact structures on $\TT^3$ by 
Giroux \cites{Giroux:plusOuMoins,Giroux:infiniteTendues,Giroux:bifurcations}
and Kanda \cite{Kanda:torus} states that these are \emph{all} of the tight
contact structures on $\TT^3$ up to contactomorphism.
We will use cylindrical contact homology to show that they are not
contactomorphic to each other.  The reader should keep Example~\ref{ex:S1S2}
in mind and try to spot the reason why the same argument cannot work
for $(S^1 \times S^2,\xi_k)$.

\begin{remark}
It has been conjectured that the converse of Exercise~\ref{EX:GTinfty}
might also hold, so every closed tight contact $3$-manifold would have finite
Giroux torsion.  This conjecture is wide open.
\end{remark}

\section{Definition of cylindrical contact homology}

\subsection{Preliminary remarks}

\begin{figure}
\includegraphics{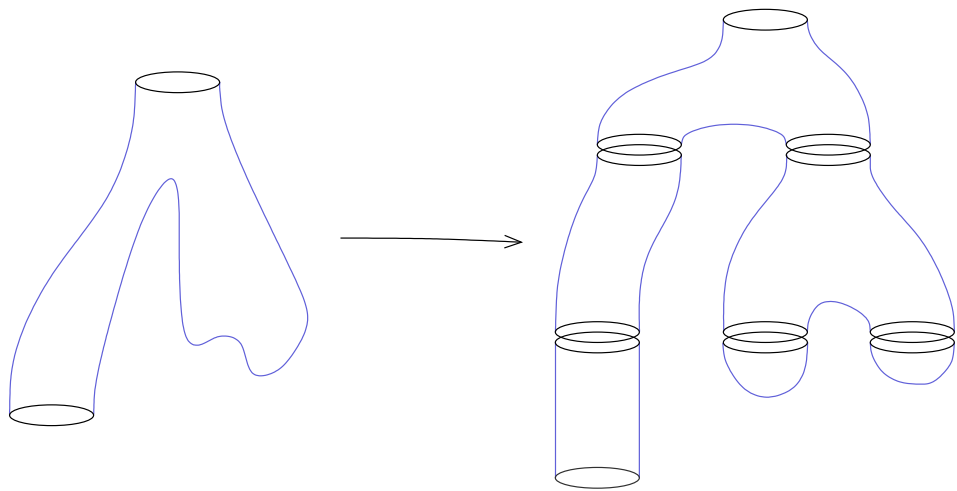}
\caption{\label{fig:cylCH} A family of holomorphic cylinders 
can converge in the SFT topology to buildings that include more
complicated curves than cylinders---this is why cylindrical contact
homology is not well defined for all contact manifolds.}
\end{figure}

Cylindrical contact homology is the natural ``first attempt'' at
using holomorphic curves in symplectizations to define a Floer-type
invariant of contact manifolds~$(M,\xi)$.  The idea is to define a chain complex
generated by Reeb orbits in $M$ and a differential $\p$ that counts holomorphic
cylinders in $\RR \times M$.  We already know some pretty good reasons why 
this idea cannot work in general: in order to prove $\p^2 = 0$, we need to
be able to identify the space of rigid ``broken'' holomorphic cylinders
(these are what is counted by~$\p^2$) with the boundary of the compactified
$1$-dimensional space of index~$2$ cylinders (up to $\RR$-translation).
But this compactified boundary has more than just broken cylinders in it,
see Figure~\ref{fig:cylCH}.  In order to define cylindrical contact
homology, one must therefore restrict to situations in which complicated
pictures like Figure~\ref{fig:cylCH} cannot occur.  The first useful remark
in this direction is that since we are working with a stable Hamiltonian 
structure of the form $(d\alpha,\alpha)$ for a contact form $\alpha$, 
a certain subset of the scenarios allowed by the SFT compactness theorem
can be excluded immediately.  Indeed:

\begin{prop}
\label{prop:onePositive}
If $J \in \jJ(\alpha)$ and $u : (\dot{\Sigma},j) \to (\RR \times M,J)$
is an asymptotically cylindrical $J$-holomorphic curve, then $u$ has
at least one positive puncture.
\end{prop}

Let us give two proofs of this result, since both contain useful ideas.
As preparation for the first proof,
recall the definition of energy for curves in symplectizations of contact
manifolds that we wrote down in Lecture~\ref{lec:intro}:
$$
E(u) := \sup_{f \in \tT} \int_{\dot{\Sigma}} u^*d(e^{f(r)}\, \alpha),
$$
where
$$
\tT := \left\{ f \in C^\infty(\RR,(-1,1))\ \big|\ f' > 0 \right\}.
$$
This formula is not identical to the definition of energy used
in Lecture~\ref{lec:compactness}, but it is equivalent in the sense that any uniform bounds on 
one imply similar uniform bounds on the other.  

\begin{proof}[First proof of Proposition~\ref{prop:onePositive}]
Denote the positive and negative punctures of $u : \dot{\Sigma} \to \RR \times M$
by $\Gamma^+$ and $\Gamma^-$ respectively, and suppose $u$ is asymptotic at
$z \in \Gamma^\pm$ to the orbit $\gamma_z$ with period $T_z > 0$.
Choose any $f \in \tT$ and denote $f_\pm := \lim_{r \to \pm\infty} f(r) \in
[-1,1]$.  
Since $d(e^{f(r)} \, \alpha)$ tames $J \in \jJ(\alpha)$, Stokes' theorem gives
\begin{equation}
\label{eqn:energyStokes}
0 \le E(u) = e^{f_+} \sum_{z \in \Gamma^+} T_z - e^{f_-} \sum_{z \in \Gamma^-} T_z,
\end{equation}
hence $\Gamma^+$ cannot be empty.
\end{proof}
\begin{remark}
The proof via Stokes' theorem works just as well if instead of $\RR \times M$,
$u$ lives in the completion of an exact symplectic cobordism $(W,\omega)$
with concave boundary $(M_-,\xi_- = \ker \alpha_-)$ and 
convex boundary $(M_+,\xi_+ = \ker\alpha_+)$.  Recall that this means
$\p W = -M_- \sqcup M_+$, and $\omega = d\lambda$ for a $1$-form $\lambda$
that restricts to positive contact forms $\lambda|_{TM_\pm} = \alpha_\pm$.
As in Lecture~\ref{lec:intro}, we will write
$$
\jJ(W,\omega,\alpha_+,\alpha_-) \subset \jJ(\widehat{W})
$$
for the space of almost complex structures $J$ on 
$\widehat{W} := \left( (-\infty,0] \times M_-\right) \cup_{M_-} W 
\cup_{M_+} \left( [0,\infty) \times M_+ \right)$ that are compatible with
$\omega$ on $W$ and belong to $\jJ(\alpha_\pm)$ on the cylindrical ends.
The energy of a $J$-holomorphic curve $u : (\dot{\Sigma},j) \to (\widehat{W},J)$
is then
$$
E(u) := \sup_{f \in \tT} \int_{\dot{\Sigma}} u^*d\lambda_f,
$$
where
$\tT := \left\{ f \in C^\infty(\RR,(-1,1))\ |\ \text{$f' > 0$ and
$f(r) = r$ near $r=0$} \right\}$ and
$$
\lambda_f := \begin{cases}
e^{f(r)} \alpha_+ & \text{ on $[0,\infty) \times M_+$},\\
\lambda & \text{ on $W$},\\
e^{f(r)} \alpha_- & \text{ on $(-\infty,0] \times M_-$}.
\end{cases}
$$
The above proof now generalizes verbatim to show that $u$ must always have
a positive puncture.  Notice that in both settings, the argument also gives a 
uniform bound for the energy in terms of the periods of the positive asymptotic
orbits.
\end{remark}
\begin{remark}
We can also prove Prop.~\ref{prop:onePositive} using the fact that
$u^*d\alpha \ge 0$ for any $u : (\dot{\Sigma},j) \to (\RR \times M,J)$
with $J \in \jJ(\alpha)$.  Indeed, Stokes' theorem then gives
\begin{equation}
\label{eqn:contactArea}
0 \le \int_{\dot{\Sigma}} u^*d\alpha = \sum_{z \in \Gamma^+} T_z -
\sum_{z \in \Gamma^-} T_z.
\end{equation}
The quantity $\int_{\dot{\Sigma}} u^*d\alpha$ is sometimes called the
\defin{contact area} of~$u$.  This version of the argument however does not 
easily generalize to arbitrary exact cobordisms.
\end{remark}

The second proof is based on the maximum principle for subharmonic functions.

\begin{prop}
\label{prop:maximum}
Suppose $J \in \jJ(\alpha)$ and $u = (u_\RR,u_M) : (\dot{\Sigma},j) \to
(\RR \times M,J)$ is $J$-holomorphic, where $\dot{\Sigma}$ has no boundary.  
Then $u_\RR : \dot{\Sigma} \to \RR$ has no local maxima.
\end{prop}
\begin{proof}
In any local holomorphic coordinates $(s,t)$ on a region in $\dot{\Sigma}$,
the nonlinear Cauchy-Riemann equation for $u$ is equivalent to the system of
equations
\begin{equation*}
\begin{split}
\p_s u_\RR - \alpha(\p_t u_M) &= 0, \\
\p_t u_\RR + \alpha(\p_s u_M) &= 0, \\
\pi_\xi \, \p_s u_M + J \pi_\xi \, \p_t u_M &= 0,
\end{split}
\end{equation*}
where $\pi_\xi : TM \to \xi$ denotes the projection along the Reeb vector
field.  This gives
\begin{equation*}
\begin{split}
- \Delta u_\RR &= -\p_s^2 u_\RR - \p_t^2 u_\RR =
-\p_s\left[ \alpha(\p_t u_M) \right] + \p_t\left[ \alpha(\p_s u_M) \right] \\
&= -d\alpha(\p_s u_M,\p_t u_M) = -d\alpha(\pi_\xi \p_s u_M, J \pi_\xi \p_s u_M) \le 0
\end{split}
\end{equation*}
since $J|_\xi$ is tamed by $d\alpha|_\xi$, hence
$u_\RR$ is subharmonic.  The result thus follows from the
maximum principle, see e.g.~\cite{Evans}.
\end{proof}
\begin{proof}[Second proof of Proposition~\ref{prop:onePositive}]
If $u = (u_\RR,u_M) : \dot{\Sigma} \to \RR \times M$ has no positive puncture
then $u_\RR : \dot{\Sigma} \to \RR$ is a proper function bounded above,
and therefore has a local maximum, contradicting Proposition~\ref{prop:maximum}.
\end{proof}
\begin{remark}
The proof via the maximum principle does not generalize to arbitrary exact
cobordisms $(W,d\lambda)$, but it does work in \emph{Stein} cobordisms, 
i.e.~if $\lambda_f$ and $J$ are related by $\lambda_f = - dF \circ J$ for some
plurisubharmonic function $F : \widehat{W} \to \RR$, then $F \circ u :
\dot{\Sigma} \to \RR$ is subharmonic (cf.~\cite{CieliebakEliashberg}).
\end{remark}

With these preliminaries understood, the next two exercises reveal one 
natural setting in which breaking of cylinders can be kept under control.
Both exercises are essentially combinatorial.

\begin{exercise}
\label{EX:onePositive}
Suppose $\mathbf{u}$ is a stable $J$-holomorphic building in a completed
symplectic cobordism $\widehat{W}$ with the following properties:
\begin{enumerate}
\item $\mathbf{u}$ has arithmetic genus~$0$ and exactly one positive puncture;
\item every connected component of $\mathbf{u}$ has at least one positive puncture.
\end{enumerate}
Show that $\mathbf{u}$ has no nodes, and all of its connected components have
\emph{exactly} one positive puncture.
\end{exercise}

\begin{exercise}
\label{EX:noPlanes}
Suppose that in addition to the conditions of Exercise~\ref{EX:onePositive},
$\mathbf{u}$ has exactly one negative puncture and
no connected component of $\mathbf{u}$ is a plane.  Show that every level
of $\mathbf{u}$ then consists of a single cylinder with one positive and
one negative end.
\end{exercise}

Exercise~\ref{EX:noPlanes} makes it reasonable to define a 
Floer-type theory counting only cylinders in any setting where planes can
be excluded, for instance because the Reeb vector field has no contractible orbits.
This is not always possible, e.g.~Hofer \cite{Hofer:weinstein} proved that on 
overtwisted contact manifolds, there is \emph{always} a plane (which is why
the Weinstein conjecture holds).  So the invariant we construct will not be
defined in such settings, but it happens to be ideally suited to the
study of~$(\TT^3,\xi_k)$.

\subsection{A compactness result for cylinders}

Fix a closed contact manifold $(M,\xi)$ of dimension $2n-1$ and a primitive
homotopy class of loops $h \in [S^1,M]$.  By \defin{primitive}, we mean
that $h$ is not equal to $N h'$ for any $h' \in [S^1,M]$ and an integer
$N > 1$, and this assumption will be crucial for technical reasons in the
following.\footnote{It is to be expected that cylindrical contact homology
can be defined also for non-primitive homotopy classes, but this would
require more sophisticated methods to address transversality problems.
The assumption that $h$ is primitive allows us to assume that all
holomorphic curves in the discussion are somewhere injective, hence they
are always regular if $J$ is generic.}
Given a contact form $\alpha$ for $\xi$, let 
$$
\pP_h(\alpha)
$$
denote the set
of closed Reeb orbits homotopic to~$h$, where two Reeb orbits are 
identified if they differ only by parametrization.

\begin{defn}
\label{defn:admis}
Given a contact manifold $(M,\xi)$ and a primitive homotopy class
$h \in [S^1,M]$, we will say that a contact form $\alpha$ for $\xi$ is
\defin{$h$-admissible} if:
\begin{enumerate}
\item All orbits in $\pP_h(\alpha)$ are nondegenerate;
\item There are no contractible closed Reeb orbits.
\end{enumerate}
Similarly, we will say that $(M,\xi)$ is \defin{$h$-admissible}
if a contact form with the above properties exists.
\end{defn}
\begin{defn}
\label{defn:regular}
Given $h \in [S^1,M]$ and an $h$-admissible contact form~$\alpha$ on
$(M,\xi)$, we will say that an almost complex structure $J \in \jJ(\alpha)$
is \defin{$h$-regular} if every $J$-holomorphic cylinder in $\RR \times M$
with a positive and a negative end both asymptotic to orbits in
$\pP_h(\alpha)$ is Fredholm regular.
\end{defn}

\begin{prop}
\label{prop:regularity}
If $h \in [S^1,M]$ is a primitive homotopy class of loops and
$\alpha$ is $h$-admissible on $(M,\xi)$, then the space of $h$-regular
almost complex structures is comeager in~$\jJ(\alpha)$.
\end{prop}
\begin{proof}
Since $h$ is primitive, the asymptotic orbits for the relevant holomorphic
cylinders cannot be multiply covered, hence all of these cylinders are
somewhere injective.  The result therefore follows from the standard
transversality results proved in Lecture~\ref{lec:Dragnev} for somewhere injective curves
in symplectizations.
\end{proof}

\begin{figure}
\includegraphics{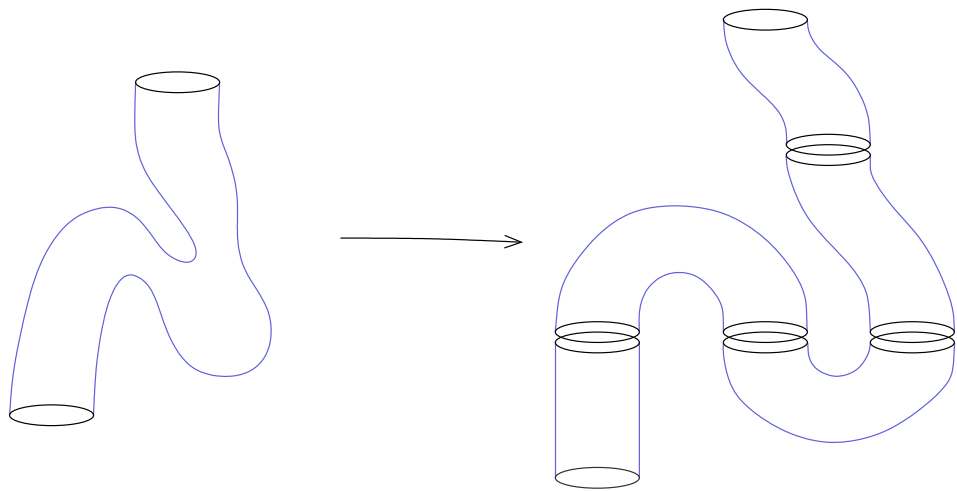}
\caption{\label{fig:curvyCyl} A degenerating sequence of holomorphic cylinders
$u_k : \RR \times S^1 \to \RR \times M$ cannot have a limiting level with
a puncture of the ``wrong'' sign unless $u_k$ violates the maximum principle
for large~$k$.}
\end{figure}

\begin{prop}
\label{prop:cylinders}
Given an $h$-admissible contact form $\alpha$, an $h$-regular almost complex
structure $J \in \jJ(\alpha)$ and an orbit
$\gamma \in \pP_h(\alpha)$, suppose $u_k$ is a sequence
of $J$-holomorphic cylinders in $\RR \times M$ with one positive puncture
at~$\gamma$ and one negative puncture.  Then $u_k$
has a subsequence convergent in the SFT topology to a broken $J$-holomorphic
cylinder, i.e.~a stable building $\mathbf{u}_\infty$ whose levels 
$u_\infty^1,\ldots,u_\infty^{N_+}$ are each cylinders with one positive and 
one negative puncture.  Moreover, each level satisfies
$\ind(u_\infty^N) \ge 1$, thus for large~$k$ in the convergent subsequence,
$$
\ind(u_k) = \sum_{N=1}^{N_+} \ind(u_\infty^N) \ge N_+.
$$
\end{prop}
\begin{proof}
Let's start with some bad news: the standard SFT compactness theorem is not 
applicable in this situation, because we have not assumed that 
$\alpha$ is nondegenerate, nor even Morse Bott---there is no assumption at all
about Reeb orbits in homotopy classes other than $h$ and~$0$.  This fairly
loose set of hypotheses is very convenient in applications, as nondegeneracy
of a contact form is generally a quite difficult condition to check.  The
price we pay is that we will have to prove compactness manually instead
of applying the big theorem (see Remark~\ref{remark:nondeg}).
Fortunately, it is not that hard: the crucial point is that in the situation at
hand, there can be no bubbling at all.

Indeed, we claim that the given sequence 
$u_k : (\RR \times S^1,i) \to (\RR \times M,J)$ must satisfy a uniform bound
$$
|d u_k| \le C
$$
with respect to any translation-invariant Riemannian metrics on
$\RR \times S^1$ and $\RR \times M$.  To see this, note first that since
all the $u_k$ have the same positive asymptotic orbit~$\gamma$, their
energies are uniformly bounded via \eqref{eqn:energyStokes}.  Thus if
$|d u_k(z_k)| \to \infty$ for some sequence $z_k \in \RR \times S^1$,
we can perform the usual rescaling trick from Lecture~\ref{lec:compactness} and deduce the
existence of a nonconstant finite-energy plane $v_\infty : \CC \to \RR \times M$.
Its singularity
at $\infty$ cannot be removable since this would produce a nonconstant
$J$-holomorphic sphere, violating Proposition~\ref{prop:onePositive}.
It follows that $v_\infty$ is asymptotic to a Reeb orbit at~$\infty$, but 
this is also impossible since $\alpha$ does not admit any contractible orbits, 
and the claim is thus proved.

Suppose now that $\gamma$ has period $T_+ > 0$, and observe
that by nondegeneracy, the set
$$
\pP_h(\alpha,T_+) := \left\{ \gamma \in \pP_h(\alpha)\ \big|\  
\text{$\gamma$ has period at most~$T_+$} \right\}
$$
is finite.  Let 
$$
\aA_h(\alpha), \aA_h(\alpha,T_+) \subset (0,\infty)
$$
denote the set of all
periods of orbits in $\pP_h(\alpha)$ and $\pP_h(\alpha,T_+)$ respectively.
By \eqref{eqn:contactArea}, the negative
asymptotic orbit of each $u_k$ is in $\pP_h(\alpha,T_+)$, so we can take a 
subsequence and assume that these are all the same orbit; call it
$\gamma_- \in \pP_h(\alpha,T_+)$ and its period
$T_- \in \aA_h(\alpha,T_+)$.  If $T_- = T_+$ then $u_k^*d\alpha \equiv 0$ for
all~$k$, implying that all $u_k$ are the trivial cylinder over $\gamma$ and
thus trivially converge.  Assume therefore $T_- < T_+$.  Then since
$u_k^*d\alpha \ge 0$, Stokes' theorem implies that for each~$k$, the function
$$
\RR \to \RR : s \mapsto \int_{S^1} u_k(s,\cdot)^*\alpha
$$
is increasing and is a surjective map onto $(T_-,T_+)$.
The uniform bound on the derivatives implies that for any sequences
$s_k, r_k \in \RR$ with $u_k(s_k,0) \in \{r_k\} \times M$, the 
sequence\footnote{Recall from Lecture~\ref{lec:compactness} that we denote the $\RR$-translation
action on $\RR \times M$ by $\tau_c(r,x) := (r + c,x)$.}
$$
v_k : \RR \times S^1 \to \RR \times M : (s,t) \mapsto
\tau_{-r_k} \circ u_k(s + s_k,t)
$$
has a subsequence convergent in $C^\infty_\loc(\RR \times S^1)$ to some
finite-energy $J$-holomorphic cylinder 
$$
v_\infty : \RR \times S^1 \to \RR \times M,
$$
which necessarily satisfies
$$
\int_{S^1} v_\infty(s,\cdot)^*\alpha = \lim_{k \to \infty}
\int_{S^1} u_k(s + s_k,\cdot)^*\alpha \in [T_-,T_+]
$$
for every $s \in \RR$.  This proves that $v_\infty$ is nonconstant,
with a positive puncture at $s=\infty$ and negative puncture at $s=-\infty$,
and both of its asymptotic orbits are in~$\pP_h(\alpha,T_+)$.\footnote{For
an alternative argument that $v_\infty$ must have a positive puncture at
$s=\infty$ and negative at $s=-\infty$, see Figure~\ref{fig:curvyCyl}.}
If $v_\infty$ is not a trivial cylinder, then it therefore satisfies
$$
\int_{\RR \times S^1} v_\infty^*d\alpha \ge \delta,
$$
where $\delta$ is any positive number less than the smallest distance
between neighboring elements of~$\aA_h(\alpha,T_+)$.

Let us call a sequence $s_k \in \RR$ \emph{nontrivial} whenever the
limiting cylinder $v_\infty$ obtained by the above procedure is not
a trivial cylinder, and call two such sequences $s_k$ and
$s_k'$ \emph{compatible} if $s_k - s_k'$ is not bounded.
We claim now that if $s_k^1,\ldots,s_k^m$
is a collection of nontrivial sequences that are all compatible with
each other, then
$$
m < \frac{2(T_+ - T_-)}{\delta}.
$$
Indeed, we can assume after ordering our collection appropriately and
restricting to a subsequence that $s_k^{N+1} - s_k^N \to \infty$ for each
$N=1,\ldots,m-1$, and let $v_\infty^N : \RR \times S^1 \to \RR \times M$ 
denote the limits of the 
corresponding convergent subsequences.  Then we can find $R > 0$ such that
$$
\int_{[-R,R] \times S^1} (v_\infty^N)^*d\alpha > \frac{\delta}{2}
$$
and thus
$$
\int_{[s_k^N - R, s_k^N + R] \times S^1} u_k^*d\alpha > \frac{\delta}{2}
$$
for each $N=1,\ldots,m$ for sufficiently large~$k$.  But these domains are
also all disjoint for sufficiently large~$k$, implying
$$
T_+ - T_- = \int_{\RR \times S^1} u_k^*d\alpha \ge 
\sum_{N=1}^m \int_{[s_k^N - R, s_k^N + R] \times S^1} u_k^*d\alpha > 
\frac{\delta m}{2}.
$$

We've shown that there exists a maximal collection of nontrivial sequences
$s_k^1,\ldots,s_k^{N_+} \in \RR$ satisfying $s_k^{N+1} - s_k^N \to \infty$
for each~$N$, such that if $u_k(s_k^N,0) \in \{r_k^N\} \times M$, then
after restricting to a subsequence, the cylinders
$$
v_k^N(s,t) := \tau_{-r_k^N} \circ u_k(s + s_k^N,t)
$$
each converge in $C^\infty_\loc(\RR \times S^1)$ as $k \to \infty$ to a 
nontrivial $J$-holomorphic cylinder 
$u_\infty^N : \RR \times S^1 \to \RR \times M$.  Let $\gamma_N^\pm$ denote
the asymptotic orbit of $u_\infty^N$ at $s=\pm\infty$.  We claim,
$$
\gamma_N^+ = \gamma_{N+1}^- \quad \text{ for each $N=1,\ldots,N_+ - 1$}.
$$
If $\gamma_N^+ \ne \gamma_{N+1}^-$ for some~$N$, choose a neighborhood
$\uU \subset M$ of the image of $\gamma_N^+$ that does not intersect
any other orbit in $\pP_h(\alpha,T_+)$.
Then since each $u_k$ is continuous, there must exist
a sequence $s_k' \in \RR$ with
$$
s_k' - s_k^N \to \infty \quad \text{ and } \quad
s_k^{N+1} - s_k' \to \infty
$$
such that $u_k(s_k',0)$ lies in $\uU$ for all~$k$ but stays a positive
distance away from the image of~$\gamma_N^+$.  A subsequence of
$(s,t) \mapsto u_k(s + s_k',t)$ then converges after suitable $\RR$-translations
to a cylinder $u_\infty' : \RR \times S^1 \to \RR \times M$ that cannot be trivial since
$u_\infty'(0,0)$ is not contained in any orbit in $\pP_h(\alpha,T_+)$.
This contradicts the assumption that our collection
$s_k^1,\ldots,s_k^{N_+}$ is maximal.  A similar argument shows
$$
\gamma_1^- = \gamma^- \quad \text{ and } \quad
\gamma_{N_+}^+ = \gamma,
$$
so the curves $u_\infty^1,\ldots,u_\infty^{N_+}$ form the levels of a
stable holomorphic building~$\mathbf{u}_\infty$.  A similar argument by
contradiction also shows that the sequence $u_k$ must converge in
the SFT topology to~$\mathbf{u}_\infty$.

Finally, note that since all the breaking orbits in $\mathbf{u}_\infty$
are homotopic to~$h$ and $J$ is $h$-regular,
the levels $u_\infty^N$ are Fredholm regular.  Since all of them also 
come in $1$-parameter
families of distinct curves related by the $\RR$-action, this implies
$\ind(u_\infty^N) \ge 1$ for each $N=1,\ldots,N_+$.
\end{proof}

\begin{remark}
\label{remark:nondeg}
Nondegeneracy or Morse-Bott conditions are required for several reasons in the
proof of SFT compactness, and indeed, the theorem is not true in general
without some such assumption.  One can see this by considering what happens
to a sequence $u_k$ of $J_k$-holomorphic curves where $J_k \to J_\infty$ 
is compatible with a sequence of nondegenerate contact forms $\alpha_k$ 
converging to one that is only Morse-Bott.  A compactness theorem for this
scenario is proved in \cite{Bourgeois:thesis}, but it requires more general
limiting objects than holomorphic buildings.  On the other hand, it is useful
for certain kinds of applications to know when one can do without
nondegeneracy assumptions and prove compactness anyway.
There are two main
advantages to knowing that all Reeb orbits are nondegenerate or belong to
Morse-Bott families:
\begin{enumerate}
\item It implies that the set of all periods of closed orbits, the
so-called \defin{action spectrum} of~$\alpha$, is a \emph{discrete}
subset of $(0,\infty)$; in fact, for any $T > 0$, the set of all periods
less than~$T$ is finite.  Using the relations \eqref{eqn:energyStokes}
and \eqref{eqn:contactArea}, this implies lower bounds on the possible
energies of limiting components and thus helps show that only finitely
many such components can arise.
\item Curves asymptotic to nondegenerate or Morse-Bott orbits also satisfy
exponential convergence estimates proved in 
\cites{HWZ:props1,HWZ:FIMpreprint,HWZ:props4,Bourgeois:thesis}, and
similar asymptotic estimates yield a result about
``long cylinders with small area'' (see \cite{HWZ:cylinders} and \cite{SFTcompactness}*{Prop.~5.7}) 
which helps in proving that neighboring levels connect to each other along 
breaking orbits.
\end{enumerate}
Our situation in Proposition~\ref{prop:cylinders} was simple enough to avoid
using the ``long cylinder'' lemma, and we did use the discreteness of the
action spectrum, but only needed it for orbits in $\pP_h(\alpha)$ since we
were able to rule out bubbling in the first step.  An alternative would
have been to assume that all orbits (in all homotopy classes) with period
up to the period of $\gamma$ are nondegenerate: then \eqref{eqn:contactArea}
implies that degenerate orbits never play any role in the main arguments
of \cite{SFTcompactness}, so the big theorem becomes safe to use.
\end{remark}

\subsection{The chain complex}

We now define a $\ZZ_2$-graded chain complex with coefficients in $\ZZ_2$
and generators $\langle \gamma \rangle$ for $\gamma \in \pP_h(\alpha)$, i.e.
$$
CC_*^h(M,\alpha) := \bigoplus_{\gamma \in \pP_h(\alpha)} \ZZ_2.
$$
The degree of each generator $\langle \gamma \rangle \in CC_*^h(M,\alpha)$ 
is defined by
$$
|\langle \gamma \rangle| = n - 3 + \muCZ(\gamma) \in \ZZ_2,
$$
where $\muCZ(\gamma) \in \ZZ_2$ denotes the parity of the Conley-Zehnder
index with respect to any choice of trivialization.  The choice to write
$n-3$ in front of this is a convention that will make no difference at
all in this lecture, but it is consistent with a $\ZZ$-grading that we will
be able to define under suitable assumptions in Lecture~\ref{lec:H}.
To define the differential on $CC_*^h(M,\alpha)$, choose an $h$-regular almost
complex structure $J \in \jJ(\alpha)$.  Given Reeb
orbits $\gamma^+,\gamma^- \in \pP_h(\alpha)$ and a number $I \in \ZZ$, let
$$
\mM^I(J,\gamma^+,\gamma^-)
$$
denote the space of all $\RR$-equivalence classes of index~$I$ holomorphic 
cylinders in $(\RR \times M,J)$ 
asymptotic to $\gamma^\pm$ at $\pm\infty$, i.e.~the union of all components 
$\mM_{0,0}(J,A,\gamma^+,\gamma^-) / \RR$ for which
$\virdim \mM_{0,0}(J,A,\gamma^+,\gamma^-) = I$.
Since $J$ is $h$-regular, all the curves in
$\mM^I(J,\gamma^+,\gamma^-)$ are Fredholm regular,
so if $I \ge 1$, $\mM^I(J,\gamma^+,\gamma^-)$ is a smooth manifold with
$$
\dim \mM^I(J,\gamma^+,\gamma^-) = I - 1.
$$
Similarly, $\mM^0(J,\gamma^+,\gamma^-)$ only contains trivial cylinders
and is thus empty unless $\gamma^+ = \gamma^-$, and
$\mM^I(J,\gamma^+,\gamma^-)$ is always empty for $I < 0$.
In particular, $\mM^1(J,\gamma^+,\gamma^-)$ is a discrete set whenever
$\gamma^+ \ne \gamma^-$,
and by Proposition~\ref{prop:cylinders}, it is also compact, hence finite.
We can therefore define
$$
\p\langle \gamma \rangle = \sum_{\gamma' \in \pP_h(\alpha)}
\#_2 \mM^1(J,\gamma,\gamma') \langle \gamma' \rangle,
$$
where for any set $X$, we denote by $\#_2 X$ the cardinality of $X$
modulo~$2$.  The operator $\p$ has odd degree with respect to the grading
since every index~$1$ holomorphic cylinder $u$ with asymptotic orbits
$\gamma^+$ and $\gamma^-$ satisfies
$$
\ind(u) = 1 = \muCZ^{\tau}(\gamma^+) - \muCZ^{\tau}(\gamma^-)
$$
for suitable choices of the trivialization~$\tau$.

\subsection{The homology}

Following the standard Floer theoretic prescription, the relation $\p^2 = 0$
should arise by viewing the compactification 
$\overline{\mM}^2(J,\gamma^+,\gamma^-)$ for each $\gamma^+,\gamma^- \in \pP_h(\alpha)$
as a compact $1$-manifold whose boundary is identified with the set of
rigid broken cylinders, as these are what is counted by~$\p^2$.
Here $\overline{\mM}^2(J,\gamma^+,\gamma^-)$ is defined as the closure of
$\mM^2(J,\gamma^+,\gamma^-)$ in the space of all $J$-holomorphic buildings
in $\RR \times M$ modulo $\RR$-translation.
Proposition~\ref{prop:cylinders} gives a natural inclusion
$$
\overline{\mM}^2(J,\gamma^+,\gamma^-) \setminus
\mM^2(J,\gamma^+,\gamma^-) \subset \bigsqcup_{\gamma_0 \in \pP_h(\alpha)}
\mM^1(J,\gamma^+,\gamma_0) \times \mM^1(J,\gamma_0,\gamma^-).
$$
We therefore need an inclusion in the other direction, and for this we need to
say a word about gluing.  We have not had time to discuss gluing in earnest
in these notes, and we will not do so now either, but the basic idea should
be familiar from Floer homology: given $u_+ \in \mM^1(J,\gamma^+,\gamma_0)$ and
$u_- \in \mM^1(J,\gamma_0,\gamma^-)$, one would like to show that there
exists a unique (up to $\RR$-translation) one-parameter family 
$\{ u_R \in \mM^2(J,\gamma^+,\gamma^-) \}_{R \in [R_0,\infty)}$ such
that $u_R$ converges as $R \to \infty$ to the building $\mathbf{u}_\infty$ 
with bottom level $u_-$ and top level~$u_+$.  One starts by constructing a
family of \emph{preglued} maps
$$
\tilde{u}_R : \RR \times S^1 \to \RR \times M,
$$
meaning a smooth family of maps
which converge in the SFT topology as $R \to \infty$ 
to $\mathbf{u}_\infty$ but are only
\emph{approximately} $J$-holomorphic.  More precisely, fix parametrizations
of $u_-$ and $u_+$ and a parametrization of the orbit
$\gamma_0 : \RR / T\ZZ \to M$ such that 
\begin{equation*}
\begin{split}
u_+(s,t) = \exp_{(Ts,\gamma_0(Tt))} h_+(s,t) &\quad \text{ for $s \ll 0$},\\
u_-(s,t) = \exp_{(Ts,\gamma_0(Tt))} h_-(s,t) &\quad \text{ for $s \gg 0$},
\end{split}
\end{equation*}
where $h_\pm$ are vector fields along the trivial cylinder satisfying
$\lim_{s \to \mp\infty} h_\pm(s,t) = 0$.  By interpolating between suitable
reparametrizations of $h_+$ and $h_-$, one can now define $\tilde{u}_R$
such that
\begin{equation*}
\begin{split}
\tilde{u}_R(s,t) = \tau_{2RT} \circ u_+(s - 2 R,t) & \quad \text{ for $s \ge R$},\\
\tilde{u}_R(s,t) \approx (Ts,\gamma_0(Tt)) & \quad \text{ for $s \in [-R,R]$,}\\
\tilde{u}_R(s,t) = \tau_{-2RT} \circ u_-(s + 2 R,t) & \quad \text{ for $s \le -R$},\\
\dbar_J \tilde{u}_R \to 0 & \quad \text{ as $R \to \infty$}.
\end{split}
\end{equation*}
Given regularity of $u_+$ and $u_-$, one can now use a quantitative version
of the implicit function theorem (cf.~\cite{McDuffSalamon:Jhol}*{\S 3.5})
to show that a distinguished $J$-holomorphic cylinder $u_R$ close to
$\tilde{u}_R$ exists for all $R$ sufficiently large.  For a more detailed
synopsis of the analysis involved, see \cite{Nelson:thesis}*{Chapter~7}, and
\cite{AudinDamian}*{Chapters~9 and~13} for the analogous story in
Floer homology.  The result is:

\begin{prop}
For an $h$-admissible $\alpha$, an $h$-regular $J \in \jJ(\alpha)$ and any two 
orbits $\gamma^+,\gamma^- \in \pP_h(\alpha)$, the space
$\overline{\mM}^2(J,\gamma^+,\gamma^-)$ admits the structure of a compact
$1$-dimensional manifold with boundary, where its boundary points can be
identified naturally with
$\bigsqcup_{\gamma_0 \in \pP_h(\alpha)}
\mM^1(J,\gamma^+,\gamma_0) \times \mM^1(J,\gamma_0,\gamma^-)$.
\qed
\end{prop}

\begin{cor}
The homomorphism $\p : CC_*^h(M,\alpha) \to CC_{*-1}^h(M,\alpha)$ satisfies $\p^2 = 0$.
\qed
\end{cor}

We shall denote the homology of this chain complex by
$$
HC_*^h(M,\alpha,J) := H_*\big( CC_*^h(M,\alpha) , \p \big).
$$
The goal of the rest of this section is to prove that up to natural
isomorphisms, $HC_*^h(M,\alpha,J)$ depends on $(M,\xi)$ and~$h$ but not
on the auxiliary data $\alpha$ and~$J$.

\subsection{Chain maps}

For any constant $c > 0$, there is an obvious bijection between the generators
of $CC_*^h(M,\alpha)$ and $CC_*^h(M,c\alpha)$, as the rescaling changes
periods of orbits but not the set of closed orbits itself.
Moreover, if $J \in \jJ(\alpha)$ and $J_c \in \jJ(c\alpha)$ are defined to
match on~$\xi$, then there is a biholomorphic diffeomorphism
$$
(\RR \times M,J) \to (\RR \times M,J_c) : (r,x) \mapsto (cr,x),
$$
thus giving a bijective correspondence between the moduli spaces of
$J$-holomorphic and $J_c$-holomorphic curves.  It follows that our
bijection of chain complexes is also a chain map and therefore defines
a canonical isomorphism
\begin{equation}
\label{eqn:canonIso}
HC_*^h(M,\alpha,J) = HC_*^h(M,c\alpha,J_c).
\end{equation}

Next suppose $\alpha_-$ and $\alpha_+$ are two distinct contact forms
for~$\xi$, hence 
$$
\alpha_\pm = e^{f_\pm} \alpha
$$
for some fixed contact
form $\alpha$ and a pair of smooth functions $f_\pm : M \to \RR$.
After rescaling $\alpha_+$ by a constant, we are free to assume
$f_+ > f_-$ everywhere.  Fix $h$-regular almost complex structures 
$J_\pm \in \jJ(\alpha_\pm)$ and let 
$$
\p_\pm : CC_*^h(M,\alpha_\pm) \to CC_{*-1}^h(M,\alpha_\pm)
$$
denote the resulting differentials on the two chain complexes.
The region
$$
W := \left\{ (r,x) \in \RR \times M\ \big|\ f_-(x) \le r \le f_+(x) \right\}
$$
now defines an exact symplectic cobordism from $(M,\xi)$ to itself: more
precisely, setting 
$$
M_\pm := \left\{ (f_\pm(x),x) \in W\ \big|\ x \in M \right\}
$$
gives $\p W = -M_- \sqcup M_+$, and the 
Liouville form $\lambda := e^r \alpha$ satisfies $\lambda|_{TM_\pm} =
\alpha_\pm$.  Choose a
generic $d\lambda$-compatible almost complex structure $J$ on the
completion $\widehat{W}$ that restricts to $J_\pm$ on the cylindrical ends.
Now given $\gamma^+ \in \pP_h(\alpha_+)$ and $\gamma^- \in \pP_h(\alpha_-)$
and a number $I \in \ZZ$, we shall denote by
$$
\mM^I(J,\gamma^+,\gamma^-)
$$
the union of all components $\mM_{0,0}(J,A,\gamma^+,\gamma^-)$ that have
virtual dimension~$I$.  Note that we are not dividing by any $\RR$-action
here since $J$ need not be $\RR$-invariant.
Since $\gamma^\pm$ are still guaranteed to be
simply covered, curves in $\mM^I(J,\gamma^+,\gamma^-)$ are again always
somewhere injective and therefore regular, hence
$\mM^I(J,\gamma^+,\gamma^-)$ is a smooth manifold with
$$
\dim \mM^I(J,\gamma^+,\gamma^-) = I
$$
if $I \ge 0$, and $\mM^I(J,\gamma^+,\gamma^-) = \emptyset$ for $I < 0$.
The compactification $\overline{\mM}^I(J,\gamma^+,\gamma^-)$ is described
via the following straightforward generalization of
Proposition~\ref{prop:cylinders}:

\begin{prop}
\label{prop:cylinders2}
For $J$ as described above, suppose $u_k$ is a sequence
of $J$-holomorphic cylinders in $\widehat{W}$ with one positive puncture
at an orbit $\gamma \in \pP_h(\alpha_+)$ and one negative puncture.  Then $u_k$
has a subsequence convergent in the SFT topology to a broken $J$-holomorphic
cylinder, i.e.~a stable building $\mathbf{u}_\infty$ whose levels 
$u_\infty^N$ for $N = -N_-,\ldots,-1,0,1,\ldots,N_+$
are each cylinders with one positive and 
one negative puncture, living in $\RR \times M^\pm$ for $\pm N > 0$ and
$\widehat{W}$ for $N=0$.  Moreover, the levels satisfy
$\ind(u_\infty^0) \ge 0$ and $\ind(u_\infty^N) \ge 1$ for $N \ne 0$,
thus for large~$k$ in the convergent subsequence,
$$
\ind(u_k) = \sum_{N=-N_-}^{N_+} \ind(u_\infty^N) \ge N_- + N_+.
$$
\qed
\end{prop}

It follows that the set $\mM^0(J,\gamma^+,\gamma^-)$ is always finite,
and we use this to define a map
$$
\Phi_J : CC_*^h(M,\alpha_+) \to CC_*^h(M,\alpha_-) : \langle \gamma \rangle 
\mapsto \sum_{\gamma' \in \pP_h(\alpha_-)} \#_2 \mM^0(J,\gamma,\gamma')
\langle \gamma' \rangle.
$$
This map preserves degrees since it counts index~$0$ curves, and
we claim that it is a chain map:
$$
\Phi_J \circ \p_+ = \p_- \circ \Phi_J.
$$
This follows from the fact that by Proposition~\ref{prop:cylinders2} (in
conjunction with a corresponding gluing theorem), $\overline{\mM}^1(J,\gamma^+,\gamma^-)$
is a compact $1$-manifold whose boundary consists of two types of broken
cylinders, depending whether the index~$1$ curve appears in an upper or
lower level:
\begin{equation*}
\begin{split}
\p \overline{\mM}^1(J,\gamma^+,\gamma^-) &=
\bigsqcup_{\gamma_0 \in \pP_h(\alpha_+)} \left( \mM^1(J_+,\gamma^+,\gamma_0) \times
\mM^0(J,\gamma_0,\gamma^-)\right)  \\
& \quad \cup
\bigsqcup_{\gamma_0 \in \pP_h(\alpha_-)} \left( \mM^0(J,\gamma^+,\gamma_0)
\times \mM^1(J_-,\gamma_0,\gamma^-) \right).
\end{split}
\end{equation*}
Counting broken cylinders of the first type gives the coefficient in front of
$\langle \gamma^- \rangle$ in $\Phi_J \circ \p_+(\langle \gamma^+ \rangle)$,
and the second type gives $\p_- \circ \Phi_J(\langle \gamma^+ \rangle)$.

It follows that $\Phi_J$ descends to a homomorphism
\begin{equation}
\label{eqn:chainMap}
\Phi_J : HC_*^h(M,\alpha_+,J_+) \to HC_*^h(M,\alpha_-,J_-).
\end{equation}

\subsection{Chain homotopies}

We claim that the map $\Phi_J$ in \eqref{eqn:chainMap} does not depend
on~$J$.  To see this, suppose $J_0$ and $J_1$ are two generic choices
of compatible almost complex structures on $\widehat{W}$ that both match
$J_\pm$ on the cylindrical ends.  The space of almost complex structures with
these properties is contractible, so we can find a smooth path
$$
\{ J_s \}_{s \in [0,1]}
$$
connecting them.  For $I \in \ZZ$, consider the parametric moduli space
$$ 
\mM^I(\{J_s\}, \gamma^+,\gamma^-) := \left\{ (s,u)\ \big|\ 
s \in [0,1],\ u \in \mM^I(J_s,\gamma^+,\gamma^-) \right\}.
$$
As we observed in Remark~\ref{remark:parametric}, a generic choice of the homotopy
$\{J_s\}$ makes $\mM^I(\{J_s\})$ a smooth manifold with
$$
\dim \mM^I(\{J_s\},\gamma^+,\gamma^-) = I + 1
$$
whenever $I \ge -1$, and $\mM^I(\{J_s\},\gamma^+,\gamma^-) = \emptyset$
when $I < -1$.  Adapting Proposition~\ref{prop:cylinders2} to allow for a
converging sequence of almost complex structures, it implies that
$\mM^{-1}(\{J_s\},\gamma^+,\gamma^-)$ is compact and thus finite, so we
can use it to define a homomorphism of odd degree by
$$
H : CC_*^h(M,\alpha_+) \to CC_{*+1}^h(M,\alpha_-) : 
\langle \gamma \rangle \mapsto \sum_{\gamma' \in \pP_h(\alpha_-)}
\# \mM^{-1}(\{J_s\},\gamma,\gamma') \langle \gamma' \rangle.
$$
We claim that this is a chain homotopy between $\Phi_{J_0}$ and $\Phi_{J_1}$,
i.e.
$$
\Phi_{J_1} - \Phi_{J_0} = \p_- \circ H + H \circ \p_+.
$$
This follows by looking at the boundary of the compactified
$1$-dimensional space
$\overline{\mM}^0(\{J_s\},\gamma^+,\gamma^-)$, which consists of four
types of objects:
\begin{enumerate}
\item Pairs $(0,u)$ with $u \in \mM^0(J_0,\gamma^+,\gamma^-)$, which are
counted by~$\Phi_{J_0}$.
\item Pairs $(1,u)$ with $u \in \mM^0(J_1,\gamma^+,\gamma^-)$, which are
counted by~$\Phi_{J_1}$.
\item Pairs $(s,\mathbf{u})$ with $\mathbf{u}$ a broken cylinder with
upper level $u_+ \in \mM^1(J_+,\gamma^+,\gamma_0)$ and main level
$u_0 \in \mM^{-1}(J_s,\gamma_0,\gamma^-)$ for some $s \in (0,1)$;
these are counted by $H \circ \p_+$.
\item Pairs $(s,\mathbf{u})$ with $\mathbf{u}$ a broken cylinder with
lower level $u_- \in \mM^1(J_-,\gamma_0,\gamma^-)$ and main level
$u_0 \in \mM^{-1}(J_s,\gamma^+,\gamma_0)$ for some $s \in (0,1)$;
these are counted by $\p_- \circ H$.
\end{enumerate}
The sum $\Phi_{J_0} + \Phi_{J_1} + \p_- \circ H + H \circ \p_+$ therefore
counts (modulo~$2$) the boundary points of a compact $1$-manifold, so it
vanishes.

Since the action of $\Phi_J$ on homology no longer depends on $J$, we will 
denote it from now on by
$$
\Phi : HC_*^h(M,\alpha_+,J_+) \to HC_*^h(M,\alpha_-,J_-).
$$
It is well defined for any pair of $h$-admissible contact forms $\alpha_\pm$
and $h$-regular $J_\pm \in \jJ(\alpha_\pm)$ since one can first rescale
$\alpha_+$ to assume $\alpha_\pm = e^{f_\pm} \alpha$ with $f_+ > f_-$,
using the canonical isomorphism \eqref{eqn:canonIso}.

\subsection{Proof of invariance}

We claim that for any $h$-admissible $\alpha$ and $h$-regular $J \in \jJ(\alpha)$,
the cobordism map
$$
\Phi : HC_*^h(M,\alpha,J) \to HC_*^h(M,\alpha,J)
$$
is the identity.  Indeed, the literal meaning of this statement is that
for any $c > 1$, the composition of the canonical isomorphism
\eqref{eqn:canonIso} with the map
$$
\Phi : HC_*^h(M,c\alpha,J_c) \to HC_*^h(M,\alpha,J)
$$
defined by counting index~$0$ cylinders in a trivial cobordism from
$(M,\alpha,J)$ to $(M,c\alpha,J_c)$ is the identity.  Writing $c = e^a$
for $a > 0$, the Liouville cobordism in question is simply
$$
(W,d\lambda) = ([0,a] \times M, d(e^r \alpha)),
$$
and one can choose a compatible almost complex structure on this which matches
$J$ and $J_c$ on~$\xi$ while taking $\p_r$ to $g(r) R_\alpha$ for a
suitable function $g$ with $g(r) = 1$ near $r=0$ and $g(r) = 1/c$ near $r=a$.
The resulting almost complex manifold is biholomorphically diffeomorphic
to the usual symplectization $(\RR \times M,J)$, so our count of index~$0$
cylinders is equivalent to the count of such cylinders in
$(\RR \times M,J)$.  The latter are simply the trivial cylinders, all of which
are Fredholm regular, so counting these defines the identity map on the
chain complex.

Finally, we need to show that for any three $h$-admissible pairs
$(\alpha_i,J_i)$ with $i=0,1,2$, the cobordism maps
$\Phi_{ij} : HC_*^h(M,\alpha_j,J_j) \to HC_*^h(M,\alpha_i,J_i)$ satisfy
\begin{equation}
\label{eqn:composition}
\Phi_{21} \circ \Phi_{10} = \Phi_{20}.
\end{equation}
We will only sketch this part: the idea is to use a stretching construction.
After rescaling, suppose without loss of generality that
$\alpha_i = e^{f_i} \alpha$ with $f_2 > f_1 > f_0$.  Then the cobordism
$$
W_{20} := \left\{ (r,x)\ \big|\ f_0(x) \le r \le f_2(x) \right\}
$$
contains a contact-type hypersurface
$$
M_1 := \left\{ (f_1(x),x)\ \big|\ x \in M \right\} \subset W_{20}.
$$
As described at the end of Lecture~\ref{lec:compactness}, one can now choose a sequence of
compatible almost complex structures $\{J_{20}^N\}_{N \in \NN}$ on $\widehat{W}_{20}$
that are fixed outside a neighborhood of~$M_1$ but degenerate in this
neighborhood as $N \to \infty$, equivalent to replacing a small tubular
neighborhood of $M_1$ with increasingly large collars 
$[-N,N] \times M$ in which $J_{20}^N$ belongs to $\jJ(\alpha_1)$.
The resulting chain maps
$$
\Phi_{J_{20}^N} : CC_*^h(M,\alpha_2,J_2) \to CC_*^h(M,\alpha_0,J_0)
$$
are chain homotopic for all~$N$, but as $N \to \infty$, the index~$0$
cylinders counted by these maps converge to buildings with two levels,
the top one an index~$0$ cylinder in the completion of a cobordism from 
$(M,\alpha_1,J_1)$ to $(M,\alpha_2,J_2)$, while the bottom one also has
index~$0$ and lives in a cobordism from $(M,\alpha_0,J_0)$ to
$(M,\alpha_1,J_1)$.  The composition $\Phi_{21} \circ \Phi_{10}$
counts these broken cylinders, so this proves \eqref{eqn:composition}.

In particular, we conclude now that each of the cobordism maps
$$
\Phi : HC_*^h(M,\alpha_+,J_+) \to HC_*^h(M,\alpha_-,J_-)
$$ 
is an isomorphism,
since composing it with a cobordism map in the opposite direction must
give the identity.  The isomorphism class of $HC_*^h(M,\alpha,J)$ is therefore
independent of the auxiliary data $(\alpha,J)$, and will be denoted by
$$
HC_*^h(M,\xi).
$$
This is the \defin{cylindrical contact homology} of $(M,\xi)$ in the homotopy
class~$h$.  It is defined for any primitive homotopy class $h \in [S^1,M]$
and closed contact manifold that is $h$-admissible in the sense of
Definition~\ref{defn:admis}.  It is also invariant under contactomorphisms
in the following sense:

\begin{prop}
\label{prop:contacto}
Suppose $\varphi : (M_0,\xi_0) \to (M_1,\xi_1)$ is a contactomorphism with
$\varphi_*h_0 = h_1$, where $h_0 \in [S^1,M]$ is a primitive homotopy
class of loops, and $(M_1,\xi_1)$ is $h_1$-admissible.  Then $(M_0,\xi_0)$
is $h_0$-admissible, and $HC_*^{h_0}(M_0,\xi_0) \cong HC_*^{h_1}(M_1,\xi_1)$.
\end{prop}
\begin{proof}
Given an $h_1$-admissible contact form $\alpha_1$ on $(M_1,\xi_1)$ and an
$h_1$-regular $J_1 \in \jJ(\alpha_1)$, the contact form 
$\alpha_0 := \varphi^*\alpha_1$ on $M_0$
is $h_0$-admissible since $\varphi$ defines a bijection from 
$\pP_{h_0}(\alpha_0)$ to $\pP_{h_1}(\alpha_1)$ and also a bijection between
the sets of contractible Reeb orbits for $\alpha_0$ and~$\alpha_1$.
Since $\varphi_*\xi_0 = \xi_1$, $\alpha_0$ is a contact form for
$(M_0,\xi_0)$, hence the latter is $h_0$-admissible.  The diffeomorphism
$\tilde{\varphi} := \Id \times \varphi : \RR \times M_0 \to \RR \times M_1$
then maps $\p_r$ to $\p_r$, $R_{\alpha_0}$ to $R_{\alpha_1}$ and
$\xi_0$ to~$\xi_1$, thus $J_0 := \tilde{\varphi}^*J_1 \in \jJ(\alpha_0)$,
so $\tilde{\varphi}$ defines a biholomorphic map 
$(\RR \times M_0,J_0) \to (\RR \times M_1,J_1)$ and thus a bijection between
the sets of holomorphic cylinders in each.  It follows that $J_0$ is
$h_0$-regular, and the bijection
$\pP_{h_0}(\alpha_0) \to \pP_{h_1}(\alpha_1)$ defines an isomorphism between
the chain complexes defining $HC_*^{h_0}(M_0,\alpha_0,J_0)$ and
$HC_*^{h_1}(M_1,\alpha_1,J_1)$.
\end{proof}

\section{Computing $HC_*(\TT^3,\xi_k)$}

\subsection{The Morse-Bott setup}
\label{sec:MorseBott}

The contact form $\alpha_k$ on $\TT^3$ defined at the beginning of this
lecture has Reeb vector field
$$
R_k(\rho,\phi,\theta) = \cos(2\pi k \rho) \, \p_\theta +
\sin(2\pi k \rho)\, \p_\phi.
$$
Its Reeb orbits therefore preserve and define linear foliations on each of
the tori $\{\rho\} \times \TT^2$.  In particular, none of the closed orbits 
are contractible, though all of them are also degenerate, as they all come in
$S^1$-parametrized families foliating $\{\text{const}\} \times \TT^2$.
For certain homotopy classes $h \in [S^1,\TT^3]$, this yields a very easy
computation of $HC_*^h(\TT^3,\xi_k)$, namely whenever $h$ contains no
periodic orbits:

\begin{thm}
\label{thm:trivial}
Suppose $h \in [S^1,\TT^3]$ is any primitive homotopy class of loops
such that the projection $p: \TT^3 \to S^1 : (\rho,\phi,\theta) \mapsto \rho$
satisfies $p_*h \ne 0 \in [S^1,S^1]$.  Then $\alpha_k$ is $h$-admissible
and the resulting contact homology $HC_*^h(\TT^3,\xi_k)$ is trivial.
\qed
\end{thm}

Now for the interesting part.  Every primitive class $h \in [S^1,\TT^3]$
not covered by Theorem~\ref{thm:trivial} contains closed orbits of~$R_k$,
all of them degenerate since they come in $S^1$-parametrized families
foliating the tori $\{\text{const}\} \times \TT^2$.  This makes it not
immediately clear whether $(\TT^3,\xi_k)$ is $h$-admissible, though the
following observation in conjunction with Proposition~\ref{prop:contacto}
shows that if $HC_*^h(\TT^3,\xi_k)$ can be
defined, it will be the same for all the homotopy classes under consideration.

\begin{lemma}
\label{lemma:contacto}
Suppose $h_0, h_1 \in [S^1,\TT^3]$ are primitive homotopy classes that are
both mapped to the trivial class under the projection
$\TT^3 \to S^1 : (\rho,\phi,\theta) \mapsto \rho$.  Then there exists a
contactomorphism $\varphi : (\TT^3,\xi_k) \to (\TT^3,\xi_k)$ satisfying
$\varphi_*h_0 = h_1$.
\end{lemma}
\begin{proof}
We can represent $h_i$ for $i=0,1$ by loops of the form 
$\gamma_i(t) = (0,\beta_i(t)) \in S^1 \times \TT^2$, where the loops
$\beta_i : S^1 \to \TT^2$ are embedded and thus represent generators
of $\pi_1(\TT^2) = \ZZ^2$.  One can thus find a matrix 
$\begin{pmatrix} m & n \\ p & q \end{pmatrix} \in \SL(2,\ZZ)$
such that the diffeomorphism
$$
\varphi : \TT^3 \to \TT^3 : (\rho,\phi,\theta) \mapsto (\rho,m\phi + n\theta,
p\phi + q\theta)
$$
satisfies $\varphi_*h_0 = h_1$.  We have
\begin{equation*}
\begin{split}
\varphi^*\alpha_k &= \left[ q \cos(2\pi k \rho) + n \sin(2\pi k\rho) \right]\, d\theta
+ \left[ p \cos(2\pi k \rho) + m \sin(2\pi k \rho) \right] \, d\phi \\
&=: F(\rho)\, d\theta + G(\rho)\, d\phi.
\end{split}
\end{equation*}
The loop $(F, G) : S^1 \to \RR^2$ satisfies
$$
\begin{pmatrix}
F(\rho) \\
G(\rho)
\end{pmatrix}
=
\begin{pmatrix}
q & n \\
p & m
\end{pmatrix}
\begin{pmatrix}
\cos(2\pi k \rho) \\
\sin(2\pi k \rho)
\end{pmatrix},
$$
where $\begin{pmatrix} q & n \\ p & m \end{pmatrix} \in \SL(2,\ZZ)$,
thus $(F,G)$ winds $k$ times about the origin.  Any choice of homotopy
from $(F,G)$ to $(\cos(2\pi k \rho),\sin(2\pi k \rho))$ through loops
$(F_s,G_s) : S^1 \to \RR^2$
winding $k$ times about the origin with positive rotational velocity
then gives rise to a homotopy from $\varphi^*\alpha_k$ to $\alpha_k$ through
contact forms $F_s(\rho)\, d\theta + G_s(\rho)\, d\phi$.  Gray's stability
theorem therefore yields a contactomorphism 
$\psi : (\TT^3,\xi_k) \to (\TT^3,\ker \varphi^*\alpha_k)$ with $\psi$
smoothly isotopic to the identity.  The map $\varphi \circ \psi$ is thus
a contactomorphism of $(\TT^3,\xi_k)$ with $(\varphi \circ \psi)_*h_0 =
\varphi_*\psi_*h_0 = \varphi_*h_0 = h_1$.
\end{proof}

In light of the lemma, we are free from now on to restrict our attention to
the particular homotopy class
$$
h := [t \mapsto (0,0,t)],
$$
which is the homotopy class of the $1$-periodic orbits foliating the
$k$ tori 
$$
T_m := \{m/k\} \times \TT^2, \qquad m=0,\ldots,k-1
$$
since $R_k(m/k,\phi,\theta) = \p_\theta$.
Though the orbits on these tori are degenerate, it is not hard to show 
that they all satisfy the
Morse-Bott condition; in fact, $\alpha_k$ is a Morse-Bott contact form.
We will explain a self-contained computation of $HC_*^h(\TT^3,\xi_k)$ in
the next two sections without using the Morse-Bott condition---but first,
it seems worthwhile to sketch how one can guess the answer using
Morse-Bott data.

Bourgeois's thesis \cite{Bourgeois:thesis} gives a prescription for
calculating contact homology in Morse-Bott settings, i.e.~for deducing what 
orbits and what holomorphic curves will appear under certain standard ways of
perturbing the Morse-Bott contact form to make it nondegenerate.
Notice first that the only orbits in $\pP_h(\alpha_k)$ are the ones that
foliate the $k$ tori $T_0,\ldots,T_{k-1}$, and they all have period~$1$.
By \eqref{eqn:contactArea}, it follows that
for any $J \in \jJ(\alpha_k)$, there can be no nontrivial $J$-holomorphic
cylinders connecting two orbits in~$\pP_h(\alpha_k)$.  This makes the
calculation of $HC_*^h(\TT^3,\xi_k)$ sound trivial, but of course there is
more to the story since $\alpha_k$ is not admissible; indeed, the chain
complex $CC_*(\TT^3,\alpha_k)$ is not even well defined.  The prescription
in \cite{Bourgeois:thesis} now gives the following.
Each of the families of orbits in $T_0,\ldots,T_{k-1}$ is
parametrized by~$S^1$, and by a standard perturbation technique,
any choice of a Morse function $f_m : S^1 \to \RR$ for $m=0,\ldots,k-1$
yields a contact form $\alpha_k'$ that is $C^\infty$-close to~$\alpha_k$, 
matches it outside a neighborhood of~$T_m$, but has a nondegenerate
Reeb orbit on $T_m$ for each critical point of~$f_m$, while
every other closed orbit in the perturbed region can be assumed to have
arbitrarily large period.  Moreover, there is a corresponding perturbation
from $J \in \jJ(\alpha_k)$ to $J' \in \jJ(\alpha_k')$ such that every
gradient flow line of the function $f_m : S^1 \to \RR$ gives rise to a
$J'$-holomorphic cylinder in $\RR \times \TT^3$ connecting the corresponding 
nondegenerate Reeb orbits along~$T_m$.  In the present situation, since 
no $J$-holomorphic
cylinders of the relevant type exist before the perturbation, the only
ones after the perturbation are those that come from gradient flow lines.

Now imagine performing a similar perturbation near every $T_0,\ldots,T_{k-1}$,
using Morse functions $f_0,\ldots,f_{k-1} : S^1 \to \RR$ that each have exactly
two critical points.  For the perturbed contact form $\alpha_k'$, 
$\pP_h(\alpha_k')$ now consists of exactly $2k$ orbits 
$$
\gamma_0^\pm,\ldots,\gamma_{k-1}^\pm \in \pP_h(\alpha_k'),
$$
where we denote by $\gamma_m^+$ and $\gamma_m^-$ the orbits on $T_m$
corresponding to the maximum and minimum of~$f_m$ respectively.  For the
obvious choice of trivialization $\tau$ for the contact bundle along
$\gamma_m^\pm$, one can relate the Conley-Zehnder indices to the Morse indices
of the corresponding critical points, giving
$$
\muCZ^\tau(\gamma_m^+) = 0, \qquad \muCZ^\tau(\gamma_m^-) = 1, \qquad
m=0,\ldots,k-1.
$$
Moreover, the two gradient flow lines connecting maximum and minimum for
each $f_m$ give rise two exactly two holomorphic cylinders in
$\mM^1(J',\gamma_m^-,\gamma_m^+)$ for each $m=0,\ldots,k-1$, and these are
all the curves that are counted for the differential on
$CC_*^h(\TT^3,\alpha_k',J')$.  Counting modulo~$2$, we thus have
$$
\p \langle \gamma_m^\pm \rangle = 0 \quad \text{ for all } \quad m=0,\ldots,k-1,
$$
implying
$$
HC_*^h(\TT^3,\alpha_k',J') = \begin{cases}
\ZZ_2^k & \ * = \text{odd},\\
\ZZ_2^k & \ * = \text{even}.
\end{cases}
$$
Let us state this as a theorem.

\begin{thm}
\label{thm:interesting}
Suppose $h \in [S^1,\TT^3]$ is a primitive homotopy class that maps to the
trivial class under the projection $\TT^3 \to S^1 : (\rho,\phi,\theta)
\mapsto \rho$.  Then $(\TT^3,\xi_k)$ is $h$-admissible and
$$
HC_*^h(\TT^3,\xi_k) \cong
\begin{cases}
\ZZ_2^k & \ * = \text{odd},\\
\ZZ_2^k & \ * = \text{even}.
\end{cases}
$$
\end{thm}

Theorem~\ref{thm:T3} is an immediate corollary of this: indeed, if $\varphi :
(\TT^3,\xi_k) \to (\TT^3,\xi_\ell)$ is a contactomorphism, choose any
$h \in [S^1,\TT^3]$ for which Theorem~\ref{thm:interesting} applies,
and let $h_0 := \varphi^*h \in [S^1,\TT^3]$.  Then $HC_*^h(\TT^3,\xi_\ell)
\cong \ZZ_2^{2\ell}$ implies via Proposition~\ref{prop:contacto} that
$HC_*^{h_0}(\TT^3,\xi_k) \cong \ZZ_2^{2\ell}$.  But Theorems~\ref{thm:trivial}
and~\ref{thm:interesting} imply that the latter is also either $0$ or
$\ZZ_2^{2k}$, hence $k=\ell$.

\subsection{A digression on the Floer equation}
\label{sec:FloerDigression}

In preparation for giving a self-contained proof of 
Theorem~\ref{thm:interesting}, we now explain a general procedure for
relating holomorphic cylinders in a symplectization to solutions of
the Floer equation.  This idea is loosely inspired by arguments in
\cite{EliashbergKimPolterovich}.

To motivate what follows, notice that on a neighborhood of
$T_0 = \{0\} \times \TT^2 \subset (\TT^3,\xi_k)$, we can write
$$
\alpha_k = \cos(2\pi k \rho) \left( d\theta + \beta \right),
$$
where $\beta := \tan(2\pi k \rho)\, d\phi$ defines a Liouville form on
the annulus $\AA := [-1/8,1/8] \times S^1$ with coordinates $(\rho,\phi)$.
This makes the neighborhood $\AA \times S^1 \subset (\TT^3,\xi_k)$ a
special case of the following general construction.

\begin{defn}
Suppose $V$ is a $2n$-dimensional manifold with an exact symplectic form 
$d\beta$.  The contact manifold $(V \times S^1,\ker(d\theta + \beta))$
is then called the \defin{contactization} of~$(V,\beta)$.\footnote{Elsewhere
in the literature, the contactization is also often defined as $V \times \RR$
instead of $V \times S^1$.  The usage here is consistent with 
\cite{MassotNiederkruegerWendl}.}  Here $\theta$ denotes the coordinate
on the $S^1$~factor.
\end{defn}
It's easy to check that $d\theta + \beta$ is indeed a contact form on
$V \times S^1$ whenever $d\beta$ is symplectic on~$V$: the latter means
$(d\beta)^n > 0$ on~$V$, so
$$
(d\theta + \beta) \wedge \left[ d(d\theta + \beta)\right]^n =
(d\theta + \beta) \wedge (d\beta)^n = d\theta \wedge (d\beta)^n > 0.
$$

Now here's a cute trick one can play with contactizations.  For the rest of
this subsection, assume 
$$
(V,d\beta)
$$
is an arbitrary compact $2n$-dimensional exact symplectic manifold with boundary.
Fix a smooth function
$$
H : V \times S^1 \to \RR,
$$
which we shall think of in the following as a time-dependent Hamiltonian
$H_\theta := H(\cdot,\theta) : V \to \RR$ on $(V,d\beta)$.  
The $2$-form on $V \times S^1$ defined by
$$
\Omega = d\beta + d\theta \wedge dH = d(\beta - H\, d\theta)
$$
is then \emph{fiberwise symplectic}, meaning its restriction to each
of the fibers of the projection map $V \times S^1 \to S^1$ is symplectic.
We claim that for every $\epsilon > 0$ sufficiently small,
$$
\lambda_\epsilon := d\theta + \epsilon(\beta - H\, d\theta)
$$
defines a contact form on $V \times S^1$.  This is a variation on the
construction that was used by Thurston and Winkelnkemper
\cite{ThurstonWinkelnkemper} to define contact forms out of open book
decompositions, and the proof is simple enough: since $d\lambda_\epsilon =
\epsilon \Omega$, we just need to check that $\lambda_\epsilon \wedge
\Omega^n > 0$ for $\epsilon > 0$ sufficiently small, and indeed,
$$
\lambda_\epsilon \wedge \Omega^n = d\theta \wedge (d\beta)^n +
\epsilon (\beta - H\, d\theta) \wedge \Omega^n > 0
$$
since the first term is a volume form and $\epsilon$ is small.  
To see the relation between
$\lambda_\epsilon$ and the contactization, we can write
$$
\lambda_\epsilon = (1 - \epsilon H) \, d\theta + \epsilon \beta =
(1 - \epsilon H) \left( d\theta + \frac{\epsilon}{1 - \epsilon H} \beta \right)
$$
and observe that $\frac{\epsilon}{1 - \epsilon H}\beta$ is also a Liouville
form on $V$ whenever $H$ is $\theta$-independent and $\epsilon > 0$ is
sufficiently small.

The Reeb vector fields $R_\epsilon$ for $\lambda_\epsilon$ vary with~$\epsilon$, 
but their directions do not, since $d\lambda_\epsilon = \epsilon \Omega$ has
the same kernel for every~$\epsilon$.  Moreover, while $\lambda_\epsilon$
ceases to be a contact form when $\epsilon \to 0$, the Reeb vector fields
still have a well-defined limit: they converge as $\epsilon \to 0$ to the
unique vector field $R_0$ satisfying
$$
d\theta(R_0) \equiv 1 \quad \text{ and } \quad \Omega(R_0,\cdot) \equiv 0.
$$
The latter can be written more explicitly as
$$
R_0 = \p_\theta + X_\theta,
$$
where $X_\theta$ is the time-dependent Hamiltonian vector field determined
by~$H_\theta$, i.e.~via the condition
$$
d\beta(X_\theta,\cdot) = - dH_\theta.
$$
As one can easily compute, the reason for this nice behavior as
$\epsilon \to 0$ is that the $R_\epsilon$ are also the
Reeb vector fields for a smooth family of stable Hamiltonian structures:

\begin{prop}
The pairs $\hH_\epsilon := (\Omega,\lambda_\epsilon)$ for $\epsilon \ge 0$
sufficiently small define a smooth family of stable Hamiltonian structures
whose Reeb vector fields are~$R_\epsilon$.
\qed
\end{prop}

We shall write the hyperplane distributions induced by $\hH_\epsilon$ as
$$
\Xi_\epsilon := \ker \lambda_\epsilon \subset T(V \times S^1).
$$
These are contact structures for $\epsilon > 0$ small, and the space
$\jJ(\hH_\epsilon)$ of $\RR$-invariant almost complex structures on
$\RR \times (V \times S^1)$ compatible with $\hH_\epsilon$ is then
identical to~$\jJ(\lambda_\epsilon)$.  On the other hand for $\epsilon=0$,
$\Xi_0 = \ker d\theta$ is a foliation, namely it is the vertical subbundle
of the trivial fibration $V \times S^1 \to S^1$.  To interpret $\hH_0$, notice
that its closed Reeb orbits in the homotopy class of
$\gamma : S^1 \to V \times S^1 : t \mapsto (\text{const},t)$ are all of
the form $\gamma(t) = (x(t),t)$ where $x : S^1 \to V$
is a contractible $1$-periodic orbit of~$X_\theta$.
Moreover, suppose $J \in \jJ(\hH_0)$, which is equivalent to a choice of
compatible complex structure on the symplectic bundle $(\Xi_0,\Omega|_{\Xi_0})$,
or in other words, an $S^1$-parametrized family of $d\beta$-compatible
almost complex structures $\{J_\theta\}_{\theta \in S^1}$ on~$V$.  Then if
$$
u = (f,v,g) : \RR \times S^1 \to \RR \times (V \times S^1)
$$
is a $J$-holomorphic cylinder asymptotic at 
$\{\pm\infty\} \times S^1$ to two orbits of the form described above,
the nonlinear Cauchy-Riemann equation for $u$ turns out to imply that
$(f,g) : \RR \times S^1 \to \RR \times S^1$ is a holomorphic map with
degree~$1$ sending $\{\pm\infty\} \times S^1$ to $\{\pm\infty\} \times S^1$,
and we can therefore choose a unique
biholomorphic reparametrization of $u$ so that
$(f,g)$ becomes the identity map.  Having done this, the equation satisfied
by $v : \RR \times S^1 \to V$ is now
$$
\p_s v + J_t(v) (\p_t v - X_t(v)) = 0,
$$
in other words, the Floer equation for the data $\{J_\theta\}_{\theta \in S^1}$
and $\{H_\theta\}_{\theta \in S^1}$.  

To complete the analogy, notice that since $\Omega$ is exact, we can write
down a natural symplectic action functional with respect to each 
$\hH_\epsilon$ as
$$
\aA_\epsilon : C^\infty(S^1,V \times S^1) \to \RR : \gamma \mapsto
\int_{S^1} \gamma^*(\beta - H\, d\theta).
$$
For loops of the form $\gamma(t) = (x(t),t)$ with $x : S^1 \to V$
contractible, this reduces (give or take a sign---see Remark~\ref{remark:horror}) to the usual formula for
the Floer action functional
\begin{equation}
\label{eqn:FloerAction}
\aA_H(\gamma) = \int_{S^1} x^*\beta - \int_{S^1} H(x(t))\, dt =
\int_{\DD} \bar{x}^*d\beta - \int_{S^1} H(x(t))\, dt,
\end{equation}
where $\bar{x} : \DD \to V$ is any map satisfying $\bar{x}|_{\p\DD} = x$.
Stokes' theorem gives an easy relation between the action and the so-called
\emph{$\Omega$-energy} if $u : \RR \times S^1 \to \RR \times (V \times S^1)$
is a $J$-holomorphic curve for $J \in \jJ(\hH_\epsilon)$ and is
positively/negatively asymptotic to orbits $\gamma^\pm : S^1 \to V \times S^1$
at $s=\pm\infty$: we have
$$
0 \le \int_{\RR \times S^1} u^*\Omega = \aA_\epsilon(\gamma^+) - \aA_\epsilon(\gamma^-).
$$
If $u(s,t) = (s,v(s,t),t)$, then the left hand side is identical to the
definition of energy in Floer homology, namely
$$
E_H(v) := \int_{\RR \times S^1} d\beta(\p_s v, \p_t v - X_t(v)) \, ds \wedge dt =
\int_{\RR \times S^1} d\beta(\p_s v,J_t(v) \p_s v)\, ds \wedge dt,
$$
thus giving the familiar relation
\begin{equation}
\label{eqn:FloerBound}
E_H(v) = \aA_H(\gamma^+) - \aA_H(\gamma^-).
\end{equation}
To relate this to the usual notion of energy with respect to a stable
Hamiltonian structure, we write the usual formula
$$
E_\epsilon(u) := \sup_{\varphi \in \tT} \int_{\dot{\Sigma}} u^*\left[
d\big(\varphi(r) \lambda_\epsilon\big) + \Omega \right],
$$
with $\tT := \left\{ \varphi \in C^\infty(\RR , (-\epsilon_0,\epsilon_0)) \ \big|\ 
\varphi' > 0 \right\}$ for some constant $\epsilon_0 > 0$ sufficiently
small.  Notice first that for any fixed $\epsilon$, Stokes' theorem gives 
a bound for $E_\epsilon(u)$ in terms of the asymptotic orbits of $u$ since
$\Omega$ is exact.  Finally, in the case $\epsilon=0$ with
$u(s,t) = (s,v(s,t),t)$, we find
$$
E_0(u) = \sup_{\varphi \in \tT} \int_{\RR \times S^1} \varphi'(s)\, 
ds \wedge dt + \int_{\RR \times S^1} u^*\Omega =
2\epsilon_0 + E_H(v),
$$
so bounds on $E_0(u)$ are equivalent to bounds on the Floer homological
energy~$E_H(v)$.  The basic fact that Floer trajectories $v : \RR \times S^1
\to V$ with $E_H(v) < \infty$ are asymptotic to contractible $1$-periodic 
Hamiltonian orbits can now be regarded as a corollary of our
Theorem~\ref{thm:asymptotics}.

The above discussion gives a one-to-one correspondence
between a certain moduli space of unparametrized $J$-holomorphic cylinders
in $\RR \times (V \times S^1)$ and the moduli space of Floer trajectories
between contractible $1$-periodic orbits in $(V,d\beta)$ with Hamiltonian 
function~$H$.  If we can adequately understand the moduli space of Floer 
trajectories---in particular if we can classify them and prove that they are 
regular---then the idea will be to extend this classification via the implicit 
function theorem to any $J_\epsilon \in \jJ(\lambda_\epsilon)$ sufficiently 
close to $J$ for $\epsilon > 0$ small.
As the reader may be aware, classifying Floer trajectories is also not easy
in general, but it does become easy under certain conditions.  Simple examples
of contractible $1$-periodic Hamiltonian orbits are furnished by the constant loops
$\gamma(t) = x$ at critical points $x \in \Crit(H)$, and for each such 
orbit, $\gamma^*\Xi_0$
has a canonical homotopy class of unitary trivializations, the so-called
\defin{constant trivialization}.  The following fundamental result is 
commonly used in proving the isomorphism from Hamiltonian Floer homology 
to singular homology.

\begin{thm}
\label{thm:Floer}
Suppose $H : V \to \RR$ is a smooth Morse function with no critical points on 
the boundary, $J$ is a fixed $d\beta$-compatible almost complex structure 
on~$V$, and the gradient flow of $H$ with respect to the metric
$d\beta(\cdot,J\cdot)$ is Morse-Smale and transverse to~$\p V$.  
Given $\delta > 0$, let $H^\delta := \delta H : V \to \RR$,
with Hamiltonian vector field $X_{H^\delta} = \delta X_H$,
and consider the stable Hamiltonian structure
$$
\hH_0^\delta := (d\beta + d\theta \wedge d H^\delta , d\theta)
$$
on $V \times S^1$ with induced Reeb vector field $R_0^\delta = \p_\theta +
X_{H^\delta}$.  Then for all $\delta > 0$ sufficiently small,
the following statements hold.
\begin{enumerate}
\item The $1$-periodic $R_0^\delta$-orbit $\gamma_x : S^1 \to V \times S^1 : 
t \mapsto (x,t)$ arising from any critical
point $x \in \Crit(H)$ is nondegenerate, and its Conley-Zehnder index
relative to the constant trivialization~$\tau$ is related to the Morse index
$\ind(x) \in \{0,\ldots,2n\}$ by
\begin{equation}
\label{eqn:CZMorse}
\muCZ^\tau(\gamma_x) = n - \ind(x).
\end{equation}
\item Any trajectory $\gamma : \RR \to V$ satisfying the negative
gradient flow question $\dot{\gamma} = -\nabla H^\delta(\gamma)$ gives rise to a
Fredholm regular solution $v : \RR \times S^1 \to V : (s,t) \mapsto \gamma(s)$ 
of the time-independent Floer equation
\begin{equation}
\label{eqn:Floer}
\p_s v + J(v)(\p_t v - X_{H^\delta}(v)) = 0,
\end{equation}
and the virtual dimensions of the spaces of Floer trajectories near $v$ and
gradient flow trajectories near~$\gamma$ are the same.
\item Every $1$-periodic orbit of $X_{H^\delta}$ in $\mathring{V}$ 
is a constant loop at a critical point of~$H$.
\item Every finite-energy solution $v : \RR \times S^1 \to \mathring{V}$ 
of \eqref{eqn:Floer} is of the form $v(s,t) = \gamma(s)$ for some negative 
gradient flow trajectory $\gamma : \RR \to V$.
\end{enumerate}
\end{thm}
\begin{proof}
The following proof is based on arguments in \cite{SalamonZehnder:Morse},
see in particular Theorem~7.3.

For the first statement, let $\gamma(t) = (x,t)$ for $x \in \Crit(H)$ and
recall from Lecture~\ref{lec:asymptotic} the formula for the asymptotic operator of a
$1$-periodic orbit,
$$
\mathbf{A}_\gamma : \Gamma(\gamma^*\Xi_0) \to \Gamma(\gamma^*\Xi_0) :
\eta \mapsto - J \left( \nabla_t \eta - \nabla_\eta R_0^\delta\right),
$$
where $\nabla$ is any symmetric connection on $V \times S^1$.
Identifying $\Gamma(\gamma^*\Xi_0)$ in the natural way with
$C^\infty(S^1,T_x V)$, using the trivial connection and writing
$R_0^\delta(z,\theta) = \p_\theta + X_{H^\delta}(z) =
\p_\theta + \delta J(z) \nabla H(z)$, $\mathbf{A}_\gamma$ becomes the operator
$$
\mathbf{A}_\gamma = -J \p_t - \delta\nabla^2 H(x)
$$
on $C^\infty(S^1,T_x V)$, where $\nabla^2 H(x) : T_x V \to T_x V$ denotes the
Hessian of~$H$ at~$x$.  Choosing a unitary basis for $T_x V$ identifies
this with $-J_0 \p_t - \delta S$ for some symmetric $2n$-by-$2n$ matrix $S$
and the standard complex structure $J_0 = \begin{pmatrix} 0 & -\1 \\ \1 & 0 \end{pmatrix}$,
so $\ker \mathbf{A}_\gamma$ corresponds to the space of $1$-periodic 
solutions to $\dot{\eta} = \delta J_0 S \eta$.  The Morse condition implies
that $S$ is nonsingular, so the eigenvalues of $\delta J_0 S$ are all
nonzero, but they are also small since $\delta$ is small.  It follows that
nontrivial solutions of $\dot{\eta} = \delta J_0 S \eta$ cannot be $1$-periodic
if $S$ is nonsingular and $\delta$ is sufficiently small, thus proving
that $\ker \mathbf{A}_\gamma$ is trivial, hence $\gamma$ is nondegenerate.

To calculate $\muCZ^\tau(\gamma)$, note that $\lambda \in \sigma(\mathbf{A}_\gamma)$
if and only if there exists a nontrivial $1$-periodic solution $\eta$ to the
equation
$$
\dot{\eta} = J_0 (\delta S + \lambda) \eta.
$$
If $\delta$ and $\lambda$ are both close to~$0$, then the same argument again
implies that no such solution exists unless $\delta S + \lambda$ is
singular, meaning $\lambda \in \sigma(-\delta S)$.  On the other hand,
any constant loop $\eta(t) \in \ker(\lambda + \delta S)$ furnishes an
element of the $\lambda$-eigenspace of $\mathbf{A}_\gamma$, so we obtain
a bijection between the spectra of $\mathbf{A}_\gamma$ and $-\delta S$
in some neighborhood of~$0$.  It follows that if $S_\pm$ denotes a pair of
nonsingular symmetric matrices defining asymptotic operators
$\mathbf{A}_\pm = -J_0 \p_t - \delta S_\pm$, then the spectral flows are
related by
$$
\muspec(\mathbf{A}_-,\mathbf{A}_+) = - \muspec(S_-,S_+)
$$
when $\delta > 0$ is sufficiently small.  Denoting the maximal negative-definite
subspace of $S_\pm$ by $E^-(S_\pm)$, this relation implies
$$
\dim E^-(S_+) - \dim E^-(S_-) = \muCZ(\mathbf{A}_-) - \muCZ(\mathbf{A}_+).
$$
Now suppose $S_+$ is a coordinate expression for the Hessian $\nabla^2 H(x)$,
hence $\dim E^-(S_+) = \ind(x)$ and $\muCZ(\mathbf{A}_+) =
\muCZ^\tau(\gamma)$.  Choosing $S_- = \begin{pmatrix} \1 & 0 \\ 0 & -\1 \end{pmatrix}$
then gives $\dim E^-(S_-) = n$ and $\muCZ(\mathbf{A}_-) = 0$ by definition,
so $\muCZ^\tau(\gamma) = n - \ind(x)$ follows.

The second statement follows in a similar manner by writing down and
comparing the linearized operators for the Floer equation and the
negative gradient flow equation.  Let's leave this as an exercise.

For the third statement, suppose we have a sequence
$\delta_k \to 0$ and a sequence of loops $x_k : S^1 \to \mathring{V}$
satisfying $\dot{x}_k = X_{H^{\delta_k}}(x_k) = \delta_k X_H(x_k)$.
Pick a number $c > 0$ small enough for part~(1) of the theorem to hold
with $\delta=c$, choose a sequence of integers $N_k \in \NN$ such that
$$
N_k \delta_k \to c,
$$
and consider the loops $y_k : S^1 \to \mathring{V} : t \mapsto x_k(N_k t)$.
These satisfy
$$
\dot{y}_k = N_k \delta_k X_H(y_k),
$$
and since $X_H$ is $C^\infty$-bounded on $V$ and $N_k \delta_k$ is also
bounded, the Arzel\`a-Ascoli theorem provides a subsequence with
$$
y_k \to y_\infty \quad \text{ in } \quad C^\infty(S^1,V),
$$
where $y_\infty : S^1 \to V$ satisfies $\dot{y}_\infty = X_{H^c}(y_\infty)$
for $H^c := cH : V \to \RR$.  But $y_\infty$ is also constant: indeed,
since $y_k(t + 1/N_k) = y_k(t)$ and $N_k \to \infty$, we can find
for any $t \in S^1$ a sequence $q_k \in \ZZ$ satisfying $q_k / N_k \to t$, so
\begin{equation}
\label{eqn:tindependent}
y_\infty(t) = \lim_{k \to \infty} y_k(q_k/N_k) = \lim_{k \to \infty} y_k(0)
= y_\infty(0).
\end{equation}
Since the constant orbit $y_\infty$ is nondegenerate by part~(1) of the
theorem, there can only be one sequence of solutions to $\dot{y}_k =
X_{H^{N_k\delta_k}}(y_k)$ converging to~$y_\infty$, and we conclude that
$y_k$ is also constant for all $k$ sufficiently large.

We will now use a similar trick to prove the fourth statement in the theorem.
We shall work under the additional assumption that
\begin{equation}
\label{eqn:MorseCondition}
|\ind(x) - \ind(y)| \le 1 \quad \text{ for all pairs } \quad
x , y \in \Crit(H),
\end{equation}
which suffices for the application in \S\ref{sec:admis} 
below.\footnote{Lifting this
assumption requires gluing, whereas we shall only need the usual implicit
function theorem for Fredholm regular solutions of the Floer equation.}

Suppose to the contrary that there exists a sequence of positive numbers
$\delta_k \to 0$ with finite-energy solutions 
$v_k : \RR \times S^1 \to \mathring{V}$ of 
the equation $\p_s v_k + J(v_k)(\p_t v_k - X_{H^{\delta_k}}(v_k)) = 0$, where
each $v_k(s,t)$ is not $t$-independent.  By part~(3) of the theorem, we can
restrict to a subsequence and assume each $v_k$ for large $k$ is asymptotic
to a fixed pair of critical points $x_\pm = \lim_{s \to \pm\infty} v_k(s,\cdot)
\in \Crit(H)$, and $x_+ \ne x_-$ since $v_k$ would otherwise by constant
and therefore $t$-independent.  Choose a sequence $N_k \in \NN$ with
$$
N_k \to \infty \quad \text{ and } \quad N_k \delta_k \to c,
$$
where $c > 0$ is chosen sufficiently small for the first three statements
in the theorem to hold with $\delta=c$.  Define $w_k : \RR \times S^1 \to V$ by
$$
w_k(s,t) = v_k(N_k s, N_k t).
$$
Then $w_k$ satisfies another time-independent Floer equation,
\begin{equation}
\label{eqn:wkEquation}
\p_s w_k + J(w_k) \left(\p_t w_k - X_{H^{N_k \delta_k }}(w_k)\right) = 0,
\end{equation}
where the Hamiltonian functions $H^{N_k \delta_k}$ converge to~$H^c$.
The standard compactness theorem for Floer trajectories should now imply
that a subsequence of $w_k$ converges to a broken Floer trajectory
whose levels will be $t$-independent.  Since the setting may seem a bit
nonstandard, here are some details.

The sequence $w_k$ is uniformly $C^0$-bounded since $V$ is compact.
We claim that it is also $C^1$-bounded.  If not, then there is a sequence
$z_k = (s_k,t_k) \in \RR \times S^1$ with $|d w_k(z_k)| =: R_k \to \infty$,
and we can use the usual rescaling trick from Lecture~\ref{lec:compactness} to define a sequence
$$
f_k : \DD_{\epsilon_k R_k} \to V : z \mapsto w_k(z_k + z / R_k)
$$
for a suitable sequence $\epsilon_k \to 0$ with $\epsilon_k R_k \to \infty$
and $|d w_k(z)| \le 2 R_k$ for all $z \in \DD_{\epsilon_k}(z_k)$.  The
latter implies that $f_k$ satisfies a local $C^1$-bound independent of~$k$,
and since
$$
\p_s f_k + J(f_k) \left( \p_t f_k - \frac{1}{R_k} J(f_k) X_{H^{N_k \delta_k}}(f_k) \right),
$$
elliptic regularity (see Remark~\ref{remark:reg} below) provides a subsequence
for which $f_k$ converges in $C^\infty_\loc(\CC,V)$ to a $J$-holomorphic plane
$f_\infty : \CC \to V$, which is nonconstant since
$$
| d f_\infty(0) | = \lim_{k \to \infty} | d f_k(0) | = 1.
$$
Since $v_k$ and therefore $w_k$ are all asymptotic to fixed constant 
orbits $x_\pm$, we have a uniform bound on the Floer energies of $w_k$,
\begin{equation}
\label{eqn:wEnergyBound}
\begin{split}
E_{H^{N_k \delta_k}}(w_k) = \aA_{H^{N_k \delta_k}}(x_+) -
\aA_{H^{N_k \delta_k}}(x_-) = N_k \delta_k \left[H(x_-) - H(x_+) \right],
\end{split}
\end{equation}
where the right hand side is bounded since $N_k \delta_k \to c$.
Using change of variables and the fact that $d\beta(\p_s f_k,J(f_k)\, \p_s f_k)
\ge 0$, this implies a uniform bound
\begin{equation*}
\begin{split}
\int_{\DD_{\epsilon_k R_k}} d\beta(\p_s f_k, & J(f_k)\, \p_s f_k)\, ds \wedge dt =
\int_{\DD_{\epsilon_k(z_k)}} d\beta(\p_s v_k,J(v_k)\, \p_s v_k)\, ds \wedge dt \\ &
\le \int_{\RR \times S^1} d\beta(\p_s v_k,J(v_k)\, \p_s v_k)\, ds \wedge dt 
= E_{H^{N_k \delta_k}}(w_k) \le C,
\end{split}
\end{equation*}
thus
$$
\int_\CC f_\infty^*d\beta = \int_\CC d\beta(\p_s f_\infty,\p_t f_\infty)\, ds \wedge dt
= \int_\CC d\beta(\p_s f_\infty,J(f_\infty)\, \p_s f_\infty)\, ds \wedge dt < \infty.
$$
The removable singularity theorem now extends $f_\infty$ to a nonconstant
$J$-holomorphic sphere $f_\infty : S^2 \to V$, but this violates Stokes'
theorem since $J$ is tamed by an exact symplectic form.

We've now shown that the sequence $w_k : \RR \times S^1 \to V$ is uniformly
$C^1$-bounded, and it has bounded energy due to \eqref{eqn:wEnergyBound}.
Pick any sequence $s_k \in \RR$ and consider the sequence of 
translated Floer trajectories
$$
\tilde{w}_k(s,t) := w_k(s + s_k,t).
$$
These are also uniformly $C^1$-bounded, so by elliptic regularity
(see Remark~\ref{remark:reg} again), a subsequence converges in 
$C^\infty_\loc(\RR \times S^1)$ to a map $w_\infty : \RR \times S^1 \to V$ 
satisfying
$$
\p_s w_\infty + J(w_\infty) \left( \p_t w_\infty - X_{H^c}(w_\infty) \right) = 0,
$$
and it has finite energy $E_{H^c}(w_\infty) < \infty$
due to \eqref{eqn:wEnergyBound}, implying that $w_\infty$ is asymptotic 
to a pair of $1$-periodic orbits of $X_{H^c}$ as $s \to \pm\infty$.
By the same argument used in \eqref{eqn:tindependent} above,
$w_\infty$ is also $t$-independent.
It follows that $w_\infty(s,t) = \gamma_\infty(s)$ for some nonconstant
gradient flow trajectory $\gamma_\infty : \RR \to \mathring{V}$.  
Depending on the choice of sequence $s_k$, this trajectory 
may or may not be constant, but we can always choose $s_k$ to guarantee that
$\gamma_\infty$ is not constant: indeed,
since each $w_k$ is asymptotic to two separate critical points at $\pm\infty$,
$s_k \in \RR$ can be chosen such that $w_k(s_k,0)$ stays a fixed
distance away from every critical point of~$H$, and then
$$
w_\infty(0,0) = \lim_{k \to \infty} w_k(s_k,0) \not\in \Crit(H^c).
$$
One can now adapt the argument of Proposition~\ref{prop:cylinders}
to find various sequences $s_k \in \RR$ that yield potentially separate 
limiting trajectories forming the levels of a broken trajectory, which is the
limit of $w_k$ in the Floer topology.  But since all the levels are
$t$-independent and the gradient flow of $H^c$ is Morse-Smale,
condition \eqref{eqn:MorseCondition} implies that the most complicated (and therefore
the only) limit possible involves a single level $w_\infty(s,t) = \gamma(s)$,
which is a gradient flow trajectory between critical points whose Morse
indices differ by~$1$.  This trajectory is Fredholm regular and has index~$1$
due to part~(2) of the theorem, thus by the implicit function theorem,
the only solutions to 
\eqref{eqn:wkEquation} that can converge to $w_\infty$ are the obvious
reparametrizations of~$\gamma$,
i.e.~they are also $t$-independent.  This is a contradiction.
\end{proof}

\begin{remark}
\label{remark:reg}
In previous lectures we've used the theorem that ``$C^1$-bounds imply
$C^\infty$-bounds'' to prove compactness for $J$-holomorphic curves, but
not for solutions of inhomogeneous Cauchy-Riemann type equations such as
the Floer trajectories $w_k$ and rescalings~$f_k$ in the above proof.
There is an easy trick to reduce these to our standard setup: as we've
already seen, solutions of the Floer equation are equivalent to honest
pseudoholomorphic curves in the symplectization of a certain stable Hamiltonian
structure, which is a manifold of two dimensions higher.  A similar trick
can be used for any inhomogeneous Cauchy-Riemann type equation
$\dbar_J f = \nu$, reducing it to an honest Cauchy-Riemann type equation at
the cost of adding two dimensions.  This trick was used already by Gromov,
see \cite{Gromov}*{1.4.C}.
\end{remark}

\begin{remark}
\label{remark:horror}
You may notice with some horror
that \eqref{eqn:CZMorse} differs by a sign from what is stated in 
\cite{SalamonZehnder:Morse}.  As far as I can tell, the discrepancy
arises from the fact that while Floer homology is traditionally defined in
terms of a negative gradient flow for the action functional, SFT is based on
a \emph{positive} gradient flow---this is also why the action functional 
in \eqref{eqn:FloerAction} differs by a sign from what we saw
in Lecture~\ref{lec:intro}.  If one takes as an axiom that the Conley-Zehnder index should
serve as a ``relative Morse index'' for the action functional, then 
changing the sign of the functional also reverses the signs
of Conley-Zehnder indices, so as a result there appear to be two parallel
sign conventions for Conley-Zehnder indices in different sectors of the
literature.  I'm sorry.  It's not my fault.
\end{remark}

Returning now to the family $\hH_\epsilon$, choose $\delta > 0$ sufficiently
small for Theorem~\ref{thm:Floer} to hold and define a modified family of
stable Hamiltonian structures on $V \times S^1$ by
$$
\hH_\epsilon^\delta = (\Omega^\delta,\lambda_\epsilon^\delta),
$$
where
$$
\Omega^\delta := d\beta + d\theta \wedge d H^\delta \quad \text{ and } \quad
\lambda_\epsilon^\delta := d\theta + \epsilon (\beta - H^\delta\, d\theta).
$$
Denote the induced hyperplane distributions and Reeb vector fields by
$\Xi_\epsilon^\delta$ and $R_\epsilon^\delta$ respectively.
We have only changed the Hamiltonian $H$ by rescaling, so all previous
statements about $\hH_\epsilon$ also apply to $\hH_\epsilon^\delta$,
in particular $\lambda_\epsilon^\delta$ is contact and $\jJ(\hH_\epsilon^\delta) =
\jJ(\lambda_\epsilon^\delta)$ for all $\epsilon > 0$ sufficiently small,
though the upper bound for the allowed range of $\epsilon$ may now depend
on~$\delta$.  Once $\delta > 0$ is fixed by the requirements of
Theorem~\ref{thm:Floer}, we are still free to take $\epsilon > 0$ is small
as we like.

\begin{thm}
\label{thm:perturbation}
Assume the same hypotheses as in Theorem~\ref{thm:Floer}, including
\eqref{eqn:MorseCondition}, and denote the unique extension of $J$ to
an $\RR$-invariant almost complex structure in $\jJ(\hH_0^\delta)$ by~$J_0$.
Given $\delta$ sufficiently small and any smooth family of compatible 
$\RR$-invariant almost 
complex structures $J_\epsilon \in \jJ(\hH_\epsilon^\delta)$ matching
$J_0$ at $\epsilon=0$, there exists $\epsilon_0 > 0$ such that
every critical point $x \in \Crit(H)$ gives rise to a smooth family of
nondegenerate closed $R_\epsilon^\delta$-orbits
$$
x^\epsilon : S^1 \to V \times S^1 \qquad \epsilon \in [0,\epsilon_0]
$$
with $x^0(t) = (x,t)$, and every gradient flow trajectory
$\gamma : \RR \to V$ for $H$ gives rise to a smooth family of Fredholm
regular $J_\epsilon$-holomorphic cylinders
$$
u_\gamma^\epsilon : \RR \times S^1 \to \RR \times (V \times S^1) \qquad
\epsilon \in [0,\epsilon_0]
$$
with $u_\gamma^0(s,t) = (s,\gamma(\delta s),t)$.
Moreover, for all $\epsilon \in [0,\epsilon_0]$, every
closed $R_\epsilon^\delta$-orbit homotopic to $t \mapsto (\operatorname{const},t)$
belongs to one of the families $x^\epsilon$ up to parametrization, and
every $J_\epsilon$-holomorphic cylinder with a positive and a negative
end asymptotic to orbits of this type belongs to one of the
families $u_\gamma^\epsilon$, up to biholomorphic parametrization.
\end{thm}
\begin{proof}
The first part is immediate from the implicit function theorem since the
orbits $x^0(t) = (x,t)$ are nondegenerate and the curves
$u_\gamma^0(s,t) = (s,\gamma(\delta s),t)$ are Fredholm regular by
Theorem~\ref{thm:Floer}.  For the uniqueness statement, observe that if
$\epsilon_k \to 0$ and $\gamma_k$ is a sequence of $R_{\epsilon_k}^\delta$-orbits
in the relevant homotopy class, then their periods are uniformly bounded,
so Arzel\`a-Ascoli gives a subsequence convergent to a closed
$R_0^\delta$-orbit, which is a nondegenerate orbit of the form
$x^0(t) = (x,t)$ for $x \in \Crit(H)$ by Theorem~\ref{thm:Floer}, so 
sequences converging to this orbit are unique by the implicit function theorem.
A similar argument proves uniqueness of $J_\epsilon$-holomorphic cylinders:
if $\epsilon_k \to 0$ and $u_k$ is a $J_{\epsilon_k}$-holomorphic sequence,
then first by the uniqueness of the orbits, we can extract a subsequence
for which all $u_k$ are asymptotic at both ends to orbits in fixed
families $x_\pm^{\epsilon_k}$ converging to~$x_\pm^0(t) = (x_\pm,t)$ as
$k \to \infty$.  Since $\Omega$ is exact, Stokes' theorem then gives a
uniform bound on the energies $E_{\epsilon_k}(u_k)$.  Since all
$R_0^\delta$-orbits in the relevant homotopy class are nondegenerate and
none are contractible, one can now prove as in Proposition~\ref{prop:cylinders}
that $u_k$ has a subsequence convergent to a finite-energy stable
$J_0$-holomorphic building~$\mathbf{u}_\infty$ consisting only of cylinders.  
Its levels are asymptotic to orbits of the form
$x(t) = (x,t)$ for $x \in \Crit(H)$, thus they can be parametrized as
$(s,t) \mapsto (s,v(s,t),t)$ for $v : \RR \times S^1 \to V$ satisfying
the $H^\delta$-Floer equation, hence $v(s,t) = \gamma(\delta s)$ by
Theorem~\ref{thm:Floer}.  Now since $\nabla H$ is Morse-Smale and
indices of critical points can only differ by at most~$1$, the building
$\mathbf{u}_\infty$ can have at most one nontrivial level
$u_\infty(s,t) = (s,\gamma(\delta s),t)$, implying $u_k \to u_\infty$.
Since $u_\infty$ is Fredholm regular, the implicit function theorem does
the rest.
\end{proof}

\subsection{Admissible data for $(\TT^3,\xi_k)$}
\label{sec:admis}

We now complete the computation of the cylindrical contact homology
$HC_*^h(\TT^3,\xi_k)$.
We can assume via Lemma~\ref{lemma:contacto}
that $h$ is the homotopy class of the orbits in the special set of tori
$$
T_m = \{m/k\} \times \TT^2 \subset \TT^3, \qquad m=0,\ldots,k-1.
$$
Let's focus for now on the case $k=1$, as the general case will simply
be a $k$-fold cover of this.  Thanks to the Morse-Bott discussion
in \S\ref{sec:MorseBott}, we know what we're looking for: we want an 
$h$-admissible contact form
$\alpha$ for $(\TT^3,\xi_1)$ such that $\pP_h(\alpha)$ contains exactly
two orbits, both in $T_0 \subset \TT^3$, along with an $h$-regular
$J \in \jJ(\alpha)$ such that the differential on $CC_*^h(\TT^3,\alpha)$
counts exactly two $J$-holomorphic cylinders that connect the two
orbits in~$T_0$.  Let $\AA$ denote the annulus
$$
\AA = [-1,1] \times S^1
$$
with coordinates $(\rho,\phi)$.  This will play the role of the Liouville
manifold $(V,d\beta)$ from the previous section, and we set
$$
\beta := \rho\, d\phi.
$$
For the Hamiltonian $H : \AA \to \RR$, choose a Morse function with the 
following properties:
\begin{enumerate}
\item $H$ has a minimum at $x_0 = (0,0)$, an index~$1$ critical point
at $x_1 = (0,1/2)$, and no other critical points;
\item $H(\rho,\phi) = |\rho|$ for $1/2 \le |\rho| \le 1$;
\item The gradient flow of $H$ with respect to the standard Euclidean
metric on $[-1,1] \times S^1$ is Morse-Smale.
\end{enumerate}
Fix a number $\delta > 0$ sufficiently small so that 
Theorem~\ref{thm:Floer} applies for Floer trajectories of 
$H^\delta := \delta H$ in~$\AA$, and since it will turn out to be useful in
Lemma~\ref{lemma:cannotEscape} below, assume without loss of generality
$$
\delta \in \QQ.
$$
Then following the prescription described above,
we consider the family of stable Hamiltonian structures
$\hH_\epsilon^\delta = (\Omega^\delta,\lambda_\epsilon^\delta)$ on 
$\AA \times S^1$ for $\epsilon \ge 0$ small, where
$$
\lambda_\epsilon^\delta = (1 - \epsilon \delta H)\, d\theta + \epsilon \rho\, d\phi,
\qquad
\Omega^\delta = d\rho \wedge d\phi + \delta \, d\theta \wedge dH,
$$
with induced Reeb vector fields $R_\epsilon^\delta$ and hyperplane
distributions $\Xi_\epsilon^\delta := \ker \lambda_\epsilon^\delta$.
Choose $J_\epsilon \in \jJ(\hH_\epsilon^\delta)$ to be any smooth family such 
that $J_0|_{\Xi_0^\delta}$ matches the standard complex structure on $\AA$
defined by $J_0 \p_\rho = \p_\phi$.  Then for all $\epsilon > 0$ sufficiently
small, Theorems~\ref{thm:Floer} and \ref{thm:perturbation} give a complete 
classification of all closed $R_\epsilon^\delta$-orbits in 
$\AA \times S^1$ homotopic to $t \mapsto (0,0,t)$, as well as a classification
of all $J_\epsilon$-holomorphic cylinders asymptotic to them.
Up to parametrization, there are exactly two such orbits,
$$
\gamma_i^\epsilon : S^1 \to \AA \times S^1, \qquad i=0,1,
$$
which correspond to the Morse critical points $x_0$ and $x_1$ and thus
by \eqref{eqn:CZMorse} have Conley-Zehnder indices
$$
\muCZ^\tau(\gamma_i^\epsilon) = 1 - \ind(x_i) = 1 - i \in \{0,1\}
$$
relative to the constant trivialization~$\tau$.  There are also exactly
two $J_\epsilon$-holomorphic cylinders
$$
u_\pm^\epsilon : \RR \times S^1 \to \RR \times (\AA \times S^1),
$$
corresponding to the two negative gradient flow lines that descend from
$x_1$ to $x_0$, thus the $u_\pm^\epsilon$ are index~$1$ curves with a
negative end approaching $\gamma_1^\epsilon$ and a positive end 
approaching~$\gamma_0^\epsilon$.  If we can suitably embed this model into
$(\TT^3,\xi_1)$ and show that all the orbits and curves needing to be
counted are contained in the model, then we will have a complete description
of $HC_*^h(\TT^3,\xi_1)$, with two generators $\langle \gamma_0^\epsilon \rangle$
and $\langle \gamma_1^\epsilon \rangle$, of even and odd degree respectively,
satisfying 
$$
\p\langle \gamma_0^\epsilon \rangle = 2 \langle \gamma_1^\epsilon \rangle = 0
\quad \text{ and } \quad
\p \langle \gamma_1^\epsilon \rangle = 0
$$
since the former counts two curves and the latter counts none.

\begin{lemma}
\label{lemma:cannotEscape}
For any $\epsilon > 0$ sufficiently small, there exists a contact embedding
of 
$$
(\AA \times S^1,\ker \lambda_\epsilon^\delta) \hookrightarrow (\TT^3,\xi_1)
$$
identifying the homotopy class of the loops
$t \mapsto (0,0,t)$ in $\AA \times S^1$ with~$h$.  Moreover, the
contact form $\lambda_\epsilon^\delta$ and almost complex structure
$J_\epsilon \in \jJ(\hH_\epsilon^\delta)$ can then be extended 
to an $h$-admissible contact form $\alpha$ on $(\TT^3,\xi_1)$ and an
$h$-regular almost complex structure $J \in \jJ(\alpha)$ such that
$\gamma_0^\epsilon$ and $\gamma_1^\epsilon$ are the only orbits in
$\pP_h(\alpha)$, and all
$J$-holomorphic cylinders with a positive and a negative end asymptotic
to either of these orbits are contained in the interior of $\AA \times S^1$.
\end{lemma}
\begin{proof}
We've chosen $\beta$ and $H$ so that in the region $1/2 \le |\rho| \le 1$,
$$
\alpha := \lambda_\epsilon^\delta = (1 - \epsilon \delta |\rho|)\, d\theta +
\epsilon \rho\, d\phi =: f(\rho)\, d\theta + g(\rho)\, d\phi,
$$
so the Reeb vector field on this region has the form
$\frac{1}{D(\rho)} (g'(\rho)\, \p_\theta - f'(\rho)\, \p_\phi)$.
Notice that
$$
\frac{f'(\rho)}{g'(\rho)} = \mp \frac{\epsilon \delta}{\epsilon} = \mp\delta,
$$
and we assumed $\delta \in \QQ$, so the Reeb orbits in this region are all
periodic.  Next, pick a large number $N \gg 1$ and extend $\alpha$ to a
contact form on $[-N,N] \times S^1 \times S^1$ via the same formula.
Now extend the path $(f,g) : [-N,N] \to \RR^2$ to 
$\RR$ such that it has period $2N + 2$ and winds once around the origin
over the interval $[-N-1,N+1]$, with positive angular velocity.  This produces
a contact form $\alpha$ on
$$
\TT^3_N := \left(\RR \Big/ (2N+2)\ZZ \right) \times S^1 \times S^1
$$
which takes the form $f(\rho)\, d\theta + g(\rho)\, d\phi$ outside of
$|\rho| \le 1/2$.  We claim in fact that $\alpha$ is homotopic through
contact forms to one that takes this form globally, where $(f,g)$ may be
assumed to be a smooth loop winding once around the origin.
To see this, one need only homotop $H$ in the region $|\rho| \le 1/2$ to
a Morse-Bott function that depends only on the $\rho$-coordinate; the
contact condition holds for all Hamiltonians in this homotopy as long as
$\epsilon > 0$ is sufficiently small.  With this understood, the obvious
diffeomorphism
$$
\TT^3_N \to \TT^3 : (\rho,\phi,\theta) \mapsto \left( \frac{\rho}{2N+2},\phi,\theta \right)
$$
pushes $\ker \alpha$ forward to a contact structure isotopic to one of the
form $F(\rho)\, d\theta + G(\rho)\, d\phi$ for a loop
$(F,G) : S^1 \to \RR^2$ winding once around the origin, so taking a homotopy of
this loop to $(\cos(2\pi\rho),\sin(2\pi\rho))$ and applying Gray's stability
theorem produces a contactomorphism
$$
(\TT^3_N,\ker\alpha) \to (\TT^3,\xi_1)
$$
that is isotopic to the above diffeomorphism.

The construction clearly guarantees that no closed Reeb orbit of $\alpha$
outside $\AA \times S^1$ is homotopic to the preferred class~$h$,
and there are also no contractible orbits, so $\alpha$ is an $h$-admissible
contact form on~$\TT^3_N$.
Choose any extension of $J_\epsilon$ to some $J \in \jJ(\alpha)$
on~$\TT^3_N$.  We claim now that if $N$ is chosen sufficiently large, then
no $J$-holomorphic cylinder in $\RR \times \TT^3_N$ with one
positive end at either of the orbits $\gamma_i^\epsilon$ can ever 
venture outside the region $\RR \times (-1/2,1/2) \times \TT^2$.
Suppose in particular that $u$ is such a curve and its image intersects
$\RR \times \{1/2\} \times \TT^2$.  Since the entire
region $[1/2,N] \times \TT^2$ is foliated by closed Reeb orbits, we can
define $\Upsilon$ to be the set of Reeb orbits $\gamma$ in that region for which
the image of $u$ intersects $\RR \times \gamma$.  This is a closed subset
of the connected topological space of all Reeb orbits in
$[1/2,N] \times \TT^2$: indeed, if $\gamma_k \in \Upsilon$ is a sequence converging
to some orbit $\gamma_\infty$, then $u(z_k) \in \RR \times \gamma_k$ for some
sequence $z_k \in \RR \times S^1$, which must be contained in a compact
subset since the asymptotic orbits of $u$ lie
outside of $[1/2,N] \times \TT^2$, hence $z_k$ has a convergent subsequence
$z_k \to z_\infty \in \RR \times S^1$ with $u(z_\infty) \in \RR \times \gamma_\infty$,
proving $\gamma_\infty \in \Upsilon$.  We claim that $\Upsilon$ is also an open subset
of the space of orbits in $[1/2,N] \times \TT^2$.  This follows from
positivity of intersections, as every $\RR \times \gamma$ is also a
$J$-holomorphic curve: if $u(z) \in \RR \times \gamma$, then for every other
closed orbit $\gamma'$ close enough to~$\gamma$, there is a point $z' \in
\RR \times S^1$ near $z$ with $u(z') \in \RR \times \gamma'$.
This proves that, in fact, $u$ passes through $\RR \times \gamma$ for
\emph{every} orbit $\gamma$ in the region $[1/2,N] \times \TT^2$.  We will
now use this to show that if $N$ is sufficiently large, the contact area
of $u$ will be larger than is allowed by Stokes' theorem.

Let us write
$$
u(s,t) = (r(s,t),\rho(s,t),\phi(s,t),\theta(s,t)) \in \RR \times \left(\RR \big/ (2N+2)\ZZ \right) \times S^1 \times S^1
$$
and choose two points $\rho_1 \in [1/2,1]$ and $\rho_2 \in [N-1,N]$ which
are both regular values of the function $\rho : \RR \times S^1 \to \RR / (2N+2)\ZZ$.
The intersections of $u$ with the orbits in $[1/2,N] \times \TT^2$ imply that
the function $\rho(s,t)$ attains every value in $[1/2,N]$, and since
the asymptotic limits of $u$ lie outside this region,
$$
\uU := \rho^{-1}([\rho_1,\rho_2]) \subset \RR \times S^1
$$
is then a nonempty and compact smooth submanifold with boundary
$$
\p \uU = -C_1 \sqcup C_2,
$$
where $C_i := \rho^{-1}(\rho_i)$ for $i=1,2$.  Restricting $u$ to the
multicurves $C_i$ then gives a pair of smooth maps
$$
w_i : C_i \to \TT^2 : (s,t) \mapsto (\phi(s,t),\theta(s,t)), \qquad i=1,2,
$$
which are homologous to each other.  Denote the generators of $H_1(\TT^2)$
corresponding to the $\phi$- and $\theta$-coordinates by $\ell_\phi$ and
$\ell_\theta$ respectively, and suppose $[w_i] = m \ell_\phi + n \ell_\theta$
for $m,n \in \ZZ$.  The key observation now is that the restriction of
$\alpha$ to each of the tori $\{\rho_i\} \times \TT^2$ is a closed
$1$-form, thus for each $i=1,2$, $\int_{C_i} u^*\alpha$ depends only on the homology
class $m \ell_\phi + n \ell_\theta \in H_1(\TT^2)$ and not any further on
the maps~$w_i$.  In particular,
$$
\int_{C_i} u^*\alpha = f(\rho_i) n + g(\rho_i) m
$$
for $i=1,2$.  We now compute,
\begin{equation*}
\begin{split}
\int_{\uU} u^*d\alpha &= \int_{C_2} u^*\alpha - \int_{C_1} u^*\alpha =
n [ f(\rho_2) - f(\rho_1) ] + m [ g(\rho_2) - g(\rho_1) ] \\
&= n [ (1 - \epsilon \delta \rho_2) - (1 - \epsilon \delta \rho_1) ]
+ m [ \epsilon \rho_2 - \epsilon \rho_1 ] \\
&= \epsilon (\rho_2 - \rho_1) (m - n \delta)
\end{split}
\end{equation*}
This integral has to be positive since $u^*d\alpha \ge 0$ and $u$ is not
a trivial cylinder, thus $m - n\delta > 0$.  Moreover, $\delta$ was
assumed rational, so if $\delta = p/q$ for some $p,q \in \NN$, we have
$$
m - n \delta = \frac{1}{q} (m q - n p) \ge \frac{1}{q},
$$
implying
$$
\int_{\RR \times S^1} u^*d\alpha \ge \int_{\uU} u^*d\alpha \ge
\frac{\epsilon}{q} (\rho_2 - \rho_1) \ge \frac{\epsilon(N-2)}{q}.
$$
Having chosen $\delta$ (which determines $q$) and $\epsilon$ in advance,
we are free to make $N$ as large as we like.  But by
\eqref{eqn:contactArea}, $\int_{\RR \times S^1} u^*d\alpha$ cannot be
any larger than the period of its positive asymptotic orbit, which does
not depend on~$N$.  So this gives a contradiction, proving that $u$
cannot touch the region $\{\rho \ge 1/2\}$.  The mirror image of this
argument shows that $u$ also cannot touch the region $\{ \rho \le -1/2\}$.
\end{proof}

With Lemma~\ref{lemma:cannotEscape} in hand, the calculation of
$HC_*^h(\TT^3_N,\alpha,J)$ for sufficiently large $N$ is straightforward: 
there is one odd generator
and one even generator, with a trivial differential, giving
$$
HC_*^h(\TT^3,\xi_1) \cong \begin{cases}
\ZZ_2 & * = \text{odd},\\
\ZZ_2 & * = \text{even}.
\end{cases}
$$
This calculation can now be extended to $(\TT^3,\xi_k)$ by a cheap trick:
using the contactomorphism $(\TT^3_N,\ker\alpha) \to (\TT^3,\xi_1)$,
let us identify $\TT^3_N$ with $\TT^3$ and write $\alpha = F \alpha_1$
for some function $F : \TT^3 \to (0,\infty)$.  Then the $k$-fold covering
map
$$
\Phi_k : \TT^3 \to \TT^3 : (\rho,\phi,\theta) \mapsto (k\rho,\phi,\theta)
$$
maps the homotopy class $h$ to itself and pulls back $\xi_1$ to~$\xi_k$,
so $\Phi_k^*\alpha$ is a contact form for~$\xi_k$.  It is also
$h$-admissible: indeed, $\Phi_k^*\alpha$ admits no contractible orbits
since they would project down to contractible orbits on $(\TT^3,\alpha)$,
and every orbit in $\pP_h(\Phi_k^*\alpha)$ projects to one in
$\pP_h(\alpha)$, hence they are all nondegenerate.  The almost complex
structure $\Phi_k^*J \in \jJ(\Phi_k^*\alpha)$ then makes the map
$\Id \times \Phi_k : (\RR \times \TT^3,\Phi_k^*J) \to (\RR \times \TT^3,J)$
holomorphic, so every $\Phi_k^*J$-holomorphic cylinder counted by
$HC_*^h(\TT^3,\Phi_k^*\alpha,\Phi_k^*J)$ projects to
a $J$-holomorphic cylinder counted by $HC_*^h(\TT^3,\alpha,J)$, and 
conversely, each orbit in $\pP_h(\alpha)$ and each $J$-holomorphic cylinder 
has exactly $k$ lifts to the cover.  The generators of
$CC_*^h(\TT^3,\Phi_k^*\alpha)$ thus consist of $2k$ orbits, $k$ odd and
$k$ even, with $2k$ connecting $\Phi_k^*J$-holomorphic cylinders
that cancel each other in pairs, giving a trivial differential.
In summary:
$$
HC_*^h(\TT^3,\xi_k) \cong \begin{cases}
\ZZ_2^k & * = \text{odd},\\
\ZZ_2^k & * = \text{even}.
\end{cases}
$$

\psfrag{uA}{$u_A$}
\psfrag{uB}{$u_B$}
\psfrag{uC}{$u_C$}
\psfrag{u}{$u$}
\psfrag{gminus}{$\boldsymbol{\gamma}^-$}
\psfrag{g1}{$\gamma_1$}
\psfrag{g2}{$\gamma_2$}
\psfrag{g3}{$\gamma_3$}
\psfrag{g4}{$\gamma_4$}
\psfrag{g5}{$\gamma_5$}
\psfrag{What}{$\widehat{W}$}
\psfrag{RtimesM+}{$\RR \times M_+$}

\chapter{Coherent orientations}
\label{lec:orientations}

\minitoc

\vspace{12pt}

\section{Gluing maps and coherence}
\label{sec:coherence}

This lecture will be concerned with orienting the moduli spaces
$$
\mM(J) := \mM_{g,m}(J,A,\boldsymbol{\gamma}^+,\boldsymbol{\gamma}^-)
$$ 
of $J$-holomorphic curves in a completed symplectic cobordism~$\widehat{W}$,
in cases where they are smooth.  We assume as usual that all Reeb orbits
are nondegenerate so that the usual linearized Cauchy-Riemann operators
are Fredholm.

For SFT and other Floer-type theories,
it is not enough to know that each component of $\mM(J)$ is 
orientable---relations like $\p^2 = 0$ rely on having certain compatibility
conditions between the orientations on different components.  The point is
that whenever a space of broken curves is meant to be interpreted as the
boundary of some other compactified moduli space, we need to make sure that
it carries the boundary orientation.  This compatibility is what is known
as \emph{coherence}, and in order to define it properly, we need to
return to the subject of gluing.

Our discussion of gluing in Lecture~\ref{lec:tight3tori} was fairly simple because it was
limited to somewhere injective holomorphic cylinders that could only
break along simply covered Reeb orbits.  Recall however that more general
holomorphic buildings carry a certain amount of extra structure that was
not relevant in that simple case.  Even in a building $\mathbf{u}$
that has only two
nontrivial levels $u_-$ and $u_+$, the breaking punctures carry 
\emph{decorations}: i.e.~if $\{z^+,z^-\}$ is a breaking pair in~$\mathbf{u}$,
then the decoration defines an orientation-reversing orthogonal map
$$
\delta_{z^+} \stackrel{\Phi}{\longrightarrow} \delta_{z^-}
$$
between the two ``circles at infinity'' $\delta_{z^\pm}$ associated to the
punctures~$z^\pm$ (see \S\ref{sec:DM}).  This extra information
is uniquely determined if the breaking orbit is simply covered, but at a
multiply covered breaking orbit there is ambiguity, and the decoration cannot
be deduced from knowledge of $u_-$ and $u_+$ alone.  We therefore need to
consider moduli spaces of curves with a bit of extra structure.

For each Reeb orbit $\gamma$ in $M_+$ or $M_-$, choose a point on its image
$$
p_\gamma \in \im \gamma \subset M_\pm.
$$
For a $J$-holomorphic curve $u : (\dot{\Sigma} = \Sigma \setminus
(\Gamma^+ \cup \Gamma^-),j) \to (\widehat{W},J)$ with
a puncture $z \in \Gamma^\pm$ asymptotic to~$\gamma$, an \defin{asymptotic
marker} is a choice of a ray $\ell \subset T_z\Sigma$ such that
$$
\lim_{t \to 0^+} u(c(t)) = (\pm\infty,p_\gamma)
$$
for any smooth path $c(t) \in \Sigma$ with $c(0) = z$ and $0 \ne \dot{c}(0) \in \ell$.
If $\gamma$ has covering multiplicity $m \in \NN$, then there are exactly
$m$ choices of asymptotic markers at~$z$, related to each other by the action
on $T_z\Sigma$ by the $m$th roots of unity.  We shall denote
$$
\mM^\$(J) := \mM_{g,m}^\$(J,A,\boldsymbol{\gamma}^+,\boldsymbol{\gamma}^-) :=
\left\{ (\Sigma,j,\Gamma^+,\Gamma^-,\Theta,u,\ell) \right\} \big/ \sim,
$$
where $(\Sigma,j,\Gamma^+,\Gamma^-,\Theta,u)$ represents an element of
$\mM_{g,m}(J,A,\boldsymbol{\gamma}^+,\boldsymbol{\gamma}^-)$,
$\ell$ denotes an assignment of asymptotic markers to every puncture $z \in \Gamma^\pm$,
and
$$
(\Sigma_0,j_0,\Gamma^+_0,\Gamma^-_0,\Theta_0,u_0,\ell_0) \sim
(\Sigma_1,j_1,\Gamma^+_1,\Gamma^-_1,\Theta_1,u_1,\ell_1)
$$
means the existence of a biholomorphic map $\psi : (\Sigma_0,j_0) \to 
(\Sigma_1,j_1)$ which defines an equivalence of
$(\Sigma_0,j_0,\Gamma^+_0,\Gamma^-_0,\Theta_0,u_0)$ with
$(\Sigma_1,j_1,\Gamma^+_1,\Gamma^-_1,\Theta_1,u_1)$ and satisfies
$\psi_*\ell_0 = \ell_1$.  There is a natural surjection
$$
\mM^\$(J) \to \mM(J)
$$
defined by forgetting the markers.  We will say that an element
$u \in \mM^\$(J)$ is Fredholm regular whenever its image under the map
to $\mM(J)$ is regular.  Let
$$
\mM^{\$,\reg}(J) = \mM^{\$,\reg}_{g,m}(J,A,\boldsymbol{\gamma}^+,\boldsymbol{\gamma}^-)
\subset \mM^\$(J)
$$
denote the open subset consisting of Fredholm regular curves with asymptotic
markers.  Note that components of $\mM(J)$ and $\mM^\$(J)$ consisting of
closed curves are identical spaces; components with punctures have the
following simple relationship to each other.

\begin{prop}
\label{prop:markerStructure}
Each component of $\mM^{\$,\reg}(J)$ consisting of curves with at least one
puncture admits the structure of a smooth
manifold, whose dimension on each connected component matches that of
$\mM^\reg(J)$.  Moreover, the natural map
$$
\mM^{\$,\reg}(J) \to \mM^\reg(J)
$$
is smooth, and the preimage of a curve $u \in \mM^\reg(J)$ with asymptotic
orbits $\{\gamma_z\}_{z \in \Gamma}$ of covering multiplicities
$\{\kappa_z\}_{z \in \Gamma}$ contains exactly
$$
\frac{\prod_{z \in \Gamma} \kappa_z}{|\Aut(u)|}
$$
distinct elements.
\end{prop}
\begin{proof}
The smooth structure of $\mM^{\$,\reg}(J)$ arises from the same argument we
used in Lecture~\ref{lec:transversality} for $\mM^\reg(J)$, supplemented by the following remarks:
first, every nontrivial automorphism $\psi \in \Aut(u)$ for
$u \in \mM(J)$ acts nontrivially on the asymptotic markers.  Indeed,
$\psi$ is required to fix each of the punctures and is a biholomorphic map
with $\psi^k \equiv \Id$ for some $k \in \NN$, thus it takes the form
$z \mapsto e^{2\pi i m / k}$ in suitable holomorphic coordinates near
each puncture for suitable integers $m , k \in \ZZ$.  If $m=0$, then
unique continuation implies $\psi \equiv \Id$, and otherwise $\psi$ changes
the asymptotic marker at every puncture.  With this understood, one can
define as in \S\ref{sec:IFT} a local identification of
$\mM^\$(J)$ with $\dbar_J^{-1}(0) / \Aut(\Sigma,j_0,\Gamma \cup \Theta)$,
where $\dbar_J^{-1}(0)$ includes information about asymptotic markers and
is a smooth manifold by the implicit function theorem, but 
$\Aut(\Sigma,j_0,\Gamma \cup \Theta)$ acts on it \emph{freely}, producing
a quotient with no isotropy.

Finally, if $(\Sigma,j,\Gamma \cup \Theta,u)$ represents an element of
$\mM(J)$ with asymptotic orbits $\{ \gamma_z \}_{z \in \Gamma}$, then
the number of possible choices of asymptotic markers is precisely
$\prod_{z \in \Gamma} \kappa_z$.  However, not all of these produce
inequivalent elements of $\mM^\$(J)$: indeed, the previous paragraph shows
that $\Aut(u)$ acts freely on the set of all choices of markers, so that
the total number of inequivalent choices is as stated.
\end{proof}

Suppose $u_+$ and $u_-$ are two (possibly disconnected and/or nodal) holomorphic 
curves, with asymptotic markers, such that the number of negative punctures
of $u_+$ equals the number of positive punctures of $u_-$, and
the asymptotic orbit of $u_+$ at its $i$th negative puncture matches
that of $u_-$ at its $i$th positive puncture for every~$i$.  Then
the pair $(u_-,u_+)$
naturally determines a holomorphic building: indeed, the breaking punctures
admit unique decorations determined by identifying the markers on $u_+$
with the markers at corresponding punctures of~$u_-$.

\begin{figure}
\includegraphics{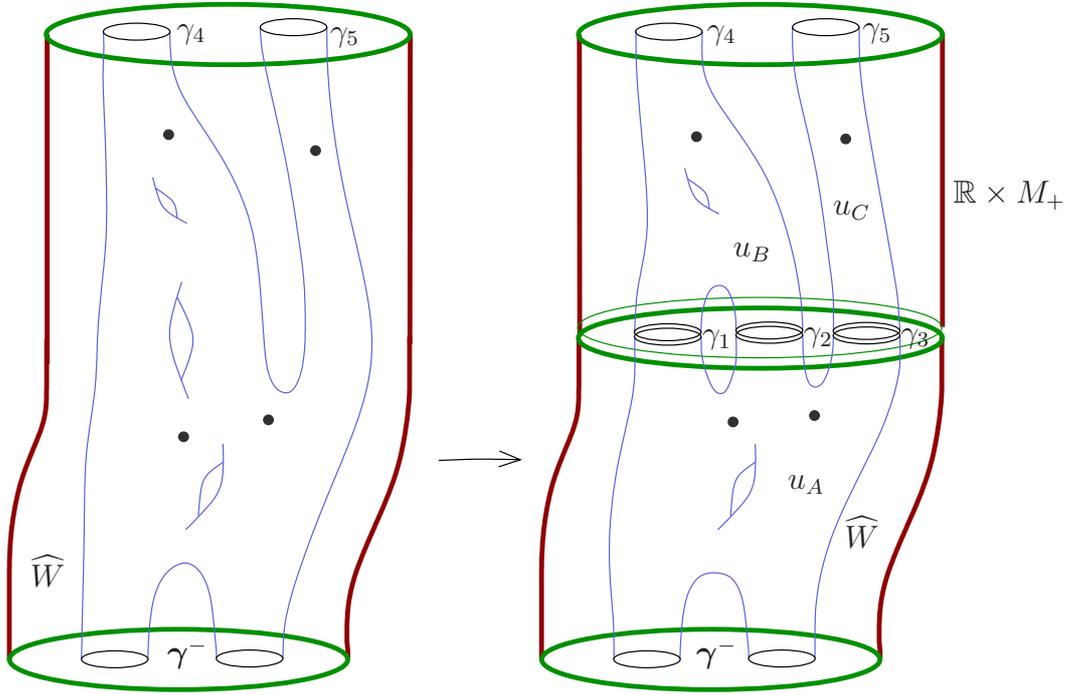}
\caption{\label{fig:gluing} The degeneration scenario behind the gluing
map \eqref{eqn:gluingMap}}
\end{figure}

Let us now consider a concrete example of a gluing scenario.
Figure~\ref{fig:gluing} shows the degeneration of a sequence of curves
in $\mM_{3,4}(J,A_k,(\gamma_4,\gamma_5),\boldsymbol{\gamma}^-)$ to a
building $\mathbf{u} \in \overline{\mM}_{3,4}(J,A+B+C,(\gamma_4,\gamma_5),
\boldsymbol{\gamma}^-)$ with one main level and one upper level.  The main 
level is a connected curve $u_A \in \mM_{1,2}(J,A,(\gamma_1,\gamma_2,\gamma_3),
\boldsymbol{\gamma}^-)$, and the upper level consists of two
connected curves
$$
u_B \in \mM_{1,1}(J_+,B,\gamma_4,(\gamma_1,\gamma_2)), \qquad
u_C \in \mM_{0,1}(J_+,C,\gamma_5,\gamma_3).
$$
One can endow each of these curves with asymptotic markers compatible
with the decoration of~$\mathbf{u}$; this is a non-unique choice, but
e.g.~if one chooses markers for $u_A$ arbitrarily, then the markers at the
negative punctures of $u_B$ and $u_C$ are uniquely determined.  Now if
all three curves are Fredholm regular, then a substantial generalization of the
gluing procedure outlined in Lecture~\ref{lec:tight3tori} provides open neighborhoods
$\uU^\$_A$ and $\uU^\$_{BC}$,
\begin{equation*}
\begin{split}
u_A \in \uU^\$_A &\subset \mM_{1,2}^\$(J,A,(\gamma_1,\gamma_2,\gamma_3)), \\
[(u_B,u_c)] \in \uU^\$_{BC} &\subset 
\left(\mM_{1,1}^\$(J_+,B,\gamma_4,(\gamma_1,\gamma_2)) \times \mM_{0,1}^\$(J_+,C,\gamma_5,\gamma_3)\right)
\Big/ \RR
\end{split}
\end{equation*}
which are smooth manifolds of dimensions 
\begin{equation*}
\begin{split}
\dim \uU^\$_A &= \virdim \mM_{1,2}(J,A,(\gamma_1,\gamma_2,\gamma_3)), \\
\dim \uU^\$_{BC} &= \virdim \mM_{1,1}(J_+,B,\gamma_4,(\gamma_1,\gamma_2))
+ \virdim \mM_{0,1}(J_+,C,\gamma_5,\gamma_3) - 1,
\end{split}
\end{equation*}
along with a smooth embedding
\begin{equation}
\label{eqn:gluingMap}
\Psi : [R_0,\infty) \times \uU^\$_A \times \uU^\$_{BC} \hookrightarrow
\mM_{3,4}^\$(J,A+B+C,(\gamma_4,\gamma_5),\boldsymbol{\gamma}^-),
\end{equation}
defined for $R_0 \gg 1$.  This is an example of a \defin{gluing map}: it
has the property that for any $u \in \uU^\$_A$ and $v \in \uU^\$_{BC}$,
$\Psi(R,u,v)$ converges in the SFT topology as $R \to \infty$ to the
unique building (with asymptotic markers) having main level $u$ and
upper level~$v$, and moreover, every sequence of smooth curves
degenerating in this way is eventually in the image of~$\Psi$.

In analogous ways one can define gluing maps for buildings with a main level
and a lower level, or more than two levels, or multiple levels in a
symplectization (always dividing symplectization levels by the
$\RR$-action).  It's important to notice that in all such scenarios, the
domain and target of the gluing map have the same dimension,
e.g.~the dimension of both sides of \eqref{eqn:gluingMap} is the sum of the
virtual dimensions of the three moduli spaces concerned.

\begin{defn}
\label{defn:coherence}
A set of orientations for the connected components of $\mM^\$(J)$ and
$\mM^\$(J_\pm)$ is called \defin{coherent} if all gluing maps are
orientation preserving.
\end{defn}

Stated in this way, this definition is based on the pretense that we 
never have to worry about non-regular curves in any components
of $\mM^\$(J)$, and that is of course false---sometimes regularity cannot
be achieved, in particular for multiply covered curves.  
As we'll see though in \S\ref{sec:det}, the question of 
orientations can be reframed in a way that completely disjoins it from the
question of regularity, thus we will later be able to state a more
general version of the above definition that is independent of
regularity (see Definition~\ref{defn:coherence2}).
The main result whose proof we will outline in the next few sections is then:

\begin{thm}
\label{thm:coherent}
Coherent orientations exist.
\end{thm}

But there is also some bad news.  The space $\mM^\$(J)$ with asymptotic
markers is not actually the space we want to orient.  In fact, even the
usual moduli space $\mM(J)$ has a certain amount of extra information
in it that we'd rather not keep track of when we don't have to, 
for instance the ordering of the punctures.  Can we forget this information
without forgetting the orientation of the moduli space?  Not always:

\begin{prop}
\label{prop:permute}
Suppose $\hat{\boldsymbol{\gamma}}^+ = (\gamma_1^+,\ldots,\gamma_{k_+}^+)$, 
and $\check{\boldsymbol{\gamma}}^+$ is a similar ordered list of Reeb
orbits obtained from $\hat{\boldsymbol{\gamma}}^+$ by exchanging
$\gamma_j^+$ with $\gamma_k^+$ for some $1 \le j < k \le k_+$.
Then for any choice of coherent orientations, the natural map
$$
\mM^\$_{g,m}(J,A,\hat{\boldsymbol{\gamma}}^+,\boldsymbol{\gamma}^-) \to
\mM^\$_{g,m}(J,A,\check{\boldsymbol{\gamma}}^+,\boldsymbol{\gamma}^-)
$$
defined by permuting the corresponding punctures $z_j^+,z_k^+ \in \Gamma^+$
along with their asymptotic markers is orientation reversing if and only if
the numbers
$$
n - 3 + \muCZ(\gamma_i^+)
$$
for $i=j,k$ are both odd.  A similar statement holds for permutations of
negative punctures.
\end{prop}

This result is the reason for the super-commutative algebra that we will
see in the next lecture.
What about forgetting the markers?  It turns out that we can sometimes
do that as well, but again not always.

\begin{prop}
\label{prop:bad}
Suppose $\mM^\$_{g,m}(J,A,\boldsymbol{\gamma}^+,\boldsymbol{\gamma}^-) \to
\mM^\$_{g,m}(J,A,\boldsymbol{\gamma}^+,\boldsymbol{\gamma}^-)$ is the map
defined by multiplying the asymptotic marker by $e^{2\pi i / m}$ at
one of the punctures for which the asymptotic orbit is an $m$-fold cover
$\gamma^m$ of a simple orbit~$\gamma$.  For any choice of coherent
orientations, this map reverse orientation if and only if $m$ is even
and $\muCZ(\gamma^m) - \muCZ(\gamma)$ is odd.
\end{prop}

Note that in both of the above propositions, only the odd/even parity
of the Conley-Zehnder indices matters, so there is no need to choose
trivializations.  Proposition~\ref{prop:bad} motivates one of the more
mysterious technical definitions in SFT.

\begin{defn}
\label{defn:bad}
A closed nondegenerate Reeb orbit $\gamma$ is called a \defin{bad} orbit
if it is an $m$-fold cover of some simple orbit $\gamma'$ where $m$ is
even and $\muCZ(\gamma) - \muCZ(\gamma')$ is odd.  Orbits that are not bad
are called \defin{good}.
\end{defn}

The upshot is that coherent orientations can be defined on the union of
all components $\mM_{g,m}(J,A,\boldsymbol{\gamma}^+,\boldsymbol{\gamma}^-)$
for which all of the orbits in the lists 
$\boldsymbol{\gamma}^+$ and $\boldsymbol{\gamma}^-$ are good.
This does not mean that moduli spaces involving bad orbits cannot be
dealt with---in fact, such moduli spaces have the convenient property that
the number of distinct choices of asymptotic markers is always even,
and every such choice can be cancelled by an alternative choice that
induces the opposite orientation.  For this reason, while bad orbits
certainly can appear in breaking of holomorphic curves, we will see that
they do not need to serve as generators of SFT.

\section{Permutations of punctures and bad orbits}
\label{sec:funPart}

Before addressing the actual construction of coherent orientations, we can
already give heuristic proofs of Propositions~\ref{prop:permute} and~\ref{prop:bad}.
They are not fully rigorous because they are based on the same pretense as
Definition~\ref{defn:coherence}, namely that all curves we ever have to worry about
(including multiple covers) are regular.  But we will be able to turn these
into precise arguments in~\S\ref{sec:funPart2}, after discussing the
determinant line bundle.

\begin{figure}
\includegraphics{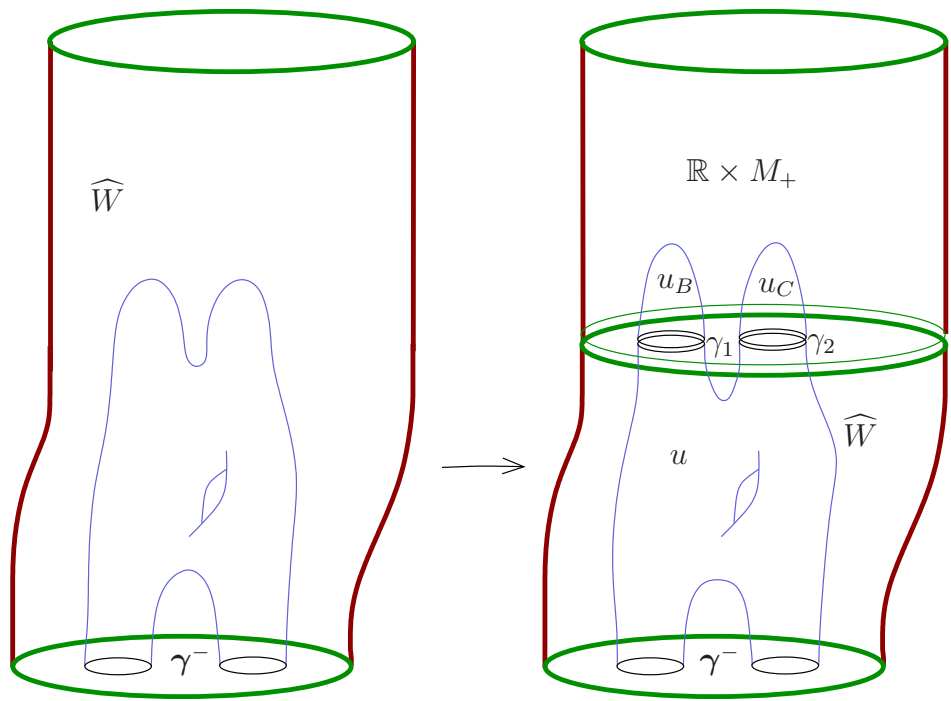}
\caption{\label{fig:permute} The gluing thought-experiment used for proving
Propositions~\ref{prop:permute} and~\ref{prop:bad}.}
\end{figure}

\begin{proof}[Heuristic proof of Proposition~\ref{prop:permute}]
To simplify the notation, suppose $\hat{\boldsymbol{\gamma}}^+$ consists of
only two orbits, so $\hat{\boldsymbol{\gamma}}^+ = (\gamma_1,\gamma_2)$ and
$\check{\boldsymbol{\gamma}}^+ = (\gamma_2,\gamma_1)$.  Consider the
gluing scenario shown in Figure~\ref{fig:permute}, where
$u \in \mM^\$_{g,m}(J,A,(\gamma_1,\gamma_2),\boldsymbol{\gamma}^-)$
needs to be glued to a disjoint union of two planes
$$
u_B \in \mM^\$_{0,0}(J_+,B,\emptyset,\gamma_1), \qquad
u_C \in \mM^\$_{0,0}(J_+,C,\emptyset,\gamma_2).
$$
You might object that there's no guarantee that such planes must exist
in $\RR \times M_+$, e.g.~the orbits $\gamma_1$ and $\gamma_2$ might not
even be contractible.  This concern is valid so far as it goes, but it
misses the point: since we're talking about gluing rather than
compactness, we do not need any seriously global information about 
$\widehat{W}$ and $M_+$, as the gluing process doesn't
depend on anything outside a small neighborhood of the curves we're
considering.  Thus we are free to change the global structure of $M_+$ elsewhere
so that the planes $u_B$ and $u_C$ will exist.\footnote{Of course by the
maximum principle, planes
with only negative ends will not exist in $\RR \times M_+$ if this is the
symplectization of a contact manifold.  But we could also change the
contact data to a stable Hamiltonian structure for which such planes are allowed.}
If you still can't
imagine how one might do this, try not to worry about it and just think of
Figure~\ref{fig:permute} as a thought-experiment: it's a situation that
certainly does sometimes happen, so when it does, let's see what it implies
about orientations.

Assuming all three curves in the picture are regular, there will be smooth
open neighborhoods
\begin{equation*}
\begin{split}
u \in \uU_{12} &\subset \mM^\$_{g,m}(J,A,(\gamma_1,\gamma_2),\boldsymbol{\gamma}^-) \\
[(u_B,u_C)] \in \uU_{BC} &\subset \left( 
\mM^\$_{0,0}(J_+,B,\emptyset,\gamma_1) \times \mM^\$_{0,0}(J_+,C,\emptyset,\gamma_2)
\right) \Big/ \RR
\end{split}
\end{equation*}
and a gluing map
$$
\Psi_{BC} : [R_0,\infty) \times \uU_{12} \times \uU_{BC} \hookrightarrow
\mM_{g,m}^\$(J,A+B+C,\emptyset,\boldsymbol{\gamma}^-),
$$
which must be orientation preserving by assumption.  But reversing the
order of the product 
$\mM^\$_{0,0}(J_+,B,\emptyset,\gamma_1) \times \mM^\$_{0,0}(J_+,C,\emptyset,\gamma_2)$
and letting $u' \in \mM^\$_{g,m}(J,A,(\gamma_2,\gamma_1),\boldsymbol{\gamma}^-)$
denote the image of $u$ under the map that switches the order of its
positive punctures, there are also smooth open neighborhoods
\begin{equation*}
\begin{split}
u' \in \uU_{21} &\subset \mM^\$_{g,m}(J,A,(\gamma_2,\gamma_1),\boldsymbol{\gamma}^-) \\
[(u_C,u_B)] \in \uU_{CB} &\subset \left( 
\mM^\$_{0,0}(J_+,C,\emptyset,\gamma_2) \times \mM^\$_{0,0}(J_+,B,\emptyset,\gamma_1) 
\right) \Big/ \RR
\end{split}
\end{equation*}
and a gluing map
$$
\Psi_{CB} : [R_0,\infty) \times \uU_{21} \times \uU_{CB} \hookrightarrow
\mM_{g,m}^\$(J,A+B+C,\emptyset,\boldsymbol{\gamma}^-).
$$
If both of these gluing maps preserve orientation, then the effect 
on orientations of the
map from $\mM^\$_{g,m}(J,A,(\gamma_1,\gamma_2),\boldsymbol{\gamma}^-)$
to $\mM^\$_{g,m}(J,A,(\gamma_2,\gamma_1),\boldsymbol{\gamma}^-)$ defined by
interchanging the positive punctures must be the same as that of the map
\begin{equation*}
\begin{split}
\mM^\$_{0,0}(J_+,B,\emptyset,\gamma_1) \times \mM^\$_{0,0}(J_+,C,\emptyset,\gamma_2) &\to
\mM^\$_{0,0}(J_+,C,\emptyset,\gamma_2) \times \mM^\$_{0,0}(J_+,B,\emptyset,\gamma_1) \\
(u_B,u_C) &\mapsto (u_C,u_B).
\end{split}
\end{equation*}
The latter is orientation reversing if and only if both moduli spaces of
planes are odd dimensional, which means $n-3 + \muCZ(\gamma_i)$ is odd
for $i=1,2$.
\end{proof}

\begin{proof}[Heuristic proof of Proposition~\ref{prop:bad}]
Let us reuse the thought-experiment of Figure~\ref{fig:permute}, but with
different details in focus.  Suppose $\gamma_1$ in the picture is an
$m$-fold covered orbit $\gamma^m$, where $\gamma$ is simply covered,
and suppose that $u_B$ is also an $m$-fold cover, taking the form
$$
u_B(z) = v(z^m)
$$
for a somewhere injective plane $v \in \mM_{0,0}(J_+,B_0,\emptyset,\gamma)$.
We're going to assume again that all curves in the discussion are regular,
including the multiple cover~$u_B$; while this doesn't sound very plausible,
we will see once the determinant line bundle enters the picture 
in \S\ref{sec:det} that it is an irrelevant detail.  Now, $u_B$ has a cyclic
automorphism group 
$$
\Aut(u_B) = \ZZ_m \subset \U(1)
$$
which acts freely on
the set of $m$ choices of asymptotic marker for~$u_B$.  Then if we act
with the same element of $\ZZ_m$ on $u_B$ and on the corresponding asymptotic 
marker for~$u$, the building is unchanged, as it has the same decoration.
Coherence therefore implies that the effect on orientations of the map from
$\mM^\$_{g,m}(J,A,(\gamma_1,\gamma_2),\boldsymbol{\gamma}^-)$ to itself
defined by acting with the canonical generator of $\ZZ_m \subset \U(1)$
on the marker at $\gamma_1$ is the same as the effect of the map
$\mM_{0,0}^\$(J_+,m B_0,\emptyset,\gamma^m) \to \mM_{0,0}^\$(J_+,m B_0,\emptyset,\gamma^m)$
defined by composing $u_B : \CC \to \RR \times M_+$ with 
$\psi(z) := e^{2\pi i/m} z$.  

The derivative of this map from $\mM_{0,0}^\$(J_+,m B_0,\emptyset,\gamma^m)$
to itself at $u_B$ defines a linear self-map
$$
\Psi : T_{u_B} \mM_{0,0}(J_+,m B_0,\emptyset,\gamma^m) \to
T_{u_B} \mM_{0,0}(J_+,m B_0,\emptyset,\gamma^m)
$$
with $\Psi^m = \1$.  The latter implies that $\Psi$ cannot reverse
orientation if $m$ is odd.  If $m$ is even, observe that the representation 
theory of $\ZZ_m$ gives a decomposition
$$
T_{u_B} \mM_{0,0}(J_+,m B_0,\emptyset,\gamma^m) 
= V_1 \oplus V_{-1} \oplus V_{\operatorname{rot}},
$$
where $\Psi$ acts on $V_{\pm 1}$ as $\pm \1$, and $V_{\operatorname{rot}}$
is a direct sum of real $2$-dimensional subspaces on which $\Psi$ acts
by rotations (and therefore preserves orientations).  Thus $\Psi$ reverses
the orientation of $T_{u_B} \mM_{0,0}(J_+,m B_0,\emptyset,\gamma^m)$ 
if and only if
$\dim V_{-1}$ is odd.  As we will review in the next section,
$T_{u_B} \mM_{0,0}(J_+,m B_0,\emptyset,\gamma^m)$ is a space of holomorphic
sections of $u_B^*T(\RR \times M_+)$ modulo a subspace defined via the
linearized automorphisms of~$\CC$, so $V_1$ consists of precisely those
sections $\eta$ that satisfy $\eta = \eta \circ \psi$, meaning they are
$m$-fold covers of sections of $v^*T(\RR \times M_+)$.  This defines a
bijective correspondence between $V_1$ and 
$T_v \mM_{0,0}(J_+,B_0,\emptyset,\gamma)$, so
$$
\dim V_{-1} = \dim \mM_{0,0}(J_+,m B_0,\emptyset,\gamma^m) -
\dim \mM_{0,0}(J_+,B_0,\emptyset,\gamma) \qquad \text{(mod 2)}.
$$
The result then comes from plugging in the dimension formulas for these
two moduli spaces.
\end{proof}

\section{Orienting moduli spaces in general}
\label{sec:moduli}

We now discuss concretely what is involved in orienting a moduli
space of $J$-holomorphic curves.

Recall from Lecture~\ref{lec:transversality} that whenever a
curve $u_0 : (\dot{\Sigma} = \Sigma \setminus \Gamma,j_0) \to (\widehat{W},J)$ 
with marked points $\Theta \subset \dot{\Sigma}$ is Fredholm regular,
a neighborhood of $u_0$ in $\mM(J)$ can be identified with
$$
\dbar_J^{-1}(0) \big/ G_0,
$$
where $G_0 = \Aut(\Sigma,j_0,\Gamma \cup \Theta)$ and
$\dbar_J$ is the smooth Fredholm section
$$
\tT \times \bB^{k,p,\delta} \to \eE^{k-1,p,\delta} : (j,u) \mapsto
Tu + J \circ Tu \circ j,
$$
defined on the product of a $G_0$-invariant Teichm\"uller slice $\tT$ 
through $j_0$ with
a Banach manifold $\bB^{k,p,\delta}$ of $W^{k,p}$-smooth maps
$\dot{\Sigma} \to \widehat{W}$ satisfying an exponential decay condition
at the cylindrical ends.  Here $G_0$ acts on $\dbar_J^{-1}(0)$ by
\begin{equation}
\label{eqn:AutAction}
G_0 \times \dbar_J^{-1}(0) \to \dbar_J^{-1}(0) : (\varphi,(j,u)) \mapsto
(\varphi^*j , u \circ \varphi).
\end{equation}
Regularity means that the linearization
$D\dbar_J(j_0,u_0) : T_{j_0}\tT \oplus T_{u_0} \bB^{k,p,\delta} \to
\eE^{k-1,p,\delta}_{(j_0,u_0)}$ is surjective, and the implicit function
theorem then gives a natural identification
$$
T_{u_0} \mM(J) = \ker D\dbar_J(j_0,u_0) \big/ \aut(\Sigma,j_0,\Gamma \cup \Theta),
$$
where $\aut(\Sigma,j_0,\Gamma \cup \Theta)$ denotes the Lie algebra of
$G_0$, which acts on $\ker D\dbar_J(j_0,u_0)$ by differentiating
\eqref{eqn:AutAction}.\footnote{The presence of $\aut(\Sigma,j_0,\Gamma \cup \Theta)$
in this discussion is only relevant in the finite set of ``non-stable'' cases
where $\chi(\dot{\Sigma} \setminus \Theta) \ge 0$, since otherwise
$G_0$ is finite and thus $\aut(\Sigma,j_0,\Gamma \cup \Theta)$ is trivial.}
This action actually defines an \emph{inclusion} of 
$\aut(\Sigma,j_0,\Gamma \cup \Theta)$ into $\ker D\dbar_J(j_0,u_0)$ whenever
$u_0$ is not constant, thus we can regard $\aut(\Sigma,j_0,\Gamma \cup \Theta)$
as a subspace of $\ker D\dbar_J(j_0,u_0)$.  

As outlined in Proposition~\ref{prop:markerStructure}, the space
$\mM^\$(J)$ with asymptotic markers admits a similar local description:
here one only needs to enhance the structure of the Banach manifold
$\bB^{k,p,\delta}$ with information about asymptotic markers at each
puncture, so the Banach manifold needed to describe $\mM^\$(J)$ is a
finite covering space of~$\bB^{k,p,\delta}$.  The rest of the discussion
is identical, except for the fact that when markers are included,
$G_0$ always acts freely on $\dbar_J^{-1}(0)$.

We now make a useful observation about the spaces
$\aut(\Sigma,j_0,\Gamma \cup \Theta)$ and $T_{j_0}\tT$: namely, they both
carry natural complex structures and are thus canonically oriented.
This follows from the fact that both the automorphism group 
$G_0$ and
the Teichm\"uller space $\tT(\Sigma,\Gamma \cup \Theta) = \jJ(\Sigma) \big/
\Diff_0(\Sigma,\Gamma \cup \Theta)$ are naturally
complex manifolds.  On the linearized level, one way to see it is via the
fact---mentioned previously in \S\ref{sec:moduliSpaces}---that 
$\aut(\Sigma,j_0,\Gamma \cup \Theta)$ and $T_{[j_0]} \tT(\Sigma,\Gamma \cup \Theta)$
can be naturally identified with the kernel and cokernel respectively of
the natural linear Cauchy-Riemann type operator on~$(\Sigma,j_0)$,
\begin{equation}
\label{eqn:linSigma}
\mathbf{D}_{\Id} : W^{k,p}_{\Gamma \cup \Theta}(T\Sigma) \to W^{k-1,p}(\overline{\End}_\CC(T\Sigma)),
\end{equation}
which is the linearization at $\Id$ of the nonlinear operator that detects
holomorphic maps $(\Sigma,j_0) \to (\Sigma,j_0)$.  This operator is equivalent
to the operator that defines the holomorphic structure of $T\Sigma$, thus
it is complex linear.  To handle the punctures and marked points, one needs
to restrict the nonlinear operator to the space of $W^{k,p}$-smooth maps
$\Sigma \to \Sigma$ that fix every point in $\Gamma \cup \Theta$, thus the
domain of the linearization becomes the finite-codimensional subspace
$$
W^{k,p}_{\Gamma \cup \Theta}(T\Sigma) := \left\{ X \in W^{k,p}(T\Sigma)\ \big|\ 
X|_{\Gamma \cup \Theta} = 0 \right\}.
$$
This subspace is still complex, thus so is \eqref{eqn:linSigma}, and its
kernel and cokernel inherit natural complex structures.

The complex structure on $\aut(\Sigma,j_0,\Gamma \cup \Theta)$ means that
defining an orientation on the tangent space $T_{u_0} \mM^\$(J)$ is equivalent 
to defining one on $\ker D\dbar_J(j_0,u_0)$.  The latter operator takes
the form
$$
D\dbar_J(j_0,u_0) : T_{j_0}\tT \oplus T_{u_0} \bB^{k,p,\delta} \to
\eE^{k-1,p,\delta}_{(j_0,u_0)} : (y,\eta) \mapsto J \circ T u_0 \circ y +
\mathbf{D}_{u_0} \eta,
$$
where $\mathbf{D}_{u_0} : W^{k,p,\delta}(u_0^*T\widehat{W}) \oplus V_\Gamma \to
W^{k-1,p,\delta}(\overline{\Hom}_\CC(T\dot{\Sigma},u_0^*T\widehat{W}))$
is the usual linearized Cauchy-Riemann operator at~$u_0$, with
$V_\Gamma$ denoting a complex $(\#\Gamma)$-dimensional space of smooth sections
that are constant near infinity..  The remarks
above and the fact that $u_0$ is $J$-holomorphic imply that the first 
term in this operator,
$$
T_{j_0} \tT \to \eE^{k-1,p,\delta}_{(j_0,u_0)} : y \mapsto
J \circ T u_0 \circ y
$$
is a complex-linear map.  Now if $\mathbf{D}_{u_0}$ happens also to be
a complex-linear map, then we are done, because $\ker D\dbar_J(j_0,u_0)$
will then be a complex vector space and inherit a natural orientation.

In general, $\mathbf{D}_{u_0}$ is not complex linear, though it does have
a \emph{complex-linear part},
$$
\mathbf{D}_{u_0}^\CC \eta := \frac{1}{2} \left( \mathbf{D}_{u_0} \eta - 
J \mathbf{D}_{u_0}(J \eta) \right),
$$
which is also a Cauchy-Riemann type operator.  The space of all
Cauchy-Riemann type operators on a fixed vector bundle is affine, so
one can interpolate from $\mathbf{D}_{u_0}$ to $\mathbf{D}_{u_0}^\CC$
through a path of Cauchy-Riemann type operators, though they may not
all be Fredholm---this depends on the asymptotic operators at the punctures.
In the special case however where there are no punctures, one can easily
imagine making use of this idea: if $\dot{\Sigma} = \Sigma$ is a
closed surface, then the obvious homotopy from $\mathbf{D}_{u_0}$ to its
complex-linear part yields a homotopy from $D\dbar_J(j_0,u_0)$ to its
complex-linear part, and if every operator along this homotopy happens to
be surjective, then the canonical orientation defined on the kernel of
the complex-linear operator determines an orientation on
$\ker D\dbar_J(j_0,u_0)$.

There are two obvious problems with the above discussion:
\begin{enumerate}
\item We have no way to ensure that every operator in the homotopy from
$D\dbar_J(j_0,u_0)$ to its complex-linear part is surjective;
\item If there are punctures, then we cannot even expect every operator
in this homotopy to be Fredholm.
\end{enumerate}

The first problem motivates the desire to define a notion of orientations
for a Fredholm operator $\mathbf{T}$ that does not require $\mathbf{T}$
to be surjective but reduces to the usual notion of orienting
$\ker \mathbf{T}$ whenever it is.  The solution to this problem is the
\emph{determinant line bundle}, which we will discuss in the next section.
With this object in hand, the above discussion for the case of closed
curves can be made rigorous, so that all smooth moduli spaces of
closed $J$-holomorphic curves inherit canonical orientations.
One of the advantages of using the determinant line bundle is that the
question of orientations becomes entirely disjoined from the question of
transversality: if one can orient the determinant line bundle then
moduli spaces of regular curves inherit orientations, but orienting the
determinant bundle does not require knowing in advance whether the curves
are regular.

The second problem is obviously significant because in the punctured case,
moduli spaces of $J$-holomorphic curves sometimes have \emph{odd} real
dimension, making it clearly impossible to homotop $D\dbar_J(j_0,u_0)$
through Fredholm operators to one that is complex linear.  The solution
in this case will be to define orientations algorithmically via the
coherence condition, and we will describe a suitable algorithm for
this in \S\ref{sec:construct}.

\section{The determinant line bundle}
\label{sec:det}

Fix real Banach spaces $X$ and $Y$ and let $\Fred_\RR(X,Y)$ denote the
space of real-linear Fredholm operators, viewed as an open subset of
the Banach space $\Lin_\RR(X,Y)$ of all bounded linear operators.
We'll use the following notation throughout: if $V$ is an $n$-dimensional
real vector space, then the top-dimensional exterior power of $V$ is denoted by
$$
\Lambda\Max V := \Lambda^n V.
$$
This $1$-dimensional real vector space is spanned by any wedge product of the 
form $v_1 \wedge \ldots \wedge v_n$ where $(v_1,\ldots,v_n)$ is a basis of~$V$.
Denoting the dual space of $V$ by $V^*$, note that there is a canonical
isomorphism $(\Lambda\Max V)^* = \Lambda\Max V^*$.  If $\dim V = 0$, then
we adopt the convention $\Lambda\Max V = \RR$.

\begin{defn}
Given $\mathbf{T} \in \Fred_\RR(X,Y)$, the \defin{determinant line} of
$\mathbf{T}$ is the real $1$-dimensional vector space
$$
\det(\mathbf{T}) = \left(\Lambda\Max \ker \mathbf{T} \right) \otimes
\left( \Lambda\Max \coker \mathbf{T} \right)^*.
$$
\end{defn}

Our main goal in this section is to prove:

\begin{thm}
\label{thm:det}
There exists a topological vector bundle $\det(X,Y) \stackrel{\pi}{\longrightarrow}
\Fred_\RR(X,Y)$ of real rank~$1$ such that $\pi^{-1}(\mathbf{T}) = 
\det(\mathbf{T})$ for each $\mathbf{T} \in \Fred_\RR(X,Y)$.
\end{thm}

Observe that whenever $\mathbf{T} \in \Fred_\RR(X,Y)$ is surjective, 
$\det(\mathbf{T}) = \Lambda\Max \ker\mathbf{T}$, so an orientation
of $\det(\mathbf{T})$ is equivalent to an orientation of $\ker\mathbf{T}$.
More generally, an orientation of $\det(\mathbf{T})$ is equivalent to an
orientation for $\ker \mathbf{T} \oplus \coker \mathbf{T}$.
If $\mathbf{T}$ is an isomorphism, then $\det(\mathbf{T})$ is simply~$\RR$,
so an orientation of $\det(\mathbf{T})$ amounts to a choice of sign $\pm 1$.

To construct local trivializations of $\det(X,Y) \to \Fred_\RR(X,Y)$, we start 
with the case where $X$ and $Y$ are both finite dimensional.  Note that in this 
case, every linear map is Fredholm, including the zero map, and its determinant
is simply $\Lambda\Max X \otimes (\Lambda\Max Y)^*$.

\begin{lemma}
\label{lemma:finDimDet}
Suppose $X$ and $Y$ are real vector spaces of finite dimensions $n$ and $m$
respectively.  Then for every $\mathbf{T} \in \Lin_\RR(X,Y)$, there exists
a canonical isomorphism
$$
\left( \Lambda\Max \ker \mathbf{T} \right) \otimes
\left( \Lambda\Max \coker \mathbf{T} \right)^* =
(\Lambda\Max V) \otimes (\Lambda\Max W)^*.
$$
\end{lemma}
\begin{proof}
Suppose $\dim \ker \mathbf{T} = k$ and $\dim \coker \mathbf{T} = \ell$,
so $\ind(\mathbf{T}) = k - \ell = n - m$, thus $n - k = m - \ell$.
We define a linear map $\Phi : (\Lambda^n X) \otimes (\Lambda^m Y)^* \to
\left( \Lambda^k \ker \mathbf{T} \right) \otimes 
\left( \Lambda^\ell \coker \mathbf{T} \right)^*$
via the following procedure.
Fix $\mathbf{x} \in \Lambda^n X$ and $\mathbf{y}^* \in (\Lambda^m Y)^*$ 
and suppose both are nontrivial.  Then for any nontrivial element 
$\mathbf{k} \in \Lambda^k \ker \mathbf{T}$, there exists a unique
element $\mathbf{v} \in \Lambda^{n-k}\left(X / \ker\mathbf{T}\right)$
such that for any subspace $V \subset X$ complementary to $\ker\mathbf{T}$,
the element $\tilde{\mathbf{v}} \in \Lambda^{n-k} V \subset \Lambda^{n-k} X$
obtained from $\mathbf{v}$ by inverting the natural isomorphism
$V \to X / \ker\mathbf{T}$ induced by the projection $X \to X / \ker\mathbf{T}$
satisfies
$$
\mathbf{k} \wedge \tilde{\mathbf{v}} = \mathbf{x}.
$$
The map $\mathbf{T}$ descends to an isomorphism
$X / \ker\mathbf{T} \to \im\mathbf{T}$ and thus induces an isomorphism
$\Lambda^{n-k}\left( X / \ker\mathbf{T}\right) \to
\Lambda^{m-\ell}\left( \im \mathbf{T} \right) \subset \Lambda^{m-\ell} Y$,
which takes $\mathbf{v}$ to a nontrivial element $\mathbf{T}\mathbf{v}$.
There is then a unique element $\mathbf{c} \in \Lambda^\ell \coker \mathbf{T}
= \Lambda^\ell\left(Y / \im\mathbf{T}\right)$ such that for any
subspace $W \subset Y$ complementary to $\im\mathbf{T}$, the element
$\tilde{\mathbf{c}} \in \Lambda^\ell W \subset \Lambda^\ell Y$ obtained
from $\mathbf{c}$ by inverting the isomorphism
$W \to Y / \im\mathbf{T}$ induced by the projection $Y \to Y / \im\mathbf{T}$
satisfies
$$
\mathbf{y}^*\left(\tilde{\mathbf{c}} \wedge \mathbf{T}\mathbf{v}\right) = 1.
$$
Now define $\Phi$ as the unique linear map such that
$$
\Phi(\mathbf{x} \otimes \mathbf{y}^*) = \mathbf{k} \otimes \mathbf{c}^*,
$$
where $\mathbf{c}^* \in (\Lambda^\ell \coker \mathbf{T})^*$ is defined by
$\mathbf{c}^*(\mathbf{c}) = 1$.  It is straightforward to check that this
definition does not depend on any choices: indeed, if we replace
$\mathbf{k}$ by $\lambda \mathbf{k}$ for some $\lambda \in \RR \setminus \{0\}$
in the above procedure, then $\mathbf{v}$ is replaced by 
$\frac{1}{\lambda} \mathbf{v}$, hence $\mathbf{T}\mathbf{v}$ becomes
$\frac{1}{\lambda} \mathbf{T}\mathbf{v}$, $\mathbf{c}$ becomes
$\lambda \mathbf{c}$ and $\mathbf{c}^*$ therefore becomes 
$\frac{1}{\lambda} \mathbf{c}^*$, so that $\mathbf{k} \otimes \mathbf{c}^*$
is replaced by
$$
(\lambda \mathbf{k}) \otimes \left( \frac{1}{\lambda} \mathbf{c}^* \right)
= \mathbf{k} \otimes \mathbf{c}^*.
$$
\end{proof}

To construct local trivializations of $\det(X,Y)$ in the infinite-dimensional
case, recall the following construction from Lecture~\ref{lec:asymptotic}.
Given $\mathbf{T}_0 \in \Fred_\RR(X,Y)$, we can write $X = V \oplus K$
and $Y = W \oplus C$ where $K = \ker \mathbf{T}_0$, $C \cong \coker \mathbf{T}_0$,
$W = \im \mathbf{T}_0$ and $\mathbf{T}_0|_V : V \to W$ is an isomorphism.
We shall use these splittings to write any other operator
$\mathbf{T} \in \Fred_\RR(X,Y)$ as
$$
\mathbf{T} = \begin{pmatrix}
\mathbf{A} & \mathbf{B} \\
\mathbf{C} & \mathbf{D}
\end{pmatrix}
$$
and let $\uU \subset \Fred_\RR(X,Y)$ denote the open neighborhood of
$\mathbf{T}_0$ for which the block $\mathbf{A} : V \to W$ is invertible.
This gives rise to a pair of smooth maps
$$
\Phi : \uU \to \Lin_\RR(K,C) : \mathbf{T} \mapsto \mathbf{D} -
\mathbf{C} \mathbf{A}^{-1} \mathbf{B}
$$
and
$$
F : \uU \to \Lin_\RR(V \oplus K) = \Lin_\RR(X) : \mathbf{T} \mapsto
\begin{pmatrix}
\1 & -\mathbf{A}^{-1} \mathbf{B} \\
0 & \1
\end{pmatrix},
$$
such that $F(\mathbf{T})$ is always invertible and maps
$\{0\} \oplus \ker \Phi(\mathbf{T})$ isomorphically to $\ker \mathbf{T}$.
Similarly, there is a smooth map
$$
G : \uU \to \Lin_\RR(W \oplus C) = \Lin_\RR(Y) : \mathbf{T} \mapsto
\begin{pmatrix}
\1 & 0 \\
- \mathbf{C} \mathbf{A}^{-1} & \1
\end{pmatrix}
$$
such that $G(\mathbf{T})$ is always invertible and maps
$\im \mathbf{T}$ isomorphically to $W \oplus \im \Phi(\mathbf{T})$, so it
descends to an isomorphism of $\coker \mathbf{T}$ to
$\coker \Phi(\mathbf{T})$.  Given the canonical isomorphism
$\det(\Phi(\mathbf{T})) = \Lambda\Max K \otimes (\Lambda\Max C)^* =
\det(\mathbf{T}_0)$ from Lemma~\ref{lemma:finDimDet}, the resulting smooth
families of isomorphisms $\ker \mathbf{T} \to \ker \Phi(\mathbf{T})$ and
$\coker \mathbf{T} \to \coker \Phi(\mathbf{T})$ determine
a local trivialization
$$
\det(X,Y)|_{\uU} \to \uU \times \det(\mathbf{T}_0).
$$
I will leave it as an exercise for the reader to check that the resulting
transition maps are continuous.\footnote{This detail should not be
underestimated, e.g.~\cite{McDuffWehrheim:trivial}*{\S 7.4} observes
that the local trivializations constructed in \cite{McDuffSalamon:Jhol}*{\S A.2}
are, unfortunately, not continuously compatible.  See \cite{Zinger:determinant}
for further discussion of this point.  If you discover that my local
trivializations are also not continuously compatible, please let me know.}

\begin{exercise}
\label{EX:complexOr}
Show that if $X$ and $Y$ are complex Banach spaces, then the restriction of
$\det(X,Y)$ to the subspace of complex-linear Fredholm operators $\Fred_\CC(X,Y)
\subset \Fred_\RR(X,Y)$
admits a canonical orientation compatible with the complex structures
of $\ker \mathbf{T}$ and $\coker \mathbf{T}$ for each
$\mathbf{T} \in \Fred_\CC(X,Y)$.  Show also that whenever $\mathbf{T} \in
\Fred_\CC(X,Y)$ is an isomorphism, the canonical orientation of
$\det(\mathbf{T})$ agrees with the standard orientation of~$\RR$.
\end{exercise}

The orientation of $\det(\mathbf{T})$ for $\mathbf{T} \in \Fred_\CC(X,Y)$
described in Exercise~\ref{EX:complexOr} is called the
\defin{complex orientation}.

\section{Determinant bundles of moduli spaces}
\label{sec:closed}

Combining ideas from the previous two sections, let
$$
\det(J) \to \mM^\$(J)
$$
denote the topological line bundle that associates to any
$u \in \mM_{g,m}^\$(J,A,\boldsymbol{\gamma}^+,\boldsymbol{\gamma}^-)$ the
determinant line of the Fredholm operator
$$
\mathbf{D}_u : W^{k,p,\delta}(u^*T\widehat{W}) \oplus V_\Gamma \to
W^{k-1,p,\delta}(\overline{\Hom}_\CC(T\dot{\Sigma},u^*T\widehat{W})).
$$
One can construct local trivializations for this bundle using
Theorem~\ref{thm:det} and any choice of local trivializations for the
Banach space bundles $T \bB^{k,p,\delta}$ and~$\eE^{k-1,p,\delta}$.

\begin{prop}
Any orientation of $\det(J) \to \mM^\$(J)$ canonically determines an orientation
of~$\mM^\reg(J)$.
\end{prop}
\begin{proof}
As explained in \S\ref{sec:moduli}, an orientation of $\mM^\reg(J)$ near
a particular curve $u_0 : (\dot{\Sigma},j_0) \to (\widehat{W},J)$ is
equivalent to a continuously varying choice of orientations for the
kernels
$$
\ker D\dbar_J(j,u) \subset T_j \tT \oplus T_u \bB^{k,p,\delta}
$$
for all $(j,u) \in \dbar_J^{-1}(0)$, where $\tT$ is a Teichm\"uller slice
through~$j_0$.  The operator $D\dbar_J(j,u)$ is of the form
$$
\mathbf{L} (y,\eta) := J \circ Tu \circ y + \mathbf{D}_u \eta
$$
and thus is homotopic through Fredholm operators to
$$
\mathbf{L}^0(y,\eta) := \mathbf{D}_u \eta,
$$
namely via the homotopy $\mathbf{L}^s(y,\eta) := s J \circ Tu \circ y +
\mathbf{D}_u \eta$ for $s \in [0,1]$.  The kernel and cokernel of
$\mathbf{L}^0$ are $T_j\tT \oplus \ker \mathbf{D}_u$ and $\coker \mathbf{D}_u$
respectively, and since $T_j\tT$ carries a complex structure, the orientation
of $\det(\mathbf{D}_u)$ naturally determines an orientation
of $\det(\mathbf{L}^0)$.
Using the homotopy $\mathbf{L}^s$, this determines orientations
of $\det(D\dbar_J(j,u))$ and thus orientations of $\ker D\dbar_J(j,u)$
for all $(j,u)$ near~$(j_0,u_0)$, and this orientation does not depend on
the choice of Teichm\"uller slice since the operators $\mathbf{D}_u$
also do not.
\end{proof}

From now on, when we speak of an \defin{orientation of $\mM^\$(J)$}, we will
actually mean an orientation of the bundle $\det(J) \to \mM^\$(J)$.
The above proposition implies that this is equivalent to what we want in
applications, but one advantage of talking about $\det(J)$ is that
there is no need to limit the discussion to curves that are regular,
i.e.~the notion of an orientation of $\mM^\$(J)$ now makes sense even
though $\mM^\$(J)$ is not globally a smooth object.

\begin{prop}
\label{prop:complexOr}
Suppose all Reeb orbits in $\boldsymbol{\gamma}^\pm$ have the property that
their asymptotic operators are complex linear.  Then 
$\mM^\$_{g,m}(J,A,\boldsymbol{\gamma}^+,\boldsymbol{\gamma}^-)$ admits a
natural orientation, known as the \defin{complex orientation}.
\end{prop}
\begin{proof}
Having complex-linear asymptotic operators implies that the obvious
homotopy from each Cauchy-Riemann operator $\mathbf{D}_u$ to its
complex-linear part does not change the asymptotic operators and is therefore
a homotopy through Fredholm operators.  We therefore have a continuously
varying homotopy of each of the relevant fibers of $\det(J)$ to the
determinant bundle over a family of complex-linear operators, which inherit
the complex orientation described in Exercise~\ref{EX:complexOr}.
\end{proof}

Proposition~\ref{prop:complexOr} applies in particular to all moduli spaces
of closed $J$-holomorphic curves, and thus solves the orientation problem
in that case.

\section{An algorithm for coherent orientations}
\label{sec:construct}

We now briefly describe the construction of coherent orientations due to
Bourgeois and Mohnke \cite{BourgeoisMohnke}.
A slightly different construction is described in \cite{SFT}, though it
appears to have minor errors in some details.

Recall from Lecture~\ref{lec:Fredholm} the notion of an \emph{asymptotically Hermitian}
vector bundle $(E,J)$ over a punctured Riemann surface $(\dot{\Sigma},j)$.
Here $(\dot{\Sigma},j)$ is endowed with the extra structure of fixed
cylindrical ends $(\dot{\uU}_z,j) \cong (Z_\pm,i)$ for each puncture
$z \in \Gamma^\pm$, which determines a choice of asymptotic markers.
Likewise, the bundle $E$ comes with an asymptotic bundle
$(E_z,J_z,\omega_z) \to S^1$ associated to each puncture, carrying compatible 
complex and symplectic structures.  We shall now endow
$E$ with a bit more structure that is always naturally present in the
case $E = u^*T\widehat{W}$: namely, assume each of the asymptotic bundles 
comes with a splitting
\begin{equation}
\label{eqn:splitting2}
(E_z,J_z,\omega_z) = (\CC \oplus \widehat{E}_z, i \oplus \hat{J}_z,
\omega_0 \oplus \hat{\omega}_z), 
\end{equation}
where $\omega_0$ is the standard symplectic structure on the trivial complex
line bundle $(\CC,i)$ over $S^1$, and 
$(\widehat{E}_z,\hat{J}_z,\hat{\omega}_z) \to S^1$
is another Hermitian bundle.  Fix a choice $\{\mathbf{A}_z\}_{z \in \Gamma}$
of nondegenerate asymptotic operators on each of the 
bundles $(\widehat{E}_z,\hat{J}_z,\hat{\omega}_z)$, and define the 
topological space
$$
\CR(E,\{\mathbf{A}_z\}_{z \in \Gamma})
$$
to consist of all Cauchy-Riemann type operators on $E$ that are asymptotic
at the punctures $z \in \Gamma$ to the asymptotic operators
$$
(-i\p_t) \oplus \mathbf{A}_z : \Gamma(\CC \oplus \widehat{E}_z) \to
\Gamma(\CC \oplus \widehat{E}_z).
$$
This is an affine space, so it is contractible, and if $\delta > 0$ is
sufficiently small and $V_\Gamma \subset \Gamma(E)$ denotes a 
complex $(\#\Gamma)$-dimensional space of smooth sections that take
constant values in $\CC \oplus \{0\} \subset E_z$ near each puncture~$z$,
then every $\mathbf{D} \in \CR(E,\{\mathbf{A}_z\}_{z \in \Gamma})$
determines a Fredholm operator
$$
\mathbf{D} : W^{k,p,\delta}(E) \oplus V_\Gamma \to W^{k-1,p,\delta}(\overline{\Hom}_\CC(T\Sigma,E)).
$$
It follows that a choice of orientation
of the determinant line for any one of these operators determines an 
orientation for all of them.  The point of this construction is that
every $u \in \mM^\$(J)$ determines an operator $\mathbf{D}_u$ belonging to
a space of this form.

We now construct a gluing operation for Cauchy-Riemann operators that
linearizes the gluing maps described in \S\ref{sec:coherence}.
Suppose $(E^i,J^i) \to (\dot{\Sigma}_i = \Sigma_i \setminus \Gamma_i,j_i)$
for $i=0,1$ is a pair of asymptotically Hermitian bundles of the same rank,
endowed with asymptotic splittings as in \eqref{eqn:splitting2} and
asymptotic operators $\{\mathbf{A}_z\}_{z \in \Gamma_i}$, and
that there exists a pair of punctures $z_0 \in \Gamma_0^+$ and 
$z_1 \in \Gamma_1^-$ such that some unitary bundle isomorphism
$$
\widehat{E}^1_{z_1} \stackrel{\cong}{\longrightarrow} \widehat{E}^0_{z_0}
$$
identifies $\mathbf{A}_{z_1}$ with $\mathbf{A}_{z_0}$.  Note that such an
isomorphism is uniquely determined up to homotopy whenever it exists.
For $R > 0$, we can define a family of glued Riemann surfaces
$$
(\dot{\Sigma}_R = \Sigma_R \setminus \Gamma_R,j_R)
$$
by cutting off the ends $(R,\infty) \times S^1 \subset \dot{\uU}_{z_0}$
and $(-\infty,-R) \times S^1 \subset \dot{\uU}_{z_1}$ and gluing 
$\{R\} \times S^1 \subset \dot{\Sigma}_0$ to $\{-R\} \times S^1 \subset
\dot{\Sigma}_1$.  The glued Riemann surface contains an annulus
biholomorphic to $([-R,R] \times S^1,i)$ in place of the infinite cylindrical
ends at the punctures $z_0$ and~$z_1$.
The unitary isomorphism $\widehat{E}^1_{z_1} \to \widehat{E}^0_{z_0}$ then 
determines an isomorphism $E^1_{z_1} \to E^0_{z_0}$ via the splitting
\eqref{eqn:splitting2} and hence an asymptotically Hermitian bundle
$$
(E^R,J^R) \to (\dot{\Sigma}_R,J_R).
$$
Using cutoff functions in the neck $[-R,R] \times S^1$, any Cauchy-Riemann
operators $\mathbf{D}_i \in \CR(E^i,\{\mathbf{A}_z\}_{z \in \Gamma_i})$ for $i=0,1$
now determine a family of operators
$$
\mathbf{D}_R \in \CR(E^R,\{\mathbf{A}_z\}_{z \in \Gamma_R})
$$
uniquely up to homotopy.  Analogously to the gluing maps in
\S\ref{sec:coherence}, one can arrange this construction so that the
operators $\mathbf{D}_R$ converge in some sense to the pair
$(\mathbf{D}_0,\mathbf{D}_1)$ as $R \to \infty$, which has the following
consequence:

\begin{lemma}[\cite{BourgeoisMohnke}*{Corollary~7}]
\label{lemma:iso}
For $R > 0$ sufficiently large, there is a natural isomorphism
$$
\det(\mathbf{D}_0) \otimes \det(\mathbf{D}_1) \to \det(\mathbf{D}_R)
$$
that is defined up to homotopy.
\qed
\end{lemma}

Up to some additional direct sums and quotients by finite-dimensional
complex vector spaces, this isomorphism should be understood as the
linearization of a gluing map between moduli spaces, generalized to a
setting in which the holomorphic curves involved need not be regular.
To orient $\mM^\$(J)$ coherently, it now suffices to choose orientations for the
operators in $\CR(E,\{\mathbf{A}_z\}_{z \in \Gamma})$ that vary
continuously under deformations of $j$ and $E$ and are preserved by the
isomorphisms of Lemma~\ref{lemma:iso}.  This motivates the following
generalization of Definition~\ref{defn:coherence}.

\begin{defn}
\label{defn:coherence2}
A system of \defin{coherent orientations} is an assignment to each
asymptotically Hermitian bundle $(E,J) \to (\dot{\Sigma},j)$ with
asymptotic splittings as in \eqref{eqn:splitting2} and asymptotic
operators $\{\mathbf{A}_z\}_{z \in \Gamma}$ of an orientation for the
determinant line of each $\mathbf{D} \in \CR(E,\{\mathbf{A}_z\})$, such that
these orientations vary continuously with $\mathbf{D}$ as well as the
data $j$ and~$J$, and such that the isomorphisms in Lemma~\ref{lemma:iso} 
are always orientation preserving.
\end{defn}

The prescription of \cite{BourgeoisMohnke} to construct such systems is
now as follows.

\begin{enumerate}
\item For any trivial bundle $E$ over $\dot{\Sigma} = \CC$ with $\infty$
as a negative puncture and any asymptotic operator $\mathbf{A}_\infty$,
choose an arbitrary continuous family of orientations for the operators
in $\CR(E,\{\mathbf{A}_\infty\})$, subject only to the requirement that these
should match the complex orientation whenever $\mathbf{A}_\infty$ is
complex linear.
\item For any trivial bundle $E_-$ over $\dot{\Sigma} = \CC$ with $\infty$
as a positive puncture, any asymptotic operator $\mathbf{A}_\infty$ and
any $\mathbf{D}_- \in \CR(E_-,\{\mathbf{A}_\infty\})$, let $E_+$ denote
the trivial bundle over $\CC$ with a negative puncture as in step~(1),
choose any $\mathbf{D}_+ \in \CR(E_+,\{\mathbf{A}_\infty\})$ and
construct the resulting family of glued operators
$$
\mathbf{D}_R \in \CR(E^R),
$$
where the $E^R$ are trivial bundles over~$S^2$.  Since $S^2$ has no
punctures, $\mathbf{D}_R$ has a natural complex orientation, so define
the orientation of $\mathbf{D}_-$ to be the one that is compatible via
Lemma~\ref{lemma:iso} with this and the orientation chosen for
$\mathbf{D}_+$ in step~(1).
\item For an arbitrary $(E,J) \to (\dot{\Sigma},j)$, glue positive and
negative planes to $\dot{\Sigma}$ to produce a bundle over a closed surface~$\widehat{\Sigma}$,
and define the orientation of any $\mathbf{D} \in \CR(E,\{\mathbf{A}_z\}_{z \in \Gamma})$
to be compatible via Lemma~\ref{lemma:iso} with the choices in
steps~(1) and~(2) and the complex orientation for operators 
over~$\widehat{\Sigma}$.
\end{enumerate}

It should be easy to convince yourself that if we now vary the bundle
$(E,J) \to (\dot{\Sigma},j)$ or the operators on this bundle 
(but \emph{not} the asymptotic operators!) continuously, the capping
procedure described in step~(3) above produces a continuous family of
Cauchy-Riemann type operators on bundles over closed Riemann surfaces.
Since these all carry the complex orientation, the resulting orientations
of the original operators vary continuously.  It is similarly clear from
the construction that any Cauchy-Riemann operator whose asymptotic operators
are all complex linear will end up with the complex orientation.
Bourgeois and Mohnke use this fact to prove that any system of orientations 
constructed in this way is compatible with \emph{all} possible linear 
gluing maps arising from Lemma~\ref{lemma:iso}.  The idea is to reduce it
to the complex-linear case by gluing cylinders to the ends of any
asymptotically Hermitian bundle so that the asymptotic operators can be
changed at will; see \cite{BourgeoisMohnke}*{Proposition~8}.

\section{Permutations and bad orbits revisited}
\label{sec:funPart2}

The heuristic proofs in \S\ref{sec:funPart} can now be made precise in the
following way.

Suppose $\mathbf{D} \in \CR(E,\{\mathbf{A}_z\}_{z \in \Gamma})$, and
$\mathbf{D}'$ is the same operator after interchanging two of the
punctures in~$\Gamma$.  Imagine gluing $(E,J) \to (\dot{\Sigma},j)$
to trivial bundles $E^1$ and $E^2$ over planes in order to cap off the two 
punctures that are being interchanged, and choose Cauchy-Riemann operators
$\mathbf{D}_1$ and $\mathbf{D}_2$ on these planes to form a glued
operator on the capped surface.
This capping procedure is done one plane
at a time, and the order of the two punctures determines which plane is
glued first.  Compatibility with the isomorphisms of
Lemma~\ref{lemma:iso} then dictates that the orientations of
$\det(\mathbf{D})$ and $\det(\mathbf{D}')$ match if and only if
the orientations of $\det(\mathbf{D}_1) \otimes \det(\mathbf{D}_2)$ and
$\det(\mathbf{D}_2) \otimes \det(\mathbf{D}_1)$ match.  Since orientations
of $\det(\mathbf{D}_i)$ for $i=1,2$ are equivalent to orientations
of $\ker \mathbf{D}_i \oplus \coker \mathbf{D}_i$, reversing the order
of the tensor product changes orientations if and only if
both of these direct sums are odd dimensional, which means
$\ind(\mathbf{D}_1)$ and $\ind(\mathbf{D}_2)$ are both odd.  If the bundles
have complex rank $n$ and the asymptotic operators are $\mathbf{A}_i$
for $k=1,2$, we have
$$
\ind(\mathbf{D}_i) = n \chi(\CC) \pm \muCZ((-i\p_t \oplus \mathbf{A}_i) \pm \delta)
= n - 1 \pm \muCZ(\mathbf{A}_i),
$$
which matches $n-3 + \muCZ(\mathbf{A}_i)$ modulo~$2$.
This proves Proposition~\ref{prop:permute}.

Similarly for Proposition~\ref{prop:bad}, we consider the action of the
generator $\psi \in \ZZ^m$ on $\det(\mathbf{D})$ where
$\psi$ rotates the cylindrical end by $1/m$ at some puncture where
the trivialized asymptotic operator $\mathbf{A}$ is of the form $-i\p_t - S(mt)$ 
for a loop of symmetric matrices~$S(t)$.  Capping off this puncture with a
plane carrying a Cauchy-Riemann operator $\mathbf{D}_\infty$, coherence
dictates that the same transformation must act the same way on the orientation 
of~$\det(\mathbf{D}_\infty)$.  Since $\psi^m = 1$, $\psi$ cannot reverse
this orientation if $m$ is odd.  To understand the case of $m$ even, note
first that we are free to choose $\mathbf{D}_\infty$
so that it is an $m$-fold cover, meaning it is related to the branched
cover $\varphi : \CC \to \CC : z \mapsto z^m$ by
$$
\mathbf{D}_\infty(\eta \circ \varphi) = \varphi^*\widehat{\mathbf{D}}_\infty \eta
$$
for some other Cauchy-Riemann operator $\widehat{\mathbf{D}}_\infty$,
which is asymptotic to $\hat{\mathbf{A}} := - i\p_t - S(t)$.
Now the group $\ZZ_m$ generated by $\psi$ acts on $\ker \mathbf{D}_\infty$
and $\coker \mathbf{D}_\infty$, so representation theory tells us
\begin{equation*}
\begin{split}
\ker \mathbf{D}_\infty &= V_1 \oplus V_{-1} \oplus V_{\operatorname{rot}} \\
\coker \mathbf{D}_\infty &= W_1 \oplus W_{-1} \oplus W_{\operatorname{rot}},
\end{split}
\end{equation*}
where $\psi$ acts on $V_{\pm 1}$ and $W_{\pm 1}$ as $\pm \1$ and 
acts as orientation-preserving rotations on $V_{\operatorname{rot}}$ 
and~$W_{\operatorname{rot}}$.  It follows that $\psi$ reverses the orientation
of $\ker\mathbf{D}_\infty \oplus \coker \mathbf{D}_\infty$ if and only if
$\dim V_{-1} - \dim W_{-1}$ is odd.  Now observe that there are natural
isomorphisms
$$
V_1 = \ker \widehat{\mathbf{D}}_\infty, \qquad
W_1 = \coker \widehat{\mathbf{D}}_\infty,
$$
hence
$$
\dim V_{-1} - \dim W_{-1} = \ind(\mathbf{D}_\infty) - \ind(\widehat{\mathbf{D}}_\infty) \qquad
\text{(mod~2)}.
$$
This difference in Fredholm indices is precisely
$\muCZ(\mathbf{A}) - \muCZ(\hat{\mathbf{A}})$ up to a sign, and this
completes the proof of Proposition~\ref{prop:bad}.

\psfrag{ua}{$u$}
\psfrag{ub}{$v$}
\psfrag{g1}{$\gamma_1$}
\psfrag{g2}{$\gamma_2$}
\psfrag{ga}{$\boldsymbol{\gamma}^-_u$}
\psfrag{gb}{$\boldsymbol{\gamma}^+_v$}
\psfrag{g3}{$\gamma_3$}
\psfrag{g4}{$\gamma_4$}
\psfrag{g5}{$\gamma_5$}
\psfrag{g6}{$\gamma_6$}
\psfrag{g7}{$\gamma_7$}
\psfrag{g8}{$\gamma_8$}
\psfrag{u}{$u$}

\chapter{The generating function of SFT}
\label{lec:H}

\minitoc

\vspace{12pt}

It is time to begin deriving algebraic consequences from the analytical
results of the previous lectures.  We saw the simplest possible example
of this in Lecture~\ref{lec:tight3tori}, where the behavior of holomorphic cylinders in 
symplectizations of contact manifolds without contractible Reeb orbits
led to a rudimentary version of cylindrical contact homology
$HC_*(M,\xi)$ with $\ZZ_2$ coefficients.  Unfortunately, the condition on
contractible orbits means that this version of
$HC_*(M,\xi)$ cannot always be defined, and even when 
it can, it only counts cylinders---we would only expect it to capture a small 
fragment of the information contained in more general moduli spaces of 
holomorphic curves.  Extracting information from these general moduli spaces
will require enlarging our algebraic notion of what a Floer-type theory can look like.

\section{Some important caveats on transversality}
\label{sec:caveat}

For this and the next lecture, we fix the following fantastically optimistic
assumption:
\begin{assumption}[science fiction]
\label{ass:optimism}
One can choose suitably compatible almost complex structures so that
all pseudoholomorphic curves are Fredholm regular.
\end{assumption}
This assumption held in Lecture~\ref{lec:tight3tori} for the curves we were interested in,
because they were all guaranteed for topological reasons to be somewhere
injective.  It can also be shown to hold under some very restrictive
conditions on Conley-Zehnder indices in dimension three, see \cites{Nelson:Abendblatt,Nelson:thesis}.
Both of those are very lucky situations, and
as we've discussed before, the assumption cannot generally be achieved merely 
by perturbing $J$ generically---it \emph{must} sometimes fail for curves
that are multiply covered, and such curves always exist (see \S\ref{sec:orbifold}
for more on this).  The only way in reality to ensure something
like Assumption~\ref{ass:optimism} is to perturb the nonlinear Cauchy-Riemann
equation more abstractly, e.g.~by replacing $\dbar_J u = 0$ with an
inhomogeneous equation of the form
$$
\dbar_J u = \nu
$$
for a generic perturbation~$\nu$.  This is the standard technique in
certain versions of Gromov-Witten theory, see 
e.g.~\cites{RuanTian,RuanTian:higherGenus}.  Alternatively, one can allow
$J$ to depend generically on points in the domain rather than just points
in the target, as in \cite{McDuffSalamon:Jhol}*{\S 7.3}.  
Both approaches eliminate the initial problem with
multiple covers, but they both also run into serious and subtle difficulties
concerning the relationship between $\mM(J)$ and the strata of its
compactification $\overline{\mM}(J)$.  As observed in 
\cite{Salamon:Floer}*{\S 5}, the possibility of symmetry in strata of
$\overline{\mM}(J)$ makes it necessary for any sufficiently general 
abstract perturbation scheme to involve \emph{multivalued} perturbations,
and it is important for these perturbations to be ``coherent'' in a sense
analogous to our discussion of orientations in the previous lecture.
These notions have not yet all been developed in a sufficiently consistent
and general way to give a rigorous definition of SFT, though there has been 
much progress: this is the main objective of the long-running
\emph{polyfold} project by Hofer-Wysocki-Zehnder \cite{Hofer:CDM}.
Recently, a quite different and much more topological approach has been
proposed by John Pardon \cite{Pardon:contact}.

For most of this lecture we will ignore these subtleties and
simply adopt Assumption~\ref{ass:optimism} as a convenient fiction, thus
pretending that all components of $\mM(J)$ are smooth orbifolds of the
correct dimension and
all gluing maps are smooth.  All ``theorems'' stated under this assumption
should be read with the caveat that they are only true in a fictional
world in which the assumption holds.  Even if it is a fiction, one can get
quite far with this point of view: it is still possible not only to deduce the
essential structure of what we assume will someday be a rigorously defined
polyfold-based SFT, but also to infer the existence of certain contact
invariants that have interesting rigorous applications requiring only
well-established techniques, e.g.~the cobordism obstructions discovered
in~\cite{LatschevWendl}.

\section{Auxiliary data, grading and supercommutativity}
\label{sec:setup3}

The goal is to define an invariant of closed $(2n-1)$-dimensional
contact manifolds $(M,\xi)$
with closed nondegenerate Reeb orbits as generators and a Floer-type
differential counting $J$-holomorphic curves in the symplectization
$(\RR \times M,d(e^r \alpha))$.  The auxiliary data we choose must obviously
therefore include a nondegenerate contact form $\alpha$ and a generic
$J \in \jJ(\alpha)$, for which we shall assume Assumption~\ref{ass:optimism}
holds.  For convenience, we will also assume throughout most of this lecture:

\begin{assumption}
\label{ass:torsion}
$H_1(M)$ is torsion free.
\end{assumption}

This is needed mainly in order to be able to define an integer grading,
though without this assumption, it is still always possible to define
a $\ZZ_2$-grading---see \S\ref{sec:torsion} for more on what to do
when Assumption~\ref{ass:torsion} does not hold.
We now supplement the auxiliary data $(\alpha,J)$ with the following
additional choices:
\begin{enumerate}
\item Coherent orientations as in Lecture~\ref{lec:orientations} for the moduli spaces
$\mM^\$(J)$ with asymptotic markers.
\item A collection of \defin{reference curves}
$$
S^1 \cong C_1,\ldots,C_r \subset M
$$
whose homology classes form a basis of~$H_1(M)$.
\item A unitary trivialization of $\xi$ along each of the reference curves
$C_1,\ldots,C_r$, denoted collectively by~$\tau$.
\item A \defin{spanning surface} $C_\gamma$ for each periodic Reeb orbit $\gamma$:
this is a smooth map of a compact and oriented surface with boundary into $M$
such that
$$
\p C_\gamma = \sum_i m_i [C_i] - [\gamma]
$$
in the sense of singular $2$-chains, where $m_i \in \ZZ$ are the unique
coefficients with
$[\gamma] = \sum_i m_i [C_i] \in H_1(M)$.
\end{enumerate}
These choices determine the following.  To any collections of Reeb orbits
$\boldsymbol{\gamma}^\pm = (\gamma_1^\pm,\ldots,\gamma_{k_\pm}^\pm)$ and
any relative homology class $A \in 
H_2(M,\bar{\boldsymbol{\gamma}}^+\cup \bar{\boldsymbol{\gamma}}^-)$ with
$\p A = \sum_i [\gamma_i^+] - \sum_j [\gamma_j^-]$, we can now associate
a cycle in absolute homology,
$$
A + \sum_i C_{\gamma_i^+} - \sum_j C_{\gamma_j^-} \in H_2(M).
$$
Indeed, the boundary of this real $2$-chain is a sum of linear combinations
of the reference curves $C_i$, which add up to zero because
$\sum_i [\gamma_i^+]$ and $\sum_j [\gamma_j^-]$ are homologous.  We shall
abuse notation and use this correspondence to associate the absolute homology class
$$
[u] \in H_2(M)
$$
to any asymptotically cylindrical holomorphic curve $u$ in $\RR \times M$.
Adapting the previous notation,
$$
\mM_{g,m}(J,A,\boldsymbol{\gamma}^+,\boldsymbol{\gamma}^-)
$$ 
for $A \in H_2(M)$ will now denote a moduli space of curves whose relative 
homology classes glue to the chosen capping surfaces to form~$A$.

Secondly, the chosen trivializations $\tau$
along the reference curves can be pulled back and extended over every
capping surface $C_\gamma$, giving trivializations of $\xi$ along
every orbit $\gamma$ uniquely up to homotopy.  We shall define
$$
\muCZ(\gamma) \in \ZZ
$$
from now on to mean the Conley-Zehnder index of $\gamma$ relative to this
trivialization.

\begin{exercise}
\label{EX:relc1}
Show that if $H_1(M)$ has no torsion and
$u : \dot{\Sigma} \to \RR \times M$ is asymptotically cylindrical,
then its relative first Chern number with respect to the trivializations
$\tau$ described above satisfies
$$
c_1^\tau(u^*T(\RR \times M)) = c_1([u]),
$$
where $c_1([u])$ denotes the evaluation of $c_1(\xi) \in H^2(M)$ on
$[u] \in H_2(M)$.
\end{exercise}

By Exercise~\ref{EX:relc1}, the index of a curve 
$u : (\dot{\Sigma} = \Sigma \setminus \Gamma,j) \to
(\RR \times M,J)$ with $[u] = A \in H_2(M)$ and asymptotic orbits
$\{\gamma_z\}_{z \in \Gamma^\pm}$ can now be written as
\begin{equation}
\label{eqn:indexAbs}
\ind(u) = -\chi(\dot{\Sigma}) + 2 c_1(A) + \sum_{z \in \Gamma^+} \muCZ(\gamma_z)
- \sum_{z \in \Gamma^-} \muCZ(\gamma_z).
\end{equation}

In order to keep track of homology classes of holomorphic curves algebraically,
we can define our theory to have coefficients in the group ring
$\QQ[H_2(M)]$, or more generally,
$$
R := \QQ[H_2(M) / G]
$$
for a given subgroup $G \subset H_2(M)$.  Elements of $R$ will be written as 
finite sums
$$
\sum_i c_i e^{A_i} \in R, \qquad c_i \in \QQ,\ A_i \in H_2(M) / G,
$$
where the multiplicative structure of the group ring is derived from the 
additive structure of $H_2(M) / G$ by $e^A e^B := e^{A+B}$.  The most
common examples of $G$ are $H_2(M)$ and the trivial subgroup, giving
$R = \QQ$ or $R = \QQ[H_2(M)]$ respectively.  We will see a geometrically
meaningful example in between these two extremes in the next lecture.

Finally, we define certain formal variables which have degrees in $\ZZ$ or
$\ZZ_{2N}$ for some $N \in \NN$, and will serve as generators in our graded
algebra.  To each closed Reeb orbit $\gamma$ we associate
two variables, $q_\gamma$, $p_\gamma$, whose integer-valued degrees are
$$
|q_\gamma| = n - 3 + \muCZ(\gamma), \qquad |p_\gamma| = n - 3 - \muCZ(\gamma).
$$
To remember these numbers, think of the index of a $J$-holomorphic plane $u$
positively or negatively asymptotic to $\gamma$, with $[u]=0$.

We also assign an integer grading to the group ring $\QQ[H_2(M)]$ such that 
rational numbers have degree $0$ and
$$
|e^A| = - 2 c_1(A), \quad \text{ for }\quad A \in H_2(M).
$$
If $c_1(A) = 0$ for every $A \in G$, in particular if $c_1(\xi) = 0$,
then this descends to an integer
grading on the ring $R = \QQ[H_2(M) / G]$.  Otherwise, $R$ inherits a
$\ZZ_{2N}$-grading, where
$$
N := \min \left\{ c_1(A) > 0 \ \big|\ A \in G \right\}.
$$
A $\ZZ_2$-grading is well defined in every case.

The algebra will include one additional formal variable $\hbar$,
which is defined to have degree
$$
|\hbar| = 2(n-3).
$$
The degrees of $\hbar$ and the $p_\gamma$ and $q_\gamma$ variables should
all be interpreted modulo $2N$ if $c_1(\xi)|_G \ne 0$.

The algebra of SFT uses monomials in the variables $p_\gamma$ and $q_\gamma$
respectively to encode sets of positive and negative asymptotic orbits of 
holomorphic curves, while the group ring $R = \QQ[H_2(M) / G]$ is used to
keep track of the homology classes of such curves, and powers of
$\hbar$ are used to keep track of their genus.
More precisely, given $g \ge 0$, $A \in H_2(M)$ and
ordered lists of Reeb orbits
$\boldsymbol{\gamma}^\pm = (\gamma_1^\pm,\ldots,\gamma_{k_\pm}^\pm)$, we
encode the moduli space 
$\mM_{g,0}(J,A,\boldsymbol{\gamma}^+,\boldsymbol{\gamma}^-)$ formally via
the product
\begin{equation}
\label{eqn:monomial}
e^A \hbar^{g-1} q^{\boldsymbol{\gamma}^-} p^{\boldsymbol{\gamma}^+} :=
e^A \hbar^{g-1} q_{\gamma_1^-}\ldots q_{\gamma_{k_-}^-}
p_{\gamma_1^+}\ldots p_{\gamma_{k_+}^+},
\end{equation}
where we are abusing notation by identifying $A$ with its equivalence class
in $H_2(M) / G$ if $G$ is nontrivial.
Notice that according to the above definitions, this expression has degree
\begin{equation}
\label{eqn:degree}
\begin{split}
| e^A \hbar^{g-1} & q^{\boldsymbol{\gamma}^-} p^{\boldsymbol{\gamma}^+} | =
|e^A| + (g-1)|\hbar| + \sum_{i=1}^{k_-} \left[ (n-3) + \muCZ(\gamma_i^-) \right] \\
&\qquad + \sum_{i=1}^{k_+} \left[ (n-3) - \muCZ(\gamma_i^+) \right] \\
&= - 2 c_1(A) + (2g-2 + k_+ + k_-)(n-3) - \sum_{i=1}^{k_+} \muCZ(\gamma_i^+)
+ \sum_{i=1}^{k_-} \muCZ(\gamma_i^-) \\
&= - \virdim \mM_{g,0}(J,A,\boldsymbol{\gamma}^+,\boldsymbol{\gamma}^-),
\end{split}
\end{equation}
interpreted modulo $2N$ if $c_1(\xi)|_G \ne 0$.
The orientation results in Lecture~\ref{lec:orientations} suggest introducing a 
\emph{supercommutativity} relation for the variables $q_\gamma$ and $p_\gamma$:
defining the graded commutator bracket by
\begin{equation}
\label{eqn:commutator}
[F,G] := FG - (-1)^{|F| |G|} GF,
\end{equation}
we define a relation on the set of all monomials of the form
$q^{\boldsymbol{\gamma}^-} p^{\boldsymbol{\gamma}^+}$ by setting
\begin{equation}
\label{eqn:nonquantum}
[q_{\gamma_1},q_{\gamma_2}] = [p_{\gamma_1},p_{\gamma_2}] = 0
\end{equation}
for all pairs of orbits $\gamma_1$ and~$\gamma_2$.  As a consequence,
permuting the orbits in the lists $\boldsymbol{\gamma}^\pm$ changes the
sign of the monomial \eqref{eqn:monomial} if and only if it changes the
orientation of the corresponding moduli space.  In particular, any
product that includes multiple copies of an odd generator $q_\gamma$ or
$p_\gamma$ is identified with~$0$.  This accounts for the fact that any 
rigid moduli space $\mM_{g,0}(J,A,\boldsymbol{\gamma}^+,\boldsymbol{\gamma}^-)$ 
with two copies of $\gamma$ among its positive or negative asymptotic orbits
contains zero curves when counted with the correct signs: every curve is
cancelled by a curve that looks identical except for a permutation of two of
its punctures.

\section{The definition of $\mathbf{H}$ and commutators}
\label{sec:H}

To write down the SFT generating function, let
$$
\mM^\sigma(J) := \mM(J) \big/ \sim
$$
denote the space of equivalence classes where two curves are considered
equivalent if they have parametrizations that differ only in the ordering
of the punctures.  This space is in some sense more geometrically natural
than $\mM(J)$ or $\mM^\$(J)$, but due to the orientation results in the
previous lecture, less convenient for technical reasons.  Given 
$u : (\dot{\Sigma},j) \to (\RR \times M,J)$ representing a nonconstant
element of $\mM^\sigma(J)$ with no marked points, it is natural to define
$$
\Aut^\sigma(u) \subset \Aut(\Sigma,j)
$$
as the (necessarily finite) group of biholomorphic transformations
$\varphi : (\Sigma,j) \to (\Sigma,j)$ satisfying $u = u \circ \varphi$;
in particular, elements of $\Aut^\sigma(u)$ are allowed to permute the
punctures, so $\Aut^\sigma(u)$ is generally a larger group than the
usual~$\Aut(u)$.  For $k \in \ZZ$, let
$$
\mM_k^\sigma(J) \subset \mM^\sigma(J)
$$
denote the subset consisting of index~$k$ curves that have no marked points
and whose asymptotic orbits are all \emph{good} (see Definition~\ref{defn:bad} in
Lecture~\ref{lec:orientations}).

We now define the \defin{SFT generating function} as a formal power series
\begin{equation}
\label{eqn:Hsigma}
\mathbf{H} = \sum_{u \in \mM_1^\sigma(J) / \RR}
\frac{\epsilon(u)}{|\Aut^\sigma(u)|}
\hbar^{g-1} e^A q^{\boldsymbol{\gamma}^-} p^{\boldsymbol{\gamma}^+},
\end{equation}
where the terms of each monomial are determined by $u \in \mM_1^\sigma(J)$
as follows:
\begin{itemize}
\item $g$ is the genus of~$u$;
\item $A$ is the equivalence class of $[u] \in H_2(M)$ in $H_2(M) / G$;
\item $\boldsymbol{\gamma}^\pm = (\gamma_1^\pm,\ldots,\gamma_{k_\pm}^\pm)$
are the asymptotic orbits of $u$ after arbitrarily fixing orderings of its
positive and negative punctures;
\item $\epsilon(u) \in \{1,-1\}$ is determined by the chosen coherent 
orientations on $\mM^\$(J)$.  Specifically, given the chosen ordering of
the punctures and an arbitrary choice of asymptotic markers at each puncture,
$u$ determines a $1$-dimensional connected component of $\mM^\$(J)$, and
we define $\epsilon(u) = +1$ if and only if the coherent orientation of
$\mM^\$(J)$ matches its tautological orientation determined by the
$\RR$-action.
\end{itemize}
Note that while both $\epsilon(u)$ and the corresponding monomial
$q^{\boldsymbol{\gamma}^-} p^{\boldsymbol{\gamma}^+}$ depend on a choice of
orderings of the punctures, their product does not depend on this choice.
Moreover, $\epsilon(u)$ does not depend on the choice of asymptotic markers
since curves with bad asymptotic orbits are excluded from $\mM_1^\sigma(J)$.
Since every monomial in $\mathbf{H}$ corresponds to a holomorphic curve
of index~$1$, \eqref{eqn:degree} implies
$$
|\mathbf{H}| = -1.
$$

There are various combinatorially more elaborate ways to rewrite~$\mathbf{H}$.
For any Reeb orbit $\gamma$, let
$$
\kappa_\gamma := \cov(\gamma) \in \NN
$$
denote its covering multiplicity, and for a finite list of orbits
$\boldsymbol{\gamma} = (\gamma_1,\ldots,\gamma_k)$, let
$$
\kappa_{\boldsymbol{\gamma}} := \prod_{i=1}^k \kappa_{\gamma_i}.
$$
Given $u \in \mM^\sigma(J)$ with $k_\pm \ge 0$ positive/negative punctures
asymptotic to the set of orbits 
$\boldsymbol{\gamma}^\pm = (\gamma_\pm^1,\ldots,\gamma_\pm^{k_\pm})$,
there are $k_+! k_-! \kappa_{\boldsymbol{\gamma}^+} \kappa_{\boldsymbol{\gamma}^-}$ 
ways to order the punctures and choose
asymptotic markers, but some of them are equivalent since (by an easy
variation on Proposition~\ref{prop:markerStructure}) the finite group
$\Aut^\sigma(u)$ acts freely
on this set of choices.  As a result, \eqref{eqn:Hsigma} is the same as
\begin{equation}
\label{eqn:Hdollar}
\mathbf{H} = \sum_{u \in \mM_1^\$(J) / \RR}
\frac{\epsilon(u)}{k_+! k_-! \kappa_{\boldsymbol{\gamma}^+} \kappa_{\boldsymbol{\gamma}^-}} 
\hbar^{g-1} e^A q^{\boldsymbol{\gamma}^-} p^{\boldsymbol{\gamma}^+},
\end{equation}
where $\mM_1^\$(J)$ denotes the space of all index~$1$ curves without marked 
points in $\mM^\$(J)$, and the rest of the mononomial is determined by
the condition that $u$ belongs to 
$\mM_{g,0}^\$(J,A,\boldsymbol{\gamma}^+,\boldsymbol{\gamma}^-)$,
with no need for any arbitrary choices.  Another way of writing this is
\begin{equation}
\label{eqn:HpowerSeries}
\mathbf{H} = \sum_{g,A,\boldsymbol{\gamma}^+,\boldsymbol{\gamma}^-}
\frac{\#\left( \mM_{g,0}^\$(J,A,\boldsymbol{\gamma}^+,\boldsymbol{\gamma}^-) \big/ \RR\right)}
{k_+! k_-! \kappa_{\boldsymbol{\gamma}^+} \kappa_{\boldsymbol{\gamma}^-}}
\hbar^{g-1} e^A q^{\boldsymbol{\gamma}^-} p^{\boldsymbol{\gamma}^+},
\end{equation}
where the sum ranges over all integers $g \ge 0$, homology classes
$A \in H_2(M)$ and ordered tuples of Reeb orbits $\boldsymbol{\gamma}^\pm =
(\gamma_1^\pm,\ldots,\gamma_{k_\pm}^\pm)$, and
$\#\left( \mM_{g,0}^\$(J,A,\boldsymbol{\gamma}^+,\boldsymbol{\gamma}^-) \big/ \RR\right)
\in \ZZ$ is the signed count of index~$1$ connected components in
$\mM_{g,0}^\$(J,A,\boldsymbol{\gamma}^+,\boldsymbol{\gamma}^-)$.  For fixed
$g$ and $\boldsymbol{\gamma}^\pm$, the union of these spaces for all
$A \in H_2(M)$ is finite due to SFT compactness, as the energy of curves in
$(\RR \times M,d(e^t\alpha))$ is computed by integrating exact symplectic forms
and thus (by Stokes) admits a uniform upper bound in terms of $\boldsymbol{\gamma}^+$.
For this reason, \eqref{eqn:HpowerSeries} defines a formal power series in the
$p$ variables and in~$\hbar$, with coefficients that are \emph{polynomials} 
in the $q$ variables and the group ring~$R$.

We played a slightly sneaky trick in writing down \eqref{eqn:Hdollar}
and \eqref{eqn:HpowerSeries}: these summations to not exclude bad orbits,
whereas \eqref{eqn:Hsigma} was a sum over curves $u$ that are not asymptotic to
any bad orbits---a necessary exclusion in that case because $\epsilon(u)$
would otherwise depend on choices of asymptotic markers.  The reason bad
orbits are allowed in \eqref{eqn:HpowerSeries} is that their total contribution
adds up to zero: indeed, bad orbits are always multiple covers with even
multiplicity, so whenever $u \in \mM^\$(J)$ has a puncture approaching
a bad orbit with multiplicity~$2m$, there are exactly $2m-1$ other elements
of $\mM^\$(J)$ that differ only by adjustment of the marker at that one
puncture, and by Proposition~\ref{prop:bad}, half of these cancel out
the other half in the signed count.  We've already seen that a similar
remark explains the harmless absence from \eqref{eqn:HpowerSeries} of 
terms with multiple factors of any odd generator $q_\gamma$ or~$p_\gamma$.

\begin{remark}
Readers famliar with Floer homology may see a resemblance between the
group ring $R = \QQ[H_2(M) / G]$ and the Novikov rings that often appear
in Floer homology, though $R$ is not a Novikov ring since it only allows
finite sums.  In Floer homology, the Novikov ring sometimes must be included
because counts of curves may fail to be finite, though they only do so if 
the energies of those curves blow up.  The situation above
is somewhat different: since the symplectization is an exact symplectic manifold,
Stokes' theorem implies that energy cannot blow up
if the positive asymptotic orbits are fixed, and one therefore obtains
well-defined curve counts no matter the choice of the coefficient
ring~$R$.  The use of the group ring is convenient however for two reasons:
first, without it one cannot always define an integer grading, 
and second, different choices of coefficients
can sometimes be used to detect different geometric phenomena via SFT.
We will see an example of the latter in Lecture~\ref{lec:SFT}.
\end{remark}

The compactness and gluing theory of SFT is encoded algebraically by
viewing $\mathbf{H}$ as an element on a noncommutative 
operator algebra determined by the commutator relations
\begin{equation}
\label{eqn:quantum}
\begin{split}
[p_\gamma,q_\gamma] &= \kappa_\gamma \hbar \\
[p_\gamma,q_{\gamma'}] &= 0 \quad \text{ if $\gamma \ne \gamma'$}.
\end{split}
\end{equation}
Here $[\ ,\ ]$ again denotes the graded commutator \eqref{eqn:commutator},
so ``commuting'' generators actually anticommute whenever they are both
odd.  The rest of the multiplicative structure of this algebra is
determined by requiring all elements of $R$ and powers of $\hbar$ (all of which
are even generators) to commute with everything, meaning all operators are
$R[[\hbar]]$-linear.

One concrete representation
of this operator algebra is as follows: let $\aA$ denote the graded
supercommutative unital algebra over $R$ generated by the set
$$
\left\{ q_\gamma\ \big|\ \text{$\gamma$ a good Reeb orbit} \right\}.
$$
The ring of formal power series $\aA[[\hbar]]$ is then an $R[[\hbar]]$-module.
Define each of the generators $q_\gamma$ to be $R[[\hbar]]$-linear operators
on $\aA[[\hbar]]$ via multiplication from the left, and define
$p_\gamma : \aA[[\hbar]] \to \aA[[\hbar]]$ by
\begin{equation}
\label{eqn:pgamma}
p_\gamma = \kappa_\gamma \hbar \frac{\p}{\p q_\gamma}.
\end{equation}
Here the $R[[\hbar]]$-linear partial derivative operator is defined via
$$
\frac{\p}{\p q_\gamma} q_{\gamma} = 1, \qquad
\frac{\p}{\p q_\gamma} q_{\gamma'} = 0 \quad \text{ for $\gamma \ne \gamma'$}
$$
and the graded Leibniz rule
$$
\frac{\p}{\p q_\gamma} (F G) = \frac{\p F}{\p q_\gamma} G +
(-1)^{|q_\gamma| |F|} F \frac{\p G}{\p q_\gamma}
$$
for all homogeneous elements $F,G \in \aA[[\hbar]]$.

\begin{exercise}
Check that the operator $p_\gamma : \aA[[\hbar]] \to \aA[[\hbar]]$ defined
above has the correct degree and satisfies the
commutation relations \eqref{eqn:nonquantum} and~\eqref{eqn:quantum}.
\end{exercise}

Notice that while $\mathbf{H}$ contains terms of order $-1$ in~$\hbar$,
every term also contains at least one $p_\gamma$ variable since all
index~$1$ holomorphic curves in $(\RR \times M,d(e^t\alpha))$ have at least
one positive puncture.  The substitution \eqref{eqn:pgamma} thus produces
a differential operator in which every term contains a nonnegative power
of~$\hbar$, giving a well-defined $R[[\hbar]]$-linear operator
$$
\mathbf{D}_\SFT : \aA[[\hbar]] \stackrel{\mathbf{H}}{\longrightarrow} \aA[[\hbar]].
$$
The following may be regarded as the fundamental theorem of SFT.

\begin{thm}
$\mathbf{H}^2 = 0$.
\end{thm}

We will discuss in \S\ref{sec:combi} how this relation follows from
the compactness and gluing theory of punctured holomorphic curves, and we 
will use it in Lecture~\ref{lec:SFT} to define various Floer-type contact invariants.
The first and most obvious of these is the homology
$$
H_*^\SFT(M,\xi) := H_*(\aA[[\hbar]],\mathbf{D}_\SFT),
$$
which will turn out to be an invariant of $(M,\xi)$ in the sense that 
any two choices of $\alpha$, $J$ and the other auxiliary data described
in \S\ref{sec:setup3} gives rise to a functorial isomorphism between the two 
graded homology
groups.  Notice that while $\aA[[\hbar]]$ is an algebra, its product
structure does not descend to $H_*^\SFT(M,\xi)$ since $\mathbf{D}_\SFT$
is not a derivation---indeed, it is a formal sum of differential operators of
all orders, not just order one.  In the next lecture we will discuss various
ways to produce homological invariants out of $\mathbf{H}$ with nicer
algebraic structures.  

\begin{figure}
\includegraphics{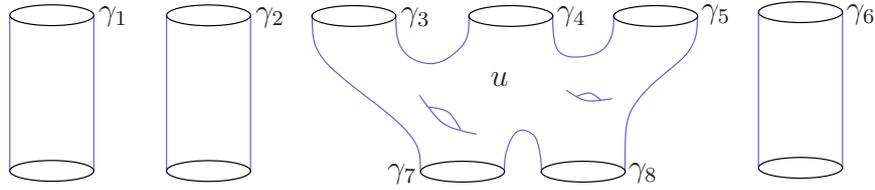}
\caption{\label{fig:DSFT} 
Counting disjoint unions of index~$1$ curves $u \in 
\mM_{2,0}(J,A,(\gamma_3,\gamma_4,\gamma_5),(\gamma_7,\gamma_8))$ with some
trivial cylinders contributes a multiple of 
$\hbar^4 e^A q_{\gamma_1} q_{\gamma_2} q_{\gamma_7} q_{\gamma_8} q_{\gamma_6}$
to $\mathbf{D}_\SFT (q_{\gamma_1} q_{\gamma_2} q_{\gamma_3} q_{\gamma_4} 
q_{\gamma_5} q_{\gamma_6} )$.}
\end{figure}

On the other hand, it is fairly easy to understand
the geometric meaning of the complex $(\aA[[\hbar]],\mathbf{D}_\SFT)$ in
Floer-theoretic terms.  Each individual curve $u \in \mM_1^\sigma(J)$ with
genus $g$, homology class $A \in H_2(M)$ and asymptotic orbits 
$\boldsymbol{\gamma}^\pm = (\gamma_1^\pm,\ldots,\gamma_{k_\pm}^\pm)$
contributes to $\mathbf{D}_\SFT$ the differential operator
$$
\frac{\epsilon(u)}{|\Aut^\sigma(u)|} \kappa_{\boldsymbol{\gamma}^+} 
\hbar^{g + k_+ - 1} e^A q_{\gamma_1^-} \ldots q_{\gamma_{k_-}^-}
\frac{\p}{\p q_{\gamma_1^+}} \ldots \frac{\p}{\p q_{\gamma_{k_+}^+}}.
$$
Applying this operator to a monomial $q_{\gamma_1}\ldots q_{\gamma_m} \in 
\aA[[\hbar]]$ that does not contain all of the generators
$q_{\gamma_1^+},\ldots,q_{\gamma_{k_+}^+}$ will produce zero, and its
effect on a product that does contain all of these generators will be to
eliminate them and multiply
$q_{\gamma_1^-} \ldots q_{\gamma_{k_-}^-}$ by whatever remains, plus some
combinatorial factors and signs that may arise from differentiating by the 
same $q_\gamma$ more than once.  Ignoring the combinatorics and signs for the
moment, this operation on $q_{\gamma_1}\ldots q_{\gamma_m}$ has a 
geometric interpretation: it counts all \emph{potentially disconnected}
$J$-holomorphic curves of index~$1$ (i.e.~disjoint unions of $u$ with
trivial cylinders) that have $\gamma_1,\ldots,\gamma_m$ as their positive
asymptotic orbits; see Figure~\ref{fig:DSFT}.  
In other words, the action of $\mathbf{D}_\SFT$ 
on each monomial $q^{\boldsymbol{\gamma}}$ for 
$\boldsymbol{\gamma} = (\gamma_1,\ldots,\gamma_m)$ is determined by a
formula of the form
\begin{equation}
\label{eqn:DSFT}
\mathbf{D}_\SFT q^{\boldsymbol{\gamma}} =
\sum_{g=0}^\infty \sum_{A \in H_2(M)} 
\sum_{\boldsymbol{\gamma}'} \sum_{k=1}^m \hbar^{g + k - 1} e^A 
n_g(\boldsymbol{\gamma},\boldsymbol{\gamma}',k) q^{\boldsymbol{\gamma}'},
\end{equation}
where $n_g(\boldsymbol{\gamma},\boldsymbol{\gamma}',k)$ is a product of
some combinatorial factors with a signed count of generally disconnected
index~$1$ holomorphic curves of genus $g$ and homology class $A$ with 
positive ends at $\boldsymbol{\gamma}$ and negative ends at
$\boldsymbol{\gamma}'$, such that the nontrivial connected component has
exactly $k$ positive ends.  The presence of the combinatorial factors
hidden in $n_g(\boldsymbol{\gamma},\boldsymbol{\gamma}',k)$ is a slightly
subtle point which we will try to clarify in the following sections.

\section{Interlude: How to count points in an orbifold}
\label{sec:orbifold}

As in all versions of Floer theory, the proof that $\mathbf{H}^2 = 0$ is
based on the fact that certain moduli spaces are compact oriented
$1$-dimensional manifolds with boundary, and the signed count of their
boundary points is therefore zero.  We must be careful of course because,
strictly speaking, $\overline{\mM}(J)$ is not a manifold even when
Assumption~\ref{ass:optimism} holds---it is an orbifold, with the
possibility of singularities at multiply covered curves with nontrivial 
automorphism groups.  On the other hand, one can show that (after excluding
curves with bad asymptotic orbits) it is an \emph{oriented} orbifold,
and oriented $1$-dimensional orbifolds happen to be very simple objects:
since smooth finite group actions on $\RR$ cannot be nontrivial without
reversing orientation, all oriented $1$-dimensional orbifolds are actually
manifolds, suggesting the simple formula
$$
\text{``$\#\p\overline{\mM}_1(J) = 0$.''}
$$
I have placed this formula in quotation marks for a reason.
The reality of the situation is somewhat more complicated.

This is in fact where it becomes important to remember that
Assumption~\ref{ass:optimism}, in the way that we stated it, really is not
just science fiction but \emph{fantasy}: transversality is sometimes 
impossible to achieve for multiple covers, and we must therefore at least
have a sensible back-up plan for such cases.  To see the problem, remember
that our local structure theorem for $\mM(J)$ was proved by identifying
it in a neighborhood of any curve $u_0 : (\dot{\Sigma},j_0) \to (\RR \times M,J)$
with a set of the form
$$
\dbar_J^{-1}(0) \big/ G,
$$
where $\dbar_J : \tT \times \bB^{k,p,\delta} \to \eE^{k-1,p,\delta}$ is a smooth section of a
Banach space bundle $\eE^{k-1,p,\delta}$ over the product of a Teichm\"uller slice $\tT$
through~$j_0$ with a Banach manifold $\bB^{k,p,\delta}$ of maps $\dot{\Sigma} \to \RR \times M$,
and $G$ is the group of automorphisms of~$j_0$,
whose action on the base\footnote{As you may know if you've ever heard
a talk about polyfolds, there are some analytical
problems with this discussion if $G$ is a Lie group of positive dimension:
its action on the infinite-dimensional manifold $\bB^{k,p,\delta}$ of 
non-smooth maps cannot then be considered smooth in any
conventional sense.  This problem leads to the introduction
of \emph{sc-smooth} structures, cf.~\cite{HWZ:Fredholm1}.  There is no
problem however if $G$ is finite, e.g.~if the underlying Riemann surface is
stable, which we may as well assume for this discussion.}
$$
G \times (\tT \times \bB^{k,p,\delta}) \to \tT \times \bB^{k,p,\delta} : (\psi,(j,u)) \mapsto
(\psi^*j, u \circ \psi)
$$
preserves $\dbar_J^{-1}(0)$.  In fact, the action of $G$ on $\tT \times \bB^{k,p,\delta}$
is covered by a natural action on the
bundle~$\eE^{k-1,p,\delta}$, and the reason for it preserving the zero-set is that
$\dbar_J$ is an equivariant section,
$$
\dbar_J(\psi^*j,u \circ \psi) = \psi^*\dbar_J(j,u).
$$
If $G$ is finite, then another way to say this is that
$\dbar_J$ is a smooth Fredholm section of the infinite-dimensional
\defin{orbibundle} $\eE^{k-1,p,\delta} / G$ 
over the orbifold~$(\tT \times \bB^{k,p,\delta}) / G$, whose isotropy group at $(j_0,u_0)$ is
$\Aut(u_0)$.  This section is transverse to the zero-section if and only if
the usual regularity condition holds, making $\dbar_J^{-1}(0) / G$ a
suborbifold of $(\tT \times \bB^{k,p,\delta}) / G$ whose isotropy group at $(j_0,u_0)$
is some quotient of $\Aut(u_0)$.

\begin{remark}
Most sensible definitions of
the term \defin{orbifold} (cf.~\cites{AdemLeidaRuan,Davis:orbifolds,FukayaOno})
require local models of the form $\uU / G$,
where $\uU$ is a $G$-invariant open subset of a vector space on which the
finite group $G$ acts smoothly and \emph{effectively}---the latter condition 
is necessary in order to have isotropy groups that are well-defined up to
isomorphism at every point.  In the above example, $G$ acts effectively on
$\tT \times \bB^{k,p,\delta}$ but might have a nontrivial subgroup $H \subset G$ of
transformations that fix every element of $\dbar^{-1}_J(0)$, in which case the 
$G$-action on $\dbar^{-1}_J(0)$ can be replaced by an effective action
of $G / H$.  The isotropy group of $(j_0,u_0) \in \dbar^{-1}_J(0) / G$ is
then $\Aut(u_0) / (\Aut(u_0) \cap H)$.
\end{remark}

Now to see just how unreasonably optimistic Assumption~\ref{ass:optimism} is,
notice that it's easy to think up examples of smooth orbibundles in which
zeroes of sections can \emph{never} be regular if they have nontrivial
isotropy.

\begin{example}
\label{ex:trivialBndl}
Let $M = \CC / \ZZ_2$ with $\ZZ_2$ acting as the antipodal map, and consider
the trivial complex line bundle $E = M \times \CC = (\CC \times \CC) / \ZZ_2$,
where the $\ZZ_2$ action on $\CC \times \CC$ identifies $(z,v)$ with $(-z,v)$.
A smooth function $f : \CC \to \CC$ then represents a section of the
orbibundle $E \to M$ if and only if $f(z) = f(-z)$ for all~$z$.  This implies
that if $f(0) = 0$, then $df(0) = 0$.  It is possible to perturb $f$
generically to a section that is transverse to the zero-section, but such a
perturbation can never have zeroes at~$0$.

Of course, we do know how to assign $\ZZ$-valued orders to degenerate
zeroes of sections, e.g.~$f(z) = z^2$ defines a section of $E \to M$
with a zero of order~$2$ at~$0$.  Notice however that if we perturb
this to $f_\epsilon(z) = z^2 + \epsilon$ for $\epsilon > 0$ small, then 
$f_\epsilon$ has two simple zeroes at points near the origin, but they
are actually \emph{the same point} in $\CC / \ZZ_2$, giving a count of
only~$1$ zero.  This means that if we give the zero of $f$ at the origin
its full weight, then we are counting wrongly---the resulting count will not
be homotopy invariant.  The correct algebraic count of zeroes is evidently
\begin{equation}
\label{eqn:orbifoldCount}
\# f^{-1}(0) := \sum_{z \in f^{-1}(0) \subset M} \frac{\ord(f;z)}{\kappa_z} \in \QQ,
\end{equation}
where $\ord(f;z) \in \ZZ$ is the order of the zero (computed in the usual way 
as a winding number, or in higher dimensions as the degree of a map of spheres,
cf.~\cite{Milnor:differentiable}), and $\kappa_z \in \NN$ denotes the
order of the isotropy group at~$z$.
\end{example}

\begin{exercise}
\label{EX:orbifoldCount}
Convince yourself that for any smooth oriented orbibundle $E \to M$ of real
rank~$m$ over a compact, smooth and oriented $m$-dimensional orbifold $M$ 
without boundary,
the count \eqref{eqn:orbifoldCount} gives the same result for any section
with isolated zeroes.\footnote{If you're still not sure what an orbibundle is, a definition can be
found in \cite{FukayaOno}*{Chapter~1}.} \\
\textsl{Hint: The space of sections of an orbibundle is still a vector space,
so any two are homotopic.  Since $M$ and $[0,1]$ are both compact, it suffices
to focus on small perturbations of a single section on a single orbifold
chart.}
\end{exercise}

For a slightly different perspective on \eqref{eqn:orbifoldCount}, consider
the special case of a closed orbifold that is the quotient of a closed
manifold $\widetilde{M}$ by an effective orientation-preserving 
finite group action,
$$
M = \widetilde{M} / G.
$$
Suppose $\widetilde{E} \to \widetilde{M}$ is an oriented vector bundle with 
rank equal to $\dim M$, and $G$ also acts on $\widetilde{E}$ by 
orientation-preserving linear bundle maps that cover its action on
$\widetilde{M}$, so the quotient
$$
E = \widetilde{E} / G \to M
$$
is an orbibundle.  A section $f : M \to E$ is then equivalent to a
$G$-equivariant section $\tilde{f} : \widetilde{M} \to \widetilde{E}$,
and the signed count of zeroes
$$
\#\tilde{f}^{-1}(0) = \sum_{z \in \tilde{f}^{-1}(0) \subset \widetilde{M}}
\ord(\tilde{f};z) \in \ZZ
$$
is of course the same for any section that has only isolated zeroes.
It can also be expressed in terms of $f$ since any $z \in f^{-1}(0) \subset M$
has exactly $|G| / \kappa_z$ lifts to points in $\tilde{f}^{-1}(0) \subset
\widetilde{M}$, implying
$$
\#\tilde{f}^{-1}(0) = \sum_{z \in f^{-1}(0) \subset M} 
\frac{|G|}{\kappa_z} \ord(f;z)
$$
and thus $\#f^{-1}(0) = \frac{1}{|G|} \#\tilde{f}^{-1}(0)$.  The
invariance of \eqref{eqn:orbifoldCount} is now an immediate consequence
of the invariance of $\#\tilde{f}^{-1}(0)$, which follows from the
standard argument as in \cite{Milnor:differentiable}.

Now, if you enjoyed reading \cite{Milnor:differentiable} as much as I did, 
then it may seem tempting to try proving invariance of \eqref{eqn:orbifoldCount}
in general by choosing a generic homotopy $H : [0,1] \times M \to E$ between two
generic sections $f_0$ and $f_1$ and
showing that $H^{-1}(0) \subset [0,1] \times M$ is a compact
oriented $1$-dimensional orbifold with boundary.  As we observed at the
beginning of this section, $H^{-1}(0)$ is then actually a manifold, so the
signed count of its boundary points should be zero.  But this would give
the wrong result: it would suggest that $\sum_{z \in f^{-1}(0) \subset M} \ord(f;z)$
should be homotopy invariant, without the rational weights, and we've already seen
that this is not true.  What is going on here?  The answer is that
the homogopy $H$ cannot in general be made transverse to the zero-section,
now matter how generically we perturb it!  It is an illustration of the
fundamental conflict between the notions of \emph{genericity} and 
\emph{equivariance}.

\begin{example}
\label{ex:nontrivialBndl}
Let $M = \CC / \ZZ_2$ as in Example~\ref{ex:trivialBndl}, but define the
complex orbibundle $E \to M$ by
$$
E = (\CC \times \CC) \big/ (z,v) \sim (-z,-v),
$$
i.e.~the $\ZZ_2$-action also acts antipodally on fibers. 
Now a smooth function $f : \CC \to \CC$ defines a section of $E$ if and only if
$f(-z) = -f(z)$, hence \emph{all} such sections have a zero at the origin.
Compare the two sections
$$
f_0(x+iy) = x + iy, \qquad f_1(x+iy) = (x^3 - x) + iy.
$$
They have qualitatively the same behavior near infinity, meaning in particular
that they are homotopic through a family of sections whose zeroes are 
confined to some compact subset, thus we expect the algebraic count of zeroes
to be the same for both.  This is true if the count is defined by
\eqref{eqn:orbifoldCount}: we have $\#f_0^{-1}(0) = \#f_1^{-1}(0) = \frac{1}{2}$,
in particular the negative zero of $f_1$ at the origin counts for $-1/2$
while the positive zero at $(1,0) \sim (-1,0)$ counts for~$1$.  We see that
the inclusion of the rational weights $\frac{1}{\kappa_x}$ is crucial
for this result.  Notice that if $H : [0,1] \times M \to E$ is a homotopy
of sections from $f_0$ to $f_1$, then $H(\tau,0) = 0$ for all $\tau$, thus
$\p_\tau H(\tau,0)$ vanishes and
$$
dH(\tau,0) = d f_\tau(0)
$$
where $f_\tau = H(\tau,\cdot)$.  But $df_\tau(0)$ cannot be an isomorphism
for all $\tau \in (0,1)$ since $df_0(0)$ preserves orientation while
$df_1(0)$ reverses it.  This is not a problem that can be fixed by
making $H$ more generic---the homotopy will never be transverse to the
zero-section, no matter what we do.
\end{example}

The need to address issues of the type raised by the above examples
leads naturally to the notion of \defin{multisections} as outlined in
\cite{Salamon:Floer}*{\S 5} and \cite{FukayaOno}, and this is a major
feature of the analysis under development by Hofer-Wysocki-Zehnder,
see for example \cite{HWZ:integration}.  In Example~\ref{ex:nontrivialBndl}
for instance, one can consider functions
$$
f : \CC \to \Sym_2(\CC) := (\CC \times \CC) \big/ (z_1,z_2) \sim (z_2,z_1),
$$
which can be regarded as doubly-valued sections of $E \to M$ if $f$ is
$\ZZ_2$-equivariant for the antipodal action of $\ZZ_2$ on the symmetric
product $\Sym_2(\CC)$.  Such a section is considered single-valued at any
point $z$ where $f(z)$ is of the form $[(v,v)]$, so one can now imagine
homotopies from $f_0$ to $f_1$ through doubly-valued sections.  One advantage 
of this generalization is that $f$ can now take nonzero values of the form
$[(v,-v)]$ at the origin, e.g.~if $g : \CC \to \CC$ is any odd function, then
$$
f(z) := [(g(z) + c,g(z) - c)]
$$
is a well-defined multisection for every $c \in \CC$.  

\begin{exercise}
\label{EX:multiHomotopy}
Find a homotopy between the sections $f_0$ and
$f_1$ of Example~\ref{ex:nontrivialBndl} through doubly-valued sections,
such that the homotopy is transverse to the zero-section.
\end{exercise}

You may notice if you work out Exercise~\ref{EX:multiHomotopy} that the
zero set of the homotopy in $[0,1] \times M$ is still not submanifold
or suborbifold.  Instead, it naturally carries the structure of a 
\emph{weighted branched manifold with boundary}.  The rational weights
attached to every point in this object can be used to explain the
weights appearing in \eqref{eqn:orbifoldCount} and thus give a
Milnor-style proof that $\#f^{-1}(0) \in \QQ$ is invariant.

We will not discuss multisections or weighted branched manifolds any further,
but the main takeaway from this discussion should be that the ``right'' way
to count $0$-dimensional orbifolds algebraically is always some version
of \eqref{eqn:orbifoldCount}, and the count in general is a
rational number, not an integer.  We've discussed this above from the
perspective of obtaining a homotopy-invariant count, but the same logic
applies to any Floer-type theory since the relation $\p^2 = 0$ is typically
based on similar arguments via $1$-dimensional moduli spaces with boundary.
While a more simplistic notion of counting
may produce well-defined homology theories in isolated cases where
Assumption~\ref{ass:optimism} holds (e.g.~in \cite{Nelson:Abendblatt}),
we cannot expect it to generalize beyond these cases, due to the
fundamental conflict between transversality and equivariance.  On the
other hand, it will be possible in our situation to remove isotropy from
the picture by lifting to moduli spaces with asymptotic markers; the
moduli space we're interested in is always the quotient of this larger
space by a finite group action, so the situation is analogous to replacing
an orbibundle $E = \widetilde{E} / G \to \widetilde{M} / G$ by an
ordinary vector bundle $\widetilde{E}$ over a manifold~$\widetilde{M}$.
In the infinite-dimensional setting, transversality is still a hard
problem, but having lifted to a manifold and thus removed the need for
equivariance, there is no longer any \emph{a priori} reason why it cannot
be solved by choosing sufficiently generic perturbations.
This makes counting curves with rational weights seem a much
more promising method for defining invariants, and we will adopt this 
perspective in the discussion to follow.

\section{Cylindrical contact homology revisited}
\label{sec:cyl}

Under an extra assumption on the complex $(\aA[[\hbar]],\mathbf{D}_\SFT)$,
we can recover from it a more general version of the cylindrical contact
homology we saw in Lecture~\ref{lec:tight3tori}.  Suppose in particular that there are no
index~$1$ holomorphic planes in $\RR \times M$, so every term in
$\hbar \mathbf{H}$ has at least one factor of either $\hbar$ or one of the
$q_\gamma$ variables.  Then 
$$
\mathbf{D}_\SFT = \sum_{\gamma,\gamma',A} \kappa_\gamma \left( \sum_{u \in \mM_{0,0}(J,A,\gamma,\gamma')/\RR}
\frac{\epsilon(u)}{|\Aut(u)|}  e^A q_{\gamma'} \frac{\p}{\p q_\gamma} \right) +
\ldots,
$$
where the first sum is over all pairs of good Reeb orbits $\gamma$ and~$\gamma'$,
and the ellipsis is a sum of terms that all include at least a positive 
power of $\hbar$ or two $q_\gamma$ variables or two partial derivatives.
Let us abbreviate the spaces $\mM_{0,0}(J,A,\gamma,\gamma') / \RR$ of 
$\RR$-equivalence classes of $J$-holomorphic
cylinders by $\mM_A(\gamma,\gamma')$, and notice that for any
$u \in \mM_A(\gamma,\gamma')$, the automorphism group is a cyclic
group of order equal to the covering multiplicity
$$
|\Aut(u)| = \kappa_u := \cov(u) \in \NN.
$$
Thus for any single generator $q_\gamma$, we have
$$
\mathbf{D}_\SFT q_\gamma = \p_\CCH q_\gamma + O(|q|^2,\hbar),
$$
where 
\begin{equation}
\label{eqn:CCH}
\p_\CCH q_\gamma := \kappa_\gamma \sum_{\gamma',A} \left(\sum_{u \in \mM_A(\gamma,\gamma')}
\frac{\epsilon(u)}{\kappa_u} \right) e^A q_{\gamma'} .
\end{equation}
The fact that $\mathbf{D}_\SFT^2 = 0$ thus implies 
$$
\p_\CCH^2 = 0,
$$
and the homology of the graded $R$-module generated by 
$\{ q_\gamma\ |\ \text{$\gamma$ good} \}$
with differential $\p_\CCH$ is an obvious generalization of the cylindrical
contact homology from Lecture~\ref{lec:tight3tori}.  What we saw there was a special case of this
where the combinatorial factor $\kappa_\gamma / \kappa_u$ did not appear
because we were restricting to a homotopy class in which all orbits were
simply covered, and all holomorphic cylinders were thus somewhere injective.

The presence of the factor $\kappa_\gamma / \kappa_u$ deserves further 
comment.  According to the above formula, we have
$$
\p_\CCH^2 q_\gamma = \sum_{\gamma',\gamma'',A,A'} \sum_{u \in \mM_A(\gamma,\gamma')}
\sum_{v \in \mM_{A'}(\gamma',\gamma'')} e^{A+A'} \frac{\kappa_\gamma \kappa_{\gamma'} \epsilon(u) \epsilon(v)}
{\kappa_u \kappa_v} q_{\gamma''},
$$
hence $\p_\CCH^2 = 0$ holds if and only if for all $A \in H_2(M)$ and all pairs of
good orbits $\gamma_+,\gamma_-$,
\begin{equation}
\label{eqn:brokenCount}
\sum_{\gamma_0} \sum_{B + C = A} \left( \sum_{(u,v) \in \mM_B(\gamma_+,\gamma_0) \times
\mM_C(\gamma_0,\gamma_-)}
\frac{\kappa_{\gamma_0}}{\kappa_u \kappa_v}  \epsilon(u) \epsilon(v) \right) = 0.
\end{equation}
If $\gamma_+$ and $\gamma_-$ happen to be simply covered orbits, then
$u$ and $v$ in this expression always have trivial automorphism groups and
it is clear what this sum means: every such pair $(u,v) \in 
\mM_B(\gamma_+,\gamma_0) \times \mM_C(\gamma_0,\gamma_-)$ corresponds to
exactly $\kappa_{\gamma_0}$ distinct holomorphic buildings obtained by
different choices of decoration, so \eqref{eqn:brokenCount} is the count
of boundary points of the compactified $1$-dimensional manifold of index~$2$ 
cylinders $\mM_A(\gamma_+,\gamma_-) / \RR$.  This sum skips over
all bad orbits~$\gamma_0$, but this is fine because whenever the breaking
orbit is bad, there are evenly many choices of decoration such that half
of these choices cancel the other half when counted with the correct signs.

To understand why this formula is still correct in the presence of 
automorphisms, let us outline two equivalent approaches.

The easiest option is to instead consider moduli spaces with asymptotic
markers, which never have automorphisms:
removing unnecessary factors of $\kappa_{\gamma_+}$ and $\kappa_{\gamma_-}$
then transforms \eqref{eqn:brokenCount} into
$$
\sum_{\gamma_0} \sum_{B+C = A} \frac{1}{\kappa_{\gamma_0}}
\# \mM_B^\$(\gamma_+,\gamma_0) \cdot \# \mM_C^\$(\gamma_0,\gamma_-) = 0.
$$
Now since each pair $(u,v) \in \mM_B^\$(\gamma_+,\gamma_0) \times
\mM_C^\$(\gamma_0,\gamma_-)$ carries a canonical decoration and thus
determines a holomorphic building, the division by $\kappa_{\gamma_0}$
accounts for the fact that 
$\# \mM_B^\$(\gamma_+,\gamma_0) \cdot \# \mM_C^\$(\gamma_0,\gamma_-)$
overcounts the set of broken cylinders from $\gamma_+$ to $\gamma_-$
with asymptotic markers at $\gamma_\pm$ by precisely this factor, as a
simultaneous adjustment of the marker at $\gamma_0$ in both 
$u \in \mM_B^\$(\gamma_+,\gamma_0)$ and $v \in \mM_C^\$(\gamma_0,\gamma_-)$
produces the same decoration and therefore the same building.

The following alternative perspective will be more useful when we generalize
beyond cylinders in the next section.  We can directly count points in
$\p\overline{\mM}_A(\gamma_+,\gamma_-)$, though as we saw in \S\ref{sec:orbifold},
rational weights should be included in the count whenever there is isotropy.
Let us write 
$$
\mM_A(\gamma_+,\gamma_-) = \mM_A^\$(\gamma_+,\gamma_-) / G,
$$ 
where $G \cong \ZZ_{\kappa_{\gamma_+}} 
\times \ZZ_{\kappa_{\gamma_-}}$ is a finite group acting by adjustment of
the asymptotic markers.  Since $\overline{\mM}_A^\$(\gamma_+,\gamma_-)$
is a compact oriented $1$-manifold with boundary under Assumption~\ref{ass:optimism}, 
the signed count of its
boundary points is~$0$.  We can ignore buildings broken along bad orbits
in this count, since these always come in cancelling pairs.
Let us now transform this into a count of buildings
$(u|\Phi|v) \in \p\overline{\mM}_A(\gamma_+,\gamma_-)$ broken along 
good orbits~$\gamma_0$: here $u \in \mM_B(\gamma_+,\gamma_0)$ and 
$v \in \mM_C(\gamma_0,\gamma_-)$ for some homology classes with $B+C=A$,
and $\Phi$ is a decoration which describes how to glue the ends
of $u$ and $v$ at~$\gamma_0$.  The automorphism group of such a building is
the subgroup
$$
\Aut(u|\Phi|v) \subset \Aut(u) \times \Aut(v)
$$
consisting of all pairs $(\varphi,\psi) \in \Aut(u) \times \Aut(v)$ that
define the same rotation at the two punctures asymptotic to~$\gamma_0$;
note that this group does not actually depend on the decoration~$\Phi$.
Since we're talking about cylinders, we can be much more specific:
we have $\Aut(u) = \ZZ_{\kappa_u}$ and $\Aut(v) = \ZZ_{\kappa_v}$, and if both
are regarded as subgroups of $\U(1)$,
$$
\Aut(u|\Phi|v) = \ZZ_{\kappa_u} \cap \ZZ_{\kappa_v} = \ZZ_{\gcd(\kappa_u,\kappa_v)},
$$
which is injected into $\Aut(u) \times \Aut(v)$ by $\psi \mapsto (\psi,\psi)$.
The boundary of $\overline{\mM}_A^\$(\gamma_+,\gamma_-)$ can be understood
likewise as a space of equivalence classes
$$
[(u,v)] \in \left( \mM_B^\$(\gamma_+,\gamma_0) \times 
\mM_C^\$(\gamma_0,\gamma_-) \right) \big/ \sim,
$$
where two such pairs are equivalent if their asymptotic markers at the ends
asymptotic to $\gamma_0$ determine the same decoration.
Now observe that the group $G \cong \ZZ_{\kappa_{\gamma_+}} \times \ZZ_{\kappa_{\gamma_-}}$
also acts on buildings in $\p\overline{\mM}_A^\$(\gamma_+,\gamma_-)$, and the
stabilizer of this action at $(u,v)$ is $\Aut(u|\Phi|v)$, hence 
each $(u|\Phi|v) \in \p\overline{\mM}_A(\gamma_+,\gamma_-)$ gives rise to
$\frac{|G|}{\gcd(\kappa_u,\kappa_v)}$ terms in the count of
$\p\overline{\mM}_A^\$(\gamma_+,\gamma_-)$, implying
\begin{equation}
\label{eqn:kappas}
\sum_{(u|\Phi|v) \in \p\overline{\mM}_A(\gamma_+,\gamma_-)}
\frac{\epsilon(u) \epsilon(v)}{\gcd(\kappa_u,\kappa_v)} = 0.
\end{equation}
Finally, notice that while each pair $(u,v) \in \mM_B(\gamma_+,\gamma_0) \times
\mM_C(\gamma_0,\gamma_-)$ determines buildings with $\kappa_{\gamma_0}$
distinct choices of decoration, some of these buildings may be equivalent:
every pair of automorphisms $(\varphi,\psi) \in \Aut(u) \times \Aut(v)$
transforms a building $(u|\Phi|v)$ by potentially changing the 
decoration~$\Phi$, thus producing an equivalent building.  This action on
buildings is trivial if and only if $(\varphi,\psi) \in \Aut(u|\Phi|v)$,
hence every pair $(u,v) \in \mM_B(\gamma_+,\gamma_0) \times
\mM_C(\gamma_0,\gamma_-)$ gives rise to exactly
$$
\frac{\kappa_{\gamma_0}}{\left| \left(\Aut(u) \times \Aut(v)\right) \big/
\Aut(u|\Phi|v) \right|} =
\frac{\kappa_{\gamma_0} \gcd(\kappa_u,\kappa_v)}{\kappa_u \kappa_v}
$$
elements of $\p\overline{\mM}_A(\gamma_+,\gamma_-)$, so that
\eqref{eqn:kappas} becomes
\begin{multline*}
\sum_{\gamma_0} \sum_{B+C=A} \left(\sum_{(u,v) \in \mM_B(\gamma_+,\gamma_0) \times
\mM_C(\gamma_0,\gamma_-)} \frac{\epsilon(u)\epsilon(v)}{\gcd(\kappa_u,\kappa_v)}
\frac{\kappa_{\gamma_0} \gcd(\kappa_u,\kappa_v)}{\kappa_u \kappa_v} \right) \\
= \sum_{\gamma_0} \sum_{B+C=A} \left( \sum_{(u,v) \in \mM_B(\gamma_+,\gamma_0) \times
\mM_C(\gamma_0,\gamma_-)} \frac{\epsilon(u)\epsilon(v) \kappa_{\gamma_0}}{\kappa_u \kappa_v}\right) = 0,
\end{multline*}
reproducing \eqref{eqn:brokenCount}.

\section{Combinatorics of gluing}
\label{sec:combi}

Now let's try to justify the formula $\mathbf{H}^2 = 0$.  The product of
$\mathbf{H}$ with itself is the formal sum over all pairs of index~$1$ curves
$u,v \in \mM_1^\sigma(J) / \RR$ of certain monomials: in particular
if these two curves respectively have genus $g_u$ and $g_v$, homology
classes $A_u$ and $A_v$, and asymptotic orbits $\boldsymbol{\gamma}_u^\pm$ and
$\boldsymbol{\gamma}_v^\pm$, then the corresponding term in $\mathbf{H}^2$ is
$$
\frac{\epsilon(u) \epsilon(v)}{|\Aut^\sigma(u)| |\Aut^\sigma(v)|}
\hbar^{g_u + g_v - 2} e^{A_u+A_v} q^{\boldsymbol{\gamma}^-_u}
p^{\boldsymbol{\gamma}^+_u} q^{\boldsymbol{\gamma}^-_v}
p^{\boldsymbol{\gamma}^+_v}.
$$
Before we can add up all monomials of this form, we need to put all the
$q$ and $p$ variables in the same order: within each of the products
$q^{\boldsymbol{\gamma}^-_u}$, $p^{\boldsymbol{\gamma}^+_u}$ and so
forth this is simply a matter of permuting the variables and changing
signs as appropriate, but the interesting part is the product
$p^{\boldsymbol{\gamma}^+_u} q^{\boldsymbol{\gamma}^-_v}$, for which we can
apply the commutation relations \eqref{eqn:quantum} to put all $q$ variables
before all $p$ variables.  Before discussing how this works in general,
let us consider a more specific example.

\begin{figure}
\includegraphics[width=6in]{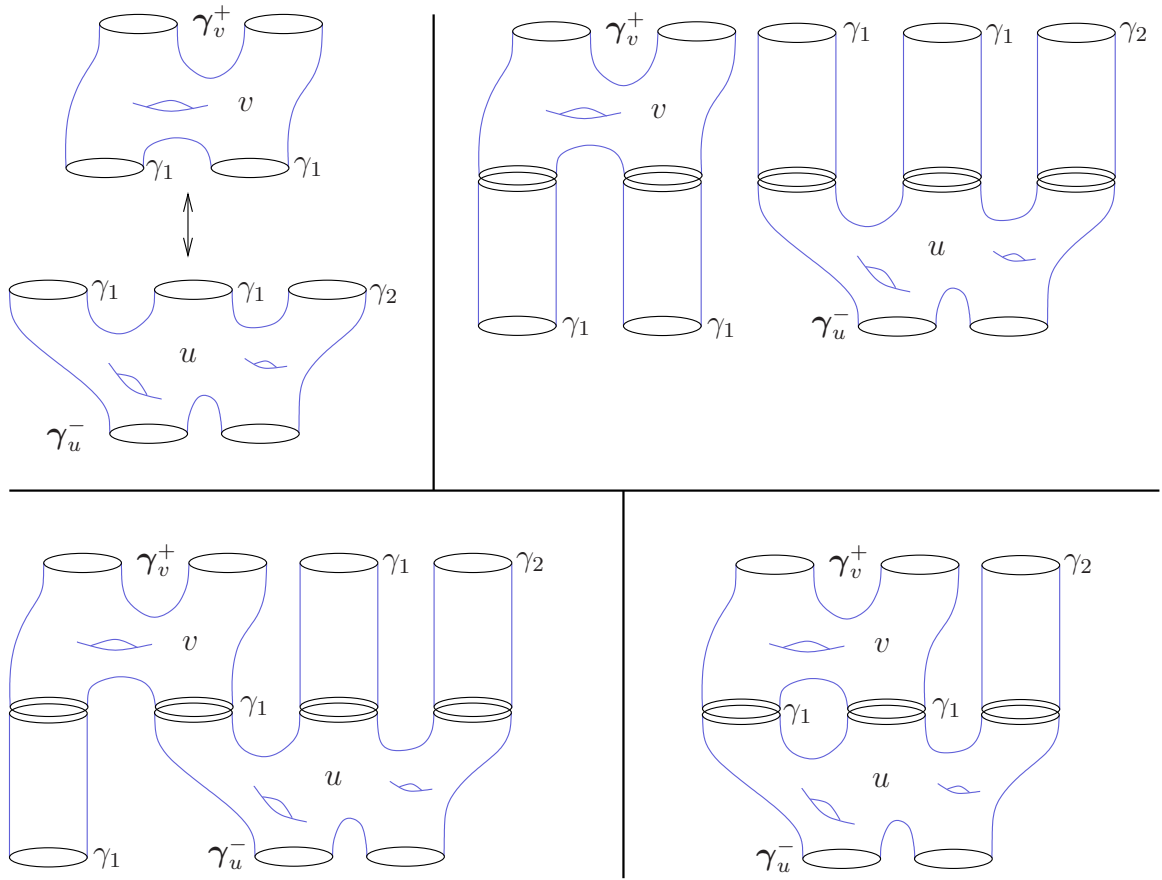}
\caption{\label{fig:glueCombi} Three possible ways of gluing the curves
$u$ and $v$ along with trivial cylinders to form index~$2$ curves.}
\end{figure}

Assume $\gamma_i$ for $i=1,2$ are two specific orbits with
$n-3 + \muCZ(\gamma_i)$ even, so the corresponding $q$ and $p$ variables
have even degree, and suppose
$$
\boldsymbol{\gamma}^+_u = (\gamma_1,\gamma_1,\gamma_2), \qquad
\boldsymbol{\gamma}^-_v = (\gamma_1,\gamma_1).
$$
After applying the relation $p_{\gamma_1} q_{\gamma_1} = q_{\gamma_1}
p_{\gamma_1} + \kappa_{\gamma_1} \hbar$ a total of five times, one
obtains the expansion
$$
p_{\gamma_1} p_{\gamma_1} p_{\gamma_2} q_{\gamma_1} q_{\gamma_1} = 
q_{\gamma_1}^2 p_{\gamma_1}^2 p_{\gamma_2} +
4 \kappa_{\gamma_1} \hbar q_{\gamma_1} p_{\gamma_1} p_{\gamma_2} +
2 \kappa_{\gamma_1}^2 \hbar^2 p_{\gamma_2},
$$
thus contributing a total of three terms to $\mathbf{H}^2$, namely the
products of the factor 
$\frac{\epsilon(u) \epsilon(v)}{|\Aut(u)| |\Aut(v)|} e^{A_u+A_v}$
with each of the expressions
\begin{align}
\label{eqn:gluing1}
& \hbar^{g_u + g_v - 2} q^{\boldsymbol{\gamma}^-_u} q_{\gamma_1}^2 p_{\gamma_1}^2
p_{\gamma_2} p^{\boldsymbol{\gamma}^+_v}, \\
\label{eqn:gluing2}
& 4 \kappa_{\gamma_1} \hbar^{g_u + g_v - 1} q^{\boldsymbol{\gamma}^-_u}
q_{\gamma_1} p_{\gamma_1} p_{\gamma_2} p^{\boldsymbol{\gamma}^+_v}, \\
\label{eqn:gluing3}
& 2 \kappa_{\gamma_1}^2 \hbar^{g_u + g_v} q^{\boldsymbol{\gamma}^-_u}
p_{\gamma_2} p^{\boldsymbol{\gamma}^+_v}.
\end{align}
As shown in Figure~\ref{fig:glueCombi}, this sum of three terms can be
interpreted as the count of all possible holomrphic buildings obtained by
gluing $v$ on top of $u$ together with a collection of trivial cylinders.
Indeed, since $\boldsymbol{\gamma}^+_u$ and $\boldsymbol{\gamma}^-_v$
include two matching orbits (which also happen to be the same one), there
are several choices to be made:
\begin{enumerate}
\item The top-right picture shows what we
might call the ``stupid gluing,'' in which no ends of $u$ are matched with
any ends of $v$, but all are instead glued to trivial cylinders, thus
producing a disconnected building.  This possibility is encoded by
\eqref{eqn:gluing1}, and we will see that in the total sum forming 
$\mathbf{H}^2$, this term gets cancelled out by a similar term for the
stupid gluing of $u$ on top of~$v$.  
\item The lower-left picture shows
the building obtained by gluing one end of $u$ to an end of $v$ along
the matching orbit~$\gamma_1$.  This option is encoded by \eqref{eqn:gluing2},
where the factor $4\kappa_{\gamma_1}$ appears because there are precisely
$4\kappa_{\gamma_1}$ distinct buildings of this type: indeed, there are
four choices of which end of $u$ should be glued to which end of~$v$,
and for each of these, a further $\kappa_{\gamma_1}$ choices of the decoration.
The arithmetic genus of the resulting building is $g_u + g_v$, as represented
by the factor $\hbar^{g_u + g_v - 1}$.
\item The lower-right picture is encoded by \eqref{eqn:gluing3}:
here there are two choices of bijections between the two pairs of punctures
asymptotic to~$\gamma_1$, and taking the choices of decoration at each
breaking orbit into account, we obtain the combinatorial 
factor~$2\kappa_{\gamma_1}^2$.  The presence of two nontrivial breaking
orbits increases the arithmetic genus to $g_u + g_v + 1$, as encoded in the
factor $\hbar^{g_u + g_v}$.
\end{enumerate}

You may now be able to extrapolate from the above example why the
commutator algebra we've defined encodes gluing of holomorphic curves in
the symplectization and thus leads to the relation $\mathbf{H}^2 = 0$.
Think of the algorithm by which you change 
$q^{\boldsymbol{\gamma}^-_u} p^{\boldsymbol{\gamma}^+_u}
q^{\boldsymbol{\gamma}^-_v} p^{\boldsymbol{\gamma}^+_u}$ 
into a sum of products with all $q$'s
appearing before $p$'s: for the first $q$ you see appearing after a~$p$, move 
it past each $p$ for different orbits (changing signs as necessary) until it
encounters a $p$ for the \emph{same} orbit.  Now you replace
$p_\gamma q_\gamma$ with $(-1)^{|p_\gamma| |q_\gamma|} q_\gamma p_\gamma +
\kappa_\gamma \hbar$, turning one product into a sum of two.
This represents a choice between two options: either you move $q_\gamma$
past $p_\gamma$ and apply the usual sign change, or you eliminate them
both but replace them with the combinatorial factor $\kappa_\gamma$ and an
extra~$\hbar$.  Then you continue this process until all $q$'s appear
before all~$p$'s.

The key point is that the process of gluing $v$ on top of $u$ in all
possible ways is governed by \emph{exactly the same algorithm}: first
consider the disjoint union of the two curves as a single disconnected
curve, with its punctures ordered in the same way in which their orbits
appear in the monomial.  Now reorder negative punctures of $v$ and
positive punctures of $u$, changing orientations as appropriate, until
you see two such punctures next to each other approaching the same 
orbit~$\gamma$.  Here you have two options:
either glue them together, or don't glue them but exchange their order.
If you exchange the order, then you may again have to change orientations
(depending on the parity of $n-3 + \muCZ(\gamma)$), but if you glue,
then you have $\kappa_\gamma$ distinct choices of decoration and will also
increase the arithmetic genus of the eventual building by~$1$.
In this way, every individual term in the final expansion of
$q^{\boldsymbol{\gamma}^-_u} p^{\boldsymbol{\gamma}^+_u}
q^{\boldsymbol{\gamma}^-_v} p^{\boldsymbol{\gamma}^+_u}$ 
represents a particular choice of which positive of ends of $u$ should
or should not be glued to which negative ends of~$v$.  Additional
factors of $\hbar$ appear to keep track of the increase in arithmetic
genus, and covering multiplicities of the breaking orbits also appear due to 
distinct choices of decorations.  At the end these must still be divided by 
orders of automorphism groups in order to avoid counting equivalent buildings
separately.  Fleshing out these details leads to the following explanation
for the relation $\mathbf{H}^2 = 0$:

\begin{prop}
\label{prop:H2}
Let $\p\overline{\mM}_2^\sigma(J)$ denote the space of two-level holomorphic
buildings in $\overline{\mM}(J)$ that have total index~$2$ and no bad 
asymptotic or breaking orbits, 
divided by the equivalence relation that forgets the order of the punctures.  Then
$$
\mathbf{H}^2 =
\sum_{\mathbf{u} \in \p\overline{\mM}^\sigma_2(J)} 
\frac{\epsilon(\mathbf{u})}{|\Aut^\sigma(\mathbf{u})|} 
\hbar^{g-1} e^A q^{\boldsymbol{\gamma}^-} p^{\boldsymbol{\gamma}^+},
$$
where the terms in each monomial are determined by
$\mathbf{u} \in \p\overline{\mM}^\sigma_2(J)$ as follows:
\begin{enumerate}
\item $g$ is the arithmetic genus of~$\mathbf{u}$;
\item $A$ is the equivalence class of $[\mathbf{u}] \in H_2(M)$ in
$H_2(M) / G$;
\item $\boldsymbol{\gamma}^\pm = (\gamma_1^\pm,\ldots,\gamma_{k_\pm}^\pm)$
are the asymptotic orbits of $\mathbf{u}$ after arbitrarily fixing orderings
of its positive and negative punctures;
\item $\epsilon(\mathbf{u}) \in \{1,-1\}$ is the boundary orientation at
$\mathbf{u}$ determined by the chosen coherent orientations on $\mM^\$(J)$.
Specifically, given the chosen ordering of the punctures and an arbitrary
choice of asymptotic markers at each puncture, $\mathbf{u}$ determines a
boundary point of a $1$-dimensional connected component of $\overline{\mM}^\$(J)$,
and we define $\epsilon(\mathbf{u}) = +1$ if and only if the orientation
of $\overline{\mM}^\$(J)$ at this point is outward.
\end{enumerate}
\end{prop}
Once again $\epsilon(\mathbf{u})$ and $q^{\boldsymbol{\gamma}^-} p^{\boldsymbol{\gamma}^+}$
change signs in the same way under any reordering of the punctures, so
their product is well defined, and there is no dependence on choices of
markers since bad orbits have been excluded.

\begin{proof}[Proof of Proposition~\ref{prop:H2}]
Our original formula for $\mathbf{H}$ gives rise to an expansion
$$
\mathbf{H}^2 = \sum_{(u,v) \in \mM_1^\sigma(J) / \RR \times \mM_1^\sigma(J) / \RR}
\frac{\epsilon(u) \epsilon(v)}{|\Aut^\sigma(u)| |\Aut^\sigma(v)|} \hbar^{g_u + g_v - 2}
e^{A_u + A_v} q^{\boldsymbol{\gamma}^-_u} p^{\boldsymbol{\gamma}^+_u}
q^{\boldsymbol{\gamma}^-_v} p^{\boldsymbol{\gamma}^+_v}.
$$
As explained in the previous paragraph, the process of reordering
$p^{\boldsymbol{\gamma}^+_u} q^{\boldsymbol{\gamma}^-_v}$ to put all $q$'s
before $p$'s produces an expansion, each term of which can be identified with
a specific choice of which positive punctures of $u$ should be glued to
which negative punctures of~$v$.  If $k$ punctures are glued, then 
the resulting power of $\hbar$ is $g_u + g_v - 2 + k$, corresponding to the
fact that the resulting building has arithmetic genus $g_u + g_v + k - 1$.
We claim that the term for $k=0$ is cancelled out by the corresponding
term of $\mathbf{H}^2$ that has the roles of $u$ and $v$ reversed.
To see this, imagine first the case where $u$ and $v$ have no asymptotic
orbits in common, hence no nontrivial gluings are possible and all the
$q$ and $p$ variables in the expression supercommute with each other.
Then since both curves have index~$1$, the monomials
$q^{\boldsymbol{\gamma}^-_u} p^{\boldsymbol{\gamma}^+_u}$ and
$q^{\boldsymbol{\gamma}^-_v} p^{\boldsymbol{\gamma}^+_v}$ must both have
odd degree, implying
$$
q^{\boldsymbol{\gamma}^-_u} p^{\boldsymbol{\gamma}^+_u}
q^{\boldsymbol{\gamma}^-_v} p^{\boldsymbol{\gamma}^+_v} =
- q^{\boldsymbol{\gamma}^-_v} p^{\boldsymbol{\gamma}^+_v} 
q^{\boldsymbol{\gamma}^-_u} p^{\boldsymbol{\gamma}^+_u}
$$
and thus the desired cancellation.  If $u$ and $v$ do have orbits in common,
then the result for the $k=0$ terms is still not any different from this:
all signs still change in the same way when applying 
$[p_\gamma,q_\gamma] = \kappa_\gamma \hbar$ to change $p_\gamma q_\gamma$
into $q_\gamma p_\gamma$, we simply ignore the extra term $\kappa_\gamma \hbar$
since it is only relevant for gluings with $k > 0$.  This proves the claim,
and consequently, that the expansion resulting from the curves $u$ and $v$
has no term containing $\hbar^{g_u + g_v - 2}$.

The combinatorial factors can be explained as follows.  The commutator
expansion for $p^{\boldsymbol{\gamma}^+_u} q^{\boldsymbol{\gamma}^-_v}$
automatically produces combinatorial factors that count the different
possible gluings, but if $u$ and $v$ have automorphisms, then not all of
these give inequivalent buildings.  This part of the discussion is a
straightforward extension of what we did for cylindrical contact homology
at the end of \S\ref{sec:cyl}.
Indeed, the actual set of inequivalent
buildings is the quotient of this larger set by an action of
$$
\left( \Aut^\sigma(u) \times \Aut^\sigma(v) \right) \big/ \Aut^\sigma(\mathbf{u}),
$$
where for a building $\mathbf{u}$ formed by endowing the pair $(u,v)$ with
decorations, $\Aut^\sigma(\mathbf{u})$ denotes the subgroup consisting of pairs
$(\varphi,\psi) \in \Aut^\sigma(u) \times \Aut^\sigma(v)$ that preserve
pairs of breaking punctures along with their decorations.  This is what
changes the factor $\frac{1}{|\Aut^\sigma(u)| |\Aut^\sigma(v)|}$ into
$\frac{1}{|\Aut^\sigma(\mathbf{u})|}$ as in the statement of the proposition.
\end{proof}

The theorem that $\mathbf{H}^2 = 0$ now follows once you believe the
propaganda from \S\ref{sec:orbifold}, arguing that
$\sum_{\mathbf{u} \in \p\overline{\mM}^\sigma_2(J)} \frac{\epsilon(\mathbf{u})}{|\Aut^\sigma(\mathbf{u})|}$
is the correct way to count the boundary points of $\overline{\mM}^\sigma_2(J)$.
As we did with cylindrical contact homology, we can use the obvious
projection $\overline{\mM}^\$(J) \to \overline{\mM}^\sigma(J)$ to reduce
this to the fact that if the $1$-dimensional components of
$\overline{\mM}^\$(J)$ are manifolds (which is true if Assumption~\ref{ass:optimism}
holds), then the integer-valued signed count of their boundary points vanishes.

\section{Some remarks on torsion, coefficients, and conventions}

\subsection{What if $H_1(M)$ has torsion?}
\label{sec:torsion}

The main consequence for SFT if $H_1(M)$ has torsion is that one cannot
define an integer grading, though there is always a canonical 
$\ZZ_2$-grading.\footnote{In fact there is a bit more than a $\ZZ_2$-grading,
see \cite{SFT}*{\S 2.9.1}.}
The setup in \S\ref{sec:setup3} must now be modified as follows.
The \defin{reference curves}
$$
C_1,\ldots,C_r \subset M
$$
are required to form a basis of $H_1(M) / \text{torsion}$, so for every
integral homology class $[\gamma]$, there is a unique collection of integers
$m_1,\ldots,m_r$ such that $[\gamma] = \sum_i m_i [C_i] \in H_1(M;\QQ)$.
Instead of spanning surfaces for each orbit, one can define
\defin{spanning chains} $C_\gamma$, which are singular $2$-chains with
rational coefficients satisfying
$$
\p C_\gamma = \sum_i m_i [C_i] - [\gamma]
$$
for the aforementioned set of integers $m_i \in \ZZ$.  Note that $C_\gamma$
must in general have nonintegral coefficients since $\sum_i m_i [C_i]$ and
$[\gamma]$ might not be homologous in $H_1(M;\ZZ)$, so $C_\gamma$ cannot
always be represented by a smooth map of a surface.  One consequence of
this is that the absolute homology class associated to an asymptotically
cylindrical holomorphic curve $u : \dot{\Sigma} \to \RR \times M$ will now
be rational,
$$
[u] \in H_2(M;\QQ),
$$
and we must therefore take $G$ to be a linear subspace
$$
G \subset H_2(M;\QQ).
$$
Another consequence is that
we cannot use capping chains to transfer trivializations from the reference
curves to the orbits, so there is no natural way to define $\muCZ(\gamma)$
as an integer.  The easiest thing to do instead is to take the mod~2
Conley-Zehnder index
$$
\muCZ(\gamma) \in \ZZ_2
$$
and define all degrees of generators as either even or odd with no further
distinction.  In particular, we now have
$$
|q_\gamma| = n - 3 + \muCZ(\gamma) \in \ZZ_2, \qquad
|p_\gamma| = n - 3 - \muCZ(\gamma) \in \ZZ_2,
$$
while $\hbar$ and all elements of $R = \QQ[H_2(M;\QQ) / G]$ are even.
With these modifications, the rest of the discussion also becomes valid
for the case where $H_1(M)$ has torsion, and leads to $\ZZ_2$-graded
contact invariants.

\subsection{Combinatorial conventions}

The combinatorial factors appearing in our definition of $\mathbf{H}$
may at first look slightly different from what appears elsewhere in the
literature.  Actually, most papers seem to agree on this detail, but various
subtle differences and ambiguities in notation mean that it sometimes requires 
intense concentration to recognize this fact.

The original propaganda paper \cite{SFT} expresses
everything in terms of moduli spaces with asymptotic markers, and the
formula for $\mathbf{H}$ in \S 2.2.3 of that paper (which is expressed
in a slightly more general form involving marked points) agrees with
our \eqref{eqn:HpowerSeries}.

Cieliebak and Latschev \cite{CieliebakLatschev:propaganda}*{\S 2} write down the same
formula in terms of moduli spaces that have no asymptotic markers but
remember the order of the punctures, thus it includes some factorials
that do not appear in \eqref{eqn:Hsigma} but is missing the $\kappa_\gamma$
terms of \eqref{eqn:HpowerSeries}.  The notation $n_g(\Gamma^-,\Gamma^+)$
used in \cite{CieliebakLatschev:propaganda} for curve counts must be understood 
implicitly to include rational weights arising from automorphisms (or
multivalued perturbations, as the case may be).

My paper with Latschev \cite{LatschevWendl} uses moduli spaces with 
asymptotic markers and
attempts to write down the same formula as in \cites{SFT,CieliebakLatschev:propaganda},
but gets it slightly wrong due to some missing $\kappa_\gamma$ terms
that should appear in front of each~$\frac{\p}{\p q_\gamma}$.
Mea culpa.

For cylindrical contact homology, the combinatorial factors in
\S\ref{sec:cyl} also agree with what appears in \cite{Bourgeois:contactHom}.
As observed by Nelson \cite{Nelson:thesis}*{Remark~8.3}, there are other
conventions for $\p_\CCH$ that appear in the literature and lead to
equivalent theories: in particular it is possible to replace \eqref{eqn:CCH}
with
$$
\p_\CCH q_\gamma := \sum_{\gamma',A} \kappa_{\gamma'} \left(\sum_{u \in \mM_A(\gamma,\gamma')}
\frac{\epsilon(u)}{\kappa_u} \right) e^A q_{\gamma'} .
$$
One can derive this from the same definition of $\mathbf{H}$ by applying a
``change of coordinates'' to the algebra $\aA[[\hbar]]$, or equivalently,
by choosing a slightly different representation of the operator algebra
defined by the $p_\gamma$ and~$q_\gamma$ variables.  To avoid confusion,
let us write the generators of $\aA$ as $x_\gamma$ instead of $q_\gamma$,
and then define the operators $q_\gamma$ and $p_\gamma$ on $\aA[[\hbar]]$ by
$$
q_\gamma = \kappa_\gamma x_\gamma , \qquad
p_\gamma = \hbar \frac{\p}{\p x_\gamma}.
$$
These operators still satisfy $[p_\gamma,q_\gamma] = \kappa_\gamma \hbar$
and thus define an equivalent theory, but the resulting differential
operator $\mathbf{D}_\SFT$ on $\aA[[\hbar]]$ now includes extra
factors of $\kappa_\gamma$ for the negative punctures instead of the positive
punctures.

\subsection{Coefficients: $\QQ$, $\ZZ$ or $\ZZ_2$?}

While we were able to use $\ZZ_2$ coefficients for cylindrical contact
homology in a primitive homotopy class in Lecture~\ref{lec:tight3tori}, a quick glance at
any version of the formula for $\mathbf{H}$ should make the reader very
skeptical about doing this for more general versions of SFT.
The existence of curves with automorphisms means that $\mathbf{H}$ always
contains terms with rational (but nonintegral) coefficients.  And this
is only what is true in the fictional world of Assumption~\ref{ass:optimism}:
in the general version of the theory, we expect to have to replace
expressions like $\sum_{u} \frac{\epsilon(u)}{|\Aut(u)|}$ with counts of
$0$-dimensional weighted branched orbifolds with rational weights, arising
as zero-sets of generic multisections.  In this case we not only
obtain rational counts but may also lose all control over the sizes of the
denominators.  

A similar phenomenon occurs in general versions of
Gromov-Witten theory.  For instance, in the approach of Cieliebak-Mohnke 
\cite{CieliebakMohnke:transversality} for the rational Gromov-Witten 
invariants of a closed symplectic manifold $(W^{2n},\omega)$ with
$[\omega] \in H^2(W;\QQ)$, the invariants are defined
by replacing the usual moduli space $\mM_{0,m}(J,A)$ by a space
$\mM_{0,m+N}(J,A;Y)$ consisting of $J$-holomorphic spheres
$u : S^2 \to W$ with some large number of auxiliary marked points
$\zeta_1,\ldots,\zeta_N$ required to satisfy the condition
$$
u(\zeta_i) \in Y, \qquad i=1,\ldots,N.
$$
Here $Y^{2n-2} \subset W^{2n}$ is a $J$-holomorphic hypersurface with
$[Y] = D \cdot \PD([\omega]) \in H_{2n-2}(W)$ for some degree $D \in \NN$, and the number
of extra marked points is determined by
$$
N = A \cdot [Y] = D \langle [\omega],A \rangle,
$$
so positivity of intersections implies that $u$ \emph{only} intersects $Y$
at the auxiliary marked points.  These auxiliary points are convenient
for technical reasons involving trans\-ver\-sality---their role is vaguely
analogous to the way that asymptotic markers get rid of isotropy in SFT---but
they are not geometrically meaningful, as we'd actually prefer to count
curves in $\mM_{0,m}(J,A)$.  Every such curve has $N$ intersections with $Y$, 
so accounting for permutations, it lifts to $N!$ distinct elements of 
$\mM_{0,m+N}(J,A;Y)$, and the correct count is therefore obtained as an
integer count of curves in the latter space divided by~$N!$.
Perturbing to achieve transversality breaks the symmetry, however, so 
there is no guarantee that counting curves in $\mM_{0,m+N}(J,A;Y)$ will
produce a multiple of~$N!$, and moreover, $N$ could be arbitrarily large
since one needs to take hypersurfaces of arbitrarily large degree in
order to show that the invariants don't depend on this choice.
For these reasons, the resulting Gromov-Witten invariants are rational
numbers rather than integers, and their denominators cannot be predicted
or bounded.

The upshot of this discussion is that there is probably no hope of defining
SFT with integer coefficients in general, much less with $\ZZ_2$
coefficients---for this reason the inclusion of orientations in the 
picture is unavoidable.  That is the bad news.

The good news however is that whenever formulas like
$\sum_{u} \frac{\epsilon(u)}{|\Aut(u)|}$ can be taken literally as a count
of curves, the chain complex $(\aA[[\hbar]],\mathbf{D}_\SFT)$ can in fact
be defined with $\ZZ$ coefficients, and one can even reduce to a $\ZZ_2$
version in order to ignore signs.  A special case of this was
observed for cylindrical contact homology in 
\cite{Nelson:Abendblatt}*{Remark~1.5}, and you may notice it already when
you look at the formula \eqref{eqn:CCH} for $\p_\CCH$: the factor
$\kappa_\gamma / \kappa_u$ is always an integer since the multiplicity
of a holomorphic cylinder always divides the covering
multiplicity of both its asymptotic orbits.  Surprisingly, something similar
turns out to be true for the much larger chain complex of SFT.  The following
result is stated under Assumption~\ref{ass:optimism} for safety's sake,
but in light of the discussion in \S\ref{sec:orbifold}, we should expect it
to hold somewhat more generally.

\begin{prop}
\label{prop:integer}
If Assumption~\ref{ass:optimism} holds then the rational coefficients 
$n_g(\boldsymbol{\gamma},\boldsymbol{\gamma}',k)$
in the formula \eqref{eqn:DSFT} for $\mathbf{D}_\SFT q^{\boldsymbol{\gamma}}$
are all integers.
\end{prop}
\begin{cor}
Under Assumption~\ref{ass:optimism}, there exist well-defined chain complexes 
$$
(\aA_\ZZ[[\hbar]],\mathbf{D}_\SFT) \quad \text{ and } \quad
(\aA_{\ZZ_2}[[\hbar]],\mathbf{D}_\SFT),
$$
where for a general commutative ring~$\rR$, $\aA_\rR$ denotes the graded
supercommutative unital algebra over $\rR[H_2(M)/G]$ generated by the
$q_\gamma$ variables for good Reeb orbits~$\gamma$.  The differentials
$\mathbf{D}_\SFT$ on $\aA_\ZZ[[\hbar]]$ and $\aA_{\ZZ_2}[[\hbar]]$ are 
defined by the same formula as on $\aA[[\hbar]]$, where in the $\ZZ_2$ case
we are free to set all signs~$\epsilon(u)$ equal to~$1$.
\end{cor}
\begin{proof}[Proof of Proposition~\ref{prop:integer}]
We need to show that expressions of the form
$$
\frac{\kappa_{\boldsymbol{\gamma}^+}}{|\Aut^\sigma(u)|}  
\frac{\p}{\p q_{\gamma_1^+}} \ldots \frac{\p}{\p q_{\gamma_{k_+}^+}}
q^{\boldsymbol{\gamma}}
$$
produce integer coefficients for every holomorphic curve
$u$ with asymptotic orbits $\boldsymbol{\gamma}^\pm =
(\gamma_1^\pm,\ldots,\gamma_{k_\pm}^\pm)$ and every
tuple $\boldsymbol{\gamma} = (\gamma_1,\ldots,\gamma_m)$.  It suffices to
consider the special case $\boldsymbol{\gamma} = \boldsymbol{\gamma}^+$,
as the derivative in question is only nonzero on monomials that are
divisible by~$q^{\boldsymbol{\gamma}^+}$.  Up to a sign change, we can
reorder the orbits and write $\boldsymbol{\gamma}^+$ in the form
$$
\boldsymbol{\gamma}^+ = ( \underbrace{\gamma_1,\ldots,\gamma_1}_{m_1},
\ldots,\underbrace{\gamma_N,\ldots,\gamma_N}_{m_N})
$$
for some finite set of distinct orbits $\gamma_1,\ldots,\gamma_N$ and numbers
$m_i \in \NN$, $i=1,\ldots,N$.  We then have
\begin{equation}
\label{eqn:fraction}
\begin{split}
\frac{\kappa_{\boldsymbol{\gamma}^+}}{|\Aut^\sigma(u)|}
\frac{\p}{\p q_{\gamma_1^+}} \ldots \frac{\p}{\p q_{\gamma_{k_+}^+}}
q^{\boldsymbol{\gamma}^+} &=
\frac{\kappa_{\gamma_1}^{m_1} \ldots \kappa_{\gamma_N}^{m_N}}{|\Aut^\sigma(u)|}
\left(\frac{\p}{\p q_{\gamma_1}}\right)^{m_1} \ldots 
\left(\frac{\p}{\p q_{\gamma_N}}\right)^{m_N} 
\left(q_{\gamma_1}^{m_1} \ldots q_{\gamma_N}^{m_N} \right) \\
&= \pm \frac{\kappa_{\gamma_1}^{m_1} \ldots \kappa_{\gamma_N}^{m_N} 
m_1! \ldots m_N!}{|\Aut^\sigma(u)|}.
\end{split}
\end{equation}
We claim that this number is always an integer.  Indeed, if $\Aut^\sigma(u)$
is nontrivial, then $u : \dot{\Sigma} \to \RR \times M$ is a multiple cover
$u = v \circ \varphi$ for some holomorphic branched cover 
$\varphi : (\Sigma,j) \to (\Sigma',j')$ and somewhere injective curve
$v : (\dot{\Sigma}' = \Sigma' \setminus \Gamma',j') \to (\RR \times M,J)$.
Automorphisms $\psi \in \Aut^\sigma(u)$ thus define biholomorphic maps on
$(\Sigma,j)$ that permute each of the sets of punctures asymptotic
to the same orbit.  Given any puncture $z \in \Gamma$
where $u$ is asymptotic to~$\gamma_i$, the $\Aut^\sigma(u)$-orbit of $z$ consists of
$\ell \le m_i$ other punctures also asymptotic to~$\gamma_i$, and its
stabilizer is a cyclic subgroup of order $k = |\Aut^\sigma(u)| / \ell$, 
acting on a neighborhood of $z$ by biholomorphic rotations.  It follows that
$\kappa_{\gamma_i}$ is divisible by~$k$, hence
$$
\frac{\kappa_{\gamma_i} \ell}{|\Aut^\sigma(u)|} \in \NN,
$$
and \eqref{eqn:fraction} is a multiple of this.
\end{proof}

\begin{remark}
Since $1 = -1$ in $\aA_{\ZZ_2}$, anticommuting elements of
$\aA_{\ZZ_2}[[\hbar]]$ actually commute, so unless one imposes extra
algebraic conditions in the case of $\ZZ_2$ coefficients,
higher powers of odd generators $p_\gamma$ and $q_\gamma$ do not vanish.  
Nonetheless,  these powers still do not appear in $\mathbf{H}$, so the complex
$(\aA_{\ZZ_2}[[\hbar]],\mathbf{D}_\SFT)$ ignores curves with multiple ends
approaching an orbit of odd degree (and also bad orbits, for that matter).
\end{remark}

\chapter{Contact invariants}
\label{lec:SFT}

\minitoc

\vspace{12pt}

In the previous lecture, we introduced an operator algebra defined via the
supercommutators $[p_\gamma,q_\gamma] = \kappa_\gamma \hbar$, then we defined
the SFT generating function
$$
\mathbf{H} = \sum_{u \in \mM^\sigma_1(J) / \RR} \frac{\epsilon(u)}{|\Aut^\sigma(u)|}
\hbar^{g-1} e^A q^{\boldsymbol{\gamma}^-} p^{\boldsymbol{\gamma}^+}
$$
and proved (modulo transversality) that $\mathbf{H}^2 = 0$.  The generating
function is a formal power series whose coefficients are rational counts of 
holomorphic
curves, and these counts are strongly dependent on the choices of contact 
form $\alpha$, almost complex structure~$J \in \jJ(\alpha)$ and further
auxiliary data such as coherent orientations.  Thus in contrast to
Gromov-Witten theory, the generating function does not define an invariant,
but one can follow the standard prescription of Floer-type theories and define
invariants via homology.  We saw that for the natural representation 
$\aA[[\hbar]]$ of the operator algebra defined by setting
$p_\gamma = \kappa_\gamma \hbar \frac{\p}{\p q_\gamma}$, $\mathbf{H}$ defines a
differential operator $\mathbf{D}_\SFT : \aA[[\hbar]] \to \aA[[\hbar]]$
with $\mathbf{D}_\SFT^2 = 0$.
One of our goals in this lecture will be to explain (again modulo transversality)
why the resulting homology
$$
H_*^\SFT(M,\xi;R) = H_*(\aA[[\hbar]],\mathbf{D}_\SFT)
$$
is an invariant of the contact structure.  We will then use it to define
simpler numerical invariants that detect symplectic fillability properties
of contact manifolds.

But first, $\aA[[\hbar]]$ is not the only possible representation of the 
operator algebra of SFT: other choices lead to different invariants
with different algebraic structures.  Let's begin by describing the
original hierarchy of contact invariants that were outlined in~\cite{SFT}.

\begin{remark}
\label{remark:torsion}
Throughout this lecture, we assume for simplicity that $H_2(M)$ has no torsion,
and the same assumption is made about cobordisms in \S\ref{sec:disconnected}.
Only minor changes are necessary if this condition is lifted, e.g.~one could 
then replace all instances of $H_2(M)$ with $H_2(M;\QQ)$ and assume always
that the grading is~$\ZZ_2$; see \S\ref{sec:torsion}.
\end{remark}

\section{The Eliashberg-Givental-Hofer package}

In the following, $(M,\xi)$ is a $(2n-1)$-dimensional closed contact
manifold with a contact form $\alpha$ and almost complex structure
$J \in \jJ(\alpha)$ for which the optimistic transversality condition
(Assumption~\ref{ass:optimism}) of Lecture~\ref{lec:H} is assumed to hold.  We fix also the
auxiliary data described in \S\ref{sec:setup3}, plus a choice of subgroup
$G \subset H_2(M)$ which determines the coefficient ring
$$
R = \QQ[H_2(M) / G].
$$
Each of the differential graded algebras described below then carries the same
grading that was described in that lecture, i.e.~there is always at least
a $\ZZ_2$-grading, and it lifts to $\ZZ$ if $H_1(M)$ is torsion free and
$c_1(\xi)|_G = 0$, or possibly $\ZZ_{2N}$ if $N \in \NN$ is the
smallest possible value for $c_1(A)$ with $A \in G$.

\subsection{Full SFT as a Weyl superalgebra}

We start with some seemingly trivial algebraic observations.  First, the 
relation $\mathbf{H}^2 = 0$ is equivalent to
$$
[\mathbf{H},\mathbf{H}] = 0.
$$
Remember that $[\ ,\ ]$ is a \emph{super}-commutator, so 
$[\mathbf{F},\mathbf{F}] = 0$ holds automatically for operators~$\mathbf{F}$ 
with even degree, but $\mathbf{H}$ is odd, and for odd operators the commutator 
is defined by $[\mathbf{F},\mathbf{G}] = \mathbf{F}\mathbf{G} + 
\mathbf{G}\mathbf{F}$, hence $[\mathbf{H},\mathbf{H}] = 2 \mathbf{H}^2$.  
Formally speaking $[\ ,\ ]$ is a \defin{super Lie bracket} and thus 
satisfies the ``super Jacobi identity'':
\begin{equation}
\label{eqn:Jacobi}
\big[\mathbf{F},[\mathbf{G},\mathbf{K}]\big] + (-1)^{|\mathbf{F}| |\mathbf{G}| +
|\mathbf{F}| |\mathbf{K}|} \big[\mathbf{G},[\mathbf{K},\mathbf{F}]\big] +
(-1)^{|\mathbf{F}| |\mathbf{K}| + |\mathbf{G}| |\mathbf{K}|} 
\big[\mathbf{K},[\mathbf{F},\mathbf{G}]\big] = 0.
\end{equation}
A consequence of this is that in order to create a homology theory out of
$\mathbf{H}$, we don't absolutely need to find a representation of the
entire operator algebra: it suffices to find a representation of the induced
super Lie algebra.  Indeed, suppose $V$ is a graded $R[[\hbar]]$-module and
$L$ is a linear grading-preserving map that associates to operators $\mathbf{F}$
(expressed as power series functions of $p$'s, $q$'s and $\hbar$ with
coefficients in~$R$) an $R[[\hbar]]$-linear map
$$
L_{\mathbf{F}} : V \to V
$$
such that
$$
L_{[\mathbf{F},\mathbf{G}]} = 
L_{\mathbf{F}} L_{\mathbf{G}} - (-1)^{|\mathbf{F}| |\mathbf{G}|}
L_{\mathbf{G}} L_{\mathbf{F}}
$$
for every pair of operators $\mathbf{F},\mathbf{G}$.  Then the 
$R[[\hbar]]$-linear map $L_{\mathbf{H}} : V \to V$ satisfies
$$
L_{\mathbf{H}}^2 = \frac{1}{2} [L_{\mathbf{H}},L_{\mathbf{H}}] = 
\frac{1}{2} L_{[\mathbf{H},\mathbf{H}]} = 0,
$$
hence $(V,L_{\mathbf{H}})$ is a chain complex.  The complex
$(\aA[[\hbar]],\mathbf{D}_\SFT)$ was a special case of this, in which we
represented the super Lie algebra via a faithful representation of the whole
operator algebra.

\begin{exercise}
Verify \eqref{eqn:Jacobi}.
\end{exercise}

\begin{remark}[supersymmetric sign rules]
To see where the signs in \eqref{eqn:Jacobi} come from, it suffices
to know the following basic rule of superalgebra: for any pair of
$\ZZ_2$-graded vector spaces $V$ and $W$, the natural ``commutation''
isomorphism $c : V \otimes W \to W \otimes V$ is defined on homogeneous 
elements by
$$
c(v \otimes w) = (-1)^{|v| |w|} w \otimes v.
$$
For any permutation of a finite tuple of $\ZZ_2$-graded vector spaces,
one can derive the appropriate isomorphism from this: in particular the
cyclic permutation isomorphism 
$\sigma : X \otimes Y \otimes Z \to Y \otimes Z \otimes X$ takes the form
$$
\sigma = (\1 \otimes c_{23}) \circ (c_{12} \otimes \1) :
x \otimes y \otimes z \mapsto (-1)^{|x| |y| + |x| |z|} y \otimes z \otimes x.
$$
Writing the Jacobi identity as $[\cdot,[\cdot,\cdot ]] \circ (\1 + \sigma + \sigma^2) = 0$
then produces \eqref{eqn:Jacobi}.  In this sense, it only differs from
the usual Jacobi identity in being based on a different definition of the
commutation isomorphism $V \otimes W \to W \otimes V$.  For more on this
perspective, see \cite{Varadarajan:susy}*{\S 3.1}.
\end{remark}

Now here is a different kind of example, where the representation does not respect
the product structure of the operator algebra but does respect its Lie bracket.
Let $\Weyl$ denote the graded unital algebra consisting of formal power series
$$
\sum_{\boldsymbol{\gamma},k} f_{\boldsymbol{\gamma},k}(q) \hbar^k p^{\boldsymbol{\gamma}},
$$
where the sum ranges over all integers $k \ge 0$ and all ordered sets 
$\boldsymbol{\gamma} = (\gamma_1,\ldots,\gamma_m)$ of good Reeb orbits for
$m \ge 0$, and the $f_{\boldsymbol{\gamma},k}$ are polynomial functions
of the $q_\gamma$ variables, with coefficients in~$R$.  Note that the
case of the empty set of orbits is included here, which means
$p^{\boldsymbol{\gamma}} = 1$.  The multiplicative structure of $\Weyl$ is
defined via the usual (super)commutation relations, and its elements can be
interpreted as operators.  If we now associate
to each $\mathbf{F} \in \Weyl$ the $R[[\hbar]]$-linear map
$$
D_{\mathbf{F}} : \Weyl \to \Weyl : \mathbf{G} \mapsto [\mathbf{F},\mathbf{G}],
$$
then the Jacobi identity \eqref{eqn:Jacobi} implies
$$
D_{[\mathbf{F},\mathbf{G}]} = D_{\mathbf{F}} D_{\mathbf{G}} -
(-1)^{|\mathbf{F}| |\mathbf{G}|} D_{\mathbf{G}} D_{\mathbf{F}}.
$$
This is just the graded version of the standard \emph{adjoint representation} of
a Lie algebra.  The only problem in applying this idea to define a differential
\begin{equation}
\label{eqn:HWeyl}
D_{\mathbf{H}} : \Weyl \to \Weyl : \mathbf{F} \mapsto 
[\mathbf{H},\mathbf{F}]
\end{equation}
is that $\mathbf{H}$ is not technically an element of~$\Weyl$: indeed,
$\mathbf{H}$ contains terms of order $-1$ in~$\hbar$, thus
$$
\mathbf{H} \in \frac{1}{\hbar} \Weyl.
$$
On the other hand, the failure of supercommutativity in $\Weyl$ is a
``phenomenon of order~$\hbar$,'' i.e.~since every nontrivial commutator
contains a factor of~$\hbar$, we have
$$
[\mathbf{F},\mathbf{G}] = \Order(\hbar) \quad \text{ for all } \quad
\mathbf{F},\mathbf{G} \in \Weyl.
$$
Here and in the following we use the symbol
$$
\Order(\hbar^k)
$$
to denote any element of the form $\hbar^k \mathbf{F}$ for
$\mathbf{F} \in \Weyl$.  As a consequence, $[\mathbf{H},\mathbf{F}]
\in \Weyl$ whenever $\mathbf{F} \in \Weyl$, hence \eqref{eqn:HWeyl} is
well defined, and the Jacobi identity now implies
$$
D_{\mathbf{H}}^2 = 0.
$$
The homology of the resulting chain complex gives another version of what
is often called \defin{full SFT},
$$
H_*^{\Weyl}(M,\xi ; R) := H_*(\Weyl,D_{\mathbf{H}}).
$$
A proof (modulo transversality) that this defines a contact invariant is 
outlined in \cite{SFT}*{\S 2}, but it is algebraically somewhat more involved
than for $H_*^\SFT(M,\xi;R)$, so I will 
skip it since I don't have any applications of 
$H_*^{\Weyl}(M,\xi;R)$ in mind.  As far as I am aware, no contact topological
applications of this invariant or computations of it (outside the trivial
case---see \S\ref{sec:OT} below) have yet appeared in the literature.
This is a pity, because $H_*^{\Weyl}(M,\xi;R)$ actually has much more
algebraic structure than $H_*^\SFT(M,\xi;R)$.  Indeed, using the identities
\begin{equation}
\label{eqn:Leibniz2}
\begin{split}
[\mathbf{F},\mathbf{G}\mathbf{K}] &= [\mathbf{F},\mathbf{G}] \mathbf{K}
+ (-1)^{|\mathbf{F}| |\mathbf{G}|} \mathbf{G} [\mathbf{F},\mathbf{K}], \\
[\mathbf{F}\mathbf{G},\mathbf{K}] &= \mathbf{F} [\mathbf{G},\mathbf{K}]
+ (-1)^{|\mathbf{G}| |\mathbf{K}|} [\mathbf{F},\mathbf{K}] \mathbf{G},
\end{split}
\end{equation}
one sees that $D_{\mathbf{H}} : \Weyl \to \Weyl$ satisfies a graded
Leibniz rule,
$$
D_{\mathbf{H}} (\mathbf{F}\mathbf{G}) = (D_{\mathbf{H}} \mathbf{F}) \mathbf{G}
+ (-1)^{|\mathbf{F}|} \mathbf{F} \, D_{\mathbf{H}} \mathbf{G}.
$$
It follows that $D_{\mathbf{H}} : \Weyl \to \Weyl$ is also a derivation with
respect to the bracket structure on~$\Weyl$, i.e.
$$
D_{\mathbf{H}} [\mathbf{F},\mathbf{G}] = [D_{\mathbf{H}} \mathbf{F} , \mathbf{G}]
+ (-1)^{|\mathbf{F}|} [\mathbf{F}, D_{\mathbf{H}} \mathbf{G}]
$$
for all $\mathbf{F}, \mathbf{G} \in \Weyl$.  As a consequence, the product
and bracket structures on 
$\Weyl$ descend to $H_*^{\Weyl}(M,\xi;R)$, giving it the structure of a
\emph{Weyl superalgebra}.

As a matter of interest, we observe that $(\Weyl,D_{\mathbf{H}})$, as with
$(\aA[[\hbar]],\mathbf{D}_\SFT)$ in the previous lecture,
can be defined with $\ZZ$ or $\ZZ_2$ coefficients whenever the 
transversality results are good enough to take the usual expression 
$\sum_u \frac{\epsilon(u)}{|\Aut^\sigma(u)|}$ literally as a
count of holomorphic curves.  This result is of limited interest since
it cannot hold in general cases where transversality for multiple covers
is impossible without multivalued perturbations---nonetheless I find it
amusing.\footnote{The same arguments used to define SFT chain complexes
over the integers can also be applied to the chain maps involved in the
proof of invariance (see~\S\ref{sec:invariance}), so the SFT invariants
\emph{should} be defined over the integers if transversality can be achieved
for multiple covers.  There are known situations however in which this
cannot hold: even if the chain complexes are well defined over~$\ZZ$,
invariance may hold only over~$\QQ$, due to the failure of transversality
in cobordisms.  See \cite{Hutchings:integer}.}

\begin{prop}
\label{prop:integers}
If Assumption~\ref{ass:optimism} in Lecture~\ref{lec:H} holds, then $D_{\mathbf{H}}$ is also well
defined if the ring $R = \QQ[H_2(M) / G]$ is replaced by
$\ZZ[H_2(M) / G]$ or $\ZZ_2[H_2(M) / G]$.
\end{prop}
\begin{proof}
Since $D_{\mathbf{H}}$ is a derivation, it suffices to check that for every
good Reeb orbit~$\gamma$, $D_{\mathbf{H}} q_\gamma$ and 
$D_{\mathbf{H}} p_\gamma$ are each sums of monomials of the form
$c e^A \hbar^k q^{\boldsymbol{\gamma}^-} p^{\boldsymbol{\gamma}^+}$ with
coefficients $c \in \ZZ$.  Suppose $u \in \mM_1(J)$ is an index~$1$
holomorphic curve with positive and/or negative asymptotic orbits
$$
\boldsymbol{\gamma}^\pm = 
(\underbrace{\gamma_1^\pm,\ldots,\gamma_1^\pm}_{m_1^\pm},\ldots,
\underbrace{\gamma_{k_\pm}^\pm,\ldots,\gamma_{k_\pm}^\pm}_{m_{k_\pm}^\pm}),
$$
where $\gamma_i^\pm \ne \gamma_j^\pm$ for $i \ne j$.  We can assume
all the orbits $\gamma_i^\pm$ are good and that $m_i^\pm = 1$ whenever
$n - 3 + \muCZ(\gamma_i^\pm)$ is odd.  Up to a sign and
factors of $e^A$ and $\hbar$ which are not relevant to this discussion,
$u$ then contributes a monomial
$$
\mathbf{H}_u := \frac{1}{|\Aut^\sigma(u)|} q_{\gamma_1^-}^{m_1^-} \ldots
q_{\gamma_{k_-}^-}^{m_{k_-}^-}  p_{\gamma_1^+}^{m_1^+} \ldots
p_{\gamma_{k_+}^+}^{m_{k_+}^+}
$$
to~$\mathbf{H}$.  The commutator $[\mathbf{H}_u,q_\gamma]$ vanishes unless
$\gamma$ is one of the orbits $\gamma_1^+,\ldots,\gamma_{k_+}^+$, so
suppose $\gamma = \gamma_{k_+}^+$.  If $n - 3 + \muCZ(\gamma)$ is odd, then
$m := m_{k_+}^+ = 1$, and \eqref{eqn:Leibniz2} with $[p_\gamma,q_\gamma] =
\kappa_\gamma \hbar$ implies
\begin{equation*}
\begin{split}
[\mathbf{H}_u,q_\gamma] &= \frac{1}{|\Aut^\sigma(u)|} 
\left[ q_{\gamma_1^-}^{m_1^-} \ldots q_{\gamma_{k_-}^-}^{m_{k_-}^-}  
p_{\gamma_1^+}^{m_1^+} \ldots p_{\gamma_{k_+ - 1}^+}^{m_{k_+ - 1}^+} p_\gamma,q_\gamma\right] \\
&= \frac{\kappa_\gamma}{|\Aut^\sigma(u)|} \hbar 
q_{\gamma_1^-}^{m_1^-} \ldots q_{\gamma_{k_-}^-}^{m_{k_-}^-}  
p_{\gamma_1^+}^{m_1^+} \ldots p_{\gamma_{k_+ - 1}^+}^{m_{k_+ - 1}^+}.
\end{split}
\end{equation*}
The fraction in front of this expression is an integer since $u$ can have only
one end asymptotic to~$\gamma$, and $\kappa_\gamma$ is thus divisible by the
covering multiplicity of~$u$.  If $n - 3 + \muCZ(\gamma)$ is even, then we
generalize this calculation by using \eqref{eqn:Leibniz2} to write
$$
[p_\gamma^m,q_\gamma] = m \kappa_\gamma \hbar p_\gamma^{m-1},
$$
so then,
\begin{equation*}
\begin{split}
[\mathbf{H}_u,q_\gamma] &= \frac{1}{|\Aut^\sigma(u)|} 
\left[ q_{\gamma_1^-}^{m_1^-} \ldots q_{\gamma_{k_-}^-}^{m_{k_-}^-}  
p_{\gamma_1^+}^{m_1^+} \ldots p_{\gamma_{k_+ - 1}^+}^{m_{k_+ - 1}^+} p_\gamma^m,q_\gamma\right] \\
&= \frac{\kappa_\gamma m}{|\Aut^\sigma(u)|} \hbar 
q_{\gamma_1^-}^{m_1^-} \ldots q_{\gamma_{k_-}^-}^{m_{k_-}^-}  
p_{\gamma_1^+}^{m_1^+} \ldots p_{\gamma_{k_+ - 1}^+}^{m_{k_+ - 1}^+} p_\gamma^{m-1}.
\end{split}
\end{equation*}
To see that $\frac{\kappa_\gamma m}{|\Aut^\sigma(u)|}$ is always an integer,
recall from our proof of Prop.~\ref{prop:integer} in the previous lecture that 
transformations in $\Aut^\sigma(u)$ permute each of the sets of punctures
that are asymptotic to the same Reeb orbit.  Suppose the set of positive punctures 
of $u$ asymptotic to~$\gamma$ is partitioned by the $\Aut^\sigma(u)$-action
into $N$ subsets, each consisting of
$\ell_1,\ldots,\ell_N$ punctures, where $\ell_1 + \ldots + \ell_N = m$.  
If $z$ is a puncture in the $i$th of these subsets,
then its stabilizer is a cyclic subgroup
of order $k_i$ acting on a neighborhood of~$z$ by biholomorphic rotations,
where $k_i \ell_i = |\Aut^\sigma(u)|$.  Each of these orders $k_i$ necessarily
divides the multiplicity $\kappa_\gamma$, so we can write 
$k_i a_i = \kappa_\gamma$ for some $a_i \in \NN$.  Putting all this together,
we have
$$
\kappa_\gamma m = \sum_{i=1}^N \kappa_\gamma \ell_i = \sum_{i=1}^N k_i a_i \ell_i
= |\Aut^\sigma(u)| \sum_{i=1}^N a_i.
$$

Following this same procedure, you should now be able to verify on your own
that the coefficient appearing in $[\mathbf{H}_u,p_\gamma]$ is also always
an integer.  The existence of a chain complex with $\ZZ_2$ coefficients 
follows from this simply by projecting $\ZZ$ to~$\ZZ_2$.
\end{proof}

\subsection{The semiclassical limit: rational SFT}

The idea of rational symplectic field theory (RSFT) is to extract as much information
as possible from genus zero holomorphic curves but ignore curves of higher
genus.  The algebra of SFT provides a fairly obvious mechanism for this:
RSFT should be what SFT becomes in the ``limit as $\hbar \to 0$,'' 
i.e.~the classical approximation to a quantum theory.  Let
$$
\Poisson := \Weyl \big/ \hbar \Weyl,
$$
so $\Poisson$ is a graded unital algebra generated by the $p_\gamma$ and 
$q_\gamma$ variables and the coefficient ring~$R$, but it does not include
$\hbar$ as a generator.  
Since all commutators in $\Weyl$ are in $\hbar \Weyl$,
the product structure of $\Poisson$ is supercommutative.  Let us use the
distinction between capital and lowercase letters to denote
the quotient projection
$$
\Weyl \to \Poisson : \mathbf{F} \mapsto \mathbf{f}.
$$
We will make an exception for the letter ``H'': recall that $\mathbf{H}$
is not an element of $\Weyl$ since its genus zero terms have order $-1$
in~$\hbar$, but $\hbar \mathbf{H} \in \Weyl$, so we will define
$$
\mathbf{h} = \sum_u \frac{\epsilon(u)}{|\Aut^\sigma(u)|} e^A
q^{\boldsymbol{\gamma}^-} p^{\boldsymbol{\gamma}^+} \in \Poisson
$$
to be the image of $\hbar \mathbf{H}$ under the projection.  The sum in
this expression ranges over all $\RR$-equivalence classes of index~$1$
curves with genus zero, so $\mathbf{h}$ will
serve as the generating function of RSFT.  To encode gluing of genus
zero terms, note first that the commutator operation would not be appropriate
since it prodcues terms for \emph{every} possible gluing of two curves,
including those which glue genus zero curves along more than one breaking
orbit to produce buildings with positive arithmetic genus.  We need instead
to have an algebraic operation on $\Poisson$ that encodes gluing along
only one breaking orbit at a time.

You already know what to expect if you've ever taken a quantum mechanics course:
in the ``classical limit,'' commutators become Poisson brackets.
To express this properly, we need to make a distinction between differential
operators operating from the left or the right: let
$$
\overrightarrow{\frac{\p}{\p q_\gamma}} : \Weyl \to \Weyl
$$
denote the usual operator $\frac{\p}{\p q_\gamma}$, which was previously
defined on $\aA[[\hbar]]$ but has an obvious extension to $\Weyl$ 
such that $\overrightarrow{\frac{\p}{\p q_\gamma}} p_{\gamma'} = 0$ for
all~$\gamma'$.  This operator satisfies the graded Leibniz rule
$$
\overrightarrow{\frac{\p}{\p q_\gamma}} (\mathbf{F} \mathbf{G})
= \left( \overrightarrow{\frac{\p}{\p q_\gamma}} \mathbf{F} \right) \mathbf{G}
+ (-1)^{|q_\gamma| |\mathbf{F}|} \mathbf{F} 
\left( \overrightarrow{\frac{\p}{\p q_\gamma}} \mathbf{G} \right).
$$
The related operator
$$
\overleftarrow{\frac{\p}{\p q_\gamma}} : \Weyl \to \Weyl :
\mathbf{F} \mapsto \mathbf{F} \overleftarrow{\frac{\p}{\p q_\gamma}}
$$
is defined exactly the same way on individual variables $p_\gamma$
and $q_\gamma$, but satisfies a slightly different Leibniz rule,
$$
(\mathbf{F} \mathbf{G}) \overleftarrow{\frac{\p}{\p q_\gamma}} =
\mathbf{F} \left( \mathbf{G} \overleftarrow{\frac{\p}{\p q_\gamma}} \right)
+ (-1)^{|q_\gamma| |\mathbf{G}|}
\left( \mathbf{F} \overleftarrow{\frac{\p}{\p q_\gamma}} \right) \mathbf{G}.
$$
The point of writing $\overleftarrow{\frac{\p}{\p q_\gamma}}$ so that it
acts from the right is to obey the usual conventions of superalgebra:
signs change whenever the order of two odd elements (or operators) is
interchanged.  Partial derivatives with respect to $p_\gamma$ can be defined
analogously on~$\Weyl$.  With this notation in hand, 
the \defin{graded Poisson bracket} on $\Weyl$ is defined by
\begin{equation}
\label{eqn:Poisson}
\{ \mathbf{F},\mathbf{G} \} = \sum_\gamma \kappa_\gamma \left(
\mathbf{F} \overleftarrow{\frac{\p}{\p p_\gamma}}
\overrightarrow{\frac{\p}{\p q_\gamma}} \mathbf{G} -
(-1)^{|\mathbf{F}| |\mathbf{G}|}
\mathbf{G} \overleftarrow{\frac{\p}{\p p_\gamma}}
\overrightarrow{\frac{\p}{\p q_\gamma}} \mathbf{F} \right),
\end{equation}
where the sum ranges over all good Reeb orbits.  In the same manner,
the differential operators and the bracket 
$\{\ ,\ \}$ can also be defined on~$\Poisson$.

It is easy to check that $\{\ ,\ \}$ on $\Weyl$ \emph{almost} satisfies 
a version of \eqref{eqn:Leibniz2}: we have
\begin{equation}
\label{eqn:LeibnizPoisson}
\begin{split}
\{\mathbf{F},\mathbf{G}\mathbf{K}\} &= \{\mathbf{F},\mathbf{G}\} \mathbf{K}
+ (-1)^{|\mathbf{F}| |\mathbf{G}|} \mathbf{G} \{\mathbf{F},\mathbf{K}\} + \Order(\hbar), \\
\{\mathbf{F}\mathbf{G},\mathbf{K}\} &= \mathbf{F} \{\mathbf{G},\mathbf{K}\}
+ (-1)^{|\mathbf{G}| |\mathbf{K}|} \{\mathbf{F},\mathbf{K}\} \mathbf{G} + \Order(\hbar)
\end{split}
\end{equation}
for all $\mathbf{F},\mathbf{G},\mathbf{K} \in \Weyl$.
The extra terms denoted by $\Order(\hbar)$ arise from the fact that in proving
\eqref{eqn:LeibnizPoisson}, we must sometimes reorder products
$\mathbf{F} \mathbf{G}$ by writing them as 
$(-1)^{|\mathbf{F}| |\mathbf{G}|} \mathbf{G} \mathbf{F} + [\mathbf{F},\mathbf{G}]$,
where $[\mathbf{F},\mathbf{G}] = \Order(\hbar)$.  Since the terms with $\hbar$
disappear in $\Poisson$, the relations become exact in $\Poisson$:
\begin{equation}
\label{eqn:LeibnizPoisson2}
\begin{split}
\{\mathbf{f},\mathbf{g}\mathbf{k}\} &= \{\mathbf{f},\mathbf{g}\} \mathbf{k}
+ (-1)^{|\mathbf{f}| |\mathbf{g}|} \mathbf{g} \{\mathbf{f},\mathbf{k}\}, \\
\{\mathbf{f}\mathbf{g},\mathbf{k}\} &= \mathbf{f} \{\mathbf{g},\mathbf{k}\}
+ (-1)^{|\mathbf{g}| |\mathbf{k}|} \{\mathbf{f},\mathbf{k}\} \mathbf{g}
\end{split}
\end{equation}
for all $\mathbf{f},\mathbf{g},\mathbf{k} \in \Poisson$.

\begin{prop}
\label{prop:classical}
For all $\mathbf{F},\mathbf{G} \in \Weyl$,
$$
[\mathbf{F},\mathbf{G}] = \hbar \{ \mathbf{f}, \mathbf{g} \} + \Order(\hbar^2),
$$
and $\{\ ,\ \}$ satisfies the conditions of a super Lie bracket
on~$\Poisson$.
\end{prop}
\begin{remark}
\label{remark:inclusion}
In formulas like the one in the above proposition, we interpret
$\{ \mathbf{f},\mathbf{g} \} \in \Poisson$ as an element of $\Weyl$ via
any choice of $R$-linear inclusion $\Poisson \hookrightarrow \Weyl$ that acts 
as the identity on the generators $p_\gamma, q_\gamma$.  There is ambiguity
in this choice due to the noncommutativity of $\Weyl$, but the ambiguity
is in $\hbar \Weyl$ and thus makes no difference to the formula.
\end{remark}
\begin{proof}[Proof of Proposition~\ref{prop:classical}]
The formula is easily checked when $\mathbf{F}$ and $\mathbf{G}$ are
individual variables of the form $p_\gamma$ or~$q_\gamma$; in fact the
extra term $\Order(\hbar^2)$ can be omitted in these cases.  The case where
$\mathbf{F}$ and $\mathbf{G}$ are general monomials follows from this 
via \eqref{eqn:Leibniz2} and \eqref{eqn:LeibnizPoisson} using
induction on the number of variables in the product.  This implies the
general case via bilinearity.

Given the formula, the condition $\{\mathbf{f},\mathbf{g}\} +
(-1)^{|\mathbf{f}| |\mathbf{g}|} \{\mathbf{g},\mathbf{f}\} = 0$ and the
Poisson version of the super Jacobi identity \eqref{eqn:Jacobi}
follow from the corresponding properties of~$[\ ,\ ]$.
\end{proof}

The proposition implies that our genus zero generating function 
$\mathbf{h} \in \Poisson$ satisfies
$0 = \hbar^2 [\mathbf{H},\mathbf{H}] = [\hbar \mathbf{H},\hbar \mathbf{H}]
= \hbar \{\mathbf{h},\mathbf{h}\} + \Order(\hbar^2)$, thus
$$
\{ \mathbf{h},\mathbf{h} \} = 0.
$$
This relation can be interpreted as the count of boundary points of all
$1$-dimensional moduli spaces of genus zero curves: indeed, any pair of 
genus two curves $u , v \in \mM_1^\sigma(J) / \RR$ constributes to
$\{\mathbf{h},\mathbf{h}\}$ a term of the form
$$
\sum_\gamma \frac{\kappa_\gamma}{|\Aut^\sigma(u)| |\Aut^\sigma(v)|}
e^{A_u + A_v} q^{\boldsymbol{\gamma}^-_u} \left( p^{\boldsymbol{\gamma}^+_u} 
\overleftarrow{\frac{\p}{\p p_\gamma}} \right)
\left( \overrightarrow{\frac{\p}{\p q_\gamma}} q^{\boldsymbol{\gamma}^-_v} \right)
p^{\boldsymbol{\gamma}^+_v},
$$
plus a corresponding term with the roles of $u$ and $v$ reversed.
This sums all the monomials that one can construct by cancelling one
$p_\gamma$ variable from $u$ with a matching $q_\gamma$ variable from~$v$,
in other words, constructing a building by gluing $v$ on top of $u$ along
one matching Reeb orbit.

The graded Jacobi identity will again imply that any representation of the
super Lie algebra $(\Poisson,\{\ ,\ \})$ gives rise to a chain complex
with $\mathbf{h}$ as its differential.  For example we can take the adjoint
representation,
$$
\Poisson \to \End_R(\Poisson) : \mathbf{f} \mapsto d_{\mathbf{f}}, \qquad
d_{\mathbf{f}} \mathbf{g} := \{ \mathbf{f},\mathbf{g} \},
$$
which satisfies $d_{\{\mathbf{f},\mathbf{g}\}} = d_{\mathbf{f}} d_{\mathbf{g}} -
(-1)^{|\mathbf{f}| |\mathbf{g}|} d_{\mathbf{g}} d_{\mathbf{f}}$ due to the
Jacobi identity.  Then $d_{\mathbf{h}}^2 = 0$ since
$\mathbf{h}$ has odd degree and $\{\mathbf{h},\mathbf{h}\} = 0$,
and the homology of \defin{rational SFT} is defined as
$$
H_*^\RSFT(M,\xi;R) := H_*(\Poisson,d_{\mathbf{h}}).
$$
We again refer to \cite{SFT} for an argument that
$H_*^\RSFT(M,\xi;R)$ is an invariant of the contact structure.
Notice that Proposition~\ref{prop:classical} yields a simple relationship
between the chain complexes $(\Weyl,D_{\mathbf{H}})$ and
$(\Poisson,d_{\mathbf{h}})$, namely
\begin{equation}
\label{eqn:Dd}
D_{\mathbf{H}} \mathbf{F} = d_{\mathbf{h}} \mathbf{f} + \Order(\hbar),
\end{equation}
where $d_{\mathbf{h}} \mathbf{f}$ is interpreted as an element of $\Weyl$ via
Remark~\ref{remark:inclusion}.  In other words, the projection
$\Weyl \to \Poisson : \mathbf{F} \to \mathbf{f}$ is a chain map.
Moreover, $d_{\mathbf{H}}$ is a derivation on $\Poisson$ with respect to
both the product and the Poisson bracket: this follows via
Proposition~\ref{prop:classical} and \eqref{eqn:Dd} from the fact that 
$D_{\mathbf{H}}$ satisfies the corresponding properties on~$\Weyl$.
We conclude that $H_*^\RSFT(M,\xi;R)$ inherits the structure of a Poisson
superalgebra, and the map
$$
H_*^{\Weyl}(M,\xi;R) \to H_*^\RSFT(M,\xi;R)
$$
induced by the chain map $(\Weyl,D_{\mathbf{H}}) \to (\Poisson,d_{\mathbf{h}})$
is both an algebra homomorphism and a homomorphism of graded super Lie algebras.

\subsection{The contact homology algebra}

Contact homology is the most popular tool in the SFT package and was probably
the first to be understood beyond the more straightforward cylindrical
theory.  In situations where cylindrical contact homology cannot be defined
due to bubbling of holomorphic planes, the next simplest thing one can do
is to define a theory that counts genus zero curves with one positive end
but \emph{arbitrary} numbers of negative ends (cf.~Exercise~\ref{EX:onePositive} in
Lecture~\ref{lec:tight3tori}).  

The proper algebraic setting for such a theory turns out to
be the algebra $\aA$ generated by the $q_\gamma$ variables, and it can be
derived from RSFT by setting all $p_\gamma$ variables to zero.
Using the obvious inclusion $\aA \hookrightarrow \Poisson$, define
$\p_\CH : \aA \to \aA$ by
$$
\p_\CH \mathbf{f} = d_{\mathbf{h}} \mathbf{f}|_{p=0}.
$$
We can thus write $d_{\mathbf{h}} \mathbf{f} = \p_\CH \mathbf{f} + \Order(p)$,
where
$$
\Order(p^k)
$$
will be used generally to denote any formal sum consisting exclusively of
terms of the form $p_{\gamma_1} \ldots p_{\gamma_k} \mathbf{f}$ for
$\mathbf{f} \in \Poisson$.  Now observe that for any good orbit~$\gamma$,
$$
d_{\mathbf{h}} p_\gamma = \{\mathbf{h},p_\gamma\} =
- (-1)^{|p_\gamma|} \sum_{\gamma'} \left( p_\gamma 
\overleftarrow{\frac{\p}{\p p_{\gamma'}}} \right)
\left( \overrightarrow{\frac{\p}{\p q_{\gamma'}}} \mathbf{h} \right) =
- (-1)^{|p_\gamma|} \frac{\p \mathbf{h}}{\p q_\gamma} = \Order(p)
$$
since every term in $\mathbf{h}$ has at least one $p$ variable.  It follows
that $d_{\mathbf{h}} \left( \Order(p) \right) = \Order(p)$, so the fact
that $d_{\mathbf{h}}^2 = 0$ implies $\p_\CH^2 = 0$, and \defin{contact homology}
is defined as
$$
HC_*(M,\xi;R) := H_*(\aA,\p_\CH).
$$
Since $d_{\mathbf{h}}$ is a derivation on~$\Poisson$, the formula
$d_{\mathbf{h}} \mathbf{f} = \p_\CH \mathbf{f} + \Order(p)$ implies that
$\p_\CH$ is likewise a derivation on~$\aA$, so $HC_*(M,\xi;R)$ has the
structure of a graded supercommutative algebra with unit.  Moreover, the
projection $\Poisson \to \aA : \mathbf{f} \mapsto \mathbf{f}|_{p=0}$ is a chain
map, giving rise to an algebra homomorphism
$$
H_*^\RSFT(M,\xi;R) \to HC_*(M,\xi;R).
$$
The invariance of $HC_*(M,\xi;R)$ will follow from the invariance of
$H_*^\SFT(M,\xi;R)$, to be discussed in \S\ref{sec:invariance} below.

To interpret $\p_\CH$, we can separate the part of $\mathbf{h}$ that is
linear in $p$ variables, writing
$$
\mathbf{h} = \sum_{\gamma} \mathbf{h}_\gamma(q) p_\gamma + \Order(p^2),
$$
where for each good Reeb orbit~$\gamma$,
$\mathbf{h}_\gamma(q)$ denotes a polynomial in $q$ variables with
coefficients in~$R$.  Since elements $\mathbf{f} \in \aA$ have no dependence
on $p$ variables, we then have
$$
d_{\mathbf{h}} \mathbf{f} = \{\mathbf{h},\mathbf{f}\} = 
\sum_\gamma \kappa_\gamma 
\left( \mathbf{h} \overleftarrow{\frac{\p}{\p p_\gamma}} \right)
\left( \overrightarrow{\frac{\p}{\p q_\gamma}} \mathbf{f} \right) =
\sum_\gamma \kappa_\gamma \mathbf{h}_\gamma \frac{\p \mathbf{f}}{\p q_\gamma} + \Order(p),
$$
hence
$$
\p_\CH \mathbf{f} = \sum_\gamma \kappa_\gamma \mathbf{h}_\gamma 
\frac{\p \mathbf{f}}{\p q_\gamma}.
$$
In particular, $\p_\CH$ acts on each generator $q_\gamma \in \aA$ as
$$
\p_\CH q_\gamma = \kappa_\gamma \mathbf{h}_\gamma = 
\sum_u \frac{\epsilon(u) \kappa_\gamma}{\Aut^\sigma(u)}
e^A q^{\boldsymbol{\gamma}^-},
$$
where the sum is over all $\RR$-equivalence classes of index~$1$ 
$J$-holomorphic curves $u$ with genus zero, one positive end at~$\gamma$,
and negative ends $\boldsymbol{\gamma}^-$, and homology class 
$A \in H_2(M) / G$.

\subsection{Algebraic overtwistedness}
\label{sec:OT}

Even the simplest of the three differential graded algebras described above
is too large to compute in most cases.  The major exception is the case of
overtwisted contact manifolds.

\begin{thm}
\label{thm:OT}
If $(M,\xi)$ is overtwisted, then $HC_*(M,\xi;R) = 0$ for all choices of the
coefficient ring~$R$.
\end{thm}
\begin{remark}
\label{remark:unit}
If $X$ is an algebra with unit, then saying $X = 0$ is equivalent to saying
that $1 = 0$ in~$X$.
\end{remark}

The notion of overtwisted contact structures in dimension three was introduced by 
Eliashberg in \cite{Eliashberg:overtwisted}, who proved that they are \emph{flexible}
in the sense that their classification up to isotopy reduces to the purely
obstruction-theoretic classification of almost contact structures up to homotopy.
This means in effect that an overtwisted contact structure carries no distinctly
contact geometric information, so it should not be surprising when ``interesting''
contact invariants such as $HC_*(M,\xi)$ vanish.  The three-dimensional
case of Theorem~\ref{thm:OT}
seems to have been among the earliest insights about SFT: its first appearance
in the literature was in \cite{Eliashberg:invariants}, and a proof later appeared
in a paper by Mei-Lin Yau \cite{Yau:overtwisted}, which includes a brief appendix
sketching Eliashberg's original proof.  We will discuss Eliashberg's proof in
detail in Lecture~\ref{lec:torsion}.

The definitive higher-dimensional notion of overtwistedness was introduced a few
years ago by Borman-Eliashberg-Murphy \cite{BEM}, following earlier steps in 
this direction by Niederkr\"uger \cite{Plastikstufe} and others.  There are now 
two known proofs
of Theorem~\ref{thm:OT} in higher dimensions: the first uses the fact that since
overtwisted contact manifolds are flexible, they always admit an embedding of
a \emph{plastikstufe}, which implies vanishing of contact homology by an
unpublished result of Bourgeois and Niederkr\"uger 
(see \cite{Bourgeois:contactSurvey}*{Theorem~4.10} for a sketch of the
argument).  The second argument
appeals to an even more recent result of Casals-Murphy-Presas \cite{CasalsMurphyPresas}
showing that $(M,\xi)$ is overtwisted if and only if it is supported
by a negatively stabilized open book, in which case $HC_*(M,\xi) = 0$ was proven
by Bourgeois and van Koert \cite{BourgeoisVanKoert}.

It is not known whether the vanishing of contact homology \emph{characterizes}
overtwistedness, i.e.~there are not yet any examples of tight contact manifolds
with $HC_*(M,\xi) = 0$.  I will go out on a limb and say that such examples
seem unlikely to exist in dimension three but are much more likely in
higher dimensions; in fact various candidates are known 
\cites{MassotNiederkruegerWendl,ChiangDingVanKoert:nonfillable},
but we do not yet have adequate methods to prove that any of them are tight.
The analogous question about Legendrian submanifolds and relative contact
homology was recently answered by Ekholm \cite{Ekholm:nonLoose}, giving 
examples of Legendrians
that are not \emph{loose} in the sense of Murphy \cite{Murphy:loose} but
have vanishing Legendrian contact homology.

Nevertheless, the lack of known counterexamples has given rise to the
following definition.

\begin{defn}
\label{defn:algOT}
A closed contact manifold $(M,\xi)$ is \defin{algebraically overtwisted}
if $HC_*(M,\xi;R) = 0$ for every choice of the coefficient ring~$R$.
\end{defn}
\begin{remark}
\label{remark:coefs}
The coefficient ring is not always mentioned in statements of the above
definition, but it should be.  We will see in \S\ref{sec:algtorsion} below
that this detail makes a difference to issues like symplectic filling
obstructions.  Note that for any nested pair of subgroups
$G \subset G' \subset H_2(M)$, the natural projection
$H_2(M) / G' \to H_2(M) / G$ induces an algebra homomorphism
$$
HC_*(M,\xi;\QQ[H_2(M) / G']) \to HC_*(M,\xi;\QQ[H_2(M) / G]).
$$
Since algebra homomorphisms necessarily map $1 \mapsto 1$ and
$0 \mapsto 0$, the target of this map must vanish whenever its
domain does, so for checking Definition~\ref{defn:algOT}, it suffices
to check the case $R = \QQ[H_2(M)]$.
\end{remark}

We've seen above that there exist algebra homomorphisms
\begin{equation}
\label{eqn:homs}
H_*^{\Weyl}(M,\xi;R) \to H_*^\RSFT(M,\xi;R) \to HC_*(M,\xi;R),
\end{equation}
thus the vanishing of either of the algebras $H_*^{\Weyl}(M,\xi;R)$
or $H_*^\RSFT(M,\xi;R)$ with all coefficient rings $R$ 
is another sufficient condition for algebraic overtwistedness.
Bourgeois and Niederkr\"uger observed that, in fact, these conditions
are also necessary:

\begin{thm}[\cite{BourgeoisNiederkrueger:algebraically}]
For any coefficient ring~$R$, the following conditions are equivalent:
\begin{enumerate}
\item $HC_*(M,\xi;R) = 0$,
\item $H_*^\RSFT(M,\xi;R) = 0$,
\item $H_*^\SFT(M,\xi;R) = 0$.
\end{enumerate}
\end{thm}
\begin{proof}
The implications (3) $\Rightarrow$ (2) $\Rightarrow$ (1) are immediate from
the algebra homomorphisms \eqref{eqn:homs}, thus it will suffice to
prove (1) $\Rightarrow$ (3).  Suppose $1 = 0 \in HC_*(M,\xi;R)$, which
means $\p_\CH \mathbf{f} = 1$ for some $\mathbf{f} \in \aA$.  Using the
obvious inclusion $\aA \hookrightarrow \Weyl$, this means
$$
D_{\mathbf{H}} \mathbf{f} = 1 - \mathbf{G},
$$
where $\mathbf{G} = \Order(p,\hbar)$, i.e.~$\mathbf{G}$
is a sum of terms that all contain at least one $p_\gamma$
variable or a power of~$\hbar$.  It follows that $\mathbf{G}^k = \Order(p^k,\hbar^k)$
for all $k \in \NN$, and the infinite sum
$$
\sum_{k=0}^\infty \mathbf{G}^k
$$
is therefore an element of $\Weyl$, as the coefficient in front of any fixed
monomial $\hbar^k p^{\boldsymbol{\gamma}}$ in this sum is a polynomial
function of the $q$ variables.  This sum is then a multiplicative inverse
of $1 - \mathbf{G}$, and since
$$
0 = D_{\mathbf{H}}^2 \mathbf{f} = 0 = -D_{\mathbf{H}} \mathbf{G},
$$
it also satisfies $D_{\mathbf{H}} \left((1 - \mathbf{G})^{-1}\right) = 0$.
Using the fact that $D_{\mathbf{H}}$ is a derivation, we therefore have
$$
D_{\mathbf{H}} \left( (1 - \mathbf{G})^{-1} \mathbf{f} \right) =
(1 - \mathbf{G})^{-1} (1 - \mathbf{G}) = 1,
$$
implying $1 = 0 \in H_*^\SFT(M,\xi;R)$.
\end{proof}

\section{SFT generating functions for cobordisms}

All invariance proofs in SFT are based on a generating function analogous
to $\mathbf{H}$ that counts index~$0$ holomorphic curves in symplectic
cobordisms.  The basic definition is a straightforward extension of
what we saw in Lecture~\ref{lec:H}, but there is an added wrinkle due to the fact
that, in general, one must include \emph{disconnected} curves in the count.  

\subsection{Weak, strong and stable cobordisms}
\label{sec:weakStrong}

First some remarks about the category we are working in.  Since the stated
purpose of SFT is to define invariants of contact structures, we have been
working since Lecture~\ref{lec:H} with symplectizations of contact manifolds rather
than more general stable Hamiltonian structures.  We've made use of this
restriction on several occasions, namely so that we can assume:
\begin{enumerate}
\item All nontrivial holomorphic curves in $\RR \times M$ have at least one
positive puncture;
\item The energy of a holomorphic curve in $\RR \times M$ can be bounded
in terms of its positive asymptotic orbits.
\end{enumerate}
It will be useful however for certain applications to permit a slightly
wider class of stable Hamiltonian structure.  Recall that a hypersurface
$V$ in an almost complex manifold $(W,J)$ is called \defin{pseudoconvex}
if the maximal complex subbundle
$$
\xi := TV \cap J(TV) \subset TV
$$
defines a contact structure on $V$ whose canonical conformal symplectic bundle 
structure tames $J|_{\xi}$.  For example, if $\alpha$ is a contact form on $M$ and
$J \in \jJ(\alpha)$, then each of the hypersurfaces $\{\text{const}\} \times M$
is pseudoconvex in $(\RR \times M,J)$.  The contact structure $\xi$
induces an orientation on the hypersurface~$V$; if $V$ comes with its own
orientation (e.g.~as a boundary component of~$W$), then we call it
\emph{pseudoconvex} if $\xi$ is a positive contact structure with respect
to this orientation, and \defin{pseudoconcave} otherwise.  For example,
if $(W,\omega)$ is a symplectic cobordism from $(M_-,\xi_-)$ to $(M_+,\xi_+)$
and $J \in \jJ(W,\omega,\alpha_+,\alpha_-)$, then $M_+$ is pseudoconvex
and $M_-$ is pseudoconcave.

\begin{defn}
\label{defn:psConvex}
Given an odd-dimensional manifold $M$, we will say that an almost complex
structure $J$ on $\RR \times M$ is \defin{pseudoconvex}
if $\{r\} \times M$ is a pseudoconvex hypersurface
in $(\RR \times M,J)$ for every $r \in \RR$, with the induced orientation
such that $\p_r$ and $\{r\} \times M$ are positively transverse.
\end{defn}

If $\hH = (\omega,\lambda)$ is a stable Hamiltonian structure on~$M$, then
pseudoconvexity of $J \in \jJ(\hH)$ imposes conditions on~$\hH$, in particular
$\lambda$ must be a contact form.  It also requires $J|_\xi$ to be tamed by 
$d\lambda|_\xi$, but unlike the case when $J \in \jJ(\lambda)$,
$J|_{\xi}$ need not be \emph{compatible} with it, i.e.~the positive bilinear
form $d\lambda(\cdot,J\cdot)|_\xi$ need not be symmetric.  As always,
$J|_\xi$ must be compatible with $\omega|_\xi$, but $\omega$ need not be an
\emph{exact} form for this to hold---the freedom to change $[\omega] \in
H^2_\dR(M)$ will be the main benefit of this generalization, particularly
when we discuss weak symplectic fillings below.

\begin{prop}
Suppose $\hH = (\omega,\lambda)$ is a stable Hamiltonian structure on a closed
manifold $M$ and $J \in \jJ(\hH)$ is pseudoconvex.
Then all nonconstant finite-energy
$J$-holomorphic curves in $\RR \times M$ have at least one positive puncture,
and their energies satisfy a uniform upper bound in terms of the periods of
their positive asymptotic orbits.
\end{prop}
\begin{proof}
It is straightforward to check that either of the two proofs of
Proposition~\ref{prop:onePositive} given in Lecture~\ref{lec:tight3tori} generalizes to any $J$ on $\RR \times M$
that is pseudoconvex.  In particular, pseudoconvexity implies that
if $u : (\dot{\Sigma},j) \to (\RR \times M,J)$ is a $J$-holomorphic curve,
then $u^*d\lambda \ge 0$, with equality only at points where $u$ is tangent
to $\p_r$ and the Reeb vector field.  Stokes' theorem thus gives
\begin{equation}
\label{eqn:Stokes}
0 \le \int_{\dot{\Sigma}} u^*d\lambda = \sum_{z \in \Gamma^+} T_z -
\sum_{z \in \Gamma^-} T_z,
\end{equation}
where $T_z > 0$ denotes the period of the asymptotic orbit at each 
positive/negative puncture $z \in \Gamma^\pm$.  Since $J|_\xi$ is also tamed
by $\omega|_\xi$ and $\omega$ annihilates the Reeb vector field, we similarly 
have $u^*\omega \ge 0$, with the same condition for equality, and the compactness
of $M$ then implies an estimate of the form
$$
0 \le u^*\omega \le c u^*d\lambda
$$
for every $J$-holomorphic curve $u : (\dot{\Sigma},j) \to (\RR \times M,J)$,
with a constant $c > 0$ that depends only on $M$, $\hH$ and~$J$.  In light of
\eqref{eqn:Stokes}, this implies an upper bound on $\int_{\dot{\Sigma}} u^*\omega$ in terms
of the periods $T_z$ for $z \in \Gamma^+$.  Writing $\omega_\varphi =
\omega + d(\varphi(r) \lambda)$ for suitable $C^0$-small increasing functions
$\varphi : \RR \to \RR$, we can then apply Stokes' theorem to the second term in
$$
E(u) = \sup_\varphi \int_{\dot{\Sigma}} u^*\omega_\varphi = 
\int_{\dot{\Sigma}} u^*\omega + \sup_\varphi \int_{\dot{\Sigma}} u^*d(\varphi(r) \lambda),
$$
implying a similar upper bound for~$E(u)$.
\end{proof}

\begin{cor}
\label{cor:chainCpxs}
For any stable Hamiltonian structure $\hH = (\omega,\lambda)$ with a nondegenerate
Reeb vector field $R_\hH$ and a pseudoconvex $J \in \jJ(\hH)$, one can use closed
$R_\hH$-orbits and count $J$-holomorphic curves in $\RR \times M$ to define
the chain complexes $(\aA[[\hbar]],\mathbf{D}_\SFT)$, 
$(\Weyl,D_{\mathbf{H}})$, $(\Poisson,d_{\mathbf{h}})$ and
$(\aA,\p_\CH)$.
\end{cor}

We shall denote the homologies of the above chain complexes with coefficients
in $R = \QQ[H_2(M) / G]$ by
$$
H_*^\SFT(M,\hH,J;R), \quad
H_*^{\Weyl}(M,\hH,J;R),\quad
H_*^\RSFT(M,\hH,J;R),\quad
HC_*(M,\hH,J;R).
$$
We make no claim at this point about these homologies being invariant.
For the examples that we actually care about, this will turn out to be an
irrelevant question due to Proposition~\ref{prop:weakVanish} and
Exercise~\ref{EX:weakVanish} below.

\begin{example}
\label{ex:perturbedJ}
Suppose $\alpha$ is a contact form on $(M,\xi)$ and $\hH = (\Omega,\alpha)$ is a
stable Hamiltonian structure.  Then for all constants $c > 0$ sufficiently
large, $\hH_c := (\Omega + c\, d\alpha,\alpha)$ is also a stable Hamiltonian
structure and there exists a pseudoconvex $J_c \in \jJ(\hH_c)$.  To see the latter,
notice that $\hH_c' := \left(\frac{1}{c}\Omega + d\alpha,\alpha\right)$ is
another family of stable Hamiltonian structures, with $\jJ(\hH_c') = \jJ(\hH_c)$
for all~$c$, and $\hH_c' \to (d\alpha,\alpha)$ as $c \to \infty$.  Thus one
can select $J_c \in \jJ(\hH_c)$ converging to some $J_\infty \in \jJ(\alpha)$ as
$c \to \infty$, and these are pseudoconvex for $c > 0$ sufficiently large
since $J_\infty$~is.
\end{example}

\begin{prop}
\label{prop:weakVanish}
In the setting of Example~\ref{ex:perturbedJ}, assume $\alpha$ is nondegenerate
and $J_\infty \in\jJ(\alpha)$ is generic.  If $HC_*(M,\xi;R) = 0$,
then $HC_*(M,\hH_c,J_c;R)$ also vanishes for all $c > 0$ sufficiently large.
\end{prop}
\begin{proof}
We will assume in the following that the usual (unrealistic) transversality
assumptions hold, but the essential idea of the argument would not change in
the presence of abstract perturbations.

Let $(\aA,\p_\CH^\infty)$ denote the contact homology chain complex generated by
closed $R_\alpha$-orbits, with $\p_\CH^\infty$ counting $J_\infty$-holomorphic 
curves in $\RR \times M$.  The assumption $HC_*(M,\xi;R) = 0$ means there 
exists an element $\mathbf{f} \in \aA$ with $\p_\CH^\infty \mathbf{f} = 1$.  
Here  $\mathbf{f}$ is a polynomial function of the $q_\gamma$ variables, and 
$\p_\CH^\infty \mathbf{f}$ counts a specific finite set of Fredholm regular 
index~$1$ curves in $(\RR \times M,J_\infty)$.  Now let $(\aA,\p_\CH^c)$ 
denote the chain complex
for $HC_*(M,\hH_c,J_c;R)$, and notice that since the stable Hamiltonian
structures $(d\alpha,\alpha)$ and $\hH_c$ define matching Reeb vector fields,
the set of generators is unchanged.  There is also no change to this complex
if we replace $\hH_c = (\Omega + c\, d\alpha,\alpha)$ by
$\hH_c' = \left(\frac{1}{c}\Omega + d\alpha,\alpha\right)$: this changes
the energies of individual $J_c$-holomorphic curves, but the sets of
finite-energy curves are still the same in both cases.
We can assume $J_c \to J_\infty$ in $C^\infty$ as 
$c \to \infty$.  The implicit function theorem then extends each of the
finitely many $J_\infty$-holomorphic curves counted by $\p^\infty \mathbf{f}$ 
uniquely to a smooth $1$-parameter family of $J_c$-holomorphic curves for $c > 0$
sufficiently large.\footnote{In case you are concerned about the parametric
moduli space being an orbifold instead of a manifold, just add asymptotic
markers so that there is no isotropy, and divide by the appropriate
combinatorial factors to count.}  We claim that these are the only curves
counted by $\p_\CH^c \mathbf{f}$ when $c > 0$ is large.  Indeed, there
would otherwise exist a sequence $c_k \to \infty$ for which additional
$J_{c_k}$-holomorphic index~$1$ curves $u_k$ contribute to 
$\p_\CH^{c_k} \mathbf{f}$, and since $\mathbf{f}$ has only finitely many
terms representing possible positive asymptotic orbits, we can find a
subsequence for which all the $u_k$ have the same positive asymptotic orbits.
A further subsequence then has all the same negative asymptotic orbits as
well since the Reeb flow is nondegenerate and the total period of the
negative orbits is bounded by the total period of the positive orbits.
Finally, since the sequence of stable Hamiltonian structures $\hH_{c_k}'$ 
converges to $(d\alpha,\alpha)$, the curves $u_k$ have uniformly bounded
energy with respect to $\hH_{c_k}'$, so that SFT compactness yields a
subsequence converging to a $J_\infty$-holomorphic building of index~$1$,
which can only be one of the curves counted by $\p_\CH^\infty \mathbf{f}$.
This contradicts the uniqueness in the implicit function theorem and thus
proves the claim.  We conclude that for all $c > 0$ sufficiently large,
$\p_\CH^c \mathbf{f} = 1$.
\end{proof}

\begin{defn}
\label{defn:psConvexCob}
Assume $(W,\omega)$ is a symplectic cobordism with stable boundary
$\p W = -M_- \sqcup M_+$, with induced stable Hamiltonian structures
$\hH_\pm = (\omega_\pm,\lambda_\pm)$ at~$M_\pm$, and suppose $J$ is an
almost complex structure on the completion $\widehat{W}$ that is
$\omega$-tame on $W$ and belongs to $\jJ(\hH_\pm)$ on the cylindrical ends.
We will say that $J$ is
\defin{pseudoconvex near infinity}\footnote{If I were being hypercorrect about
use of language, I might insist on saying that $J$ is ``pseudoconvex near 
$+\infty$ and \emph{pseudoconcave} near~$-\infty$,'' as the orientation
reversal at the negative boundary makes $M_-$ technically
a pseudoconcave hypersurface in $(\widehat{W},J)$, not pseudoconvex.
But this definition will only be useful to us in 
cases where $M_- = \emptyset$, so my linguistic guilt is limited.}
if the $\RR$-invariant almost complex
structures $J_\pm$ defined by restricting $J$ to $[0,\infty) \times M_+$
and $(-\infty,0] \times M_-$ are both pseudoconvex.
\end{defn}

Note that the condition on $J$ in the above definition can only be satisfied
if $\lambda_\pm$ are both positive contact forms on~$M_\pm$, but the
$2$-forms $\omega_\pm$ need not be exact.

Proving contact invariance of SFT requires counting curves in trivial exact
symplectic cobordisms, but it is also natural to try to say things about
non-exact \defin{strong} symplectic cobordisms using SFT.\footnote{By 
\emph{strong cobordism}, we mean the usual notion of a 
compact symplectic manifold with convex and/or concave boundary components 
(see~\S\ref{sec:cobordisms}).  The word ``strong'' is included in order to contrast 
this notion with its weaker cousin described in Definition~\ref{defn:weak}.}
These fit naturally into our previously established picture since every
strong cobordism has collar neighborhoods near the boundary in which it
matches the symplectization of a contact manifold.  The following more
general notion of cobordism is also natural from a contact topological
perspective, but fits less easily into the SFT picture.

\begin{defn}[\cite{MassotNiederkruegerWendl}]
\label{defn:weak}
Given closed contact manifolds $(M_+,\xi_+)$ and $(M_-,\xi_-)$ of dimension
$2n-1$, a \defin{weak symplectic cobordism} from $(M_-,\xi_-)$ to $(M_+,\xi_+)$
is a compact symplectic manifold $(W,\omega)$ with $\p W = -M_- \sqcup M_+$
admitting an $\omega$-tame almost complex structure $J$ for which the
almost complex manifold $(W,J)$ is pseudoconvex at $M_+$ and pseudoconcave
at~$M_-$, with
$$
\xi_\pm = T M_\pm \cap J(T M_\pm).
$$
\end{defn}

Weak cobordisms are characterized by the existence of a tame almost 
complex structure $J$ whose restriction to $\xi_\pm$ is tamed by \emph{two}
symplectic bundle structures, $\omega|_{\xi_\pm}$ and $d\alpha_\pm|_{\xi_\pm}$
(for any choices of contact forms $\alpha_\pm$ defining~$\xi_\pm$).
Notice that in dimension~$4$, the second condition is mostly vacuous, and 
the weak cobordism condition just reduces to
$$
\omega|_{\xi_\pm} > 0.
$$
In this form, the low-dimensional case of Definition~\ref{defn:weak} has
been around since the late 1980's, and there are many interesting results
about it, e.g.~examples of contact $3$-manifolds that are weakly but not
strongly fillable \cites{Giroux:plusOuMoins,Eliashberg:fillableTorus}.
We will see in \S\ref{sec:algtorsion} that this distinction is detectable via SFT.
Higher-dimensional examples of this phenomenon were found in 
\cite{MassotNiederkruegerWendl}.

One major difference between weak and strong cobordisms is that the latter are
always exact near the boundary, as the Liouville vector field is dual to a
primitive of~$\omega$.  It turns out that up to deformation, weak fillings
that are exact at the boundary are the same thing as strong fillings---this 
was first observed by Eliashberg in dimension three
\cite{Eliashberg:contactProperties}*{Prop.~3.1}, and was extended to higher
dimensions in \cite{MassotNiederkruegerWendl}:

\begin{prop}
\label{prop:weakExact}
Suppose $(W,\omega)$ is a weak filling of a $(2n-1)$-dimensional contact manifold 
$(M,\xi)$ such that $\omega|_{TM}$ is exact.  Then after a homotopy of $\omega$
through a family of symplectic forms that vary only in a collar
neighborhood of~$\p W$ and define weak fillings of $(M,\xi)$, 
$(W,\omega)$ is a strong filling of $(M,\xi)$.
\end{prop}
\begin{proof}
Choose any contact form $\alpha$ for~$\xi$, denote its Reeb vector field
by~$R_\alpha$, and let $\Omega = \omega|_{TM}$.
Identify a collar neighborhood of $\p W$ in $W$ smoothly with
$(-\epsilon,0] \times M$, with the coordinate on $(-\epsilon,0]$ denoted
by~$r$, such that $\p_r$ and $R_\alpha$ span the symplectic complement
of $\xi$ at $\p W$ and satisfy $\omega(\p_r,R_\alpha) = 1$.
Then $\omega$ and $\Omega + d(r\alpha)$ are cohomologous symplectic
forms on $(-\epsilon,0] \times M$ that match at $r=0$, hence a Moser
deformation argument implies they are isotopic.  We can therefore assume
without loss of generality that $\omega = \Omega + d(r \alpha)$ on the collar
near~$\p W$.

By assumption, $\Omega = d\eta$ for some $1$-form $\eta$ on~$M$, and since
$(W,\omega)$ is a weak filling of $(M,\xi = \ker\alpha)$, we can choose a
complex structure $J_\xi$ on $\xi$ that is tamed by both 
$d\alpha|_\xi$ and~$d\eta|_\xi$.
Now choose a smooth cutoff function $\beta : [0,\infty) \to [0,1]$ that
has compact support and equals $1$ near~$0$.  We claim that
$$
\omega := d(\beta(r) \eta) + d(r\alpha)
$$
is a symplectic form on $[0,\infty) \times M$ if $|\beta'|$ is sufficiently
small.  Indeed, writing $\omega = dr \wedge (\alpha + \beta'(r)\, \eta) +
[\beta(r)\, d\eta + r\, d\alpha]$, we have
$$
\omega^n = n\, dr \wedge \alpha \wedge [\beta(r)\, d\eta + r\, d\alpha]^{n-1} +
n \beta'(r) \, dr \wedge \eta \wedge [\beta(r)\, d\eta + r\, d\alpha]^{n-1}.
$$
The first term is positive and bounded away from zero since 
$d\eta|_\xi$ and $d\alpha|_\xi$ both tame $J_\xi$, hence do does 
$\beta\, d\eta + r\, d\alpha|_{\xi}$.  The second term is then harmless if
$|\beta'|$ is sufficiently small, proving $\omega^n > 0$.  

This defines an
extension of the original weak filling to a symplectic completion
$\widehat{W} = W \cup_M \left( [0,\infty) \times M\right)$, and for each $r_0 \ge 0$, the 
compact subdomains defined by $r \le r_0$ define weak fillings 
of $(\{r_0\} \times M,\xi)$ since $\omega|_\xi = (\beta(r_0)\, d\eta + r_0\, d\alpha)|_\xi$
also tames~$J_\xi$.  Notice that for $r_0$ sufficiently large, the
$d\eta$ term disappears, so $\omega$ has a primitive that restricts to
$\{r_0\} \times M$ as a contact form for~$\xi$, meaning we have a \emph{strong} 
filling of this hypersurface.  The desired deformation of $\omega$ can
therefore be defined by pulling back via a smooth family of diffeomorphisms
$(-\epsilon,0] \to (-\epsilon,r_0]$, where $r_0$ varies from $0$ to a
sufficiently large constant.
\end{proof}

Unlike strong cobordisms, being a weak cobordism is an open
condition: if $(W,\omega)$ is a weak cobordism, then so is
$(W,\omega + \epsilon \sigma)$ for any $\epsilon > 0$ sufficiently small and
a closed $2$-form $\sigma$, which need not be exact at~$\p W$.
As a consequence, the cylindrical ends of a completed weak cobordism cannot
always be deformed to look like the symplectization of a contact manifold.  
This is where Definition~\ref{defn:psConvexCob} comes in useful.
The proof of the next lemma is very much analogous to
Proposition~\ref{prop:weakExact}.

\begin{lemma}[\cite{MassotNiederkruegerWendl}*{Lemma~2.10}]
\label{lemma:weakNonExact}
Suppose $(W,\omega)$ is a weak filling of a $(2n-1)$-dimensional contact manifold 
$(M,\xi)$, $\alpha$ is a contact form for~$\xi$ and $\Omega$ is a closed 
$2$-form on $M$ with $[\Omega] = [\omega|_{TM}] \in H^2_\dR(M)$.  
Then for any constant $c > 0$ sufficiently large,
after a homotopy of $\omega$
through a family of symplectic forms that vary only in a collar
neighborhood of~$\p W$ and define weak fillings of $(M,\xi)$, 
$\omega|_{TM} = \Omega + c\, d\alpha$.
\qed
\end{lemma}

The following result then provides a suitable model that can be used as 
$\Omega$ in the above lemma when $\omega|_{TM}$ is nonexact.  The statement
below is restricted to the case where $[\omega|_{TM}]$ is a rational cohomology
class; the reason for this is that it relies on a Donaldson-type existence
result for contact submanifolds obtained as zero-sets of approximately
holomorphic sections, due to Ibort, Marti\'inez-Torres and Presas
\cite{IbortMTPresas}.  It seems likely that the rationality condition could
be lifted with more work, and in dimension three this is known to be true;
see \cite{NiederkruegerWendl}*{Prop.~2.6}.

\begin{lemma}[\cite{CieliebakVolkov}*{Prop.~2.18}]
\label{lemma:weakSHS}
For any rational cohomology class $\eta \in H^2(M;\QQ)$ on a closed
$(2n-1)$-dimensional contact manifold $(M,\xi)$, there exists a closed
$2$-form $\Omega$ and a nondegenerate contact form $\alpha$ for $\xi$
such that $(\Omega,\alpha)$ is a stable Hamiltonian structure.
\qed
\end{lemma}

Combining all of the above results (including Example~\ref{ex:perturbedJ})
proves:

\begin{prop}
\label{prop:weakFilling}
Suppose $(W,\omega)$ is a weak filling of a $(2n-1)$-dimensional contact
manifold $(M,\xi)$ such that $[\omega|_{TM}] \in H^2_\dR(M)$ is rational
or $n=2$.  Fix a nondegenerate contact form $\alpha$ for~$\xi$.  Then there exists
a closed $2$-form $\Omega$ cohomologous to $\omega|_{TM}$ such that 
$\hH := (\Omega,\alpha)$ is a stable Hamiltonian structure, and for all
$c > 0$ sufficiently large, $\omega$ can be deformed in a collar neighborhood
of~$\p W$, through a family of symplectic forms defining weak fillings
of~$(M,\xi)$, to a new weak filling for which $\p W$ is also stable and
inherits the stable Hamiltonian structure 
$\hH_c := (\Omega + c\, d\alpha,\alpha)$.  In particular, after this deformation,
the completed stable filling admits a tame almost complex structure
that is pseudoconvex near infinity and may be assumed $C^\infty$-close to
any given $J \in \jJ(\alpha)$.
\qed
\end{prop}

We will use this in \S\ref{sec:algtorsion} to define obstructions to weak
fillability via SFT.

\begin{remark}
There is apparently no analogue of Propositions~\ref{prop:weakExact}
and~\ref{prop:weakFilling} for negative boundary components of weak cobordisms,
and this is one of a few reasons why they are not often discussed.
For example, if $L$ is a Lagrangian torus in the standard symplectic 
$4$-ball~$\DD^4$, then the complement of a neighborhood of $L$ in $B^4$
defines a strong cobordism from the standard contact $\TT^3$ to~$S^3$.
The symplectic form on this cobordism is obviously exact, but if any
result analogous to Proposition~\ref{prop:weakExact} were to hold at the
concave boundary, then we could deform it to a Liouville cobordism.
No such Liouville cobordism exists---it would imply that the Lagrangian
$L \subset B^4$ is exact, thus violating Gromov's famous theorem \cite{Gromov}
on exact Lagrangians.
\end{remark}

\subsection{Counting disconnected index~$0$ curves}
\label{sec:disconnected}

Fix a symplectic cobordism $(W,\omega)$ with stable boundary $\p W = -M_- \sqcup M_+$
carrying stable Hamiltonian structures $\hH_\pm = (\omega_\pm,\lambda_\pm)$,
along with a generic almost complex structure $J$ that is $\omega$-tame
on~$W$, belongs to $\jJ(\hH_\pm)$ on the cylindrical ends, and is
pseudoconvex near infinity.  This implies that the stabilizing 
$1$-forms $\lambda_\pm$ are both contact forms.  Let us also assume that
the $\lambda_\pm$ are both nondegenerate, and that the induced $\RR$-invariant
almost complex structures $J_\pm \in \jJ(\hH_\pm)$ are sufficiently generic
to achieve regularity for all holomorphic curves under consideration.
In particular, these assumptions mean that all the usual SFT chain complexes
are well defined for $(M_\pm,\hH_\pm,J_\pm;R_\pm)$ with any choice of
coefficient ring $R_\pm = \QQ[H_2(M_\pm) / G_\pm]$.  Denote the corresponding
SFT generating functions by $\mathbf{H}_\pm$.

Recall from Lecture~\ref{lec:H} that the auxiliary data on $M_+$ and $M_-$ includes a 
choice of capping surface $C_\gamma$ for each closed Reeb orbit~$\gamma$
(or a capping \emph{chain} with rational coefficients if $H_1(M_\pm)$ has 
torsion).  These surfaces satisfy
$$
\p C_\gamma = \sum_i m_i [C_i^\pm] - [\gamma],
$$
where the $m_i$ are integers and $C_i^\pm \subset M_\pm$ are fixed curves
forming a basis of~$H_1(M_\pm)$.  Assume $H_1(W)$ is torsion free, in which 
case the same is true of $H_1(M_+)$ and $H_1(M_-)$.  (Only minor modifications
are needed if this assumption fails to hold, see Remark~\ref{remark:torsion}.)
We can then fix the following additional auxiliary data:
\begin{enumerate}
\item A collection of \defin{reference curves}
$$
S^1 \cong C_1,\ldots,C_r \subset W
$$
whose homology classes from a basis of $H_1(W)$.
\item A unitary trivialization of $TW$ along each of the reference
curves $C_1,\ldots,C_r$, denoted collectively by~$\tau$.
\item A \defin{spanning surface} $S_i^\pm$ for each of the positive/negative
reference curves $C_i^\pm \subset M_\pm$, i.e.~a smooth map of a compact
and oriented surface with boundary into $W$ such that
$$
\p S_i^\pm = \sum_j m_{ji} [C_j] - [C_i^\pm]
$$
in the sense of singular $2$-chains, where $m_{ji} \in \ZZ$ are the unique
coefficients with $[C_i^\pm] = \sum_j m_{ji} [C_j] \in H_1(W)$.
\end{enumerate}
Now to any collections
of orbits $\boldsymbol{\gamma}^\pm = (\gamma_1^\pm,\ldots,\gamma_{k_\pm}^\pm)$
in $M_\pm$ and a relative homology class $A \in H_2(W,\bar{\boldsymbol{\gamma}}^+
\cup \bar{\boldsymbol{\gamma}}^-)$ with $\p A = \sum_i [\gamma_i^+] -
\sum_j [\gamma_j^-]$, we can associate an absolute homology class in two steps:
first add $A$ to suitable sums of the capping surfaces $C_{\gamma_i^\pm}$
producing a $2$-chain whose boundary is a linear combination of positive
and negative reference curves, then add a suitable linear combination of
the $S_i^\pm$ so that the boundary becomes the \emph{trivial} linear
combination of $C_1,\ldots,C_r$.  With this understood, 
we can now associate an absolute homology class
$$
[u] \in H_2(W)
$$
to any asymptotically cylindrical $J$-holomorphic curve $u :
(\dot{\Sigma},j) \to (\widehat{W},J)$, and this defines the notation
$\mM_{g,m}(J,A,\boldsymbol{\gamma}^+,\boldsymbol{\gamma}^-)$ with
$A \in H_2(W)$.  We now require the trivializations of $\xi_\pm$ along
each $C_i^\pm$ to be compatible with $\tau$ in the sense that they extend
to trivializations of $TW$ along the capping surfaces $S_i^\pm$.  With this
convention, the Fredholm index formula takes the expected form
$$
\ind(u) = (n-3) \chi(\dot{\Sigma}) + 2 c_1([u]) + \sum_{i=1}^{k_+} 
\muCZ(\gamma_i) - \sum_{j=1}^{k_-} \muCZ(\gamma_j).
$$
If $H_1(W)$ has torsion, then this whole discussion
can be adapted as in \S\ref{sec:torsion} by replacing integral homology
with rational homology and capping surfaces with capping chains, and the
Conley-Zehnder indices can be defined modulo~$2$.

We will also need to impose a compatibility condition relating the
coefficient rings $R_\pm = \QQ[H_2(M_\pm) / G_\pm]$ to a corresponding
choice on the cobordism~$W$.  Choose a subgroup $G \subset H_2(W)$ such that
\begin{equation}
\label{eqn:keromega}
\langle [\omega],A \rangle = 0 \quad \text{ for all } \quad A \in G,
\end{equation}
and such that the maps $H_2(M_\pm) \to H_2(W)$ induced by the inclusions
$M_\pm \hookrightarrow W$ send $G_\pm$ into~$G$.  If $[\omega] \ne 0
\in H^2_\dR(W)$, then we will have to deal with noncompact sequences of
$J$-holomorphic curves that have unbounded energy, so it becomes necessary
to ``complete'' $R$ to a \defin{Novikov ring} $\overline{R}$, which contains
$R$ but also includes infinite formal sums
$$
\sum_{i=1}^\infty c_i e^{A_i} \quad \text{ such that } \quad
\langle [\omega],A_i \rangle \to +\infty \text{ as } i \to \infty.
$$
Note that the evaluation $\langle [\omega],A \rangle \in \RR$ is well defined
for $A \in H_2(W) / G$ due to \eqref{eqn:keromega}.

Analogously to our definition of $\mathbf{H}$ in Lecture~\ref{lec:H}, the generating
function for index~$0$ curves in $\widehat{W}$ is defined as a formal
power series in the variables $\hbar$, $q_\gamma$ (for orbits in~$M_-$),
and $p_\gamma$ (for orbits in~$M_+$), with coefficients
in~$\overline{R}$:
\begin{equation}
\label{eqn:Fsigma}
\mathbf{F} = \sum_{u \in \mM_0^\sigma(J)}
\frac{\epsilon(u)}{|\Aut^\sigma(u)|}
\hbar^{g-1} e^A q^{\boldsymbol{\gamma}^-} p^{\boldsymbol{\gamma}^+},
\end{equation}
where $\mM_0^\sigma(J)$ denotes the moduli space of connected
$J$-holomorphic curves $u$ in $\widehat{W}$ with $\ind(u) = 0$ and only
good asymptotic orbits, modulo permutations of the punctures, and for each~$u$:
\begin{itemize}
\item $g$ is the genus of~$u$;
\item $A$ is the equivalence class of $[u] \in H_2(W)$ in $H_2(W) / G$;
\item $\boldsymbol{\gamma}^\pm = (\gamma_1^\pm,\ldots,\gamma_{k_\pm}^\pm)$
are the asymptotic orbits of $u$ after arbitrarily fixing orderings of its
positive and negative punctures;
\item $\epsilon(u) \in \{1,-1\}$ is the sign of $u$ as a point in the
$0$-dimensional component of $\mM^\$(J)$ (after choosing an ordering of the
punctures and asymptotic markers), relative to a choice of coherent 
orientations on~$\mM^\$(J)$.
\end{itemize}
As usual, the product 
$\epsilon(u) q^{\boldsymbol{\gamma}^-} p^{\boldsymbol{\gamma}^+}$ is
independent of choices.  We shall regard $\mathbf{F}$ as an element in an
enlarged operator algebra that includes $q$ and $p$ variables for good orbits
in both $M_+$ and $M_-$, related to each other by the supercommutation relations
$$
[p_{\gamma_-},q_{\gamma_+}] = [p_{\gamma_+},q_{\gamma_-}] = 
[q_{\gamma_-},q_{\gamma_+}] = [p_{\gamma_-},p_{\gamma_+}] = 0
$$
whenever $\gamma_-$ is an orbit in $M_-$ and $\gamma_+$ is an orbit in~$M_+$.
Since all curves counted by $\mathbf{F}$ have index~$0$, $\mathbf{F}$ is
homogeneous with degree
$$
|\mathbf{F}| = 0.
$$
Notice that for any fixed monomial $q^{\boldsymbol{\gamma}^-} p^{\boldsymbol{\gamma}^+}$,
the corresponding set of curves in $\mM_0^\sigma(J)$ may be infinite if
$\omega$ is nonexact, but SFT compactness implies that the set of such curves
with any given bound on $\int_{\dot{\Sigma}} u^*\omega$ is bounded.
As a consequence, the coefficient of 
$q^{\boldsymbol{\gamma}^-} p^{\boldsymbol{\gamma}^+}$ in $\mathbf{F}$ belongs
to the Novikov ring~$\overline{R}$.

Consider next the series
$$
\exp(\mathbf{F}) := \sum_{k=0}^\infty \frac{1}{k!} \mathbf{F}^k.
$$
We will be able to view this as a formal power series in $q$ and $p$ variables 
and a formal Laurent series in~$\hbar$ with coefficients in~$\overline{R}$, 
though it is not obvious at first glance whether its coefficients are
in any sense finite.  We will deduce this after interpreting it as
a count of \emph{disconnected} index~$0$ curves: first, write
\begin{equation*}
\begin{split}
\exp(\mathbf{F}) = \sum_{k=0}^\infty \frac{1}{k!} \Bigg(\sum_{(u_1,\ldots,u_k) \in (\mM_0^\sigma(J))^k}
& \frac{\epsilon(u_1)\ldots \epsilon(u_k)}{|\Aut^\sigma(u_1)| \ldots |\Aut^\sigma(u_k)|}
\hbar^{g_1 + \ldots + g_k - k} e^{A_1 + \ldots A_k} \\
& \cdot q^{\boldsymbol{\gamma}^-_1} p^{\boldsymbol{\gamma}^+_1} \ldots 
q^{\boldsymbol{\gamma}^-_k} p^{\boldsymbol{\gamma}^+_k} \Bigg).
\end{split}
\end{equation*}
Observe that since each of the curves $u_i \in \mM_0^\sigma(J)$ in this expansion
has index~$0$, the monomials $q^{\boldsymbol{\gamma}^-_i} p^{\boldsymbol{\gamma}^+_i}$
all have even degree and thus the order in which they are written does not matter.
Now for a given collection of distinct curves $v_1,\ldots,v_N$ and integers
$k_1,\ldots,k_N \in \NN$ with $k_1 + \ldots + k_N = k$, the various permutations
of
$$
(u_1,\ldots,u_k) := (\underbrace{v_1,\ldots,v_1}_{k_1},\ldots,\underbrace{v_N,\ldots,v_N}_{k_N}) \in
(\mM_0^\sigma(J))^k
$$
occur $\frac{k!}{k_1! \ldots k_N!}$ times in the above sum, so if we forget the
ordering, then the contribution of this particular $k$-tuple of curves 
to $\exp(\mathbf{F})$ is
$$
\frac{\epsilon(u_1) \ldots \epsilon(u_k)}{k_1! \ldots k_N! |\Aut^\sigma(u_1)| \ldots |\Aut^\sigma(u_k)|}
\hbar^{g_1 + \ldots + g_k - k} e^{A_1 + \ldots + A_k}
q^{\boldsymbol{\gamma}^-_1} p^{\boldsymbol{\gamma}^+_1} \ldots 
q^{\boldsymbol{\gamma}^-_k} p^{\boldsymbol{\gamma}^+_k}.
$$
Notice next that the denominator $k_1! \ldots k_N! |\Aut^\sigma(u_1)| \ldots |\Aut^\sigma(u_k)|$
is the order of the automorphism group of the \emph{disconnected} curve formed by the
disjoint union of $u_1,\ldots,u_k$: the extra factors $k_i!$ come from automorphisms that
permute connected components of the domain.  Thus $\exp(\mathbf{F})$ can also be written as
in \eqref{eqn:Fsigma}, but with $\mM_0^\sigma(J)$ replaced by the moduli space of
\emph{potentially disconnected} index~$0$ curves with unordered punctures, and $g - 1$ generalized
to $g_1 + \ldots + g_k - k$ for any curve that has $k$ connected components of genera 
$g_1,\ldots,g_k$.  One subtlety that was glossed over in the above discussion:
the sum also includes the unique curve with \emph{zero} components, i.e.~the ``empty''
$J$-holomorphic curve, which appears as the initial $1$ in the series expansion of $\exp(\mathbf{F})$.

With this interpretation of $\exp(\mathbf{F})$ understood, we can now address
the possibility that the infinite sum defining $\exp(\mathbf{F})$ might include
infinitely many terms for a given monomial 
$\hbar^m q^{\boldsymbol{\gamma}^-} p^{\boldsymbol{\gamma}^+}$, i.e.~that there
are infinitely many disconnected index~$0$ curves with fixed asymptotic
orbits and a fixed sum of the genera minus the number of connected components.
We claim that this can indeed, happen, but only if the curves belong to a
sequence of homology classes $A_i \in H_2(M) / G$ with $\langle [\omega],A_i \rangle
\to \infty$, hence the coefficient of 
$\hbar^m q^{\boldsymbol{\gamma}^-} p^{\boldsymbol{\gamma}^+}$ in
$\exp(\mathbf{F})$ belongs to the Novikov ring~$\overline{R}$.  The danger
here comes only from \emph{closed} curves, since a disjoint union of two curves
with punctures always has strictly more punctures.  Notice also that 
for any given tuples of orbits $\boldsymbol{\gamma}^\pm$, there exists a
number $c \in \RR$ depending only on these orbits and the chosen capping 
surfaces such that every (possibly disconnected) $J$-holomorphic curve
$u : \dot{\Sigma} \to \widehat{W}$ asymptotic to $\boldsymbol{\gamma}^\pm$
satisfies
$$
\langle [\omega],[u] \rangle \ge c.
$$
This follows from the fact that the integral of $\omega$ over the relative
homology class of $u$ always has a nonnegative integrand.

\begin{lemma}
\label{lemma:Novikov}
Given constants $C \in \RR$ and $k \in \ZZ$, there exists a number $N \in \NN$
such that if $u : (\Sigma,j) \to (\widehat{W},J)$ is a closed $J$-holomorphic 
curve satisfying $\int_{\Sigma} u^*\omega \le C$, with $m$ connected components 
of genera $g_1,\ldots,g_m$ satisfying $g_1 + \ldots + g_m - m = k$, then
$m \le N$.
\end{lemma}
\begin{proof}
Note first that for each integer $g \ge 0$, there is an \defin{energy threshold},
i.e.~a constant $c_g > 0$ such that every nonconstant closed and
connected $J$-holomorphic curve $u : \Sigma \to \widehat{W}$ of genus $g$ has
$$
\int_{\Sigma} u^*\omega \ge c_g.
$$
This is an easy consequence of SFT compactness: indeed, if there were no such
constant, then we would find a sequence $u_k : \Sigma \to \widehat{W}$ of
connected closed curves with genus~$g$ such that
$$
E(u_k) = \int_\Sigma u^*\omega \to 0;
$$
here we have used the fact that $\Sigma$ is closed and 
$\int_\Sigma u^*\omega_\varphi$ depends only on the homology class of $u$
in order to simplify the usual definition of energy for asymptotically
cylindrical curves.  SFT compactness then gives a subsequence of $u_k$ that
converges to a stable holomorphic building in which every component has zero 
energy and is therefore constant.  Since there are no marked points in the
picture, no such building exists, so this is a contradiction.

Now if $u$ is a disconnected curve satisfying the stated conditions, the bound 
on $\int_\Sigma u^*\omega$ combines with the energy threshold to give a
bound for each $g \ge 0$ on the number of connected components of $u$ with
genus~$g$.  In particular, there is a bound on the number of components
with genus $0$ or~$1$.  All other components contribute positively to the
left hand side of the relation $\sum_{i=1}^m (g_i - 1) = k$, so this
implies a universal bound on~$m$.
\end{proof}

\begin{cor}
Fix constants $C \in \RR$ and $k \in \ZZ$, and tuples of Reeb orbits
$\boldsymbol{\gamma}^\pm$, and assume that the usual transversality conditions
hold.  Then there exist at most finitely many potentially disconnected
$J$-holomorphic curves $u : \dot{\Sigma} \to \widehat{W}$ with index~$0$
such that the number of connected components $m$ and the genera
$g_1,\ldots,g_m$ of its components satisfy $g_1 + \ldots + g_m - m = k$.
\end{cor}
\begin{cor}
The expression $\exp(\mathbf{F})$ is a formal power series in $q$ and $p$
variables and a formal Laurent series in~$\hbar$, with coefficients in
the Novikov ring~$\overline{R}$.
\end{cor}

The necessity of considering disconnected curves becomes clear when one
tries to translate the compactness and gluing theory of 
$J$-holomorphic curves in $\widehat{W}$ into algebraic relations.
In particular, consider the $1$-dimensional moduli space of connected
index~$1$ curves in $\widehat{W}$ with genus~$g$.  The boundary points of
the compactification of this space consist of two types of buildings:
\begin{itemize}
\item[\textsc{Type~1}]: A main level of index~$0$ and an upper level of index~$1$;
\item[\textsc{Type~2}]: A main level of index~$0$ and a lower level of index~$1$.
\end{itemize}
This is clear under the usual transversality assumptions since regular curves
in $\widehat{W}$ must have index at least~$0$, while regular curves in the
symplectizations $\RR \times M_\pm$ have index at least~$1$ unless they
are trivial cylinders.  The building must also be connected and have
arithmetic genus~$g$, but there is nothing to guarantee that each individual
level is connected.  In fact, we already saw this issue in Lecture~\ref{lec:H} when
proving $\mathbf{H}^2 = 0$, but it was simpler to deal with there, because
disconnected regular curves of index~$1$ in a symplectization always have
a unique nontrivial component, while the rest are trivial cylinders.
In the cobordism $\widehat{W}$, on the other hand, a disconnected index~$0$
curve can be formed by any disjoint union of index~$0$ curves, all of which
are nontrivial.  Exponentiation provides a convenient way to encode
all data about disconnected curves in terms of connected curves.

Since the union of all buildings of types~1 and~2 described above forms the
boundary of a compact oriented $1$-manifold, the count of these buildings
is zero, and this fact is encoded in the so-called \defin{master equation}
\begin{equation}
\label{eqn:masterEqn}
\mathbf{H}_- \exp(\mathbf{F})|_{p_- = 0} - 
\exp(\mathbf{F}) \mathbf{H}_+|_{q_+ = 0} = 0,
\end{equation}
where the expressions ``$p_- = 0$'' and ``$q_+ = 0$'' mean that we discard
all terms in $\mathbf{H}_- \exp(\mathbf{F}) - \exp(\mathbf{F}) \mathbf{H}_+$
containing any variables $p_\gamma$ for orbits in $M_-$ or $q_\gamma$ for
orbits in~$M_+$.  The resulting expression is therefore a formal power
series in $q$ variables for orbits in $M_-$ and $p$ variables for
orbits in~$M_+$, representing a count of generally disconnected index~$1$ 
holomorphic buildings in $\widehat{W}$ with the specified asymptotics.
The various ways to form such buildings by choices of gluings is again
encoded by the commutator algebra.  The master equation \eqref{eqn:masterEqn}
can be used to prove the chain map property for counts of curves in
cobordisms, thus it is an essential piece of the invariance proof for
each of the homology theories introduced above.

\begin{exercise}
Fill in the details of the proof of \eqref{eqn:masterEqn}.
\end{exercise}

\section{Full SFT as a $BV_\infty$-algebra}

In this section we discuss the specific theory $H_*^\SFT(M,\xi;R)$, defined
as the homology of the chain complex $(\aA[[\hbar]],\mathbf{D}_\SFT)$.
The case $G = H_2(M)$ with trivial group ring coefficients 
$\QQ[H_2(M) / G] = \QQ$ will be abbreviated as
$$
H_*^\SFT(M,\xi) := H_*^\SFT(M,\xi;\QQ).
$$
As we defined it, $\mathbf{D}_\SFT$ acts on $\aA[[\hbar]]$ by treating the
generating function $\mathbf{H}$ as a differential operator via the substitution
\begin{equation}
\label{eqn:substitution}
p_\gamma = \kappa_\gamma \hbar \frac{\p}{\p q_\gamma}.
\end{equation}
According to \cite{CieliebakLatschev:propaganda}, this makes 
$(\aA[[\hbar]],\mathbf{D}_\SFT)$ into
a $BV_\infty$-algebra; we'll have no particular need to discuss here what
that means, but one convenient feature is the expansion
\begin{equation}
\label{eqn:expansion}
\mathbf{D}_\SFT = \frac{1}{\hbar} \sum_{k=1}^\infty 
\mathbf{D}_\SFT^{(k)} \hbar^k,
\end{equation}
in which each $\mathbf{D}_\SFT^{(k)} : \aA \to \aA$ is a differential operator
\emph{of order~$\le k$} (see \cite{CieliebakLatschev:propaganda}*{\S 5}).
For each $k \in \NN$, $\mathbf{D}_\SFT^{(k)}$ is a count of all index~$1$
holomorphic curves that have genus $g \ge 0$ and $m \ge 1$ positive
punctures such that $g + m = k$.  In particular, $\mathbf{D}_\SFT^{(1)}$
is simply the contact homology differential~$\p_\CH$, and the
expansion \eqref{eqn:expansion} implies together with
$\mathbf{D}_\SFT^2 = 0$ that $(\mathbf{D}_\SFT^{(1)})^2 = 0$, hence we
again see the chain complex for contact homology hidden inside a version
of the ``full'' SFT complex.

\subsection{Cobordism maps and invariance}
\label{sec:invariance}

One can use the master equation \eqref{eqn:masterEqn} to prove invariance
of $H_*^\SFT(M,\xi;R)$ by a straightforward generalization of the usual
Floer-theoretic argument.  Suppose $(W,d\lambda)$ is an exact symplectic
cobordism from $(M_-,\xi_-)$ to $(M_+,\xi_+)$ with $\lambda|_{TM_\pm} =
\alpha_\pm$, and choose a generic almost complex structure $J$ on
$\widehat{W}$ that is $d\lambda$-compatible on $W$ and restricts to the
cylindrical ends as generic elements
$J_\pm \in \jJ(\alpha_\pm)$.  Let $(\aA^\pm[[\hbar]],\mathbf{D}_\SFT^\pm)$
denote the chain complexes associated to the data $(\alpha_\pm,J_\pm)$,
and for simplicity in this initial discussion, choose the trivial
coefficient ring $R = \QQ$ for both.  We then define a map
$$
\boldsymbol{\Phi} : \aA^+[[\hbar]] \to \aA^-[[\hbar]] :
\mathbf{f} \mapsto \exp(\mathbf{F}) \mathbf{f}|_{q_+ = 0},
$$
where the generating function $\exp(\mathbf{F})$ is regarded as a
differential operator via the substitution \eqref{eqn:substitution},
with $e^A := 1$ for all $A \in H_2(W)$ since we are using trivial coefficients,
and ``$q_+ = 0$'' means that after applying $\exp(\mathbf{F})$ to change
$\mathbf{f}$ into a function of $q$ variables for orbits in both $M_+$ and
$M_-$, we discard all terms that involve orbits in~$M_+$.  
The exactness of the cobordism implies that negative powers of $\hbar$ do not 
appear in $\boldsymbol{\Phi}\mathbf{f}$, thus producing an element
of $\aA^-[[\hbar]]$: indeed, since there are no holomorphic curves
in $\widehat{W}$ without positive punctures, every term in $\mathbf{F}$
contains at least one $p$ variable, so that negative powers of $\hbar$ do
not appear in $\exp(\mathbf{F})$ after applying \eqref{eqn:substitution}.

The master equation for $\mathbf{F}$ now translates into the fact that
$\boldsymbol{\Phi}$ is a chain map,
$$
\mathbf{D}_\SFT^- \circ \boldsymbol{\Phi} - \boldsymbol{\Phi} \circ
\mathbf{D}_\SFT^+,
$$
thus it descends to homology.  The geometric meaning of
$\boldsymbol{\Phi}$ is straightforward to describe: analogous to
\eqref{eqn:DSFT} in Lecture~\ref{lec:H}, we can write
\begin{equation}
\label{eqn:PhiSFT}
\boldsymbol{\Phi} q^{\boldsymbol{\gamma}} =
\sum_{g=0}^\infty 
\sum_{\boldsymbol{\gamma}'} \hbar^{g + k - 1}
n_g(\boldsymbol{\gamma},\boldsymbol{\gamma}',k) q^{\boldsymbol{\gamma}'},
\end{equation}
where $n_g(\boldsymbol{\gamma},\boldsymbol{\gamma}',k)$ is a product of
some combinatorial factors with a signed count of disconnected
index~$0$ holomorphic curves with connected components of genus 
$g_1,\ldots,g_m$ satisfying $g_1 + \ldots + g_m - m = g - 1$, and with
positive ends at $\boldsymbol{\gamma}$ and negative ends at
$\boldsymbol{\gamma}'$, where $k$ is the number of positive ends.  

Let's discuss two applications of the cobordism map~$\boldsymbol{\Phi}$.
First, note that if $W$ is a \emph{trivial} symplectic cobordism
$[0,1] \times M$, then the above discussion can easily be generalized
with $(\aA^\pm,\mathbf{D}_\SFT^\pm)$ both defined over the same group ring
$R = \QQ[H_2(M) / G]$ for any choice of $G \subset H_2(M)$.  There is no
need to consider a Novikov ring in defining $\mathbf{F}$ here since the
cobordism is exact.  We therefore obtain a chain map with arbitrary group
ring coefficients, and extending this discussion along standard
Floer-theoretic principles will imply that the chain map is an isomorphism:
this can be used in particular to prove that $H_*^\SFT(M,\xi;R)$ does not
depend on the choices of contact form and almost complex structure.
There are two additional steps involved in this argument: first, one needs
to use a chain homotopy to prove that $\boldsymbol{\Phi}$ does not depend
on the choice of almost complex structure $J$ on~$\widehat{W}$.  Given a 
generic homotopy $\{J_s\}_{s \in [0,1]}$, the chain homotopy map
$$
\boldsymbol{\Psi} : \aA^+[[\hbar]] \to \aA^-[[\hbar]]
$$
is defined as a differential operator in the same manner as $\boldsymbol{\Phi}$,
but counting pairs $(s,u)$ where $s \in [0,1]$ is a parameter value for which
$J_s$ is nongeneric and $u$ is a disconnected $J_s$-holomorphic curve
in $\widehat{W}$ with index~$-1$.  We saw how this works for cylindrical
contact homology in Lecture~\ref{lec:tight3tori}, but there is a new subtlety now that should
be mentioned: in principle, a \emph{disconnected} index~$-1$ curve in
$\widehat{W}$ could have arbitrarily many components, including perhaps many
with index~$-1$ and others with arbitrarily large index.  Even worse,
the compactified $1$-dimensional space of pairs $(s,u)$ for $J_s$-holomorphic
curves $u$ of index~$0$ may include buildings that have symplectization levels
of index greater than $1$, balanced by disjoint unions of many index~$-1$
curves in the main level.  This sounds horrible, but it can actually be
ignored, for the following reason: first, since there are only finitely many pairs
$(s,u)$ where $u$ is a \emph{connected} $J_s$-holomorphic curve with
index~$-1$, one can (if transversality is achievable at all) use a genericity
argument to assume without loss of generality that for any given $s \in [0,1]$,
at most \emph{one} connected index~$-1$ curve exists.  This means that in any
building that has multiple index~$-1$ components, those components are just
multiple copies of the same curve.  Now, since that curve has odd index,
it is represented by a monomial $q^{\boldsymbol{\gamma}^-} p^{\boldsymbol{\gamma}^+}$
that contains an odd number of odd generators, and any nontrivial product of
such generators therefore \emph{disappears} in $\aA$ since odd generators 
anticommute with themselves.  This algebraic miracle encodes a convenient
fact about coherent orientations: whenever one of the horrible buildings
described above appears, one can reorder two of the index~$-1$ components to
produce from it a different building that lives in a moduli space with the
opposite orientation.  Gluing this building back together then produces a
continuation of the $1$-dimensional moduli space, so that the horrible
building can actually be interpreted as an ``interior'' point of the
$1$-dimensional space, rather than boundary.  The actual count of boundary
points is then exactly what we want it to be: it is represented algebraically
by the chain homotopy relation!

Finally, compositions of cobordism maps can be understood via a stretching
argument that is not substantially different from the case of cylindrical
contact homology.  Since the trivial cobordism with $\RR$-invariant data
gives a cobordism map that just counts trivial cylinders and is therefore
the identity, it follows that cobordism maps relating different pairs of
data $(\alpha_\pm,J_\pm)$ are always invertible, and this proves the
invariance of $H_*^\SFT(M,\xi;R)$.

The second application concerns nontrivial exact cobordisms, and it is
immediate from the fact that $\boldsymbol{\Phi}$ is a chain map:

\begin{thm}
\label{thm:cobordismMap}
Any exact cobordism $(W,d\lambda)$ from $(M_-,\xi_-)$ to $(M_+,\xi_+)$ gives
rise to a $\QQ[[\hbar]]$-linear map
$$
H_*^\SFT(M_+,\xi_+) \to H_*^\SFT(M_-,\xi_-).
$$
\qed
\end{thm}

It is much more complicated to say what happens in the event of a nonexact
cobordism, but slightly easier if we restrict our attention to fillings,
i.e.~the case with $M_- = \emptyset$.  Assume $(W,\omega)$ is a compact
symplectic manifold with stable boundary~$M$, inheriting a stable
Hamiltonian structure $\hH = (\Omega,\alpha)$ for which $\alpha$ is a
nondegenerate contact form, and assume also that the completion $\widehat{W}$ 
admits an almost complex structure $J$ that is $\omega$-tame on $W$ and has a
pseudoconvex restriction $J_+ \in \jJ(\hH)$ to the cylindrical end.
We saw in Proposition~\ref{prop:weakFilling} that these conditions can
always be achieved for a weak filling after deforming the symplectic
structure.  Let
$$
G := \ker[\omega] := \left\{ A \in H_2(W)\ |\ \langle [\omega],A \rangle = 0 \right\},
$$
and choose $G_+ \subset H_2(M)$ to be any subgroup such that the
map $H_2(M) \to H_2(W)$ induced by the inclusion $M \hookrightarrow W$
sends $G_+$ into~$G$.  In other words, $G_+$ can be any subgroup of
$\ker[\Omega] \subset H_2(M)$.  Define the group rings
$$
R_+ = \QQ[H_2(M) / G_+], \qquad R = \QQ[H_2(W) / \ker[\omega]],
$$
with the Novikov completion of $R$ denoted by~$\overline{R}$.  
The map $H_2(M) / G_+ \to H_2(W) / G$ induced by $M \hookrightarrow W$ then
gives a natural ring homomorphism 
\begin{equation}
\label{eqn:ringHom}
R_+ \to \overline{R}.
\end{equation}
If $\omega$ is not exact, then it may no longer be true that every term in
$\mathbf{F}$ has at least one $p$ variable.  Let us write
$$
\mathbf{F} = \mathbf{F}_0 + \mathbf{F}_1,
$$
where $\mathbf{F}_0$ contains no $p$ variables and $\mathbf{F}_1 = \Order(p)$,
i.e.~$\mathbf{F}_0$ counts all closed curves in~$\widehat{W}$, and
$\mathbf{F}_1$ counts everything else.  Since $\mathbf{F}_0$ and $\mathbf{F}_1$
have even degree, they commute, and thus
$$
\exp(\mathbf{F}) = \exp(\mathbf{F}_0) \exp(\mathbf{F}_1).
$$
where $\exp(\mathbf{F}_0)$ is an invertible element of 
$\overline{R}[[\hbar,\hbar^{-1}]]$ since $\exp(-\mathbf{F}_0) \exp(\mathbf{F}_0) = 1$.
By the master equation,
$$
\exp(\mathbf{F}_0) \exp(\mathbf{F}_1) \mathbf{H} = \Order(q),
$$
hence $\exp(\mathbf{F}_1) \mathbf{H} = \exp(-\mathbf{F}_0) \Order(q) = \Order(q)$
since $\exp(-\mathbf{F}_0)$ contains no $p$ variables.  Using the
substitution \eqref{eqn:substitution}, and using \eqref{eqn:ringHom} to map
coefficients in $R_+$ to~$\overline{R}$, it follows that $\exp(\mathbf{F}_1)$
gives rise to a differential operator
$$
\boldsymbol{\Phi} : \aA[[\hbar]] \to \overline{R}[[\hbar]] :
\mathbf{f} \mapsto \exp(\mathbf{F}_1) \mathbf{f}|_{q=0},
$$
which is a chain map to the SFT of the empty set with Novikov coefficients, 
meaning
$$
\boldsymbol{\Phi} \circ \mathbf{D}_\SFT = 0.
$$
This chain map counts the disconnected index~$0$ curves in $\widehat{W}$
whose connected components all have at least one positive puncture.

\begin{thm}
\label{thm:weakFillingMap}
Suppose $(W,\omega)$ is a compact symplectic manifold with stable boundary
$(M,\hH = (\Omega,\alpha))$, where $\alpha$ is a nondegenerate contact form, 
and its completion $\widehat{W}$ admits an almost complex structure that is
$\omega$-tame on $W$ and has a generic and pseudoconvex restriction $J_+ \in \jJ(\hH)$
to the cylindrical end.  Let $\overline{R}$ denote the Novikov completion
of $\QQ[H_2(W) / \ker[\omega]]$, and let $R_+ = \QQ[H_2(M) / G_+]$, where
$G_+ \subset H_2(M)$ is any subgroup on which the evaluation of
$[\Omega] \in H^2_\dR(M)$ vanishes.  
Then there exists an $\overline{R}[[\hbar]]$-linear map
$H_*^\SFT(M,\hH,J_+;R_+) \to \overline{R}[[\hbar]]$.
\qed
\end{thm}

\subsection{Algebraic torsion}
\label{sec:algtorsion}

We can now generalize the notion of algebraic overtwistedness.  Notice
that since every term in $\mathbf{D}_\SFT$ is a differential operator of
order at least~$1$, 
$$
\mathbf{D}_\SFT \mathbf{f} = 0 \quad \text{ for all } \quad
\mathbf{f} \in R[[\hbar]],
$$
hence every element of the extended coefficient ring $R[[\hbar]]$ represents
an element of $H_*^\SFT(M,\xi;R)$ that may or may not be trivial.
Since $\mathbf{D}_\SFT$ commutes with all elements of $R[[\hbar]]$, the subset
consisting of elements that are trivial in homology forms an ideal.
The following definition originates in \cite{LatschevWendl}.

\begin{defn}
\label{defn:torsion}
We say that a closed contact manifold $(M,\xi)$ has \defin{algebraic torsion
of order~$k$} (or \emph{$k$-torsion} for short) with coefficients in~$R$ if
$$
[\hbar^k] = 0 \in H_*^\SFT(M,\xi;R).
$$
The numerical invariant 
$$
\AT(M,\xi;R) \in \NN \cup \{0,\infty\}
$$ 
is defined to be the smallest integer $k$ such that $(M,\xi)$ has algebraic 
$k$-torsion but no $(k-1)$-torsion, or $\infty$ if there is no algebraic 
torsion of any order.
\end{defn}

Several consequences of algebraic torsion can be read off quickly from the
properties of SFT cobordism maps.  Consider first the case of trivial
coefficients $R = \QQ$, which we shall refer to as \defin{untwisted}
algebraic torsion and abbreviate
$$
\AT(M,\xi) := \AT(M,\xi;\QQ).
$$
If $(W,\omega)$ is a strong filling of $(M,\xi)$,
then the hypotheses of Theorem~\ref{thm:weakFillingMap} are fulfilled even
with $G_+ = H_2(M)$ since $\omega$ is exact at the boundary, thus we obtain
a $\QQ[[\hbar]]$-linear map $H_*^\SFT(M,\xi) \to \overline{R}[[\hbar]]$,
with $\overline{R}$ denoting the Novikov completion of
$\QQ[H_2(W) / \ker[\omega]]$.  If $[\hbar^k] = 0 \in H_*^\SFT(M,\xi)$,
then the cobordism map implies a contradiction since $\hbar^k$ does not
equal $0$ in $\overline{R}[[\hbar]]$.  Similarly, if $(W,d\lambda)$ is an
exact cobordism from $(M_-,\xi_-)$ to $(M_+,\xi_+)$, then the cobordism map
$H_*^\SFT(M_+,\xi_+) \to H_*^\SFT(M_-,\xi_-)$ of Theorem~\ref{thm:cobordismMap}
is also $\QQ[[\hbar]]$-linear, and thus any algebraic $k$-torsion in
$(M_+,\xi_+)$ is inherited by $(M_-,\xi_-)$.  This proves:

\begin{thm}
Contact manifolds with $\AT(M,\xi) < \infty$ are not strongly fillable.  
Moreover, if there exists an exact symplectic cobordism from $(M_-,\xi_-)$
to $(M_+,\xi_+)$, then $\AT(M_-,\xi_-) \le \AT(M_+,\xi_+)$.
\qed
\end{thm}

It is known (see \cite{Wendl:cobordisms}) that the second part of the above 
theorem does not hold for strong symplectic cobordisms in general, so
exactness of cobordisms is a meaningful symplectic topological condition, 
not just a technical hypothesis.  It is also known thanks to a construction of Ghiggini
\cite{Ghiggini:strongNotStein}
that strong and exact fillability are not equivalent conditions, but
Ghiggini's proof of this uses Heegaard Floer homology; thus far it is not
known whether this phenomenon can be detected via SFT or other holomorphic
curve techniques.

There are also many known examples of contact manifolds that have untwisted
algebraic torsion but are weakly fillable.  The simplest are the tight
tori $(\TT^3,\xi_k)$ for $k \ge 2$, for which weak fillings were first
constructed by Giroux \cite{Giroux:plusOuMoins}, but Eliashberg \cite{Eliashberg:fillableTorus}
showed that strong fillings do not exist, and we will see in Lecture~\ref{lec:torsion}
that $\AT(\TT^3,\xi_k) = 1$.  The weak/strong distinction can
often be detected via the choice of coefficients in SFT.
We saw in \S\ref{sec:weakStrong} that a weak filling of a contact manifold
$(M,\xi)$ can always be deformed
so as to have stable boundary with data $(\hH = (\Omega,\alpha),J_+)$ for 
which $\alpha$ is a nondegenerate contact form and $J_+$ is $C^\infty$-close to
any given element of $\jJ(\alpha)$.  Proposition~\ref{prop:weakVanish}
showed that if $(M,\xi)$ is algebraically overtwisted, then the contact
homology for the stable Hamiltonian data $(\hH,J_+)$ can also be made to
vanish.

\begin{exercise}
\label{EX:weakVanish}
Generalize the proof of Prop.~\ref{prop:weakVanish} to show that
if $(M,\xi)$ has algebraic $k$-torsion with coefficients in~$R$, then
also $[\hbar^k] = 0 \in H_*^\SFT(M,\hH_c,J_c;R)$ for sufficiently
large $c > 0$.
\end{exercise}

It then follows using Theorem~\ref{thm:weakFillingMap} that algebraic torsion
with suitably twisted coefficients also gives an obstruction to weak filling.
Let us say that $(M,\xi)$ has \defin{fully twisted} algebraic $k$-torsion
whenever $[\hbar^k] = 0 \in H_*^\SFT(M,\xi;\QQ[H_2(M)])$.  Note that in
parallel with Remark~\ref{remark:coefs}, any nested pair of subgroups
$G \subset G' \subset H_2(M)$ gives rise to a map
$$
H_*^\SFT(M,\xi;\QQ[H_2(M) / G']) \to H_*^\SFT(M,\xi;\QQ[H_2(M) / G]),
$$
which is a morphism in the sense that it maps the unit and all powers of $\hbar$
to themselves.  This implies that $(M,\xi)$ has fully twisted $k$-torsion
if and only if it has $k$-torsion for every choice of coefficients.

\begin{thm}
\label{thm:weakFillings}
If $(M,\xi)$ is a closed contact manifold with a finite order of algebraic
torsion with coefficients in $R = \QQ[H_2(M) / G]$ for some subgroup~$G$,
then $(M,\xi)$ does not admit any weak symplectic filling $(W,\omega)$ 
for which $[\omega|_{TM}] \in H^2_\dR(M)$ is rational and annihilates
all elements of~$G$.  In particular, if $(M,\xi)$ has fully twisted algebraic
torsion of some finite order, then it is not weakly fillable.
\end{thm}

\begin{remark}
The rationality condition in Theorem~\ref{thm:weakFillings} can probably
be lifted, and is known to be unnecessary at least in dimension three.
It is clear in any case that if $(M,\xi)$ admits a weak filling $(W,\omega)$,
then one can always make a small perturbation of $\omega$ to produce a
weak filling for which $[\omega|_{TM}] \in H^2(M;\QQ)$.
\end{remark}

We will see some concrete examples of algebraic torsion computations
in Lecture~\ref{lec:torsion}.  Let us conclude this discussion for now with the observation
that algebraic torsion of order \emph{zero} is a notion we've seen before:

\begin{prop}
For any closed contact manifold $(M,\xi)$ and group ring
$R = \QQ[H_2(M) / G]$, the following conditions are equivalent:
\begin{enumerate}
\item $(M,\xi)$ has algebraic $0$-torsion (with coefficients in~$R$);
\label{item:0torsion}
\item $(M,\xi)$ is algebraically overtwisted (with coefficients in~$R$);
\label{item:OT}
\item $H_*^\SFT(M,\xi;R) = 0$.
\label{item:vanish}
\end{enumerate}
\end{prop}
\begin{proof}
It is obvious that \eqref{item:vanish} implies \eqref{item:0torsion}.
Since $\mathbf{D}_\SFT \mathbf{f} = \p_\CH \mathbf{f} + \Order(\hbar)$
for $\mathbf{f} \in \aA$, the $R[[\hbar]]$-linear map
$$
\aA[[\hbar]] \to \aA : \mathbf{F} \mapsto \mathbf{F}|_{\hbar = 0}
$$
defines a chain map $(\aA[[\hbar]],\mathbf{D}) \to (\aA,\p_\CH)$ and thus
descends to an $R[[\hbar]]$-linear map
$H_*^\SFT(M,\xi;R) \to HC_*(M,\xi;R)$.
The existence of this map proves that \eqref{item:0torsion} implies
\eqref{item:OT}.

To prove that \eqref{item:OT} implies \eqref{item:vanish}, recall first that
if there exists $\mathbf{f} \in \aA$ with $\p_\CH \mathbf{f} = 1$, then 
the fact that $HC_*(M,\xi;R) = 0$ follows easily
since for any $\mathbf{g} \in \aA$ with $\p_\CH \mathbf{g} = 0$, the
graded Leibniz rule implies
$\p_\CH(\mathbf{f} \mathbf{g}) = (\p_\CH \mathbf{f}) \mathbf{g} -
\mathbf{f} (\p_\CH \mathbf{g}) = \mathbf{g}$.
This works because $\p_\CH$ is a derivation---but $\mathbf{D}_\SFT$ is not one,
so the same trick will not quite work for~$\mathbf{D}_\SFT$.
The trick in proving $H_*^\SFT(M,\xi;R) = 0$ will be to quantify the failure of
$\mathbf{D}_\SFT$ to be a derivation.  For our purposes, it suffices to know that
\begin{equation}
\label{eqn:almostDerivation}
\mathbf{D}_\SFT(\mathbf{F}\mathbf{G}) = (\mathbf{D}_\SFT \mathbf{F}) \mathbf{G}
+ (-1)^{|\mathbf{F}|} \mathbf{F} (\mathbf{D}_\SFT \mathbf{G}) + \Order(\hbar)
\end{equation}
holds for all $\mathbf{F},\mathbf{G} \in \aA[[\hbar]]$, which follows from
the fact that $\p_\CH$ is a derivation.

With this remark out of the way, suppose $\mathbf{f} \in \aA$ satisfies
$\p_\CH \mathbf{f} = 1$, in which case
\begin{equation}
\label{eqn:almostExact}
\mathbf{D}_\SFT \mathbf{f} = 1 + \hbar \mathbf{G}
\end{equation}
for some $\mathbf{G} \in \aA[[\hbar]]$.  We claim then that for any
$\mathbf{Q} \in \aA[[\hbar]]$ with $\mathbf{D}_\SFT \mathbf{Q} = 0$, there
exists $\mathbf{Q}_1 \in \aA[[\hbar]]$ with
\begin{equation}
\label{eqn:approximation}
\mathbf{D}_\SFT(\mathbf{f} \mathbf{Q}) = \mathbf{Q} + \hbar \mathbf{Q}_1
\end{equation}
and $\mathbf{D}_\SFT \mathbf{Q}_1 = 0$.  Indeed, \eqref{eqn:approximation}
follows from \eqref{eqn:almostDerivation} and \eqref{eqn:almostExact}
since $\mathbf{D}_\SFT \mathbf{Q} = 0$, and $\mathbf{D}_\SFT \mathbf{Q}_1 = 0$
then follows by applying $\mathbf{D}_\SFT$ to \eqref{eqn:approximation} and
using $\mathbf{D}_\SFT^2 = 0$.  Fixing $\mathbf{Q}_0 := \mathbf{Q} \in \aA[[\hbar]]$,
we can now define a sequence $\mathbf{Q}_k \in \aA[[\hbar]]$ satisfying
$\mathbf{D}_\SFT \mathbf{Q}_k = 0$ for all integers $k \ge 0$ via the inductive
condition
$$
\mathbf{D}_\SFT(\mathbf{f} \mathbf{Q}_k) = \mathbf{Q}_k + \hbar \mathbf{Q}_{k+1}.
$$
Then $\sum_{k=0}^\infty (-1)^k \hbar^k \mathbf{Q}_k \in \aA[[\hbar]]$, and
$$
\mathbf{D}_\SFT \left( \mathbf{f} \sum_{k=0}^\infty (-1)^k \hbar^k \mathbf{Q}_k \right)
= \mathbf{Q}.
$$
\end{proof}

\chapter{Transversality and embedding controls in dimension four}
\label{lec:automatic}

\minitoc

\vspace{12pt}

The final three lectures will be included in the published version of this book.
For updates on publication, see the author's website
\begin{center}
\url{https://www.mathematik.hu-berlin.de/~wendl/publications.html#notes}
\end{center}




\chapter{Intersection theory for punctured holomorphic curves}
\label{lec:intersections}

\minitoc

\vspace{12pt}

The final three lectures will be included in the published version of this book.
For updates on publication, see the author's website
\begin{center}
\url{https://www.mathematik.hu-berlin.de/~wendl/publications.html#notes}
\end{center}



\chapter{Torsion computations and applications}
\label{lec:torsion}

\minitoc

\vspace{12pt}

The final three lectures will be included in the published version of this book.
For updates on publication, see the author's website
\begin{center}
\url{https://www.mathematik.hu-berlin.de/~wendl/publications.html#notes}
\end{center}





\appendix

\renewcommand{\thesection}{\Alph{chapter}.\arabic{section}}
\renewcommand{\thefigure}{\Alph{chapter}.\arabic{figure}}

\chapter{Sobolev spaces}
\label{app:Sobolev}

\minitoc
\vspace{12pt}

In this appendix, we review some of the standard properties of Sobolev
spaces, in particular using them to prove Propositions~\ref{prop:BanachAlg},
\ref{prop:CkContinuity} and~\ref{prop:SobolevRescaling} 
from~\S\ref{sec:Sobolev}, and elucidating the construction of Sobolev spaces
of sections on vector bundles.  A good reference for the necessary background
material is \cite{AdamsFournier}.

\section{Approximation, extension and embedding theorems}
\label{sec:SobolevThms}

Unless otherwise noted, all functions in the following are assumed to be
defined on a nonempty open subset
$$
\uU \subset \RR^n
$$
with its standard Lebesgue measure,
and taking values in a finite-dimensional normed vector space
that will usually not need to be specified, though occasionally
we will assume it is $\RR$ or $\CC$ so that one can define products of
functions.  The domain $\uU$ will also sometimes have additional conditions
specified such as boundedness or regularity at the boundary, though we
will try not to add too many more restrictions than are really needed.  
The most useful assumption to impose on $\uU$ is known as the
\defin{strong local Lipschitz condition}: if $\uU$ is bounded, then it
means simply that near every boundary point of~$\uU$, one can find
smooth local coordinates in which $\uU$ looks like the region bounded by
the graph of a Lipschitz-continuous function, and in this case we call
$\uU$ a \defin{bounded Lipschitz domain}.  If $\uU$ is unbounded,
then one needs to impose extra conditions guaranteeing e.g.~uniformity
of Lipschitz constants, and the precise definition becomes a bit
lengthy (see \cite{AdamsFournier}*{\S 4.9}).  For our purposes, all we
really need to know about the strong local Lipschitz condition is that
that it is satisfied both by bounded Lipschitz domains and
by relatively tame unbounded domains such as 
$(0,1) \times (0,\infty) \subset \RR^2$ which 
have smooth boundary with finitely many corners.
We will repeatedly need to use the generalized version of
\defin{H\"older's inequality}, which states that for any finite collection
of measurable functions $f_1,\ldots,f_m$,
\begin{equation}
\label{eqn:generalizedHoelder}
\left\| \prod_{i=1}^m |f_i| \right\|_{L^p} \le \prod_{i=1}^m \| f_i \|_{L^{p_i}} \quad
\text{ for $1 \le p \le p_1,\ldots,p_m \le \infty$ with
$\frac{1}{p} = \sum_{i=1}^m \frac{1}{p_i}$}.
\end{equation}
This is an easy corollary of the standard version,
$$
\big\| |f| \cdot |g| \big\|_{L^1} \le \| f \|_{L^p} \cdot \| g \|_{L^q} \quad
\text{ whenever $1 \le p,q \le \infty$ and $1 = \frac{1}{p} + \frac{1}{q}$}.
$$

For an integer $k \ge 0$ and real number $p \in [1,\infty]$ we define
$W^{k,p}(\uU)$ as in \S\ref{sec:Sobolev} to be the Banach space of all
$f \in L^p(\uU)$ which have weak partial derivatives $\p^\alpha f \in L^p(\uU)$
for all $|\alpha| \le k$.  For $p=2$, these spaces are also 
often denoted by
$$
H^k(\uU) := W^{k,2}(\uU),
$$
and they admit Hilbert space structures with inner product
$$
\langle f , g \rangle_{H^k} = \sum_{|\beta| \le k}
\langle \p^\alpha f, \p^\alpha g \rangle_{L^2}.
$$
We denote by
$$
W^{k,p}_0(\uU) \subset W^{k,p}(\uU), \qquad
H^k_0(\uU) \subset H^k(\uU)
$$
the closed subspaces defined as the closures of $C_0^\infty(\uU)$ with
respect to the relevant norms.  Since $C_0^\infty(\uU)$ is dense in
$L^p(\uU)$ for $1 \le p < \infty$ (see e.g.~\cite{LiebLoss}*{\S 2.19}),
there is no difference between $W^{0,p}(\uU)$ and $W^{0,p}_0(\uU)$ for
$p < \infty$, but in general $W^{k,p}_0(\uU) \ne W^{k,p}(\uU)$ for $k \ge 1$,
with a few notable exceptions such as the case $\uU = \RR^n$
(cf.~Corollary~\ref{cor:C0inftyDense} below).  Let
\begin{equation*}
\begin{split}
W^{k,p}_\loc(\uU) := \big\{ \text{functions $f$ on $\uU$}\ \big|\ 
&\text{$f \in W^{k,p}(\vV)$ for all open subsets $\vV \subset \uU$} \\
&\text{with compact closure $\overline{\vV} \subset \uU$} \big\},
\end{split}
\end{equation*}
and we say that a sequence $f_j \in W^{k,p}_\loc(\uU)$ converges in
$W^{k,p}_\loc$ to $f \in W^{k,p}_\loc(\uU)$ if the restrictions to all
precompact open subsets $\vV \subset \overline{\vV} \subset \uU$ converge
in~$W^{k,p}(\vV)$.  Recall that for $k \in \{0,1,2,\ldots,\infty\}$,
$C^k(\uU)$ denotes the space of functions on $\uU$ with continuous derivatives
up to order~$k$, while
$$
C^k(\overline{\uU}) \subset C^k(\uU)
$$ 
is the space of $f \in C^k(\uU)$ such that for all $|\alpha| \le k$,
$\p^\alpha f$ is bounded and uniformly continuous.

\begin{thm}[\cite{AdamsFournier}*{\S 3.17, 3.22}]
\label{thm:CinftyDense}
For any open subset $\uU \subset \RR^n$, and any $k \ge 0$,
$1 \le p < \infty$, the subspace
$$
C^\infty(\uU) \cap W^{k,p}(\uU) \subset W^{k,p}(\uU)
$$
is dense.  Moreover, if $\uU \subset \RR^n$ satisfies the
strong local Lipschitz condition, then the space
$$
\left\{ f \in C^\infty(\uU)\ \Big|\ 
\text{$f = \tilde{f}|_{\uU}$ for some $\tilde{f} \in C_0^\infty(\RR^n)$} \right\}
$$
is also dense in $W^{k,p}(\uU)$, so in particular,
$$
C^\infty(\overline{\uU}) \cap W^{k,p}(\uU) \subset W^{k,p}(\uU)
$$
is dense.  \qed
\end{thm}
\begin{cor}
\label{cor:C0inftyDense}
The space $C_0^\infty(\RR^n)$ is dense in $W^{k,p}(\RR^n)$ for every
$k \ge 0$ and $p \in [1,\infty)$.   \qed
\end{cor}

Here is another useful characterization of $W^{k,p}_0(\uU)$:

\begin{thm}[\cite{AdamsFournier}*{\S 5.29}]
\label{thm:Wkp0}
Assume $\uU \subset \RR^n$ is an open subset satisfying the
strong local Lipschitz condition.  Then a function $f \in W^{k,p}(\uU)$
belongs to $W^{k,p}_0(\uU)$ if and only if the function
$\tilde{f}$ on $\RR^n$ defined to match $f$ on $\uU$ and $0$ everywhere else
belongs to $W^{k,p}(\RR^n)$.  \qed
\end{thm}

While it is obvious from the definitions that functions in $W^{k,p}_0(\uU)$
always admit extensions of class $W^{k,p}$ over~$\RR^n$, this is much
less obvious for functions in $W^{k,p}(\uU)$ in general, and it is not
true without sufficient assumptions about the regularity of~$\p \uU$.
For our purposes it suffices to consider the following case.

\begin{thm}[\cite{AdamsFournier}*{\S 5.22}]
\label{thm:SobolevExtension}
Assume $\uU \subset \RR^n$ is a bounded open subset such that
$\p\overline{\uU}$ is a submanifold of class~$C^m$ for some
$m \in \{1,2,3,\ldots,\infty\}$.  Then there exists a linear operator
$E$ that maps functions defined almost everywhere on $\uU$ to functions
defined almost everywhere on $\RR^n$ and has the following properties:
\begin{itemize}
\item For every function $f$ on~$\uU$, $Ef|_{\uU} \equiv f$ almost everywhere;
\item For every nonnegative integer $k \le m$ and every $p \in [1,\infty)$,
$E$ defines a bounded linear operator $W^{k,p}(\uU) \to W^{k,p}(\RR^n)$.
\end{itemize}  \qed
\end{thm}

\begin{cor}
\label{cor:compactClosureExtension}
Suppose $\uU , \uU' \subset \RR^n$ are open subsets such that
$\uU$ has compact closure contained in~$\uU'$.  If $\uU$
satisfies the hypothesis of Theorem~\ref{thm:SobolevExtension},
then the resulting extension operator $E$ can be chosen such that it
maps each $W^{k,p}(\uU)$ for $k \le m$ and $1 \le p < \infty$ into
$W^{k,p}_0(\uU')$.
\end{cor}
\begin{proof}
Choose a smooth function $\rho : \uU' \to [0,1]$ that has compact support
and equals~$1$ on~$\overline{\uU}$, then replace the operator $E$
given by Theorem~\ref{thm:SobolevExtension} with the operator
$f \mapsto \rho \cdot Ef$.
\end{proof}

To state the Sobolev embedding theorem in its proper generality, recall
that for $0 < \alpha \le 1$, the \defin{H\"older seminorm}
of a function $f$ on $\uU$ is defined by
$$
|f|_{C^\alpha} := |f|_{C^\alpha(\uU)} := \sup_{x \ne y \in \uU} \frac{|f(x) - f(y)|}{|x-y|^\alpha},
$$
and $C^{k,\alpha}(\uU)$ is then defined as the Banach space of
functions $f \in C^k(\overline{\uU})$ for which the norm
$$
\| f \|_{C^{k,\alpha}} := \| f \|_{C^k} +
\max_{|\beta|=k} | \p^\beta f |_{C^\alpha}
$$
is finite.  In reading the following statement, it is important to remember
that elements of $W^{k,p}(\uU)$ are technically not functions, but
rather \emph{equivalence classes} of functions defined almost everywhere.
Thus when we say e.g.~that there is an inclusion
$W^{k,p}(\uU) \hookrightarrow C^{m,\alpha}(\uU)$, the literal meaning is
that for every function $f$ representing an element of $W^{k,p}(\uU)$, 
one can change the values of $f$ in a unique way on some set of measure
zero in $\uU$ so that after this change, $f \in C^{m,\alpha}(\uU)$.
Continuity of the inclusion means that there is a bound of the form
$$
\| f \|_{C^{m,\alpha}} \le c \| f \|_{W^{k,p}}
$$
for all $f \in W^{k,p}(\uU)$,
where $c > 0$ is a constant which may in general depend on $m$, $\alpha$,
$k$, $p$ and~$\uU$, but not on~$f$.

\begin{thm}[\cite{AdamsFournier}*{\S 4.12}]
\label{thm:SobolevEmbedding}
Assume $\uU \subset \RR^n$ is an open subset satisfying the strong local
Lipschitz condition, $k \ge 1$ is an integer and $1 \le p < \infty$.
\begin{enumerate}
\item If $kp > n$ and $k - n/p < 1$, then there exist continuous 
inclusions
\begin{equation*}
\begin{split}
W^{k,p}(\uU) &\hookrightarrow C^{0,\alpha}(\uU) \quad
\text{ for each $\alpha \in (0,k - n/p]$},\\
W^{k,p}(\uU) &\hookrightarrow L^q(\uU) \quad
\text{ for each $q \in [p,\infty]$}.
\end{split}
\end{equation*}
\item If $kp < n$ and $p^* > p$ is defined by the condition
$$
\frac{1}{p^*} = \frac{1}{p} - \frac{k}{n},
$$
then there exist continuous inclusions
$$
W^{k,p}(\uU) \hookrightarrow L^q(\uU), \qquad
\text{ for each $q \in [p,p^*]$}.
$$
\item If $kp = n$, then there exist continuous inclusions
$$
W^{k,p}(\uU) \hookrightarrow L^q(\uU), \qquad
\text{ for each $q \in [p,\infty)$}.
$$
\end{enumerate}
Moreover, the spaces $W^{k,p}_0(\uU)$ admit similar inclusions
under no assumption on the open subset $\uU \subset \RR^n$.
\qed
\end{thm}
Under the same assumption on the domain~$\uU$, one can apply
Theorem~\ref{thm:SobolevEmbedding} to successive derivatives of
functions in $W^{k,p}(\uU)$ and thus obtain the following inclusions
for any integer $d \ge 0$:
\begin{equation}
\label{eqn:SobolevGood}
W^{k+d,p}(\uU) \hookrightarrow C^{d,\alpha}(\uU) \quad
\text{ if $kp > n$ and $0 < \alpha \le k - n/p < 1$},
\end{equation}
\begin{equation}
\label{eqn:SobolevGood2}
W^{k+d,p}(\uU) \hookrightarrow W^{d,q}(\uU) \quad
\text{ if $kp > n$ and $p \le q \le \infty$},
\end{equation}
\begin{equation}
\label{eqn:SobolevBad}
W^{k+d,p}(\uU) \hookrightarrow W^{d,q}(\uU) \quad
\text{ if $kp < n$ and $p \le q \le p^*$, with $\frac{1}{p^*} = \frac{1}{p} - \frac{k}{n}$},
\end{equation}
\begin{equation}
\label{eqn:SobolevBorderline1}
W^{k+d,p}(\uU) \hookrightarrow W^{d,q}(\uU) \quad
\text{ if $kp = n$ and $p \le q < \infty$}.
\end{equation}
This last inclusion can then be composed with \eqref{eqn:SobolevGood} for an 
arbitrarily large choice of~$q$, giving another inclusion
\begin{equation}
\label{eqn:SobolevBorderline2}
W^{k+d,p}(\uU) \hookrightarrow C^{d-1,\alpha}(\uU) \quad
\text{ if $kp = n$ and $0 < \alpha < 1$}.
\end{equation}

\begin{remark}
\label{remark:howManyDerivatives}
The embedding theorem suggests that one should intuitively think of
$W^{k,p}(\uU)$ as consisting of functions with ``$k - n/p$ continuous
derivatives,'' where the number $k - n/p$ may in general be a non-integer
and/or negative.  This provides a useful mnemonic for results about
embeddings of one Sobolev space into another, such as the following.
\end{remark}

\begin{cor}
\label{cor:kpmq}
Assume $\uU \subset \RR^n$ is an open subset satisfying the strong local
Lipschitz condition, $1 \le p,q < \infty$, and $k, m \ge 0$ are integers 
satisfying
$$
k \ge m, \qquad p \le q, \qquad\text{ and }\qquad
k - \frac{n}{p} \ge m - \frac{n}{q}.
$$
Then there exists a continuous inclusion
$W^{k,p}(\uU) \hookrightarrow W^{m,q}(\uU)$.  \qed
\end{cor}

By the Arzel\`a-Ascoli theorem, the natural inclusion
$$
C^{k,\alpha'}(\uU) \hookrightarrow C^{k,\alpha}(\uU)
$$
for $\alpha < \alpha'$ is a compact operator whenever $\uU \subset \RR^n$
is bounded.  It follows
that if $\uU \subset \RR^n$ in \eqref{eqn:SobolevGood} is bounded and
$\alpha$ is \emph{strictly} less than the extremal value
$k - n/p$, then the inclusion \eqref{eqn:SobolevGood} is also compact.  
A similar statement holds for the inclusion \eqref{eqn:SobolevBad} when
$p \le q < p^*$, and this is known as the \defin{Rellich-Kondrachov
compactness theorem}.  We summarize these as follows:

\begin{thm}[\cite{AdamsFournier}*{\S 6.3}]
\label{thm:Rellich}
Assume $\uU \subset \RR^n$ is a bounded Lipschitz domain,
$k \ge 1$ and $d \ge 0$ are integers and $1 \le p < \infty$.
\begin{enumerate}
\item
If $kp > n$ and $k - n/p < 1$, then the inclusions
\begin{equation*}
\begin{split}
W^{k+d,p}(\uU) &\hookrightarrow C^{d,\alpha}(\uU) \quad
\text{ for $\alpha \in (0,k - n/p)$},\\
W^{k+d,p}(\uU) &\hookrightarrow W^{d,q}(\uU) \quad
\text{ for $q \in [p,\infty)$}
\end{split}
\end{equation*}
are compact.
\item
If $kp \le n$ and $p^* \in (p,\infty]$ is defined by the condition
$1/p^* = 1/p - k/n$, then the inclusions
$$
W^{k+d,p}(\uU) \hookrightarrow W^{d,q}(\uU) \quad
\text{ for $q \in [p,p^*)$}
$$
are compact.
\end{enumerate}  
In particular, the continuous inclusion $W^{k,p}(\uU) \hookrightarrow
W^{m,q}(\uU)$ in Corollary~\ref{cor:kpmq} is compact whenever the
inequality $k - n/p \ge m - n/q$ is strict.  \qed
\end{thm}

\section{Products, compositions, and rescaling}
\label{sec:SobolevProducts}

We now restate and prove Propositions~\ref{prop:BanachAlg},
\ref{prop:CkContinuity} and~\ref{prop:SobolevRescaling} 
from~\S\ref{sec:Sobolev}.  These are all corollaries of
the Sobolev embedding theorem, so in particular they hold for the same
class of domains $\uU \subset \RR^n$, and the restrictions on $\uU$
can be dropped at the cost of replacing each space $W^{k,p}$ by~$W^{k,p}_0$.

We begin by generalizing Prop.~\ref{prop:BanachAlg}, hence we consider
Sobolev spaces of functions valued in $\RR$ or~$\CC$
so that pointwise products of functions are well defined almost everywhere.
We say that there is a \defin{continuous product map},
$$
W^{k_1,p_1}(\uU) \times \ldots \times W^{k_m,p_m}(\uU) \to W^{k,p}(\uU),
$$
or a continuous product \defin{pairing} in the case $m=2$,
if for every set of functions $f_i \in W^{k_i,p_i}(\uU)$ with $i=1,\ldots,m$, 
the pointwise product function
$f_1 \cdot \ldots \cdot f_m$ is in $W^{k,p}(\uU)$ and there is an estimate of the form
$$
\| f_1 \cdot \ldots \cdot f_m \|_{W^{k,p}} \le c \| f_1 \|_{W^{k_1,p_1}} \cdot
\ldots \cdot \| f_m \|_{W^{k_m,p_m}}
$$
for some constant $c > 0$ not depending on $f_1,\ldots,f_m$.  The case
$m=2$, $k_1 = k_2 = k$ and $p_1 = p_2 = p$ is especially interesting, as
the space $W^{k,p}(\uU)$ is then a \defin{Banach algebra}.  More generally,
one can ask under what circumstances multiplication by functions of class
$W^{k,p}$ defines a bounded linear operator on functions of class~$W^{m,q}$.
A hint about this comes from the world of classically differentiable
functions: multiplication by $C^k$-smooth functions defines a continuous 
map $C^m \to C^m$ if and only if $k \ge m$.  The corresponding answer in
Sobolev spaces turns out to be that functions of class $W^{k,p}$ need to have
strictly more than zero derivatives in the sense of 
Remark~\ref{remark:howManyDerivatives}, and at least as many derivatives
as functions of class~$W^{m,q}$.

\begin{thm}
\label{thm:productEstimates}
Assume $\uU \subset \RR^n$ is an open subset satisfying the strong local 
Lipschitz condition, $k$, $p$, $m$ and $q$ satisfy the same numerical 
hypotheses as in
Corollary~\ref{cor:kpmq} (so in particular $W^{k,p}(\uU)$ embeds continuously
into $W^{m,q}(\uU)$), and $kp > n$.  Then there exists a continuous product pairing
$$
W^{k,p}(\uU,\CC) \times W^{m,q}(\uU,\CC) \to W^{m,q}(\uU,\CC) : (f,g) \mapsto fg.
$$
\end{thm}

The following preparatory lemma will be useful both for proving the product 
estimate and for further results below.
It is an easy consequence of Theorem~\ref{thm:SobolevEmbedding}
and H\"older's inequality.

\begin{lemma}
\label{lemma:SobolevProduct}
Assume $\uU \subset \RR^n$ is an open subset satisfying the strong local
Lipschitz condition, $m \ge 2$ is an integer, and we are given
positive numbers $p_1,\ldots,p_m \ge 1$ and integers
$k_1,\ldots,k_m \ge 0$.  Let $I := \left\{ i \in \{1,\ldots,m\}\ \big|\ 
k_i p_i \le n \right\}$.  Then for any $q \ge 1$ satisfying
$$
\sum_{i \in I} \left( \frac{1}{p_i} - \frac{k_i}{n} \right) < \frac{1}{q} 
\le \sum_{i=1}^m \frac{1}{p_i},
$$
there is a continuous product map
$$
W^{k_1,p_1}(\uU) \times \ldots \times W^{k_m,p_m}(\uU) \to L^q(\uU).
$$
\end{lemma}
\begin{proof}
By the generalized H\"older inequality \eqref{eqn:generalizedHoelder},
it suffices to show that for any $q \ge 1$ in the stated range,
one can find numbers $q_1,\ldots,q_m \in [q,\infty]$ satisfying
$1/q = 1/q_1 + \ldots + 1/q_m$ for which
Theorem~\ref{thm:SobolevEmbedding} provides continuous inclusions
$$
W^{k_i,p_i}(\uU) \hookrightarrow L^{q_i}(\uU)
$$
for each $i=1,\ldots,m$.  Whenever $k_i p_i > n$,
this inclusion is valid with $q_i$ chosen freely from the 
interval $[p_i,\infty]$, so $1/q_i$ can then take any value subject to the 
constraint
$$
0 \le \frac{1}{q_i} \le \frac{1}{p_i}.
$$
If on the other hand $k_i p_i \le n$, then we can arrange $1/q_i$ to take
any value in the range
$$
\frac{1}{p_i} - \frac{k_i}{n} < \frac{1}{q_i} \le \frac{1}{p_i}.
$$
Adding these up, the range of values for 
$\sum_i \frac{1}{q_i}$ that we can achieve in this way covers the
stated interval.
\end{proof}

\begin{proof}[Proof of Theorem~\ref{thm:productEstimates}]
By density of smooth functions, it suffices to prove that an estimate of the 
form
$$
\| f g \|_{W^{m,q}} \le c \| f \|_{W^{k,p}} \| g \|_{W^{m,q}}
$$
holds for all $f \in C^\infty(\uU) \cap W^{k,p}(\uU)$ and
$g \in C^\infty(\uU) \cap W^{m,q}(\uU)$.  Equivalently, we need to show that for
all $f$ and $g$ of this type and every multiindex $\alpha$ of degree
$|\alpha| \le m$, there is a constant $c > 0$ independent of $f$ and $g$
such that
$$
\| \p^\alpha(f g) \|_{L^q} \le c \| f \|_{W^{k,p}} \| g \|_{W^{m,q}}.
$$
Since $f$ and $g$ are smooth, we are
free to use the product rule in computing $\p^\alpha(fg)$, which will then
be a linear combination of terms of the form $\p^\beta f \cdot \p^\gamma g$
where $|\alpha| = |\beta| + |\gamma|$, hence we have reduced the problem
to proving a bound
$$
\| \p^\beta f \cdot \p^\gamma g \|_{L^q} \le c \| f \|_{W^{k,p}} \| g \|_{W^{m,q}}
$$
for every pair of multiindices $\beta$, $\gamma$ with $|\beta| + |\gamma| \le m$.
Since $\p^\beta f \in W^{k - |\beta|,p}(\uU)$ and
$\p^\gamma f \in W^{m - |\gamma|,q}(\uU)$, the result follows if we can
assume that for every pair of integers $a , b \ge 0$ satisfying
$a + b \le m$, there exists a continuous product pairing
\begin{equation}
\label{eqn:desiredPairing}
W^{k-a,p}(\uU) \times W^{m-b,q}(\uU) \to L^q(\uU).
\end{equation}
If $(k-a) p > n$, then $W^{k-a,p} \hookrightarrow L^\infty$ and
\eqref{eqn:desiredPairing} is immediate since $W^{m-b,q} \hookrightarrow
L^q(\uU)$.  For the remaining cases, we shall apply
Lemma~\ref{lemma:SobolevProduct}, noting that the condition
$1/q \le 1/p + 1/q$ is trivially satisfied.  

If $(m-b)q > n$ but
$(k-a)p \le n$, then the hypotheses of the lemma are satisfied if and
only if
$$
\frac{1}{p} - \frac{k-a}{n} < \frac{1}{q}.
$$
Since $\frac{1}{p} - \frac{k}{n} \le \frac{1}{q} - \frac{m}{n}$ by assumption,
we have
$$
\frac{1}{p} - \frac{k-a}{n} = \frac{1}{p} - \frac{k}{n} + \frac{a}{n}
\le \frac{1}{q} - \frac{m}{n} + \frac{a}{n} \le \frac{1}{q}
$$
since $a \le m$, and equality holds only if $a=m$, $b=0$ and
$k - n/p = m - n/q$, which implies $mq > n$.  In this case
$W^{m-b,q} = W^{m,q} \hookrightarrow L^\infty$, and the pairing
\eqref{eqn:desiredPairing} follows because $W^{k-a,p} = W^{k-m,p}$ embeds
continuously into~$L^q$: the latter follows from Theorem~\ref{thm:SobolevEmbedding}
since $\frac{1}{p} - \frac{k-m}{n} = \frac{1}{q}$.

Finally, when $(k-a)p \le n$ and $(m-b)q \le n$, the hypotheses of the lemma
are satisfied since
$$
\left(\frac{1}{p} - \frac{k-a}{n}\right) + \left(\frac{1}{q} - \frac{m-b}{n}\right)
\le \frac{1}{p} - \frac{k}{n} + \frac{1}{q} - \frac{m}{n} + \frac{m}{n} =
\left( \frac{1}{p} - \frac{k}{n} \right) + \frac{1}{q} < \frac{1}{q},
$$
where we've used the assumption $kp > n$ and the fact that $a + b \le m$.
\end{proof}

The next result generalizes Proposition~\ref{prop:CkContinuity} and
concerns the following question: if $f : \uU \to \RR^m$
is a function of class $W^{k,p}$ whose graph lies in some open subset
$\vV \subset \uU \times \RR^m$, and $\Psi : \vV \to \RR^N$ is another
function, under what conditions can we conclude that the function
$$
\uU \to \RR^N : x \mapsto \Psi(x,f(x))
$$
is in $W^{k,p}(\uU,\RR^N)$?  We will abbreviate this function in the following
by $\Psi \circ (\Id \times f)$, and we would also like to know whether it
depends continuously (in the $W^{k,p}$-topology) on $f$ and~$\Psi$.
The following theorem is stated rather generally, but on first reading
you may prefer to assume $\uU \subset \RR^n$ is bounded, in which case
some of the hypotheses become vacuous.  We will
say that an open subset $\vV \subset \uU \times \RR^m$ is a
\defin{star-shaped neighborhood of~$f : \uU \to \RR^m$} if it contains the 
graph of $f_0$ and
$$
\text{$(x,v) \in \vV$} \quad \Rightarrow \quad \text{$(x, t v + (1-t) f_0(x)) \in \vV$ for all
$t \in [0,1]$}.
$$

\begin{thm}
\label{thm:precomposition}
Assume $\uU \subset \RR^n$ is an open subset satisfying the strong local
Lipschitz condition, $p \in [1,\infty)$ and $k \in \NN$ satisfy $kp > n$,
and $\vV \subset \uU \times \RR^m$ is a
star-shaped neighborhood of some function
$f_0 \in W^{k,p}(\uU,\RR^m)$.  Assume also $\oO^{k,p}(\uU ; \vV) \subset
W^{k,p}(\uU,\RR^m)$ is an open neighborhood of $f_0$ such that
$$
(x,f(x)) \in \vV \quad \text{ for all $x \in \uU$ and $f \in \oO^{k,p}(\uU ; \vV)$},
$$
and $\oO^k(\overline{\vV},\RR^N) \subset C^k(\overline{\vV},\RR^N)$ 
is a subset such that all $\Psi \in \oO^k(\overline{\vV},\RR^n)$ have 
the following properties:\footnote{Both of the conditions on 
$\Psi \in \oO^k(\overline{\vV},\RR^n)$ are vacuous if $\uU \subset \RR^n$
is bounded.}
\begin{enumerate}
\item
There exists a bounded subset $\kK \subset \uU$ such that
$\Psi(x,v)$ is independent of $x$ for all $x \in \uU \setminus \kK$;
\item
$\Psi \circ (\Id \times f_0) \in L^p(\uU,\RR^N)$.
\end{enumerate}
Then there is a well-defined and continuous map
$$
\oO^k(\overline{\vV},\RR^N) \times \oO^{k,p}(\uU;\vV) \to W^{k,p}(\uU,\RR^N) :
(\Psi,f) \mapsto \Psi \circ (\Id \times f).
$$
\end{thm}
\begin{proof}
We will show first that if $f \in \oO^{k,p}(\uU;\vV)$ is smooth, then
$\Psi \circ (\Id \times f)$ belongs to $W^{k,p}(\uU,\RR^N)$ for every 
$\Psi \in \oO^k(\overline{\vV},\RR^N)$.  Since $\vV$ is a star-shaped
neighborhood of~$f_0$, we have
\begin{equation*}
\begin{split}
|\Psi(x,f(x)) - &\Psi(x,f_0(x))| = \left| \int_0^1 
\frac{d}{dt} \Psi\big(x,t f(x) + (1-t) f_0(x)\big) \, dt \right| \\
&\le \left( \int_0^1 | D_2 \Psi\big(x, t f(x) + (1-t) f_0(x)\big) | \, dt \right) \cdot 
| f(x) - f_0(x) | \\
&\le \| \Psi \|_{C^1(\vV)} \cdot | f(x) - f_0(x) |
\end{split}
\end{equation*}
for all $x \in \uU$, implying
$$
\| \Psi \circ (\Id \times f) - \Psi \circ (\Id \times f_0) \|_{L^p} 
\le \| \Psi \|_{C^1(\vV)} \cdot \| f - f_0 \|_{L^p},
$$
hence $\Psi \circ (\Id \times f) \in L^p(\uU,\RR^N)$.

For $\ell=1,\ldots,k$, we can regard the $\ell$th derivative of $\Psi$ 
with respect to variables in $\RR^m$ as a 
bounded and uniformly continuous map from $\vV$ into the
vector space of symmetric $\ell$-multilinear maps from $\RR^m$ to $\RR^N$,
denoting this by
$$
D_2^\ell \Psi : \vV \to \Hom((\RR^m)^{\otimes \ell},\RR^N).
$$
Denote the partial derivatives with respect to variables in
$\uU \subset \RR^n$ by
$$
D^\beta_1 \Psi : \vV \to \RR^N,
$$
where $\beta$ is a multiindex in $n$ variables.
Now for any multiindex $\alpha$ with $|\alpha| \le k$, the derivative
$\p^\alpha(\Psi \circ (\Id \times f))$ is a linear combination of product
functions of the form
\begin{equation}
\label{eqn:derivProduct}
(D_1^\gamma D_2^\ell \Psi \circ (\Id \times f))(\p^{\beta_1} f,\ldots,\p^{\beta_\ell} f)
 : \uU \to \RR^N,
\end{equation}
where $\ell + |\gamma| \in \{1,\ldots,|\alpha|\}$ and $|\beta_1| + \ldots + 
|\beta_\ell| = |\alpha| - |\gamma|$.  If $\ell=0$ but $|\gamma| > 0$, then
this expression is clearly in $L^p(\uU,\RR^N)$ since it is continuous and
$D^\gamma_1 \Psi(x,v) = 0$ for $x \in \uU \setminus \kK$, where
$\kK$ is bounded.  For $\ell \ge 1$, it satisfies
$$
\left\| (D_1^\gamma D_2^\ell \Psi \circ (\Id \times f))(\p^{\beta_1} f,\ldots,\p^{\beta_\ell} f) \right\|_{L^p(\uU)} \le
\| D_1^\gamma D_2^\ell \Psi \|_{C^0(\vV)} \cdot \left\| \prod_{j=1}^\ell 
|\p^{\beta_j} f| \right\|_{L^p(\uU)}
$$
if the product on the right hand side has finite $L^p$-norm.  The latter is 
trivially true if $\ell=1$. 
To deal with the $\ell \ge 2$ case, note that
$\p^{\beta_j} f \in W^{k-|\beta_j|,p}(\uU)$ for each $j=1,\ldots,\ell$,
so the necessary bound 
will follow from the existence of a continuous product map
$$
W^{k-m_1,p}(\uU) \times \ldots \times W^{k-m_\ell,p}(\uU) \to L^p(\uU)
$$
for $m_j := |\beta_j|$, and we claim that such a product map does exist
whenever $kp > n$ and $m_1,\ldots,m_\ell \ge 0$ are integers satisfying
$m_1 + \ldots + m_\ell \le k$.  To see this, note
first that since $W^{k-m_j,p} \hookrightarrow L^\infty$ whenever
$(k-m_j)p > n$, it suffices to prove the claim under the assumption that
$(k-m_j)p \le n$ for every $j=1,\ldots,\ell$.  In this case,
Lemma~\ref{lemma:SobolevProduct} provides the desired product map if
the condition
$$
\sum_{j=1}^\ell \left( \frac{1}{p} - \frac{k - m_j}{n} \right) <
\frac{1}{p} \le \sum_{j=1}^\ell \frac{1}{p}
$$
is satisfied.  And it is: using $kp > n$, $\ell \ge 2$ 
and $m_1 + \ldots + m_\ell \le k$, we find
\begin{equation*}
\begin{split}
\sum_{j=1}^\ell \left( \frac{1}{p} - \frac{k - m_j}{n} \right) 
&= \ell \left( \frac{1}{p} - \frac{k}{n} \right) + \frac{m_1 + \ldots
+ m_\ell}{n} \\
&\le \frac{1}{p} + (\ell - 1) \left(\frac{1}{p} - \frac{k}{n} \right)
< \frac{1}{p}.
\end{split}
\end{equation*}
This proves that $\Psi \circ (\Id \times f) \in W^{k,p}(\uU,\RR^N)$.

Next, suppose $f \in \oO^{k,p}(\uU;\vV)$ is not necessarily smooth but
$f_i \in \oO^{k,p}(\uU;\vV)$ is a sequence of smooth functions converging
to $f$ in~$W^{k,p}$, while $\Psi_i \in \oO^k(\overline{\vV},\RR^N)$
converges to $\Psi \in \oO^k(\overline{\vV},\RR^N)$ in~$C^k$.  
Then the same argument we used to estimate
$\| \Psi \circ (\Id \times f) - \Psi \circ (\Id \times f_0) \|_{L^p}$
shows that $\Psi \circ (\Id \times f_i) \to \Psi \circ (\Id \times f)$
in~$L^p$, and since $f_i$ is also $C^0$-convergent, 
the compactly supported functions
$D^\gamma_1 \Psi \circ (\Id \times f_i)$ converge to
$D^\gamma_1 \Psi \circ (\Id \times f)$ in $L^p$ for each multiindex
with $1 \le |\gamma| \le k$.  For $\ell \ge 1$ and $|\gamma| + \ell \le k$,
$D_1^\gamma D_2^\ell \Psi_i \circ (\Id \times f_i)$
converges to $D_1^\gamma D_2^\ell \Psi \circ (\Id \times f)$ in 
$C^0(\overline{\uU},\RR^N)$,
and each of the derivatives $\p^{\beta_j} f_i$ appearing in
\eqref{eqn:derivProduct} also converges in~$L^p(\uU)$.  In light of the
continuous product maps discussed above, it follows that
each derivative $\p^\alpha (\Psi_i \circ (\Id \times f_i))$ for $|\alpha| \le k$ is
$L^p$-convergent, and its limit is necessarily the corresponding weak
derivative $\p^\alpha (\Psi \circ (\Id \times f))$, hence (see Exercise~\ref{EX:derivWeakly} 
below) $\Psi \circ (\Id \times f) \in W^{k,p}(\uU,\RR^N)$ and $\Psi_i 
\circ (\Id \times f_i) \stackrel{W^{k,p}}{\longrightarrow} \Psi \circ
(\Id \times f)$.
\end{proof}

\begin{exercise}
\label{EX:derivWeakly}
Show that if $f_i$ is a sequence of smooth functions on an open set
$\uU \subset \RR^n$ with $f_i \stackrel{L^p}{\to} f$ and
$\p^\alpha f_i \stackrel{L^p}{\to} g$ for some multiindex $\alpha$ and
functions $f , g \in L^p(\uU)$, then $\p^\alpha f = g$ in the sense of
distributions.
\end{exercise}

The following result on coordinate transformations of the domain can be
proved in an analogous way to Theorem~\ref{thm:precomposition}, though it
is considerably easier since there is no need to worry about Sobolev
product maps (and thus no need to assume $kp > n$ or impose regularity
conditions on the domain).

\begin{thm}[\cite{AdamsFournier}*{\S 3.41}]
\label{thm:SobolevCoordinate}
Assume $k \in \NN$, $1 \le p \le \infty$, and $\uU, \uU' \subset \RR^n$
are open subsets with a $C^k$-smooth diffeomorphism
$\varphi : \uU \to \uU'$ such that all derivatives of $\varphi$ and
$\varphi^{-1}$ up to order~$k$ are bounded and uniformly continuous.
Then there is a well-defined Banach space isomorphism
$$
W^{k,p}(\uU') \to W^{k,p}(\uU) : f \mapsto f \circ \varphi.
$$  \qed
\end{thm}

We now restate and prove Proposition~\ref{prop:SobolevRescaling}.
We denote by $\intDD^n, \intDD^n_\epsilon \subset \RR^n$
the open balls of radius~$1$ and $\epsilon$ respectively about the origin.

\begin{thm}
\label{thm:SobolevRescaling}
Assume $p \in [1,\infty)$ and $k \in \NN$ satisfy $kp > n$, and for each
$f \in W^{k,p}(\intDD^n)$ and $\epsilon \in (0,1]$, define
$f_\epsilon \in W^{k,p}(\intDD^n)$ by
$$
f_\epsilon(x) := f(\epsilon x).
$$
Then there exist constants $C > 0$ and $r > 0$ such that
for every $f \in W^{k,p}(\intDD^n)$,
$$
\| f_\epsilon - f(0) \|_{W^{k,p}(\intDD^n)} \le C \epsilon^r \| f - f(0) \|_{W^{k,p}(\intDD^n)}
\quad\text{ for all $\epsilon \in (0,1]$}.
$$
\end{thm}
\begin{proof}
Let $\beta$ denote a multiindex of order $|\beta| = k$.  Then using a change
of variables, we have
\begin{equation*}
\begin{split}
\| \p^\beta (f_\epsilon - f(0)) \|_{L^p(\intDD^n)}^p &= 
\epsilon^{kp} \int_{\DD^n} | \p^\beta f(\epsilon x) |^p 
= \epsilon^{kp - n} \int_{\DD^n_\epsilon} | \p^\beta f(x) |^p \\
&\le \epsilon^{kp - n} \| \p^\beta f \|_{L^p(\intDD^n)}^p
\le \epsilon^{kp - n} \| f - f(0) \|_{W^{k,p}(\intDD^n)}^p,
\end{split}
\end{equation*}
and $\epsilon^{kp - n} \to 0$ as $\epsilon \to 0$ since $kp - n > 0$.

Next, suppose $|\beta| = m \in \{1,\ldots,k-1\}$.  Then
$\p^\beta f$ and $\p^\beta f_\epsilon$ are in $W^{k-m,p}(\intDD^n)$, and if
$(k-m)p < n$, Theorem~\ref{thm:SobolevEmbedding} gives a continuous
inclusion
\begin{equation}
\label{eqn:goodInclusion}
W^{k-m,p}(\intDD^n) \hookrightarrow L^q(\intDD^n)
\end{equation}
with $q > p$ satisfying $1/q + (k-m)/n = 1/p$.  Likewise, if $(k-m)p \ge n$,
then \eqref{eqn:goodInclusion} is a continuous inclusion for arbitrarily
large choices of $q \ge p$.  We will
therefore assume in general that \eqref{eqn:goodInclusion} holds with
$q \in (p,\infty)$ satisfying
$$
\frac{1}{q} + \frac{1}{r} = \frac{1}{p},
$$
where $r = \frac{n}{k-m}$ if $(k-m)p < n$ and otherwise
$r = p + \delta$ for some $\delta > 0$ which may be chosen arbitrarily small.
Given this, we use H\"older's inequality and find
\begin{equation*}
\begin{split}
\| \p^\beta (f_\epsilon - f(0)) \|_{L^p(\intDD^n)}^p &=
\epsilon^{mp} \int_{\DD^n} | \p^\beta f(\epsilon x) |^p 
= \epsilon^{mp - n} \int_{\DD^n_\epsilon} | \p^\beta f(x) |^p \\
&\le \epsilon^{mp - n} \| \p^\beta f \|^p_{L^q(\intDD^n_\epsilon)} 
\| 1 \|^p_{L^r(\intDD^n_\epsilon)} \\
&\le \epsilon^{mp - n} \left[ \Vol(\DD^n_\epsilon) \right]^{p/r} 
\| \p^\beta f \|^p_{L^q(\intDD^n)} \\
&\le c \epsilon^{mp - n} \left[ \Vol(\DD^n_\epsilon) \right]^{p/r} 
\| \p^\beta f \|^p_{W^{k-m,p}(\intDD^n)} \\
&\le c \epsilon^{mp - n} \left[ \Vol(\DD^n_\epsilon) \right]^{p/r} 
\| f - f(0) \|^p_{W^{k,p}(\intDD^n)}
\end{split}
\end{equation*}
for some constant $c > 0$.  Writing $\Vol(\DD^n_\epsilon) = C \epsilon^n$
for a suitable constant $C > 0$, the exponent on $\epsilon$ in this
expression becomes
$$
mp - n + \frac{np}{r},
$$
which is positive whenever $r = p + \delta$ with $\delta > 0$ sufficiently
small since $m \ge 1$, and in the case $r = n/(k-m)$, it becomes simply
$kp - n > 0$.

Finally, to bound the $L^p$-norm of $f_\epsilon - f(0)$ itself, we can
use the fact that $f \in W^{k,p}$ is H\"older continuous, i.e.~it satisfies
$$
|f(x) - f(0)| \le c \| f - f(0) \|_{W^{k,p}(\intDD^n)} 
|x|^\alpha \quad \text{ for all $x \in \intDD^n$}
$$
for suitable constants $c > 0$ and $\alpha \in (0,1)$.  Thus
\begin{equation*}
\begin{split}
\| f_\epsilon - f(0) \|_{L^p(\intDD^n)}^p &= \int_{\DD^n} | f(\epsilon x) - f(0) |^p 
\le c^p \| f - f(0) \|_{W^{k,p}}^p \int_{\DD^n} | \epsilon x |^{\alpha p} \\
&= c^p \| f - f(0) \|_{W^{k,p}}^p \epsilon^{\alpha p} \int_{\DD^n} |x|^{\alpha p} \\
&= \epsilon^{\alpha p} \frac{c^p \Vol(S^{n-1})}{\alpha p + n} 
\| f - f(0) \|_{W^{k,p}}^p .
\end{split}
\end{equation*}
\end{proof}

\section{Spaces of sections of vector bundles}
\label{sec:SobolevMfd}

In this section, fix a field
$$
\FF := \RR \text{ or } \CC,
$$
assume $M$ is a smooth $n$-dimensional manifold,
possibly with boundary, and $\pi : E \to M$ is a smooth vector bundle of 
rank~$m$ over~$\FF$.  This comes with a ``bundle atlas'' 
$\aA(\pi)$, a set whose elements $\alpha \in \aA(\pi)$ 
each consist of the following data:
\begin{enumerate}
\item An open subset $\uU_\alpha \subset M$;
\item A smooth local coordinate chart $\varphi_\alpha : \uU_\alpha 
\stackrel{\cong}{\longrightarrow} \Omega_\alpha$, where $\Omega_\alpha$ is
an open subset of $\RR^n_+ := \{ (x_1,\ldots,x_n) \in \RR^n
\ |\ x_n \ge 0 \}$;
\item A smooth local trivialization $\Phi_\alpha : E|_{\uU_\alpha} 
\stackrel{\cong}{\longrightarrow} \uU_\alpha \times \FF^m$.
\end{enumerate}
Smoothness of $\varphi_\alpha$ and $\Phi_\alpha$ means as usual that 
for every pair $\alpha,\beta \in \aA(\pi)$,
the coordinate transformations 
$$
\varphi_{\beta\alpha} := \varphi_\beta^{-1} \circ \varphi_\alpha :
\Omega_{\alpha\beta} \stackrel{\cong}{\longrightarrow} \Omega_{\beta\alpha},
\qquad
\Omega_{\alpha\beta} := \varphi_\alpha(\uU_\alpha \cap \uU_\beta)
$$
and transition maps
\begin{equation*}
\begin{split}
g_{\beta\alpha} : \uU_\alpha \cap \uU_\beta \to \GL(m,\FF) \quad
\text{ such that } \quad & \text{$\Phi_\beta \circ \Phi_\alpha^{-1}(x,v) = 
(x,g_{\beta\alpha}(x) v)$} \\
&\text{for $x \in \uU_\alpha \cap \uU_\beta$, $v \in \FF^m$}
\end{split}
\end{equation*}
are smooth, and we shall assume the bundle atlas is maximal
in the sense that any triple $(\uU,\varphi,\Phi)$ that is smoothly
compatible with every $\alpha \in \aA(\pi)$ also belongs to~$\aA(\pi)$.

Any $\alpha \in \aA(\pi)$ now associates to sections $\eta : M \to E$ their
local coordinate representatives
$$
\eta^\alpha := \pr_2 \circ \Phi_\alpha \circ \eta \circ \varphi_\alpha^{-1} :
\Omega_\alpha \to \FF^m,
$$
where $\pr_2 : \uU_\alpha \times \FF^m \to \FF^m$ is the projection, 
and the representatives
with respect to two distinct $\alpha, \beta \in \aA(\pi)$ are related by
$$
\eta^\beta = (g_{\beta\alpha} \circ \varphi_\beta^{-1}) 
(\eta^\alpha \circ \varphi_{\alpha\beta}) \quad \text{ on $\Omega_{\beta\alpha}
\subset \Omega_\beta$}.
$$
For $p \in [1,\infty]$ and each integer $k \ge 0$, we then define the 
topological vector space of sections of class $W^{k,p}_\loc$ by
\begin{equation*}
\begin{split}
W^{k,p}_\loc(E) := \big\{ \eta : M \to E\ \big|\ &\text{sections such that
$\eta^\alpha \in W^{k,p}_\loc(\mathring{\Omega}_\alpha,\FF^m)$} \\
&\text{for all $\alpha \in \aA(\pi)$} \big\},
\end{split}
\end{equation*}
where convergence $\eta_j \to \eta$ in $W^{k,p}_\loc(E)$ means that 
$\eta^\alpha_j \to \eta^\alpha$ in $W^{k,p}_\loc(\mathring{\Omega}_\alpha,\FF^m)$
for all $\alpha \in \aA(\pi)$.  Note that $\Omega_\alpha$ is not necessarily
an open subset of $\RR^n$ since it may contain points in
$\p \RR^n_+ = \RR^{n-1} \times \{0\}$, but its interior 
$\mathring{\Omega}_\alpha$ is open in~$\RR^n$,
and $W^{k,p}_\loc(\mathring{\Omega}_\alpha)$ is thus defined as in
\S\ref{sec:SobolevThms}.  Strictly speaking, elements of
$\eta \in W^{k,p}_\loc(E)$ are not sections but \emph{equivalence classes}
of sections defined almost everywhere---the latter notion is defined with
respect to any measure arising from a smooth volume element on~$M$,
and it does not depend on this choice.

It turns out that $W^{k,p}_\loc(E)$ can be given the structure of a
Banach space if $M$ is compact.  This follows from the fact that $M$ can
then be covered by a finite subset of the atlas $\aA(\pi)$, but we must
be a little bit careful: not all charts in $\aA(\pi)$ are equally suitable
for defining $W^{k,p}$-norms on sections, because e.g.~even a nice smooth
section $\eta \in \Gamma(E)$ may have 
$\| \eta^\alpha \|_{W^{k,p}(\mathring{\Omega}_\alpha)} = \infty$ if
$\Omega_\alpha \subset \RR^n_+$ is unbounded.  One way to deal with this
is as follows: we will say that $\alpha \in \aA(\pi)$ is a
\defin{precompact chart}
if there exists $\alpha' \in \aA(\pi)$ and a compact subset $\kK \subset M$
such that
$$
\uU_\alpha \subset \kK \subset \uU_{\alpha'}.
$$
When this is the case, $\Omega_\alpha \subset \RR^n_+$ is necessarily
bounded, and the transition maps between two precompact charts necessarily
have bounded derivatives of all orders, as they are restrictions to
precompact subsets of maps that are smooth on larger domains.
If $M$ is compact, then one can always find a finite subset
$I \subset \aA(\pi)$ consisting of precompact charts
such that $M = \bigcup_{\alpha \in I} \uU_\alpha$.

\begin{defn}
\label{defn:SobolevMfd}
Suppose $E \to M$ is a smooth vector bundle over a compact manifold $M$,
and $I \subset \aA(\pi)$ is a finite set of precompact charts such that
$\{ \uU_\alpha \}_{\alpha \in I}$ is an open cover of~$M$.  We then
define $W^{k,p}(E)$ as the vector space 
of all sections $\eta : M \to E$ for which the norm
$$
\| \eta \|_{W^{k,p}} := \| \eta \|_{W^{k,p}(E)} :=
\sum_{\alpha \in I} \| \eta^\alpha \|_{W^{k,p}(\mathring{\Omega}_\alpha)}
$$
is finite.
\end{defn}

The norm in the above definition depends on auxiliary choices, but it is
easy to see that the resulting definition of the space $W^{k,p}(E)$ and
its topology do not.  In fact:

\begin{prop}
\label{prop:differentCharts}
If $M$ is compact, then $W^{k,p}(E) = W^{k,p}_\loc(E)$, and a sequence
$\eta_j$ converges to $\eta$ in $W^{k,p}_\loc(E)$ if and only if
the norm given in Definition~\ref{defn:SobolevMfd} satisfies
$\| \eta_j - \eta \|_{W^{k,p}(E)} \to 0$.
\end{prop}

The proposition is an immediate consequence of the following.

\begin{lemma}
\label{lemma:differentCharts}
Suppose $M$ is a smooth manifold, $\pi : E \to M$ is a smooth vector
bundle, $\{\beta\} \cup J \subset \aA(\pi)$ is a finite collection of charts
such that $M = \bigcup_{\alpha \in J} \uU_\alpha$ and
all coordinate transformations and transition maps relating any
two charts in the collection $\{\beta\} \cup J$ have bounded derivatives of 
all orders (e.g.~it suffices to assume all are precompact).  Then there 
exists a constant $c > 0$ such that 
$$
\| \eta^\beta \|_{W^{k,p}(\mathring{\Omega}_\beta)} \le c
\sum_{\alpha \in J} \| \eta^\alpha \|_{W^{k,p}(\mathring{\Omega}_\alpha)}
$$
for all sections $\eta : M \to E$ with $\eta^\alpha \in
W^{k,p}(\mathring{\Omega}_\alpha)$ for every $\alpha \in J$.
\end{lemma}
\begin{proof}
Choose a partition of unity 
$\{ \rho_\alpha : M \to [0,1] \}_{\alpha \in J}$ subordinate to the
finite open cover $\{ \uU_\alpha \}_{\alpha \in J}$.  Now
$\eta = \sum_{\alpha \in J} \rho_\alpha \eta$, and each $\rho_\alpha \eta$
is supported in~$\uU_\alpha$, so $(\rho_\alpha \eta)^\beta$ has support
in $\Omega_{\beta\alpha} = \varphi_\beta(\uU_\alpha \cap \uU_\beta)$.  Thus
using Theorem~\ref{thm:SobolevCoordinate} with the fact that
$g_{\beta\alpha}$, $\varphi_\beta^{-1}$, $\varphi_{\alpha\beta}$ and
$\varphi_{\beta\alpha} = \varphi_{\alpha\beta}^{-1}$ are all
smooth functions with bounded derivatives of all orders on the domains
in question, we find
\begin{equation*}
\begin{split}
\| \eta^\beta \|_{W^{k,p}(\mathring{\Omega}_\beta)} &=
\left\| \sum_{\alpha \in J} (\rho_\alpha \eta)^\beta \right\|_{W^{k,p}(\mathring{\Omega}_\beta)}
\le \sum_{\alpha \in J} \| (\rho_\alpha \eta)^\beta \|_{W^{k,p}(\mathring{\Omega}_{\beta\alpha})} \\
&= \sum_{\alpha \in J} \| (\rho_\alpha \circ \varphi_\beta^{-1}) 
(g_{\beta\alpha} \circ \varphi_\beta^{-1})(\eta^\alpha \circ \varphi_{\alpha\beta})
\|_{W^{k,p}(\mathring{\Omega}_{\beta\alpha})} \\
&\le c \sum_{\alpha \in J} \| \eta^\alpha \|_{W^{k,p}(\mathring{\Omega}_{\alpha\beta})} \le
c \sum_{\alpha \in J} \| \eta^\alpha \|_{W^{k,p}(\mathring{\Omega}_\alpha)}.
\end{split}
\end{equation*}
\end{proof}

\begin{cor}
If $M$ is compact, then the norm on $W^{k,p}(E)$ given by
Definition~\ref{defn:SobolevMfd} is independent of all auxiliary choices up to
equivalence of norms.   \qed
\end{cor}

\begin{thm}
\label{thm:mfdComplete}
For any smooth vector bundle $\pi : E \to M$ over a compact manifold~$M$,
$W^{k,p}(E)$ is a Banach space.
\end{thm}
\begin{proof}
If $\eta_j \in W^{k,p}(E)$
is a Cauchy sequence, then for some chosen finite collection 
$I \subset \aA(\pi)$ of precompact charts covering~$M$,
the sequences $\eta_j^\alpha$ for $\alpha \in I$
are Cauchy in $W^{k,p}(\mathring{\Omega}_\alpha)$ and thus have limits
$\xi^{(\alpha)} \in W^{k,p}(\mathring{\Omega}_\alpha,\FF^m)$.
Choosing a partition of unity $\{ \rho_\alpha : M \to [0,1] \}_{\alpha \in I}$
subordinate to $\{ \uU_\alpha \}_{\alpha \in I}$, we can now associate to
each $\alpha \in I$ a section
$\eta_{\infty,\alpha} \in W^{k,p}(E)$ characterized uniquely by the condition
that it vanishes outside of $\uU_\alpha$ and is represented in the
trivialization on $\uU_\alpha$ by
$$
\eta_{\infty,\alpha}^\alpha = (\rho_\alpha \circ \varphi_\alpha^{-1})
\xi^{(\alpha)}.
$$
We claim that $\rho_\alpha \eta_j \to \eta_{\infty,\alpha}$ 
in $W^{k,p}(E)$ for each $\alpha \in I$.  Indeed, we have
$$
(\rho_\alpha \eta_j)^\alpha = (\rho_\alpha \circ \varphi_\alpha^{-1})
\eta_j^\alpha \to (\rho_\alpha \circ \varphi_\alpha^{-1}) \xi^{(\alpha)}
= \eta_{\infty,\alpha}^\alpha \quad
\text{ in $W^{k,p}(\mathring{\Omega}_\alpha)$}
$$
since $\eta_j^\alpha \to \xi^{(\alpha)}$.  For all other $\beta \in I$
not equal to~$\alpha$, 
$(\rho_\alpha \eta_j)^\beta - \eta_{\infty,\alpha}^\beta \in
W^{k,p}(\mathring{\Omega}_\beta,\FF^m)$
has support in $\Omega_{\beta\alpha} = \varphi_\beta(\uU_\alpha \cap \uU_\beta)$,
thus
$$
\| (\rho_\alpha \eta_j)^\beta - \eta_{\infty,\alpha}^\beta \|_{W^{k,p}(\mathring{\Omega}_\beta)}
=
\| (\rho_\alpha \eta_j)^\beta - \eta_{\infty,\alpha}^\beta \|_{W^{k,p}(\mathring{\Omega}_{\beta\alpha})}
\le
c \| (\rho_\alpha \eta_j)^\alpha - \eta_{\infty,\alpha}^\alpha \|_{W^{k,p}(\mathring{\Omega}_{\alpha})},
$$
where the inequality comes from Lemma~\ref{lemma:differentCharts} after 
replacing $M$ with $\uU_\alpha$, and $\uU_\beta$ with $\uU_\beta \cap \uU_\alpha$
(note that the lemma does not require $M$ to be compact).  With the claim
established, we have
$$
\eta_j = \sum_{\alpha \in I} \rho_\alpha \eta_j \to \sum_{\alpha \in I}
\eta_{\infty,\alpha} \quad
\text{ in $W^{k,p}(E)$}.
$$
\end{proof}

\begin{exercise}
\label{EX:locEquiv}
For $\uU \subset \RR^n$ an open subset, the space $W^{k,p}_\loc(\uU)$ was
defined in \S\ref{sec:SobolevThms}, but one can give it an alternative
definition in the present context by viewing functions on $\uU$ as sections
of a trivial vector bundle over $\uU$, with the latter viewed as a
noncompact smooth $n$-manifold.  Show that these two definitions
of $W^{k,p}_\loc(\uU)$ are equivalent.
\end{exercise}
\begin{exercise}
\label{EX:compactDomain}
Suppose $\uU \subset \RR^n$ is a bounded open subset with smooth boundary,
so its closure
$\overline{\uU} \subset \RR^n$ is a smooth compact submanifold with boundary, 
and let $E \to \overline{\uU}$ be a trivial vector bundle.  Show that
there is a canonical Banach space isomorphism between
$W^{k,p}(\uU)$ as defined in \S\ref{sec:SobolevThms} and
$W^{k,p}(E)$ as defined in the present section.
\textsl{Hint: Recall that sections in $W^{k,p}(E)$ are only required to be
defined almost everywhere, so in particular if the domain $M$ is a manifold
with boundary, they need not be well defined on~$\p M$.}
\end{exercise}

In light of Exercise~\ref{EX:compactDomain}, the natural
generalization of $W^{k,p}_0(\uU)$ in the present setting is
$$
W^{k,p}_0(E) := \overline{C_0^\infty(E|_{M \setminus \p M})},
$$
i.e.~it is the closure in the $W^{k,p}$-norm of the space of smooth
sections that vanish near the boundary.  Density of smooth sections will
imply that this is the same as $W^{k,p}(E)$ if $M$ is closed, but in
general $W^{k,p}_0(E)$ is a closed subspace of~$W^{k,p}(E)$.

The partition of unity argument in Theorem~\ref{thm:mfdComplete} contains
all the essential ideas needed to generalize results about Sobolev spaces
on domains in $\RR^n$ to compact manifolds.  We now state the essential
results, leaving the proofs as exercises.

\begin{thm}
\label{thm:SobolevMfd}
Assume $M$ is a smooth compact $n$-manifold, possibly with boundary,
$\pi : E \to M$ is a smooth vector bundle of finite rank,
$k \ge 0$ is an integer and $1 \le p < \infty$.  Then the Banach space
$W^{k,p}(E)$ has the following properties.
\begin{enumerate}
\item The space $\Gamma(E)$ of smooth sections is dense
in~$W^{k,p}(E)$.
\item If $N \subset M$ is a smooth compact
$n$-dimensional submanifold with boundary, then there exists a bounded
linear operator
$$
E : W^{k,p}(E|_N) \to W^{k,p}_0(E)
$$
which is an \emph{extension} operator in the sense that
$E\eta|_N = \eta$ for all $\eta \in W^{k,p}(E|_N)$.  Moreover, a section
$\eta \in W^{k,p}(E|_N)$ belongs to $W^{k,p}_0(E|_N)$ if and only if
the section $\tilde{\eta}$ defined to match $\eta$ on $N$ and to vanish
on $M \setminus N$ belongs to $W^{k,p}(E)$.
\item 
If $kp > n$, then for each integer $d \ge 0$, there exists a
continuous and compact inclusion
$$
W^{k+d,p}(E) \hookrightarrow C^d(E).
$$
\item The natural inclusion
$$
W^{k+1,p}(E) \hookrightarrow W^{k,p}(E)
$$
is compact.
\item Suppose $F, G \to M$ are smooth vector bundles such that there exists a
smooth bundle map
$$
E \otimes F \to G : \eta \otimes \xi \mapsto \eta \cdot \xi.
$$
Then if $kp > n$ and $0 \le m \le k$, there exists a continuous product pairing
$$
W^{k,p}(E) \times W^{m,p}(F) \to W^{m,p}(G) : (\eta,\xi) \mapsto \eta \cdot \xi.
$$
In particular, products of $W^{k,p}$ sections give $W^{k,p}$ sections
whenever $kp > n$.
\item Suppose $F \to M$ is another smooth vector bundle,
$\vV \subset E$ is an open subset that intersects every fiber of~$E$,
and we consider the spaces
$$
W^{k,p}(\vV) := \left\{ \eta \in W^{k,p}(E)\ \big|\ {\eta(M)} \subset \vV \right\}
$$
and
$$
C^k_M(\vV,F) := \left\{ \Phi : \vV \to F\ |\ \text{fiber-preserving maps
of class~$C^k$} \right\},
$$
where the latter is assigned the topology of $C^k$-convergence on compact
subsets.  If $kp > n$, then $W^{k,p}(\vV)$ is an open subset of
$W^{k,p}(E)$, and the map
$$
C^k_M(\vV,F) \times W^{k,p}(\vV) \to W^{k,p}(F) : 
(\Phi,\eta) \mapsto \Phi \circ \eta
$$
is well defined and continuous.
\item If $N$ is another smooth compact manifold and $\varphi : N \to M$
is a smooth diffeomorphism, then there is a Banach space isomorphism
$$
W^{k,p}(E) \to W^{k,p}(\varphi^*E) : \eta \mapsto \eta \circ \varphi.
$$
\end{enumerate}
\qed
\end{thm}

\section{Some remarks on domains with cylindrical ends}
\label{sec:SobolevEnds}

For bundles $\pi : E \to M$ with $M$ noncompact, $W^{k,p}(E)$ is not generally
well defined without making additional choices.  When $M = \dot{\Sigma} =
\Sigma \setminus \Gamma$ is a punctured Riemann surface
and $\pi : E \to \dot{\Sigma}$ is
equipped with an asymptotically Hermitian structure
$\{ (E_z,J_z,\omega_z) \}_{z \in \Gamma}$,
one nice way to define $W^{k,p}(E)$ was introduced in
\S\ref{sec:setup}: one takes it to be the space of sections in
$W^{k,p}_\loc(E)$ whose $W^{k,p}$-norms on each cylindrical end are finite
with respect to a choice of asymptotic trivialization.  This definition
requires the convenient fact that complex vector bundles over $S^1$ are
always trivial, though one can also do without this by using
the ideas in the previous section.  Indeed, any collection of local
trivializations on the asymptotic bundle $E_z \to S^1$ covering $S^1$
gives rise via the asymptotically Hermitian structure to a collection of
trivializations on $E$ covering the corresponding cylindrical end~$\dot{\uU}_z$.
The key fact is then that $S^1$ is compact, hence one can always choose such
a covering to be finite: combining this with a finite covering of 
$\dot{\Sigma}$ in the complement of its cylindrical ends by precompact charts,
we obtain a covering of $\dot{\Sigma}$ by 
a finite collection of bundle charts that are not all precompact,
but nonetheless have the property that all transition maps have bounded
derivatives of all orders.  This is enough to define a $W^{k,p}$-norm
for sections of $E \to \dot{\Sigma}$ as in
Definition~\ref{defn:SobolevMfd} and to prove that it does not depend on
the choices of charts or local trivializations, though it does depend on the
asymptotically Hermitian structure.

With this definition understood, one can easily generalize the Sobolev
embedding theorem and other important statements in Theorem~\ref{thm:SobolevMfd}
to the setting of an asymptotically Hermitian bundle over a punctured
Riemann surface.  We shall leave the details of this generalization as an
exercise, but take the opportunity to point out a few important differences
from the compact case.

First, since
$\dot{\Sigma}$ is not compact, neither are the inclusions
$$
W^{k+d,p}(E) \hookrightarrow C^d(E), \qquad
W^{k+1,p}(E) \hookrightarrow W^{k,p}(E).
$$
The proof of compactness fails due to the fact that cylindrical ends require
local trivializations over unbounded domains of the form 
$(0,\infty) \times (0,1) \subset \RR^2$, for which
Theorem~\ref{thm:Rellich} does not hold.  And indeed, considering unbounded
shifts on the infinite cylinder $\dot{\Sigma} = \RR \times S^1$, it is
easy to find a sequence of $W^{k,p}$-bounded functions with $kp > 2$ that
do not have a $C^0$-convergent subsequence.  That is the bad news.

The good news is that if $\eta \in W^{k+d,p}(E)$ for $kp > 2$, then one can
say considerably more about $\eta$ than just that it is $C^d$-smooth.
Indeed, restricting to one of the cylindrical ends $[0,\infty) \times S^1
\subset \dot{\Sigma}$, notice that finiteness of the $W^{k+d,p}$-norm 
over $\dot{\Sigma}$ implies
$$
\| \eta \|_{W^{k+d,p}((R,\infty) \times S^1)} \to 0 \quad \text{ as } \quad
R \to \infty.
$$
Since these domains are all naturally diffeomorphic for different values of
$R$, the $C^d$-norm of $\eta$ over $(R,\infty) \times S^1$ is bounded by
the $W^{k+d,p}$-norm via a constant that does not depend on~$R$, so this
implies an asymptotic decay condition
$$
\| \eta \|_{C^d([R,\infty) \times S^1)} \to 0 \quad \text{ as } \quad
R \to \infty
$$
for every $\eta \in W^{k+d,p}(E)$.

Here is another useful piece of good news: since $\dot{\Sigma}$ does not have 
boundary, $W^{k,p}(E) = W^{k,p}_0(E)$.

\begin{thm}
\label{thm:C0inftyDense}
Given an asymptotically Hermitian bundle $E$ over a punctured Riemann 
surface~$\dot{\Sigma}$, the space $C_0^\infty(E)$ of smooth sections with
compact support is dense in $W^{k,p}(E)$ for all $k \ge 0$ and $1 \le p < \infty$.
\end{thm}
\begin{proof}
We can assume as in Definition~\ref{defn:SobolevMfd} that the
$W^{k,p}$-norm for sections $\eta$ of $E$ is given by
$$
\| \eta \|_{W^{k,p}} = \sum_{\alpha \in I} \| \eta^\alpha \|_{W^{k,p}(\Omega_\alpha)},
$$
where $I \subset \aA(\pi)$ is a finite collection of bundle charts 
$$
\alpha = \left(\varphi_\alpha : \uU_\alpha \stackrel{\cong}{\longrightarrow} \Omega_\alpha \,,\,
\Phi_\alpha : E|_{\uU_\alpha} \stackrel{\cong}{\longrightarrow} \uU_\alpha \times \CC^n \right)
$$
such that each of the open sets $\Omega_\alpha \subset \CC$ is either bounded
or (for charts over the cylindrical ends) of the form $$
\Omega_\alpha = (0,\infty) \times \omega_\alpha \subset \RR^2 = \CC
$$
for some bounded open subset $\omega_\alpha \subset \RR$.  Now given
$\eta \in W^{k,p}(E)$, Theorem~\ref{thm:CinftyDense} provides for
each $\alpha \in I$ a sequence $\eta^\alpha_j \in W^{k,p}(\Omega_\alpha)$
of smooth functions with bounded support such that
$\eta^\alpha_j \to \eta^\alpha$ in $W^{k,p}(\Omega_\alpha)$.  
Choose a partition of unity $\{ \rho_\alpha : \dot{\Sigma} \to [0,1] \}_{\alpha \in I}$
subordinate to the open cover $\{ \uU_\alpha \}_{\alpha \in I}$ and let
$$
\eta_j := \sum_{\alpha \in I} \rho_\alpha (\eta^\alpha_j \circ \varphi_\alpha)
\in W^{k,p}(E).
$$
These sections are smooth and have compact support since the $\eta^\alpha_j$
have bounded support in~$\Omega_\alpha$, and they converge in
$W^{k,p}$ to~$\eta$.
\end{proof}

\chapter{The Floer $C_\varepsilon$ space}
\label{app:Cepsilon}

\minitoc
\vspace{12pt}

The $C_\varepsilon$-topology for functions was introduced by Floer 
\cite{Floer:action} to provide a Banach manifold of perturbed geometric
structures without departing from the smooth category: it is
a way to circumvent the annoying fact that spaces of smooth functions
which arise naturally in geometric settings are not Banach spaces.
The construction of $C_\varepsilon$ spaces generally depends on several
arbitrary choices and is thus far from canonical, but this detail is
unimportant since the $C_\varepsilon$ space itself is never the main object of
interest.  What is important is merely the properties that it has, namely
that it not only embeds continuously into $C^\infty$
and contains an abundance of nontrivial functions, but also
is a separable Banach space
and can therefore be used in the Sard-Smale theorem for genericity arguments.
We shall prove these facts in this appendix.

Fix a smooth finite-rank vector bundle $\pi : E \to M$ over a finite-dimensional
compact manifold $M$, possibly with boundary.  For each integer $k \ge 0$, 
we denote by $C^k(E)$ the Banach space of $C^k$-smooth sections of~$E$;
note that the norm on $C^k(E)$ depends on various auxiliary choices but
is well defined up to equivalence of norms since $M$ is compact.
Now if $(\varepsilon_k)_{k=0}^\infty$ is a sequence of positive numbers with
$\varepsilon_k \to 0$, set
$$
C_\varepsilon(E) = \left\{ \eta \in \Gamma(E)\ \big|\ \|\eta\|_{C_\varepsilon} < \infty \right\},
$$
where the $C_\varepsilon$-norm is defined by
\begin{equation}
\label{eqn:Cepsilon}
\| \eta \|_{C_\varepsilon} = \sum_{k=0}^\infty \varepsilon_k \| \eta \|_{C^k}.
\end{equation}
The norm for $C_\varepsilon(E)$ is somewhat more delicate than for
$C^k(E)$, e.g.~its equivalence class is not obviously
independent of auxiliary choices.  This remark is meant as a sanity check, but 
it should not cause extra concern since, in practice, the space 
$C_\varepsilon(E)$ is typically regarded as an auxiliary choice in itself.
In many applications, one fixes an open subset $\uU \subset M$ and considers
the closed subspace 
$$
C_\varepsilon(E;\uU) = \left\{ \eta \in C_\varepsilon(E)\ \big|\ 
\eta|_{M \setminus \uU} \equiv 0 \right\}.
$$

\begin{remark}
\label{remark:epsNoncompact}
The requirement for $M$ to be compact can be relaxed as long as
$\uU \subset M$ has compact closure: e.g.~in one situation of frequent
interest in this book, we take $M$ to be the noncompact completion of
a symplectic cobordism.  In this case $C_\varepsilon(E;\uU)$ can be
defined as a closed subspace of $C_\varepsilon(E|_{M_0})$ where
$M_0 \subset M$ is any compact manifold with boundary that contains
the closure of~$\uU$.  For this reason, we lose no generality in continuing
under the assumption that $M$ is compact.
\end{remark}

In order to prove things about $C_\varepsilon(E)$, we will need
to specify a more precise definition of the $C^k$-norms.  To this end,
define a sequence of vector bundles $E^{(k)} \to M$ for integers $k \ge 0$
inductively by
$$
E^{(0)} := E, \qquad E^{(k+1)} := \Hom(TM,E^{(k)}).
$$
Choose connections and bundle metrics on both $TM$ and~$E$; these induce
connections and bundle metrics on each of the $E^{(k)}$, so that for any
section $\xi \in \Gamma(E^{(k)})$, the covariant derivative $\nabla \xi$ is
now a section of $E^{(k+1)}$.  In particular for $\eta \in \Gamma(E)$, we can
define the ``$k$th covariant derivative'' of $\eta$ as a section
$$
\nabla^k \eta \in \Gamma(E^{(k)}).
$$
Using the bundle metrics to define $C^0$-norms for sections of $E^{(k)}$,
we can then define
$$
\| \eta \|_{C^k(E)} = \sum_{m=0}^k \| \nabla^m \eta \|_{C^0(E^{(m)})},
$$
where by convention $\nabla^0 \eta := \eta$.  We will assume throughout
the following that the $C^k$-norms appearing in \eqref{eqn:Cepsilon}
are defined in this way.

\begin{thm}
\label{thm:CepsilonComplete}
$C_\varepsilon(E)$ is a Banach space.
\end{thm}
\begin{proof}
We need to show that $C_\varepsilon$-Cauchy sequences converge in the
$C_\varepsilon$-norm.  It is clear from the definitions that if
$\eta_j \in C_\varepsilon(E)$ is Cauchy, then $\eta_j$ is also
$C^k$-Cauchy for every $k \ge 0$, hence its derivatives
$\nabla^k \eta_j$ for every $k$ are $C^0$-convergent to continuous sections
$\xi^k$ of~$E^{(k)}$.  This convergence implies that $\xi^{k+1} = \nabla \xi^k$
in the sense of distributions, hence by the equivalence of classical and
distributional derivatives (see e.g.~\cite{LiebLoss}*{\S 6.10}),
$\eta_\infty := \xi^0$ is smooth with $\nabla^k \eta_\infty = \xi^k$, so that
$\nabla^k \eta_j \to \nabla^k \eta_\infty$ in $C^0(E^{(k)})$
for all~$k$.

We claim $\eta_\infty \in C_\varepsilon(E)$.  Choose $N > 0$ such that
$\| \eta_i - \eta_j \|_{C_\varepsilon} < 1$ for all $i,j \ge N$.  Then for
every $m \in \NN$ and every $i \ge N$,
\begin{equation*}
\begin{split}
\sum_{k=0}^m \varepsilon_k \| \eta_i \|_{C^k} &\le
\sum_{k=0}^m \varepsilon_k \| \eta_i - \eta_N \|_{C^k} + \sum_{k=0}^m
\varepsilon_k \| \eta_N \|_{C^k} \\
&\le \| \eta_i - \eta_N \|_{C_\varepsilon} + \| \eta_N \|_{C_\varepsilon}
< 1 + \| \eta_N \|_{C_\varepsilon}.
\end{split}
\end{equation*}
Fixing $m$ and letting $i \to \infty$, we then have
$$
\sum_{k=0}^m \varepsilon_k \| \eta_\infty \|_{C^k} \le 1 + \| \eta_N \|_{C_\varepsilon}
$$
for all $m$, so we can now let $m \to \infty$ and conclude
$\| \eta_\infty \|_{C_\varepsilon} \le 1 + \| \eta_N \|_{C_\varepsilon}< \infty$.

The argument that $\|\eta_j \to \eta_\infty\|_{C_\varepsilon} \to 0$ as
$j \to \infty$ is similar: pick $\epsilon > 0$ and $N$ such that
$\| \eta_i - \eta_j \|_{C_\varepsilon} < \epsilon$ for all $i,j \ge N$.
Then for a fixed $m \in \NN$, we can let $i \to \infty$ in the expression
$\sum_{k=0}^m \varepsilon_k \| \eta_i - \eta_j \|_{C^k} < \epsilon$, giving
$$
\sum_{k=0}^m \varepsilon_k \| \eta_\infty - \eta_j \|_{C^k} \le \epsilon.
$$
This is true for every~$m$, so we can take $m \to \infty$ and conclude
$\| \eta_\infty - \eta_j \|_{C_\varepsilon} \le \epsilon$ for all $j \ge N$.
\end{proof}

To show that $C_\varepsilon(E)$ is also separable, we will follow a
hint\footnote{Thanks to Sam Lisi for explaining to me what the hint
in \cite{HoferSalamon:Novikov} was referring to.}
from \cite{HoferSalamon:Novikov} and embed it isometrically into another 
Banach space that can be more easily shown to be separable.
For each integer $k \ge 0$, define the vector bundle
$$
F^{(k)} = E^{(0)} \oplus \ldots \oplus E^{(k)},
$$
and let $X_\varepsilon$ denote the vector space of all sequences
$$
\xi := (\xi^0,\xi^1,\xi^2,\ldots) \in \prod_{k=0}^\infty C^0(F^{(k)})
$$
such that
$$
\| \xi \|_{X_\varepsilon} := \sum_{k=0}^\infty \varepsilon_k \| \xi^k \|_{C^0} < \infty.
$$

\begin{exercise}
\label{EX:Xinfty}
Adapt the proof of Theorem~\ref{thm:CepsilonComplete} to show that 
$X_\varepsilon$ is also a Banach space.
\end{exercise}

\begin{lemma}
\label{lemma:XepsSeparable}
$X_\varepsilon$ is separable.
\end{lemma}
\begin{proof}
Since $C^0(F^{(k)})$ is
separable for each $k \ge 0$, we can fix countable dense subsets
$P^k \subset C^0(F^{(k)})$.  The set
$$
P := \left\{ (\xi^0,\ldots,\xi^N,0,0,\ldots) \in X_\varepsilon\ \big|\ 
\text{$N \ge 0$ and $\xi^k \in P^k$ for all $k=0,\ldots,N$} \right\}
$$
is then countable and dense in $X_\varepsilon$.
\end{proof}

\begin{thm}
\label{thm:CepsilonSeparable}
$C_\varepsilon(E)$ is separable.
\end{thm}
\begin{proof}
Consider the injective linear map
$$
C_\varepsilon(E) \hookrightarrow X_\varepsilon :
\eta \mapsto \left( \eta,(\eta,\nabla \eta),(\eta,\nabla\eta,\nabla^2\eta),\ldots\right).
$$
This is an isometric embedding and thus presents $C_\varepsilon(E)$
as a closed linear subspace of $X_\varepsilon$, hence the theorem follows
from Lemma~\ref{lemma:XepsSeparable} and the fact that subspaces of separable
metric spaces are always separable.
\end{proof}

Note that given any open subset $\uU \subset M$, Theorems~\ref{thm:CepsilonComplete}
and~\ref{thm:CepsilonSeparable} also hold for $C_\varepsilon(E;\uU)$, as a
closed subspace of~$C_\varepsilon(E)$.  So far in this discussion, however,
there has been no guarantee that $C_\varepsilon(E)$ or $C_\varepsilon(E;\uU)$
contains anything other than the zero-section, though it is clear that 
in theory, one should always be able to enlarge the space by choosing new
sequences $\varepsilon_k$ that converge to zero faster.  The following
result says that $C_\varepsilon(E;\uU)$ can always be made large enough to be
useful in applications.

\begin{thm}
\label{thm:bump}
Given an open subset $\uU \subset M$, the sequence $\varepsilon_k$ can be 
chosen to have the following properties:
\begin{enumerate}
\item $C_\varepsilon(E;\uU)$ is dense in the space of continuous
sections vanishing outside~$\uU$.
\item
Given any point $p \in \uU$, a neighborhood $\nN_p \subset \uU$ of~$p$,
a number $\delta > 0$ and a continuous section $\eta_0$ of~$E$, there exists
a section $\eta \in \Gamma(E)$ and a smooth compactly supported
function $\beta : \nN_p \to [0,1]$ such that
$$
\beta \eta \in C_\varepsilon(E;\uU), \qquad
\beta(p)\eta(p) = \eta_0(p),\quad \text{ and } \quad
\| \eta - \eta_0 \|_{C^0} < \delta.
$$
\end{enumerate}
\end{thm}
\begin{proof}
Note first that it suffices to find two separate sequences
$\varepsilon_k$ and $\varepsilon_k'$ that
have the first and second
property respectively, as the sequence of minima $\min(\varepsilon_k,\varepsilon_k')$
will then have both properties.

The following construction for the first property is based on a suggestion by
Barney Bramham.  Observe first that the space $C^0(E;\uU)$ of continuous
sections vanishing outside $\uU$ is a closed subspace of $C^0(E)$ and is
thus separable, so we can choose a countable $C^0$-dense subset 
$P \subset C^0(E;\uU)$.  Moreover, the space of \emph{smooth} sections
vanishing outside $\uU$ is dense in $C^0(E;\uU)$, hence we can assume
without loss of generality that the sections in $P$ are smooth.
Now write $P = \{\eta_1,\eta_2,\eta_3,\ldots\}$ and
define $\varepsilon_k > 0$ for every integer $k \ge 0$ to have the property
$$
\varepsilon_k < \frac{1}{2^k} \min\left\{ \frac{1}{\| \eta_1 \|_{C^k}} ,
\ldots , \frac{1}{\| \eta_k \|_{C^k}} \right\}.
$$
Then every $\eta_j$ is in $C_\varepsilon(E;\uU)$, as
$$
\| \eta_j \|_{C_\varepsilon} < \sum_{k=0}^{j-1} \varepsilon_k \| \eta_j \|_{C^k}
+ \sum_{k=j}^\infty \frac{1}{2^k} < \infty.
$$

The second property is essentially local, so it can be deduced from
Lemma~\ref{lemma:bump} below.
\end{proof}

\begin{lemma}
\label{lemma:bump}
Suppose $\beta : \intDD^n \to [0,1]$ is a smooth function with compact support
on the open unit ball $\intDD^n \subset \RR^n$ and $\beta(0)=1$.
One can choose a sequence of positive numbers $\varepsilon_k \to 0$ such that 
for every $\eta_0 \in \RR^m$ and $r > 0$, the function $\eta : \RR^n \to \RR^m$
defined by
$$
\eta(p) := \beta(p/r) \eta_0
$$
satisfies $\sum_{k=0}^\infty \varepsilon_k \| \eta \|_{C^k} < \infty$.
\end{lemma}
\begin{proof}
Define $\varepsilon_k > 0$ so that for $k \ge 1$,
$$
\varepsilon_k = \frac{1}{k^k \| \beta \|_{C^k}}.
$$
Then
$$
\sum_{k=1}^\infty \varepsilon_k \| \eta \|_{C^k} \le
\sum_{k=1}^\infty \frac{1}{k^k \| \beta \|_{C^k}} \frac{\| \beta \|_{C^k}}{r^k}
= \sum_{k=1}^\infty \left( \frac{1 / r}{k} \right)^k < \infty.
$$
\end{proof}

\chapter{Genericity in the space of asymptotic operators}
\label{app:spectral}

\minitoc
\vspace{12pt}

The purpose of this appendix is to prove Lemma~\ref{lemma:generic}, which was
needed for our definition of spectral flow in~\S\ref{sec:spectral}.  
The proof combines some ideas from that section
with the technique used in Lecture~\ref{lec:transversality} to prove
generic transversality of moduli spaces via the Sard-Smale theorem.
Some knowledge of that technique should thus be considered a prerequisite
for this appendix; if you have never seen it before and were directed
here after reading the statement of Lemma~\ref{lemma:generic}, you might
want to skip this for now and come back after you've read as far as
Lecture~\ref{lec:transversality}.

Recalling the notation from Lecture~\ref{lec:asymptotic}, we fix the 
real Hilbert spaces
$$
\hH = L^2(S^1,\RR^{2n}), \qquad \dD = H^1(S^1,\RR^{2n}),
$$
the symmetric index~$0$ Fredholm operator
$$
\Tref = -J_0\, \p_t : \dD \to \hH
$$
and, given a smooth loop of symmetric matrices $S : S^1 \to \End_\RR^\sym(\RR^{2n})$,
refer to any operator of the form
$$
\mathbf{A} = -J_0\, \p_t - S : \dD \to \hH
$$
as an \defin{asymptotic operator}.  Such operators belong to the space of
symmetric compact perturbations of~$\Tref$,
$$
\Fred_\RR^\sym(\dD,\hH,\Tref) = \left\{ \Tref + \mathbf{K} : \dD \to \hH \ 
\big|\ \mathbf{K} \in \Lin_\RR^\sym(\hH) \right\},
$$
which we regard as a smooth Banach manifold via its obvious identification
with the space $\Lin_\RR^\sym(\hH)$ of symmetric bounded linear operators
on~$\hH$.  For $k \in \NN$, we denote by
$$
\Fred_\RR^{\sym,k}(\dD,\hH,\Tref) \subset \Fred_\RR^\sym(\dD,\hH,\Tref)
$$
the finite-codimensional submanifold determined by the condition
$\dim_\RR \ker \mathbf{A} = \dim_\RR \coker \mathbf{A} = k$.

Here is the statement of Lemma~\ref{lemma:generic} again.

\begin{lemmu}
Fix a smooth map $S : [-1,1] \times S^1 \to \End_\RR^\sym(\RR^{2n})$ and
consider the $1$-parameter family of operators
$$
\mathbf{A}_s := -J_0 \, \p_t - S(s,\cdot) 
\in \Fred_\RR^\sym(\dD,\hH,\Tref)
$$
for $s \in [-1,1]$.  Then after a $C^\infty$-small perturbation of 
$S$ fixed at $s= \pm 1$, one can assume the following:
\begin{enumerate}
\item For every $s \in (-1,1)$, all
eigenvalues of $\mathbf{A}_s$ (regarded as an unbounded operator on~$\hH$)
are simple.
\item All intersections of the path
$$
(-1,1) \to \Fred_\RR^\sym(\dD,\hH,\Tref) : s \mapsto \mathbf{A}_s
$$
with $\Fred_\RR^{\sym,1}(\dD,\hH,\Tref)$ are transverse.
\end{enumerate}
\end{lemmu}

We shall now prove this by constructing a Floer-type space of
$C_\varepsilon$-smooth (see Appendix~\ref{app:Cepsilon}) perturbed families 
of asymptotic operators, and using the Sard-Smale theorem to find a countable
collection of comeager subsets whose intersection contains perturbations
achieving the desired conditions.

Choose a sequence of positive numbers 
$(\varepsilon)_{k=0}^\infty$ with
$\varepsilon_k \to 0$ to define a separable Banach space
$$
\aA_\varepsilon := \left\{ B \in C^\infty([-1,1] \times S^1,\End_\RR^\sym(\RR^{2n}))\ \big|\ 
\| B \|_{C_\varepsilon} < \infty \text{ and } B(\pm 1,\cdot) \equiv 0 \right\},
$$
and assume via Theorem~\ref{thm:bump} that $\aA_\varepsilon$ is dense in the
Banach space of continuous functions $[-1,1] \times S^1 \to \End_\RR^\sym(\RR^{2n})$
vanishing at $\{\pm 1\} \times S^1$.
We then consider perturbed $1$-parameter families of asymptotic operators of the form
$$
\mathbf{A}_s^B := \mathbf{A}_s + B(s,\cdot) : \dD \to \hH
$$ 
for $B \in \aA_\varepsilon$, $s \in [-1,1]$.  For each $k \in \NN$ and
$B \in \aA_\varepsilon$, define the set
$$
\vV^k(B) = \left\{ (s,\lambda) \in (-1,1) \times \RR\ \big|\ 
\dim_\RR \ker \left( \mathbf{A}_s^B - \lambda \right) = k \right\}.
$$
To show that eigenvalues are generically simple,
we need to show that for a comeager set of choices of $B \in \aA_\varepsilon$,
$\vV^k(B)$ is empty for all $k \ge 2$.
Given $(s_0,\lambda_0) \in \vV^k(B)$, recall from \S\ref{sec:spectral}
that there exist decompositions
$$
\dD = V \oplus K, \qquad \hH = W \oplus K
$$
where $K = \ker \left( \mathbf{A}_{s_0}^B - \lambda_0 \right)$,
$W = \im \left( \mathbf{A}_{s_0}^B - \lambda_0 \right)$ is the $L^2$-orthogonal
complement of~$K$, and
$V = W \cap \dD$, so that any symmetric bounded linear operator
$\mathbf{T}$ in a sufficiently small neighborhood $\oO \subset \Lin_\RR^\sym(\dD,\hH)$ of
$\mathbf{A}_{s_0}^B - \lambda_0$ can be written in block form
$$
\mathbf{T} = \begin{pmatrix} \mathbf{A} & \mathbf{B} \\ \mathbf{C} & \mathbf{D} \end{pmatrix}
$$
with $\mathbf{A} : V \to W$ invertible, giving rise to a smooth map
$$
\Phi : \oO \to \End_\RR^\sym(K) : \mathbf{T} \mapsto \mathbf{D} - 
\mathbf{C} \mathbf{A}^{-1} \mathbf{B}
$$
whose zero-set is precisely the set of nearby symmetric operators with
$k$-dimensional kernel.  A neighborhood of $(s_0,\lambda_0)$ 
in $\vV^k(B)$ can thus be identified with the zero-set of the map
$$
\Psi_B(s,\lambda) := \Phi(\mathbf{A}_s^B - \lambda) \in \End_\RR^\sym(K),
$$
defined for $(s,\lambda) \in (-1,1) \times \RR$ sufficiently close to 
$(s_0,\lambda_0)$.  Notice that the derivative 
$d\Psi_B(s,\lambda) : \RR \oplus \RR \to \End_\RR^\sym(K)$ is
Fredholm since its domain and target are both finite dimensional, and
it can only ever be surjective when $k = \dim_\RR K = 1$.

The following space will now play the role of a ``universal moduli space''
as in Lecture~\ref{lec:transversality}: let
$$
\vV^k = \left\{ (s,\lambda,B) \in (-1,1) \times \RR \times
\aA_\varepsilon \ \big|\ (s,\lambda) \in \vV^k(B) \right\}.
$$
The proof that this is a smooth Banach manifold depends on the following
algebraic lemma.

\begin{lemma}
\label{lemma:symAlg}
Fix an asymptotic operator $\mathbf{A} = -J_0\, \p_t - S$ and a linear
transformation
$$
\Upsilon : \ker\mathbf{A} \to \ker\mathbf{A}
$$
that is symmetric with respect to the $L^2$-product.  Then there exists a
smooth loop $B : S^1 \to \End^\sym(\RR^{2n})$ such that
$$
\langle \eta , B \xi \rangle_{L^2} = \langle \eta , \Upsilon \xi \rangle_{L^2}
$$
for all $\eta , \xi \in \ker \mathbf{A}$.
\end{lemma}
\begin{proof}
Note first that every nontrivial loop $\eta \in \ker \mathbf{A} \subset
H^1(S^1,\RR^{2n})$ is smooth and nowhere zero since it satisfies a linear
first-order ODE with smooth coefficients.  It follows that if we fix a basis
$(\eta_1,\ldots,\eta_k)$ for $\ker \mathbf{A}$, then the
vectors $\eta_1(t),\ldots,\eta_k(t) \in \RR^{2n}$ are also linearly independent
for all $t \in S^1$ and thus span a smooth $S^1$-family of $k$-dimensional
subspaces $V_t \subset \RR^{2n}$, each equipped with a distinguished basis.
It follows that there exists a unique smooth $S^1$-family of linear transformations
$\widehat{B}(t) : V_t \to V_t$ such that for every $\eta \in \ker\mathbf{A}$,
$\widehat{B}(t) \eta(t) = (\Upsilon \eta)(t)$ for all~$t$.  Extend $\widehat{B}(t)$
arbitrarily to a smooth family of linear maps on $\RR^{2n}$.

The matrices $\widehat{B}(t) \in \End_\RR(\RR^{2n})$ need not be symmetric,
but they do satisfy
$$
\langle \eta , \widehat{B} \xi \rangle_{L^2} = \langle \eta , \Upsilon \xi \rangle_{L^2}
\quad \text{ for all $\eta,\xi \in \ker \mathbf{A}$}.
$$
Since $\Upsilon$ is symmetric, this implies moreover that for all
$\eta,\xi \in \ker\mathbf{A}$,
$$
\langle \eta,\Upsilon \xi \rangle_{L^2} = \langle \xi,\Upsilon \eta \rangle_{L^2}
= \langle \xi , \widehat{B} \eta \rangle_{L^2}
= \langle \eta , \widehat{B}^\transpose \xi \rangle_{L^2}.
$$
The loop $B := \frac{1}{2}( \widehat{B} + \widehat{B}^\transpose )$ thus has
the desired properties.
\end{proof}

Now using the previously described construction in the space of symmetric
Fredholm operators, a neighborhood of any point
$(s_0,\lambda_0,B_0)$ in $\vV^k$ can be identified with the zero-set of a
smooth map of the form
$$
\Psi(s,\lambda,B) := \Psi_{B}(s,\lambda) \in \End_\FF^\sym(K),
$$
defined for all $(s,\lambda,B)$ sufficiently close to $(s_0,\lambda_0,B_0)$
in $(-1,1) \times \RR \times \aA_\varepsilon$, where
$K = \ker \left( \mathbf{A}_{s_0}^{B_0} - \lambda_0 \right)$.
The partial derivative of $\Psi$ with respect to the third variable
at $(s_0,\lambda_0,B_0)$ is then a linear map
$$
\mathbf{L} := D_3 \Psi(s_0,\lambda_0,B_0) : \aA_\varepsilon \to \End_\RR^\sym(K)
$$
of the form
\begin{equation}
\label{eqn:defL}
\mathbf{L}B : K \to K : \eta \mapsto \pi_K(B(s_0,\cdot)\eta),
\end{equation}
where $\pi_K : W \oplus K \to K$ is the orthogonal projection.  We claim that
$\mathbf{L}$ is surjective.  Indeed, for any $\Upsilon \in \End_\RR^\sym(K)$,
Lemma~\ref{lemma:symAlg} provides a smooth loop $C_0 : S^1 \to \End_\RR^\sym(K)$
such that
$$
\pi_K(C_0\eta) = \Upsilon \eta \quad \text{ for all } \eta \in K,
$$
and this can be extended to a smooth function $C : [-1,1] \times S^1 \to \End_\RR^\sym(K)$
satisfying $C(s_0,\cdot) \equiv C_0$ and $C(\pm 1,\cdot) \equiv 0$ since
$s_0 \ne \pm 1$.  The function
$C$ might fail to be of class $C_\varepsilon$, but since it can be approximated
arbitrarily well in the $C^0$-norm by functions in~$\aA_\varepsilon$, we
conclude that the image of $\mathbf{L}$ is dense in $\End_\RR^\sym(K)$.
Since the latter is finite dimensional, the claim follows.

The implicit function theorem now gives $\vV^k$ the structure of a smooth Banach
submanifold of $(-1,1) \times \RR \times \aA_\varepsilon$, and it is
separable since the latter is also separable.  Consider the projection
\begin{equation}
\label{eqn:spectralProj}
\pi : \vV^k \to \aA_\varepsilon : (s,\lambda,B) \mapsto B,
\end{equation}
which is a smooth map of separable Banach manifolds whose fibers
$\pi^{-1}(B)$ are the spaces $\vV^k(B)$.
Using Lemma~\ref{lemma:linearAlg}, the fact that each map $\Psi_{B}$
is Fredholm implies that $\pi$ is also a Fredholm map, so
the Sard-Smale theorem implies that the regular values of $\pi$ form a comeager
subset 
$$
\aA_\varepsilon^{\reg,k} \subset \aA_\varepsilon.
$$
The intersection
$$
\aA_\varepsilon^\reg := \bigcap_{k \in \NN} \aA_\varepsilon^{\reg,k}
$$
is then another comeager subset of~$\aA_\varepsilon$, with the property that
for each $B \in \aA_\varepsilon^\reg$ and every $(s,\lambda) \in \vV^k(B)$,
$d\Psi_B(s,\lambda)$ is (by Lemma~\ref{lemma:linearAlg}) surjective.
As was observed previously, this is impossible for dimensional reasons 
if $k \ge 2$, implying that $\vV^k(B)$ is then empty.

To find perturbations that also achieve the transversality condition, we
use a similar argument: define for each $B \in \aA_\varepsilon$ the subset
$$
\vV^0(B) = \left\{ s \in (-1,1) \ \big|\ 
\dim_\RR \ker \mathbf{A}_s^B = 1 \right\},
$$
along with the corresponding universal set
$$
\vV^0 = \left\{ (s,B) \in (-1,1) \times \aA_\varepsilon \ \big|\ 
s \in \vV^0(B) \right\}.
$$
A neighborhood of any $(s_0,B_0)$ in $\vV^0$ is then the
zero-set of a smooth map of the form
$$
\Psi(s,B) = \Phi(\mathbf{A}_s^B) \in \End_\RR^\sym(\ker \mathbf{A}_{s_0}^{B_0}),
$$
defined for all $(s,B) \in (-1,1) \times \aA_\varepsilon$ close enough to
$(s_0,B_0)$.  For a fixed $B \in \aA_\varepsilon$ near~$B_0$ and
$s_1 \in \vV^0(B)$ near~$s_0$, a neighborhood of $s_1$ in $\vV^0(B)$ is then 
the zero-set of $\Psi_B(s) := \Psi(s,B)$,
and the intersection of the path $s \mapsto \mathbf{A}_s \in
\Fred_\RR^\sym(\dD,\hH,\Tref)$ with $\Fred_\RR^{\sym,1}(\dD,\hH,\Tref)$
at $s=s_1$ is transverse if and only if 
$$
d\Psi_B(s_1) : \RR \to \End_\RR^\sym(\ker \mathbf{A}_{s_0}^{B_0})
$$
is surjective.  At $(s_0,B_0)$, the partial derivative of
$\Psi$ with respect to $B$ is again the same operator
$$
\mathbf{L} = D_2 \Psi(s_0,B_0) : \aA_\varepsilon \to \End_\RR^\sym(\ker \mathbf{A}_{s_0}^{B_0})
$$
as in \eqref{eqn:defL},
which we've already seen is surjective due to Lemma~\ref{lemma:symAlg}.
Thus one can apply the Sard-Smale theorem to the projection
$$
\vV^0 \to \aA_\varepsilon : (s,B) \mapsto B,
$$
obtaining a comeager subset $\aA_\varepsilon^{\reg,0} \subset \aA_\varepsilon$
such that all paths $\mathbf{A}_s + B(s,\cdot)$ for $B \in \aA_\varepsilon^{\reg,0}$
satisfy the required transversality condition.  The comeager subset
$\aA_\varepsilon^{\reg,0} \cap \aA_\varepsilon^\reg \subset \aA_\varepsilon$ 
thus consists of perturbed families of operators for which all desired
conditions are satisfied, and it contains a sequence converging in
the $C^\infty$-topology to~$0$.  This concludes the proof of
Lemma~\ref{lemma:generic}.

\backmatter

\end{document}